\documentclass[numbers=enddot,12pt,final,onecolumn,notitlepage]{scrartcl}%
\usepackage[headsepline,footsepline,manualmark]{scrlayer-scrpage}
\usepackage[all,cmtip]{xy}
\usepackage{amssymb}
\usepackage{amsmath}
\usepackage{amsthm}
\usepackage{framed}
\usepackage{comment}
\usepackage{cancel}
\usepackage{color}
\usepackage[breaklinks=True]{hyperref}
\usepackage[sc]{mathpazo}
\usepackage[T1]{fontenc}
\usepackage{needspace}
\usepackage{tabls}
\usepackage{tikz}
\usepackage[type={CC}, modifier={zero}, version={1.0},]{doclicense}
\providecommand{\U}[1]{\protect\rule{.1in}{.1in}}
\newcounter{exer}

\numberwithin{exer}{subsection}
\theoremstyle{definition}
\newtheorem{theo}{Theorem}[subsection]
\newenvironment{theorem}[1][]
{\begin{theo}[#1]\begin{leftbar}}
{\end{leftbar}\end{theo}}
\newtheorem{lem}[theo]{Lemma}
\newenvironment{lemma}[1][]
{\begin{lem}[#1]\begin{leftbar}}
{\end{leftbar}\end{lem}}
\newtheorem{prop}[theo]{Proposition}
\newenvironment{proposition}[1][]
{\begin{prop}[#1]\begin{leftbar}}
{\end{leftbar}\end{prop}}
\newtheorem{defi}[theo]{Definition}
\newenvironment{definition}[1][]
{\begin{defi}[#1]\begin{leftbar}}
{\end{leftbar}\end{defi}}
\newtheorem{remk}[theo]{Remark}
\newenvironment{remark}[1][]
{\begin{remk}[#1]\begin{leftbar}}
{\end{leftbar}\end{remk}}
\newtheorem{coro}[theo]{Corollary}
\newenvironment{corollary}[1][]
{\begin{coro}[#1]\begin{leftbar}}
{\end{leftbar}\end{coro}}
\newtheorem{conv}[theo]{Convention}
\newenvironment{convention}[1][]
{\begin{conv}[#1]\begin{leftbar}}
{\end{leftbar}\end{conv}}
\newtheorem{quest}[theo]{Question}
\newenvironment{question}[1][]
{\begin{quest}[#1]\begin{leftbar}}
{\end{leftbar}\end{quest}}
\newtheorem{warn}[theo]{Warning}
\newenvironment{warning}[1][]
{\begin{warn}[#1]\begin{leftbar}}
{\end{leftbar}\end{warn}}
\newtheorem{conj}[theo]{Conjecture}

\newtheorem{exam}[theo]{Example}
\newenvironment{example}[1][]
{\begin{exam}[#1]\begin{leftbar}}
{\end{leftbar}\end{exam}}
\newtheorem{exmp}[exer]{Exercise}
\newenvironment{exercise}[1][]
{\begin{exmp}[#1]\begin{leftbar}}
{\end{leftbar}\end{exmp}}
\newtheorem{faketh}[theo]{Fake Theorem}
\newenvironment{faketheorem}[1][]
{\begin{faketh}[#1]\begin{leftbar}}
{\end{leftbar}\end{faketh}}
\newenvironment{statement}{\begin{quote}}{\end{quote}}
\newenvironment{fineprint}{\begin{small}}{\end{small}}
\let\sumnonlimits\sum
\let\prodnonlimits\prod
\let\cupnonlimits\bigcup
\let\capnonlimits\bigcap
\renewcommand{\sum}{\sumnonlimits\limits}
\renewcommand{\prod}{\prodnonlimits\limits}
\renewcommand{\bigcup}{\cupnonlimits\limits}
\renewcommand{\bigcap}{\capnonlimits\limits}
\newenvironment{noncompile}{}{}
\excludecomment{verlong}
\includecomment{vershort}
\excludecomment{noncompile}

\newcommand{\set}[1]{\left\{ #1 \right\}}

\newcommand{\arinj}{\ar@{_{(}->}}
\newcommand{\arinjrev}{\ar@{^{(}->}}
\newcommand{\arsurj}{\ar@{->>}}
\newcommand{\arelem}{\ar@{|->}}
\newcommand{\arback}{\ar@{<-}}
\newcommand{\symd}{\mathbin{\bigtriangleup}}

\usetikzlibrary{arrows.meta}
\usetikzlibrary{calc}
\usetikzlibrary{chains}
\usetikzlibrary{shapes}
\usetikzlibrary{decorations.pathmorphing}
\usetikzlibrary{lindenmayersystems}
\definecolor{darkgreen}{rgb}{0,.5,0}
\newtheoremstyle{plainsl}
{8pt plus 2pt minus 4pt}
{8pt plus 2pt minus 4pt}
{\slshape}
{0pt}
{\bfseries}
{.}
{5pt plus 1pt minus 1pt}
{}
\theoremstyle{plainsl}
\ihead{Introduction to Discrete Mathematics, version November 27, 2026}
\ohead{page \thepage}
\begin{document}

\title{An Introduction to Discrete Mathematics\\{\normalsize (Text for Math 221 Winter 2024 at Drexel University)}}
\author{Darij Grinberg}
\date{November 27, 2026}
\maketitle

\begin{abstract}
\textbf{Abstract.} This is a rigorous undergraduate-level introduction to
discrete mathematics, covering mathematical induction, finite sums and
products, elementary number theory, maps and basic enumerative combinatorics.
A familiarity with logic and proofs is assumed of the reader.

The text is written for a quarter-long course; it will likely be extended with
further topics (e.g., equivalence relations, games) as Drexel University
introduces the semester schedule.

Over 90 exercises (without solutions) are scattered through the text.

\end{abstract}
\tableofcontents

\doclicenseThis

\newpage


This is a set of lecture notes for my Math 221 (Discrete Mathematics) course
at Drexel University in the Winter quarters of 2024 and 2025 and in the Spring
quarter of 2025.
It will likely be extended when Drexel switches to a semester schedule (adding
chapters on equivalence relations and combinatorial game theory).

\setcounter{section}{-1}

\section{Preface}

\subsection{What is this?}

This is a course on \textbf{discrete mathematics}. To us, discrete mathematics
means the mathematics of finite, discrete objects: integers, finite sets,
occasionally some more complex creatures such as graphs and polynomials.
Integer sequences, while theoretically infinite, are also included since one
usually makes statements about finite pieces of the sequence. Much of linear
algebra logically belongs to discrete mathematics, but there are separate
courses entirely devoted to it, so we won't touch on it here.

Discrete mathematics stands in contrast to \textbf{continuous mathematics},
which studies real numbers, continuous functions and infinite sets; the latter
mostly begins with analysis (or calculus, which is its less rigorous variant).
The two types of mathematics are intertwined and often borrow ideas and
results from each other; however, at the elementary level of this course, we
will encounter almost nothing from the continuous side.

So this course will introduce you to some of the major topics of discrete mathematics:

\begin{itemize}
\item \textbf{mathematical induction} and \textbf{recursion};

\item \textbf{elementary number theory} (the properties of divisibility, prime
numbers, coprimality, possibly applications like the RSA cryptosystem);

\item basic \textbf{enumerative combinatorics} (counting and binomial coefficients).


\end{itemize}

We will neither go very deep nor be fully rigorous about everything. There are
deeper, more specific classes on most of these subjects:

\begin{itemize}
\item Math 222 (notes: \cite{19fco}, \cite{22fco}) is a quarter-length
introduction to enumerative combinatorics.

\item Math 530 (notes: \cite{23s}) is an introduction to graph theory.

\item CS 303 is a course on cryptography.

\item Math 235 (notes: \cite{20f} and \cite{23f}) is an introduction to
mathematical problem-solving. This can be viewed as a continuation of this
course, leading into more advanced techniques and more exotic results.
\end{itemize}

The reader is assumed to have learned the concept of a mathematical proof, the
language of sets, and fundamental logical rules and techniques (such as proof
by contradiction). This material can be found in various texts, such as
\cite{BecGeo10}, \cite[Chapters 1--5]{LeLeMe16}, \cite{Day-proofs},
\cite{DonPan25}, \cite{Hammac15}, \cite[Chapters 1--2]{KeeGui23}, \cite[Part I
and Appendices A--B]{Newste23} or \cite[Chapters 1--3]{Swanso20}.

I do not intend to give a Bourbaki-style axiomatic treatment of the subject
here, nor to spell out each proof in maximum possible rigor (though I am more
rigorous than many other undergraduate texts). My goal is merely to give a
taste of each of several important and (if I dare say so) interesting topics,
each time veering deep enough to see some substance but not to get lost in the
jungle. Other introductions to discrete mathematics are \cite{Levin21},
\cite{LeLeMe16}, \cite{Newste23} and \cite{GKP}, just to name a few. There is
no \textquotedblleft standard\textquotedblright\ choice of material for such a
text; each author goes their own way through the vast landscape. So do I in
these notes. I have deliberately avoided anything analytic or geometric in
order to stick to the subject declared (\textbf{discrete} mathematics!), but
otherwise I have picked from different topics and fields. Some topics such as
graphs, posets and the construction of the number systems are nevertheless
missing from this introduction, as the lack of days in an academic quarter has
forced choices upon me.

Some homework exercises to accompany this text can be found on the websites of
the courses that it was written for:
\begin{align*}
&  \text{\url{https://darijgrinberg.gitlab.io/t/24wd/}}\\
&  \text{\url{https://darijgrinberg.gitlab.io/t/25wd/}}\\
&  \text{\url{https://darijgrinberg.gitlab.io/t/25sd/}}%
\end{align*}

\subsection{Plan}

This text is split into $6$ chapters:

\begin{enumerate}
\item \textbf{Induction and recursion.} Induction is one of the foundational
principles of mathematics; in a sense, it is the essence of the notion of an
integer. We explore two of its versions and many of its uses. Along the way,
we introduce some concepts (e.g., Fibonacci numbers and prime factorization)
that we will later revisit.

\item \textbf{Sums and products.} Here we define and study finite sums, finite
products, factorials and binomial coefficients. These basic algebraic concepts
appear all over mathematics, and also offer us some more practice in using induction.

\item \textbf{Elementary number theory.} Here we explore the
divisibility-related properties of integers: congruence modulo $n$; division
with remainder; prime numbers; greatest common divisors and least common
multiples. One of the most famous algorithms in mathematics -- the Euclidean
algorithm -- is encountered here, as are some more curious results such as the
$xa+yb$ theorem. As an application, we briefly discuss two cryptographical
algorithms (methods for encrypting data): Caesar ciphers and RSA.

\item \textbf{An informal introduction to enumeration.} Enumeration is another
word for counting -- specifically, counting objects satisfying some
properties, such as $3$-element subsets of a given $7$-element set. We take a
first dip into this subject here, without formally defining what counting
means; this is to be done in a later chapter.

\item \textbf{Maps (aka functions).} The notion of a map (or function, which
is synonymous) is absolutely fundamental to modern mathematics. It encompasses
the functions known from calculus, but is more general, as it allows any kinds
of input and output. We define it rigorously and introduce its basic features:
composition, inverses, injectivity, surjectivity, bijectivity.

\item \textbf{Enumeration revisited.} Now that we have learned the language of
maps, we define the size of a finite set, which allows us to rigorously speak
about counting. We reprove the results previously shown informally, and dig
deeper, answering several classes of counting questions: e.g., subsets, (weak)
compositions, or anagrams of a word.
\end{enumerate}

\begin{noncompile}
\textbf{A bit of combinatorial games.} \textit{If time allows:} We introduce
the notion of a game -- more precisely, a combinatorial game for two players,
with complete information and no chance (i.e., no randomness). We study a few
classical examples of such games, such as the game of Nim. \textbf{This
chapter has not been written yet.}
\end{noncompile}

\subsection{Notations}

We shall use the following notations:

\begin{itemize}
\item We let $\mathbb{N}$ denote the set of all nonnegative integers, that is,
$\left\{  0,1,2,\ldots\right\}  $.

\item We let $\mathbb{Z}$ denote the set of all integers (both nonnegative and negative).

\item The notation $\left\vert S\right\vert $ denotes the size (i.e., the
number of elements) of a set $S$.

\item The symbol \# means \textquotedblleft number\textquotedblright. For
example, \textquotedblleft\# of positive integers that have two
digits\textquotedblright\ means \textquotedblleft number of positive integers
that have two digits\textquotedblright.

\item The abbreviation \textquotedblleft LHS\textquotedblright\ means
\textquotedblleft left hand side\textquotedblright\ (of an equation). The
abbreviation \textquotedblleft RHS\textquotedblright\ means \textquotedblleft
right hand side\textquotedblright.

\item The symbol $\varnothing$ means the empty set.

\item The notations $\left[  n\right]  $ and $\left[  a,b\right]  $ are
defined in Definition \ref{def.sets.ive-n}.

\item The notation \textquotedblleft$:=$\textquotedblright\ means
\textquotedblleft is defined to be\textquotedblright. For example,
\textquotedblleft$s_{n}:=1+2+\cdots+n$\textquotedblright\ means that we define
$s_{n}$ to be $1+2+\cdots+n$.
\end{itemize}

\subsection{Acknowledgments}

I thank Keith Conrad, Karen Edwards, Andy Hicks and Tom Roby for helpful
advice and conversations about what a course on discrete mathematics should
contain. (Needless to say, I did not heed all of this advice in these notes.)
Furthermore, I thank Mikhail Botchkarev for corrections and comments.

Your name could stand here: Please send corrections and comments to
\newline\texttt{darijgrinberg@gmail.com} .

\newpage

\section{\label{chp.ind}Induction and recursion}

\subsection{\label{sec.ind.hanoi}The Tower of Hanoi}

\subsubsection{The puzzle}

Let me start with a puzzle called
\href{https://en.wikipedia.org/wiki/Tower_of_Hanoi}{the \textbf{Tower of
Hanoi}}.

You have $3$ pegs (or rods). The first peg has $n$ disks stacked on it. The
$n$ disks have $n$ different sizes, and they are stacked in the order of their
size, with the smallest one on top. Here is what this looks like for $n=8$
(with the $3$ pegs numbered $1,2,3$ from left to right):%
\[
\includegraphics[scale=0.3]{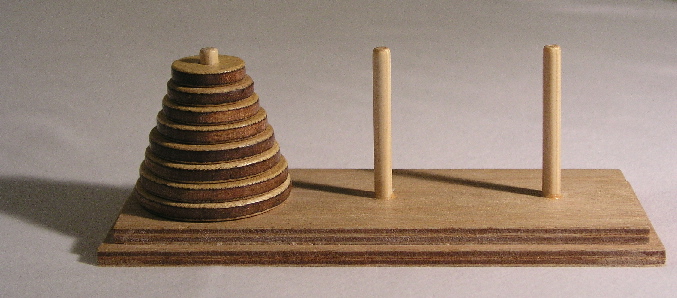}
\]
(image by \href{https://en.wikipedia.org/wiki/User:Evanherk}{User:Evanherk} on
\href{https://en.wikipedia.org/wiki/File:Tower_of_Hanoi.jpeg}{Wikipedia},
licensed under
\href{https://creativecommons.org/licenses/by-sa/3.0/deed.en}{the CC
Attribution-Share Alike 3.0 Unported License}).

You can make only one kind of moves (the so-called \textquotedblleft%
\textbf{Hanoi moves}\textquotedblright): You can take the topmost disk from
one peg and move it on top of another peg. However, you are only allowed to do
this if this disk is smaller than the other disks currently on the latter peg;
in other words, you must never stack a larger disk atop a smaller disk.

Your \textbf{goal} is to move all $n$ disks onto the third peg.

\bigskip

This game can actually be played online, e.g., at
\url{https://codepen.io/eliortabeka/pen/yOrrxG}. (Be warned that this site has
$n=7$ hardcoded into it. But you can easily fix this by modifying
\textquotedblleft disksNum = 7\textquotedblright\ and changing
\textquotedblleft minMoves = 127\textquotedblright\ to \textquotedblleft
minMoves = 0\textquotedblright. Also note that the game allows you to win by
moving all disks to peg $2$ as well, but this is clearly not a significant difference.)

\subsubsection{Some explorations}

Let us analyze the case $n=3$. In this case, one strategy to win the game
(i.e., achieve the goal) is as follows:

\begin{enumerate}
\item Move the smallest disk from peg $1$ to peg $3$.

\item Move the middle disk from peg $1$ to peg $2$.

\item Move the smallest disk from peg $3$ to peg $2$.

\item Move the largest disk from peg $1$ to peg $3$.

\item Move the smallest disk from peg $2$ to peg $1$.

\item Move the middle disk from peg $2$ to peg $3$.

\item Move the smallest disk from peg $1$ to peg $3$.
\end{enumerate}

So we can win in $7$ moves for $n=3$.

\bigskip

What about other values of $n$ ? The questions we can ask are the following:

\begin{question}
\label{quest.ind.hanoi.qn}\textbf{(a)} Can we always win the game? \medskip

\textbf{(b)} If so, then what is the smallest \# of moves\footnotemark\ we
need to make?
\end{question}

\footnotetext{The symbol \textquotedblleft\#\textquotedblright\ means
\textquotedblleft number\textquotedblright.}Let us record the answers for
small values of $n$:

\begin{itemize}
\item For $n=0$, we win in $0$ moves (since all disks -- of which there are
none -- are on peg $3$ already). This sounds very pedantic and pointless, but
it's not a bad start.

\item For $n=1$, we win in $1$ move (just move the single disk directly).

\item For $n=2$, we win in $3$ moves. Fewer moves are not enough, for fairly
simple logical reasons: We need $1$ move to free the largest disk, then $1$
move to move it to peg $3$, then $1$ more move to get the other disk on top of it.

\item For $n=3$, we win in $7$ moves (as we have seen above). But do we need
$7$ moves, or can we do with fewer?

\item For $n=4$, what happens?
\end{itemize}

Solving the problem by brute force gets harder and harder as $n$ grows. But we
can try to analyze our strategy for $n=3$ and see if there is a pattern behind it.

We observe that the largest disk moves only once, and its move is right in the
middle of the strategy. So our strategy for $n=3$ can be summarized as follows:

\begin{enumerate}
\item[1.--3.] Move the two smaller disks from peg $1$ onto peg $2$.

\item[4.] Move the largest disk from peg $1$ onto peg $3$.

\item[5.--7.] Move the two smaller disks from peg $2$ onto peg $3$.
\end{enumerate}

Moreover, the moves 1--3 in this strategy are essentially a Tower of Hanoi
game played only with the two smaller disks, except that the goal is not to
move them to peg $3$ but to move them to peg $2$ (but this doesn't matter,
because the two games are clearly \textquotedblleft
isomorphic\textquotedblright\ -- i.e., the roles of pegs $2$ and $3$ are
swapped but otherwise everything is the same). The largest disk stays at the
bottom of peg $1$ all the time and thus does not prevent any of the moves
(since all the other disks are smaller than it and thus can fit on top of it).

Move 4 moves the newly liberated largest disk from peg $1$ onto peg $3$.

Moves 5--7 are again a little Tower of Hanoi game for the two smaller disks,
except that now they have to be moved from peg $2$ to peg $3$. Again, the
largest disk (which is now on the bottom of peg $3$) does not interfere with
any of the moves.

Now the logic behind the above strategy has become clear (and also easier to memorize).

\bigskip

Does this help us solve the $n=4$ case?

Yes! We can win in $15$ moves by a strategy that has the same structure:

\begin{enumerate}
\item[1.--7.] Move the three smaller disks from peg $1$ onto peg $2$. (This is
a little Tower of Hanoi game for these three smaller disks. The largest disk
rests at the bottom of peg $1$ and does not interfere.)

\item[8.] Move the largest disk from peg $1$ onto peg $3$.

\item[9.--15.] Move the three smaller disks from peg $2$ onto peg $3$. (This
is again a little Tower of Hanoi game for these three smaller disks. The
largest disk rests at the bottom of peg $3$ and does not interfere.)
\end{enumerate}

Thus, we don't just have a strategy for $n=3$ and one for $n=4$, but actually
a \textquotedblleft meta-strategy\textquotedblright\ that lets us win the game
for $n$ disks if we know how to win it for $n-1$ disks. In a nutshell, it says
\textquotedblleft first move the $n-1$ smaller disks onto peg $2$; then move
the largest disk onto peg $3$; then move the $n-1$ smaller disks onto peg
$3$\textquotedblright. We will still call this \textquotedblleft
meta-strategy\textquotedblright\ a strategy.

\subsubsection{The numbers $m_{n}$}

Let us summarize what we gain from this strategy.

\begin{definition}
\label{def.ind.hanoi.mn}For any integer $n\geq0$, we let $m_{n}$ denote the
smallest \# of moves needed to win the Tower of Hanoi game with $n$ disks. If
the game cannot be won with $n$ disks, then we set $m_{n}=\infty$ (where
$\infty$ is not a number but just a symbol).
\end{definition}

Thus, both Question \ref{quest.ind.hanoi.qn} \textbf{(a)} and Question
\ref{quest.ind.hanoi.qn} \textbf{(b)} boil down to computing $m_{n}$.

Here is a table of small values of $m_{n}$ obtained using our strategy:%
\[%
\begin{tabular}
[c]{|c||c|c|c|c|c|c|c|c|c|}\hline
$n$ & $0$ & $1$ & $2$ & $3$ & $4$ & $5$ & $6$ & $7$ & $8$\\\hline
$m_{n}$ & $0$ & $1$ & $3$ & $7$ & $15$ & $31$ & $63$ & $127$ & $255$\\\hline
\end{tabular}
\ \
\]
Note that these values are easily computed using our strategy, because in
order to win the game for a given $n$, we have to win it for $n-1$, then make
one extra move, then win it for $n-1$ again. So we get $m_{n}=m_{n-1}%
+1+m_{n-1}=2m_{n-1}+1$ (for $n\geq1$).

Right?

\bigskip

Not so fast! We have proved that, e.g., the game can be won in $127$ moves for
$n=7$. We have not proved that it cannot be won in fewer moves. So the formula
$m_{n}=2m_{n-1}+1$ has been proved not for the smallest \# of moves needed to
win, but rather for the \# of moves needed to win \textbf{using our strategy}.
Maybe there is a better strategy that wins for $n=7$ in (say) $109$ moves?

\bigskip

So what we really have proved is the following:

\begin{proposition}
\label{prop.hanoi.1}Let $n$ be a positive integer. If $m_{n-1}$ is an integer
(i.e., if $m_{n-1}\neq\infty$), then $m_{n}\leq2m_{n-1}+1$.
\end{proposition}

To gain some writing experience, let us write out the proof in detail:

\begin{proof}
Assume that $m_{n-1}$ is an integer. Thus, we can win the game for $n-1$ disks
in $m_{n-1}$ moves. Let $S$ be the strategy (i.e., the sequence of moves)
needed to do this. So the strategy $S$ moves $n-1$ disks from peg $1$ onto peg
$3$ in $m_{n-1}$ moves.

Let $S_{23}$ be the same strategy as $S$, but with the roles of pegs $2$ and
$3$ swapped. Thus, $S_{23}$ moves $n-1$ disks from peg $1$ onto peg $2$ in
$m_{n-1}$ moves.

Let $S_{12}$ be the same strategy as $S$, but with the roles of pegs $1$ and
$2$ swapped. Thus, $S_{12}$ moves $n-1$ disks from peg $2$ onto peg $3$ in
$m_{n-1}$ moves.

Now, we proceed as follows to win the game with $n$ disks:

\begin{enumerate}
\item[A.] We use strategy $S_{23}$ to move the $n-1$ smaller disks from peg
$1$ onto peg $2$. (This is allowed because the largest disk rests at the
bottom of peg $1$ and does not interfere with the movement of smaller disks.)

\item[B.] We move the largest disk from peg $1$ onto peg $3$. (This is allowed
because this disk is free (i.e., there are no disks on top of it) and because
peg $3$ is empty, since all the other disks are on peg $2$.)

\item[C.] We use strategy $S_{12}$ to move the $n-1$ smaller disks from peg
$2$ onto peg $3$. (Again, this is allowed since the largest disk rests at the
bottom of peg $3$ and does not interfere.)
\end{enumerate}

This strategy wins the game (for $n$ disks) in $m_{n-1}+1+m_{n-1}=2m_{n-1}+1$
many moves. So the game for $n$ disks can be won in $2m_{n-1}+1$ many moves.
In other words, $m_{n}\leq2m_{n-1}+1$. This proves Proposition
\ref{prop.hanoi.1}.
\end{proof}

\bigskip

Now, let us see if the inequality $m_{n}\leq2m_{n-1}+1$ that we have proved is
an equality or just an inequality -- i.e., whether the above strategy is
optimal or there is a faster one. I claim it is the former:

\begin{proposition}
\label{prop.hanoi.2}Let $n$ be a positive integer. If $m_{n-1}$ is an integer
(i.e., if $m_{n-1}\neq\infty$), then $m_{n}=2m_{n-1}+1$.
\end{proposition}

\begin{proof}
Again, assume that $m_{n-1}$ is an integer.

We need to show that $m_{n}=2m_{n-1}+1$. It suffices to show that $m_{n}%
\geq2m_{n-1}+1$ (since Proposition \ref{prop.hanoi.1} yields $m_{n}%
\leq2m_{n-1}+1$, and we can combine these two inequalities to get
$m_{n}=2m_{n-1}+1$). In other words, it suffices to show that any winning
strategy for $n$ disks has at least $2m_{n-1}+1$ many moves.

So let us consider a winning strategy $T$ for $n$ disks. Somewhere during the
strategy $T$, the largest disk has to move (since it starts out on peg $1$ but
has to end up on peg $3$). Let us refer to these moves (the ones that move the
largest disk) as the \textbf{special moves}. There may be several special
moves or just one, but as we just said, there has to be \textbf{at least} one.

\textbf{Before the first special move} can happen, the smallest $n-1$ disks
have to be moved away from peg $1$ (since they would otherwise block the
largest disk from moving). Moreover, these smallest $n-1$ disks must all be
moved onto the same peg (since otherwise, both pegs $2$ and $3$ would be
occupied, and then the largest disk would have nowhere to move). Thus, before
the first special move can happen, we must have won the Tower of Hanoi game
for $n-1$ disks (albeit perhaps with pegs $2$ and $3$ swapped\footnote{The
goal of the Tower of Hanoi game is to move all disks onto peg $3$. If we
instead have landed them on peg $2$, then we have technically won a different
game; but this game is clearly \textquotedblleft isomorphic\textquotedblright%
\ to the Tower of Hanoi game (i.e., identical in all but how the pegs are
named), and thus requires the same number of moves.}). Hence, before the first
special move can happen, we already need to have made $m_{n-1}$ moves (since
$m_{n-1}$ is the smallest \# of moves that can win the game for $n-1$ disks).

Now, consider what happens \textbf{after the last special move}. This last
special move necessarily moves the largest disk to peg $3$ (since that's where
this disk has to come to rest). After that, we still need to move all the
other disks onto peg $3$. At the time we are making the last special move,
these other disks must all be on the same peg (since they can be neither on
the peg from which the largest disk is moving, nor on the peg to which it is
moving\footnote{because in either case, they would block the move of the
largest disk}). Therefore, after the last special move, we still need to move
all the remaining $n-1$ disks from one peg to another. And this is again
tantamount to winning the game for $n-1$ disks. So this again needs at least
$m_{n-1}$ moves.

So in total, we know that our strategy $T$ needs to have

\begin{enumerate}
\item at least $m_{n-1}$ moves before the first special move,

\item at least one special move, and

\item at least $m_{n-1}$ moves after the last special move.
\end{enumerate}

Thus, it needs to have at least $m_{n-1}+1+m_{n-1}=2m_{n-1}+1$ many moves in
total. This proves $m_{n}\geq2m_{n-1}+1$. As explained above, this completes
the proof of Proposition \ref{prop.hanoi.2}.
\end{proof}

\bigskip

Proposition \ref{prop.hanoi.2} confirms the table we have carelessly made
before:%
\[%
\begin{tabular}
[c]{|c||c|c|c|c|c|c|c|c|c|}\hline
$n$ & $0$ & $1$ & $2$ & $3$ & $4$ & $5$ & $6$ & $7$ & $8$\\\hline
$m_{n}$ & $0$ & $1$ & $3$ & $7$ & $15$ & $31$ & $63$ & $127$ & $255$\\\hline
\end{tabular}
\ \
\]
Obviously, you can keep using Proposition \ref{prop.hanoi.2} to compute
$m_{9},m_{10},m_{11},\ldots$. Indeed, the equation%
\begin{equation}
m_{n}=2m_{n-1}+1 \label{eq.hanoi.rec}%
\end{equation}
is what is called a \textbf{recursive formula} for the numbers $m_{n}$. This
means a formula that allows you to compute $m_{n}$ using the previous values
$m_{0},m_{1},\ldots,m_{n-1}$. In our case, we only need the direct predecessor
$m_{n-1}$, so this is a particularly convenient recursive formula.

\subsubsection{In search of an explicit formula}

Still, can we perhaps do better? Can we find an \textbf{explicit formula} --
i.e., one that gives us $m_{n}$ directly?

You might have guessed a formula from our table of numbers already:%
\[
m_{n}=2^{n}-1.
\]

Is there a way to see this without guessing? Let's try applying the recursive
formula (\ref{eq.hanoi.rec}) again and again, simplifying each time:%
\begin{align*}
m_{n}  &  =2m_{n-1}+1\ \ \ \ \ \ \ \ \ \ \left(  \text{by (\ref{eq.hanoi.rec}%
)}\right) \\
&  =2\left(  2m_{n-2}+1\right)  +1\ \ \ \ \ \ \ \ \ \ \left(  \text{by
(\ref{eq.hanoi.rec}), applied to }n-1\right) \\
&  =4m_{n-2}+2+1\\
&  =4\left(  2m_{n-3}+1\right)  +2+1\ \ \ \ \ \ \ \ \ \ \left(  \text{by
(\ref{eq.hanoi.rec}), applied to }n-2\right) \\
&  =8m_{n-3}+4+2+1\\
&  =8\left(  2m_{n-4}+1\right)  +4+2+1\ \ \ \ \ \ \ \ \ \ \left(  \text{by
(\ref{eq.hanoi.rec}), applied to }n-3\right) \\
&  =16m_{n-4}+8+4+2+1\\
&  =\cdots\ \ \ \ \ \ \ \ \ \ \left(  \text{keep going until you reach }%
m_{0}\right) \\
&  =2^{n}\underbrace{m_{0}}_{=0}+\,2^{n-1}+2^{n-2}+\cdots+2^{0}\\
&  =2^{n-1}+2^{n-2}+\cdots+2^{0}\\
&  =2^{0}+2^{1}+2^{2}+\cdots+2^{n-1}.
\end{align*}
I claim that the right hand side is $2^{n}-1$. Why?

\subsection{The Principle of Mathematical Induction}

At this place, I could explain why. But I prefer not to, since there is an
easier way to prove that $m_{n}=2^{n}-1$ (and anyway, the above proof of
$m_{n}=2^{0}+2^{1}+2^{2}+\cdots+2^{n-1}$ through a long computation was rather
messy and untrustworthy, so I would rather avoid relying on it).

This easier way uses one of the fundamental proof techniques in mathematics.
This technique is called \textbf{proof by induction}, and it relies on the
following principle:

\begin{theorem}
[Principle of Mathematical Induction]Let $b$ be an integer.

Let $P\left(  n\right)  $ be a mathematical statement defined for each integer
$n\geq b$.

(For example, $P\left(  n\right)  $ can be \textquotedblleft$n+1>n$%
\textquotedblright\ or \textquotedblleft$n$ is even\textquotedblright\ or
\textquotedblleft$n$ is prime\textquotedblright\ or \textquotedblleft there
exists a prime number larger than $n$\textquotedblright. Note that not every
statement needs to be true (for example, \textquotedblleft$n$ is
even\textquotedblright\ is true for some $n$'s and false for others). So
$P\left(  n\right)  $ is a statement that depends on $n$; in logic, such a
statement is called a \textbf{predicate}.)

Assume the following:

\begin{enumerate}
\item The statement $P\left(  b\right)  $ holds (i.e., the statement $P\left(
n\right)  $ holds for $n=b$).

\item For each integer $n\geq b$, the implication $P\left(  n\right)
\Longrightarrow P\left(  n+1\right)  $ holds (i.e., if $P\left(  n\right)  $
holds, then $P\left(  n+1\right)  $ does as well)\footnotemark.
\end{enumerate}

Then, the statement $P\left(  n\right)  $ holds for every integer $n\geq b$.
\end{theorem}

\footnotetext{Let me recall the meaning of the \textquotedblleft%
$\Longrightarrow$\textquotedblright\ symbol:
\par
If $A$ and $B$ are two statements, then \textquotedblleft$A\Longrightarrow
B$\textquotedblright\ means the statement \textquotedblleft if $A$, then
$B$\textquotedblright. This statement is true whenever $B$ is true, but also
true whenever $A$ is false; only in the remaining case (i.e., when $A$ is true
but $B$ is false) is it false. In other words, its truth table is as follows:%
\[%
\begin{tabular}
[c]{|c|c||c|}\hline
$A$ & $B$ & $A\Longrightarrow B$\\\hline
true & true & true\\\hline
true & false & false\\\hline
false & true & true\\\hline
false & false & true\\\hline
\end{tabular}
\ \ .
\]
You can think of it as a contract: \textquotedblleft If you make $A$ true,
then I make $B$ true\textquotedblright. If you don't make $A$ true, then this
contract places no obligation on me, since you haven't done your part! The
only way for me to violate the contract is if you make $A$ true but I don't
make $B$ true. In other words, $A\Longrightarrow B$ is a \textquotedblleft
relative\textquotedblright\ statement, which is true by default if $A$ is not.
\par
Usually, if you want to prove an implication $A\Longrightarrow B$, you start
by assuming that $A$ holds, and you need to show that $B$ holds (under this
assumption).}Before we discuss the true meaning of this principle, let me show
how to use it to prove our $m_{n}=2^{n}-1$ claim. We state this claim as a theorem:

\begin{theorem}
[explicit answer to Tower of Hanoi]\label{thm.hanoi.mn=2n-1}For each integer
$n\geq0$, we let $m_{n}$ be the smallest \# of moves needed to win the Tower
of Hanoi game with $n$ disks (or $\infty$ if it cannot be won).

Then,%
\[
m_{n}=2^{n}-1\ \ \ \ \ \ \ \ \ \ \text{for each integer }n\geq0.
\]

\end{theorem}

\begin{proof}
We denote the statement \textquotedblleft$m_{n}=2^{n}-1$\textquotedblright\ by
$P\left(  n\right)  $. So we must prove that $P\left(  n\right)  $ holds for
each integer $n\geq0$.

According to the Principle of Mathematical Induction (applied to $b=0$), it
suffices (for this purpose) to show that

\begin{enumerate}
\item the statement $P\left(  0\right)  $ holds;

\item for each integer $n\geq0$, the implication $P\left(  n\right)
\Longrightarrow P\left(  n+1\right)  $ holds.
\end{enumerate}

Proving these two claims will be our two goals; we call them Goal 1 and Goal
2. Let us see if we can achieve them.

Goal 1 is easy: The statement $P\left(  0\right)  $ is just saying that
$m_{0}=2^{0}-1$, which is true since both sides are $0$.

We now start working towards Goal 2. Let $n\geq0$ be an integer. We must prove
the implication $P\left(  n\right)  \Longrightarrow P\left(  n+1\right)  $. To
prove this, we assume that $P\left(  n\right)  $ holds, and we set out to
prove that $P\left(  n+1\right)  $ holds.

Our assumption says that $P\left(  n\right)  $ holds, i.e., that
\[
m_{n}=2^{n}-1.
\]
In particular, $m_{n}$ is an integer, so that the Tower of Hanoi game for $n$
disks is winnable.

We need to prove that $P\left(  n+1\right)  $ holds, i.e., that
\[
m_{n+1}\overset{?}{=}2^{n+1}-1.
\]
(The question mark above the equality sign just serves to remind us that we
have not proved this equality yet.)

Proposition \ref{prop.hanoi.2} yields that $m_{n}=2m_{n-1}+1$ if $n\geq1$ (and
if $m_{n-1}$ is not $\infty$). But this is not very helpful, since we are
looking for $m_{n+1}$, not for $m_{n}$.

However, we can also apply Proposition \ref{prop.hanoi.2} to $n+1$ instead of
$n$ (since $n$ is just an arbitrary integer $\geq1$ in that proposition; it is
not bound to be our current $n$). This gives us%
\[
m_{n+1}=2m_{n}+1.
\]
Thus,%
\begin{align*}
m_{n+1}  &  =2\underbrace{m_{n}}_{=2^{n}-1}+\,1=2\cdot\left(  2^{n}-1\right)
+1=2\cdot2^{n}-2+1=\underbrace{2\cdot2^{n}}_{\substack{=2^{n+1}\\\text{(by one
of the laws}\\\text{of exponents)}}}-\,1\\
&  =2^{n+1}-1.
\end{align*}
But this is precisely the statement $P\left(  n+1\right)  $. So we have shown
that $P\left(  n+1\right)  $ holds.

More precisely, we have shown that $P\left(  n+1\right)  $ holds under the
assumption that $P\left(  n\right)  $ holds. In other words, we have proved
the implication $P\left(  n\right)  \Longrightarrow P\left(  n+1\right)  $.
This achieves Goal 2.

So we have achieved both goals, and thus the Principle of Mathematical
Induction yields that $P\left(  n\right)  $ holds for every integer $n\geq0$.
In other words, $m_{n}=2^{n}-1$ holds for every integer $n\geq0$. This proves
the theorem.
\end{proof}

\bigskip

What have we really done here? How did this proof work? What is the logic
underlying the Principle of Mathematical Induction?

Let us take a look at the structure of our above proof.

Our goal was to prove that $P\left(  n\right)  $ holds for every $n\geq0$.

In other words, our goal was to prove the statements%
\[
P\left(  0\right)  ,\ P\left(  1\right)  ,\ P\left(  2\right)  ,\ P\left(
3\right)  ,\ \ldots.
\]
This is an infinite sequence of statements.

We have proved that $P\left(  0\right)  $ holds; that was our Goal 1.

We have then proved that $P\left(  n\right)  \Longrightarrow P\left(
n+1\right)  $ for each $n$. In other words, we have proved that each statement
in our sequence implies the next. In particular, $P\left(  0\right)
\Longrightarrow P\left(  1\right)  $ and $P\left(  1\right)  \Longrightarrow
P\left(  2\right)  $ and $P\left(  2\right)  \Longrightarrow P\left(
3\right)  $ and so on.

Combining $P\left(  0\right)  $ with $P\left(  0\right)  \Longrightarrow
P\left(  1\right)  $, we obtain $P\left(  1\right)  $.

Combining $P\left(  1\right)  $ with $P\left(  1\right)  \Longrightarrow
P\left(  2\right)  $, we obtain $P\left(  2\right)  $.

Combining $P\left(  2\right)  $ with $P\left(  2\right)  \Longrightarrow
P\left(  3\right)  $, we obtain $P\left(  3\right)  $.

And so on. Continuing this logic, you obtain $P\left(  4\right)  $, then
$P\left(  5\right)  $, then $P\left(  6\right)  $, and so on. By common sense,
it is clear that if you keep going on like this, you will eventually reach
each statement in our infinite sequence; i.e., you will obtain $P\left(
n\right)  $ for any given integer $n\geq0$. Of course, this reasoning is
informal (\textquotedblleft common sense\textquotedblright\ is not a
mathematical concept, nor are the words \textquotedblleft and so
on\textquotedblright).

Thus, if we want to use this kind of reasoning in a mathematical proof, we
need to state it as a precise principle and we need this principle to be true.
The Principle of Mathematical Induction is doing precisely that.

\bigskip

\begin{remark}
\label{rmk.ind.domino-metaphor}You can metaphorically think of our proof (or
any proof using the Principle of Mathematical Induction) as an infinite daisy
chain of lamps, which stand for the statements $P\left(  0\right)  ,\ P\left(
1\right)  ,\ P\left(  2\right)  ,\ \ldots$: Goal 1 turns the first lamp on,
whereas Goal 2 ensures that each lamp turns the next on when it is turned on itself.

Or, to use a more commonplace illustration, you have an infinite sequence of
dominos arranged in a row, at sufficiently close distances so that tipping
over one domino will tip over the next. After you tip over the first domino,
all the dominos will eventually fall down. (The dominos here stand for the
statements $P\left(  0\right)  ,\ P\left(  1\right)  ,\ P\left(  2\right)
,\ \ldots$.)
\end{remark}

\bigskip

I called the Principle of Mathematical Induction a theorem, but I will not
prove it, since it is one of the fundamental axioms of mathematics. You can at
best replace it by a different axiom, but this doesn't change much; you need
some kind of axiom that allows you to \textquotedblleft chain
together\textquotedblright\ arbitrarily many little implications.

\subsection{Some more proofs by induction}

A proof that uses the Principle of Mathematical Induction is called a
\textbf{proof by induction} (or an \textbf{induction proof}, or an
\textbf{inductive proof}\footnote{This has barely anything to do with
\textquotedblleft inductive reasoning\textquotedblright\ as understood by
philosophers (known to mathematics as \textquotedblleft
generalization\textquotedblright, and not considered as a method of proof per
se).}). So our above proof of Theorem \ref{thm.hanoi.mn=2n-1} was a proof by induction.

\subsubsection{The sum of the first $n$ positive integers}

Let us see another (simpler) example of a proof by induction. We will prove
the following result:

\begin{theorem}
[\textquotedblleft Little Gauss formula\textquotedblright]%
\label{thm.ind.1+2+...+n}For every integer $n\geq0$, we have%
\[
1+2+\cdots+n=\dfrac{n\left(  n+1\right)  }{2}.
\]

\end{theorem}

The LHS (= left hand side) here is understood to be the sum of the first $n$
positive integers. For $n=0$, this sum is an empty sum (i.e., it has no
addends at all), so its value is $0$ by definition.

\begin{proof}
[First proof of Theorem \ref{thm.ind.1+2+...+n}.]We set
\[
s_{n}:=1+2+\cdots+n
\]
for each $n\geq0$. Thus, we must prove that $s_{n}=\dfrac{n\left(  n+1\right)
}{2}$ for each $n\geq0$.

Let us denote the statement \textquotedblleft$s_{n}=\dfrac{n\left(
n+1\right)  }{2}$\textquotedblright\ by $P\left(  n\right)  $. So we need to
prove that $P\left(  n\right)  $ holds for every $n\geq0$.

According to the Principle of Mathematical Induction, it suffices to show that

\begin{enumerate}
\item the statement $P\left(  0\right)  $ holds;

\item for each $n\geq0$, the implication $P\left(  n\right)  \Longrightarrow
P\left(  n+1\right)  $ holds.
\end{enumerate}

Goal 1 is easy: To prove $P\left(  0\right)  $, we must show that
$s_{0}=\dfrac{0\left(  0+1\right)  }{2}$, but this is true because both sides
equal $0$.

Now to Goal 2. We let $n\geq0$ be an integer, and we want to prove the
implication $P\left(  n\right)  \Longrightarrow P\left(  n+1\right)  $. So we
assume that $P\left(  n\right)  $ holds, and we set out to prove $P\left(
n+1\right)  $.

By assumption, $P\left(  n\right)  $ holds, so that we have%
\[
s_{n}=\dfrac{n\left(  n+1\right)  }{2}.
\]
We must prove $P\left(  n+1\right)  $; in other words, we must prove that%
\[
s_{n+1}\overset{?}{=}\dfrac{\left(  n+1\right)  \left(  \left(  n+1\right)
+1\right)  }{2}.
\]

To do so, we observe that%
\begin{align*}
s_{n+1}  &  =1+2+\cdots+\left(  n+1\right)  =\underbrace{\left(
1+2+\cdots+n\right)  }_{=s_{n}}+\left(  n+1\right) \\
&  =s_{n}+\left(  n+1\right)  =\dfrac{n\left(  n+1\right)  }{2}+\left(
n+1\right)  \ \ \ \ \ \ \ \ \ \ \left(  \text{since }s_{n}=\dfrac{n\left(
n+1\right)  }{2}\right) \\
&  =\dfrac{n\left(  n+1\right)  }{2}+\dfrac{2\left(  n+1\right)  }{2}%
=\dfrac{\left(  n+2\right)  \left(  n+1\right)  }{2}=\dfrac{\left(
n+1\right)  \left(  n+2\right)  }{2}\\
&  =\dfrac{\left(  n+1\right)  \left(  \left(  n+1\right)  +1\right)  }{2}.
\end{align*}
In other words, $P\left(  n+1\right)  $ holds. Thus, we have proved the
implication $P\left(  n\right)  \Longrightarrow P\left(  n+1\right)  $.

We have now achieved both goals, so the Principle of Mathematical Induction
yields that $P\left(  n\right)  $ holds for every $n\geq0$. This proves the theorem.
\end{proof}

There is also a non-inductive proof; this is how Gauss supposedly did it:

\begin{proof}
[Second proof of Theorem \ref{thm.ind.1+2+...+n}.]We have%
\begin{align*}
&  2\cdot\left(  1+2+\cdots+n\right) \\
&  =\left(  1+2+\cdots+n\right)  +\left(  1+2+\cdots+n\right) \\
&  =\left(  1+2+\cdots+n\right)  +\left(  n+\left(  n-1\right)  +\cdots
+1\right) \\
&  \ \ \ \ \ \ \ \ \ \ \ \ \ \ \ \ \ \ \ \ \left(
\begin{array}
[c]{c}%
\text{here, we turned the second sum upside-down, i.e.,}\\
\text{we reversed the order of its addends}%
\end{array}
\right) \\
&  =\underbrace{\left(  1+n\right)  }_{=n+1}+\underbrace{\left(  2+\left(
n-1\right)  \right)  }_{=n+1}+\cdots+\underbrace{\left(  n+1\right)  }%
_{=n+1}\\
&  \ \ \ \ \ \ \ \ \ \ \ \ \ \ \ \ \ \ \ \ \left(
\begin{array}
[c]{c}%
\text{here, we rearranged the sum by matching}\\
\text{up each addend inside the first pair of}\\
\text{parentheses with the corresponding addend}\\
\text{inside the second pair of parentheses}%
\end{array}
\right) \\
&  =\underbrace{\left(  n+1\right)  +\left(  n+1\right)  +\cdots+\left(
n+1\right)  }_{n\text{ addends}}\\
&  =n\cdot\left(  n+1\right)  .
\end{align*}
Dividing this by $2$, we find%
\[
1+2+\cdots+n=\dfrac{n\cdot\left(  n+1\right)  }{2},
\]
and thus Theorem \ref{thm.ind.1+2+...+n} is proved again.
\end{proof}

\subsubsection{The sum of the squares of the first $n$ positive integers}

Here is a similar theorem:

\begin{theorem}
\label{thm.ind.12+22+...+n2}For every integer $n\geq0$, we have%
\[
1^{2}+2^{2}+\cdots+n^{2}=\dfrac{n\left(  n+1\right)  \left(  2n+1\right)  }%
{6}.
\]

\end{theorem}

\begin{proof}
The following proof is almost a word-by-word copy of the first proof of
Theorem \ref{thm.ind.1+2+...+n}. The structure is the same; only the
calculations change.

We set
\[
s_{n}:=1^{2}+2^{2}+\cdots+n^{2}.
\]
Thus, we must prove that $s_{n}=\dfrac{n\left(  n+1\right)  \left(
2n+1\right)  }{6}$ for each $n\geq0$.

Let us denote the statement \textquotedblleft$s_{n}=\dfrac{n\left(
n+1\right)  \left(  2n+1\right)  }{6}$\textquotedblright\ by $P\left(
n\right)  $. So we need to prove that $P\left(  n\right)  $ holds for every
$n\geq0$.

According to the Principle of Mathematical Induction, it suffices to show that

\begin{enumerate}
\item the statement $P\left(  0\right)  $ holds;

\item for each $n\geq0$, the implication $P\left(  n\right)  \Longrightarrow
P\left(  n+1\right)  $ holds.
\end{enumerate}

Goal 1 is easy: To prove $P\left(  0\right)  $, we must show that
$s_{0}=\dfrac{0\left(  0+1\right)  \left(  2\cdot0+1\right)  }{6}$, but this
is true because both sides equal $0$.

Now to Goal 2. We let $n\geq0$ be an integer, and we want to prove the
implication $P\left(  n\right)  \Longrightarrow P\left(  n+1\right)  $. So we
assume that $P\left(  n\right)  $ holds, and we set out to prove $P\left(
n+1\right)  $.

By assumption, $P\left(  n\right)  $ holds, so that we have%
\[
s_{n}=\dfrac{n\left(  n+1\right)  \left(  2n+1\right)  }{6}.
\]
We must prove $P\left(  n+1\right)  $; in other words, we must prove that
\[
s_{n+1}\overset{?}{=}\dfrac{\left(  n+1\right)  \left(  \left(  n+1\right)
+1\right)  \left(  2\left(  n+1\right)  +1\right)  }{6}.
\]

To do so, we observe that%
\begin{align*}
s_{n+1}  &  =1^{2}+2^{2}+\cdots+\left(  n+1\right)  ^{2}\\
&  =\underbrace{\left(  1^{2}+2^{2}+\cdots+n^{2}\right)  }_{=s_{n}}+\left(
n+1\right)  ^{2}\\
&  =s_{n}+\left(  n+1\right)  ^{2}\\
&  =\dfrac{n\left(  n+1\right)  \left(  2n+1\right)  }{6}+\left(  n+1\right)
^{2}\ \ \ \ \ \ \ \ \ \ \left(  \text{since }s_{n}=\dfrac{n\left(  n+1\right)
\left(  2n+1\right)  }{6}\right) \\
&  =\left(  n+1\right)  \cdot\left(  \dfrac{n\left(  2n+1\right)  }{6}+\left(
n+1\right)  \right) \\
&  =\left(  n+1\right)  \cdot\dfrac{2n^{2}+7n+6}{6}\\
&  =\dfrac{\left(  n+1\right)  \left(  2n^{2}+7n+6\right)  }{6}\\
&  =\dfrac{\left(  n+1\right)  \left(  n+2\right)  \left(  2n+3\right)  }%
{6}\ \ \ \ \ \ \ \ \ \ \left(
\begin{array}
[c]{c}%
\text{since }2n^{2}+7n+6\text{ can be}\\
\text{factored as }\left(  n+2\right)  \left(  2n+3\right)
\end{array}
\right) \\
&  =\dfrac{\left(  n+1\right)  \left(  \left(  n+1\right)  +1\right)  \left(
2\left(  n+1\right)  +1\right)  }{6}.
\end{align*}
In other words, $P\left(  n+1\right)  $ holds. Thus, we have proved the
implication $P\left(  n\right)  \Longrightarrow P\left(  n+1\right)  $.

We have now achieved both goals, so the Principle of Mathematical Induction
yields that $P\left(  n\right)  $ holds for every $n\geq0$. This proves the theorem.
\end{proof}

\bigskip

As we said, our above proof of Theorem \ref{thm.ind.12+22+...+n2} was an
almost verbatim copy of our first proof of Theorem \ref{thm.ind.1+2+...+n}; we
only needed to make the obvious changes and calculate a little bit harder.
Both proofs were more or less determined by the idea to use induction. In
contrast, the slick second proof of Theorem \ref{thm.ind.1+2+...+n} is not as
easy to adapt to Theorem \ref{thm.ind.12+22+...+n2}\footnote{But it can be
adapted; see \cite{Levrie15} for such a proof of Theorem
\ref{thm.ind.12+22+...+n2}.}. So the induction proof has the advantage of
being easier to generalize.

However, it has the disadvantage that it can only be used to \textbf{prove} a
formula (in our case, $1+2+\cdots+n=\dfrac{n\left(  n+1\right)  }{2}$ or
$1^{2}+2^{2}+\cdots+n^{2}=\dfrac{n\left(  n+1\right)  \left(  2n+1\right)
}{6}$), not to \textbf{find} this formula in the first place. We could not
have used induction to answer the question \textquotedblleft what is
$1+2+\cdots+n$?\textquotedblright; we could only use it to prove the answer
after guessing it in some way.

\begin{exercise}
\label{exe.23hw1.1a}Prove that
\[
1^{3}+2^{3}+\cdots+n^{3}=\left(  \dfrac{n\left(  n+1\right)  }{2}\right)  ^{2}%
\]
for each nonnegative integer $n$. (The left hand side here is the sum of the
\textbf{cubes} of the first $n$ positive integers.)
\end{exercise}

\subsection{Notations for an induction proof}

Here is some standard terminology that is commonly used in proofs by
induction. Let's say that you are proving a statement of the form $P\left(
n\right)  $ for every integer $n\geq b$ (where $b$ is some fixed integer).

\begin{itemize}
\item The $n$ is called the \textbf{induction variable}; you say that you
\textbf{induct on }$n$. It does not have to be called $n$. Your statement
might just as well be \textquotedblleft for every integer $a\geq0$, we have
$1+2+\cdots+a=\dfrac{a\left(  a+1\right)  }{2}$\textquotedblright, and then
you can prove it by inducting on $a$.

\item The proof of $P\left(  b\right)  $ (that is, Goal 1 in our above proofs)
is called the \textbf{induction base} or the \textbf{base case}. In our above
examples, this was always the proof of $P\left(  0\right)  $, but in general
$b$ can be another integer. (For example, if you are proving the statement
\textquotedblleft every integer $n\geq4$ satisfies $2^{n}\geq n^{2}%
$\textquotedblright, then $b$ will have to be $4$, so your induction base
consists in proving that $2^{4}\geq4^{2}$.)

\item The proof of \textquotedblleft$P\left(  n\right)  \Longrightarrow
P\left(  n+1\right)  $ for every $n\geq b$\textquotedblright\ (that is, Goal 2
in our above proofs) is called the \textbf{induction step}. For example, in
the proof of Theorem \ref{thm.ind.12+22+...+n2}, this was the part where we
assumed that $s_{n}=\dfrac{n\left(  n+1\right)  \left(  2n+1\right)  }{6}$ and
proved that $s_{n+1}=\dfrac{\left(  n+1\right)  \left(  \left(  n+1\right)
+1\right)  \left(  2\left(  n+1\right)  +1\right)  }{6}$.

In the induction step, the assumption that $P\left(  n\right)  $ holds is
called the \textbf{induction hypothesis} or the \textbf{induction assumption},
and the claim that $P\left(  n+1\right)  $ holds (this is the claim that you
are trying to prove) is called the \textbf{induction goal}. The induction step
is complete when the induction goal is reached (i.e., proved).
\end{itemize}

As an example, let us rewrite our above proof of Theorem
\ref{thm.hanoi.mn=2n-1} using this language:

\begin{proof}
[Proof of Theorem \ref{thm.hanoi.mn=2n-1}, rewritten.]We induct on $n$.

\textit{Base case:} The theorem\footnote{i.e., Theorem \ref{thm.hanoi.mn=2n-1}%
} holds for $n=0$, since both $m_{0}$ and $2^{0}-1$ equal $0$.

\textit{Induction step:} Let $n\geq0$ be an integer. We assume that the
theorem holds for $n$ (this is what we previously called $P\left(  n\right)
$). We will now show that the theorem holds for $n+1$ as well (this is what we
previously called $P\left(  n+1\right)  $).

We have assumed that the theorem holds for $n$. In other words, $m_{n}%
=2^{n}-1$. This is our induction hypothesis.

We must prove that the theorem holds for $n+1$. In other words, we must prove
that $m_{n+1}\overset{?}{=}2^{n+1}-1$.

To prove this, we apply Proposition \ref{prop.hanoi.2} to $n+1$ instead of $n$
(we can do this, since $m_{n}=2^{n}-1$ is not $\infty$). This gives us%
\begin{align*}
m_{n+1}  &  =2\underbrace{m_{n}}_{\substack{=2^{n}-1\\\text{(by the induction
hypothesis)}}}+\,1=2\cdot\left(  2^{n}-1\right)  +1\\
&  =2\cdot2^{n}-2+1=\underbrace{2\cdot2^{n}}_{=2^{n+1}}-\,1=2^{n+1}-1.
\end{align*}
Thus, the induction goal is reached, and the induction is complete. Hence, the
theorem is proved.
\end{proof}

\subsection{The Fibonacci numbers}

\subsubsection{Definition}

Our next applications of induction will be some properties of the
\textbf{Fibonacci sequence}. The Fibonacci sequence is defined
\textbf{recursively} -- i.e., a given entry is not defined directly, but
rather defined in terms of the previous entries. Specifically, it is defined
as follows:

\begin{definition}
\label{def.ind.fib}The \textbf{Fibonacci sequence} is the sequence $\left(
f_{0},f_{1},f_{2},\ldots\right)  $ of nonnegative integers defined recursively
by setting%
\begin{align*}
f_{0}  &  =0,\ \ \ \ \ \ \ \ \ \ f_{1}=1,\ \ \ \ \ \ \ \ \ \ \text{and}\\
f_{n}  &  =f_{n-1}+f_{n-2}\ \ \ \ \ \ \ \ \ \ \text{for each }n\geq2.
\end{align*}

\end{definition}

In other words, the Fibonacci sequence starts with the two entries $0$ and
$1$, and then every next entry is the sum of the two previous entries.

The entries of the Fibonacci sequence are called the \textbf{Fibonacci
numbers}. Let us compute the first fourteen of them:%
\[%
\begin{tabular}
[c]{|c||c|c|c|c|c|c|c|c|c|c|c|c|c|c|}\hline
$n$ & $0$ & $1$ & $2$ & $3$ & $4$ & $5$ & $6$ & $7$ & $8$ & $9$ & $10$ & $11$
& $12$ & $13$\\\hline
$f_{n}$ & $0$ & $1$ & $1$ & $2$ & $3$ & $5$ & $8$ & $13$ & $21$ & $34$ & $55$
& $89$ & $144$ & $233$\\\hline
\end{tabular}
\ \ \ .
\]

As we see, a recursive definition is a perfectly valid way to define (e.g.) a
sequence of numbers. It allows you to compute each entry of the sequence
eventually, as long as you compute the entries in order (i.e., first $f_{0}$,
then $f_{1}$, then $f_{2}$, and so on). In a sense, the reason why this works
is the same as the reason why induction works: You can get to any integer
$n\geq0$ if you start at $0$ and keep adding $1$.

Note that it is important that our recursive definition of $f_{n}$ uses only
previous entries of the sequence (in our case, $f_{n-1}$ and $f_{n-2}$). If we
had instead defined the Fibonacci sequence by%
\[
f_{n}=f_{n+1}-f_{n-2},
\]
then we could not even compute $f_{2}$, since this would require knowing
$f_{3}$, which would in turn require knowing $f_{4}$, and so on.

\subsubsection{The sum of the first $n$ positive Fibonacci numbers}

The Fibonacci sequence is famous for its many properties and
patterns\footnote{There is an entire book about it (Vorobiev's \cite{Vorobi02}%
, which I can recommend).}. Here is a first one:

\begin{theorem}
\label{thm.fib.sum-first-n}For any integer $n\geq0$, we have%
\[
f_{1}+f_{2}+\cdots+f_{n}=f_{n+2}-1.
\]

\end{theorem}

For example, for $n=8$, this is saying that%
\[
1+1+2+3+5+8+13+21=55-1.
\]

\begin{proof}
[Proof of Theorem \ref{thm.fib.sum-first-n}.]We induct on $n$.

\textit{Base case:} For $n=0$, the theorem claims that $f_{1}+f_{2}%
+\cdots+f_{0}=f_{0+2}-1$. This is true, since the LHS\footnote{I remind that
the abbreviations \textquotedblleft LHS\textquotedblright\ and
\textquotedblleft RHS\textquotedblright\ mean \textquotedblleft left hand
side\textquotedblright\ and \textquotedblleft right hand
side\textquotedblright, respectively.} is an empty sum (thus $=0$) whereas the
RHS is $f_{2}-1=1-1=0$.

\textit{Induction step:} Let $n\geq0$ be an integer. Assume that the theorem
holds for $n$.

We must prove that the theorem holds for $n+1$.

So we assumed that
\[
f_{1}+f_{2}+\cdots+f_{n}=f_{n+2}-1.
\]
We must prove that%
\[
f_{1}+f_{2}+\cdots+f_{n+1}\overset{?}{=}f_{\left(  n+1\right)  +2}-1.
\]

We have%
\begin{align*}
f_{1}+f_{2}+\cdots+f_{n+1}  &  =\underbrace{\left(  f_{1}+f_{2}+\cdots
+f_{n}\right)  }_{\substack{=f_{n+2}-1\\\text{(by our induction hypothesis)}%
}}+\,f_{n+1}=f_{n+2}-1+f_{n+1}\\
&  =\underbrace{f_{n+2}+f_{n+1}}_{\substack{=f_{n+3}\\\text{(since the
recursive definition}\\\text{of the Fibonacci sequence}\\\text{yields }%
f_{n+3}=f_{n+2}+f_{n+1}\text{)}}}-\,1=f_{n+3}-1\\
&  =f_{\left(  n+1\right)  +2}-1\ \ \ \ \ \ \ \ \ \ \left(  \text{since
}n+3=\left(  n+1\right)  +2\right)  .
\end{align*}
This is precisely what we wanted to prove -- i.e., it says that the theorem
holds for $n+1$. This completes the induction step. Thus, the theorem is proved.
\end{proof}

\subsection{\label{sec.ind.geosums}Some more examples of induction}

Let us see some more examples of proofs by induction. The following theorem I
have already mentioned at the end of Section \ref{sec.ind.hanoi}:

\begin{theorem}
\label{thm.geo-sum.2}For any integer $n\geq0$, we have
\[
2^{0}+2^{1}+2^{2}+\cdots+2^{n-1}=2^{n}-1.
\]

\end{theorem}

\begin{proof}
We induct on $n$.

\textit{Base case:} For $n=0$, the equality $2^{0}+2^{1}+2^{2}+\cdots
+2^{n-1}=2^{n}-1$ is true, because the LHS\footnote{\textquotedblleft
LHS\textquotedblright\ means \textquotedblleft left-hand
side\textquotedblright. Likewise, \textquotedblleft RHS\textquotedblright%
\ means \textquotedblleft right-hand side\textquotedblright.} is an empty sum
and thus equals $0$, whereas the RHS is $2^{0}-1=1-1=0$.

\textit{Induction step:} Let $n$ be an integer $\geq0$. Assume that Theorem
\ref{thm.geo-sum.2} holds for $n$, i.e., that we have%
\[
2^{0}+2^{1}+2^{2}+\cdots+2^{n-1}=2^{n}-1.
\]
We must prove that Theorem \ref{thm.geo-sum.2} holds for $n+1$ as well, i.e.,
that we have%
\[
2^{0}+2^{1}+2^{2}+\cdots+2^{\left(  n+1\right)  -1}=2^{n+1}-1.
\]

However,%
\begin{align*}
2^{0}+2^{1}+2^{2}+\cdots+2^{\left(  n+1\right)  -1}  &  =2^{0}+2^{1}%
+2^{2}+\cdots+2^{n}\\
&  =\underbrace{\left(  2^{0}+2^{1}+2^{2}+\cdots+2^{n-1}\right)
}_{\substack{=2^{n}-1\\\text{(by the induction hypothesis)}}}+\,2^{n}\\
&  =2^{n}-1+2^{n}=\underbrace{2\cdot2^{n}}_{=2^{n+1}}-\,1=2^{n+1}-1,
\end{align*}
which is precisely what we want: This shows that Theorem \ref{thm.geo-sum.2}
holds for $n+1$. Thus, our induction step is complete, and Theorem
\ref{thm.geo-sum.2} is proved.
\end{proof}

Theorem \ref{thm.geo-sum.2} can be generalized:

\begin{theorem}
\label{thm.geo-sum.xy}Let $x$ and $y$ be any two numbers. Then, for any
integer $n\geq0$, we have%
\[
\left(  x-y\right)  \left(  x^{n-1}+x^{n-2}y+x^{n-3}y^{2}+\cdots+x^{2}%
y^{n-3}+xy^{n-2}+y^{n-1}\right)  =x^{n}-y^{n}.
\]
Here, the big sum in the parentheses is the sum of all products $x^{i}y^{j}$
where $i$ and $j$ are nonnegative integers with $i+j=n-1$.
\end{theorem}

Before we prove this, let us give some examples for what this theorem actually says:

\begin{itemize}
\item For $n=2$, Theorem \ref{thm.geo-sum.xy} says that%
\[
\left(  x-y\right)  \left(  x+y\right)  =x^{2}-y^{2}.
\]

\item For $n=3$, Theorem \ref{thm.geo-sum.xy} says that%
\[
\left(  x-y\right)  \left(  x^{2}+xy+y^{2}\right)  =x^{3}-y^{3}.
\]

\item For $n=4$, Theorem \ref{thm.geo-sum.xy} says that%
\[
\left(  x-y\right)  \left(  x^{3}+x^{2}y+xy^{2}+y^{3}\right)  =x^{4}-y^{4}.
\]

\item For $x=2$ and $y=1$, Theorem \ref{thm.geo-sum.xy} says that%
\[
\left(  2-1\right)  \left(  2^{n-1}+2^{n-2}1+2^{n-3}1^{2}+\cdots+2^{2}%
1^{n-3}+2\cdot1^{n-2}+1^{n-1}\right)  =2^{n}-1^{n}.
\]
Since any power of $1$ is $1$ (and since the $2-1$ factor also equals $1$),
this simplifies to%
\[
2^{n-1}+2^{n-2}+2^{n-3}+\cdots+2^{2}+2+1=2^{n}-1,
\]
which is precisely Theorem \ref{thm.geo-sum.2}. Thus, Theorem
\ref{thm.geo-sum.xy} generalizes Theorem \ref{thm.geo-sum.2}.
\end{itemize}

Let us now prove Theorem \ref{thm.geo-sum.xy}:

\begin{proof}
[Proof of Theorem \ref{thm.geo-sum.xy}.]We induct on $n$.

\textit{Base case:} For $n=0$, the claim%
\[
\left(  x-y\right)  \left(  x^{n-1}+x^{n-2}y+x^{n-3}y^{2}+\cdots+x^{2}%
y^{n-3}+xy^{n-2}+y^{n-1}\right)  =x^{n}-y^{n}%
\]
is true, since the LHS is $0$ (because the second factor is an empty sum),
while the RHS is $x^{0}-y^{0}=1-1=0$ as well.

\textit{Induction step:} Let $n\geq0$ be an integer. Assume that Theorem
\ref{thm.geo-sum.xy} is true for $n$. That is, assume that%
\[
\left(  x-y\right)  \left(  x^{n-1}+x^{n-2}y+x^{n-3}y^{2}+\cdots+x^{2}%
y^{n-3}+xy^{n-2}+y^{n-1}\right)  =x^{n}-y^{n}.
\]
We must prove that Theorem \ref{thm.geo-sum.xy} is also true for $n+1$. That
is, we must prove that%
\[
\left(  x-y\right)  \left(  x^{n}+x^{n-1}y+x^{n-2}y^{2}+\cdots+x^{3}%
y^{n-3}+x^{2}y^{n-2}+xy^{n-1}+y^{n}\right)  =x^{n+1}-y^{n+1}.
\]

We begin by extracting the $y^{n}$ addend from the long sum in the second pair
of parentheses in this equation. We thus obtain
\begin{align*}
&  \left(  x-y\right)  \left(  x^{n}+x^{n-1}y+x^{n-2}y^{2}+\cdots+x^{3}%
y^{n-3}+x^{2}y^{n-2}+xy^{n-1}+y^{n}\right) \\
&  =\left(  x-y\right)  \underbrace{\left(  x^{n}+x^{n-1}y+x^{n-2}y^{2}%
+\cdots+x^{3}y^{n-3}+x^{2}y^{n-2}+xy^{n-1}\right)  }_{\substack{=\left(
x^{n-1}+x^{n-2}y+x^{n-3}y^{2}+\cdots+x^{2}y^{n-3}+xy^{n-2}+y^{n-1}\right)
x\\\text{(here, we have factored out an }x\text{ from the sum)}}}+\left(
x-y\right)  y^{n}\\
&  =\underbrace{\left(  x-y\right)  \left(  x^{n-1}+x^{n-2}y+x^{n-3}%
y^{2}+\cdots+x^{2}y^{n-3}+xy^{n-2}+y^{n-1}\right)  }_{\substack{=x^{n}%
-y^{n}\\\text{(by the induction hypothesis)}}}x+\left(  x-y\right)  y^{n}\\
&  =\left(  x^{n}-y^{n}\right)  x+\left(  x-y\right)  y^{n}=x^{n+1}%
-xy^{n}+xy^{n}-y^{n+1}=x^{n+1}-y^{n+1}.
\end{align*}
This means precisely that Theorem \ref{thm.geo-sum.xy} is also true for $n+1$.
Thus, the induction step is complete, and the theorem is proved.
\end{proof}

Another useful particular case of Theorem \ref{thm.geo-sum.xy} is the
following equality:\footnote{A \textquotedblleft corollary\textquotedblright%
\ means a theorem that follows easily from another theorem.}

\begin{corollary}
\label{cor.geo-sum.q}Let $q$ be a number distinct from $1$. Let $n\geq0$ be an
integer. Then,%
\[
q^{0}+q^{1}+q^{2}+\cdots+q^{n-1}=\dfrac{q^{n}-1}{q-1}.
\]

\end{corollary}

\begin{proof}
Apply Theorem \ref{thm.geo-sum.xy} to $x=q$ and $y=1$. We obtain%
\[
\left(  q-1\right)  \left(  q^{n-1}+q^{n-2}1+q^{n-3}1^{2}+\cdots+q^{2}%
1^{n-3}+q\cdot1^{n-2}+1^{n-1}\right)  =q^{n}-1^{n}.
\]
Simplifying this, we obtain%
\[
\left(  q-1\right)  \left(  q^{n-1}+q^{n-2}+q^{n-3}+\cdots+q^{2}+q+1\right)
=q^{n}-1.
\]
Thus,%
\[
q^{n-1}+q^{n-2}+q^{n-3}+\cdots+q^{2}+q+1=\dfrac{q^{n}-1}{q-1}.
\]
In other words,%
\[
q^{0}+q^{1}+q^{2}+\cdots+q^{n-1}=\dfrac{q^{n}-1}{q-1}%
\]
(since the sum on the left hand side can be rearranged in any order). This
proves Corollary \ref{cor.geo-sum.q}.
\end{proof}

Corollary \ref{cor.geo-sum.q} is often called the \emph{geometric sum
formula}, since it helps compute the sum of a geometric sequence or a
geometric series.

\begin{exercise}
\label{exe.23hw1.2}Let $\mathbb{N}$ denote the set of all nonnegative integers
(that is, $\left\{  0,1,2,\ldots\right\}  $). Let $q$ and $d$ be two real
numbers such that $q\neq1$. Let $\left(  a_{0},a_{1},a_{2},\ldots\right)  $ be
a sequence of real numbers. Assume that%
\begin{equation}
a_{n+1}=qa_{n}+d\ \ \ \ \ \ \ \ \ \ \text{for each }n\in\mathbb{N}.
\label{eq.prop.ind.ari-geo.ass}%
\end{equation}
Prove that%
\begin{equation}
a_{n}=q^{n}a_{0}+\dfrac{q^{n}-1}{q-1}d\ \ \ \ \ \ \ \ \ \ \text{for each }%
n\in\mathbb{N}. \label{eq.prop.ind.ari-geo.claim}%
\end{equation}

\end{exercise}

\subsection{How not to use induction}

Induction proofs can be slippery:

\begin{faketheorem}
\label{fakethm.ind.horses}In any set of $n\geq1$ horses, all the horses are
the same color.
\end{faketheorem}

\begin{proof}
We induct on $n$.

\textit{Base case:} This is clearly true for $n=1$, since a single horse
always has the same color as itself.

\textit{Induction step:} Let $n\geq1$ be an integer. We assume that the fake
theorem holds for $n$, i.e., that any $n$ horses are the same color.

We must prove that it also holds for $n+1$, i.e., that any $n+1$ horses are
the same color.

So let $H_{1},H_{2},\ldots,H_{n+1}$ be $n+1$ horses.

By our induction hypothesis, the first $n$ horses $H_{1},H_{2},\ldots,H_{n}$
are the same color.

Again by our induction hypothesis, the last $n$ horses $H_{2},H_{3}%
,\ldots,H_{n+1}$ are the same color.

Now, consider the first horse $H_{1}$ and the last horse $H_{n+1}$. They both
have the same color as the \textquotedblleft middle horses\textquotedblright%
\ $H_{2},H_{3},\ldots,H_{n}$ (according to the preceding two paragraphs).
Thus, all the $n+1$ horses have the same color, right? \medskip

When a claim is as obviously wrong as this one, there is an easy way to find
the mistake in the proof: You just look at some example in which the claim is
wrong, and you trace the proof on this example. The first time you see a wrong
conclusion, that's where the error probably is.

Fake Theorem \ref{fakethm.ind.horses} is wrong for $n=2$ already, i.e., for
two horses. So let us see where the induction step goes wrong when $n=1$ (that
is, going from $1$ horse to $2$ horses). In this induction step, we claim that
$H_{1}$ and $H_{n+1}=H_{2}$ both have the same color as the \textquotedblleft
middle horses\textquotedblright\ $H_{2},H_{3},\ldots,H_{1}$. But there are no
\textquotedblleft middle horses\textquotedblright, so it makes no sense to
have the same color as these \textquotedblleft middle horses\textquotedblright%
. So the argument doesn't work.

Thus, our mistake was to implicitly treat the \textquotedblleft middle
horses\textquotedblright\ as if they existed. They do exist for any $n>1$, but
not for $n=1$, and thus our induction step breaks down for $n=1$.
\end{proof}

Note how one little mistake has brought down the entire proof! For an
induction proof to work, the induction step needs to work for all $n$; that
is, we need the implication $P\left(  n\right)  \Longrightarrow P\left(
n+1\right)  $ to hold for every $n$. If even one of these implications breaks
down, the whole chain is disconnected, and all the statements $P\left(
n\right)  $ \textquotedblleft to the right of\textquotedblright\ this breaking
point are no longer guaranteed to hold. For example, if we have a statement
$P\left(  n\right)  $ for each $n\geq0$, and we have proved the base case
$P\left(  0\right)  $ and the implication $P\left(  n\right)  \Longrightarrow
P\left(  n+1\right)  $ for all $n\neq4$, then we can conclude that $P\left(
0\right)  ,\ P\left(  1\right)  ,\ P\left(  2\right)  ,\ P\left(  3\right)  $
and $P\left(  4\right)  $ hold, but we cannot guarantee that any of $P\left(
5\right)  ,\ P\left(  6\right)  ,\ P\left(  7\right)  ,\ \ldots$ hold. As so
often, a chain is only as strong as its weakest link.

\subsection{More on the Fibonacci numbers}

Recall the Fibonacci sequence, which we defined in Definition
\ref{def.ind.fib}. We recall that it is the sequence $\left(  f_{0}%
,f_{1},f_{2},\ldots\right)  $ of nonnegative integers defined recursively by
setting $f_{0}=0$ and $f_{1}=1$ and $f_{n}=f_{n-1}+f_{n-2}$ for each $n\geq2$.

The entries of the Fibonacci sequence are called the \textbf{Fibonacci
numbers}. Here are the first few:%
\[%
\begin{tabular}
[c]{|c||c|c|c|c|c|c|c|c|c|c|c|c|c|c|}\hline
$n$ & $0$ & $1$ & $2$ & $3$ & $4$ & $5$ & $6$ & $7$ & $8$ & $9$ & $10$ & $11$
& $12$ & $13$\\\hline
$f_{n}$ & $0$ & $1$ & $1$ & $2$ & $3$ & $5$ & $8$ & $13$ & $21$ & $34$ & $55$
& $89$ & $144$ & $233$\\\hline
\end{tabular}
\ \ \ \ .
\]

We proved a first property of Fibonacci numbers (Theorem
\ref{thm.fib.sum-first-n}) a while ago. In this section, we shall prove some
deeper properties of the Fibonacci sequence.

As a warm-up, we begin with two (inconsequential but neat) identities:

\begin{exercise}
\textbf{(a)} Prove that every nonnegative integer $n$ satisfies
\[
f_{1}+f_{3}+f_{5}+\cdots+f_{2n-1}=f_{2n}.
\]
(The left hand side is the sum of all $f_{2i-1}$ with $i\in\left\{
1,2,\ldots,n\right\}  $.) \medskip

\textbf{(b)} Prove that every nonnegative integer $n$ satisfies
\[
f_{0}+f_{2}+f_{4}+\cdots+f_{2n}=f_{2n+1}-1.
\]
(The left hand side is the sum of all $f_{2i}$ with $i\in\left\{
0,1,\ldots,n\right\}  $.)
\end{exercise}

\subsubsection{The addition theorem}

The next theorem is one of the most important properties of the Fibonacci sequence.

\begin{theorem}
[addition theorem for Fibonacci numbers]\label{thm.fib.add}We have%
\[
f_{n+m+1}=f_{n}f_{m}+f_{n+1}f_{m+1}\ \ \ \ \ \ \ \ \ \ \text{for all integers
}n,m\geq0.
\]

\end{theorem}

\begin{proof}
Can you induct on two variables at the same time? Not directly (although you
can induct on $n$ and then induct on $m$ in the induction step, so that you
have one induction proof inside another). Fortunately, we don't need to do
this here. It suffices to induct on one of the variables.

To be specific, let us induct on $n$. To that purpose, for every integer
$n\geq0$, we define the statement $P\left(  n\right)  $ to say%
\[
\text{\textquotedblleft for all integers }m\geq0\text{, we have }%
f_{n+m+1}=f_{n}f_{m}+f_{n+1}f_{m+1}\text{\textquotedblright.}%
\]
(Don't forget the \textquotedblleft for all integers $m\geq0$%
\textquotedblright\ part! The statement $P\left(  n\right)  $ is not just a
single equality $f_{n+m+1}=f_{n}f_{m}+f_{n+1}f_{m+1}$ for some specific value
of $m$, but rather combines infinitely many such equalities, one for each
integer $m\geq0$. If we fixed a value of $m$ and defined $P\left(  n\right)  $
to be just the single equality $f_{n+m+1}=f_{n}f_{m}+f_{n+1}f_{m+1}$ for this
particular value of $m$, then the induction proof below would not work,
because we are going to apply the induction hypothesis to a different $m$ than
we start with.)

We shall now prove this statement $P\left(  n\right)  $ for all $n\geq0$ by
induction on $n$.

\textit{Base case:} We must prove $P\left(  0\right)  $. In other words, we
must prove that%
\[
\text{\textquotedblleft for all integers }m\geq0\text{, we have }%
f_{0+m+1}=f_{0}f_{m}+f_{0+1}f_{m+1}\text{\textquotedblright.}%
\]
This is easy to show: For all integers $m\geq0$, we have $f_{0+m+1}=f_{m+1}$
and $\underbrace{f_{0}}_{=0}f_{m}+\underbrace{f_{0+1}}_{=f_{1}=1}%
f_{m+1}=0f_{m}+1f_{m+1}=f_{m+1}$, so the two sides are equal.

\textit{Induction step:} Let $n\geq0$ be an integer. We assume that $P\left(
n\right)  $ holds. We must show that $P\left(  n+1\right)  $ holds.

Our induction hypothesis says that $P\left(  n\right)  $ holds, i.e., that%
\[
\text{\textquotedblleft for all integers }m\geq0\text{, we have }%
f_{n+m+1}=f_{n}f_{m}+f_{n+1}f_{m+1}\text{\textquotedblright\ holds.}%
\]

We must prove that $P\left(  n+1\right)  $ holds, i.e., that%
\[
\text{\textquotedblleft for all integers }m\geq0\text{, we have }%
f_{n+1+m+1}=f_{n+1}f_{m}+f_{n+1+1}f_{m+1}\text{\textquotedblright\ holds.}%
\]

To prove this, we let $m\geq0$ be an integer. Then,%
\begin{align}
&  f_{n+1}f_{m}+\underbrace{f_{n+1+1}}_{\substack{=f_{n+2}\\=f_{n+1}%
+f_{n}\\\text{(by the recursive}\\\text{definition of the}\\\text{Fibonacci
numbers)}}}f_{m+1}\nonumber\\
&  =f_{n+1}f_{m}+\left(  f_{n+1}+f_{n}\right)  f_{m+1}\nonumber\\
&  =f_{n+1}f_{m}+f_{n+1}f_{m+1}+f_{n}f_{m+1}\nonumber\\
&  =f_{n+1}\underbrace{\left(  f_{m}+f_{m+1}\right)  }_{\substack{=f_{m+1}%
+f_{m}\\=f_{m+2}\\\text{(by the recursive}\\\text{definition of the}%
\\\text{Fibonacci numbers)}}}+\,f_{n}f_{m+1}\nonumber\\
&  =f_{n+1}f_{m+2}+f_{n}f_{m+1}=f_{n}f_{m+1}+f_{n+1}f_{m+2}.
\label{pf.thm.fib.add.maincalc}%
\end{align}

Now, recall that the induction hypothesis says that $P\left(  n\right)  $
holds, i.e., that%
\[
\text{\textquotedblleft for all integers }m\geq0\text{, we have }%
f_{n+m+1}=f_{n}f_{m}+f_{n+1}f_{m+1}\text{\textquotedblright\ holds.}%
\]
Note that the $m$ in this statement is a bound variable, i.e., it has nothing
to do with the $m$ that we have fixed; it just happens to have the same name.
Thus, we are free to apply our induction hypothesis $P\left(  n\right)  $ not
to the current $m$, but to any other $m$ as well. In particular, we can apply
it to $m+1$ instead of $m$. Thus, we obtain%
\[
f_{n+\left(  m+1\right)  +1}=f_{n}f_{m+1}+f_{n+1}f_{\left(  m+1\right)  +1}.
\]
\footnote{Let me explain this again in a slightly clearer (if longer) way.
\par
Our induction hypothesis tells us that
\[
\text{\textquotedblleft for all integers }m\geq0\text{, we have }%
f_{n+m+1}=f_{n}f_{m}+f_{n+1}f_{m+1}\text{\textquotedblright\ holds.}%
\]
We can rename the variable $m$ as $p$ in this statement (since it is just a
bound variable). Thus, we obtain that%
\[
\text{\textquotedblleft for all integers }p\geq0\text{, we have }%
f_{n+p+1}=f_{n}f_{p}+f_{n+1}f_{p+1}\text{\textquotedblright\ holds.}%
\]
Now, applying this latter statement to $p=m+1$ (where $m$ is the $m$ that we
fixed), we obtain%
\[
f_{n+\left(  m+1\right)  +1}=f_{n}f_{m+1}+f_{n+1}f_{\left(  m+1\right)  +1}.
\]
} This can be trivially simplified to%
\[
f_{n+m+2}=f_{n}f_{m+1}+f_{n+1}f_{m+2}.
\]

This equality has the same right hand side as (\ref{pf.thm.fib.add.maincalc}).
Thus, the left hand sides of the two equalities must be equal as well. In
other words, we must have%
\[
f_{n+m+2}=f_{n+1}f_{m}+f_{n+1+1}f_{m+1}.
\]
Since $n+m+2=n+1+m+1$, we can rewrite this as%
\[
f_{n+1+m+1}=f_{n+1}f_{m}+f_{n+1+1}f_{m+1}.
\]

Thus, we have proved that for all integers $m\geq0$, we have $f_{n+1+m+1}%
=f_{n+1}f_{m}+f_{n+1+1}f_{m+1}$. In other words, we have proved that $P\left(
n+1\right)  $ holds.

So the induction step is complete, and Theorem \ref{thm.fib.add} is proved.
\end{proof}

The next exercise gives two further properties of the Fibonacci sequence:

\begin{exercise}
\label{exe.23hw1.3}\textbf{(a)} Show that every positive integer $n$
satisfies
\[
f_{n+1}f_{n-1}-f_{n}^{2}=\left(  -1\right)  ^{n}.
\]
(The word \textquotedblleft Show\textquotedblright\ is a synonym for
\textquotedblleft Prove\textquotedblright.) \medskip

\textbf{(b)} Show that every nonnegative integer $n$ satisfies
\[
f_{1}^{2}+f_{2}^{2}+\cdots+f_{n}^{2}=f_{n}f_{n+1}.
\]
(The left hand side here is the sum of the squares of the first $n$ positive
Fibonacci numbers.)
\end{exercise}

The following exercise generalizes Theorem \ref{thm.fib.add} to a more general
class of recursively defined sequences:

\begin{exercise}
\label{exe.23hw1.4}Let $u$ and $v$ be two real numbers. Let $\left(
x_{0},x_{1},x_{2},\ldots\right)  $ be a sequence of real numbers such that
$x_{0}=0$ and $x_{1}=1$ and%
\[
x_{n}=ux_{n-1}+vx_{n-2}\ \ \ \ \ \ \ \ \ \ \text{for each }n\geq2.
\]
(When $u=1$ and $v=1$, this is the Fibonacci sequence.) Prove that%
\[
x_{n+m+1}=vx_{n}x_{m}+x_{n+1}x_{m+1}\ \ \ \ \ \ \ \ \ \ \text{for all integers
}n,m\geq0.
\]

\end{exercise}

\subsubsection{Divisibility of Fibonacci numbers}

Our next theorem involves divisibility of integers. We will study this in more
detail in Section \ref{sec.ent.div} (it is the fundamental concept of number
theory), but for now let me give its definition:

\begin{definition}
\label{def.divides}Let $a$ and $b$ be two integers. We say that $a$
\textbf{divides }$b$ (and we write $a\mid b$) if there exists an integer $c$
such that $b=ac$. Equivalently, we say that $b$ is \textbf{divisible by }$a$
in this case.
\end{definition}

For example, we have $2\mid4$ and $3\mid12$ and $10\mid30$ and $0\mid0$ and
$5\mid0$. But we don't have $2\mid3$ or $0\mid1$. The integer $0$ is divisible
by every integer, but only divides itself.

Now we can state a divisibility property of Fibonacci numbers:

\begin{theorem}
\label{thm.fib.div}If $a,b\geq0$ are two integers that satisfy $a\mid b$, then
$f_{a}\mid f_{b}$.
\end{theorem}

In other words, in our above table of Fibonacci numbers, if some entry of the
first row divides some other entry of the first row, then the same holds for
the corresponding entries of the second row. For example, $6\mid12$ implies
$f_{6}\mid f_{12}$ (which is saying that $8\mid144$).

\begin{proof}
[Proof of Theorem \ref{thm.fib.div}.]It is reasonable to try induction.
However, inducting on $a$ does not lead anywhere: The base case is easy, but
in the induction step it is completely unclear how to reach the goal, since
the condition $a\mid b$ in the induction hypothesis usually has nothing to do
with the condition $a+1\mid b$ in the induction goal.

Similar problems appear if you try to induct on $b$. So neither of the two
variables in the theorem is suitable for being inducted on.

What can we do? Give up on induction?

Not so fast. One thing we haven't tried is to introduce a new variable and
then induct on that new variable.

To do so, we observe that two integers $a,b\geq0$ satisfy $a\mid b$ if and
only if there exists an integer $c$ such that $b=ac$ (by the definition of
\textquotedblleft divides\textquotedblright). Moreover, if this integer $c$
exists, then it can be chosen to be $\geq0$ (this is automatic when $b\neq0$,
because $c=\dfrac{b}{a}>0$ in this case; but otherwise we can achieve this by
simply choosing $c=0$). Thus, two integers $a,b\geq0$ satisfy $a\mid b$ if and
only if there exists an integer $c\geq0$ such that $b=ac$.

Hence, a pair of integers $a,b\geq0$ satisfying $a\mid b$ is nothing but a
pair of the form $a,ac$ where $a,c\geq0$ are integers. This allows us to
restate Theorem \ref{thm.fib.div} as follows:

\begin{statement}
\textit{Restated theorem:} \textquotedblleft For any integers $a,c\geq0$, we
have $f_{a}\mid f_{ac}$.\textquotedblright
\end{statement}

Now, we shall prove this restated theorem by induction on $c$. In other words,
for each $c\geq0$, we shall prove the statement%
\[
P\left(  c\right)  :=\left(  \text{\textquotedblleft for any integer }%
a\geq0\text{, we have }f_{a}\mid f_{ac}\text{\textquotedblright}\right)  .
\]

\textit{Base case:} We must prove $P\left(  0\right)  $. In other words, we
must prove that
\[
\text{\textquotedblleft for any integer }a\geq0\text{, we have }f_{a}\mid
f_{a\cdot0}\text{\textquotedblright.}%
\]
But this is easy, because for any integer $a\geq0$, we have $f_{a\cdot0}%
=f_{0}=0$, which is divisible by any integer (thus in particular by $f_{a}$).

\textit{Induction step:} Let $c\geq0$ be an integer. We assume that $P\left(
c\right)  $ holds, i.e., that
\[
\text{\textquotedblleft for any integer }a\geq0\text{, we have }f_{a}\mid
f_{ac}\text{\textquotedblright\ holds.}%
\]
We must prove that $P\left(  c+1\right)  $ holds, i.e., that%
\[
\text{\textquotedblleft for any integer }a\geq0\text{, we have }f_{a}\mid
f_{a\left(  c+1\right)  }\text{\textquotedblright\ holds.}%
\]

Let $a\geq0$ be any integer. Then, the induction hypothesis (i.e., our
assumption that $P\left(  c\right)  $ holds) yields that $f_{a}\mid f_{ac}$.
In other words, $f_{ac}=f_{a}p$ for some integer $p$. Now,%
\begin{align*}
f_{a\left(  c+1\right)  }  &  =f_{ac+a}=f_{ac+\left(  a-1\right)  +1}\\
&  =\underbrace{f_{ac}}_{=f_{a}p}f_{a-1}+f_{ac+1}f_{a}%
\ \ \ \ \ \ \ \ \ \ \left(
\begin{array}
[c]{c}%
\text{by Theorem \ref{thm.fib.add},}\\
\text{applied to }n=ac\text{ and }m=a-1
\end{array}
\right) \\
&  =f_{a}pf_{a-1}+f_{ac+1}f_{a}=f_{a}\cdot\underbrace{\left(  pf_{a-1}%
+f_{ac+1}\right)  }_{\text{an integer}}.
\end{align*}
This immediately yields that $f_{a}\mid f_{a\left(  c+1\right)  }$. Thus, we
have shown that for any integer $a\geq0$, we have $f_{a}\mid f_{a\left(
c+1\right)  }$. In other words, we have proved that $P\left(  c+1\right)  $
holds. This completes the induction step, and thus the restated theorem is
proved. Therefore, the original Theorem \ref{thm.fib.div} is also proved.
\medskip

$\ldots\ldots\ldots\ldots\ldots$

Is it? There is a subtle gap in our above argument. Can you find it? \medskip

$\ldots\ldots\ldots\ldots\ldots$

Can you? Don't look down just yet. The gap is somewhere above! \medskip

$\ldots\ldots\ldots\ldots\ldots$

This time, the theorem itself is correct, so you can't find the gap by tracing
the proof through a case where the theorem is false. Though an example might
be useful... \medskip

$\ldots\ldots\ldots\ldots\ldots$

No, we didn't misuse the principle of induction. The structure of the proof is
fine. (Actually, we could have made our statements a bit shorter by fixing
$a\geq0$, but this wouldn't have made much of a difference.) \medskip

$\ldots\ldots\ldots\ldots\ldots$

The base case was fine, too. \medskip

$\ldots\ldots\ldots\ldots\ldots$

A computer, of course, would spot the problem.

If you tried to formalize the above proof in a computer language (e.g.,
\href{https://en.wikipedia.org/wiki/Coq}{Coq} or
\href{https://en.wikipedia.org/wiki/Lean_(proof_assistant)}{Lean}), you would
run into a type mismatch error. Some statement has been proved for variables
of a certain type, but is being used for variables of a different type. Very
slightly different. \medskip

$\ldots\ldots\ldots\ldots\ldots$

The statement in question is Theorem \ref{thm.fib.add}. It is stated for one
kind of variables, but we have used it for a slightly different kind. \medskip

$\ldots\ldots\ldots\ldots\ldots$

OK, I am spelling it out: Theorem \ref{thm.fib.add} (i.e., the addition
formula $f_{n+m+1}=f_{n}f_{m}+f_{n+1}f_{m+1}$) has been stated and proved for
all integers $n,m\geq0$, but we have applied it to $n=ac$ and $m=a-1$. For
this to work, we need $ac\geq0$ and $a-1\geq0$. Now, $ac\geq0$ is indeed
satisfied (since $a\geq0$ and $c\geq0$), but $a-1\geq0$ holds only if $a\geq
1$, which is not guaranteed. Thus, our use of Theorem \ref{thm.fib.add} was
illegal when $a=0$. And indeed, if we apply Theorem \ref{thm.fib.add} for
$a=0$, then we end up with an $f_{-1}$ term, which is undefined. Even if you
define $f_{-1}$ appropriately (and there is a good definition; see Subsection
\ref{subsec.ind.moreex.f-n}), we have not proved Theorem \ref{thm.fib.add} for
negative $n,m$. So there is a gap in our proof. Can we fix it?\medskip

$\ldots\ldots\ldots\ldots\ldots$

Fortunately, we can: Our argument breaks down only in the case when $a=0$, and
we can just treat this case $a=0$ manually, since it is an easy case. So we
build a case distinction into our above induction step. Thus, the induction
step takes the following form: \medskip

\textit{Induction step (corrected):} Let $c\geq0$ be an integer. We assume
that $P\left(  c\right)  $ holds, i.e., that
\[
\text{\textquotedblleft for any integer }a\geq0\text{, we have }f_{a}\mid
f_{ac}\text{\textquotedblright\ holds.}%
\]
We must prove that $P\left(  c+1\right)  $ holds, i.e., that%
\[
\text{\textquotedblleft for any integer }a\geq0\text{, we have }f_{a}\mid
f_{a\left(  c+1\right)  }\text{\textquotedblright\ holds.}%
\]

Let $a\geq0$ be any integer. We must show that $f_{a}\mid f_{a\left(
c+1\right)  }$. We are in one of the following two cases:

\textit{Case 1:} We have $a=0$.

\textit{Case 2:} We have $a\neq0$.

In Case 1, we have $a=0$, so that both $f_{a}$ and $f_{a\left(  c+1\right)  }$
equal $f_{0}=0$, and thus $f_{a}\mid f_{a\left(  c+1\right)  }$ holds (since
$0\mid0$). Thus, the divisibility $f_{a}\mid f_{a\left(  c+1\right)  }$ is
proved in Case 1.

Now, consider Case 2. In this case, $a\neq0$, so that $a\geq1$ (because $a$ is
an integer and $\geq0$).\ Hence, $a-1\geq0$. This will allow us to apply
Theorem \ref{thm.fib.add} to $n=ac$ and $m=a-1$ in a few moments. The
induction hypothesis (i.e., our assumption that $P\left(  c\right)  $ holds)
yields that $f_{a}\mid f_{ac}$. In other words, $f_{ac}=f_{a}p$ for some
integer $p$. Now,%
\begin{align*}
f_{a\left(  c+1\right)  }  &  =f_{ac+a}=f_{ac+\left(  a-1\right)  +1}\\
&  =\underbrace{f_{ac}}_{=f_{a}p}f_{a-1}+f_{ac+1}f_{a}%
\ \ \ \ \ \ \ \ \ \ \left(
\begin{array}
[c]{c}%
\text{by Theorem \ref{thm.fib.add},}\\
\text{applied to }n=ac\text{ and }m=a-1
\end{array}
\right) \\
&  =f_{a}pf_{a-1}+f_{ac+1}f_{a}=f_{a}\cdot\underbrace{\left(  pf_{a-1}%
+f_{ac+1}\right)  }_{\text{an integer}}.
\end{align*}
This immediately yields that $f_{a}\mid f_{a\left(  c+1\right)  }$.

So we have proved $f_{a}\mid f_{a\left(  c+1\right)  }$ in both Cases 1 and 2.
Therefore, $f_{a}\mid f_{a\left(  c+1\right)  }$ always holds.

Thus, $P\left(  c+1\right)  $ is proved. This completes the induction step,
and thus the restated theorem is proved. Therefore, Theorem \ref{thm.fib.div}
is proved -- correctly this time!
\end{proof}

\subsubsection{Binet's formula}

Is there an explicit formula for $f_{n}$, that is, a formula that does not
rely on the previous entries of the Fibonacci sequence?

Yes, there is one; it is known as \textbf{Binet's formula}:

\begin{theorem}
[Binet's formula]\label{thm.fib.binet}Let%
\[
\varphi=\dfrac{1+\sqrt{5}}{2}\approx1.618\ldots\ \ \ \ \ \ \ \ \ \ \text{and}%
\ \ \ \ \ \ \ \ \ \ \psi=\dfrac{1-\sqrt{5}}{2}\approx-0.618\ldots.
\]
Then,%
\[
f_{n}=\dfrac{\varphi^{n}-\psi^{n}}{\sqrt{5}}\ \ \ \ \ \ \ \ \ \ \text{for
every integer }n\geq0.
\]

\end{theorem}

Some remarks:

\begin{itemize}
\item The number $\varphi$ is called the
\textbf{\href{https://en.wikipedia.org/wiki/Golden_ratio}{\textbf{golden
ratio}}}, and is famous for many properties, including the fact that
$\varphi^{2}=\varphi+1$ (which you can easily check by expanding both
sides\footnote{Namely: From $\varphi=\dfrac{1+\sqrt{5}}{2}$, we obtain
\[
\varphi^{2}=\left(  \dfrac{1+\sqrt{5}}{2}\right)  ^{2}=\dfrac{1+2\sqrt{5}%
+5}{4}=\dfrac{6+2\sqrt{5}}{4}=\dfrac{3+\sqrt{5}}{2}=\dfrac{1+\sqrt{5}}%
{2}+1=\varphi+1.
\]
}). The number $\psi$ is its so-called conjugate and also satisfies $\psi
^{2}=\psi+1$.

\item The numbers $f_{n}$ are integers, but Binet's formula expresses them in
terms of two irrational numbers $\varphi$ and $\psi$. This should be rather unexpected.

\item As $n$ grows large, $\psi^{n}$ approaches $0$ (since $-1<\psi<1$),
whereas $\varphi^{n}$ grows exponentially (since $\varphi>1$). So $f_{n}$ also
grows exponentially (according to Binet's formula), with growth rate
$\varphi\approx1.618\ldots$.
\end{itemize}

Two questions arise:

\begin{enumerate}
\item How do we prove Binet's formula?

\item How could we find Binet's formula if we didn't already know it?
\end{enumerate}

We will answer Question 1 soon. Question 2 is significantly trickier and will
not be answered in this course\footnote{Answers at different levels of
generality can be found in:
\par
\begin{itemize}
\item \cite[Subsection 4.9.2]{20f} (which solves any linear recurrence of the
form $x_{n}=ax_{n-1}+bx_{n-2}$ for constant numbers $a$ and $b$ in an explicit
and elementary way);
\par
\item \cite{Melian16} and \cite{Ivanov14} (which solve the more general
version $x_{n}=a_{1}x_{n-1}+a_{2}x_{n-2}+\cdots+a_{k}x_{n-k}$ in terms of the
eigenvalues of a matrix).
\end{itemize}
\par
Textbooks on combinatorics or advanced linear algebra also tend to discuss
such sequences (called \textbf{linearly recurrent sequences}).}.

\bigskip

Let us try to prove Binet's formula by induction on $n$:

\begin{proof}
[Attempted proof of Binet's formula.]We induct on $n$:

\textit{Base case:} For $n=0$, we have $f_{n}=f_{0}=0$ and
\[
\dfrac{\varphi^{n}-\psi^{n}}{\sqrt{5}}=\dfrac{\varphi^{0}-\psi^{0}}{\sqrt{5}%
}=\dfrac{1-1}{\sqrt{5}}=0.
\]
Thus, Binet's formula holds for $n=0$.

\textit{Induction step:} Let $n\geq0$ be an integer.

Assume (as induction hypothesis) that Binet's formula holds for $n$; we must
prove that it holds for $n+1$.

So we must prove that%
\[
f_{n+1}=\dfrac{\varphi^{n+1}-\psi^{n+1}}{\sqrt{5}}.
\]
The recursive definition of the Fibonacci sequence yields%
\[
f_{n+1}=f_{n}+f_{n-1}=\dfrac{\varphi^{n}-\psi^{n}}{\sqrt{5}}+f_{n-1}%
\ \ \ \ \ \ \ \ \ \ \left(  \text{by the induction hypothesis}\right)  .
\]
So far so good, but how can we simplify $f_{n-1}$ ? Our induction hypothesis
only tells us that $f_{n}=\dfrac{\varphi^{n}-\psi^{n}}{\sqrt{5}}$, but it says
nothing about $f_{n-1}$.
\end{proof}

So this induction proof does not work.\footnote{There is also one more little
(fixable) gap in the above attempted proof. Do you see it?}

Let us see how to fix this by introducing a more advanced version of induction.

\subsection{\label{sec.ind.strind}Strong induction}

\subsubsection{Reminder on regular induction}

Recall the (original) principle of mathematical induction:

\begin{theorem}
[Principle of Mathematical Induction]Let $b$ be an integer.

Let $P\left(  n\right)  $ be a mathematical statement defined for each integer
$n\geq b$.

Assume the following:

\begin{enumerate}
\item \textquotedblleft\textbf{Base case}\textquotedblright: The statement
$P\left(  b\right)  $ holds.

\item \textquotedblleft\textbf{Induction step}\textquotedblright: For each
integer $n\geq b$, the implication $P\left(  n\right)  \Longrightarrow
P\left(  n+1\right)  $ holds.
\end{enumerate}

Then, the statement $P\left(  n\right)  $ holds for every integer $n\geq b$.
\end{theorem}

We can restate this principle slightly by renaming the $n$ in the induction
step as $n-1$ (so that the implication $P\left(  n\right)  \Longrightarrow
P\left(  n+1\right)  $ turns into $P\left(  n-1\right)  \Longrightarrow
P\left(  n\right)  $). Thus, it takes the following form:

\begin{theorem}
[Principle of Mathematical Induction, restated]Let $b$ be an integer.

Let $P\left(  n\right)  $ be a mathematical statement defined for each integer
$n\geq b$.

Assume the following:

\begin{enumerate}
\item \textquotedblleft\textbf{Base case}\textquotedblright: The statement
$P\left(  b\right)  $ holds.

\item \textquotedblleft\textbf{Induction step}\textquotedblright: For each
integer $n>b$, the implication $P\left(  n-1\right)  \Longrightarrow P\left(
n\right)  $ holds.
\end{enumerate}

Then, the statement $P\left(  n\right)  $ holds for every integer $n\geq b$.
\end{theorem}

The idea behind the principle (in either form) is that the base case gives us
$P\left(  b\right)  $ whereas the induction step gives us the implications%
\begin{align*}
P\left(  b\right)   &  \Longrightarrow P\left(  b+1\right)  ,\\
P\left(  b+1\right)   &  \Longrightarrow P\left(  b+2\right)  ,\\
P\left(  b+2\right)   &  \Longrightarrow P\left(  b+3\right)  ,\\
&  \ldots.
\end{align*}
In the domino metaphor (see Remark \ref{rmk.ind.domino-metaphor}), the base
case tips over the first domino, and the induction step ensures that each
domino falls from the impact of the previous domino's falling.

\subsubsection{\label{subsec.ind.strind.strind}Strong induction}

Now, assume that the $\left(  b+2\right)  $-domino (i.e., $P\left(
b+2\right)  $) falls not from the impact of the previous domino $P\left(
b+1\right)  $, but rather from the combined force of the dominos $P\left(
b\right)  $ and $P\left(  b+1\right)  $. This would still suffice, because the
latter two dominos have already fallen. In other words, instead of the
implication $P\left(  b+1\right)  \Longrightarrow P\left(  b+2\right)  $, we
could just as well prove the implication%
\[
\left(  P\left(  b\right)  \text{ AND }P\left(  b+1\right)  \right)
\Longrightarrow P\left(  b+2\right)  ,
\]
which is somewhat weaker (since it assumes more to get to the same conclusion)
but nevertheless gives the same result. Likewise, we could just as well
replace the implication $P\left(  b+2\right)  \Longrightarrow P\left(
b+3\right)  $ by the weaker implication%
\[
\left(  P\left(  b\right)  \text{ AND }P\left(  b+1\right)  \text{ AND
}P\left(  b+2\right)  \right)  \Longrightarrow P\left(  b+3\right)  .
\]
More generally, for each $n>b$, instead of proving the implication $P\left(
n-1\right)  \Longrightarrow P\left(  n\right)  $, it will suffice to prove the
weaker implication%
\[
\underbrace{\left(  P\left(  b\right)  \text{ AND }P\left(  b+1\right)  \text{
AND }P\left(  b+2\right)  \text{ AND }\cdots\text{ AND }P\left(  n-1\right)
\right)  }_{\text{i.e., the statement }P\left(  k\right)  \text{ holds for
each }k\in\left\{  b,b+1,\ldots,n-1\right\}  }\Longrightarrow P\left(
n\right)
\]
(so that the domino $P\left(  n\right)  $ is tipped over by the combined force
of all the preceding dominos, not just the one domino directly to its left).

This induction principle is called \textbf{strong induction}. Explicitly, it
says the following:

\begin{theorem}
[Principle of Strong Induction]Let $b$ be an integer.

Let $P\left(  n\right)  $ be a mathematical statement defined for each integer
$n\geq b$.

Assume the following:

\begin{enumerate}
\item \textquotedblleft\textbf{Base case}\textquotedblright: The statement
$P\left(  b\right)  $ holds.

\item \textquotedblleft\textbf{Induction step}\textquotedblright: For each
integer $n>b$, the implication
\[
\left(  P\left(  b\right)  \text{ AND }P\left(  b+1\right)  \text{ AND
}P\left(  b+2\right)  \text{ AND }\cdots\text{ AND }P\left(  n-1\right)
\right)  \Longrightarrow P\left(  n\right)
\]
holds.
\end{enumerate}

Then, the statement $P\left(  n\right)  $ holds for every integer $n\geq b$.
\end{theorem}

Proofs using this principle are called \textbf{proofs by strong induction} (or
\textbf{strong induction proofs}). They differ from proofs by (regular)
induction as follows: In the induction step of a strong induction proof, you
can use not just the preceding statement $P\left(  n-1\right)  $, but also all
the statements before it ($P\left(  n-2\right)  $ and $P\left(  n-3\right)  $
and so on, all the way down to $P\left(  b\right)  $). In other words, the
induction hypothesis is now stronger (thus the name \textquotedblleft strong
induction\textquotedblright). Roughly speaking, strong induction is
\textquotedblleft induction with a long memory\textquotedblright\ (as opposed
to regular induction, whose memory only is $1$ step long).\footnote{A remark
for the logically inclined:
\par
Surprisingly, the Principle of Strong Induction is logically equivalent to the
regular Principle of Mathematical Induction (i.e., each of the two principles
can be derived from the other). Thus, we don't need to assume the former as an
extra axiom (once we have assumed the latter). See \cite[\S 2.8.1]{detnotes}
for how the former can be derived from the latter.
\par
In essence, this means that strong induction is just a \textquotedblleft more
convenient user interface\textquotedblright\ for regular induction; everything
that can be proved using strong induction can still be proved using regular
induction. (But it requires a little trick: If you can prove $P\left(
n\right)  $ by strong induction on $n$, then you can prove the statement%
\[
Q\left(  n\right)  :=\left(  P\left(  b\right)  \text{ AND }P\left(
b+1\right)  \text{ AND }P\left(  b+2\right)  \text{ AND }\cdots\text{
AND\ }P\left(  n\right)  \right)
\]
by regular induction on $n$, and then you can derive $P\left(  n\right)  $
from $Q\left(  n\right)  $.)}

(We will later see a slightly nicer form of strong induction, in which the
base case is incorporated in the induction step.)

\bigskip

Before we see an example of a strong induction proof, let me explain why it
works. Let's say you have proved a statement $P\left(  n\right)  $ for all
$n\geq0$ by strong induction. Thus,

\begin{itemize}
\item you have proved $P\left(  0\right)  $ (this is the base case);

\item you have proved the implication $P\left(  0\right)  \Longrightarrow
P\left(  1\right)  $ (this is the induction step for $n=1$), so you conclude
that $P\left(  1\right)  $ holds (since $P\left(  0\right)  $ holds);

\item you have proved the implication $\left(  P\left(  0\right)  \text{ AND
}P\left(  1\right)  \right)  \Longrightarrow P\left(  2\right)  $ (this is the
induction step for $n=2$), so you conclude that $P\left(  2\right)  $ holds
(since $P\left(  0\right)  $ and $P\left(  1\right)  $ hold);

\item you have proved the implication $\left(  P\left(  0\right)  \text{ AND
}P\left(  1\right)  \text{ AND }P\left(  2\right)  \right)  \Longrightarrow
P\left(  3\right)  $ (this is the induction step for $n=3$), so you can
conclude that $P\left(  3\right)  $ holds (since $P\left(  0\right)  $ and
$P\left(  1\right)  $ and $P\left(  2\right)  $ hold);

\item and so on.
\end{itemize}

\subsubsection{\label{subsec.ind.strind.binet}Example: Proof of Binet's
formula}

Let us now prove Binet's formula by strong induction:

\begin{proof}
[Proof of Theorem \ref{thm.fib.binet} (i.e., of Binet's formula).]We strongly
induct on $n$ (i.e., we use strong induction on $n$). That is, we let
$P\left(  n\right)  $ denote the statement
\[
\left(  \text{\textquotedblleft}f_{n}=\dfrac{\varphi^{n}-\psi^{n}}{\sqrt{5}%
}\text{\textquotedblright}\right)
\]
for each $n\geq0$, and we apply the Principle of Strong Induction (for $b=0$)
to prove this statement $P\left(  n\right)  $ for each $n\geq0$.

\textit{Base case:} As above, we check that Binet's formula (i.e., the
statement $P\left(  n\right)  $) holds for $n=0$.

\textit{Induction step:} Let $n>0$ be an integer. We must prove the
implication%
\[
\left(  P\left(  0\right)  \text{ AND }P\left(  1\right)  \text{ AND }P\left(
2\right)  \text{ AND }\cdots\text{ AND }P\left(  n-1\right)  \right)
\Longrightarrow P\left(  n\right)  .
\]
Thus, we assume that $P\left(  0\right)  $ AND $P\left(  1\right)  $ AND
$P\left(  2\right)  $ AND $\cdots$ AND $P\left(  n-1\right)  $ holds. In other
words, we assume that Binet's formula holds for $0$, for $1$, for $2$, and so
on, all the way up to $n-1$. (In other words, we assume that $f_{k}%
=\dfrac{\varphi^{k}-\psi^{k}}{\sqrt{5}}$ for each $k\in\left\{  0,1,\ldots
,n-1\right\}  $.)

We have to prove $P\left(  n\right)  $. In other words, we have to prove that
Binet's formula also holds for $n$. In other words, we have to prove that
$f_{n}=\dfrac{\varphi^{n}-\psi^{n}}{\sqrt{5}}$.

We assumed that Binet's formula holds for $n-1$. That is, we have
$f_{n-1}=\dfrac{\varphi^{n-1}-\psi^{n-1}}{\sqrt{5}}$.

We assumed that Binet's formula holds for $n-2$. That is, we have
$f_{n-2}=\dfrac{\varphi^{n-2}-\psi^{n-2}}{\sqrt{5}}$.

As we have seen above, we have $\varphi^{2}=\varphi+1$ and $\psi^{2}=\psi+1$.

But the recursive definition of the Fibonacci sequence yields%
\begin{align*}
f_{n}  &  =f_{n-1}+f_{n-2}=\dfrac{\varphi^{n-1}-\psi^{n-1}}{\sqrt{5}}%
+\dfrac{\varphi^{n-2}-\psi^{n-2}}{\sqrt{5}}\\
&  \ \ \ \ \ \ \ \ \ \ \ \ \ \ \ \ \ \ \ \ \left(  \text{since }f_{n-1}%
=\dfrac{\varphi^{n-1}-\psi^{n-1}}{\sqrt{5}}\text{ and }f_{n-2}=\dfrac
{\varphi^{n-2}-\psi^{n-2}}{\sqrt{5}}\right) \\
&  =\dfrac{1}{\sqrt{5}}\left(  \varphi^{n-1}-\psi^{n-1}+\varphi^{n-2}%
-\psi^{n-2}\right) \\
&  =\dfrac{1}{\sqrt{5}}\left(  \underbrace{\varphi^{n-1}+\varphi^{n-2}%
}_{=\varphi^{n-2}\left(  \varphi+1\right)  }-\underbrace{\left(  \psi
^{n-1}+\psi^{n-2}\right)  }_{=\psi^{n-2}\left(  \psi+1\right)  }\right) \\
&  =\dfrac{1}{\sqrt{5}}\left(  \varphi^{n-2}\underbrace{\left(  \varphi
+1\right)  }_{=\varphi^{2}}-\,\psi^{n-2}\underbrace{\left(  \psi+1\right)
}_{=\psi^{2}}\right) \\
&  =\dfrac{1}{\sqrt{5}}\left(  \underbrace{\varphi^{n-2}\varphi^{2}}%
_{=\varphi^{n}}-\underbrace{\psi^{n-2}\psi^{2}}_{=\psi^{n}}\right)  =\dfrac
{1}{\sqrt{5}}\left(  \varphi^{n}-\psi^{n}\right)  =\dfrac{\varphi^{n}-\psi
^{n}}{\sqrt{5}}.
\end{align*}
So we have proved Binet's formula for $n$. Right?

$\ldots\ldots\ldots\ldots\ldots\ldots\ldots\ldots\ldots\ldots$

Wait a moment! We have assumed (as the induction hypothesis) that Binet's
formula holds for each of the numbers $0,1,\ldots,n-1$. But then we have used
it for $n-2$ and for $n-1$. This tacitly relied on the fact that $n-2$ and
$n-1$ are among the numbers $0,1,\ldots,n-1$. However, this fact is only true
if $n\geq2$. If $n=1$, then $n-2$ is not among the numbers $0,1,\ldots,n-1$
(because it is negative).\footnote{The particularly astute reader might notice
that the recursion $f_{n}=f_{n-1}+f_{n-2}$ also doesn't work for $n=1$, since
$f_{-1}$ is undefined -- at least with the definition of Fibonacci numbers we
gave above. Later on, in Subsection \ref{subsec.ind.moreex.f-n}, we shall
\textquotedblleft fix\textquotedblright\ this by extending the Fibonacci
sequence to the left (defining $f_{-1},f_{-2},f_{-3},\ldots$). But this still
does not help our induction step.}

So our induction step worked for $n=2,3,4,\ldots$ but not for $n=1$. What can
we do?

We can fix this by just proving the claim for $n=1$ by hand. So we must prove
that $f_{1}=\dfrac{\varphi^{1}-\psi^{1}}{\sqrt{5}}$. This can be checked by a
direct computation:%
\[
\dfrac{\varphi^{1}-\psi^{1}}{\sqrt{5}}=\dfrac{\varphi-\psi}{\sqrt{5}}%
=\dfrac{\dfrac{1+\sqrt{5}}{2}-\dfrac{1-\sqrt{5}}{2}}{\sqrt{5}}=\dfrac{\sqrt
{5}}{\sqrt{5}}=1=f_{1}.
\]

Now our induction step is really complete, and Binet's formula is proved.
\end{proof}

\bigskip

Let us summarize: We have used strong induction in our above proof of Theorem
\ref{thm.fib.binet}, because the \textquotedblleft extra
memory\textquotedblright\ in a strong induction step allowed us to express not
just $f_{n-1}$ but also $f_{n-2}$ via the induction hypothesis.

Note that we have had to handle the two cases $n=0$ and $n=1$ by hand in our
above proof, because we had to reach \textquotedblleft$2$ steps
back\textquotedblright\ in memory in the induction step (i.e., we had to apply
the induction hypothesis both to $n-1$ and to $n-2$).\ \ \ \ \footnote{Had we
reached further back, we would have needed extra cases (e.g., if we had
applied the induction hypothesis to $n-5$, then we would have to handle all
the cases $n=0,1,2,3,4$ by hand).} The case $n=0$ was our base case, whereas
the case $n=1$ was part of the induction step, but nevertheless had to be
singled out for special treatment (since $n-2$ is negative for $n=1$).
Nevertheless, it makes sense to think of the $n=1$ case as a \textquotedblleft
second base case\textquotedblright, even if it is de-jure part of the
induction step.

\subsubsection{\label{subsec.ind.strind.baseless}Baseless strong induction}

You can actually reformulate the principle of strong induction in a form that
does not have a de-jure base case at all:

\begin{theorem}
[Principle of Strong Induction, restated]Let $b$ be an integer.

Let $P\left(  n\right)  $ be a mathematical statement defined for each integer
$n\geq b$.

Assume the following:

\begin{itemize}
\item \textquotedblleft\textbf{Induction step\textquotedblright:}\ For each
integer $n\geq b$, the implication
\[
\left(  P\left(  b\right)  \text{ AND }P\left(  b+1\right)  \text{ AND
}P\left(  b+2\right)  \text{ AND }\cdots\text{ AND }P\left(  n-1\right)
\right)  \Longrightarrow P\left(  n\right)
\]
holds.
\end{itemize}

Then, the statement $P\left(  n\right)  $ holds for every integer $n\geq b$.
\end{theorem}

How does this restated principle work without a base case? Easy: We have just
repackaged the base case into the induction step. Indeed, note that the
induction step now says \textquotedblleft$n\geq b$\textquotedblright, not
\textquotedblleft$n>b$\textquotedblright. In particular, this means that the
implication%
\[
\left(  P\left(  b\right)  \text{ AND }P\left(  b+1\right)  \text{ AND
}P\left(  b+2\right)  \text{ AND }\cdots\text{ AND }P\left(  n-1\right)
\right)  \Longrightarrow P\left(  n\right)
\]
has to hold for $n=b$. However, for $n=b$, the antecedent (= if-part) of this
implication is a tautology (i.e., is an empty statement that is automatically
true by dint of its emptiness\footnote{Don't believe it? Observe that this
antecedent%
\[
\left(  P\left(  b\right)  \text{ AND }P\left(  b+1\right)  \text{ AND
}P\left(  b+2\right)  \text{ AND }\cdots\text{ AND }P\left(  n-1\right)
\right)
\]
is a conjunction of $n-b$ statements (since there are $n-b$ numbers between
$b$ and $n-1$ inclusive). If $n=b$, this means that it is a conjunction of
$b-b=0$ statements, i.e., of no statements whatsoever. So it is an empty
statement, automatically true.}), and thus proving this implication is
tantamount to just unconditionally proving $P\left(  b\right)  $, which was
what we previously viewed as our base case. So we have not magically removed
the need for a base case; we just have merged it into the induction step.
Nevertheless, this makes for a slightly cleaner version of strong induction.

\begin{noncompile}
So a strong induction needs no base case.

In practice, in the induction step, some values of $n$ will require extra
treatment, because there are not enough induction hypotheses in these cases;
thus, we might end up doing de-facto base cases (often more than one) inside
our induction step.
\end{noncompile}

\subsubsection{Example: Prime factorizations exist}

Another example of a strong induction proof comes from elementary number
theory. We recall two basic definitions (more on this later, when we cover
number theory):

\begin{definition}
Let $b$ be an integer. A \textbf{divisor} of $b$ means an integer $a$
satisfying $a\mid b$.
\end{definition}

For example, the divisors of $6$ are $1,2,3,6,-1,-2,-3,-6$.

\begin{definition}
A \textbf{prime} (or \textbf{prime number}) means an integer $p>1$ whose only
positive divisors are $1$ and $p$.
\end{definition}

So the primes (in increasing order) are%
\[
2,\ \ 3,\ \ 5,\ \ 7,\ \ 11,\ \ 13,\ \ 17,\ \ 19,\ \ 23,\ \ 29,\ \ \ldots.
\]
There are infinitely many primes, as we will show later.

\begin{theorem}
\label{thm.ent.primefac.exists}Every positive integer is a product of finitely
many primes.
\end{theorem}

Here and in the following, I understand an empty product (i.e., a product of
no numbers whatsoever) to be $1$. Thus, Theorem \ref{thm.ent.primefac.exists}
does hold for $1$, since $1$ is a product of no primes.

Here are more interesting examples:

\begin{itemize}
\item $2023=7\cdot17\cdot17$ is a product of three primes.

\item $2024=2\cdot2\cdot2\cdot11\cdot23$ is a product of five primes.

\item $2=2$ is a product of one prime (namely, $2$ itself).
\end{itemize}

How do we prove Theorem \ref{thm.ent.primefac.exists} in general?

\begin{proof}
[Proof of Theorem \ref{thm.ent.primefac.exists}.]We must prove the statement%
\[
P\left(  n\right)  =\left(  \text{\textquotedblleft}n\text{ is a product of
finitely many primes\textquotedblright}\right)
\]
for each integer $n\geq1$.

We shall prove this by strong induction on $n$. (We use the original variant
of strong induction, with a base case.)

\textit{Base case:} $P\left(  1\right)  $ is true, since $1$ is a product of
finitely many primes (specifically, of $0$ primes, as we saw).

\textit{Induction step:} Let $n>1$. We must prove the implication%
\[
\left(  P\left(  1\right)  \text{ AND }P\left(  2\right)  \text{ AND }%
\cdots\text{ AND }P\left(  n-1\right)  \right)  \Longrightarrow P\left(
n\right)  .
\]

So we assume that $P\left(  1\right)  $ AND $P\left(  2\right)  $
AND\ $\cdots$ AND $P\left(  n-1\right)  $ holds. We must prove that $P\left(
n\right)  $ holds.

In other words, we must prove that $n$ is a product of finitely many primes.

We are in one of the following two cases:

\textit{Case 1:} The only positive divisors of $n$ are $1$ and $n$.

\textit{Case 2:} There is a positive divisor $d$ of $n$ that is neither $1$
nor $n$.

(Other cases are not possible, since $1$ and $n$ always are positive divisors
of $n$.)

Consider Case 1 first. In this case, $n$ itself is a prime (by the definition
of a prime), and thus is a product of finitely many primes (namely, of just
$1$ prime: itself). Thus, $P\left(  n\right)  $ holds in Case 1.

Now, consider Case 2. In this case, there is a positive divisor $d$ of $n$
that is neither $1$ nor $n$. Consider such a $d$ (you might have to choose
one, but any choice is fine). Since $d$ is a positive divisor of $n$, we have
$1\leq d\leq n$ (strictly speaking, this needs to be proved, but we take this
for granted here\footnote{Actually, the inequality $1\leq d$ is obvious (since
$d$ is a positive integer), whereas the inequality $d\leq n$ follows from
Proposition \ref{prop.div.abs} \textbf{(b)}.}). Therefore, $1<d<n$ (since $d$
is neither $1$ nor $n$). Hence, $d$ is one of the numbers $1,2,\ldots,n-1$
(actually $2,3,\ldots,n-1$, but we don't care).

Furthermore, $\dfrac{n}{d}$ is an integer (since $d$ is a divisor of $n$) and
positive (since $n$ and $d$ are positive). Multiplying the inequality $1<d$ by
$\dfrac{n}{d}$, we obtain $1\cdot\dfrac{n}{d}<d\cdot\dfrac{n}{d}$ (since we
can always multiply an inequality by a positive number\footnote{This is a
basic fact that we are taking for granted.}). In other words, $\dfrac{n}{d}%
<n$. Since $\dfrac{n}{d}$ is a positive integer, we thus conclude that
$\dfrac{n}{d}$ is one of the numbers $1,2,\ldots,n-1$.

Now, our induction hypothesis says that $P\left(  1\right)  $ AND $P\left(
2\right)  $ AND\ $\cdots$ AND $P\left(  n-1\right)  $ holds. In particular,
$P\left(  d\right)  $ holds (since $d$ is one of the numbers $1,2,\ldots
,n-1$). In other words, $d$ is a product of primes. That is, we can write $d$
as%
\[
d=p_{1}p_{2}\cdots p_{k}\ \ \ \ \ \ \ \ \ \ \text{for some primes }p_{1}%
,p_{2},\ldots,p_{k}.
\]
Consider these primes $p_{1},p_{2},\ldots,p_{k}$.

Again, our induction hypothesis says that $P\left(  1\right)  $ AND $P\left(
2\right)  $ AND\ $\cdots$ AND $P\left(  n-1\right)  $ holds. In particular,
$P\left(  \dfrac{n}{d}\right)  $ holds (since $\dfrac{n}{d}$ is one of the
numbers $1,2,\ldots,n-1$). In other words, $\dfrac{n}{d}$ is a product of
primes. That is, we can write $\dfrac{n}{d}$ as%
\[
\dfrac{n}{d}=q_{1}q_{2}\cdots q_{\ell}\ \ \ \ \ \ \ \ \ \ \text{for some
primes }q_{1},q_{2},\ldots,q_{\ell}.
\]
Consider these primes $q_{1},q_{2},\ldots,q_{\ell}$.

Now,%
\[
n=d\cdot\dfrac{n}{d}=p_{1}p_{2}\cdots p_{k}\cdot q_{1}q_{2}\cdots q_{\ell}%
\]
(since $d=p_{1}p_{2}\cdots p_{k}$ and $\dfrac{n}{d}=q_{1}q_{2}\cdots q_{\ell}%
$). This shows that $n$ is a product of primes (since $p_{1},p_{2}%
,\ldots,p_{k}$ as well as $q_{1},q_{2},\ldots,q_{\ell}$ are primes). In other
words, $P\left(  n\right)  $ holds. Thus, we have proved $P\left(  n\right)  $
in Case 2.

Now, we have proved $P\left(  n\right)  $ both in Case 1 and Case 2.
Therefore, $P\left(  n\right)  $ always holds. Thus, the induction step is
complete, and Theorem \ref{thm.ent.primefac.exists} is proven.
\end{proof}

The above proof is just reflecting the elementary recursive algorithm for
factoring an integer $n$ into a product of primes: We search for a positive
divisor $d$ of $n$ that is neither $1$ nor $n$. If such a $d$ does not exist,
then $n$ itself is a prime (or $1$, in which case the empty factorization
works). If it does, then we are reduced to the simpler problems of factoring
$d$ and $\dfrac{n}{d}$, and just have to multiply the resulting factorizations
at the end.

\begin{noncompile}
The name \textquotedblleft strong induction\textquotedblright\ refers to the
fact that it \textbf{feels} stronger than regular induction: After all, in the
induction step of a strong induction, we can use not just the single previous
claim $P\left(  n-1\right)  $, but all previously proved claims $P\left(
0\right)  ,\ P\left(  1\right)  ,\ \ldots,\ P\left(  n-1\right)  $.
Nevertheless, \textbf{logically} the principle of strong induction is
\textbf{equivalent} to the principle of regular induction. In other words,
anything you can prove by regular induction you can also prove by strong
induction, and vice versa. The first part of this is obvious, but the second
might appear surprising.

For example, let me explain how to prove Binet's formula by regular induction.
This is despite our first attempt at doing so being unsuccessful!
Nevertheless, this can be done. The trick is to define $P\left(  n\right)  $
to be \textbf{not} the statement \textquotedblleft$f_{n}=\dfrac{\varphi
^{n}-\psi^{n}}{\sqrt{5}}$\textquotedblright\ (for a specific $n$), but rather
the statement \textquotedblleft$f_{k}=\dfrac{\varphi^{k}-\psi^{k}}{\sqrt{5}}$
for all $k\in\left\{  0,1,\ldots,n\right\}  $\textquotedblright. In other
words, we \textbf{incorporate} all the statements $P\left(  0\right)
,P\left(  1\right)  ,\ldots,P\left(  n-1\right)  $ into our new statement
$P\left(  n\right)  $. Then, the induction hypothesis will include Binet's
formula not just for $n-1$ but also for $n-2$ (and for all the other numbers
before it). So we get all the advantages of strong induction from regular
induction by just pushing all the necessary memory into the claim that we are proving.

For another example, instead of proving the claim%
\[
\text{\textquotedblleft}n\text{ is a product of primes\textquotedblright}%
\]
by strong induction, we can prove the claim%
\[
\text{\textquotedblleft each integer between }1\text{ and }n\text{ is a
product of primes\textquotedblright}%
\]
by regular induction.

This all is an aside; it's perfectly fine to use strong induction.
\end{noncompile}

\subsubsection{\label{subsec.ind.strind.coin-35}Example: Paying with $3$-cent
and $5$-cent coins}

Here is another example of how strong induction can be used:

\begin{exercise}
\label{exe.ind.strind.coins-35}Assume that you have $3$-cent coins and
$5$-cent coins (each in infinite supply). What denominations can you pay with
these coins?
\end{exercise}

Let's make a table (\textquotedblleft yes\textquotedblright\ means that you
can pay it; \textquotedblleft no\textquotedblright\ means that you can't):%
\[%
\begin{tabular}
[c]{|c|c|}\hline
$0$ cents & yes\\\hline
$1$ cents & no\\\hline
$2$ cents & no\\\hline
$3$ cents & yes\\\hline
$4$ cents & no\\\hline
$5$ cents & yes\\\hline
$6$ cents & yes: $2\cdot3$\\\hline
$7$ cents & no\\\hline
$8$ cents & yes: $3+5$\\\hline
$9$ cents & yes: $3\cdot3$\\\hline
$10$ cents & yes: $2\cdot5$\\\hline
$11$ cents & yes: $2\cdot3+5$\\\hline
$12$ cents & yes: $4\cdot3$\\\hline
$13$ cents & yes: $3+2\cdot5$\\\hline
$\cdots$ & $\cdots$\\\hline
\end{tabular}
\ \ .
\]
Experimentally, we seem to observe that any denomination $\geq8$ cents can be
paid. Why?

We can notice that if a denomination $k$ (that is, $k$ cents) can be paid,
then so can $k+3$ (just add a $3$-cent coin). Thus, because we can pay $8$
cents, we can also pay $11,14,17,\ldots$ cents. Because we can pay $9$ cents,
we can also pay $12,15,18,\ldots$ cents. Because we can pay $10$ cents, we can
also pay $13,16,19,\ldots$ cents. Together, these three sequences account for
all the integers $\geq8$. Thus, any denomination of $\geq8$ cents can be paid.

Let us formalize this argument as an induction proof.

We define $\mathbb{N}$ to be the set of all nonnegative integers:%
\[
\mathbb{N}=\left\{  0,1,2,\ldots\right\}  .
\]

\begin{proposition}
\label{prop.int.coins-35}For any integer $n\geq8$, we can pay $n$ cents with
$3$-cent and $5$-cent coins. In other words, any integer $n\geq8$ can be
written as $n=3a+5b$ with $a,b\in\mathbb{N}$.
\end{proposition}

\begin{proof}
We proceed by strong induction on $n$:

\textit{Base case:} For $n=8$, the claim is true, since $8=3\cdot1+5\cdot1$.

\textit{Induction step:} Fix an integer $n>8$. Assume that the proposition is
already proved for all the integers $8,9,\ldots,n-1$. We must prove that it
also holds for $n$. In other words, we must prove that we can pay $n$ cents
with $3$-cent and $5$-cent coins.

We are in one of the following three cases (since $n>8$):

\textit{Case 1:} We have $n=9$.

\textit{Case 2:} We have $n=10$.

\textit{Case 3:} We have $n\geq11$.

In Case 1, we are done, since $n=9=3\cdot3+5\cdot0$ (that is, $n$ cents can be
paid with three $3$-cent coins).

In Case 2, we are done, since $n=10=3\cdot0+5\cdot2$ (that is, $n$ cents can
be paid with two $5$-cent coins).

Now, consider Case 3. In this case, we have $n\geq11$. Hence, $n-3\geq8$. This
shows that $n-3$ is one of the numbers $8,9,\ldots,n-1$.

Thus, we can apply the induction hypothesis to $n-3$. We conclude that $n-3$
cents can be paid with $3$-cent and $5$-cent coins, i.e., we can write $n-3$
as $n-3=3c+5d$ with $c,d\in\mathbb{N}$. Using these $c,d\in\mathbb{N}$, we
therefore have%
\begin{align*}
n  &  =3+3c+5d\ \ \ \ \ \ \ \ \ \ \left(  \text{since }n-3=3c+5d\right) \\
&  =3\left(  c+1\right)  +5d,
\end{align*}
which shows that $n$ cents can also be paid with $3$-cent and $5$-cent coins.
This shows that the proposition is true for $n$, and thus the induction step
is complete.

The proposition is thus proved.
\end{proof}

Note that the above proof had one \textquotedblleft de-jure base
case\textquotedblright\ (the case $n=8$) and two \textquotedblleft de-facto
base cases\textquotedblright\ (the cases $n=9$ and $n=10$, which were formally
part of the induction step but had to be treated separately because $n-3$
would be smaller than $8$ in these cases). We could have just as well used the
baseless form of strong induction, in which case we would have to treat all
three of these cases as \textquotedblleft de-facto base
cases\textquotedblright. This would be a bit more uniform, although this is
entirely a matter of taste.

More examples of strong induction can be found in \cite[Appendix A]{AlNoWo19}
(where it is called \textquotedblleft top-down induction\textquotedblright,
and generalized beyond integers).

\subsection{\label{sec.ind.moreex}More exercises}

Let us finish this chapter with some further exercises on induction.

\subsubsection{A fake proof}

\begin{exercise}
\label{exe.21finhw3.4}Find the error(s) in the following fake proof:

We claim that $3^{n}=1$ for each $n\in\mathbb{N}$.

\textit{\textquotedblleft Proof:\textquotedblright} We proceed by strong
induction on $n$. So we let $n\in\mathbb{N}$ be arbitrary, and we assume (as
the induction hypothesis) that $3^{k}=1$ for each $k<n$. We must now prove
that $3^{n}=1$.

By our induction hypothesis, we have $3^{n-1}=1$ (since $n-1<n$) and
$3^{n-2}=1$ (since $n-2<n$). Now, $3^{n}=\dfrac{\left(  3^{n-1}\right)  ^{2}%
}{3^{n-2}}$ (since the laws of exponents yield $\dfrac{\left(  3^{n-1}\right)
^{2}}{3^{n-2}}=3^{2\cdot\left(  n-1\right)  -\left(  n-2\right)  }=3^{n}$). In
view of $3^{n-1}=1$ and $3^{n-2}=1$, this rewrites as $3^{n}=\dfrac{1^{2}}%
{1}=1$. This completes the induction step, and thus the claim is proved.
\end{exercise}

\subsubsection{\label{subsec.ind.moreex.f-n}Negative Fibonacci numbers}

Recall again the Fibonacci sequence $\left(  f_{0},f_{1},f_{2},\ldots\right)
$ from Definition \ref{def.ind.fib}. Let us now extend this sequence
\textquotedblleft to the left\textquotedblright\ by defining $f_{n}$ not only
for nonnegative integers $n$, but also for negative integers $n$. To do so, we
simply rewrite the equation $f_{n}=f_{n-1}+f_{n-2}$ (which we used to
recursively define the Fibonacci sequence) as $f_{n-2}=f_{n}-f_{n-1}$. This
allows us to compute $f_{n-2}$ from $f_{n}$ and $f_{n-1}$. Thus, we can
compute $f_{-1}$ from $f_{1}$ and $f_{0}$, then compute $f_{-2}$ from $f_{0}$
and $f_{-1}$, and so on:%
\begin{align*}
f_{-1}  &  =f_{1}-f_{0}=1-0=1;\\
f_{-2}  &  =f_{0}-f_{-1}=0-1=-1;\\
f_{-3}  &  =f_{-1}-f_{-2}=1-\left(  -1\right)  =2;\\
f_{-4}  &  =f_{-2}-f_{-3}=\left(  -1\right)  -2=-3;\\
&  \ldots.
\end{align*}
Thus, we gradually extend the Fibonacci sequence to the left, obtaining a
\textquotedblleft two-sided sequence\textquotedblright\ $\left(  \ldots
,f_{-2},f_{-1},f_{0},f_{1},f_{2},\ldots\right)  $ that is \textquotedblleft
infinite in both directions\textquotedblright. By virtue of its construction,
it satisfies $f_{n}=f_{n-1}+f_{n-2}$ not only for all $n\geq2$, but also for
all integers $n$. However, a quick look at the first (say) $7$
\textquotedblleft extended\textquotedblright\ Fibonacci numbers to the left of
$f_{0}$ reveals that they are not as new as they might seem: They are just
copies of the positive Fibonacci numbers with signs. More precisely, it looks
like we have%
\begin{equation}
f_{-n}=\left(  -1\right)  ^{n-1}f_{n}\ \ \ \ \ \ \ \ \ \ \text{for each }%
n\geq0. \label{eq.fib.neg.f-n=}%
\end{equation}

\begin{exercise}
\label{exe.ind.fib.neg.f-n=}\textbf{(a)} Try to prove (\ref{eq.fib.neg.f-n=})
directly by induction on $n$. (So the induction step involves assuming that
$f_{-n}=\left(  -1\right)  ^{n-1}f_{n}$ and proving that $f_{-\left(
n+1\right)  }=\left(  -1\right)  ^{n}f_{n+1}$. Don't use strong induction
yet!) Does this work? \medskip

\textbf{(b)} Now, instead, try to prove the \textbf{stronger} claim that
\textquotedblleft$f_{-n}=\left(  -1\right)  ^{n-1}f_{n}$ and $f_{-n+1}=\left(
-1\right)  ^{n-2}f_{n-1}$ for each $n\geq0$\textquotedblright\ by induction on
$n$. Does this work? \medskip

\textbf{(c)} Now, prove (\ref{eq.fib.neg.f-n=}) by \textbf{strong} induction
on $n$.
\end{exercise}

\subsubsection{More on the Hanoi tower}

\begin{exercise}
\label{exe.23hw1.6}Let $n\geq0$ be an integer, and let $k\in\left\{
1,2,\ldots,n\right\}  $. In the proof of Proposition \ref{prop.hanoi.1}, we
presented a certain strategy for solving the Tower of Hanoi puzzle with $n$ disks.

Prove that the $k$-th largest disk is moved exactly $2^{k-1}$ many times
during this strategy.
\end{exercise}

\subsubsection{More on recursively defined sequences}

\begin{exercise}
\label{exe.23hw2.4}Let $\left(  a_{0},a_{1},a_{2},\ldots\right)  $ be a
sequence of integers defined recursively by%
\begin{align*}
a_{0}  &  =2,\ \ \ \ \ \ \ \ \ \ a_{1}=3,\\
a_{n}  &  =3a_{n-1}-2a_{n-2}\ \ \ \ \ \ \ \ \ \ \text{for all }n\geq2.
\end{align*}
Prove that $a_{n}=2^{n}+1$ for each integer $n\geq0$.
\end{exercise}

\begin{exercise}
\label{exe.23mt1.2}Let $\left(  a_{0},a_{1},a_{2},\ldots\right)  $ be a
sequence of integers defined recursively by%
\begin{align*}
a_{0}  &  =2,\ \ \ \ \ \ \ \ \ \ a_{1}=1,\\
a_{n}  &  =a_{n-1}+6a_{n-2}\ \ \ \ \ \ \ \ \ \ \text{for all }n\geq2.
\end{align*}
Prove that $a_{n}=3^{n}+\left(  -2\right)  ^{n}$ for each $n\in\mathbb{N}$.
\end{exercise}

\begin{exercise}
\label{exe.23hw2.5}Recall the Fibonacci sequence (Definition \ref{def.ind.fib}%
) again. \medskip

\textbf{(a)} Let $k$ be a nonnegative integer. Show that
\[
f_{n}^{2}-f_{n+k}f_{n-k}=\left(  -1\right)  ^{n-k}f_{k}^{2}%
\ \ \ \ \ \ \ \ \ \ \text{for every integer }n\geq k.
\]

\textbf{(b)} Which of the previously posed exercises does this generalize?
\end{exercise}

\begin{exercise}
Recall the Fibonacci sequence (Definition \ref{def.ind.fib}) again. Let
$n\geq0$. Prove that%
\begin{align*}
&  f_{3n}\text{ is even;}\\
&  f_{3n+1}\text{ is odd;}\\
&  f_{3n+2}\text{ is odd.}%
\end{align*}

(In this exercise, you can freely use basic properties of even and odd numbers
-- such as Proposition \ref{prop.ent.odd=2k+1}.)
\end{exercise}

\begin{exercise}
\label{exe.21finhw3.5}Define a sequence $\left(  t_{0},t_{1},t_{2}%
,\ldots\right)  $ of positive rational numbers recursively by setting%
\begin{align*}
t_{0}  &  =1,\ \ \ \ \ \ \ \ \ \ t_{1}=1,\ \ \ \ \ \ \ \ \ \ t_{2}%
=1,\ \ \ \ \ \ \ \ \ \ \text{and}\\
t_{n}  &  =\dfrac{1+t_{n-1}t_{n-2}}{t_{n-3}}\ \ \ \ \ \ \ \ \ \ \text{for each
}n\geq3.
\end{align*}
(So its next entries after $t_{2}$ are%
\begin{align*}
t_{3}  &  =\dfrac{1+t_{2}t_{1}}{t_{0}}=\dfrac{1+1\cdot1}{1}=2;\\
t_{4}  &  =\dfrac{1+t_{3}t_{2}}{t_{1}}=\dfrac{1+2\cdot1}{1}=3;\\
t_{5}  &  =\dfrac{1+t_{4}t_{3}}{t_{2}}=\dfrac{1+3\cdot2}{1}=7;\\
t_{6}  &  =\dfrac{1+t_{5}t_{4}}{t_{3}}=\dfrac{1+7\cdot3}{2}=11,
\end{align*}
and so on.) \medskip

\textbf{(a)} Prove that $t_{n+2}=4t_{n}-t_{n-2}$ for each $n\geq2$. \medskip

\textbf{(b)} Prove that $t_{n}$ is a positive integer for each integer
$n\geq0$. \medskip

[\textbf{Hint:} Use regular induction for part \textbf{(a)} and strong
induction for part \textbf{(b)}. Note that the \textquotedblleft
positive\textquotedblright\ part is clear from the definition, so you only
need to prove the \textquotedblleft integer\textquotedblright\ part in
\textbf{(b)}.]
\end{exercise}

\subsubsection{More coin problems}

\begin{exercise}
\textbf{(a)} Prove the following: For any integer $n\geq12$, we can pay $n$
cents with $3$-cent and $7$-cent coins. In other words, any integer $n\geq12$
can be written as $n=3a+7b$ with $a,b\in\mathbb{N}$. (Here, again,
$\mathbb{N}=\left\{  0,1,2,\ldots\right\}  $.) \medskip

\textbf{(b)} Find the largest integer $k$ such that $k$ cents cannot be paid
with $2$-cent and $13$-cent coins. Prove that for every integer $n>k$, we can
pay $n$ cents with these kinds of coins. \medskip

\textbf{(c)} Is there a largest integer $k$ such that $k$ cents cannot be paid
with $2$-cent and $6$-cent coins?
\end{exercise}

\subsubsection{A bit of matrix algebra}

The next two exercises are about matrix multiplication. For an introduction to
matrix multiplication, see any textbook on linear algebra (e.g.,
\cite[\S 10.1]{BoyVan18} covers it in detail). However, all we need for these
exercises will be $2\times2$-matrices, so let us recall how matrix
multiplication works for them:

\begin{itemize}
\item The product $AX$ of two $2\times2$-matrices $A=\left(
\begin{array}
[c]{cc}%
a & b\\
c & d
\end{array}
\right)  $ and $X=\left(
\begin{array}
[c]{cc}%
x & y\\
z & w
\end{array}
\right)  $ is defined to be $\left(
\begin{array}
[c]{cc}%
ax+bz & ay+bw\\
cx+dz & cy+dw
\end{array}
\right)  $.

\item The $n$-th power $A^{n}$ of a $2\times2$-matrix $A$ is defined to be the
product $\underbrace{AA\cdots A}_{n\text{ factors}}$.
\end{itemize}

\begin{exercise}
\label{exe.21finhw3.1}\textbf{(a)} Prove that $\left(
\begin{array}
[c]{cc}%
1 & 1\\
0 & 1
\end{array}
\right)  ^{n}=\left(
\begin{array}
[c]{cc}%
1 & n\\
0 & 1
\end{array}
\right)  $ for each positive integer $n$. \medskip

\textbf{(b)} Find a formula for $\left(
\begin{array}
[c]{cc}%
a & b\\
0 & c
\end{array}
\right)  ^{n}$, where $a,b,c$ are real numbers and $n$ is a positive integer.
\end{exercise}

\begin{exercise}
Recall the Fibonacci sequence $\left(  f_{0},f_{1},f_{2},\ldots\right)
=\left(  0,1,1,2,3,5,\ldots\right)  $. Prove that%
\[
\left(
\begin{array}
[c]{cc}%
1 & 1\\
1 & 0
\end{array}
\right)  ^{n}=\left(
\begin{array}
[c]{cc}%
f_{n+1} & f_{n}\\
f_{n} & f_{n-1}%
\end{array}
\right)  \ \ \ \ \ \ \ \ \ \ \text{for each positive integer }n.
\]

\end{exercise}

\subsubsection{More induction proofs}

\begin{exercise}
\label{exe.21finhw3.6}Let $m\in\mathbb{N}$. Prove that there exists a way to
arrange the first $m$ positive integers ($1,2,\ldots,m$) in a row in such a
way that the average of two numbers never stands between these two numbers.

(For example, for $m=8$, one such arrangement is $1,5,3,7,2,6,4,8$. The
arrangement $1,3,2,7,8,5,6,4$ is invalid because the average of $1$ and $5$ is
$3$, which stands between $1$ and $5$.) \medskip

[\textbf{Hint:} First show that there is such an arrangement when $m$ is a
power of $2$ (that is, when $m=2^{n}$ for some $n\in\mathbb{N}$). Then, choose
a sufficiently large power of $2$ and remove all entries larger than $m$.]
\end{exercise}

More advanced and creative uses of induction can be found in \cite[Chapter
2]{20f}, \cite[Lecture 1]{23f}, \cite{AndCri17}, \cite{Gunder10} and
\cite{Weintr17}.

\newpage

\section{\label{chp.sums}Sums and products}

\subsection{\label{sec.sums.sums}Finite sums}

Previously, we have encountered sums such as%
\[
x^{n-1}+x^{n-2}y+x^{n-3}y^{2}+\cdots+x^{2}y^{n-3}+xy^{n-2}+y^{n-1}%
\]
(in Section \ref{sec.ind.geosums}). Such sums can be tricky to decipher: You
need to guess the pattern of the addends to understand what the
\textquotedblleft$\cdots$\textquotedblright\ means. There is a notation that
makes such sums both shorter and easier to understand. This is the
\textbf{finite sum notation} (also known as the \textbf{sigma notation}). In
its simplest form, it is defined as follows:

\begin{definition}
\label{def.finsum.k=u-v}Let $u$ and $v$ be two integers. Let $a_{u}%
,a_{u+1},\ldots,a_{v}$ be some numbers. Then,%
\[
\sum_{k=u}^{v}a_{k}%
\]
is defined to be the sum%
\[
a_{u}+a_{u+1}+\cdots+a_{v}%
\]
(in more detail: $a_{u}+a_{u+1}+a_{u+2}+a_{u+3}+\cdots+a_{v-1}+a_{v}$). It is
called the \textbf{sum of the numbers }$a_{k}$ \textbf{where }$k$
\textbf{ranges from }$u$ \textbf{to }$v$. When $v<u$, this sum is called
\textbf{empty} and defined to be $0$. When $v=u$, this sum contains only one
addend (namely, $a_{u}$) and equals this addend.
\end{definition}

For example:%
\begin{align*}
\sum_{k=5}^{10}k  &  =5+6+7+8+9+10=45;\\
\sum_{k=5}^{10}\dfrac{1}{k}  &  =\dfrac{1}{5}+\dfrac{1}{6}+\dfrac{1}{7}%
+\dfrac{1}{8}+\dfrac{1}{9}+\dfrac{1}{10}=\dfrac{2131}{2520};\\
\sum_{k=5}^{10}k^{k}  &  =5^{5}+6^{6}+7^{7}+8^{8}+9^{9}+10^{10};\\
\sum_{k=5}^{5}k  &  =5;\\
\sum_{k=5}^{4}k  &  =0\ \ \ \ \ \ \ \ \ \ \left(  \text{an empty sum}\right)
;\\
\sum_{k=5}^{3}k  &  =0\ \ \ \ \ \ \ \ \ \ \left(  \text{an empty sum}\right)
;\\
\sum_{k=5}^{8}3  &  =3+3+3+3=12\ \ \ \ \ \ \ \ \ \ \left(  \text{a sum of four
equal terms}\right)  ;\\
\sum_{k=0}^{n-1}q^{k}  &  =q^{0}+q^{1}+\cdots+q^{n-1}%
\ \ \ \ \ \ \ \ \ \ \text{for any }n\in\mathbb{N}\text{ and any number }q;\\
\sum_{k=0}^{n-1}x^{k}y^{n-1-k}  &  =x^{0}y^{n-1}+x^{1}y^{n-2}+x^{2}%
y^{n-3}+\cdots+x^{n-3}y^{2}+x^{n-2}y^{1}+x^{n-1}y^{0}\\
&  =y^{n-1}+xy^{n-2}+x^{2}y^{n-3}+\cdots+x^{n-3}y^{2}+x^{n-2}y+x^{n-1}\\
&  =x^{n-1}+x^{n-2}y+x^{n-3}y^{2}+\cdots+x^{2}y^{n-3}+xy^{n-2}+y^{n-1}\\
&  \ \ \ \ \ \ \ \ \ \ \ \ \ \ \ \ \ \ \ \ \text{for any }n\in\mathbb{N}\text{
and any numbers }x\text{ and }y.
\end{align*}
Thus, Theorem \ref{thm.geo-sum.xy} is saying that%
\[
\left(  x-y\right)  \left(  \sum_{k=0}^{n-1}x^{k}y^{n-1-k}\right)
=x^{n}-y^{n}%
\]
for any numbers $x$ and $y$ and any $n\in\mathbb{N}$.

The variable $k$ is not set in stone; you can replace it by any other variable
(unless this other variable already stands for something else). For example,%
\[
\sum_{k=u}^{v}a_{k}=\sum_{i=u}^{v}a_{i}=\sum_{\circledS=u}^{v}a_{\circledS
}=\sum_{\spadesuit=u}^{v}a_{\spadesuit}.
\]
Just don't make it $\sum\limits_{u=u}^{v}a_{u}$.

Here are a couple more examples: For any $n\in\mathbb{N}$, we have%
\begin{align*}
\sum_{k=1}^{n}k  &  =1+2+\cdots+n=\dfrac{n\left(  n+1\right)  }{2}%
\ \ \ \ \ \ \ \ \ \ \left(  \text{by Theorem \ref{thm.ind.1+2+...+n}}\right)
;\\
\sum_{k=1}^{n}k^{2}  &  =1^{2}+2^{2}+\cdots+n^{2}\\
&  =\dfrac{n\left(  n+1\right)  \left(  2n+1\right)  }{6}%
\ \ \ \ \ \ \ \ \ \ \left(  \text{by Theorem \ref{thm.ind.12+22+...+n2}%
}\right)  ;\\
\sum_{k=1}^{n}1  &  =\underbrace{1+1+\cdots+1}_{n\text{ times}}=n\cdot1=n;\\
\sum_{k=1}^{n}\left(  2k-1\right)   &  =\left(  2\cdot1-1\right)  +\left(
2\cdot2-1\right)  +\left(  2\cdot3-1\right)  +\cdots+\left(  2n-1\right) \\
&  =1+3+5+\cdots+\left(  2n-1\right) \\
&  =\left(  \text{the sum of the first }n\text{ odd positive integers}\right)
.
\end{align*}

We have not computed this last sum, so let us do this. I will use the
following \textquotedblleft laws of summation\textquotedblright:

\begin{itemize}
\item We have
\begin{equation}
\sum_{k=u}^{v}\left(  a_{k}-b_{k}\right)  =\sum_{k=u}^{v}a_{k}-\sum_{k=u}%
^{v}b_{k} \label{eq.sums.sum(a-b)}%
\end{equation}
for any integers $u,v$ and any numbers $a_{k},b_{k}$. Indeed, if you rewrite
this without finite sum notation, it takes the form%
\begin{align*}
&  \left(  a_{u}-b_{u}\right)  +\left(  a_{u+1}-b_{u+1}\right)  +\cdots
+\left(  a_{v}-b_{v}\right) \\
&  =\left(  a_{u}+a_{u+1}+\cdots+a_{v}\right)  -\left(  b_{u}+b_{u+1}%
+\cdots+b_{v}\right)  ,
\end{align*}
which is rather clear. (A formal proof can be given by induction on $v$.)

\item We have%
\begin{equation}
\sum_{k=u}^{v}\lambda a_{k}=\lambda\sum_{k=u}^{v}a_{k} \label{eq.sums.sum(la)}%
\end{equation}
for any integers $u,v$ and any numbers $\lambda,a_{k}$. Indeed, rewritten
without the use of finite sum notation, this is just saying that%
\[
\lambda a_{u}+\lambda a_{u+1}+\cdots+\lambda a_{v}=\lambda\left(
a_{u}+a_{u+1}+\cdots+a_{v}\right)  ,
\]
which is again clear (and can be proved by induction on $v$).
\end{itemize}

Rules like these are a dime a dozen, and you should be able to come up with
them on the spot when you need them. (See \cite[\S 1.4.2]{detnotes} for these
and several others.)

Let us now compute our sum:%
\begin{align*}
\sum_{k=1}^{n}\left(  2k-1\right)   &  =\sum_{k=1}^{n}2k-\sum_{k=1}%
^{n}1\ \ \ \ \ \ \ \ \ \ \left(  \text{by (\ref{eq.sums.sum(a-b)})}\right) \\
&  =2\underbrace{\sum_{k=1}^{n}k}_{=\dfrac{n\left(  n+1\right)  }{2}%
}-\underbrace{\sum_{k=1}^{n}1}_{=n}\ \ \ \ \ \ \ \ \ \ \left(  \text{by
(\ref{eq.sums.sum(la)})}\right) \\
&  =2\cdot\dfrac{n\left(  n+1\right)  }{2}-n=n\left(  n+1\right)  -n=n^{2}.
\end{align*}

As another illustration of the use of our notation, we can rewrite Gauss's
proof of the equality
\begin{equation}
\sum_{k=1}^{n}k=1+2+\cdots+n=\dfrac{n\left(  n+1\right)  }{2}
\label{eq.sums.little-gauss}%
\end{equation}
(Theorem \ref{thm.ind.1+2+...+n}) using finite sum notation. We will need
three new rules this time:

\begin{itemize}
\item We have
\begin{equation}
\sum_{k=u}^{v}a_{k}+\sum_{k=u}^{v}b_{k}=\sum_{k=u}^{v}\left(  a_{k}%
+b_{k}\right)  \label{eq.sums.suma+sumb}%
\end{equation}
for any integers $u,v$ and any numbers $a_{k},b_{k}$. Indeed, if you rewrite
this without finite sum notation, it takes the form%
\begin{align*}
&  \left(  a_{u}+a_{u+1}+\cdots+a_{v}\right)  +\left(  b_{u}+b_{u+1}%
+\cdots+b_{v}\right) \\
&  =\left(  a_{u}+b_{u}\right)  +\left(  a_{u+1}+b_{u+1}\right)
+\cdots+\left(  a_{v}+b_{v}\right)  .
\end{align*}

\item We have%
\begin{equation}
\sum_{k=u}^{v}a_{k}=\sum_{k=u}^{v}a_{u+v-k} \label{eq.sums.revert}%
\end{equation}
for any integers $u,v$ and any numbers $a_{k}$. This is called
\textquotedblleft substituting $u+v-k$ for $k$ in the sum\textquotedblright%
\ or just \textquotedblleft turning the sum upside-down\textquotedblright, as
it amounts to reversing the order of the addends; restated without finite sum
notation, this is just saying that%
\[
a_{u}+a_{u+1}+\cdots+a_{v}=a_{v}+a_{v-1}+\cdots+a_{u},
\]
which is saying that a sum of a bunch of numbers does not change if we add its
addends together in reverse order.

\item For any integers $u\leq v$ and any number $\lambda$, we have
\begin{equation}
\sum_{k=u}^{v}\lambda=\left(  v-u+1\right)  \lambda. \label{eq.sums.lam}%
\end{equation}
(This is just saying that a sum of $v-u+1$ many equal addends $\lambda$ is
$\left(  v-u+1\right)  \lambda$. Note that the sum on the left hand side has
$v-u+1$ addends, because there are $v-u+1$ numbers in the set $\left\{
u,u+1,\ldots,v\right\}  $.)

Actually, the \textquotedblleft right\textquotedblright\ condition for
(\ref{eq.sums.lam}) to hold is not $u\leq v$ but $u\leq v+1$. If $u=v+1$, then
the left hand side of (\ref{eq.sums.lam}) is an empty sum, thus $0$ by
definition, while the right hand side is $\underbrace{\left(  v-u+1\right)
}_{\substack{=0\\\text{(since }u=v+1\text{)}}}\lambda=0$ as well. Only once
$u$ surpasses $v+1$ does (\ref{eq.sums.lam}) start failing.
\end{itemize}

Now, Gauss's proof of (\ref{eq.sums.little-gauss}) takes the following shape:%
\begin{align*}
2\sum_{k=1}^{n}k  &  =\sum_{k=1}^{n}k+\sum_{k=1}^{n}k\\
&  =\sum_{k=1}^{n}k+\sum_{k=1}^{n}\left(  n+1-k\right) \\
&  \ \ \ \ \ \ \ \ \ \ \ \ \ \ \ \ \ \ \ \ \left(
\begin{array}
[c]{c}%
\text{here, we substituted }n+1-k\text{ for }k\text{ in the second}\\
\text{sum (i.e., rewrote it using (\ref{eq.sums.revert}))}%
\end{array}
\right) \\
&  =\sum_{k=1}^{n}\underbrace{\left(  k+\left(  n+1-k\right)  \right)
}_{=n+1}\ \ \ \ \ \ \ \ \ \ \left(  \text{by (\ref{eq.sums.suma+sumb})}\right)
\\
&  =\sum_{k=1}^{n}\left(  n+1\right)  =n\cdot\left(  n+1\right)
\ \ \ \ \ \ \ \ \ \ \left(  \text{by (\ref{eq.sums.lam})}\right)  .
\end{align*}
Dividing both sides by $2$, we recover (\ref{eq.sums.little-gauss}) again.

\bigskip

We have found closed-form expressions (i.e., expressions without $\sum$ signs
or \textquotedblleft$\cdots$\textquotedblright s) for several sums. Not every
sum has a closed-form expression\footnote{It is actually not clear what
\textquotedblleft closed-form expression\textquotedblright\ even means: the
more you think about it, the more nebulous the concept becomes. For instance,
we view $\dfrac{n\left(  n+1\right)  }{2}$ as a closed form because we have
learnt division in elementary school, even though the expression
$1+2+\cdots+n$ is more elementary (relying only on addition). Are roots
allowed in a closed form? What if these roots are imaginary? What about
logarithms, factorials, trigonometric functions? The internet is awash in
tongue-in-cheek \textquotedblleft closed forms\textquotedblright\ for the
$n$-th prime number (e.g.,
\href{https://mathworld.wolfram.com/PrimeFormulas.html}{on MathWorld} or
\href{https://hsm.stackexchange.com/questions/13353/who-discovered-this-closed-form-formula-for-the-n-th-prime-number}{on
StackExchange}); these formulas usually take longer to evaluate than simply
testing the integers $1,2,3,\ldots$ successively until one has found $n$
primes. So \textquotedblleft closed form\textquotedblright\ is more of a vibe
than a mathematical notion; in particular, it is impossible to \textbf{prove}
that a closed form for a given sum (or function, or sequence) does not exist,
unless one specifies precisely what such a closed form is allowed to consist
of.
\par
Nevertheless, I shall talk about closed forms several times in this course;
just keep in mind that such claims (about something being a closed form, or
something not having a closed form) are not meant to be precise mathematical
statements.}. For instance, there is no closed form for%
\[
\sum_{k=1}^{n}\dfrac{1}{k}\ \ \ \ \ \ \ \ \ \ \text{or for}%
\ \ \ \ \ \ \ \ \ \ \sum_{k=1}^{n}k^{k}.
\]

\bigskip

Some more terminology:

The notation $\sum\limits_{k=u}^{v}a_{k}$ is called \textbf{sigma notation} or
\textbf{finite sum notation}. The symbol $\sum$ itself is called the
\textbf{summation sign}. The numbers $u$ and $v$ are called the \textbf{lower
limit} and the \textbf{upper limit} of the summation\footnote{This use of the
word \textquotedblleft limit\textquotedblright\ is totally unrelated to the
way this word is used in analysis/calculus.}. The variable $k$ is called the
\textbf{summation index} or the \textbf{running index}, and is said to
\textbf{range} (or \textbf{run}) from $u$ to $v$. The numbers $a_{k}$ are
called the \textbf{addends} of the finite sum.

\bigskip

There are many similarities between finite sums $\sum_{k=u}^{v}a_{k}$ and
integrals $\int\limits_{u}^{v}f\left(  x\right)  dx$, but the analogy should
not be taken too far (e.g., an integral $\int\limits_{u}^{u}f\left(  x\right)
dx$ whose upper and lower limit are equal will always be $0$, but an
\textquotedblleft analogous\textquotedblright\ finite sum $\sum_{k=u}^{u}%
a_{k}$ will be $a_{u}$).

\bigskip

We note two more rules for finite sums:

\begin{itemize}
\item The \textquotedblleft splitting-off rule\textquotedblright: For any
integers $u\leq v$ and any numbers $a_{u},a_{u+1},\ldots,a_{v}$, we have%
\[
\sum_{k=u}^{v}a_{k}=\sum_{k=u}^{v-1}a_{k}+a_{v}=a_{u}+\sum_{k=u+1}^{v}a_{k}.
\]
This is just saying that%
\begin{align*}
a_{u}+a_{u+1}+\cdots+a_{v}  &  =\left(  a_{u}+a_{u+1}+\cdots+a_{v-1}\right)
+a_{v}\\
&  =a_{u}+\left(  a_{u+1}+a_{u+2}+\cdots+a_{v}\right)  .
\end{align*}
This rule allows us to split the first or the last addend out of a finite sum.
This is important for proofs by induction.

\item More generally, any finite sum $\sum_{k=u}^{v}a_{k}$ can be split at any
point: We have%
\[
\sum_{k=u}^{v}a_{k}=\sum_{k=u}^{w}a_{k}+\sum_{k=w+1}^{v}a_{k}%
\]
for any integers $u\leq w\leq v$ and any numbers $a_{k}$. This is just saying
that%
\[
a_{u}+a_{u+1}+\cdots+a_{v}=\left(  a_{u}+a_{u+1}+\cdots+a_{w}\right)  +\left(
a_{w+1}+a_{w+2}+\cdots+a_{v}\right)  .
\]
(Strictly speaking, this is true not just for $u\leq w\leq v$ but more
generally for $u-1\leq w\leq v$. If you find this confusing, recall that an
empty sum equals $0$ by definition.)
\end{itemize}

\bigskip

Finite sum notation, in the form defined above, is helpful when the summation
index is running over an integer interval (i.e., a set of consecutive
integers). For more general situations, there is a more general version of
finite sum notation, e.g.:%
\[
\sum_{k\in\left\{  1,2,\ldots,n\right\}  \text{ is even}}k=2+4+6+\cdots+m,
\]
where $m$ is the largest even element of $\left\{  1,2,\ldots,n\right\}  $ (or
$0$ if no such even element exists). We won't use it much, but it is fairly
self-explanatory; essentially, the writing under the summation sign explains
what $k$'s the sum is ranging over. See \cite[\S 1.4.1]{detnotes} for a more
precise explanation.

\begin{exercise}
Let $n\in\mathbb{N}$. Prove that%
\[
\sum_{k=0}^{n}k\left(  n-k\right)  =\dfrac{\left(  n-1\right)  n\left(
n+1\right)  }{6}.
\]

\end{exercise}

\begin{exercise}
Let $n\in\mathbb{N}$. \medskip

\textbf{(a)} Prove that%
\[
\sum_{k=1}^{n}\dfrac{1}{k\left(  k+1\right)  }=\dfrac{n}{n+1}.
\]

\textbf{(b)} More generally: Let $b$ and $d$ be two numbers, and let
\[
a_{i}:=b+id\ \ \ \ \ \ \ \ \ \ \text{for each }i\in\left\{  1,2,\ldots
,n+1\right\}  .
\]
(Thus, $\left(  a_{1},a_{2},\ldots,a_{n+1}\right)  $ is what is called an
\textbf{arithmetic progression} -- i.e., a sequence of numbers that increase
from each to the next by the same amount $d$.) Assume that all the $n+1$
numbers $a_{1},a_{2},\ldots,a_{n+1}$ are nonzero. Prove that%
\[
\sum_{k=1}^{n}\dfrac{1}{a_{k}a_{k+1}}=\dfrac{n}{a_{1}a_{n+1}}.
\]

\end{exercise}

\begin{exercise}
\label{exe.23hw2.3}The \textbf{floor} $\left\lfloor x\right\rfloor $ of a real
number $x$ means the largest integer that is smaller than or equal to $x$. For
instance, $\left\lfloor 6.2\right\rfloor =6$ and $\left\lfloor
7.7\right\rfloor =7$ and $\left\lfloor 8\right\rfloor =8$. (In other words,
$\left\lfloor x\right\rfloor $ is what you get if you round $x$ down. Beware:
$\left\lfloor -1.3\right\rfloor $ is $-2$, not $-1$.)

Let $n\in\mathbb{N}$. Prove that%
\[
\sum_{k=1}^{n}\left\lfloor \dfrac{k}{2}\right\rfloor =\left\lfloor \dfrac
{n}{2}\right\rfloor \cdot\left\lfloor \dfrac{n+1}{2}\right\rfloor .
\]

(In this exercise, you can freely use basic properties of even and odd numbers
-- such as Proposition \ref{prop.ent.odd=2k+1}.)
\end{exercise}

\begin{exercise}
\label{exe.23mt1.1}Let $n\in\mathbb{N}$, and let $q$ be any number distinct
from $1$. Prove that%
\[
\sum_{k=1}^{n}kq^{k}=q\cdot\dfrac{nq^{n+1}-\left(  n+1\right)  q^{n}%
+1}{\left(  q-1\right)  ^{2}}.
\]

\end{exercise}

\begin{exercise}
\label{exe.23mt1.3}Let $n\in\mathbb{N}$. Prove that%
\[
\underbrace{1+2+\cdots+n}_{=\sum_{k=1}^{n}k}=\underbrace{n^{2}-\left(
n-1\right)  ^{2}+\left(  n-2\right)  ^{2}-\left(  n-3\right)  ^{2}\pm
\cdots+\left(  -1\right)  ^{n-1}1^{2}}_{=\sum_{k=1}^{n}\left(  -1\right)
^{n-k}k^{2}}.
\]

\end{exercise}

\subsection{\label{sec.sums.prods}Finite products}

Finite products are analogous to finite sums, just using multiplication
instead of addition:

\begin{definition}
\label{def.finprod.k=u-v}Let $u$ and $v$ be two integers. Let $a_{u}%
,a_{u+1},\ldots,a_{v}$ be some numbers. Then,%
\[
\prod_{k=u}^{v}a_{k}%
\]
is defined to be the product%
\[
a_{u}a_{u+1}\cdots a_{v}.
\]
It is called the \textbf{product of the numbers }$a_{k}$ \textbf{where }$k$
\textbf{ranges from }$u$ \textbf{to }$v$. When $v<u$, this product is called
\textbf{empty} and defined to be $1$.
\end{definition}

For example:%
\begin{align*}
\prod_{k=5}^{10}k  &  =5\cdot6\cdot7\cdot8\cdot9\cdot10=151\ 200;\\
\prod_{k=1}^{5}\dfrac{1}{k}  &  =\dfrac{1}{1}\cdot\dfrac{1}{2}\cdot\dfrac
{1}{3}\cdot\dfrac{1}{4}\cdot\dfrac{1}{5}=\dfrac{1}{120};\\
\prod_{k=5}^{5}\dfrac{1}{k}  &  =\dfrac{1}{5};\\
\prod_{k=6}^{5}\dfrac{1}{k}  &  =1\ \ \ \ \ \ \ \ \ \ \left(  \text{an empty
product}\right)  ;\\
\prod_{k=1}^{n}a  &  =\underbrace{aa\cdots a}_{n\text{ times}}=a^{n}%
\ \ \ \ \ \ \ \ \ \ \text{for any fixed number }a\text{ and any }%
n\in\mathbb{N};\\
\prod_{k=1}^{n}a^{k}  &  =a^{1}a^{2}\cdots a^{n}\\
&  =a^{1+2+\cdots+n}\ \ \ \ \ \ \ \ \ \ \left(
\begin{array}
[c]{c}%
\text{by one of the laws of exponents:}\\
\text{namely, the law }a^{i_{1}}a^{i_{2}}\cdots a^{i_{n}}=a^{i_{1}%
+i_{2}+\cdots+i_{n}}%
\end{array}
\right) \\
&  =a^{n\left(  n+1\right)  /2}\ \ \ \ \ \ \ \ \ \ \text{for any fixed number
}a\text{ and any }n\in\mathbb{N}\text{.}%
\end{align*}
(Yes, dear reader: The $0$-th power of any number is $1$, because it is an
empty product. This entails $0^{0}=1$ in particular. This is entirely
uncontroversial in advanced mathematics, despite whatever confusion is
typically sown by calculus textbooks.)

\bigskip

In a finite product $\prod\limits_{k=u}^{v}a_{k}$, the $k$ is called the
\textbf{product index} or the \textbf{running index}\footnote{And just like in
a sum, you can use any letter for it (unless it already stands for something
different).}, and the symbol $\prod$ is called the \textbf{product sign}. The
numbers $a_{k}$ are called the \textbf{factors} of the product. Other
terminology is analogous to the case of a finite sum (e.g., lower limit, upper
limit). Almost all rules for finite sums have analogues for finite products.
Let me only state the analogues of the \textquotedblleft splitting-off
rule\textquotedblright\ and of the rule (\ref{eq.sums.sum(a-b)}):

\begin{itemize}
\item The \textquotedblleft splitting-off rule\textquotedblright\ for
products: For any integers $u\leq v$ and any numbers $a_{u},a_{u+1}%
,\ldots,a_{v}$, we have%
\[
\prod_{k=u}^{v}a_{k}=\left(  \prod_{k=u}^{v-1}a_{k}\right)  a_{v}=a_{u}%
\prod_{k=u+1}^{v}a_{k}.
\]
This is just saying that%
\[
a_{u}a_{u+1}\cdots a_{v}=\left(  a_{u}a_{u+1}\cdots a_{v-1}\right)
a_{v}=a_{u}\left(  a_{u+1}a_{u+2}\cdots a_{v}\right)  .
\]
This rule allows us to split the first or the last factor out of a finite
product. This is important for proofs by induction.

\item The analogue of the rule (\ref{eq.sums.sum(a-b)}) for products: We have
\begin{equation}
\prod_{k=u}^{v}\left(  a_{k}/b_{k}\right)  =\left(  \prod_{k=u}^{v}%
a_{k}\right)  /\left(  \prod_{k=u}^{v}b_{k}\right)  \label{eq.prods.prod(a/b)}%
\end{equation}
for any integers $u,v$ and any numbers $a_{k},b_{k}$, as long as the numbers
$b_{k}$ are nonzero. This is an analogue of (\ref{eq.sums.sum(a-b)}), since
the multiplicative counterpart to subtraction is division. (We had to assume
that the $b_{k}$ are nonzero in order for the fractions $a_{k}/b_{k}$ to be well-defined.)
\end{itemize}

\subsection{\label{sec.sums.facs}Factorials}

Now, we define a sequence of integers that appears all over mathematics.
Recall that $\mathbb{N}=\left\{  0,1,2,\ldots\right\}  $.

\begin{definition}
\label{def.fac}For any $n\in\mathbb{N}$, we define the positive integer $n!$
(called the \textbf{factorial} of $n$, and often pronounced \textquotedblleft%
$n$ \textbf{factorial}\textquotedblright) by%
\[
n!=\prod_{k=1}^{n}k=1\cdot2\cdot\cdots\cdot n.
\]
This is the product of the first $n$ positive integers.
\end{definition}

For example,%
\begin{align*}
0!  &  =\left(  \text{empty product}\right)  =1;\\
1!  &  =1=1;\\
2!  &  =1\cdot2=2;\\
3!  &  =1\cdot2\cdot3=6;\\
4!  &  =1\cdot2\cdot3\cdot4=24;\\
5!  &  =1\cdot2\cdot3\cdot4\cdot5=120;\\
6!  &  =1\cdot2\cdot3\cdot4\cdot5\cdot6=720;\\
7!  &  =5\ 040;\\
8!  &  =40\ 320;\\
9!  &  =362\ 880;\\
10!  &  =3\ 628\,800.
\end{align*}
Note the following:

\begin{proposition}
[recursion of the factorials]\label{prop.fac.rec}For any positive integer $n$,
we have%
\[
n!=\left(  n-1\right)  !\cdot n.
\]

\end{proposition}

\begin{proof}
Let $n$ be a positive integer. Then,%
\[
n!=1\cdot2\cdot\cdots\cdot n=\underbrace{\left(  1\cdot2\cdot\cdots
\cdot\left(  n-1\right)  \right)  }_{=\left(  n-1\right)  !}\cdot\,n=\left(
n-1\right)  !\cdot n.
\]

\end{proof}

\begin{exercise}
\label{exe.23hw1.1b}Prove that%
\[
1\cdot1!+2\cdot2!+3\cdot3!+\cdots+n\cdot n!=\left(  n+1\right)  !-1
\]
for each $n\in\mathbb{N}$.
\end{exercise}

(Meanwhile, there is no such simple formula for $1!+2!+3!+\cdots+n!$. Not
every sum can be simplified!)

\begin{exercise}
\label{exe.23hw2.1}\textbf{(a)} Prove that
\[
\prod_{i=2}^{n}\left(  1-\dfrac{1}{i^{2}}\right)  =\dfrac{n+1}{2n}%
\]
for each positive integer $n$. \medskip

\textbf{(b)} Find and prove a closed-form expression (i.e., no $\prod$ or
$\sum$ signs) for%
\[
\prod_{i=2}^{n}\left(  1-\dfrac{1}{i}\right)  .
\]

\end{exercise}

\begin{exercise}
\label{exe.23finfin.2a}Prove that%
\[
\prod_{k=0}^{n}k!=\prod_{k=1}^{n}k!=\prod_{k=1}^{n}k^{n-k+1}%
\ \ \ \ \ \ \ \ \ \ \text{for each }n\in\mathbb{N}.
\]

\end{exercise}

\begin{exercise}
\label{exe.23hw2.2}Prove that%
\[
\prod_{i=1}^{n}\left(  i!\cdot i^{i}\right)  =n!^{n+1}%
\ \ \ \ \ \ \ \ \ \ \text{for each }n\in\mathbb{N}.
\]

\end{exercise}

\begin{exercise}
\label{exe.23mt1.5}Let $\left(  a_{0},a_{1},a_{2},\ldots\right)  $ be a
sequence of integers defined recursively by%
\[
a_{n}=1+a_{0}a_{1}\cdots a_{n-1}\ \ \ \ \ \ \ \ \ \ \text{for all }n\geq0.
\]
(In particular, $a_{0}=1+\underbrace{a_{0}a_{1}\cdots a_{0-1}}_{=\left(
\text{empty product}\right)  =1}=1+1=2$.) Here are the first few entries of
this sequence:%
\[%
\begin{tabular}
[c]{|c||c|c|c|c|c|c|c|}\hline
$n$ & $0$ & $1$ & $2$ & $3$ & $4$ & $5$ & $6$\\\hline
$a_{n}$ & $2$ & $3$ & $7$ & $43$ & $1807$ & $3263443$ & $10650056950807$%
\\\hline
\end{tabular}
\ \ \
\]
(notice the astronomical growth!). \medskip

\textbf{(a)} Prove that%
\[
a_{n+1}=a_{n}^{2}-a_{n}+1\ \ \ \ \ \ \ \ \ \ \text{for each }n\geq0.
\]

\textbf{(b)} Prove that
\[
\dfrac{1}{a_{0}}+\dfrac{1}{a_{1}}+\cdots+\dfrac{1}{a_{n-1}}=1-\dfrac{1}%
{a_{n}-1}\ \ \ \ \ \ \ \ \ \ \text{for each }n\geq0.
\]

\end{exercise}

\begin{exercise}
Define a sequence $\left(  s_{0},s_{1},s_{2},\ldots\right)  $ of integers
recursively by%
\[
s_{0}=1\ \ \ \ \ \ \ \ \ \ \text{and}\ \ \ \ \ \ \ \ \ \ s_{n}=2+s_{0}%
s_{1}\cdots s_{n-1}\ \ \ \ \ \ \ \ \ \ \text{for all }n\geq1.
\]
(Thus, $s_{1}=3$ and $s_{2}=5$ and $s_{3}=17$.) \medskip

\textbf{(a)} Prove that $s_{n}=s_{n-1}^{2}-2s_{n-1}+2$ for all $n\geq2$.
\medskip

\textbf{(b)} Prove that%
\[
s_{n}=2^{2^{n-1}}+1\ \ \ \ \ \ \ \ \ \ \text{for all }n\geq1.
\]
(Keep in mind that a \textquotedblleft power tower\textquotedblright\ of the
form $a^{b^{c}}$ has to be understood as $a^{\left(  b^{c}\right)  }$, not as
$\left(  a^{b}\right)  ^{c}$.) \medskip

\textbf{(c)} Does this equality $s_{n}=2^{2^{n-1}}+1$ also hold for $n=0$ ?
\end{exercise}

\subsection{\label{sec.sums.bin-def}Binomial coefficients: Definition}

We shall now define one of the most important families of numbers in mathematics:

\begin{definition}
\label{def.binom}Let $n$ and $k$ be any numbers. Then, we define a number
$\dbinom{n}{k}$ as follows:

\begin{itemize}
\item If $k\in\mathbb{N}$, then we set%
\[
\dbinom{n}{k}:=\dfrac{n\left(  n-1\right)  \left(  n-2\right)  \cdots\left(
n-k+1\right)  }{k!}.
\]
(The numerator here is the product of $k$ factors, where the first factor is
$n$ and each further factor is $1$ smaller than the previous. You can also
write this product as $\prod\limits_{i=0}^{k-1}\left(  n-i\right)  $.)

\item If $k\notin\mathbb{N}$, then we set%
\[
\dbinom{n}{k}:=0.
\]

\end{itemize}

The number $\dbinom{n}{k}$ is called \textquotedblleft$n$ \textbf{choose }%
$k$\textquotedblright, and is known as the \textbf{binomial coefficient} of
$n$ and $k$. Do not mistake the notation $\dbinom{n}{k}$ for a vector $\left(
\begin{array}
[c]{c}%
n\\
k
\end{array}
\right)  $.
\end{definition}

\begin{example}
\label{exa.binom.nC3-etc}For any number $n$, we have%
\begin{align*}
\dbinom{n}{3}  &  =\dfrac{n\left(  n-1\right)  \left(  n-2\right)  }%
{3!}=\dfrac{n\left(  n-1\right)  \left(  n-2\right)  }{6};\\
\dbinom{n}{2}  &  =\dfrac{n\left(  n-1\right)  }{2!}=\dfrac{n\left(
n-1\right)  }{2};\\
\dbinom{n}{1}  &  =\dfrac{n}{1!}=n;\\
\dbinom{n}{0}  &  =\dfrac{\left(  \text{empty product}\right)  }{0!}=\dfrac
{1}{1}=1;\\
\dbinom{n}{2.5}  &  =0\ \ \ \ \ \ \ \ \ \ \left(  \text{since }2.5\notin%
\mathbb{N}\right)  ;\\
\dbinom{n}{-1}  &  =0\ \ \ \ \ \ \ \ \ \ \left(  \text{since }-1\notin%
\mathbb{N}\right)  .
\end{align*}
For any $k\in\mathbb{N}$, we have%
\begin{align*}
\dbinom{0}{k}  &  =\dfrac{0\left(  0-1\right)  \left(  0-2\right)
\cdots\left(  0-k+1\right)  }{k!}=%
\begin{cases}
1, & \text{if }k=0;\\
0, & \text{if }k\neq0
\end{cases}
\\
&  \ \ \ \ \ \ \ \ \ \ \ \ \ \ \ \ \ \ \ \ \left(
\begin{array}
[c]{c}%
\text{since the product }0\left(  0-1\right)  \left(  0-2\right)
\cdots\left(  0-k+1\right) \\
\text{is empty for }k=0\text{, and otherwise has a}\\
\text{factor equal to }0\text{ and thus must be }0
\end{array}
\right) \\
\dbinom{-1}{k}  &  =\dfrac{\left(  -1\right)  \left(  -2\right)  \left(
-3\right)  \cdots\left(  -k\right)  }{k!}=\left(  -1\right)  ^{k}%
\cdot\underbrace{\dfrac{1\cdot2\cdot\cdots\cdot k}{k!}}_{=1}=\left(
-1\right)  ^{k}.
\end{align*}

\end{example}

Let us tabulate the values of $\dbinom{n}{k}$ for nonnegative integers $n$ and
$k$:%
\[%
\begin{tabular}
[c]{|c||c|c|c|c|c|c|c|}\hline
& $k=0$ & $k=1$ & $k=2$ & $k=3$ & $k=4$ & $k=5$ & $k=6$\\\hline\hline
$n=0$ & $1$ & $0$ & $0$ & $0$ & $0$ & $0$ & $0$\\\hline
$n=1$ & $1$ & $1$ & $0$ & $0$ & $0$ & $0$ & $0$\\\hline
$n=2$ & $1$ & $2$ & $1$ & $0$ & $0$ & $0$ & $0$\\\hline
$n=3$ & $1$ & $3$ & $3$ & $1$ & $0$ & $0$ & $0$\\\hline
$n=4$ & $1$ & $4$ & $6$ & $4$ & $1$ & $0$ & $0$\\\hline
$n=5$ & $1$ & $5$ & $10$ & $10$ & $5$ & $1$ & $0$\\\hline
$n=6$ & $1$ & $6$ & $15$ & $20$ & $15$ & $6$ & $1$\\\hline
\end{tabular}
\ \ \ \ .
\]

What patterns can we spot in this table? (We are ignoring negative and
non-integer $n$'s for now.)

The following is probably the most visible one:

\begin{proposition}
\label{prop.binom.0}Let $n\in\mathbb{N}$ and $k>n$. Then, $\dbinom{n}{k}=0$.
\end{proposition}

\begin{proof}
If $k\notin\mathbb{N}$, then this is clear by definition. Otherwise, again by
definition, we have%
\[
\dbinom{n}{k}=\dfrac{n\left(  n-1\right)  \left(  n-2\right)  \cdots\left(
n-k+1\right)  }{k!}=\dfrac{0}{k!}%
\]
(since the product $n\left(  n-1\right)  \left(  n-2\right)  \cdots\left(
n-k+1\right)  $ has a factor of $n-n=0$, and thus is $0$). For example, for
$n=3$ and $k=6$, we have%
\[
\dbinom{3}{6}=\dfrac{3\cdot2\cdot1\cdot0\cdot\left(  -1\right)  \cdot\left(
-2\right)  }{6!}=\dfrac{0}{6!}=0.
\]
In the general case, we thus find $\dbinom{n}{k}=\dfrac{0}{k!}=0$.
\end{proof}

\begin{remark}
Note that Proposition \ref{prop.binom.0} would \textbf{not} hold without the
$n\in\mathbb{N}$ assumption. For example,%
\[
\dbinom{1.5}{3}=\dfrac{1.5\cdot0.5\cdot\left(  -0.5\right)  }{3!}%
\neq0\ \ \ \ \ \ \ \ \ \ \text{even though }3>1.5.
\]
The product in the numerator is not $0$, since it \textquotedblleft jumps
over\textquotedblright\ the $0$ factor.
\end{remark}

Proposition \ref{prop.binom.0} explains why our above table of $\dbinom{n}{k}$
has so many zeroes in it. More precisely, it tells us that all entries above
the main diagonal of the table are zeroes (no matter how many more rows and
columns we add). Thus, we can redraw our table as a triangular table (and fill
in a few more rows while at that):%
\[
\hspace{-1.0038pc}\hspace{-0.3cm}%
\begin{tabular}
[c]{ll|ccccccccccccccccc}%
\vphantom{$\overset{k=0}{\swarrow}$} &  & $\phantom{20}$ & $\phantom{20}$ &
$\phantom{20}$ & $\phantom{20}$ & $\phantom{20}$ & $\phantom{20}$ &  &  &  &
$\overset{k=0}{\swarrow}$ & $\phantom{20}$ & $\phantom{20}$ & $\phantom{20}$ &
$\phantom{20}$ & $\phantom{20}$ & $\phantom{20}$ & $\phantom{20}$\\
\vphantom{$\overset{k=0}{\swarrow}$}$n=0$ & $\rightarrow$ &  &  &  &  &  &  &
&  & $1$ &  & $\overset{k=1}{\swarrow}$ &  &  &  &  &  & \\
\vphantom{$\overset{k=0}{\swarrow}$}$n=1$ & $\rightarrow$ &  &  &  &  &  &  &
& $1$ &  & $1$ &  & $\overset{k=2}{\swarrow}$ &  &  &  &  & \\
\vphantom{$\overset{k=0}{\swarrow}$}$n=2$ & $\rightarrow$ &  &  &  &  &  &  &
$1$ &  & $2$ &  & $1$ &  & $\overset{k=3}{\swarrow}$ &  &  &  & \\
\vphantom{$\overset{k=0}{\swarrow}$}$n=3$ & $\rightarrow$ &  &  &  &  &  & $1$
&  & $3$ &  & $3$ &  & $1$ &  & $\overset{k=4}{\swarrow}$ &  &  & \\
\vphantom{$\overset{k=0}{\swarrow}$}$n=4$ & $\rightarrow$ &  &  &  &  & $1$ &
& $4$ &  & $6$ &  & $4$ &  & $1$ &  & $\overset{k=5}{\swarrow}$ &  & \\
\vphantom{$\overset{k=0}{\swarrow}$}$n=5$ & $\rightarrow$ &  &  &  & $1$ &  &
$5$ &  & $10$ &  & $10$ &  & $5$ &  & $1$ &  & $\overset{k=6}{\swarrow}$ & \\
\vphantom{$\overset{k=0}{\swarrow}$}$n=6$ & $\rightarrow$ &  &  & $1$ &  & $6$
&  & $15$ &  & $20$ &  & $15$ &  & $6$ &  & $1$ &  & $\overset{k=7}{\swarrow}%
$\\
\vphantom{$\overset{k=0}{\swarrow}$}$n=7$ & $\rightarrow$ &  & $1$ &  & $7$ &
& $21$ &  & $35$ &  & $35$ &  & $21$ &  & $7$ &  & $1$ & \\
\vphantom{$\overset{k=0}{\swarrow}$}$n=8$ & $\rightarrow$ & $1$ &  & $8$ &  &
$28$ &  & $56$ &  & $70$ &  & $56$ &  & $28$ &  & $8$ &  & $1$%
\end{tabular}
\ \ \ \ \ \ \ \ \ \ \
\]
This table is known as \textbf{Pascal's triangle}, and has a variety of
wonderful properties. Here are just a few:

\begin{itemize}
\item \textbf{Pascal's identity, aka the recurrence of the binomial
coefficients:} For any numbers $n$ and $k$, we have%
\[
\dbinom{n}{k}=\dbinom{n-1}{k-1}+\dbinom{n-1}{k}.
\]

\item \textbf{Symmetry of binomial coefficients:} For any $n\in\mathbb{N}$ and
any $k$, we have $\dbinom{n}{k}=\dbinom{n}{n-k}$.

\item We have $\dbinom{n}{n}=1$ for each $n\in\mathbb{N}$.

\item \textbf{Integrality of binomial coefficients:} For any $n\in\mathbb{Z}$
and any $k$, we have $\dbinom{n}{k}\in\mathbb{Z}$.
\end{itemize}

In the next section, we will prove these four propositions and more.

\subsection{\label{sec.ind.bin-props}Binomial coefficients: Properties}

\subsubsection{Pascal's identity}

We begin with the most important property of binomial coefficients:

\begin{theorem}
[Pascal's identity, aka the recurrence of the binomial coefficients]%
\label{thm.binom.rec}For any numbers $n$ and $k$, we have%
\[
\dbinom{n}{k}=\dbinom{n-1}{k-1}+\dbinom{n-1}{k}.
\]

\end{theorem}

\begin{example}
For $n=7$ and $k=3$, this is claiming that $\dbinom{7}{3}=\dbinom{6}%
{2}+\dbinom{6}{3}$, which explicitly is saying that $35=15+20$.

But note that Theorem \ref{thm.binom.rec} also can be applied when $n$ or $k$
is negative or non-integer.
\end{example}

\begin{proof}
[Proof of Theorem \ref{thm.binom.rec}.]Let $n$ and $k$ be two numbers. We are
in one of the following three cases:

\textit{Case 1:} The number $k$ is a positive integer.

\textit{Case 2:} We have $k=0$.

\textit{Case 3:} None of the above.

Let us first consider Case 1 (this is the interesting case). Here, $k$ is a
positive integer, so that both $k$ and $k-1$ belong to $\mathbb{N}$. The
definition of binomial coefficients therefore yields the three formulas%
\begin{align*}
\dbinom{n}{k}  &  =\dfrac{n\left(  n-1\right)  \left(  n-2\right)
\cdots\left(  n-k+1\right)  }{k!};\\
\dbinom{n-1}{k-1}  &  =\dfrac{\left(  n-1\right)  \left(  n-2\right)  \left(
n-3\right)  \cdots\left(  \left(  n-1\right)  -\left(  k-1\right)  +1\right)
}{\left(  k-1\right)  !}\\
&  =\dfrac{\left(  n-1\right)  \left(  n-2\right)  \left(  n-3\right)
\cdots\left(  n-k+1\right)  }{\left(  k-1\right)  !};\\
\dbinom{n-1}{k}  &  =\dfrac{\left(  n-1\right)  \left(  n-2\right)  \left(
n-3\right)  \cdots\left(  \left(  n-1\right)  -k+1\right)  }{k!}\\
&  =\dfrac{\left(  n-1\right)  \left(  n-2\right)  \left(  n-3\right)
\cdots\left(  n-k\right)  }{k!}.
\end{align*}
Let us set $a:=\left(  n-1\right)  \left(  n-2\right)  \left(  n-3\right)
\cdots\left(  n-k+1\right)  $ (this is the common factor in the numerators of
all these three formulas). Then, these three formulas can be rewritten as%
\begin{align}
\dbinom{n}{k}  &  =\dfrac{na}{k!};\label{pf.thm.binom.rec.4}\\
\dbinom{n-1}{k-1}  &  =\dfrac{a}{\left(  k-1\right)  !}%
;\label{pf.thm.binom.rec.5}\\
\dbinom{n-1}{k}  &  =\dfrac{a\left(  n-k\right)  }{k!}.
\label{pf.thm.binom.rec.6}%
\end{align}

But the claim that we are trying to prove is
\[
\dbinom{n}{k}=\dbinom{n-1}{k-1}+\dbinom{n-1}{k}.
\]
Using the formulas (\ref{pf.thm.binom.rec.4}), (\ref{pf.thm.binom.rec.5}) and
(\ref{pf.thm.binom.rec.6}), this can be rewritten as%
\[
\dfrac{na}{k!}=\dfrac{a}{\left(  k-1\right)  !}+\dfrac{a\left(  n-k\right)
}{k!}.
\]
Multiplying both sides by $k!$, we can transform this into
\[
na=a\cdot\dfrac{k!}{\left(  k-1\right)  !}+a\left(  n-k\right)  .
\]
Since $\dfrac{k!}{\left(  k-1\right)  !}=k$ (because the recursion of the
factorials (i.e., Proposition \ref{prop.fac.rec}) yields $k!=\left(
k-1\right)  !\cdot k$), we can simplify this further to%
\[
na=a\cdot k+a\left(  n-k\right)  ,
\]
which is obviously true. Thus, our claim is proved in Case 1.

Now, we consider Case 2. In this case, $k=0$. Our claim%
\[
\dbinom{n}{k}=\dbinom{n-1}{k-1}+\dbinom{n-1}{k}%
\]
thus rewrites as%
\[
\dbinom{n}{0}=\dbinom{n-1}{0-1}+\dbinom{n-1}{0},
\]
which again is true (because Example \ref{exa.binom.nC3-etc} shows that
$\dbinom{n}{0}=1$ and $\dbinom{n-1}{0}=1$ and $\dbinom{n-1}{0-1}=\dbinom
{n-1}{-1}=0$).

Finally, we consider Case 3. In this case, $k$ is neither a positive integer
nor $0$. Hence, $k\notin\mathbb{N}$. Thus, $k-1\notin\mathbb{N}$ as well.
Hence, in our claim%
\[
\dbinom{n}{k}=\dbinom{n-1}{k-1}+\dbinom{n-1}{k},
\]
all three binomial coefficients are $0$ (since a binomial coefficient
$\dbinom{m}{\ell}$ is $0$ by definition when $\ell\notin\mathbb{N}$). Thus,
again, the claim is true (since $0=0+0$).

We have now proved Theorem \ref{thm.binom.rec} in all three cases; thus, it is
always true.
\end{proof}

Pascal's identity is highly useful for proving properties of binomial
coefficients $\dbinom{n}{k}$ by induction on $n$. (We will see an example of
this very soon, in the proof of Theorem \ref{thm.binom.N}.)

Pascal's identity shows that every entry of Pascal's triangle (except the $1$
at the apex) equals the sum of the two entries directly above it (i.e., of the
entry one step northwest of it and the entry one step northeast of it). But it
also applies to binomial coefficients that are not (commonly) considered to be
part of Pascal's triangle, such as $\dbinom{-3}{5}=\dbinom{-4}{4}+\dbinom
{-4}{5}$ and $\dbinom{3.2}{2}=\dbinom{2.2}{1}+\dbinom{2.2}{2}$.

\subsubsection{The factorial formula}

Binomial coefficients $\dbinom{n}{k}$ make sense for arbitrary numbers $n$ and
$k$. However, when $n$ and $k$ are nonnegative integers with $k\leq n$ (that
is, when $n\in\mathbb{N}$ and $k\in\left\{  0,1,\ldots,n\right\}  $), there is
a particularly simple formula for $\dbinom{n}{k}$, known as the
\textbf{factorial formula}:

\begin{theorem}
[factorial formula]\label{thm.binom.fac}Let $n\in\mathbb{N}$ and $k\in\left\{
0,1,\ldots,n\right\}  $. Then,%
\[
\dbinom{n}{k}=\dfrac{n!}{k!\cdot\left(  n-k\right)  !}.
\]

\end{theorem}

\begin{proof}
The definition of $\dbinom{n}{k}$ yields%
\[
\dbinom{n}{k}=\dfrac{n\left(  n-1\right)  \left(  n-2\right)  \cdots\left(
n-k+1\right)  }{k!}.
\]
Multiplying both sides by $k!$, we obtain%
\begin{align*}
k!\cdot\dbinom{n}{k}  &  =n\left(  n-1\right)  \left(  n-2\right)
\cdots\left(  n-k+1\right) \\
&  =\left(  n-k+1\right)  \left(  n-k+2\right)  \left(  n-k+3\right)  \cdots
n\\
&  =\dfrac{1\cdot2\cdot\cdots\cdot n}{1\cdot2\cdot\cdots\cdot\left(
n-k\right)  }\\
&  \ \ \ \ \ \ \ \ \ \ \ \ \ \ \ \ \ \ \ \ \left(
\begin{array}
[c]{c}%
\text{since }n-k+1,\ n-k+2,\ n-k+3,\ \ldots,\ n\text{ are}\\
\text{precisely the factors of the product }1\cdot2\cdot\cdots\cdot n\\
\text{that do not appear in the product }1\cdot2\cdot\cdots\cdot\left(
n-k\right)
\end{array}
\right) \\
&  =\dfrac{n!}{\left(  n-k\right)  !}.
\end{align*}
Dividing this by $k!$, we obtain%
\[
\dbinom{n}{k}=\dfrac{n!}{\left(  n-k\right)  !}/k!=\dfrac{n!}{k!\cdot\left(
n-k\right)  !}.
\]
This proves the factorial formula.
\end{proof}

\begin{warning}
The factorial formula can be used to compute $\dbinom{10}{4}$ for example, but
it cannot be used to compute $\dbinom{-1}{3}$ or $\dbinom{1.2}{2}$ (because
the \textquotedblleft$n\in\mathbb{N}$ and $k\in\left\{  0,1,\ldots,n\right\}
$\textquotedblright\ conditions in the factorial formula are not satisfied
here). It is thus not as general as the definition of binomial coefficients!
\end{warning}

\begin{exercise}
\label{exe.23mt1.4}Let $n\in\mathbb{N}$. Prove that%
\[
\prod_{i=0}^{n}\dbinom{2i}{i}=2^{n}\prod_{i=0}^{n}\dbinom{n+i}{n-i}.
\]

\end{exercise}

\subsubsection{The symmetry of binomial coefficients}

Here is another property of Pascal's triangle: It has a vertical axis of
symmetry, meaning that the entries to the left of this axis equal the
corresponding entries to the right of the axis. Let us state this more precisely:

\begin{theorem}
[symmetry of Pascal's triangle]\label{thm.binom.sym}Let $n\in\mathbb{N}$, and
let $k$ be any number. Then,%
\[
\dbinom{n}{k}=\dbinom{n}{n-k}.
\]

\end{theorem}

\begin{proof}
We are in one of the following four cases:

\textit{Case 1:} We have $k\in\left\{  0,1,\ldots,n\right\}  $.

\textit{Case 2:} We have $k<0$.

\textit{Case 3:} We have $k>n$.

\textit{Case 4:} The number $k$ is not an integer.

Let us first consider Case 1. In this case, we have $k\in\left\{
0,1,\ldots,n\right\}  $ and thus also $n-k\in\left\{  0,1,\ldots,n\right\}  $.
Since $k\in\left\{  0,1,\ldots,n\right\}  $, we can apply the factorial
formula to obtain%
\[
\dbinom{n}{k}=\dfrac{n!}{k!\cdot\left(  n-k\right)  !}.
\]
Since $n-k\in\left\{  0,1,\ldots,n\right\}  $, we can also apply the factorial
formula to $n-k$ instead of $k$, and thus we find%
\[
\dbinom{n}{n-k}=\dfrac{n!}{\left(  n-k\right)  !\cdot\left(  n-\left(
n-k\right)  \right)  !}=\dfrac{n!}{\left(  n-k\right)  !\cdot k!}=\dfrac
{n!}{k!\cdot\left(  n-k\right)  !}.
\]
The right hand sides of these two equalities are equal. Thus, the left hand
sides are equal as well. This proves $\dbinom{n}{k}=\dbinom{n}{n-k}$ in Case 1.

In Case 2, we have $\dbinom{n}{k}=0$ by definition (since $k<0$ entails
$k\notin\mathbb{N}$), whereas $\dbinom{n}{n-k}=0$ by Proposition
\ref{prop.binom.0} (since $n\in\mathbb{N}$ and $n-\underbrace{k}_{<0}>n$).
This proves $\dbinom{n}{k}=\dbinom{n}{n-k}$ in Case 2.

Case 3 is analogous to Case 2, except that $k$ and $n-k$ trade places.

In Case 4, both $\dbinom{n}{k}$ and $\dbinom{n}{n-k}$ are $0$ by definition
(since neither $k$ nor $n-k$ belongs to $\mathbb{N}$).

Thus, $\dbinom{n}{k}=\dbinom{n}{n-k}$ is proved in all four cases, so that
Theorem \ref{thm.binom.sym} follows.
\end{proof}

Alternatively, Theorem \ref{thm.binom.sym} could have been proved by induction
on $n$.

\begin{warning}
Theorem \ref{thm.binom.sym} holds only for $n\in\mathbb{N}$. For $n=-1$ and
$k=0$, it is false (since $\dbinom{-1}{0}=1$ but $\dbinom{-1}{-1-0}=0$).
\end{warning}

One corollary of Theorem \ref{thm.binom.sym} is the fact that the
\textquotedblleft right border\textquotedblright\ of Pascal's triangle is
filled with $1$'s:

\begin{corollary}
\label{cor.binom.nn}For any $n\in\mathbb{N}$, we have $\dbinom{n}{n}=1$.
\end{corollary}

\begin{proof}
For any $n\in\mathbb{N}$, Theorem \ref{thm.binom.sym} (applied to $k=n$)
yields%
\[
\dbinom{n}{n}=\dbinom{n}{n-n}=\dbinom{n}{0}=1.
\]

\end{proof}

\begin{warning}
Corollary \ref{cor.binom.nn} does not hold for negative (or non-integer) $n$.
\end{warning}

\subsubsection{Pascal's triangle consists of integers}

The perhaps most surprising pattern in Pascal's triangle is that all its
entries are integers! It is tempting to take this for granted, but this is not
at all obvious from our definition of $\dbinom{n}{k}$ as a fraction.
Nevertheless, we can now prove it without much trouble:

\begin{theorem}
\label{thm.binom.N}For any $n\in\mathbb{N}$ and any number $k$, we have
$\dbinom{n}{k}\in\mathbb{N}$.
\end{theorem}

\begin{proof}
We induct on $n$.

\textit{Base case:} Theorem \ref{thm.binom.N} holds for $n=0$, since any
number $k$ satisfies%
\begin{align*}
\dbinom{0}{k}  &  =%
\begin{cases}
1, & \text{if }k=0;\\
0, & \text{if }k\neq0
\end{cases}
\ \ \ \ \ \ \ \ \ \ \left(  \text{easy to see from the definition}\right) \\
&  \in\mathbb{N}.
\end{align*}

\textit{Induction step:} We will make an induction step from $n-1$ to $n$
(instead of the more conventional step from $n$ to $n+1$). So we fix a
positive integer $n$, and we assume (as the induction hypothesis) that Theorem
\ref{thm.binom.N} holds for $n-1$ instead of $n$. In other words, we assume
that%
\begin{equation}
\dbinom{n-1}{k}\in\mathbb{N}\ \ \ \ \ \ \ \ \ \ \text{for all numbers }k.
\label{pf.thm.binom.N.IH}%
\end{equation}

Our goal is to prove that Theorem \ref{thm.binom.N} also holds for $n$. In
other words, we must prove that%
\[
\dbinom{n}{k}\in\mathbb{N}\ \ \ \ \ \ \ \ \ \ \text{for all numbers }k.
\]

But this is easy: Pascal's identity yields%
\[
\dbinom{n}{k}=\underbrace{\dbinom{n-1}{k-1}}_{\substack{\in\mathbb{N}%
\\\text{(by (\ref{pf.thm.binom.N.IH}),}\\\text{applied to }k-1\\\text{instead
of }k\text{)}}}+\underbrace{\dbinom{n-1}{k}}_{\substack{\in\mathbb{N}%
\\\text{(by (\ref{pf.thm.binom.N.IH}))}}}\in\mathbb{N}%
\ \ \ \ \ \ \ \ \ \ \text{for all numbers }k\text{.}%
\]
So the induction step is complete, and the theorem is proved.
\end{proof}

\bigskip

Theorem \ref{thm.binom.N} is crying for a better explanation: Certainly, a
number shouldn't belong to $\mathbb{N}$ for no reason! (Actually, it can, but
let's be optimistic.) Such an explanation does indeed exist:

\begin{theorem}
[combinatorial interpretation of binomial coefficients]%
\label{thm.binom.comb-1}Let $n\in\mathbb{N}$, and let $k$ be any number. Let
$A$ be any $n$-element set. (Here, \textquotedblleft$n$-element
set\textquotedblright\ means a set that has exactly $n$ distinct elements. For
example, $\left\{  2,6,11\right\}  $ is a $3$-element set, and this does not
change if I rewrite this set as $\left\{  2,6,2,11\right\}  $. Note that the
sets $\left\{  2,3\right\}  $ and $\left\{  3,2\right\}  $ are identical,
since a set doesn't care how its elements are ordered.)

Then,%
\[
\dbinom{n}{k}\text{ is the number of }k\text{-element subsets of }A.
\]

\end{theorem}

\begin{example}
Let $n=4$ and $k=2$ and $A=\left\{  1,2,3,4\right\}  $. Then, the $2$-element
subsets of $A$ are%
\[
\left\{  1,2\right\}  ,\ \ \left\{  1,3\right\}  ,\ \ \left\{  1,4\right\}
,\ \ \left\{  2,3\right\}  ,\ \ \left\{  2,4\right\}  ,\ \ \left\{
3,4\right\}  .
\]
So their number is $6$. And this agrees with Theorem \ref{thm.binom.comb-1},
since $\dbinom{n}{k}=\dbinom{4}{2}=6$.

Another example: The $3$-element subsets of $\left\{  1,2,3,4,5\right\}  $ are%
\begin{align*}
&  \left\{  1,2,3\right\}  ,\ \left\{  1,2,4\right\}  ,\ \left\{
1,2,5\right\}  ,\ \left\{  1,3,4\right\}  ,\ \left\{  1,3,5\right\}
,\ \left\{  1,4,5\right\}  ,\\
&  \left\{  2,3,4\right\}  ,\ \left\{  2,3,5\right\}  ,\ \left\{
2,4,5\right\}  ,\ \left\{  3,4,5\right\}  .
\end{align*}
There are $10$ of them, just as Theorem \ref{thm.binom.comb-1} predicts (since
$\dbinom{5}{3}=10$).
\end{example}

We will prove Theorem \ref{thm.binom.comb-1} later in this course (see Theorem
\ref{thm.count.subsets.size.S}), as we learn more about finite sets and their
sizes. Note that the $k$-element subsets of $A$ are also known as
\textbf{combinations without replacement}. Theorem \ref{thm.binom.comb-1} also
explains why $\dbinom{n}{k}$ is called \textquotedblleft$n$ choose
$k$\textquotedblright: After all, a $k$-element subset of $A$ is a
\textquotedblleft choice\textquotedblright\ of $k$ distinct elements (without
regard for order) from $A$.

Note again that Theorem \ref{thm.binom.comb-1} says nothing about binomial
coefficients $\dbinom{n}{k}$ with $n\notin\mathbb{N}$, since a number
$n\notin\mathbb{N}$ cannot be the size of a set. So Theorem
\ref{thm.binom.comb-1} explains why $\dbinom{5}{2}$ is an integer, but does
not explain why $\dbinom{-5}{2}$ is an integer.

\subsubsection{Upper negation}

Here is another property of binomial coefficients, called the \textbf{upper
negation formula}:

\begin{theorem}
[upper negation formula]\label{thm.binom.upneg}For any numbers $n$ and
$k\in\mathbb{Z}$, we have%
\[
\dbinom{-n}{k}=\left(  -1\right)  ^{k}\dbinom{n+k-1}{k}.
\]

\end{theorem}

\begin{proof}
If $k\notin\mathbb{N}$, then this is clear because both binomial coefficients
are $0$ by definition.

Thus, we only need to prove the theorem in the case when $k\in\mathbb{N}$.

In this case, the definition of binomial coefficients yields%
\begin{align*}
\dbinom{-n}{k}  &  =\dfrac{\left(  -n\right)  \left(  -n-1\right)  \left(
-n-2\right)  \cdots\left(  -n-k+1\right)  }{k!}\\
&  =\left(  -1\right)  ^{k}\cdot\dfrac{n\left(  n+1\right)  \left(
n+2\right)  \cdots\left(  n+k-1\right)  }{k!}%
\end{align*}
(here, we factored out all the minus signs from the numerator) and%
\begin{align*}
\dbinom{n+k-1}{k}  &  =\dfrac{\left(  n+k-1\right)  \left(  n+k-2\right)
\left(  n+k-3\right)  \cdots n}{k!}\\
&  =\dfrac{n\left(  n+1\right)  \left(  n+2\right)  \cdots\left(
n+k-1\right)  }{k!}.
\end{align*}
Comparing these equalities, we find $\dbinom{-n}{k}=\left(  -1\right)
^{k}\dbinom{n+k-1}{k}$. This proves the theorem.
\end{proof}

\begin{corollary}
\label{cor.binom.Z}For any $n\in\mathbb{Z}$ and any number $k$, we have
$\dbinom{n}{k}\in\mathbb{Z}$.
\end{corollary}

\begin{proof}
If $n\geq0$, then this has already been proved in Theorem \ref{thm.binom.N}.

If $k\notin\mathbb{N}$, then this is clear because $\dbinom{n}{k}=0$.

In the remaining case, use the upper negation formula. Details are left to the
reader (see \cite[Theorem 1.3.16]{19fco}).
\end{proof}

Note that (as we said above) \textquotedblleft negative\textquotedblright%
\ binomial coefficients such as $\dbinom{-3}{5}=-21$ have no immediate
combinatorial meaning, because there is no such thing as a $\left(  -3\right)
$-element set. Nevertheless, they have uses in algebra and elsewhere.

\subsubsection{Finding Fibonacci numbers in Pascal's triangle}

The binomial coefficients are related to the Fibonacci numbers:

\begin{theorem}
\label{thm.binom.fib}For any $n\in\mathbb{N}$, the Fibonacci number $f_{n+1}$
is%
\begin{align*}
f_{n+1}  &  =\dbinom{n-0}{0}+\dbinom{n-1}{1}+\dbinom{n-2}{2}+\cdots
+\dbinom{n-n}{n}\\
&  =\sum_{k=0}^{n}\dbinom{n-k}{k}.
\end{align*}

\end{theorem}

For example, for $n=7$, this is saying that%
\begin{align*}
f_{8}  &  =\dbinom{7-0}{0}+\dbinom{7-1}{1}+\dbinom{7-2}{2}+\cdots+\dbinom
{7-7}{7}\\
&  =1+6+10+4+0+0+0+0=21.
\end{align*}

We will prove Theorem \ref{thm.binom.fib} in Chapter \ref{chp.enum2} (as
Corollary \ref{cor.fib.sum-binom}) using enumerative combinatorics. You can
find proofs of Theorem \ref{thm.binom.fib} in \cite[\S 15]{Vorobi02} and in
\cite[\S 1.4.5, proof of Proposition 1.3.32]{19fco} as well.

\subsection{The binomial formula}

One of the most important properties of binomial coefficients (which,
incidentally, explains their name) is the \textbf{binomial formula}:

\begin{theorem}
[binomial formula, aka binomial theorem]\label{thm.binom.bin}Let $a$ and $b$
be any numbers, and let $n\in\mathbb{N}$. Then,%
\begin{equation}
\left(  a+b\right)  ^{n}=\sum_{k=0}^{n}\dbinom{n}{k}a^{k}b^{n-k}.
\label{eq.thm.binom.bin.1}%
\end{equation}
Restating this without the summation sign:%
\[
\left(  a+b\right)  ^{n}=\dbinom{n}{0}a^{0}b^{n}+\dbinom{n}{1}a^{1}%
b^{n-1}+\dbinom{n}{2}a^{2}b^{n-2}+\cdots+\dbinom{n}{n}a^{n}b^{0}.
\]

Equivalently:%
\begin{equation}
\left(  a+b\right)  ^{n}=\sum_{k=0}^{n}\dbinom{n}{k}a^{n-k}b^{k}.
\label{eq.thm.binom.bin.3}%
\end{equation}

\end{theorem}

\begin{example}
For $n=5$, the formula (\ref{eq.thm.binom.bin.1}) is saying that%
\begin{align*}
&  \left(  a+b\right)  ^{5}\\
&  =\sum_{k=0}^{5}\dbinom{5}{k}a^{k}b^{5-k}\\
&  =\dbinom{5}{0}a^{0}b^{5}+\dbinom{5}{1}a^{1}b^{4}+\dbinom{5}{2}a^{2}%
b^{3}+\dbinom{5}{3}a^{3}b^{2}+\dbinom{5}{4}a^{4}b^{1}+\dbinom{5}{5}a^{5}%
b^{0}\\
&  =1b^{5}+5ab^{4}+10a^{2}b^{3}+10a^{3}b^{2}+5a^{4}b+1a^{5}\\
&  =b^{5}+5ab^{4}+10a^{2}b^{3}+10a^{3}b^{2}+5a^{4}b+a^{5}.
\end{align*}
For a more familiar example, for $n=2$, the formula (\ref{eq.thm.binom.bin.1})
becomes%
\[
\left(  a+b\right)  ^{2}=b^{2}+2ab+a^{2}.
\]

\end{example}

\begin{proof}
[Proof of Theorem \ref{thm.binom.bin}.]Clearly, the formula
(\ref{eq.thm.binom.bin.3}) is just the formula (\ref{eq.thm.binom.bin.1}) with
the variables $a$ and $b$ swapped (since $b+a=a+b$). Thus, it will suffice to
prove (\ref{eq.thm.binom.bin.1}).

We will prove (\ref{eq.thm.binom.bin.1}) by induction on $n$:

\textit{Base case:} For $n=0$, this formula (\ref{eq.thm.binom.bin.1}) is
true, since%
\[
\left(  a+b\right)  ^{0}=1\ \ \ \ \ \ \ \ \ \ \text{and}%
\ \ \ \ \ \ \ \ \ \ \sum_{k=0}^{0}\dbinom{0}{k}a^{k}b^{0-k}%
=\underbrace{\dbinom{0}{0}}_{=1}\underbrace{a^{0}}_{=1}\underbrace{b^{0-0}%
}_{=b^{0}=1}=1.
\]

\textit{Induction step:} Let $n\in\mathbb{N}$. We assume (as the induction
hypothesis) that the formula (\ref{eq.thm.binom.bin.1}) holds for $n$. In
other words, we assume that
\begin{equation}
\left(  a+b\right)  ^{n}=\sum_{k=0}^{n}\dbinom{n}{k}a^{k}b^{n-k}.
\label{pf.thm.binom.bin.IH}%
\end{equation}

We must show that the formula (\ref{eq.thm.binom.bin.1}) also holds for $n+1$.
In other words, we must prove that%
\begin{equation}
\left(  a+b\right)  ^{n+1}=\sum_{k=0}^{n+1}\dbinom{n+1}{k}a^{k}b^{n+1-k}.
\label{pf.thm.binom.bin.IG}%
\end{equation}

Indeed, we have%
\begin{align}
&  \left(  a+b\right)  ^{n+1}\nonumber\\
&  =\left(  a+b\right)  ^{n}\cdot\left(  a+b\right) \nonumber\\
&  =\left(  \sum_{k=0}^{n}\dbinom{n}{k}a^{k}b^{n-k}\right)  \cdot\left(
a+b\right)  \ \ \ \ \ \ \ \ \ \ \left(  \text{by (\ref{pf.thm.binom.bin.IH}%
)}\right) \nonumber\\
&  =\left(  \sum_{k=0}^{n}\dbinom{n}{k}a^{k}b^{n-k}\right)  \cdot a+\left(
\sum_{k=0}^{n}\dbinom{n}{k}a^{k}b^{n-k}\right)  \cdot b\nonumber\\
&  =\sum_{k=0}^{n}\dbinom{n}{k}\underbrace{a^{k}b^{n-k}a}_{=a^{k+1}b^{n-k}%
}+\sum_{k=0}^{n}\dbinom{n}{k}a^{k}\underbrace{b^{n-k}b}_{=b^{n-k+1}%
}\nonumber\\
&  \ \ \ \ \ \ \ \ \ \ \ \ \ \ \ \ \ \ \ \ \left(  \text{by distributivity for
finite sums, i.e., by the rule }\left(  \sum_{s=u}^{v}a_{s}\right)
c=\sum_{s=u}^{v}a_{s}c\right) \nonumber\\
&  =\sum_{k=0}^{n}\dbinom{n}{k}a^{k+1}b^{n-k}+\sum_{k=0}^{n}\dbinom{n}{k}%
a^{k}b^{n-k+1}. \label{pf.thm.binom.bin.1}%
\end{align}

On the other hand, for each $k$, we have%
\[
\dbinom{n+1}{k}=\dbinom{n}{k-1}+\dbinom{n}{k}%
\]
(indeed, this is just Theorem \ref{thm.binom.rec}, applied to $n+1$ instead of
$n$). Hence,%
\begin{align}
&  \sum_{k=0}^{n+1}\dbinom{n+1}{k}a^{k}b^{n+1-k}\nonumber\\
&  =\sum_{k=0}^{n+1}\left(  \dbinom{n}{k-1}+\dbinom{n}{k}\right)
a^{k}b^{n+1-k}\nonumber\\
&  =\sum_{k=0}^{n+1}\left(  \dbinom{n}{k-1}a^{k}b^{n+1-k}+\dbinom{n}{k}%
a^{k}b^{n+1-k}\right) \nonumber\\
&  =\sum_{k=0}^{n+1}\dbinom{n}{k-1}a^{k}b^{n+1-k}+\sum_{k=0}^{n+1}\dbinom
{n}{k}a^{k}b^{n+1-k}\nonumber\\
&  =\left(  \underbrace{\dbinom{n}{0-1}}_{\substack{=0\\\text{(by
definition,}\\\text{since }0-1\notin\mathbb{N}\text{)}}}a^{0}b^{n+1-0}%
+\sum_{k=1}^{n+1}\dbinom{n}{k-1}a^{k}b^{n+1-k}\right) \nonumber\\
&  \ \ \ \ \ \ \ \ \ \ +\left(  \sum_{k=0}^{n}\dbinom{n}{k}a^{k}%
b^{n+1-k}+\underbrace{\dbinom{n}{n+1}}_{\substack{=0\\\text{(by Proposition
\ref{prop.binom.0},}\\\text{since }n+1>n\text{)}}}a^{n+1}b^{n+1-\left(
n+1\right)  }\right) \nonumber\\
&  \ \ \ \ \ \ \ \ \ \ \ \ \ \ \ \ \ \ \ \ \left(
\begin{array}
[c]{c}%
\text{here, we have split off the }k=0\text{ addend from the first sum,}\\
\text{and the }k=n+1\text{ addend from the second sum}%
\end{array}
\right) \nonumber\\
&  =\sum_{k=1}^{n+1}\dbinom{n}{k-1}a^{k}b^{n+1-k}+\sum_{k=0}^{n}\dbinom{n}%
{k}a^{k}b^{n+1-k}. \label{pf.thm.binom.bin.2}%
\end{align}

Let us now compare the two equalities (\ref{pf.thm.binom.bin.1}) and
(\ref{pf.thm.binom.bin.2}). Our goal is to prove that their left hand sides
are equal (because this equality will be precisely (\ref{pf.thm.binom.bin.IG}%
)). Let us look at the right hand sides instead. The right hand side of
(\ref{pf.thm.binom.bin.1}) consists of two finite sums, and so does the right
hand side of (\ref{pf.thm.binom.bin.2}). The second sums of both right hand
sides are equal, since $n-k+1=n+1-k$ for each $k$. If we can also show that
the respective first sums are equal, then we will conclude that the right hand
sides of (\ref{pf.thm.binom.bin.1}) and (\ref{pf.thm.binom.bin.2}) are equal,
and therefore the left hand sides are also equal, and thus we will conclude
that%
\[
\left(  a+b\right)  ^{n+1}=\sum_{k=0}^{n+1}\dbinom{n+1}{k}a^{k}b^{n+1-k},
\]
which is precisely our goal.

So it remains to prove that the first sums on the right hand sides of
(\ref{pf.thm.binom.bin.1}) and (\ref{pf.thm.binom.bin.2}) are equal. In other
words, it remains to prove that%
\begin{equation}
\sum_{k=0}^{n}\dbinom{n}{k}a^{k+1}b^{n-k}=\sum_{k=1}^{n+1}\dbinom{n}{k-1}%
a^{k}b^{n+1-k}. \label{pf.thm.binom.bin.sum=sum}%
\end{equation}
But this becomes clear if we observe that these two sums contain the exact
same addends: Indeed, written out without using summation signs, both sums
become%
\[
\dbinom{n}{0}a^{1}b^{n}+\dbinom{n}{1}a^{2}b^{n-1}+\dbinom{n}{2}a^{3}%
b^{n-2}+\cdots+\dbinom{n}{n}a^{n+1}b^{0}.
\]

This argument can be made more rigorously using an important summation rule,
known as \textbf{substitution}. In its simplest form, this rule says that%
\begin{equation}
\sum_{k=u}^{v}c_{k}=\sum_{k=u+\delta}^{v+\delta}c_{k-\delta}
\label{pf.thm.binom.bin.subrule}%
\end{equation}
for any integers $u,v,\delta$ and any numbers $c_{u},c_{u+1},\ldots,c_{v}$.
This is the discrete analogue of the formula%
\[
\int_{u}^{v}f\left(  x\right)  dx=\int_{u+\delta}^{v+\delta}f\left(
x-\delta\right)  dx
\]
from real analysis. A formal proof of (\ref{pf.thm.binom.bin.subrule}) can
easily be given by induction on $v$, but intuitively
(\ref{pf.thm.binom.bin.subrule}) should be obvious (since both sides are
$c_{u}+c_{u+1}+\cdots+c_{v}$).

When we use (\ref{pf.thm.binom.bin.subrule}) to rewrite a sum of the form
$\sum_{k=u}^{v}c_{k}$ as $\sum_{k=u+\delta}^{v+\delta}c_{k-\delta}$, we say
that we are \textbf{substituting} $k-\delta$ for $k$ in the sum. For example,
taking $u=4$ and $v=9$ and $c_{k}=k^{k}$ and $\delta=-2$, we see that
\[
\sum_{k=4}^{9}k^{k}=\sum_{k=4+\left(  -2\right)  }^{9+\left(  -2\right)
}\left(  k-\left(  -2\right)  \right)  ^{k-\left(  -2\right)  }=\sum_{k=2}%
^{7}\left(  k+2\right)  ^{k+2}.
\]

Now, substituting $k-1$ for $k$ in the sum $\sum_{k=0}^{n}\dbinom{n}{k}%
a^{k+1}b^{n-k}$, we obtain%
\[
\sum_{k=0}^{n}\dbinom{n}{k}a^{k+1}b^{n-k}=\sum_{k=1}^{n+1}\dbinom{n}%
{k-1}\underbrace{a^{\left(  k-1\right)  +1}}_{=a^{k}}%
\ \ \underbrace{b^{n-\left(  k-1\right)  }}_{=b^{n+1-k}}=\sum_{k=1}%
^{n+1}\dbinom{n}{k-1}a^{k}b^{n+1-k}.
\]
This proves (\ref{pf.thm.binom.bin.sum=sum}) rigorously.

Having proved (\ref{pf.thm.binom.bin.sum=sum}), we have shown that the first
sums on the right hand sides of (\ref{pf.thm.binom.bin.1}) and
(\ref{pf.thm.binom.bin.2}) are equal. As we explained, this yields
(\ref{pf.thm.binom.bin.IG}), and thus completes the induction step. This
proves (\ref{eq.thm.binom.bin.1}), thus concluding the proof of Theorem
\ref{thm.binom.bin}.
\end{proof}

\begin{exercise}
\label{exe.binom.2^n}Let $n\in\mathbb{N}$. \medskip

\textbf{(a)} Prove that
\[
\sum_{k=0}^{n}\dbinom{n}{k}=2^{n}\ \ \ \ \ \ \ \ \ \ \text{and}%
\ \ \ \ \ \ \ \ \ \ \sum_{k=0}^{n}\left(  -1\right)  ^{k}\dbinom{n}{k}=%
\begin{cases}
0, & \text{if }n>0;\\
1, & \text{if }n=0.
\end{cases}
\]
(For example, for $n=4$, this is saying that $\dbinom{4}{0}+\dbinom{4}%
{1}+\dbinom{4}{2}+\dbinom{4}{3}+\dbinom{4}{4}=2^{4}$ and that $\dbinom{4}%
{0}-\dbinom{4}{1}+\dbinom{4}{2}-\dbinom{4}{3}+\dbinom{4}{4}=0$.) \medskip

\textbf{(b)} Assume that $n$ is positive. Show that%
\[
\sum_{\substack{k\in\left\{  0,1,\ldots,n\right\}  \\\text{is even}}%
}\dbinom{n}{k}=2^{n-1}.
\]
(The left hand side can be explicitly written as $\dbinom{n}{0}+\dbinom{n}%
{2}+\dbinom{n}{4}+\cdots$, ending at $\dbinom{n}{n-1}$ or $\dbinom{n}{n}$
depending on whether $n$ is odd or even.)

[You can use Proposition \ref{prop.ent.odd=2k+1} here.]
\end{exercise}

\begin{exercise}
\label{exe.23hw2.6}Recall the Fibonacci sequence (Definition \ref{def.ind.fib}%
). Prove that every $n\in\mathbb{N}$ and $m\in\mathbb{N}$ satisfy%
\[
\sum_{k=0}^{n}\dbinom{n}{k}f_{m+k}=f_{m+2n}.
\]

\end{exercise}

\begin{exercise}
\label{exe.23hw3.5}\textbf{(a)} Prove that $k\dbinom{n}{k}=n\dbinom{n-1}{k-1}$
for any two numbers $n$ and $k$. \medskip

\textbf{(b)} Prove that $\sum_{k=0}^{n}k\dbinom{n}{k}x^{k}=nx\left(
x+1\right)  ^{n-1}$ for any positive integer $n$ and any number $x$. \medskip

[\textbf{Hint:} In part \textbf{(a)}, don't forget about cases like $k=0$ and
$k\notin\mathbb{N}$.]
\end{exercise}

\begin{exercise}
\label{exe.23hw3.6}Let $\left(  f_{0},f_{1},f_{2},\ldots\right)  $ be the
Fibonacci sequence. Let $n\in\mathbb{N}$. Prove that
\[
2^{n-1}\cdot f_{n}=\sum_{k=0}^{n}\dbinom{n}{2k+1}\cdot5^{k}.
\]

[\textbf{Hint:} The $5$ on the right hand side looks suspiciously like the $5$
in $\dfrac{1+\sqrt{5}}{2}$, whereas the binomial coefficients look like the
binomial formula...]
\end{exercise}

\subsection{More properties of binomial coefficients}

\begin{exercise}
\label{exe.21finhw3.2}Prove that every $n\in\mathbb{N}$ and every number $a$
satisfy%
\[
\sum_{k=0}^{n}\dbinom{k+a}{k}=\dbinom{n+a+1}{n}.
\]

(For example, for $n=4$ and $a=2$, this is saying that
\[
\dbinom{2}{0}+\dbinom{3}{1}+\dbinom{4}{2}+\dbinom{5}{3}+\dbinom{6}{4}%
=\dbinom{7}{4}.
\]
Keep in mind that $a$ doesn't have to be an integer in general!)
\end{exercise}

\begin{exercise}
Let $n\in\mathbb{N}$. \medskip

\textbf{(a)} Prove that $1\cdot3\cdot5\cdot\cdots\cdot\left(  2n-1\right)  $
(that is, the product of the first $n$ odd positive integers) is
$\dfrac{\left(  2n\right)  !}{2^{n}\cdot n!}$. \medskip

\textbf{(b)} Prove that $\dbinom{-1/2}{n}=\left(  \dfrac{-1}{4}\right)
^{n}\dbinom{2n}{n}$.
\end{exercise}

\begin{exercise}
Let $m,n\in\mathbb{N}$ be such that $m>0$. \medskip

\textbf{(a)} Prove that $\dbinom{mn}{m}=n\dbinom{mn-1}{m-1}$. \medskip

\textbf{(b)} Prove that $\dfrac{\left(  mn\right)  !}{m!^{n}\cdot n!}$ is an
integer. \medskip

[\textbf{Hint:} For part \textbf{(b)}, induct on $n$.]
\end{exercise}

The equalities in the above exercises (as well as Theorem \ref{thm.binom.rec},
Theorem \ref{thm.binom.bin}, Theorem \ref{thm.binom.fib}, Exercise
\ref{exe.23mt1.4} and similar results) are known as \emph{binomial identities}
-- which simply means identities that involve binomial coefficients. There are
many more. A reader interested in going deeper may want to consult
\cite[Chapter 5]{GKP}, \cite[Chapter 2]{19fco}, \cite[\S 7.6, \S 7.7 and other
places]{20f} and \cite{WilfA=B}.

\newpage

\section{\label{chp.ent}Elementary number theory}

Number theory is commonly understood to be the study of integers, and
particularly of those properties and features of integers that do not make
much sense for rational, real or complex numbers. Divisibility is one such
property; prime numbers are another. In this course, we will only cover the
very basics of elementary number theory; there is no shortage of texts that go
much deeper (some freely available ones are \cite{Mileti22}, \cite{Hackma09},
\cite{Stein08}, \cite{Shoup08} and \cite{Martin17}).

\subsection{\label{sec.ent.div}Divisibility}

\subsubsection{Definition}

We begin by defining the one most important concept in number theory:

\begin{definition}
Let $a$ and $b$ be two integers.

We write $a\mid b$ (and we say that \textquotedblleft$a$ \textbf{divides }%
$b$\textquotedblright, or \textquotedblleft$b$ is \textbf{divisible} by
$a$\textquotedblright, or \textquotedblleft$b$ is a \textbf{multiple} of
$a$\textquotedblright, or \textquotedblleft$a$ is a \textbf{divisor} of
$b$\textquotedblright; yes, all these statements are equivalent) if there
exists an integer $c$ such that $b=ac$.

We write $a\nmid b$ if we don't have $a\mid b$.
\end{definition}

\begin{example}
\textbf{(a)} We have $4\mid12$, because $12=4\cdot3$. \medskip

\textbf{(b)} We have $4\nmid11$, because there exists no integer $c$ such that
$11=4c$. \medskip

\textbf{(c)} We have $1\mid b$ for every integer $b$, since $b=1\cdot b$.
\medskip

\textbf{(d)} We have $a\mid a$ for every integer $a$, since $a=a\cdot1$. In
particular, $0\mid0$, which is somewhat controversial (but true in our
opinion). (Some authors deliberately exclude $0$ as a divisor on the grounds
that $\dfrac{0}{0}$ is not well-defined, but I believe that making this an
exception is more trouble than it is worth.) \medskip

\textbf{(e)} We have $a\mid0$ for every integer $a$, since $0=a\cdot0$.
\medskip

\textbf{(f)} An integer $b$ satisfies $0\mid b$ if and only if $b=0$.
\end{example}

The well-known concepts of even and odd integers are instances of divisibility:

\begin{definition}
\textbf{(a)} An integer $n$ is said to be \textbf{even} if $2\mid n$. \medskip

\textbf{(b)} An integer $n$ is said to be \textbf{odd} if $2\nmid n$.
\end{definition}

You probably know a few things about even and odd numbers already: e.g.,

\begin{enumerate}
\item The sum of two even numbers is even.

\item The sum of an even with an odd number is odd.

\item The sum of two odd numbers is even.
\end{enumerate}

Strictly speaking, these claims (particularly the third one) are not at all
obvious. So we need to understand divisibility better to even convince
ourselves that such fundamental statements are true. We will do this soon
(Corollary \ref{cor.ent.odd+odd}). First, let us prove some basic facts about divisibility.

\subsubsection{Basic properties}

In the next proposition, we shall let $\operatorname*{abs}x$ denote the
absolute value of a real number $x$. Thus,%
\[
\operatorname*{abs}x=%
\begin{cases}
x, & \text{if }x\geq0;\\
-x, & \text{if }x<0.
\end{cases}
\]
This absolute value $\operatorname*{abs}x$ is normally called $\left\vert
x\right\vert $, but I believe that writing \textquotedblleft%
$\operatorname*{abs}a\mid\operatorname*{abs}b$\textquotedblright\ is less
confusing than writing \textquotedblleft$\left\vert a\right\vert
\mid\left\vert b\right\vert $\textquotedblright\ (where four of the bars stand
for absolute values, while the middle bar stands for divisibility).

\begin{proposition}
\label{prop.div.abs}Let $a$ and $b$ be two integers. Then: \medskip

\textbf{(a)} We have $a\mid b$ if and only if $\operatorname*{abs}%
a\mid\operatorname*{abs}b$. \medskip

\textbf{(b)} If $a\mid b$ and $b\neq0$, then $\operatorname*{abs}%
a\leq\operatorname*{abs}b$. \medskip

\textbf{(c)} If $a\mid b$ and $b\mid a$, then $\operatorname*{abs}%
a=\operatorname*{abs}b$. \medskip

\textbf{(d)} Assume that $a\neq0$. Then, $a\mid b$ if and only if $\dfrac
{b}{a}\in\mathbb{Z}$.
\end{proposition}

\begin{proof}
\textbf{(a)} Proposition \ref{prop.div.abs} \textbf{(a)} says that the
divisibility $a\mid b$ does not depend on the signs of $a$ and $b$; in other
words, it says that we can replace the numbers $a$ and $b$ by their absolute
values without changing the truth (or falsity) of $a\mid b$.

Clearly, in order to prove this, it suffices to show the following two statements:

\begin{enumerate}
\item \textit{Statement 1:} We can replace $a$ by $-a$ without changing the
truth (or falsity) of $a\mid b$;

\item \textit{Statement 2:} We can replace $b$ by $-b$ without changing the
truth (or falsity) of $a\mid b$.
\end{enumerate}

In fact, if these two statements are proved, the claim of Proposition
\ref{prop.div.abs} \textbf{(a)} will follow, because the absolute value of an
integer $c$ is always either $c$ or $-c$.

But both of these statements are easy to prove:

To prove Statement 1, we assume that $a\mid b$. Thus, $b=ac$ for some integer
$c$ (by the definition of \textquotedblleft$a\mid b$\textquotedblright).
Hence, for this integer $c$, we have $b=ac=\left(  -a\right)  \left(
-c\right)  $, which allows us to conclude that $-a\mid b$ (since $-c$ is an
integer, too). Thus, we have shown that $a\mid b$ implies $-a\mid b$.
Conversely, a similar argument shows that $-a\mid b$ implies $a\mid b$
(indeed, it is the same argument with the roles of $a$ and $-a$ swapped,
because $-\left(  -a\right)  =a$). Thus, the statements $a\mid b$ and $-a\mid
b$ are equivalent. In other words, we can replace $a$ by $-a$ without changing
the truth (or falsity) of $a\mid b$. This proves Statement 1.

The proof of Statement 2 is similar. (This time, you need to argue that $a\mid
b$ implies $a\mid-b$. Again, write $b$ as $b=ac$, and conclude that
$-b=-ac=a\left(  -c\right)  $, so that $a\mid-b$.)

Thus, both Statements 1 and 2 are proved, so that the proof of Proposition
\ref{prop.div.abs} \textbf{(a)} is complete. \medskip

\textbf{(b)} Assume that $a\mid b$ and $b\neq0$. We must show that
$\operatorname*{abs}a\leq\operatorname*{abs}b$.

Let $x=\operatorname*{abs}a$ and $y=\operatorname*{abs}b$. Thus, $x$ is a
nonnegative integer and $y$ is a positive integer (since $b\neq0$). Thus,
$x\geq0$ and $y>0$.

Proposition \ref{prop.div.abs} \textbf{(a)} yields that $\operatorname*{abs}%
a\mid\operatorname*{abs}b$ (since $a\mid b$). In other words, $x\mid y$ (since
$x=\operatorname*{abs}a$ and $y=\operatorname*{abs}b$). In other words, $y=xz$
for some integer $z$. Consider this $z$.

If we had $z\leq0$, then we would have $y=\underbrace{x}_{\geq0}%
\ \ \underbrace{z}_{\leq0}\leq0$ (by the standard rules for inequalities),
which would contradict $y>0$. Hence, we cannot have $z\leq0$. Thus, $z>0$, so
that $z\geq1$ (since $z$ is an integer). Hence, $xz\geq x1$ (since $x\geq0$
allows us to multiply any inequality by $x$ without having to flip the sign).
Therefore, $y=xz\geq x1=x$. In other words, $x\leq y$. In other words,
$\operatorname*{abs}a\leq\operatorname*{abs}b$ (since $x=\operatorname*{abs}a$
and $y=\operatorname*{abs}b$). This proves Proposition \ref{prop.div.abs}
\textbf{(b)}. \medskip

\textbf{(c)} Let $a\mid b$ and $b\mid a$. We must prove that
$\operatorname*{abs}a=\operatorname*{abs}b$.

If $a=0$, then this is easily done (because if $a=0$, then $0=a\mid b$ quickly
leads to $b=0$, and therefore $a=0=b$, so that $\operatorname*{abs}%
a=\operatorname*{abs}b$).

Likewise, this is easily done if $b=0$.

It remains to handle the third possible case, which is when both $a$ and $b$
are $\neq0$. Consider this case. In this case, Proposition \ref{prop.div.abs}
\textbf{(b)} yields $\operatorname*{abs}a\leq\operatorname*{abs}b$ (since
$a\mid b$ and $b\neq0$). However, we can also apply Proposition
\ref{prop.div.abs} \textbf{(b)} with the roles of $a$ and $b$ interchanged
(since $b\mid a$ and $a\neq0$), and thus obtain $\operatorname*{abs}%
b\leq\operatorname*{abs}a$. Combining this with $\operatorname*{abs}%
a\leq\operatorname*{abs}b$, we find $\operatorname*{abs}a=\operatorname*{abs}%
b$. Proposition \ref{prop.div.abs} \textbf{(c)} is thus proved. \medskip

\textbf{(d)} This is quite straightforward:

Assume that $a\mid b$. Thus, there exists some integer $c$ such that $b=ac$
(by the definition of \textquotedblleft$a\mid b$\textquotedblright). This $c$
must then be $\dfrac{b}{a}$ (since $b=ac$ implies $c=\dfrac{b}{a}$ in view of
$a\neq0$). Hence, $\dfrac{b}{a}$ is an integer, i.e., we have $\dfrac{b}{a}%
\in\mathbb{Z}$.

Forget that we assumed $a\mid b$. We thus have shown that $\dfrac{b}{a}%
\in\mathbb{Z}$ if $a\mid b$. The same argument (done in reverse) yields that
conversely, if $\dfrac{b}{a}\in\mathbb{Z}$, then $a\mid b$. Combining these
two facts, we conclude that $a\mid b$ if and only if $\dfrac{b}{a}%
\in\mathbb{Z}$. This proves Proposition \ref{prop.div.abs} \textbf{(d)}.
\end{proof}

This was a warm-up (if somewhat laborious to write up). Here are some slightly
more substantial properties of divisibility:

\begin{theorem}
[rules for divisibility]\label{thm.div.rules}\textbf{(a)} We have $a\mid a$
for each $a\in\mathbb{Z}$. (This is called \textbf{reflexivity of
divisibility}.) \medskip

\textbf{(b)} If $a,b,c\in\mathbb{Z}$ satisfy $a\mid b$ and $b\mid c$, then
$a\mid c$. (This is called \textbf{transitivity of divisibility}.) \medskip

\textbf{(c)} If $a_{1},a_{2},b_{1},b_{2}\in\mathbb{Z}$ satisfy $a_{1}\mid
b_{1}$ and $a_{2}\mid b_{2}$, then $a_{1}a_{2}\mid b_{1}b_{2}$. (This is
called \textbf{multiplying two divisibilities}.) \medskip

\textbf{(d)} If $d,a,b\in\mathbb{Z}$ satisfy $d\mid a$ and $d\mid b$, then
$d\mid a+b$. (This is often restated as \textquotedblleft a sum of two
multiples of $d$ is again a multiple of $d$\textquotedblright.)
\end{theorem}

\begin{proof}
\textbf{(a)} Let $a\in\mathbb{Z}$. Then, $a=a\cdot1$, so that $a\mid a$ (since
$1$ is an integer). This proves Theorem \ref{thm.div.rules} \textbf{(a)}.
\medskip

\textbf{(b)} Let $a,b,c\in\mathbb{Z}$ satisfy $a\mid b$ and $b\mid c$.

From $a\mid b$, we see that there exists an integer $x$ such that $b=ax$.

From $b\mid c$, we see that there exists an integer $y$ such that $c=by$.

Consider these integers $x$ and $y$. Now,%
\[
c=\underbrace{b}_{=ax}y=axy.
\]
Hence, there exists some integer $z$ such that $c=az$ (namely, $z=xy$). This
shows that $a\mid c$. Theorem \ref{thm.div.rules} \textbf{(b)} is thus proven.
\medskip

\textbf{(c)} Let $a_{1},a_{2},b_{1},b_{2}\in\mathbb{Z}$ satisfy $a_{1}\mid
b_{1}$ and $a_{2}\mid b_{2}$.

From $a_{1}\mid b_{1}$, we see that $b_{1}=a_{1}c_{1}$ for some integer
$c_{1}$.

From $a_{2}\mid b_{2}$, we see that $b_{2}=a_{2}c_{2}$ for some integer
$c_{2}$.

Consider these integers $c_{1}$ and $c_{2}$. Now,%
\[
\underbrace{b_{1}}_{=a_{1}c_{1}}\ \ \underbrace{b_{2}}_{=a_{2}c_{2}}%
=a_{1}c_{1}a_{2}c_{2}=\left(  a_{1}a_{2}\right)  \underbrace{\left(
c_{1}c_{2}\right)  }_{\text{an integer}}.
\]
Thus, $a_{1}a_{2}\mid b_{1}b_{2}$. This proves Theorem \ref{thm.div.rules}
\textbf{(c)}. \medskip

\textbf{(d)} Let $d,a,b\in\mathbb{Z}$ satisfy $d\mid a$ and $d\mid b$.

From $d\mid a$, we see that $a=dx$ for some integer $x$.

From $d\mid b$, we see that $b=dy$ for some integer $y$.

Consider these integers $x$ and $y$. Now,%
\[
a+b=dx+dy=d\underbrace{\left(  x+y\right)  }_{\text{an integer}}.
\]
Thus, $d\mid a+b$. This proves Theorem \ref{thm.div.rules} \textbf{(d)}.
\end{proof}

\bigskip

Theorem \ref{thm.div.rules} \textbf{(b)} tells us that divisibilities can be
chained together: If $a\mid b$ and $b\mid c$, then $a\mid c$. Therefore, you
will often see a statement of the form \textquotedblleft$a\mid b$ and $b\mid
c$\textquotedblright\ rewritten as \textquotedblleft$a\mid b\mid
c$\textquotedblright, just like two inequalities $a\leq b$ and $b\leq c$ can
be chained together to form $a\leq b\leq c$. More generally, the statement%
\[
\text{\textquotedblleft}a_{1}\mid a_{2}\mid\cdots\mid a_{k}%
\text{\textquotedblright}%
\]
shall mean that each of the numbers $a_{1},a_{2},\ldots,a_{k}$ divides the
next (i.e., that $a_{1}\mid a_{2}$ and $a_{2}\mid a_{3}$ and so on, ending
with $a_{k-1}\mid a_{k}$). By induction on $k$, it is easy to see that such a
chain of divisibilities always entails $a_{1}\mid a_{k}$. For example,
$3\mid6\mid18\mid36$ entails $3\mid36$.

\begin{exercise}
\label{exe.23hw3.1a}Let $a,b\in\mathbb{Z}$ satisfy $a\mid b$. Prove that
$a^{k}\mid b^{k}$ for each $k\in\mathbb{N}$.
\end{exercise}

\subsubsection{Divisibility criteria}

How can you spot divisibilities between actual numbers? For small values of
$a$, there are several known \textbf{divisibility criteria} (also known as
\textbf{divisibility rules}), which give simple methods to check whether a
given integer $b$ is divisible by $a$ (without computing $\dfrac{b}{a}$). Here
are some:

\begin{theorem}
\label{thm.div.criteria}Let $b\in\mathbb{N}$. Write $b$ in decimal notation.
Then: \medskip

\textbf{(a)} We have $2\mid b$ if and only if the last digit of $b$ is $0$ or
$2$ or $4$ or $6$ or $8$. \medskip

\textbf{(b)} We have $5\mid b$ if and only if the last digit of $b$ is $0$ or
$5$. \medskip

\textbf{(c)} We have $10\mid b$ if and only if the last digit of $b$ is $0$.
\medskip

\textbf{(d)} We have $3\mid b$ if and only if the sum of the digits of $b$ is
divisible by $3$. \medskip

\textbf{(e)} We have $9\mid b$ if and only if the sum of the digits of $b$ is
divisible by $9$.
\end{theorem}

\begin{example}
Let $b=10835$. Then, $2\nmid b$, since the last digit of $b$ is neither $0$
nor $2$ nor $4$ nor $6$ nor $8$ (but $5$). However, $5\mid b$, since the last
digit of $b$ is $0$ or $5$. Do we have $3\mid b$ ? The sum of the digits of
$b$ is $1+0+8+3+5=17$, which is not divisible by $3$. Thus, $b$ is not
divisible by $3$. Hence, $b$ is not divisible by $9$ either, because if we had
$9\mid b$, then we would get $3\mid9\mid b$ (by Theorem \ref{thm.div.rules}
\textbf{(b)}), which would contradict the previous sentence.
\end{example}

How do we prove Theorem \ref{thm.div.criteria}?

The easiest part is part \textbf{(c)}: If you multiply a number (written in
decimal) by $10$, then its decimal representation just grows a new digit $0$
at the end. Thus, if $10\mid b$, then the last digit of $b$ is $0$.
Conversely, if the last digit of $b$ is $0$, then $b=10b^{\prime}$, where
$b^{\prime}$ is the number $b$ with its last digit removed. For example,
$390=10\cdot39$.

Parts \textbf{(a)} and \textbf{(b)} of Theorem \ref{thm.div.criteria} are
somewhat trickier, and parts \textbf{(d)} and \textbf{(e)} are even more so.
To get simple proofs for these parts, we will now introduce another type of
relation between integers, known as \textbf{congruence modulo }$n$.

\subsection{\label{sec.ent.cong}Congruence modulo $n$}

\subsubsection{Definition}

\begin{definition}
Let $n,a,b\in\mathbb{Z}$. We say that $a$ is \textbf{congruent to }$b$
\textbf{modulo }$n$ if and only if $n\mid a-b$.

We shall use the notation \textquotedblleft$a\equiv b\operatorname{mod}%
n$\textquotedblright\ for \textquotedblleft$a$ is congruent to $b$ modulo
$n$\textquotedblright.

We shall use the notation \textquotedblleft$a\not \equiv b\operatorname{mod}%
n$\textquotedblright\ for \textquotedblleft$a$ is not congruent to $b$ modulo
$n$\textquotedblright.
\end{definition}

\begin{example}
\textbf{(a)} Is $3\equiv7\operatorname{mod}2$ ? This would mean that
$2\mid3-7$, which is true (since $3-7=-4=2\cdot\left(  -2\right)  $). So yes,
we do have $3\equiv7\operatorname{mod}2$. \medskip

\textbf{(b)} Is $3\equiv6\operatorname{mod}2$ ? This would mean that
$2\mid3-6$, which is false (since $3-6=-3$ is not divisible by $2$). So we
have $3\not \equiv 6\operatorname{mod}2$. \medskip

\textbf{(c)} We have $a\equiv b\operatorname{mod}1$ for any integers $a$ and
$b$. This is because $1\mid a-b$ (since $1$ divides every integer). \medskip

\textbf{(d)} Two integers $a$ and $b$ satisfy $a\equiv b\operatorname{mod}0$
if and only if $a=b$ (since $0$ divides only $0$ itself). \medskip

\textbf{(e)} For any two integers $a$ and $b$, we have $a+b\equiv
a-b\operatorname{mod}2$, since $\left(  a+b\right)  -\left(  a-b\right)  =2b$
is clearly divisible by $2$.
\end{example}

The word \textquotedblleft modulo\textquotedblright\ in the phrase
\textquotedblleft$a$ is congruent to $b$ modulo $n$\textquotedblright\ was
coined by Gauss, and should be read as something like \textquotedblleft with
respect to\textquotedblright. You can translate the statement
\textquotedblleft$a$ is congruent to $b$ modulo $n$\textquotedblright\ as
\textquotedblleft$a$ equals $b$ up to a multiple of $n$\textquotedblright.
Indeed, the definition of congruence can be restated as follows:%
\[
a\equiv b\operatorname{mod}n\ \ \ \ \ \ \ \ \ \ \text{if and only
if}\ \ \ \ \ \ \ \ \ \ a=b+nc\text{ for some }c\in\mathbb{Z}.
\]

As we will soon see, congruence modulo $2$ is essentially parity:

\begin{itemize}
\item Two even numbers are always congruent (to each other) modulo $2$.

\item Two odd numbers are always congruent (to each other) modulo $2$.

\item An even number is never congruent to an odd number modulo $2$.
\end{itemize}

We will prove this in Corollary \ref{cor.ent.parity-mod2}.

\subsubsection{Basic properties}

First, we shall establish some fundamental properties of congruence.

\begin{proposition}
\label{prop.mod.0}Let $n,a\in\mathbb{Z}$. Then, $a\equiv0\operatorname{mod}n$
if and only if $n\mid a$.
\end{proposition}

\begin{proof}
By the definition of congruence, we have the following equivalences:%
\[
\left(  a\equiv0\operatorname{mod}n\right)  \ \Longleftrightarrow\ \left(
n\mid a-0\right)  \ \Longleftrightarrow\ \left(  n\mid a\right)  .
\]
Proposition \ref{prop.mod.0} thus follows.
\end{proof}

\begin{proposition}
\label{prop.mod.rules}Let $n\in\mathbb{Z}$. Then: \medskip

\textbf{(a)} We have $a\equiv a\operatorname{mod}n$ for every $a\in\mathbb{Z}%
$. (This is called the \textbf{reflexivity of congruence}.) \medskip

\textbf{(b)} If $a,b\in\mathbb{Z}$ satisfy $a\equiv b\operatorname{mod}n$,
then $b\equiv a\operatorname{mod}n$. (This is called the \textbf{symmetry of
congruence}.) \medskip

\textbf{(c)} If $a,b,c\in\mathbb{Z}$ satisfy $a\equiv b\operatorname{mod}n$
and $b\equiv c\operatorname{mod}n$, then $a\equiv c\operatorname{mod}n$. (This
is called the \textbf{transitivity of congruence}.) \medskip

\textbf{(d)} If $a_{1},a_{2},b_{1},b_{2}\in\mathbb{Z}$ satisfy%
\[
a_{1}\equiv b_{1}\operatorname{mod}n\ \ \ \ \ \ \ \ \ \ \text{and}%
\ \ \ \ \ \ \ \ \ \ a_{2}\equiv b_{2}\operatorname{mod}n,
\]
then%
\begin{align}
a_{1}+a_{2}  &  \equiv b_{1}+b_{2}\operatorname{mod}%
n;\label{eq.prop.mod.rules.d.+}\\
a_{1}-a_{2}  &  \equiv b_{1}-b_{2}\operatorname{mod}%
n;\label{eq.prop.mod.rules.d.-}\\
a_{1}a_{2}  &  \equiv b_{1}b_{2}\operatorname{mod}n.
\label{eq.prop.mod.rules.d.*}%
\end{align}
(In other words, two congruences modulo $n$ can be added, subtracted or
multiplied.) \medskip

\textbf{(e)} Let $m\in\mathbb{Z}$ be such that $m\mid n$. If $a,b\in
\mathbb{Z}$ satisfy $a\equiv b\operatorname{mod}n$, then $a\equiv
b\operatorname{mod}m$.
\end{proposition}

\begin{proof}
\textbf{(a)} Let $a\in\mathbb{Z}$. Then, $n\mid a-a$ because $a-a=0=n\cdot0$.
But this means that $a\equiv a\operatorname{mod}n$. Thus, Proposition
\ref{prop.mod.rules} \textbf{(a)} follows. \medskip

\textbf{(b)} Let $a,b\in\mathbb{Z}$ be such that $a\equiv b\operatorname{mod}%
n$. Thus, $n\mid a-b$.

We must prove that $b\equiv a\operatorname{mod}n$, i.e., that $n\mid b-a$.

However, $b-a=\left(  a-b\right)  \cdot\left(  -1\right)  $, so that $a-b\mid
b-a$. Hence, $n\mid a-b\mid b-a$. Therefore, by the transitivity of
divisibility, $n\mid b-a$. But this means precisely that $b\equiv
a\operatorname{mod}n$. Thus, Proposition \ref{prop.mod.rules} \textbf{(b)} is
proved. \medskip

\textbf{(c)} Let $a,b,c\in\mathbb{Z}$ be such that $a\equiv
b\operatorname{mod}n$ and $b\equiv c\operatorname{mod}n$.

From $a\equiv b\operatorname{mod}n$, we obtain $n\mid a-b$.

From $b\equiv c\operatorname{mod}n$, we obtain $n\mid b-c$.

Recall that a sum of two multiples of $n$ is again a multiple of $n$ (this is
Theorem \ref{thm.div.rules} \textbf{(d)}). Thus, from $n\mid a-b$ and $n\mid
b-c$, we obtain $n\mid\left(  a-b\right)  +\left(  b-c\right)  $. Since
$\left(  a-b\right)  +\left(  b-c\right)  =a-c$, we can rewrite this as $n\mid
a-c$. In other words, $a\equiv c\operatorname{mod}n$. This proves Proposition
\ref{prop.mod.rules} \textbf{(c)}. \medskip

\textbf{(d)} Let $a_{1},a_{2},b_{1},b_{2}\in\mathbb{Z}$ satisfy%
\[
a_{1}\equiv b_{1}\operatorname{mod}n\ \ \ \ \ \ \ \ \ \ \text{and}%
\ \ \ \ \ \ \ \ \ \ a_{2}\equiv b_{2}\operatorname{mod}n.
\]
Thus, $n\mid a_{1}-b_{1}$ and $n\mid a_{2}-b_{2}$.

From $n\mid a_{1}-b_{1}$, we see that $a_{1}-b_{1}=nc_{1}$ for some integer
$c_{1}$.

From $n\mid a_{2}-b_{2}$, we see that $a_{2}-b_{2}=nc_{2}$ for some integer
$c_{2}$.

Consider these integers $c_{1}$ and $c_{2}$.

From $a_{1}-b_{1}=nc_{1}$, we obtain $a_{1}=b_{1}+nc_{1}$. Similarly,
$a_{2}=b_{2}+nc_{2}$.

Adding the equalities $a_{1}=b_{1}+nc_{1}$ and $a_{2}=b_{2}+nc_{2}$ together,
we find
\[
a_{1}+a_{2}=\left(  b_{1}+nc_{1}\right)  +\left(  b_{2}+nc_{2}\right)
=b_{1}+b_{2}+n\left(  c_{1}+c_{2}\right)  .
\]
Thus, $a_{1}+a_{2}$ differs from $b_{1}+b_{2}$ by a multiple of $n$ (namely,
by $n\left(  c_{1}+c_{2}\right)  $). In other words, $n\mid\left(  a_{1}%
+a_{2}\right)  -\left(  b_{1}+b_{2}\right)  $. Hence,%
\[
a_{1}+a_{2}\equiv b_{1}+b_{2}\operatorname{mod}n.
\]

Subtracting the equalities $a_{1}=b_{1}+nc_{1}$ and $a_{2}=b_{2}+nc_{2}$ from
one another, we obtain%
\[
a_{1}-a_{2}=\left(  b_{1}+nc_{1}\right)  -\left(  b_{2}+nc_{2}\right)
=b_{1}-b_{2}+n\left(  c_{1}-c_{2}\right)  .
\]
Thus, $a_{1}-a_{2}$ differs from $b_{1}-b_{2}$ by a multiple of $n$ (namely,
by $n\left(  c_{1}-c_{2}\right)  $). Hence,%
\[
a_{1}-a_{2}\equiv b_{1}-b_{2}\operatorname{mod}n.
\]

Multiplying the equalities $a_{1}=b_{1}+nc_{1}$ and $a_{2}=b_{2}+nc_{2}$
together, we find%
\begin{align*}
a_{1}a_{2}  &  =\left(  b_{1}+nc_{1}\right)  \left(  b_{2}+nc_{2}\right)
=b_{1}b_{2}+b_{1}nc_{2}+nc_{1}b_{2}+nc_{1}nc_{2}\\
&  =b_{1}b_{2}+n\left(  b_{1}c_{2}+c_{1}b_{2}+nc_{1}c_{2}\right)  .
\end{align*}
Thus, $a_{1}a_{2}$ differs from $b_{1}b_{2}$ by a multiple of $n$ (namely, by
$n\left(  b_{1}c_{2}+c_{1}b_{2}+nc_{1}c_{2}\right)  $). Therefore,%
\[
a_{1}a_{2}\equiv b_{1}b_{2}\operatorname{mod}n.
\]
Altogether, we have proved all claims of Proposition \ref{prop.mod.rules}
\textbf{(d)} now. \medskip

\textbf{(e)} Let $m\in\mathbb{Z}$ be such that $m\mid n$. Let $a,b\in
\mathbb{Z}$ satisfy $a\equiv b\operatorname{mod}n$.

Thus, $n\mid a-b$. Hence, $m\mid n\mid a-b$, so that $m\mid a-b$ (by the
transitivity of divisibility). But this means that $a\equiv
b\operatorname{mod}m$. Thus, Proposition \ref{prop.mod.rules} \textbf{(e)} follows.
\end{proof}

\bigskip

Proposition \ref{prop.mod.rules} \textbf{(b)} says that congruences can be
turned around: From $a\equiv b\operatorname{mod}n$, we can always obtain
$b\equiv a\operatorname{mod}n$. (This is very different from divisibilities,
for which $a\mid b$ almost never implies $b\mid a$.)

Proposition \ref{prop.mod.rules} \textbf{(c)} says that congruences can be
chained together: From $a\equiv b\operatorname{mod}n$ and $b\equiv
c\operatorname{mod}n$, we can always obtain $a\equiv c\operatorname{mod}n$.
This is analogous to Theorem \ref{thm.div.rules} \textbf{(b)}, and leads to a
similar convention: Instead of writing \textquotedblleft$a\equiv
b\operatorname{mod}n$ and $b\equiv c\operatorname{mod}n$\textquotedblright, we
will often just write \textquotedblleft$a\equiv b\equiv c\operatorname{mod}%
n$\textquotedblright, understanding that (by Proposition \ref{prop.mod.rules}
\textbf{(c)}) this chain of congruences automatically implies $a\equiv
c\operatorname{mod}n$. More generally, the statement%
\[
\text{\textquotedblleft}a_{1}\equiv a_{2}\equiv\cdots\equiv a_{k}%
\operatorname{mod}n\text{\textquotedblright}%
\]
shall mean that each of the numbers $a_{1},a_{2},\ldots,a_{k}$ is congruent to
the next modulo $n$ (i.e., that $a_{i}\equiv a_{i+1}\operatorname{mod}n$ for
each $i\in\left\{  1,2,\ldots,k-1\right\}  $). By induction on $k$, it is easy
to see that such a chain of congruences always entails $a_{1}\equiv
a_{k}\operatorname{mod}n$ (and, better yet: $a_{i}\equiv a_{j}%
\operatorname{mod}n$ for all $i$ and $j$).

Note that we can only chain together two congruences modulo the same $n$, not
two congruences modulo two different $n$'s. For example, if we know that
$a\equiv b\operatorname{mod}2$ and $b\equiv c\operatorname{mod}3$, then we
cannot conclude any congruence between $a$ and $c$.

\bigskip

Proposition \ref{prop.mod.rules} \textbf{(d)} says that congruences modulo $n$
(for a fixed integer $n$) can be added, subtracted and multiplied together
(just like equalities). Before you get over-enthusiastic, keep in mind that

\begin{itemize}
\item they cannot be divided by one another: We have $2\equiv
0\operatorname{mod}2$ and $2\equiv2\operatorname{mod}2$ but $2/2\not \equiv
0/2\operatorname{mod}2$.

\item they cannot be taken to each other's power: We have $2\equiv
2\operatorname{mod}2$ and $2\equiv0\operatorname{mod}2$ but $2^{2}%
\not \equiv 2^{0}\operatorname{mod}2$.
\end{itemize}

However, we can take a congruence to a $k$-th power for a fixed $k\in
\mathbb{N}$:

\begin{exercise}
\label{exe.mod.k-th-power}Let $n,a,b\in\mathbb{Z}$ be such that $a\equiv
b\operatorname{mod}n$. Let $k\in\mathbb{N}$. Prove that $a^{k}\equiv
b^{k}\operatorname{mod}n$.
\end{exercise}

\begin{noncompile}
The preceding exercise is Exercise 1 \textbf{(b)} on homework set \#3.
\end{noncompile}

Proposition \ref{prop.mod.rules} \textbf{(e)} shows that the $n$ in a
congruence $a\equiv b\operatorname{mod}n$ can be replaced by any divisor of
$n$. For example, if two integers $a$ and $b$ satisfy $a\equiv
b\operatorname{mod}15$, then $a\equiv b\operatorname{mod}3$, since $3$ is a
divisor of $15$. \medskip

The next exercise shows that we can divide a congruence $a\equiv
b\operatorname{mod}n$ by a nonzero integer $d$ as long as we divide all three
numbers in it ($a$, $b$ and $n$) by $d$ (rather than just $a$ and $b$):

\begin{exercise}
\label{exe.23hw3.2}Let $n,d,a,b\in\mathbb{Z}$, and assume that $d\neq0$ and
$da\equiv db\operatorname{mod}dn$. \medskip

\textbf{(a)} Prove that $a\equiv b\operatorname{mod}n$. \medskip

\textbf{(b)} Show by an example that $a\equiv b\operatorname{mod}dn$ is not
necessarily true (i.e., we cannot simply cancel the $d$ from $da$ and $db$
while leaving the $dn$ unchanged).
\end{exercise}

\subsubsection{Proving the divisibility criteria}

Now, let us prove Theorem \ref{thm.div.criteria} \textbf{(e)}, restating it as follows:

\begin{proposition}
Let $m\in\mathbb{N}$. Let $s$ be the sum of the digits of $m$ written in
decimal. (For instance, if $m=302$, then $s=3+0+2=5$.)

Then, $9\mid m$ if and only if $9\mid s$.
\end{proposition}

\begin{proof}
Let the integer $m$ have decimal representation $m_{d}m_{d-1}\cdots m_{0}$
(where $m_{d}$ is the leading digit). Thus,%
\begin{align*}
m  &  =m_{d}\cdot10^{d}+m_{d-1}\cdot10^{d-1}+\cdots+m_{0}\cdot10^{0}%
\ \ \ \ \ \ \ \ \ \ \text{and}\\
s  &  =m_{d}+m_{d-1}+\cdots+m_{0}.
\end{align*}

However, $10\equiv1\operatorname{mod}9$ (since $10-1=9$ is divisible by $9$).
Hence, by Exercise \ref{exe.mod.k-th-power}, we have $10^{k}\equiv
1^{k}\operatorname{mod}9$ for every $k\in\left\{  0,1,\ldots,d\right\}  $.
Multiplying this congruence with the obvious congruence $m_{k}\equiv
m_{k}\operatorname{mod}9$, we obtain\footnote{The reason why we can multiply
two congruences together is Proposition \ref{prop.mod.rules} \textbf{(d)}
(specifically, (\ref{eq.prop.mod.rules.d.*})).}%
\[
m_{k}\cdot10^{k}\equiv m_{k}\cdot1^{k}\operatorname{mod}%
9\ \ \ \ \ \ \ \ \ \ \text{for every }k\in\left\{  0,1,\ldots,d\right\}  .
\]
In other words,%
\[
m_{k}\cdot10^{k}\equiv m_{k}\operatorname{mod}9\ \ \ \ \ \ \ \ \ \ \text{for
every }k\in\left\{  0,1,\ldots,d\right\}
\]
(since $m_{k}\cdot\underbrace{1^{k}}_{=1}=m_{k}$). In other words, we have%
\begin{align*}
m_{d}\cdot10^{d}  &  \equiv m_{d}\operatorname{mod}9;\\
m_{d-1}\cdot10^{d-1}  &  \equiv m_{d-1}\operatorname{mod}9;\\
m_{d-2}\cdot10^{d-2}  &  \equiv m_{d-2}\operatorname{mod}9;\\
&  \ldots;\\
m_{0}\cdot10^{0}  &  \equiv m_{0}\operatorname{mod}9.
\end{align*}
Adding these $d+1$ many congruences together, we obtain\footnote{The reason
why we can add two congruences together is Proposition \ref{prop.mod.rules}
\textbf{(d)} (specifically, (\ref{eq.prop.mod.rules.d.+})). To be very
pedantic, we have to apply (\ref{eq.prop.mod.rules.d.+}) several times, since
we are adding not two but $d+1$ many congruences together.}%
\[
m_{d}\cdot10^{d}+m_{d-1}\cdot10^{d-1}+\cdots+m_{0}\cdot10^{0}\equiv
m_{d}+m_{d-1}+\cdots+m_{0}\operatorname{mod}9.
\]
In other words,%
\[
m\equiv s\operatorname{mod}9
\]
(since $m=m_{d}\cdot10^{d}+m_{d-1}\cdot10^{d-1}+\cdots+m_{0}\cdot10^{0}$ and
$s=m_{d}+m_{d-1}+\cdots+m_{0}$). Turning this congruence around (i.e.,
applying Proposition \ref{prop.mod.rules} \textbf{(b)}), we obtain $s\equiv
m\operatorname{mod}9$.

Now, if $9\mid m$, then $m\equiv0\operatorname{mod}9$ (by Proposition
\ref{prop.mod.0}), whence $s\equiv m\equiv0\operatorname{mod}9$ (here we are
tacitly using Proposition \ref{prop.mod.rules} \textbf{(c)}), which entails
$9\mid s$ (again by Proposition \ref{prop.mod.0}). Thus, we have shown that if
$9\mid m$, then $9\mid s$.

Conversely, if $9\mid s$, then $s\equiv0\operatorname{mod}9$ (by Proposition
\ref{prop.mod.0}), whence $m\equiv s\equiv0\operatorname{mod}9$, which in turn
entails $9\mid m$ (by Proposition \ref{prop.mod.0}). Thus, we have shown that
if $9\mid s$, then $9\mid m$.

Now we have proved that each of the statements $9\mid m$ and $9\mid s$ implies
the other. In other words, we have $9\mid m$ if and only if $9\mid s$. This
proves the proposition.
\end{proof}

In other words, Theorem \ref{thm.div.criteria} \textbf{(e)} is proven. A
similar argument (with $9$ replaced by $3$) can be used to prove Theorem
\ref{thm.div.criteria} \textbf{(d)}. In fact, $s\equiv m\operatorname{mod}9$
entails $s\equiv m\operatorname{mod}3$ by Proposition \ref{prop.mod.rules}
\textbf{(e)}, because $3\mid9$.

Parts \textbf{(a)} and \textbf{(b)} of Theorem \ref{thm.div.criteria} can be
proved along similar lines, but are in fact easier. Indeed, if $m\in
\mathbb{N}$ has decimal representation $m_{d}m_{d-1}\cdots m_{0}$, then
$m\equiv m_{0}\operatorname{mod}10$ (since the number $m-m_{0}$ has decimal
representation $m_{d}m_{d-1}\cdots m_{1}0$ and thus is divisible by $10$), and
therefore (by Proposition \ref{prop.mod.rules} \textbf{(e)}) we have $m\equiv
m_{0}\operatorname{mod}2$ and $m\equiv m_{0}\operatorname{mod}5$ as well.

\begin{exercise}
\label{exe.21finfin.3}Let $m$ be a positive integer, and let $m_{d}%
m_{d-1}\cdots m_{0}$ be its decimal representation, so that $m_{0}%
,m_{1},\ldots,m_{d}$ are digits satisfying%
\[
m=m_{d}\cdot10^{d}+m_{d-1}\cdot10^{d-1}+\cdots+m_{0}\cdot10^{0}=\sum_{k=0}%
^{d}m_{k}\cdot10^{k}.
\]
Let $a$ be the alternating sum of digits of $m$; this is defined by%
\[
a:=m_{0}-m_{1}+m_{2}-m_{3}\pm\cdots+\left(  -1\right)  ^{d}m_{d}=\sum
_{k=0}^{d}\left(  -1\right)  ^{k}m_{k}.
\]
Prove that $11\mid m$ if and only if $11\mid a$. (This is the classical
divisibility test for divisibility by $11$.)
\end{exercise}

\subsection{Division with remainder}

\subsubsection{The theorem}

What comes next is the most fundamental theorem of number theory:

\begin{theorem}
[division-with-remainder theorem]\label{thm.ent.quo-rem}Let $n$ be an integer.
Let $d$ be a positive integer. Then, there exists a \textbf{unique} pair
$\left(  q,r\right)  $ of integers%
\[
q\in\mathbb{Z}\ \ \ \ \ \ \ \ \ \ \text{and}\ \ \ \ \ \ \ \ \ \ r\in\left\{
0,1,\ldots,d-1\right\}
\]
such that%
\[
n=qd+r.
\]

\end{theorem}

We will prove this soon. First, let us introduce some notations:

\begin{definition}
\label{def.ent.quo-rem}Let $n$ be an integer. Let $d$ be a positive integer.
Let $\left(  q,r\right)  $ be the pair whose existence and uniqueness is
claimed in Theorem \ref{thm.ent.quo-rem}. Then:

\begin{itemize}
\item The number $q$ is called the \textbf{quotient} of the division of $n$ by
$d$, and will be denoted by $n//d$.

\item The number $r$ is called the \textbf{remainder} of the division of $n$
by $d$, and will be denoted by $n\%d$.

\item The pair $\left(  q,r\right)  $ is called the \textbf{quo-rem pair} of
$n$ and $d$.
\end{itemize}
\end{definition}

For now, of course, we do not yet know that these $q$ and $r$ exist and are
unique (because we haven't proved the theorem yet). Thus, we will take care to
speak of \textquotedblleft\textbf{a} quotient\textquotedblright,
\textquotedblleft\textbf{a} remainder\textquotedblright\ and \textquotedblleft%
\textbf{a} quo-rem pair\textquotedblright, never taking their existence and
uniqueness for granted until we have proved it.

\begin{example}
What are $8//5$ and $8\%5$ ? We have%
\[
\underbrace{8}_{=n}=\underbrace{1}_{=q}\cdot\underbrace{5}_{=d}+\underbrace{3}%
_{=r\in\left\{  0,1,2,3,4\right\}  },
\]
so $8//5=1$ and $8\%5=3$. (This is taking the uniqueness of $8//5$ and $8\%5$
for granted, but we will prove this soon.)
\end{example}

\begin{example}
What are $19//5$ and $19\%5$ ? We have $19=3\cdot5+4$, so $19//5=3$ and
$19\%5=4$.
\end{example}

\begin{example}
What are $\left(  -7\right)  //5$ and $\left(  -7\right)  \%5$ ? We have
\[
\underbrace{-7}_{=n}=\underbrace{\left(  -2\right)  }_{=q}\cdot\underbrace{5}%
_{=d}+\underbrace{3}_{=r\in\left\{  0,1,2,3,4\right\}  },
\]
so $\left(  -7\right)  //5=-2$ and $\left(  -7\right)  \%5=3$.
\end{example}

So Theorem \ref{thm.ent.quo-rem} is saying that for any integer $n$ and any
positive integer $d$, there is a unique quo-rem pair of $n$ and $d$. Let us
now prove this.

\subsubsection{The proof}

\begin{proof}
[Proof of Theorem \ref{thm.ent.quo-rem}.]We need to prove two things: that a
quo-rem pair of $n$ and $d$ exists, and that it is unique. Let me prove the
uniqueness part first. \medskip

\textit{Proof of the uniqueness part:} Fix an integer $n$ and a positive
integer $d$. We must show that there is \textbf{at most one} quo-rem pair
$\left(  q,r\right)  $ of $n$ and $d$. In other words, we must show that there
are no two distinct quo-rem pairs of $n$ and $d$.

We shall prove this by contradiction. So we assume that $\left(  q_{1}%
,r_{1}\right)  $ and $\left(  q_{2},r_{2}\right)  $ are two distinct quo-rem
pairs of $n$ and $d$. We want to derive a contradiction.

Since $\left(  q_{1},r_{1}\right)  $ is a quo-rem pair of $n$ and $d$, we have%
\[
q_{1}\in\mathbb{Z}\ \ \ \ \ \ \ \ \ \ \text{and}\ \ \ \ \ \ \ \ \ \ r_{1}%
\in\left\{  0,1,\ldots,d-1\right\}  \ \ \ \ \ \ \ \ \ \ \text{and}%
\ \ \ \ \ \ \ \ \ \ n=q_{1}d+r_{1}.
\]

Since $\left(  q_{2},r_{2}\right)  $ is a quo-rem pair of $n$ and $d$, we have%
\[
q_{2}\in\mathbb{Z}\ \ \ \ \ \ \ \ \ \ \text{and}\ \ \ \ \ \ \ \ \ \ r_{2}%
\in\left\{  0,1,\ldots,d-1\right\}  \ \ \ \ \ \ \ \ \ \ \text{and}%
\ \ \ \ \ \ \ \ \ \ n=q_{2}d+r_{2}.
\]

Subtracting the equation $n=q_{2}d+r_{2}$ from $n=q_{1}d+r_{1}$, we find%
\[
0=\left(  q_{1}d+r_{1}\right)  -\left(  q_{2}d+r_{2}\right)  =\left(
r_{1}-r_{2}\right)  -\left(  q_{2}d-q_{1}d\right)  =\left(  r_{1}%
-r_{2}\right)  -\left(  q_{2}-q_{1}\right)  d.
\]
In other words,%
\begin{equation}
r_{1}-r_{2}=\left(  q_{2}-q_{1}\right)  d. \label{pf.thm.ent.quo-rem.4}%
\end{equation}

We are in one of the following three cases:

\textit{Case 1:} We have $q_{1}<q_{2}$.

\textit{Case 2:} We have $q_{1}=q_{2}$.

\textit{Case 3:} We have $q_{1}>q_{2}$.

Let us first consider Case 1. In this case, we have $q_{1}<q_{2}$, so that
$q_{2}-q_{1}>0$. Since $q_{2}-q_{1}$ is an integer, this entails that
$q_{2}-q_{1}\geq1$. We can multiply this inequality by $d$ (since $d>0$), thus
obtaining $\left(  q_{2}-q_{1}\right)  d\geq1d=d$. In view of
(\ref{pf.thm.ent.quo-rem.4}), we can rewrite this as $r_{1}-r_{2}\geq d$.
However, $r_{1}\leq d-1$ (since $r_{1}\in\left\{  0,1,\ldots,d-1\right\}  $)
and $r_{2}\geq0$ (since $r_{2}\in\left\{  0,1,\ldots,d-1\right\}  $). Hence,
$r_{1}-\underbrace{r_{2}}_{\geq0}\leq r_{1}\leq d-1<d$. This contradicts
$r_{1}-r_{2}\geq d$. Thus, we have found a contradiction in Case 1.

Let us next consider Case 2. In this case, we have $q_{1}=q_{2}$. Hence, we
can rewrite (\ref{pf.thm.ent.quo-rem.4}) as $r_{1}-r_{2}=\underbrace{\left(
q_{2}-q_{2}\right)  }_{=0}d=0$, so that $r_{1}=r_{2}$. Combining $q_{1}=q_{2}$
with $r_{1}=r_{2}$, we obtain $\left(  q_{1},r_{1}\right)  =\left(
q_{2},r_{2}\right)  $, which contradicts our assumption that the two quo-rem
pairs $\left(  q_{1},r_{1}\right)  $ and $\left(  q_{2},r_{2}\right)  $ are
distinct. Thus, we have found a contradiction in Case 2.

Finally, in Case 3, we have $q_{1}>q_{2}$ and therefore $q_{2}<q_{1}$. Thus,
Case 3 is just a copy of Case 1 with the roles of the two pairs $\left(
q_{1},r_{1}\right)  $ and $\left(  q_{2},r_{2}\right)  $ switched (since the
two quo-rem pairs $\left(  q_{1},r_{1}\right)  $ and $\left(  q_{2}%
,r_{2}\right)  $ are playing identical roles). Hence, we obtain a
contradiction in Case 3 (since we obtained one in Case 1).

We have now obtained contradictions in all three Cases 1, 2 and 3. Thus, we
always have a contradiction. Hence, our assumption was wrong. This completes
our proof of the uniqueness of the quo-rem pair of $n$ and $d$. \medskip

Now, let us come to the existence part. It is reasonable to try induction, but
there is a hurdle: Induction on $d$ does not work (there is no good way to use
the induction hypothesis), whereas induction on $n$ cannot be used as long as
$n$ can be negative. Fortunately, the latter hurdle is surmountable. One way
around it is to \textbf{first} prove the existence of a quo-rem pair in the
case when $n\in\mathbb{N}$ (that is, $n\geq0$), and \textbf{afterwards}
generalize this result to arbitrary integers $n$.

So let us prove the $n\in\mathbb{N}$ case:

\begin{lemma}
\label{lem.ent.quo-rem.existN}Let $n\in\mathbb{N}$, and let $d$ be a positive
integer. Then, there exists a quo-rem pair of $n$ and $d$.
\end{lemma}

\begin{proof}
[Proof of Lemma \ref{lem.ent.quo-rem.existN}.]Fix $d$. We apply strong
induction on $n$:

\textit{Induction step:}\footnote{Recall that a strong induction needs no base
case (see Subsection \ref{subsec.ind.strind.baseless}).} Let $n\in\mathbb{N}$.
Assume (as the induction hypothesis) that Lemma \ref{lem.ent.quo-rem.existN}
is proved for all nonnegative integers smaller than $n$ instead of $n$. In
other words, assume that for each nonnegative integer $k<n$, there exists a
quo-rem pair of $k$ and $d$. We must prove that Lemma
\ref{lem.ent.quo-rem.existN} also holds for $n$, i.e., that there exists a
quo-rem pair of $n$ and $d$.

If $n<d$, then such a pair can be explicitly constructed: it is $\left(
0,n\right)  $. (Indeed, the pair $\left(  0,n\right)  $ satisfies the
requirements for a quo-rem pair, since $n=0d+n$ and $n\in\left\{
0,1,\ldots,d-1\right\}  $.)

Otherwise, we have $n\geq d$, so that $n-d\in\mathbb{N}$. Thus, we can apply
the induction hypothesis to $n-d$ instead of $n$ (since $n-d<n$). We conclude
that there exists a quo-rem pair of $n-d$ and $d$. We denote this pair by
$\left(  q,r\right)  $. Thus,%
\[
q\in\mathbb{Z}\ \ \ \ \ \ \ \ \ \ \text{and}\ \ \ \ \ \ \ \ \ \ r\in\left\{
0,1,\ldots,d-1\right\}  \ \ \ \ \ \ \ \ \ \ \text{and}%
\ \ \ \ \ \ \ \ \ \ n-d=qd+r.
\]
Now, I claim that $\left(  q+1,r\right)  $ is a quo-rem pair of $n$ and $d$.
Indeed, from $n-d=qd+r$, we conclude that%
\[
n=\left(  qd+r\right)  +d=qd+d+r=\left(  q+1\right)  d+r,
\]
which shows that $\left(  q+1,r\right)  $ is a quo-rem pair of $n$ and $d$
(since $r\in\left\{  0,1,\ldots,d-1\right\}  $). Thus, there exists a quo-rem
pair of $n$ and $d$. This completes our induction step, and thus Lemma
\ref{lem.ent.quo-rem.existN} is proved.
\end{proof}

We now return to proving Theorem \ref{thm.ent.quo-rem}. We have shown that

\begin{itemize}
\item there is always \textbf{at most one} quo-rem pair of $n$ and $d$, and

\item there is \textbf{at least one} quo-rem pair of $n$ and $d$ if
$n\in\mathbb{N}$.
\end{itemize}

What remains to be done is proving that there is \textbf{at least one} quo-rem
pair of $n$ and $d$ if $n<0$.

This can be done in several ways. One way is to proceed similarly to the proof
of Lemma \ref{lem.ent.quo-rem.existN}, but using strong induction on $-n$.

Alternatively, there is a slicker argument: We can reduce the
\textquotedblleft negative $n$\textquotedblright\ case to the
\textquotedblleft nonnegative $n$\textquotedblright\ case (which is already
covered by Lemma \ref{lem.ent.quo-rem.existN}). Namely, let $n\in\mathbb{Z}$
be negative. Then, the product $\left(  1-d\right)  n$ is nonnegative (since
both factors $1-d$ and $n$ are $\leq0$), so we can apply Lemma
\ref{lem.ent.quo-rem.existN} to $\left(  1-d\right)  n$ instead of $n$. Thus,
we conclude that there exists a quo-rem pair $\left(  q,r\right)  $ of
$\left(  1-d\right)  n$ and $d$. This pair $\left(  q,r\right)  $ satisfies%
\[
\left(  1-d\right)  n=qd+r
\]
(by the definition of a quo-rem pair). In other words,%
\[
n-dn=qd+r.
\]
Hence,%
\[
n=dn+qd+r=\left(  n+q\right)  d+r.
\]
This shows that $\left(  n+q,r\right)  $ is a quo-rem pair of $n$ and $d$.
Hence, such a quo-rem pair exists. Hence, we have proved the existence of a
quo-rem pair in the case when $n$ is negative. This completes our proof of
Theorem \ref{thm.ent.quo-rem}.
\end{proof}

\subsubsection{An application: even and odd integers}

We shall now use this theorem to derive some basic properties of even and odd
numbers. Recall what these words mean:

\begin{definition}
\textbf{(a)} An integer $n$ is said to be \textbf{even} if $2\mid n$. \medskip

\textbf{(b)} An integer $n$ is said to be \textbf{odd} if $2\nmid n$.
\end{definition}

In other words, an integer is called \textbf{even} if it is divisible by $2$,
and is called \textbf{odd} if it is not even.

Now we shall show the following:

\begin{proposition}
\label{prop.ent.odd=2k+1}Let $n$ be an integer. \medskip

\textbf{(a)} The integer $n$ is even if and only if there exists some
$k\in\mathbb{Z}$ such that $n=2k$. \medskip

\textbf{(b)} The integer $n$ is odd if and only if there exists some
$k\in\mathbb{Z}$ such that $n=2k+1$.
\end{proposition}

\begin{proof}
Part \textbf{(a)} is a direct consequence of the definition of divisibility.
But part \textbf{(b)} is not!

So let us prove part \textbf{(b)}. This is an \textquotedblleft if and only
if\textquotedblright\ statement, so we need to prove both directions:%
\[
\left(  n\text{ is odd}\right)  \Longrightarrow\left(  \text{there exists some
}k\in\mathbb{Z}\text{ such that }n=2k+1\right)
\]
and%
\[
\left(  \text{there exists some }k\in\mathbb{Z}\text{ such that }%
n=2k+1\right)  \Longrightarrow\left(  n\text{ is odd}\right)  .
\]

For the sake of brevity, I shall refer to these two directions as the
\textquotedblleft$\Longrightarrow$\textquotedblright\ and \textquotedblleft%
$\Longleftarrow$\textquotedblright\ directions (respectively).

\textit{Proof of the \textquotedblleft}$\Longrightarrow$%
\textit{\textquotedblright\ direction:} Assume that $n$ is odd. By Theorem
\ref{thm.ent.quo-rem}, there exists a quo-rem pair $\left(  q,r\right)  $ of
$n$ and $2$. Consider this $\left(  q,r\right)  $. By the definition of a
quo-rem pair, this pair satisfies%
\[
q\in\mathbb{Z}\ \ \ \ \ \ \ \ \ \ \text{and}\ \ \ \ \ \ \ \ \ \ r\in\left\{
0,1\right\}  \ \ \ \ \ \ \ \ \ \ \text{and}\ \ \ \ \ \ \ \ \ \ n=2q+r.
\]
If $r$ were $0$, then we would thus get $n=2q+\underbrace{r}_{=0}=2q$, which
would show that $n$ is even; but this is impossible because $n$ is odd.
Therefore, we must have $r\neq0$, so that $r=1$ (since $r\in\left\{
0,1\right\}  $). Thus, $n=2q+\underbrace{r}_{=1}=2q+1$. Hence, there exists
some $k\in\mathbb{Z}$ such that $n=2k+1$ (namely, $k=q$). Thus we have shown
the \textquotedblleft$\Longrightarrow$\textquotedblright\ direction.

\textit{Proof of the \textquotedblleft}$\Longleftarrow$%
\textit{\textquotedblright\ direction:} Assume that there exists some
$k\in\mathbb{Z}$ such that $n=2k+1$. Consider this $k$.

We must show that $n$ is odd. This means showing that $2\nmid n$. This means
proving that $n$ cannot be written as $2c$ for an integer $c$.

To prove this, we assume the contrary. That is, we assume that $n=2c$ for some
integer $c$. Consider this $c$.

Now, the two pairs $\left(  k,1\right)  $ and $\left(  c,0\right)  $ both are
quo-rem pairs of $n$ and $2$, because we have $n=2k+1$ and $n=2c=2c+0$ (and
$1$ and $0$ belong to $\left\{  0,1\right\}  $). However, Theorem
\ref{thm.ent.quo-rem} says that the quo-rem pair of $n$ and $2$ is unique, so
these two pairs $\left(  k,1\right)  $ and $\left(  c,0\right)  $ must be
identical. But this is absurd, since their second entries $1$ and $0$ are
different. So we find a contradiction. This concludes our proof that $n$ is
odd. Thus, we have shown the \textquotedblleft$\Longleftarrow$%
\textquotedblright\ direction of Proposition \ref{prop.ent.odd=2k+1}
\textbf{(b)}.

This completes the proof of Proposition \ref{prop.ent.odd=2k+1} \textbf{(b)}
(since both directions are proved).
\end{proof}

\begin{corollary}
\label{cor.ent.odd+odd}\textbf{(a)} The sum of any two even integers is even.
\medskip

\textbf{(b)} The sum of any even integer with any odd integer is odd. \medskip

\textbf{(c)} The sum of any two odd integers is even.
\end{corollary}

\begin{proof}
We will only prove part \textbf{(c)}, since the other two parts are analogous
(and even simpler).

\textbf{(c)} Let $a$ and $b$ be two odd integers. We must prove that $a+b$ is even.

The integer $a$ is odd. Hence, Proposition \ref{prop.ent.odd=2k+1}
\textbf{(b)} shows that we can write $a$ as $a=2k+1$ for some integer $k$.

Similarly, we can write $b$ as $b=2\ell+1$ for some integer $\ell$.

Consider these $k$ and $\ell$. Now, from $a=2k+1$ and $b=2\ell+1$, we obtain%
\[
a+b=\left(  2k+1\right)  +\left(  2\ell+1\right)  =2k+2\ell+2=2\left(
k+\ell+1\right)  ,
\]
which is clearly even. This proves Corollary \ref{cor.ent.odd+odd}
\textbf{(c)}.
\end{proof}

\begin{remark}
Corollary \ref{cor.ent.odd+odd} \textbf{(c)} is a property specific to the
number $2$. For example, it is not true that the sum of any two integers not
divisible by $3$ is divisible by $3$.
\end{remark}

\subsubsection{Basic properties of quotients and remainders}

Here are some elementary facts about quotients and remainders:

\begin{proposition}
\label{prop.ent.quo-rem.basics}Let $n\in\mathbb{Z}$, and let $d$ be a positive
integer. Then: \medskip

\textbf{(a)} We have $n\%d\in\left\{  0,1,\ldots,d-1\right\}  $ and
$n\%d\equiv n\operatorname{mod}d$. \medskip

\textbf{(b)} We have $d\mid n$ if and only if $n\%d=0$. \medskip

\textbf{(c)} If $c\in\left\{  0,1,\ldots,d-1\right\}  $ satisfies $c\equiv
n\operatorname{mod}d$, then $c=n\%d$. \medskip

\textbf{(d)} We have $n=\left(  n//d\right)  d+\left(  n\%d\right)  $.
\medskip

\textbf{(e)} If $n\in\mathbb{N}$, then $n//d\in\mathbb{N}$.
\end{proposition}

Note that part \textbf{(a)} of this proposition can be restated as follows:
The remainder $n\%d$ is an element of $\left\{  0,1,\ldots,d-1\right\}  $ that
is congruent to $n$ modulo $d$. Part \textbf{(c)} says that, conversely, any
element $c$ of $\left\{  0,1,\ldots,d-1\right\}  $ that is congruent to $n$
modulo $d$ must be this remainder $n\%d$. Thus, together, these two parts
uniquely characterize the remainder $n\%d$ as the only element of $\left\{
0,1,\ldots,d-1\right\}  $ that is congruent to $n$ modulo $d$. This
characterization is good to keep in mind, as it describes the remainder
independently of the quotient.

\begin{proof}
[Proof of Proposition \ref{prop.ent.quo-rem.basics}.]We set%
\[
q:=n//d\ \ \ \ \ \ \ \ \ \ \text{and}\ \ \ \ \ \ \ \ \ \ r:=n\%d.
\]
Thus, $\left(  q,r\right)  $ is a quo-rem pair of $n$ and $d$ (by the
definition of a quo-rem pair). In other words, we have $n=qd+r$ and
$q\in\mathbb{Z}$ and $r\in\left\{  0,1,\ldots,d-1\right\}  $. We can now prove
all five parts of the proposition: \medskip

\textbf{(d)} We have $n=\underbrace{q}_{=n//d}d+\underbrace{r}_{=n\%d}=\left(
n//d\right)  d+\left(  n\%d\right)  $. This proves Proposition
\ref{prop.ent.quo-rem.basics} \textbf{(d)}. \medskip

\textbf{(a)} We have $n\%d=r\in\left\{  0,1,\ldots,d-1\right\}  $. Moreover,
from $n=qd+r$, we obtain $r-n=r-\left(  qd+r\right)  =-qd$, which is clearly
divisible by $d$. Hence, $d\mid r-n$. Equivalently, $r\equiv
n\operatorname{mod}d$. In other words, $n\%d\equiv n\operatorname{mod}d$
(since $r=n\%d$). Thus, Proposition \ref{prop.ent.quo-rem.basics} \textbf{(a)}
is proved (since we have shown that $n\%d\in\left\{  0,1,\ldots,d-1\right\}  $
as well). \medskip

\textbf{(c)} Let $c\in\left\{  0,1,\ldots,d-1\right\}  $ satisfy $c\equiv
n\operatorname{mod}d$. We must show that $c=n\%d$.

From $c\equiv n\operatorname{mod}d$, we obtain $d\mid c-n$. In other words,
$c-n=de$ for some $e\in\mathbb{Z}$. Consider this $e$. From $c-n=de$, we
obtain $c=n+de$, so that $n=c-de=\left(  -e\right)  d+c$. This (combined with
$c\in\left\{  0,1,\ldots,d-1\right\}  $) shows that $\left(  -e,c\right)  $ is
a quo-rem pair of $n$ and $d$. However, $\left(  q,r\right)  $ is also a
quo-rem pair of $n$ and $d$ (by its definition). Since there is only one
quo-rem pair of $n$ and $d$ (by Theorem \ref{thm.ent.quo-rem}), this shows
that $\left(  -e,c\right)  =\left(  q,r\right)  $. Hence, $c=r=n\%d$. This
proves Proposition \ref{prop.ent.quo-rem.basics} \textbf{(c)}. \medskip

\textbf{(b)} Again, this is an \textquotedblleft if and only
if\textquotedblright\ statement, and we shall prove its \textquotedblleft%
$\Longrightarrow$\textquotedblright\ and \textquotedblleft$\Longleftarrow
$\textquotedblright\ directions separately:

$\Longrightarrow:$ Assume that $d\mid n$. We must prove that $n\%d=0$. In
other words, we must prove that $r=0$.

Indeed, $d\mid n$ yields that $n\equiv0\operatorname{mod}d$ (by Proposition
\ref{prop.mod.0}). In other words, $0\equiv n\operatorname{mod}d$. Since we
furthermore have $0\in\left\{  0,1,\ldots,d-1\right\}  $, we can thus apply
Proposition \ref{prop.ent.quo-rem.basics} \textbf{(c)} to $c=0$, and conclude
that $0=n\%d$. In other words, $n\%d=0$. This proves the \textquotedblleft%
$\Longrightarrow$\textquotedblright\ direction (i.e., it proves that if $d\mid
n$, then $n\%d=0$).

$\Longleftarrow:$ If $n\%d=0$, then $d\mid n$ because%
\[
n=qd+\underbrace{r}_{=n\%d=0}=qd.
\]
This proves the \textquotedblleft$\Longleftarrow$\textquotedblright%
\ direction. Thus, both directions are proved, so that Proposition
\ref{prop.ent.quo-rem.basics} \textbf{(b)} holds. \medskip

\textbf{(e)} Assume that $n\in\mathbb{N}$. Recall that $r\in\left\{
0,1,\ldots,d-1\right\}  $, so that $r\leq d-1<d$. But $n=qd+r$, so that
$qd+r=n\geq0$ (since $n\in\mathbb{N}$). In other words, $qd\geq-r>-d$ (since
$r<d$).

If we had $q<0$, then we would have $q\leq-1$ (since $q$ is an integer) and
therefore $qd\leq\left(  -1\right)  d$ (since we can multiply the inequality
$q\leq-1$ by the positive number $d$); but this would contradict
$qd>-d=\left(  -1\right)  d$. Hence, we cannot have $q<0$. Thus, $q\geq0$, so
that $q\in\mathbb{N}$. In other words, $n//d\in\mathbb{N}$ (since $q=n//d$).
This proves Proposition \ref{prop.ent.quo-rem.basics} \textbf{(e)}.
\end{proof}

\begin{corollary}
\label{cor.ent.parity-rem2}Let $n\in\mathbb{Z}$. Then: \medskip

\textbf{(a)} The integer $n$ is even if and only if $n\%2=0$. \medskip

\textbf{(b)} The integer $n$ is odd if and only if $n\%2=1$.
\end{corollary}

\begin{proof}
\textbf{(a)} We have the following chain of logical equivalences:%
\begin{align*}
\left(  n\text{ is even}\right)  \  &  \Longleftrightarrow\ \left(  2\mid
n\right)  \ \ \ \ \ \ \ \ \ \ \left(  \text{by the definition of
\textquotedblleft even\textquotedblright}\right) \\
&  \Longleftrightarrow\ \left(  n\%2=0\right)  \ \ \ \ \ \ \ \ \ \ \left(
\text{by Proposition \ref{prop.ent.quo-rem.basics} \textbf{(b)}, applied to
}d=2\right)  .
\end{align*}
Hence, Corollary \ref{cor.ent.parity-rem2} \textbf{(a)} is proved. \medskip

\textbf{(b)} Proposition \ref{prop.ent.quo-rem.basics} \textbf{(a)} (applied
to $d=2$) yields $n\%2\in\left\{  0,1,\ldots,2-1\right\}  $ and $n\%2\equiv
n\operatorname{mod}2$. Thus, $n\%2\in\left\{  0,1,\ldots,2-1\right\}
=\left\{  0,1\right\}  $. Hence, $n\%2$ is either $0$ or $1$. Thus,
$n\%2\neq0$ holds if and only if $n\%2=1$. In other words, we have the logical
equivalence $\left(  n\%2\neq0\right)  \ \Longleftrightarrow\ \left(
n\%2=1\right)  $.

However, we have the following chain of logical equivalences:%
\begin{align*}
\left(  n\text{ is odd}\right)  \  &  \Longleftrightarrow\ \left(  2\nmid
n\right)  \ \ \ \ \ \ \ \ \ \ \left(  \text{by the definition of
\textquotedblleft odd\textquotedblright}\right) \\
&  \Longleftrightarrow\ \left(  \text{not }2\mid n\right) \\
&  \Longleftrightarrow\ \left(  \text{not }n\%2=0\right) \\
&  \ \ \ \ \ \ \ \ \ \ \ \ \ \ \ \ \ \ \ \ \left(
\begin{array}
[c]{c}%
\text{since Proposition \ref{prop.ent.quo-rem.basics} \textbf{(b)} (applied to
}d=2\text{)}\\
\text{yields that }2\mid n\text{ holds if and only if }n\%2=0
\end{array}
\right) \\
&  \Longleftrightarrow\ \left(  n\%2\neq0\right) \\
&  \Longleftrightarrow\ \left(  n\%2=1\right)  .
\end{align*}
This proves Corollary \ref{cor.ent.parity-rem2} \textbf{(b)}.
\end{proof}

Quotients and remainders are closely connected to the so-called
\href{https://en.wikipedia.org/wiki/Floor_and_ceiling_functions}{floor
function}:

\begin{definition}
\label{def.ent.floor}The \textbf{integer part} (aka \textbf{floor}) of a real
number $x$ is defined to be the largest integer that is $\leq x$. It is
denoted by $\left\lfloor x\right\rfloor $.
\end{definition}

For example,%
\begin{align*}
\left\lfloor 3.8\right\rfloor  &  =3,\ \ \ \ \ \ \ \ \ \ \left\lfloor
4.2\right\rfloor =4,\ \ \ \ \ \ \ \ \ \ \left\lfloor 5\right\rfloor
=5,\ \ \ \ \ \ \ \ \ \ \left\lfloor \sqrt{2}\right\rfloor =1,\\
\left\lfloor \pi\right\rfloor  &  =3,\ \ \ \ \ \ \ \ \ \ \left\lfloor
0.5\right\rfloor =0,\ \ \ \ \ \ \ \ \ \ \left\lfloor -1.2\right\rfloor =-2
\end{align*}
(make sure you understand the last example! $-1$ is not $\leq-1.2$, but $-2$ is).

Now, here is the connection to quotients and remainders:

\begin{proposition}
[\textquotedblleft explicit formulas\textquotedblright\ for quotient and
remainder]\label{prop.ent.quo-rem-floor}Let $n\in\mathbb{Z}$, and let $d$ be a
positive integer. Then,%
\[
n//d=\left\lfloor \dfrac{n}{d}\right\rfloor \ \ \ \ \ \ \ \ \ \ \text{and}%
\ \ \ \ \ \ \ \ \ \ n\%d=n-d\cdot\left\lfloor \dfrac{n}{d}\right\rfloor .
\]

\end{proposition}

\begin{proof}
Proposition \ref{prop.ent.quo-rem.basics} \textbf{(a)} yields $n\%d\in\left\{
0,1,\ldots,d-1\right\}  $. Hence, $n\%d\geq0$ and $n\%d\leq d-1<d$.

Proposition \ref{prop.ent.quo-rem.basics} \textbf{(d)} yields $n=\left(
n//d\right)  d+\left(  n\%d\right)  $. Thus,%
\[
n=\left(  n//d\right)  d+\underbrace{\left(  n\%d\right)  }_{<d}<\left(
n//d\right)  d+d=\left(  \left(  n//d\right)  +1\right)  d.
\]
Dividing both sides of this inequality by $d$ (we can do this, since $d>0$),
we obtain $\dfrac{n}{d}<\left(  n//d\right)  +1$. In other words, $\left(
n//d\right)  +1>\dfrac{n}{d}$.

On the other hand,
\[
n=\left(  n//d\right)  d+\underbrace{\left(  n\%d\right)  }_{\geq0}\geq\left(
n//d\right)  d.
\]
Dividing both sides of this inequality by $d$ (we can do this, since $d>0$),
we obtain $\dfrac{n}{d}\geq n//d$.

Now, the integer $n//d$ is $\leq\dfrac{n}{d}$ (since $\dfrac{n}{d}\geq n//d$),
but the next-larger integer $\left(  n//d\right)  +1$ is not (since $\left(
n//d\right)  +1>\dfrac{n}{d}$). Thus, $n//d$ is the largest integer that is
$\leq\dfrac{n}{d}$. In other words, $n//d=\left\lfloor \dfrac{n}%
{d}\right\rfloor $ (by the definition of the floor $\left\lfloor \dfrac{n}%
{d}\right\rfloor $).

Solving the equation $n=\left(  n//d\right)  d+\left(  n\%d\right)  $ for
$n\%d$, we find%
\[
n\%d=n-\underbrace{\left(  n//d\right)  }_{=\left\lfloor \dfrac{n}%
{d}\right\rfloor }d=n-\left\lfloor \dfrac{n}{d}\right\rfloor d=n-d\cdot
\left\lfloor \dfrac{n}{d}\right\rfloor .
\]
Thus, Proposition \ref{prop.ent.quo-rem-floor} is proved.
\end{proof}

Division with remainder is one of the most fundamental facts about integers;
almost all of number theory is downstream of it. Here are some applications:

\begin{exercise}
\label{exe.23hw3.3}Let $n$ be any integer. Prove the following: \medskip

\textbf{(a)} If $n$ is odd, then $8\mid n^{2}-1$. \medskip

\textbf{(b)} If $3\nmid n$, then $3\mid n^{2}-1$. \medskip

[\textbf{Hint:} In part \textbf{(a)}, write $n$ as $2k+1$. In part
\textbf{(b)}, write $n$ as $q\cdot3+r$ and consider the possible values for
$r$.]
\end{exercise}

\begin{exercise}
\label{exe.23hw3.4}Let $p$ be a positive integer.

Assume that you are given $p$-cent coins and $\left(  p+1\right)  $-cent coins
(each in infinite supply).

Prove that you can pay $n$ cents using these coins for every integer $n\geq
p^{2}-p$.

In other words, prove that each integer $n\geq p^{2}-p$ can be written as
$a\left(  p+1\right)  +bp$ with $a,b\in\mathbb{N}$.
\end{exercise}

\subsubsection{Base-$b$ representation of nonnegative integers}

Division with remainder is the main ingredient in a feature of integers that
you may well be taking for granted, but actually needs to be proved: the fact
that every integer can be uniquely expressed in decimal notation, or, more
generally, in base-$b$ notation for any given integer $b>1$.

What does this mean? For example,%
\begin{align*}
3401  &  =3\cdot1000+4\cdot100+0\cdot10+1\cdot1\\
&  =3\cdot10^{3}+4\cdot10^{2}+0\cdot10^{1}+1\cdot10^{0}.
\end{align*}
Thus, we have written the fairly large number $3401$ as a pretty short sum of
powers of $10$, with the coefficients being integers between $0$ and $9$
(commonly known as \textquotedblleft digits\textquotedblright).

This can be done for any nonnegative integer $n$ instead of $3401$. This can
also be done with any fixed integer $b>1$ instead of $10$, except that the
coefficients (\textquotedblleft generalized digits\textquotedblright) will
then be integers between $0$ and $b-1$. This is called the \textquotedblleft
base-$b$ representation\textquotedblright\ of the integer $n$.

For instance, let us find the base-$4$ representation of the integer $3401$:
This will be a representation of $3401$ in the form%
\[
3401=r_{6}4^{6}+r_{5}4^{5}+r_{4}4^{4}+r_{3}4^{3}+r_{2}4^{2}+r_{1}4^{1}%
+r_{0}4^{0},
\]
where each $r_{i}$ is a \textquotedblleft base-$4$ digit\textquotedblright%
\ (i.e., an element of $\left\{  0,1,2,3\right\}  $). Here, we are tacitly
assuming that $4^{6}$ is the highest power of $4$ that we need; but we don't
actually know this yet, so we must be prepared to add higher powers
($4^{7},4^{8},4^{9},\ldots$) if needed.

How do we find these base-$4$ digits $r_{0},r_{1},\ldots,r_{6}$?

We start by identifying $r_{0}$. Indeed, on the RHS\footnote{\textquotedblleft
RHS\textquotedblright\ means \textquotedblleft right hand
side\textquotedblright.} of the equation%
\[
3401=r_{6}4^{6}+r_{5}4^{5}+r_{4}4^{4}+r_{3}4^{3}+r_{2}4^{2}+r_{1}4^{1}%
+r_{0}4^{0},
\]
all but the last addends are multiples of $4$, whereas the last addend is
$r_{0}4^{0}=r_{0}$. Hence, we can rewrite this equation as follows (factoring
out the $4$):%
\[
3401=4\cdot\left(  r_{6}4^{5}+r_{5}4^{4}+r_{4}4^{3}+r_{3}4^{2}+r_{2}%
4^{1}+r_{1}4^{0}\right)  +r_{0}.
\]
Since $r_{0}\in\left\{  0,1,2,3\right\}  $, this equation reveals that the
pair%
\[
\left(  r_{6}4^{5}+r_{5}4^{4}+r_{4}4^{3}+r_{3}4^{2}+r_{2}4^{1}+r_{1}%
4^{0},\ r_{0}\right)
\]
is a quo-rem pair of $3401$ and $4$. In particular, we must have%
\begin{align*}
r_{0}  &  =3401\%4=1\ \ \ \ \ \ \ \ \ \ \text{and}\\
r_{6}4^{5}+r_{5}4^{4}+r_{4}4^{3}+r_{3}4^{2}+r_{2}4^{1}+r_{1}4^{0}  &
=3401//4=850.
\end{align*}

Thus, we have identified the last base-$4$ digit $r_{0}$ as $1$. In order to
find the remaining digits, we analyze the latter equation%
\[
850=r_{6}4^{5}+r_{5}4^{4}+r_{4}4^{3}+r_{3}4^{2}+r_{2}4^{1}+r_{1}4^{0}.
\]
In this equation, all but the last addends are multiples of $4$, whereas the
last addend is $r_{1}4^{0}=r_{1}$. Hence, we can rewrite this equation as%
\[
850=4\cdot\left(  r_{6}4^{4}+r_{5}4^{3}+r_{4}4^{2}+r_{3}4^{1}+r_{2}%
4^{0}\right)  +r_{1},
\]
and thus conclude that
\begin{align*}
r_{1}  &  =850\%4=2\ \ \ \ \ \ \ \ \ \ \text{and}\\
r_{6}4^{4}+r_{5}4^{3}+r_{4}4^{2}+r_{3}4^{1}+r_{2}4^{0}  &  =850//4=212.
\end{align*}

Thus, we have identified the base-$4$ digit $r_{1}$ as $2$. In order to find
the remaining digits, we analyze the latter equation%
\[
212=r_{6}4^{4}+r_{5}4^{3}+r_{4}4^{2}+r_{3}4^{1}+r_{2}4^{0}.
\]
In this equation, all but the last addends are multiples of $4$, whereas the
last addend is $r_{2}4^{0}=r_{2}$. Hence, we can rewrite this equation as%
\[
212=4\cdot\left(  r_{6}4^{3}+r_{5}4^{2}+r_{4}4^{1}+r_{3}4^{0}\right)  +r_{2},
\]
and thus conclude that
\begin{align*}
r_{2}  &  =212\%4=0\ \ \ \ \ \ \ \ \ \ \text{and}\\
r_{6}4^{3}+r_{5}4^{2}+r_{4}4^{1}+r_{3}4^{0}  &  =212//4=53.
\end{align*}

Thus, we have identified the base-$4$ digit $r_{2}$ as $0$. In order to find
the remaining digits, we analyze the latter equation%
\[
53=r_{6}4^{3}+r_{5}4^{2}+r_{4}4^{1}+r_{3}4^{0}.
\]
In this equation, all but the last addends are multiples of $4$, whereas the
last addend is $r_{3}4^{0}=r_{3}$. Hence, we can rewrite this equation as%
\[
53=4\cdot\left(  r_{6}4^{2}+r_{5}4^{1}+r_{4}4^{0}\right)  +r_{3},
\]
and thus conclude that
\begin{align*}
r_{3}  &  =53\%4=1\ \ \ \ \ \ \ \ \ \ \text{and}\\
r_{6}4^{2}+r_{5}4^{1}+r_{4}4^{0}  &  =53//4=13.
\end{align*}

Thus, we have identified the base-$4$ digit $r_{3}$ as $1$. In order to find
the remaining digits, we analyze the latter equation%
\[
13=r_{6}4^{2}+r_{5}4^{1}+r_{4}4^{0}.
\]
In this equation, all but the last addends are multiples of $4$, whereas the
last addend is $r_{4}4^{0}=r_{4}$. Hence, we can rewrite this equation as%
\[
13=4\cdot\left(  r_{6}4^{1}+r_{5}4^{0}\right)  +r_{4},
\]
and thus conclude that
\begin{align*}
r_{4}  &  =13\%4=1\ \ \ \ \ \ \ \ \ \ \text{and}\\
r_{6}4^{1}+r_{5}4^{0}  &  =13//4=3.
\end{align*}

Thus, we have identified the base-$4$ digit $r_{4}$ as $1$. In order to find
the remaining digits, we analyze the latter equation%
\[
3=r_{6}4^{1}+r_{5}4^{0}.
\]
In this equation, all but the last addends are multiples of $4$, whereas the
last addend is $r_{5}4^{0}=r_{5}$. Hence, we can rewrite this equation as%
\[
3=4\cdot\left(  r_{6}4^{0}\right)  +r_{5},
\]
and thus conclude that
\begin{align*}
r_{5}  &  =3\%4=3\ \ \ \ \ \ \ \ \ \ \text{and}\\
r_{6}4^{0}  &  =3//4=0.
\end{align*}

Thus, we have identified the base-$4$ digit $r_{5}$ as $3$. Moreover, the
equation $r_{6}4^{0}=0$ shows that $r_{6}=0$.

Thus, altogether, we have found the representation of $3401$ we were looking
for:%
\[
3401=\underbrace{r_{6}}_{=0}4^{6}+\underbrace{r_{5}}_{=3}4^{5}%
+\underbrace{r_{4}}_{=1}4^{4}+\underbrace{r_{3}}_{=1}4^{3}+\underbrace{r_{2}%
}_{=0}4^{2}+\underbrace{r_{1}}_{=2}4^{1}+\underbrace{r_{0}}_{=1}4^{0}.
\]

In analogy to the decimal system, we can state this as \textquotedblleft the
number $3401$ written in base-$4$ is $0311021$\textquotedblright\ (since the
base-$4$ digits $r_{6},r_{5},\ldots,r_{0}$ have been identified as
$0,3,1,1,0,2,1$). Commonly, one would omit the leading zeroes, so this would
become $311021$.

The method we just used can be used for any given integer $b>1$ instead of $4$
and any nonnegative integer $n\in\mathbb{N}$ instead of $3401$: To find the
\textquotedblleft base-$b$ digits\textquotedblright\ of a nonnegative integer
$n$, we first divide $n$ by $b$ with remainder, then divide the resulting
quotient again by $b$ with remainder, then divide the resulting quotient again
by $b$ with remainder, and so on, until we are left with the quotient $0$. The
remainders obtained in the process will then be the base-$b$ digits of $n$
(from right to left). This process must eventually come to an end because
(since $b>1$) each quotient will be smaller than the preceding one.

We can summarize this as a theorem:

\begin{theorem}
\label{thm.ent.base-b.main}Let $b>1$ be an integer. Let $n\in\mathbb{N}$.
Then: \medskip

\textbf{(a)} We can write $n$ in the form%
\[
n=r_{k}\cdot b^{k}+r_{k-1}\cdot b^{k-1}+\cdots+r_{1}\cdot b^{1}+r_{0}\cdot
b^{0}%
\]
with%
\[
k\in\mathbb{N}\ \ \ \ \ \ \ \ \ \ \text{and}\ \ \ \ \ \ \ \ \ \ r_{0}%
,r_{1},\ldots,r_{k}\in\left\{  0,1,\ldots,b-1\right\}  .
\]

\textbf{(b)} If $n<b^{k+1}$ for some $k\in\mathbb{N}$, then we can write $n$
in the form%
\[
n=r_{k}\cdot b^{k}+r_{k-1}\cdot b^{k-1}+\cdots+r_{1}\cdot b^{1}+r_{0}\cdot
b^{0}%
\]
with%
\[
r_{0},r_{1},\ldots,r_{k}\in\left\{  0,1,\ldots,b-1\right\}  .
\]

\textbf{(c)} These $r_{0},r_{1},\ldots,r_{k}$ are unique (when $k$ is given).
Moreover, they can be explicitly computed by the formula%
\[
r_{i}=\left(  n//b^{i}\right)  \%b\ \ \ \ \ \ \ \ \ \ \text{for each }%
i\in\left\{  0,1,\ldots,k\right\}  .
\]
That is, they can be explicitly computed by%
\begin{align*}
r_{0}  &  =n\%b,\\
r_{1}  &  =\left(  n//b\right)  \%b,\\
r_{2}  &  =\left(  n//b^{2}\right)  \%b,\\
r_{3}  &  =\left(  n//b^{3}\right)  \%b,\\
&  \ldots,\\
r_{k}  &  =\left(  n//b^{k}\right)  \%b.
\end{align*}

\end{theorem}

\begin{fineprint}
\begin{proof}
Forget that $n$ was fixed (but keep $b$ fixed). We shall prove the following
two claims:

\begin{statement}
\textit{Claim 1:} Let $n\in\mathbb{N}$ and $k\in\mathbb{N}$ be such that
$n<b^{k+1}$. Then, we can write $n$ in the form%
\[
n=r_{k}\cdot b^{k}+r_{k-1}\cdot b^{k-1}+\cdots+r_{1}\cdot b^{1}+r_{0}\cdot
b^{0}%
\]
with%
\[
r_{0},r_{1},\ldots,r_{k}\in\left\{  0,1,\ldots,b-1\right\}  .
\]

\end{statement}

\begin{statement}
\textit{Claim 2:} Let $n\in\mathbb{N}$ and $k\in\mathbb{N}$. Assume that $n$
has been written in the form%
\[
n=r_{k}\cdot b^{k}+r_{k-1}\cdot b^{k-1}+\cdots+r_{1}\cdot b^{1}+r_{0}\cdot
b^{0}%
\]
with%
\[
r_{0},r_{1},\ldots,r_{k}\in\left\{  0,1,\ldots,b-1\right\}  .
\]
Then,
\[
r_{i}=\left(  n//b^{i}\right)  \%b\ \ \ \ \ \ \ \ \ \ \text{for each }%
i\in\left\{  0,1,\ldots,k\right\}  .
\]

\end{statement}

Once these two claims are proved, Theorem \ref{thm.ent.base-b.main} will
follow, because

\begin{itemize}
\item Theorem \ref{thm.ent.base-b.main} \textbf{(b)} follows directly from
Claim 1.

\item Theorem \ref{thm.ent.base-b.main} \textbf{(c)} follows directly from
Claim 2.

\item Theorem \ref{thm.ent.base-b.main} \textbf{(a)} follows from Claim 1
(since we can pick $k\in\mathbb{N}$ high enough that $n<b^{k+1}$
holds\footnote{Indeed, the assumption $b>1$ ensures that the sequence $\left(
b^{0},b^{1},b^{2},\ldots\right)  $ is strictly increasing and thus eventually
outgrows any given integer, including our $n$. Or we can argue this directly:
An easy induction (on $n$) shows that $n<b^{n+1}$, and thus we can simply take
$k=n$.}).
\end{itemize}

Hence, it remains to prove Claim 1 and Claim 2.

\begin{proof}
[Proof of Claim 1.]We proceed by induction on $k$:

\textit{Base case:} For $k=0$, Claim 1 is saying that every $n\in\mathbb{N}$
satisfying $n<b$ can be written in the form $n=r_{0}\cdot b^{0}$ with
$r_{0}\in\left\{  0,1,\ldots,b-1\right\}  $. But this is obvious: Since
$n\in\mathbb{N}$ and $n<b$, we have $n\in\left\{  0,1,\ldots,b-1\right\}  $,
and thus we can just pick $r_{0}=n$ and have $n=r_{0}\cdot b^{0}$ (since
$r_{0}\cdot\underbrace{b^{0}}_{=1}=r_{0}=n$). Hence, Claim 1 is proved for
$k=0$.

\textit{Induction step:} We make a step from $k-1$ to $k$. Thus, we let $k$ be
a positive integer. Assume (as the induction hypothesis) that Claim 1 holds
for $k-1$ instead of $k$. We must now show that Claim 1 holds for $k$ as well.

So let $n\in\mathbb{N}$ be such that $n<b^{k+1}$. Then, Proposition
\ref{prop.ent.quo-rem.basics} \textbf{(e)} (applied to $d=b$) yields
$n//b\in\mathbb{N}$. Moreover, $n\%b\in\left\{  0,1,\ldots,b-1\right\}  $ (by
the definition of a remainder). Hence, $n\%b\geq0$. Now, Proposition
\ref{prop.ent.quo-rem.basics} \textbf{(d)} (applied to $d=b$) yields
\[
n=\left(  n//b\right)  b+\underbrace{\left(  n\%b\right)  }_{\substack{\geq
0}}\geq\left(  n//b\right)  b.
\]
Hence, $\left(  n//b\right)  b\leq n<b^{k+1}$. Dividing this inequality by the
positive number $b$, we obtain $n//b<b^{k+1}/b=b^{k}$.

Now, recall our induction hypothesis, which says that Claim 1 holds for $k-1$
instead of $k$. In other words, if $m\in\mathbb{N}$ is such that $m<b^{\left(
k-1\right)  +1}$, then we can write $m$ in the form\footnote{We are
deliberately using the letters $m$ and $s_{i}$ instead of $n$ and $r_{i}$
here, since the letter $n$ is already taken (and the letters $r_{i}$ will be
needed for something different).}%
\[
m=s_{k-1}\cdot b^{k-1}+s_{k-2}\cdot b^{k-2}+\cdots+s_{1}\cdot b^{1}+s_{0}\cdot
b^{0}%
\]
with%
\[
s_{0},s_{1},\ldots,s_{k-1}\in\left\{  0,1,\ldots,b-1\right\}  .
\]
We can apply this to $m=n//b$ (since $n//b\in\mathbb{N}$ and $n//b<b^{k}%
=b^{\left(  k-1\right)  +1}$), and conclude that we can write $n//b$ in the
form%
\[
n//b=s_{k-1}\cdot b^{k-1}+s_{k-2}\cdot b^{k-2}+\cdots+s_{1}\cdot b^{1}%
+s_{0}\cdot b^{0}%
\]
with%
\[
s_{0},s_{1},\ldots,s_{k-1}\in\left\{  0,1,\ldots,b-1\right\}  .
\]
Let us do this. Thus,%
\begin{align*}
n  &  =\underbrace{\left(  n//b\right)  }_{=s_{k-1}\cdot b^{k-1}+s_{k-2}\cdot
b^{k-2}+\cdots+s_{1}\cdot b^{1}+s_{0}\cdot b^{0}}b+\left(  n\%b\right) \\
&  =\left(  s_{k-1}\cdot b^{k-1}+s_{k-2}\cdot b^{k-2}+\cdots+s_{1}\cdot
b^{1}+s_{0}\cdot b^{0}\right)  b+\left(  n\%b\right) \\
&  =s_{k-1}\cdot b^{k}+s_{k-2}\cdot b^{k-1}+\cdots+s_{1}\cdot b^{2}+s_{0}\cdot
b^{1}+\underbrace{\left(  n\%b\right)  }_{=\left(  n\%b\right)  \cdot b^{0}}\\
&  =s_{k-1}\cdot b^{k}+s_{k-2}\cdot b^{k-1}+\cdots+s_{1}\cdot b^{2}+s_{0}\cdot
b^{1}+\left(  n\%b\right)  \cdot b^{0}.
\end{align*}
Note that the coefficients $n\%b,s_{0},s_{1},\ldots,s_{k-1}$ on the right hand
side here all belong to $\left\{  0,1,\ldots,b-1\right\}  $ (as we know).
Thus, through this equality, we have written $n$ in the form%
\[
n=r_{k}\cdot b^{k}+r_{k-1}\cdot b^{k-1}+\cdots+r_{1}\cdot b^{1}+r_{0}\cdot
b^{0}%
\]
with%
\[
r_{0},r_{1},\ldots,r_{k}\in\left\{  0,1,\ldots,b-1\right\}
\]
(namely, with $r_{0}=n\%b$ and $r_{1}=s_{0}$ and $r_{2}=s_{1}$ and $\ldots$
and $r_{k-1}=s_{k-2}$ and $r_{k}=s_{k-1}$). Hence, $n$ can be written in this form.

We have thus proved that if $n\in\mathbb{N}$ is such that $n<b^{k+1}$, then we
can write $n$ in the form%
\[
n=r_{k}\cdot b^{k}+r_{k-1}\cdot b^{k-1}+\cdots+r_{1}\cdot b^{1}+r_{0}\cdot
b^{0}%
\]
with%
\[
r_{0},r_{1},\ldots,r_{k}\in\left\{  0,1,\ldots,b-1\right\}  .
\]
In other words, we have proved Claim 1 for our $k$. This completes the
induction step. Thus, Claim 1 is proved by induction.
\end{proof}

\begin{proof}
[Proof of Claim 2.]We could prove this by induction as well, but let us
instead go for a direct proof.

By assumption, we have%
\[
n=r_{k}\cdot b^{k}+r_{k-1}\cdot b^{k-1}+\cdots+r_{1}\cdot b^{1}+r_{0}\cdot
b^{0}=\sum_{j=0}^{k}r_{j}\cdot b^{j}=\sum_{j=0}^{k}r_{j}b^{j}.
\]

Now, we must prove that $r_{i}=\left(  n//b^{i}\right)  \%b$ for each
$i\in\left\{  0,1,\ldots,k\right\}  $. So let us fix an $i\in\left\{
0,1,\ldots,k\right\}  $.

We have
\begin{equation}
n=\sum_{j=0}^{k}r_{j}b^{j}=\sum_{j=0}^{i-1}r_{j}b^{j}+\sum_{j=i}^{k}r_{j}b^{j}
\label{pf.thm.ent.base-b.main.c2.pf.1}%
\end{equation}
(here, we have split our sum into two parts: one part which contains the
addends for $j\in\left\{  0,1,\ldots,i-1\right\}  $, and one part which
contains the addends for $j\in\left\{  i,i+1,\ldots,k\right\}  $). We can
rewrite the second sum as follows:%
\[
\sum_{j=i}^{k}r_{j}\underbrace{b^{j}}_{=b^{i}b^{j-i}}=\sum_{j=i}^{k}r_{j}%
b^{i}b^{j-i}=b^{i}\sum_{j=i}^{k}r_{j}b^{j-i}.
\]
Thus, we can rewrite (\ref{pf.thm.ent.base-b.main.c2.pf.1}) as%
\begin{equation}
n=\sum_{j=0}^{i-1}r_{j}b^{j}+b^{i}\sum_{j=i}^{k}r_{j}b^{j-i}.
\label{pf.thm.ent.base-b.main.c2.pf.2}%
\end{equation}

Let us set%
\[
q^{\prime}:=\sum_{j=i}^{k}r_{j}b^{j-i}\ \ \ \ \ \ \ \ \ \ \text{and}%
\ \ \ \ \ \ \ \ \ \ r^{\prime}:=\sum_{j=0}^{i-1}r_{j}b^{j}.
\]
With these notations, we can rewrite (\ref{pf.thm.ent.base-b.main.c2.pf.2})
as
\begin{equation}
n=r^{\prime}+b^{i}q^{\prime}=q^{\prime}b^{i}+r^{\prime}.
\label{pf.thm.ent.base-b.main.c2.pf.3}%
\end{equation}

Note that both sums $q^{\prime}=\sum_{j=i}^{k}r_{j}b^{j-i}$ and $r^{\prime
}=\sum_{j=0}^{i-1}r_{j}b^{j}$ are integers (indeed, $b^{j-i}$ is always an
integer in the first sum, since $j\geq i$ entails $j-i\in\mathbb{N}$).

We have assumed that $r_{0},r_{1},\ldots,r_{k}\in\left\{  0,1,\ldots
,b-1\right\}  $. In particular, the integers $r_{0},r_{1},\ldots,r_{k}$ are
all $\geq0$ and $\leq b-1$. In other words, each $j\in\left\{  0,1,\ldots
,k\right\}  $ satisfies $r_{j}\geq0$ and $r_{j}\leq b-1$. Hence, $r^{\prime
}=\sum_{j=0}^{i-1}r_{j}b^{j}\geq0$ (since all the integers $r_{j}$ are $\geq
0$, and so is $b$) and%
\[
r^{\prime}=\sum_{j=0}^{i-1}\underbrace{r_{j}}_{\leq b-1}b^{j}\leq\sum
_{j=0}^{i-1}\left(  b-1\right)  b^{j}=\left(  b-1\right)  \underbrace{\sum
_{j=0}^{i-1}b^{j}}_{\substack{=b^{0}+b^{1}+\cdots+b^{i-1}\\=\dfrac{b^{i}%
-1}{b-1}\\\text{(by Corollary \ref{cor.geo-sum.q},}\\\text{applied to }b\text{
and }i\\\text{instead of }q\text{ and }n\text{)}}}=\left(  b-1\right)
\cdot\dfrac{b^{i}-1}{b-1}=b^{i}-1.
\]
Thus, $r^{\prime}\in\left\{  0,1,\ldots,b^{i}-1\right\}  $.

The equality (\ref{pf.thm.ent.base-b.main.c2.pf.3}) says that $n=q^{\prime
}b^{i}+r^{\prime}$. In light of $q^{\prime}\in\mathbb{Z}$ and $r^{\prime}%
\in\left\{  0,1,\ldots,b^{i}-1\right\}  $, this shows that $\left(  q^{\prime
},r^{\prime}\right)  $ is a quo-rem pair of $n$ and $b^{i}$. Therefore, in
particular, $q^{\prime}$ is the quotient of the division of $n$ by $b^{i}$. In
other words,%
\[
q^{\prime}=n//b^{i}.
\]

However,%
\begin{align*}
q^{\prime}  &  =\sum_{j=i}^{k}r_{j}b^{j-i}=r_{i}b^{0}+r_{i+1}b^{1}%
+r_{i+2}b^{2}+\cdots+r_{k}b^{k-i}\\
&  =r_{i}\underbrace{b^{0}}_{=1}+\underbrace{\left(  r_{i+1}b^{1}+r_{i+2}%
b^{2}+\cdots+r_{k}b^{k-i}\right)  }_{=\left(  r_{i+1}b^{0}+r_{i+2}b^{1}%
+\cdots+r_{k}b^{k-i-1}\right)  b}\\
&  =r_{i}+\left(  r_{i+1}b^{0}+r_{i+2}b^{1}+\cdots+r_{k}b^{k-i-1}\right)  b.
\end{align*}
Thus, $q^{\prime}-r_{i}=\left(  r_{i+1}b^{0}+r_{i+2}b^{1}+\cdots
+r_{k}b^{k-i-1}\right)  b$, which is clearly divisible by $b$. That is, $b\mid
q^{\prime}-r_{i}$. In other words, $q^{\prime}\equiv r_{i}\operatorname{mod}%
b$. In other words, $r_{i}\equiv q^{\prime}\operatorname{mod}b$. Since we
furthermore have $r_{i}\in\left\{  0,1,\ldots,b-1\right\}  $ (because
$r_{0},r_{1},\ldots,r_{k}\in\left\{  0,1,\ldots,b-1\right\}  $), we thus
conclude that $r_{i}=q^{\prime}\%b$ (by Proposition
\ref{prop.ent.quo-rem.basics} \textbf{(c)}, applied to $q^{\prime}$, $b$ and
$r_{i}$ instead of $n$, $d$ and $c$). In view of $q^{\prime}=n//b^{i}$, we can
rewrite this as $r_{i}=\left(  n//b^{i}\right)  \%b$.

Forget that we fixed $i$. We thus have shown that $r_{i}=\left(
n//b^{i}\right)  \%b$ for each $i\in\left\{  0,1,\ldots,k\right\}  $. This
proves Claim 2.
\end{proof}

Now, both Claims 1 and 2 are proved. As explained above, this completes the
proof of Theorem \ref{thm.ent.base-b.main}.
\end{proof}
\end{fineprint}

\begin{fineprint}
The inductive proof of Claim 1 in the above proof is just a formal avatar of
the algorithm for writing a nonnegative integer $n$ in base $b$ that we
demonstrated on an example before the theorem. The formula $r_{i}=\left(
n//b^{i}\right)  \%b$ from Claim 2, on the other hand, gives an alternative
way of computing each base-$b$ digit of $n$ directly.
\end{fineprint}

\subsubsection{Congruence in terms of remainders}

Here is one more application of division with remainder: a new criterion for
congruence. Specifically, two integers $a$ and $b$ are congruent modulo a
given positive integer $d$ if and only if they leave the same remainder when
divided by $d$ (that is, satisfy $a\%d=b\%d$). In other words:

\begin{proposition}
\label{prop.ent.quo-rem.modd}Let $d$ be a positive integer. Let $a$ and $b$ be
two integers. Then, $a\equiv b\operatorname{mod}d$ if and only if $a\%d=b\%d$.
\end{proposition}

\begin{proof}
Proposition \ref{prop.ent.quo-rem.basics} \textbf{(a)} (applied to $n=a$)
yields that $a\%d\in\left\{  0,1,\ldots,d-1\right\}  $ and $a\%d\equiv
a\operatorname{mod}d$. Similarly, $b\%d\in\left\{  0,1,\ldots,d-1\right\}  $
and $b\%d\equiv b\operatorname{mod}d$.

We must prove the logical equivalence $\left(  a\equiv b\operatorname{mod}%
d\right)  \ \Longleftrightarrow\ \left(  a\%d=b\%d\right)  $. In other words,
we must prove the two implications%
\[
\left(  a\equiv b\operatorname{mod}d\right)  \ \Longrightarrow\ \left(
a\%d=b\%d\right)
\]
and%
\[
\left(  a\%d=b\%d\right)  \ \Longrightarrow\ \left(  a\equiv
b\operatorname{mod}d\right)  .
\]
Let us prove these implications separately: \medskip

\textit{Proof of }$\left(  a\equiv b\operatorname{mod}d\right)
\ \Longrightarrow\ \left(  a\%d=b\%d\right)  $\textit{:} Assume that $a\equiv
b\operatorname{mod}d$. Thus, $b\equiv a\operatorname{mod}d$ (by symmetry of
congruence -- i.e., by Proposition \ref{prop.mod.rules} \textbf{(b)}).
Combining $b\%d\equiv b\operatorname{mod}d$ with $b\equiv a\operatorname{mod}%
d$, we obtain $b\%d\equiv a\operatorname{mod}d$ (by transitivity of congruence
-- i.e., by Proposition \ref{prop.mod.rules} \textbf{(c)}).

Thus, we know that $b\%d\in\left\{  0,1,\ldots,d-1\right\}  $ and $b\%d\equiv
a\operatorname{mod}d$. Hence, Proposition \ref{prop.ent.quo-rem.basics}
\textbf{(c)} (applied to $n=a$ and $c=b\%d$) yields $b\%d=a\%d$. In other
words, $a\%d=b\%d$. Thus, we have proved the implication $\left(  a\equiv
b\operatorname{mod}d\right)  \ \Longrightarrow\ \left(  a\%d=b\%d\right)  $.
\medskip

\textit{Proof of }$\left(  a\%d=b\%d\right)  \ \Longrightarrow\ \left(
a\equiv b\operatorname{mod}d\right)  $\textit{:} Assume that $a\%d=b\%d$.
However, we know that $a\%d\equiv a\operatorname{mod}d$, so that $a\equiv
a\%d\operatorname{mod}d$ (by symmetry of congruence). In view of $a\%d=b\%d$,
we can rewrite this as $a\equiv b\%d\operatorname{mod}d$. Combining this with
$b\%d\equiv b\operatorname{mod}d$, we obtain $a\equiv b\operatorname{mod}d$
(by transitivity of congruence -- i.e., by Proposition \ref{prop.mod.rules}
\textbf{(c)}). Thus, we have proved the implication $\left(  a\%d=b\%d\right)
\ \Longrightarrow\ \left(  a\equiv b\operatorname{mod}d\right)  $. \medskip

Now, both implications are proved, so that Proposition
\ref{prop.ent.quo-rem.modd} is proved.
\end{proof}

\begin{corollary}
\label{cor.ent.parity-mod2}Let $a$ and $b$ be two integers. Then, $a\equiv
b\operatorname{mod}2$ holds if and only if the numbers $a$ and $b$ are either
both even or both odd.
\end{corollary}

\begin{proof}
$\Longrightarrow:$ Assume that $a\equiv b\operatorname{mod}2$. We must show
that the numbers $a$ and $b$ are either both even or both odd.

Proposition \ref{prop.ent.quo-rem.modd} (applied to $d=2$) shows that $a\equiv
b\operatorname{mod}2$ if and only if $a\%2=b\%2$. Thus, $a\%2=b\%2$ (since
$a\equiv b\operatorname{mod}2$). However, $a\%2\in\left\{  0,1\right\}  $ (by
Proposition \ref{prop.ent.quo-rem.basics} \textbf{(a)}, applied to $n=a$ and
$d=2$). In other words, $a\%2$ is either $0$ or $1$. If $a\%2=0$, then
$b\%2=0$ as well (since $a\%2=b\%2$), and therefore both $a$ and $b$ are even
(by Corollary \ref{cor.ent.parity-rem2} \textbf{(a)}). If $a\%2=1$, then
$b\%2=1$ as well (since $a\%2=b\%2$), and therefore both $a$ and $b$ are odd
(by Corollary \ref{cor.ent.parity-rem2} \textbf{(b)}). Other cases cannot
occur, since we know that $a\%2$ is either $0$ or $1$. Hence, in every
possible case, the numbers $a$ and $b$ are either both even or both odd. This
proves the \textquotedblleft$\Longrightarrow$\textquotedblright\ direction of
Corollary \ref{cor.ent.parity-mod2}. \medskip

$\Longleftarrow:$ Assume that the numbers $a$ and $b$ are either both even or
both odd. Thus, the numbers $a\%2$ and $b\%2$ are either both $0$ (this
happens when $a$ and $b$ are both even, by Corollary \ref{cor.ent.parity-rem2}
\textbf{(a)}) or both $1$ (this happens when $a$ and $b$ are both odd, by
Corollary \ref{cor.ent.parity-rem2} \textbf{(b)}). In either case, we have
$a\%2=b\%2$. By Proposition \ref{prop.ent.quo-rem.modd} (applied to $d=2$),
this entails $a\equiv b\operatorname{mod}2$. Hence, the \textquotedblleft%
$\Longleftarrow$\textquotedblright\ direction of Corollary
\ref{cor.ent.parity-mod2} is proved.
\end{proof}

\subsubsection{The birthday lemma}

If you have lived for exactly $n$ days, then you are $n//365$ years and
$n\%365$ days old (assuming, for simplicity, that every year has exactly $365$
days; leap years would complicate this a lot). On any \textquotedblleft
normal\textquotedblright\ day, the latter number (that is, $n\%365$) increases
by $1$ while the former number (that is, $n//365$) stays unchanged. On a
birthday, however, the latter number gets reset to $0$ while the former number
increases by $1$. This simple and intuitive observation is not specific to
$365$, and is worth stating as a proposition:

\begin{proposition}
[birthday lemma]\label{prop.ent.quo-rem.bday}Let $n\in\mathbb{Z}$, and let $d$
be a positive integer. Then: \medskip

\textbf{(a)} If $d\mid n$, then
\begin{align*}
n//d  &  =\left(  \left(  n-1\right)  //d\right)
+1\ \ \ \ \ \ \ \ \ \ \text{and}\\
n\%d  &  =0\ \ \ \ \ \ \ \ \ \ \text{and}\ \ \ \ \ \ \ \ \ \ \left(
n-1\right)  \%d=d-1.
\end{align*}

\textbf{(b)} If $d\nmid n$, then
\[
n//d=\left(  n-1\right)  //d\ \ \ \ \ \ \ \ \ \ \text{and}%
\ \ \ \ \ \ \ \ \ \ n\%d=\left(  \left(  n-1\right)  \%d\right)  +1.
\]

\end{proposition}

It should be easy to prove both parts of this lemma, but we give a proof for
the sake of completeness.

\begin{fineprint}
\begin{proof}
[Proof of Proposition \ref{prop.ent.quo-rem.bday}.]\textbf{(a)} Assume that
$d\mid n$. Thus, $n=dq$ for some $q\in\mathbb{Z}$. Consider this $q$.

Recall Definition \ref{def.ent.quo-rem}. We have $q\in\mathbb{Z}$ and
$0\in\left\{  0,1,\ldots,d-1\right\}  $ and $n=qd+0$ (since $qd+0=qd=dq=n$).
In other words, $\left(  q,0\right)  $ is a quo-rem pair of $n$ and $d$ (by
the definition of a quo-rem pair). Hence, Definition \ref{def.ent.quo-rem}
shows that $n//d=q$ and $n\%d=0$.

On the other hand, from $n=dq$, we obtain%
\begin{align*}
n-1  &  =dq-1\\
&  =\left(  q-1\right)  d+\left(  d-1\right)  \ \ \ \ \ \ \ \ \ \ \left(
\text{since }\left(  q-1\right)  d+\left(  d-1\right)  =qd-d+d-1=qd-1\right)
.
\end{align*}
Thus, we have $q-1\in\mathbb{Z}$ and $d-1\in\left\{  0,1,\ldots,d-1\right\}  $
and $n-1=\left(  q-1\right)  d+\left(  d-1\right)  $. In other words, the pair
$\left(  q-1,d-1\right)  $ is a quo-rem pair of $n-1$ and $d$ (by the
definition of a quo-rem pair). Hence, Definition \ref{def.ent.quo-rem} shows
that $\left(  n-1\right)  //d=q-1$ and $\left(  n-1\right)  \%d=d-1$.

Now, from $\left(  n-1\right)  //d=q-1$, we obtain $\left(  \left(
n-1\right)  //d\right)  +1=q=n//d$. In other words, $n//d=\left(  \left(
n-1\right)  //d\right)  +1$. Combining this with $n\%d=0$ and $\left(
n-1\right)  \%d=d-1$, we see that Proposition \ref{prop.ent.quo-rem.bday}
\textbf{(a)} has been proved. \medskip

\textbf{(b)} Assume that $d\nmid n$. Let $q=n//d$ and $r=n\%d$. Then, by the
definition of quotient and remainder, we have%
\[
q\in\mathbb{Z}\ \ \ \ \ \ \ \ \ \ \text{and}\ \ \ \ \ \ \ \ \ \ r\in\left\{
0,1,\ldots,d-1\right\}  \ \ \ \ \ \ \ \ \ \ \text{and}%
\ \ \ \ \ \ \ \ \ \ n=qd+r.
\]
If we had $r=0$, then we would have $n=qd+\underbrace{r}_{=0}=qd=dq$, which
would entail $d\mid n$; but this would contradict $d\nmid n$. Hence, we cannot
have $r=0$. In other words, $r$ is not $0$.

So $r$ is an element of the set $\left\{  0,1,\ldots,d-1\right\}  $ but is not
$0$. Therefore, $r$ is one of the remaining elements $1,2,\ldots,d-1$.
Therefore, $r-1$ is one of the elements $0,1,\ldots,d-2$. Thus, $r-1\in
\left\{  0,1,\ldots,d-1\right\}  $.

Also, from $n=qd+r$, we obtain $n-1=\left(  qd+r\right)  -1=qd+\left(
r-1\right)  $. So we know that $q\in\mathbb{Z}$ and $r-1\in\left\{
0,1,\ldots,d-1\right\}  $ and $n-1=qd+\left(  r-1\right)  $. In other words,
the pair $\left(  q,r-1\right)  $ is a quo-rem pair of $n-1$ and $d$ (by the
definition of a quo-rem pair). Hence, Definition \ref{def.ent.quo-rem} shows
that $\left(  n-1\right)  //d=q$ and $\left(  n-1\right)  \%d=r-1$.

Thus, $\left(  n-1\right)  //d=q=n//d$, so that $n//d=\left(  n-1\right)
//d$. Also, from $\left(  n-1\right)  \%d=r-1$, we obtain $\left(  \left(
n-1\right)  \%d\right)  +1=r=n\%d$, so that $n\%d=\left(  \left(  n-1\right)
\%d\right)  +1$. Thus, we have proved Proposition \ref{prop.ent.quo-rem.bday}
\textbf{(b)}.
\end{proof}
\end{fineprint}

Part of Proposition \ref{prop.ent.quo-rem.bday} can be restated using the
floor notation:

\begin{corollary}
\label{cor.ent.quo-rem.floor-bday}Let $n\in\mathbb{Z}$, and let $d$ be a
positive integer. Then: \medskip

\textbf{(a)} If $d\mid n$, then
\[
\left\lfloor \dfrac{n}{d}\right\rfloor =\left\lfloor \dfrac{n-1}%
{d}\right\rfloor +1.
\]

\textbf{(b)} If $d\nmid n$, then
\[
\left\lfloor \dfrac{n}{d}\right\rfloor =\left\lfloor \dfrac{n-1}%
{d}\right\rfloor .
\]

\end{corollary}

\begin{fineprint}
\begin{proof}
Proposition \ref{prop.ent.quo-rem-floor} yields $n//d=\left\lfloor \dfrac
{n}{d}\right\rfloor $. The same argument (applied to $n-1$ instead of $n$)
yields $\left(  n-1\right)  //d=\left\lfloor \dfrac{n-1}{d}\right\rfloor $.
\medskip

\textbf{(a)} Assume that $d\mid n$. Then, Proposition
\ref{prop.ent.quo-rem.bday} \textbf{(a)} yields $n//d=\left(  \left(
n-1\right)  //d\right)  +1$. In view of $n//d=\left\lfloor \dfrac{n}%
{d}\right\rfloor $ and $\left(  n-1\right)  //d=\left\lfloor \dfrac{n-1}%
{d}\right\rfloor $, we can rewrite this as $\left\lfloor \dfrac{n}%
{d}\right\rfloor =\left\lfloor \dfrac{n-1}{d}\right\rfloor +1$. This proves
Corollary \ref{cor.ent.quo-rem.floor-bday} \textbf{(a)}. \medskip

\textbf{(b)} Assume that $d\nmid n$. Then, Proposition
\ref{prop.ent.quo-rem.bday} \textbf{(b)} yields $n//d=\left(  n-1\right)
//d$. In view of $n//d=\left\lfloor \dfrac{n}{d}\right\rfloor $ and $\left(
n-1\right)  //d=\left\lfloor \dfrac{n-1}{d}\right\rfloor $, we can rewrite
this as $\left\lfloor \dfrac{n}{d}\right\rfloor =\left\lfloor \dfrac{n-1}%
{d}\right\rfloor $. This proves Corollary \ref{cor.ent.quo-rem.floor-bday}
\textbf{(b)}.
\end{proof}
\end{fineprint}

\subsection{\label{sec.ent.gcd}Greatest common divisors}

\subsubsection{Definition}

The following definition plays a crucial role in number theory, particularly
in the study of prime numbers that will be the topic of Section
\ref{sec.ent.pri}.

\begin{definition}
\label{def.ent.gcd}Let $a$ and $b$ be two integers. \medskip

\textbf{(a)} The \textbf{common divisors} of $a$ and $b$ are the integers that
divide $a$ and simultaneously divide $b$. \medskip

\textbf{(b)} The \textbf{greatest common divisor} of $a$ and $b$ is the
largest among the common divisors of $a$ and $b$, unless $a=b=0$. In the case
$a=b=0$, it is defined to be $0$ instead.

We denote the greatest common divisor of $a$ and $b$ as $\gcd\left(
a,b\right)  $, and we refer to it as the \textbf{gcd} of $a$ and $b$.
\end{definition}

We will soon see that this greatest common divisor is well-defined (see Remark
\ref{rmk.ent.gcd.wd} below). But first, some examples:

\begin{itemize}
\item What is $\gcd\left(  4,6\right)  $ ?

The divisors of $4$ are $-4,-2,-1,1,2,4$.

The divisors of $6$ are $-6,-3,-2,-1,1,2,3,6$.

Thus, the common divisors of $4$ and $6$ are $-2,-1,1,2$.

So the greatest common divisor of $4$ and $6$ is $2$. That is, $\gcd\left(
4,6\right)  =2$.

\item What is $\gcd\left(  0,5\right)  $ ?

The divisors of $0$ are all the integers (you cannot list them all).

The divisors of $5$ are $-5,-1,1,5$.

Thus, the common divisors of $0$ and $5$ are just the divisors of $5$, which
are $-5,-1,1,5$.

So the gcd is $5$. That is, $\gcd\left(  0,5\right)  =5$.

\item What is $\gcd\left(  0,0\right)  $ ?

The common divisors of $0$ and $0$ are all the integers, so there is no
greatest one among them, but we have defined $\gcd\left(  0,0\right)  $ to be
$0$. (This is the reason why we had to make an exception for the $a=b=0$ case
in Definition \ref{def.ent.gcd} \textbf{(b)}.)
\end{itemize}

Let us now convince ourselves that $\gcd\left(  a,b\right)  $ is well-defined:

\begin{remark}
\label{rmk.ent.gcd.wd}Let $a,b\in\mathbb{Z}$. We want to show that
$\gcd\left(  a,b\right)  $ is well-defined in Definition \ref{def.ent.gcd}
\textbf{(b)}.

If $a=b=0$, then this is clear, since we defined this gcd to be $0$.

Consider the remaining case -- i.e., the case when $a\neq0$ or $b\neq0$ (or both).

For instance, let us assume that $a\neq0$. Then, the divisors $d$ of $a$ all
satisfy $\left\vert d\right\vert \leq\left\vert a\right\vert $ (since
Proposition \ref{prop.div.abs} \textbf{(b)} shows that they satisfy
$\operatorname*{abs}d\leq\operatorname*{abs}a$, which in our present notations
means $\left\vert d\right\vert \leq\left\vert a\right\vert $). In other words,
all these divisors are integers in the interval $\left[  -\left\vert
a\right\vert ,\ \left\vert a\right\vert \right]  $. Hence, there are finitely
many of them. Therefore, there are finitely many common divisors of $a$ and
$b$ (since any common divisor of $a$ and $b$ is a divisor of $a$). On the
other hand, there is at least one common divisor of $a$ and $b$ (namely, $1$).
Therefore, the set of all common divisors of $a$ and $b$ is nonempty and
finite, and thus has a maximum element. In other words, there is a (literally)
largest among the common divisors of $a$ and $b$. This shows that $\gcd\left(
a,b\right)  $ is well-defined when $a\neq0$.

An analogous argument leads to the same conclusion when $b\neq0$. Thus, we
have shown that $\gcd\left(  a,b\right)  $ is always well-defined.
\end{remark}

This argument also gives us a slow and stupid algorithm to compute
$\gcd\left(  a,b\right)  $ when $a\neq0$: We just go through all integers in
the interval $\left[  -\left\vert a\right\vert ,\ \left\vert a\right\vert
\right]  $, and check which of them are common divisors of $a$ and $b$. But
there is a much faster algorithm.

\subsubsection{Basic properties}

To find this algorithm, we first collect some basic properties of gcds:

\begin{proposition}
\label{prop.ent.gcd.basics}\textbf{(a)} We have $\gcd\left(  a,b\right)
\in\mathbb{N}$ for any $a,b\in\mathbb{Z}$. \medskip

\textbf{(b)} We have $\gcd\left(  a,0\right)  =\gcd\left(  0,a\right)
=\left\vert a\right\vert $ for any $a\in\mathbb{Z}$. \medskip

\textbf{(c)} We have $\gcd\left(  a,b\right)  =\gcd\left(  b,a\right)  $ for
any $a,b\in\mathbb{Z}$. \medskip

\textbf{(d)} If $a,b,c\in\mathbb{Z}$ satisfy $b\equiv c\operatorname{mod}a$,
then $\gcd\left(  a,b\right)  =\gcd\left(  a,c\right)  $. \medskip

\textbf{(e)} We have $\gcd\left(  a,b\right)  =\gcd\left(  a,\ ua+b\right)  $
for any $a,b,u\in\mathbb{Z}$. \medskip

\textbf{(f)} We have $\gcd\left(  a,b\right)  =\gcd\left(  a,\ b\%a\right)  $
for any positive integer $a$ and any $b\in\mathbb{Z}$. \medskip

\textbf{(g)} We have $\gcd\left(  a,b\right)  \mid a$ and $\gcd\left(
a,b\right)  \mid b$ for any $a,b\in\mathbb{Z}$. \medskip

\textbf{(h)} We have $\gcd\left(  -a,b\right)  =\gcd\left(  a,b\right)  $ and
$\gcd\left(  a,-b\right)  =\gcd\left(  a,b\right)  $ for any $a,b\in
\mathbb{Z}$. \medskip

\textbf{(i)} If $a,b\in\mathbb{Z}$ satisfy $a\mid b$, then $\gcd\left(
a,b\right)  =\left\vert a\right\vert $.
\end{proposition}

\begin{proof}
\textbf{(a)} Let $a,b\in\mathbb{Z}$. We must prove that $\gcd\left(
a,b\right)  \in\mathbb{N}$.

If $a=b=0$, then this follows from $\gcd\left(  0,0\right)  =0\in\mathbb{N}$.

Thus, let us assume that not both of $a$ and $b$ are $0$. Then, $\gcd\left(
a,b\right)  $ is literally the greatest common divisor of $a$ and $b$. If
$\gcd\left(  a,b\right)  $ was negative, then $-\gcd\left(  a,b\right)  $
would be an even greater common divisor of $a$ and $b$ (since $-\gcd\left(
a,b\right)  $ divides whatever $\gcd\left(  a,b\right)  $ divides, but the
negativity of $\gcd\left(  a,b\right)  $ implies $-\gcd\left(  a,b\right)
>\gcd\left(  a,b\right)  $), which would contradict the previous sentence.
Hence, $\gcd\left(  a,b\right)  $ cannot be negative. Thus, $\gcd\left(
a,b\right)  \in\mathbb{N}$. This proves Proposition \ref{prop.ent.gcd.basics}
\textbf{(a)}. \medskip

\textbf{(b)} Let $a\in\mathbb{Z}$. Every integer is a divisor of $0$. Thus,
the common divisors of $a$ and $0$ are just the divisors of $a$. However, the
largest divisor of $a$ is $\left\vert a\right\vert $ (unless $a=0$, which case
can be easily handled separately)\footnote{This fact is a consequence of
Proposition \ref{prop.div.abs} \textbf{(b)} (recalling that $\left\vert
a\right\vert $ was called $\operatorname*{abs}a$ back in that proposition).}.
Hence, the greatest common divisor of $a$ and $0$ is $\left\vert a\right\vert
$. In other words, we have $\gcd\left(  a,0\right)  =\left\vert a\right\vert
$. Similarly, we can see that $\gcd\left(  0,a\right)  =\left\vert
a\right\vert $. Thus, Proposition \ref{prop.ent.gcd.basics} \textbf{(b)} is
proved. \medskip

\textbf{(c)} Proposition \ref{prop.ent.gcd.basics} \textbf{(c)} follows from
observing that $a$ and $b$ play equal roles in Definition \ref{def.ent.gcd}.
\medskip

\textbf{(d)} Let $a,b,c\in\mathbb{Z}$ satisfy $b\equiv c\operatorname{mod}a$.
We must prove that $\gcd\left(  a,b\right)  =\gcd\left(  a,c\right)  $.

If $a=0$, then this is clearly true (because in this case, $b\equiv
c\operatorname{mod}a$ becomes $b\equiv c\operatorname{mod}0$, which entails
$b=c$).

It thus remains to consider the case $a\neq0$ only. In this case, $\gcd\left(
a,b\right)  $ is literally the greatest common divisor of $a$ and $b$, whereas
$\gcd\left(  a,c\right)  $ is literally the greatest common divisor of $a$ and
$c$. Hence, in order to prove that these two gcds are equal, it will suffice
to show that the common divisors of $a$ and $b$ are precisely the common
divisors of $a$ and $c$. To do this, in turn, it suffices to prove the
following two claims:

\begin{statement}
\textit{Claim 1:} Each common divisor of $a$ and $b$ is a common divisor of
$a$ and $c$.
\end{statement}

\begin{statement}
\textit{Claim 2:} Each common divisor of $a$ and $c$ is a common divisor of
$a$ and $b$.
\end{statement}

Before we prove these two claims, let us recall that $b\equiv
c\operatorname{mod}a$; in other words, $c\equiv b\operatorname{mod}a$ (by the
symmetry of congruence). Hence, the numbers $b$ and $c$ play equal roles in
our setting. Thus, Claims 1 and 2 are analogous, so that any proof of one of
the two will also prove the other (once the roles of $b$ and $c$ are switched).

\begin{proof}
[Proof of Claim 1.]Let $d$ be a common divisor of $a$ and $b$. Thus, $d\mid a$
and $d\mid b$ (by the definition of a common divisor). In other words, we have
$a=dx$ and $b=dy$ for some integers $x$ and $y$. Consider these $x$ and $y$.

But $b\equiv c\operatorname{mod}a$. In other words, $a\mid b-c$. Hence, $d\mid
a\mid b-c$ (by the transitivity of divisibility). In other words, $b-c=dz$ for
some integer $z$. Consider this $z$.

Now, $b-\left(  b-c\right)  =c$, so that
\[
c=\underbrace{b}_{=dy}-\underbrace{\left(  b-c\right)  }_{=dz}%
=dy-dz=d\underbrace{\left(  y-z\right)  }_{\text{an integer}}.
\]
Therefore, $d\mid c$. From $d\mid a$ and $d\mid c$, we conclude that $d$ is a
common divisor of $a$ and $c$.

So we have shown that if $d$ is a common divisor of $a$ and $b$, then $d$ is a
common divisor of $a$ and $c$. In other words, each common divisor of $a$ and
$b$ is a common divisor of $a$ and $c$. This proves Claim 1.
\end{proof}

\begin{proof}
[Proof of Claim 2.]As we said, we can obtain a proof of Claim 2 by switching
the roles of $b$ and $c$ in the above proof of Claim 1 (because we have
$c\equiv b\operatorname{mod}a$).
\end{proof}

Combining Claim 1 with Claim 2, we see that the common divisors of $a$ and $b$
are precisely the common divisors of $a$ and $c$. Therefore, the greatest
common divisor of $a$ and $b$ equals the greatest common divisor of $a$ and
$c$. In other words, $\gcd\left(  a,b\right)  =\gcd\left(  a,c\right)  $. This
proves Proposition \ref{prop.ent.gcd.basics} \textbf{(d)}. \medskip

\textbf{(e)} Proposition \ref{prop.ent.gcd.basics} \textbf{(e)} follows from
Proposition \ref{prop.ent.gcd.basics} \textbf{(d)} (applied to $c=ua+b$),
since $b\equiv ua+b\operatorname{mod}a$ (because $b-\left(  ua+b\right)  =-ua$
is divisible by $a$). \medskip

\textbf{(f)} Proposition \ref{prop.ent.gcd.basics} \textbf{(f)} follows from
Proposition \ref{prop.ent.gcd.basics} \textbf{(d)} (applied to $c=b\%a$),
since $b\equiv b\%a\operatorname{mod}a$ (because Proposition
\ref{prop.ent.quo-rem.basics} \textbf{(a)} yields $b\%a\equiv
b\operatorname{mod}a$). \medskip

\textbf{(g)} is obvious when $a=b=0$ (since $0\mid0$), and otherwise follows
from the definition of $\gcd\left(  a,b\right)  $. \medskip

\textbf{(h)} The divisors of $a$ are precisely the divisors of $-a$. The
divisors of $b$ are precisely the divisors of $-b$. Thus, the common divisors
of $a$ and $b$ remain unchanged if we replace $a$ by $-a$ or replace $b$ by
$-b$. Therefore, Proposition \ref{prop.ent.gcd.basics} \textbf{(h)} follows
from the definition of the gcd (except in the case when $a=b=0$, but this case
is again obvious). \medskip

\textbf{(i)} Let $a,b\in\mathbb{Z}$ satisfy $a\mid b$. Then, $b\equiv
0\operatorname{mod}a$. Hence, Proposition \ref{prop.ent.gcd.basics}
\textbf{(d)} (applied to $c=0$) yields $\gcd\left(  a,b\right)  =\gcd\left(
a,0\right)  =\left\vert a\right\vert $ (by Proposition
\ref{prop.ent.gcd.basics} \textbf{(b)}). This proves Proposition
\ref{prop.ent.gcd.basics} \textbf{(i)}.
\end{proof}

\begin{corollary}
[Euclidean recursion for the gcd]\label{cor.ent.gcd.euclid}Let $a\in
\mathbb{Z}$, and let $b$ be a positive integer. Then,%
\[
\gcd\left(  a,b\right)  =\gcd\left(  b,\ a\%b\right)  .
\]

\end{corollary}

\begin{proof}
Proposition \ref{prop.ent.gcd.basics} \textbf{(c)} yields%
\[
\gcd\left(  a,b\right)  =\gcd\left(  b,a\right)  =\gcd\left(  b,\ a\%b\right)
\]
(by Proposition \ref{prop.ent.gcd.basics} \textbf{(f)}, applied to $b$ and $a$
instead of $a$ and $b$). This proves Corollary \ref{cor.ent.gcd.euclid}.
\end{proof}

\subsubsection{The Euclidean algorithm}

By applying Corollary \ref{cor.ent.gcd.euclid} repeatedly, we can compute gcds
rather quickly: For example,%
\begin{align*}
\gcd\left(  93,\ 18\right)   &  =\gcd\left(  18,\ \underbrace{93\%18}%
_{=3}\right)  \ \ \ \ \ \ \ \ \ \ \left(  \text{by Corollary
\ref{cor.ent.gcd.euclid}}\right) \\
&  =\gcd\left(  18,\ 3\right) \\
&  =\gcd\left(  3,\ \underbrace{18\%3}_{=0}\right)
\ \ \ \ \ \ \ \ \ \ \left(  \text{by Corollary \ref{cor.ent.gcd.euclid}%
}\right) \\
&  =\gcd\left(  3,\ 0\right)  =\left\vert 3\right\vert
\ \ \ \ \ \ \ \ \ \ \left(  \text{by Proposition \ref{prop.ent.gcd.basics}
\textbf{(b)}}\right) \\
&  =3
\end{align*}
and%
\begin{align*}
\gcd\left(  1145,\ 739\right)   &  =\gcd\left(  739,\ \underbrace{1145\%739}%
_{=406}\right)  \ \ \ \ \ \ \ \ \ \ \left(  \text{by Corollary
\ref{cor.ent.gcd.euclid}}\right) \\
&  =\gcd\left(  739,\ 406\right) \\
&  =\gcd\left(  406,\ \underbrace{739\%406}_{=333}\right)
\ \ \ \ \ \ \ \ \ \ \left(  \text{by Corollary \ref{cor.ent.gcd.euclid}%
}\right) \\
&  =\gcd\left(  406,\ 333\right) \\
&  =\gcd\left(  333,\ \underbrace{406\%333}_{=73}\right)
\ \ \ \ \ \ \ \ \ \ \left(  \text{by Corollary \ref{cor.ent.gcd.euclid}%
}\right) \\
&  =\gcd\left(  333,\ 73\right) \\
&  =\gcd\left(  73,\ 333\%73\right)  \ \ \ \ \ \ \ \ \ \ \left(  \text{by
Corollary \ref{cor.ent.gcd.euclid}}\right) \\
&  =\gcd\left(  73,\ 41\right) \\
&  =\gcd\left(  41,\ 73\%41\right)  \ \ \ \ \ \ \ \ \ \ \left(  \text{by
Corollary \ref{cor.ent.gcd.euclid}}\right) \\
&  =\gcd\left(  41,\ 32\right) \\
&  =\gcd\left(  32,\ 41\%32\right)  \ \ \ \ \ \ \ \ \ \ \left(  \text{by
Corollary \ref{cor.ent.gcd.euclid}}\right) \\
&  =\gcd\left(  32,\ 9\right) \\
&  =\gcd\left(  9,\ 32\%9\right)  \ \ \ \ \ \ \ \ \ \ \left(  \text{by
Corollary \ref{cor.ent.gcd.euclid}}\right) \\
&  =\gcd\left(  9,\ 5\right) \\
&  =\gcd\left(  5,\ 9\%5\right)  \ \ \ \ \ \ \ \ \ \ \left(  \text{by
Corollary \ref{cor.ent.gcd.euclid}}\right) \\
&  =\gcd\left(  5,\ 4\right) \\
&  =\gcd\left(  4,\ 5\%4\right)  \ \ \ \ \ \ \ \ \ \ \left(  \text{by
Corollary \ref{cor.ent.gcd.euclid}}\right) \\
&  =\gcd\left(  4,\ 1\right) \\
&  =\gcd\left(  1,\ 4\%1\right)  \ \ \ \ \ \ \ \ \ \ \left(  \text{by
Corollary \ref{cor.ent.gcd.euclid}}\right) \\
&  =\gcd\left(  1,\ 0\right)  =\left\vert 1\right\vert
\ \ \ \ \ \ \ \ \ \ \left(  \text{by Proposition \ref{prop.ent.gcd.basics}
\textbf{(b)}}\right) \\
&  =1.
\end{align*}

These two computations are instances of a general algorithm for computing
$\gcd\left(  a,b\right)  $ for any two numbers $a\in\mathbb{Z}$ and
$b\in\mathbb{N}$. This algorithm proceeds as follows:

\begin{itemize}
\item If $b=0$, then the gcd is $\left\vert a\right\vert $.

\item If $b>0$, then we replace $a$ and $b$ by $b$ and $a\%b$ and recurse
(i.e., we apply the method again to $b$ and $a\%b$ instead of $a$ and $b$).
\end{itemize}

In Python code\footnote{I am using the Python programming language because of
its ease of use and abundance of inbuilt fundamental mathematical tools. All
the algorithms can be implemented in any other language as well, but the code
looks best in Python.}, this algorithm looks as follows:\medskip\Needspace{6pc}

\texttt{\qquad def gcd(a, b): \# for b nonnegative }

\texttt{\qquad\qquad if b == 0:}

\texttt{\qquad\qquad\qquad return abs(a) \# This is the absolute value of a.}

\texttt{\qquad\qquad return gcd(b, a\%b)} \medskip

This algorithm is called the \textbf{Euclidean algorithm}. Let us convince
ourselves that it really terminates (rather than getting stuck in an endless loop):

\begin{proposition}
\label{prop.ent.gcd.euclid-terminates}Let $a\in\mathbb{Z}$ and $b\in
\mathbb{N}$. Then, the Euclidean algorithm terminates after at most $b$ steps.
(Here, we count each time that the algorithm replaces $a$ and $b$ by $b$ and
$a\%b$ as a \textquotedblleft step\textquotedblright.)
\end{proposition}

\begin{proof}
In each step of the Euclidean algorithm, the second argument $b$ gets replaced
by $a\%b$. This has the consequence that $b$ decreases by at least $1$ (since
the definition of a remainder yields $a\%b\in\left\{  0,1,\ldots,b-1\right\}
$ and thus $a\%b\leq b-1$). But $b$ remains nonnegative throughout the
algorithm. Thus, $b$ cannot decrease (by at least $1$) more than $b_{0}$ times
in succession, where $b_{0}$ is the original value of $b$ (as it was fed into
the algorithm). Hence, the algorithm cannot have more than $b_{0}$ steps. In
other words, the algorithm must terminate after at most $b_{0}$ steps. This
proves Proposition \ref{prop.ent.gcd.euclid-terminates} (since $b_{0}$ is
precisely the original value of $b$).
\end{proof}

Proposition \ref{prop.ent.gcd.euclid-terminates} greatly overestimates the
actual time that the Euclidean algorithm needs to terminate: In truth, it
terminates after at most $\log_{2}\left(  ab\right)  +2$ steps (if $a$ and $b$
are positive)\footnote{\textit{Hints to the proof.} Recall that each step of
the algorithm replaces the numbers $a$ and $b$ by $b$ and $a\%b$. Since
$b>a\%b$ (because $a\%b\in\left\{  0,1,\ldots,b-1\right\}  $ entails
$a\%b<b$), this yields that after each step of the algorithm, the
\textquotedblleft current\textquotedblright\ numbers $a$ and $b$ satisfy
$a>b$.
\par
Now, consider the product $ab$ of the two numbers $a$ and $b$. We claim that
each step of the algorithm, except perhaps the first one, decreases this
number by a factor of at least $2$. In order to see this, you need to show
that $b\left(  a\%b\right)  \leq\dfrac{ab}{2}$ whenever $a>b$. But this
follows from $a\%b\leq\dfrac{a}{2}$, which in turn follows easily from $a>b$
(why?).
\par
Now you know that the product $ab$ decreases by a factor of at least $2$ at
each step of the algorithm except for the first one. In other words, its
binary logarithm $\log_{2}\left(  ab\right)  $ decreases by at least $1$ at
each step of the algorithm except for the first one (here, we use the
convention $\log_{2}0:=-\infty<0$). At the first step, it also decreases or
stays unchanged. From this, it follows easily that the algorithm cannot have
more than $\log_{2}\left(  ab\right)  +1$ steps until it reaches a situation
in which $\log_{2}\left(  ab\right)  \leq0$. But in such a situation, we must
have $a=b=1$ or $a=0$ or $b=0$, and it will only take one more step to reach
the end of the algorithm.}, which is usually much fewer than $b$. Some
variants of the Euclidean algorithm get to the goal even faster. This
speediness is part of the reason why the Euclidean algorithm (and greatest
common divisors) is so useful in practical applications of number theory.
\medskip

The Euclidean algorithm can be easily adapted to arbitrary $b\in\mathbb{Z}$
instead of just $b\in\mathbb{N}$ (by adding a first step in which we replace
$b$ by $-b$ if $b$ is negative):\medskip\Needspace{6pc}

\texttt{\qquad def gcd(a, b): \# for b arbitrary }

\texttt{\qquad\qquad if b
$<$
0:}

\texttt{\qquad\qquad\qquad return gcd(a, -b) \# replace b by -b.}

\texttt{\qquad\qquad if b == 0:}

\texttt{\qquad\qquad\qquad return abs(a) \# This is the absolute value of a.}

\texttt{\qquad\qquad return gcd(b, a\%b)}

\subsubsection{\label{subsec.ent.gcd.bezout}Bezout's theorem and the extended
Euclidean algorithm}

The Euclidean algorithm can be adapted so that it doesn't only compute
$\gcd\left(  a,b\right)  $, but also expresses $\gcd\left(  a,b\right)  $ as
an \textquotedblleft integer linear combination\textquotedblright\ of $a$ and
$b$ (that is, as a multiple of $a$ plus a multiple of $b$). This allows us to
prove the following theorem:

\begin{theorem}
[Bezout's theorem for integers]\label{thm.ent.gcd.bezout}Let $a$ and $b$ be
two integers. Then, there exist two integers $x$ and $y$ such that
\[
\gcd\left(  a,b\right)  =xa+yb.
\]

\end{theorem}

We will soon prove this theorem. First, we introduce a notation and give a few examples:

\begin{definition}
\label{def.ent.gcd.bezout-pair}Let $a$ and $b$ be two integers. Then, a
\textbf{Bezout pair} for $\left(  a,b\right)  $ means a pair $\left(
x,y\right)  $ of two integers satisfying $\gcd\left(  a,b\right)  =xa+yb$.
\end{definition}

For instance, a Bezout pair for $\left(  4,7\right)  $ is a pair $\left(
x,y\right)  $ of integers satisfying $\gcd\left(  4,7\right)  =x\cdot
4+y\cdot7$. In view of $\gcd\left(  4,7\right)  =1$, this latter equation
simplifies to $1=4x+7y$. So a Bezout pair for $\left(  4,7\right)  $ is a
solution to this equation $1=4x+7y$ in \textbf{integers} $x$ and $y$. This is
similar to the coin problem from Subsection \ref{subsec.ind.strind.coin-35},
in the sense that you can think of such a Bezout pair $\left(  x,y\right)  $
as a way to pay $1$ cent with $x$ many $4$-cent coins and $y$ many $7$-cent
coins, assuming that you are allowed to get change (because $x$ and $y$ are
allowed to be negative). Without change, of course, you could not pay $1$ cent
using $4$-cent coins and $7$-cent coins. But with change, it works: You pay
two $4$-cent coins and get one $7$-cent coin in return, and thus end up paying
$2\cdot4+\left(  -1\right)  \cdot7=1$ cent, which is what you wanted. In other
words, the pair $\left(  x,y\right)  =\left(  2,-1\right)  $ satisfies
$1=4x+7y$. In other words, $\left(  2,-1\right)  $ is a Bezout pair for
$\left(  4,7\right)  $. There are also other Bezout pairs for $\left(
4,7\right)  $, for example $\left(  -5,3\right)  $ (since $4\left(  -5\right)
+7\cdot3=1$). So a Bezout pair is usually not unique.

So Bezout's theorem can be restated as follows: For any two integers $a$ and
$b$, you can pay $\gcd\left(  a,b\right)  $ cents with $a$-cent coins and
$b$-cent coins, if you can get change\footnote{more precisely: if you can get
change in $a$-cent coins and $b$-cent coins (and there are infinitely many
coins of either denomination available)}. What amounts can be paid
\textbf{without} change is a more complicated story, and we will return to
this in Section \ref{sec.ent.sylv}.

Here is another example: A Bezout pair for $\left(  6,16\right)  $ is $\left(
3,-1\right)  $, since $\gcd\left(  6,16\right)  =2=6x+16y$ for $\left(
x,y\right)  =\left(  3,-1\right)  $.

\bigskip

So Bezout's theorem (Theorem \ref{thm.ent.gcd.bezout}) is saying that for any
two integers $a,b\in\mathbb{Z}$, there exists a Bezout pair for $\left(
a,b\right)  $.

How can we prove this theorem? Induction (particularly strong induction)
appears to be a reasonable method. Unfortunately, induction can only be used
to prove a statement about elements of a set of the form $\left\{
k,k+1,k+2,\ldots\right\}  $ for a given integer $k$ (that is, a statement
about integers from a given lower bound onwards). To put it differently,
induction can only prove a statement that \textquotedblleft starts
somewhere\textquotedblright\ (even if it is presented as a strong induction
with no base case). Meanwhile, in Bezout's theorem, both $a$ and $b$ are just
arbitrary integers, so they can be arbitrarily low.

This hurdle can be surmounted: While we cannot prove Bezout's theorem by
induction directly, we can first restrict it to the case when $b\in\mathbb{N}%
$, and prove this restriction by induction. In other words, we shall use
induction to prove the following particular case of Bezout's theorem:

\begin{lemma}
[restricted Bezout's theorem]\label{lem.ent.gcd.bezout-ZN}Let $a\in\mathbb{Z}$
and $b\in\mathbb{N}$. Then, there exists a Bezout pair for $\left(
a,b\right)  $.
\end{lemma}

Once this lemma is proved, we will quickly deduce Bezout's theorem in full
generality from it. So let us prove this lemma.

\begin{proof}
[Proof of Lemma \ref{lem.ent.gcd.bezout-ZN}.]We shall use strong induction on
$b$. Here, we do not consider $a$ to be fixed. Thus, the statement that we
will be proving for all $b\in\mathbb{N}$ is%
\[
P\left(  b\right)  :=\left(  \text{for each }a\in\mathbb{Z}\text{, there
exists a Bezout pair for }\left(  a,b\right)  \right)  .
\]
Our goal is to prove this statement $P\left(  b\right)  $ for all
$b\in\mathbb{N}$. We shall do this by strong induction on $b$:

\textit{Base case:} Let us prove the statement $P\left(  0\right)  $. Indeed,
for each $a\in\mathbb{Z}$, let us set
\[
\operatorname*{sign}a:=%
\begin{cases}
1, & \text{if }a>0;\\
0, & \text{if }a=0;\\
-1, & \text{if }a<0.
\end{cases}
\]
Then, for each $a\in\mathbb{Z}$, the pair $\left(  \operatorname*{sign}%
a,0\right)  $ is a Bezout pair for $\left(  a,0\right)  $, since%
\begin{align*}
\gcd\left(  a,0\right)   &  =\left\vert a\right\vert
\ \ \ \ \ \ \ \ \ \ \left(  \text{by Proposition \ref{prop.ent.gcd.basics}
\textbf{(b)}}\right) \\
&  =\left(  \operatorname*{sign}a\right)  \cdot a\ \ \ \ \ \ \ \ \ \ \left(
\begin{array}
[c]{c}%
\text{this is a general fact that holds for any real}\\
\text{number }a\text{, and can be easily verified by}\\
\text{checking the cases }a>0\text{, }a=0\text{ and }a<0
\end{array}
\right) \\
&  =\left(  \operatorname*{sign}a\right)  \cdot a+0\cdot0.
\end{align*}
Hence, for each $a\in\mathbb{Z}$, there exists a Bezout pair for $\left(
a,0\right)  $. In other words, the statement $P\left(  0\right)  $ holds.

\textit{Induction step:} Fix a positive integer $b$. We must prove the
implication%
\[
\left(  P\left(  0\right)  \ \text{AND }P\left(  1\right)  \text{ AND
}P\left(  2\right)  \text{ AND }\cdots\text{ AND }P\left(  b-1\right)
\right)  \Longrightarrow P\left(  b\right)  .
\]
Thus, we assume (as the induction hypothesis) that $P\left(  0\right)  \ $AND
$P\left(  1\right)  $ AND $P\left(  2\right)  $ AND $\cdots$ AND $P\left(
b-1\right)  $ holds. In other words, we assume that the $b$ statements
$P\left(  0\right)  ,\ P\left(  1\right)  ,\ P\left(  2\right)  ,\ \ldots
,\ P\left(  b-1\right)  $ all hold. In other words, we assume that%
\begin{align*}
&  \left(  \text{for each }a\in\mathbb{Z}\text{, there exists a Bezout pair
for }\left(  a,0\right)  \right)  \text{ and}\\
&  \left(  \text{for each }a\in\mathbb{Z}\text{, there exists a Bezout pair
for }\left(  a,1\right)  \right)  \text{ and}\\
&  \left(  \text{for each }a\in\mathbb{Z}\text{, there exists a Bezout pair
for }\left(  a,2\right)  \right)  \text{ and}\\
&  \cdots\ \ \ \ \ \ \ \ \ \ \text{and}\\
&  \left(  \text{for each }a\in\mathbb{Z}\text{, there exists a Bezout pair
for }\left(  a,b-1\right)  \right)  .
\end{align*}
In other words, we assume that for each $a\in\mathbb{Z}$ and each
$d\in\left\{  0,1,\ldots,b-1\right\}  $, there exists a Bezout pair for
$\left(  a,d\right)  $. Renaming $a$ as $c$ here, we can restate this as
follows: We assume that for each $c\in\mathbb{Z}$ and each $d\in\left\{
0,1,\ldots,b-1\right\}  $, there exists a Bezout pair for $\left(  c,d\right)
$. So this is our induction hypothesis (brought to its most convenient form).

Our goal is now to prove $P\left(  b\right)  $. In other words, we must prove
that for each $a\in\mathbb{Z}$, there exists a Bezout pair for $\left(
a,b\right)  $.

So we fix an $a\in\mathbb{Z}$, and we set out to find a Bezout pair for
$\left(  a,b\right)  $.

The Euclidean recursion (Corollary \ref{cor.ent.gcd.euclid}) yields%
\begin{equation}
\gcd\left(  a,b\right)  =\gcd\left(  b,\ a\%b\right)  .
\label{pf.lem.ent.gcd.bezout-ZN.4}%
\end{equation}
However, $a\%b\in\left\{  0,1,\ldots,b-1\right\}  $ (by Proposition
\ref{prop.ent.quo-rem.basics} \textbf{(a)}, applied to $n=a$ and $d=b$).

Recall our induction hypothesis, which says that for each $c\in\mathbb{Z}$ and
each $d\in\left\{  0,1,\ldots,b-1\right\}  $, there exists a Bezout pair for
$\left(  c,d\right)  $. We can apply this to $c=b$ and $d=a\%b$ (because
$b\in\mathbb{Z}$ and $a\%b\in\left\{  0,1,\ldots,b-1\right\}  $), and thus
conclude that there exists a Bezout pair for $\left(  b,\ a\%b\right)  $. Let
us denote this Bezout pair by $\left(  u,v\right)  $. Thus, by the definition
of a Bezout pair, $u$ and $v$ are integers and satisfy
\begin{equation}
\gcd\left(  b,\ a\%b\right)  =ub+v\left(  a\%b\right)  .
\label{pf.lem.ent.gcd.bezout-ZN.5}%
\end{equation}

However, Proposition \ref{prop.ent.quo-rem.basics} \textbf{(d)} (applied to
$n=a$ and $d=b$) yields%
\[
a=\left(  a//b\right)  b+\left(  a\%b\right)  .
\]
Solving this for $a\%b$, we obtain%
\begin{equation}
a\%b=a-\left(  a//b\right)  b. \label{pf.lem.ent.gcd.bezout-ZN.7}%
\end{equation}

Now, (\ref{pf.lem.ent.gcd.bezout-ZN.4}) becomes%
\begin{align*}
\gcd\left(  a,b\right)   &  =\gcd\left(  b,\ a\%b\right)
=ub+v\underbrace{\left(  a\%b\right)  }_{\substack{=a-\left(  a//b\right)
b\\\text{(by (\ref{pf.lem.ent.gcd.bezout-ZN.7}))}}}\ \ \ \ \ \ \ \ \ \ \left(
\text{by (\ref{pf.lem.ent.gcd.bezout-ZN.5})}\right) \\
&  =ub+v\left(  a-\left(  a//b\right)  b\right) \\
&  =ub+va-v\left(  a//b\right)  b\\
&  =\underbrace{v}_{\text{an integer}}a+\underbrace{\left(  u-v\left(
a//b\right)  \right)  }_{\text{an integer}}b.
\end{align*}
Thus, we have written $\gcd\left(  a,b\right)  $ as a multiple of $a$ plus a
multiple of $b$. More specifically, the pair%
\[
\left(  v,\ u-v\left(  a//b\right)  \right)
\]
is a Bezout pair for $\left(  a,b\right)  $. And so we conclude that there
exists a Bezout pair for $\left(  a,b\right)  $ (because we just found one).
This proves the statement $P\left(  b\right)  $ for our $b$, and thus
completes the induction step.

Hence, by induction, we have shown that $P\left(  b\right)  $ holds for all
$b\in\mathbb{N}$. But this is saying precisely that there exists a Bezout pair
for $\left(  a,b\right)  $ whenever $a\in\mathbb{Z}$ and $b\in\mathbb{N}$.
Thus, Lemma \ref{lem.ent.gcd.bezout-ZN} is proved.
\end{proof}

This inductive proof contains a recursive algorithm for finding a Bezout pair
for $\left(  a,b\right)  $ whenever $a\in\mathbb{Z}$ and $b\in\mathbb{N}$.
Written in the Python programming language, this algorithm looks as
follows:\footnote{Here, \texttt{sign(a)} is what was called
$\operatorname*{sign}a$ in the above proof. In Python, this can be defined as
follows:\medskip
\par
\texttt{\qquad def sign(a):}
\par
\texttt{\qquad\qquad if a
$<$
0:}
\par
\texttt{\qquad\qquad\qquad return -1}
\par
\texttt{\qquad\qquad if a == 0:}
\par
\texttt{\qquad\qquad\qquad return 0}
\par
\texttt{\qquad\qquad if a
$>$
0:}
\par
\texttt{\qquad\qquad\qquad return 1}}\medskip\Needspace{6pc}

\texttt{\qquad def bezout\_pair(a, b): \# for b nonnegative }

\texttt{\qquad\qquad if b == 0:}

\texttt{\qquad\qquad\qquad return (sign(a), 0)}

\texttt{\qquad\qquad(u, v) = bezout\_pair(b, a\%b)}

\texttt{\qquad\qquad return (v, u - v * (a//b))} \medskip

This algorithm is known as the \textbf{extended Euclidean algorithm}.

\bigskip

Now that Lemma \ref{lem.ent.gcd.bezout-ZN} has been proven, Bezout's theorem
in the general case (Theorem \ref{thm.ent.gcd.bezout}) easily follows:

\begin{proof}
[Proof of Theorem \ref{thm.ent.gcd.bezout}.]We are in one of the following two cases:

\textit{Case 1:} We have $b\geq0$.

\textit{Case 2:} We have $b<0$.

Let us first consider Case 1. In this case, $b\geq0$. Hence, $b\in\mathbb{N}$.
Thus, Lemma \ref{lem.ent.gcd.bezout-ZN} yields that there exists a Bezout pair
for $\left(  a,b\right)  $. In other words, there exists a pair $\left(
x,y\right)  $ of two integers satisfying $\gcd\left(  a,b\right)  =xa+yb$ (by
the definition of a Bezout pair). But this is precisely what Theorem
\ref{thm.ent.gcd.bezout} is claiming. Thus, Theorem \ref{thm.ent.gcd.bezout}
is proved in Case 1.

Let us now consider Case 2. In this case, $b<0$. Hence, $-b>0$, so that
$-b\in\mathbb{N}$. Hence, Lemma \ref{lem.ent.gcd.bezout-ZN} (applied to $-b$
instead of $b$) yields that there exists a Bezout pair for $\left(
a,-b\right)  $. Let $\left(  u,v\right)  $ be this Bezout pair. Then, by the
definition of a Bezout pair, $u$ and $v$ are integers and satisfy $\gcd\left(
a,-b\right)  =ua+v\left(  -b\right)  $.

However, Proposition \ref{prop.ent.gcd.basics} \textbf{(h)} yields
$\gcd\left(  a,-b\right)  =\gcd\left(  a,b\right)  $. Thus,%
\[
\gcd\left(  a,b\right)  =\gcd\left(  a,-b\right)  =ua+\underbrace{v\left(
-b\right)  }_{=\left(  -v\right)  b}=ua+\left(  -v\right)  b.
\]
Thus, there exist two integers $x$ and $y$ such that $\gcd\left(  a,b\right)
=xa+yb$ (namely, $x=u$ and $y=-v$). This proves Theorem
\ref{thm.ent.gcd.bezout} in Case 2.

We have now proved Theorem \ref{thm.ent.gcd.bezout} in both Cases 1 and 2, so
that the theorem always holds.
\end{proof}

\begin{exercise}
\label{exe.23hw4.1}Recall the \texttt{bezout\_pair} function defined above.
This function outputs a Bezout pair for any given pair $\left(  a,b\right)  $
with $a\in\mathbb{Z}$ and $b\in\mathbb{N}$. Tweak it so that it works for
arbitrary $b\in\mathbb{Z}$ (not just for $b\in\mathbb{N}$). \medskip

[Feel free to use your favorite programming language instead of Python, but do
not change the logic in the case when $b\geq0$.]
\end{exercise}

\begin{exercise}
\label{exe.23hw4.2andmore}\textbf{(a)} Prove that $\gcd\left(
2n+3,\ 3n+4\right)  =1$ for each $n\in\mathbb{Z}$. \medskip

\textbf{(b)} Prove that $\gcd\left(  15n+4,\ 12n+5\right)  =1$ for each
$n\in\mathbb{Z}$.
\end{exercise}

\begin{exercise}
\label{exe.23hw4.3}Let $a\in\mathbb{Z}$ and $b\in\mathbb{Z}$ be nonzero
integers. Let $\left(  x,y\right)  $ be some Bezout pair for $\left(
a,b\right)  $.

Let $g=\gcd\left(  a,b\right)  $. Let $a^{\prime}=a/g$ and $b^{\prime}=b/g$.

Prove that each Bezout pair for $\left(  a,b\right)  $ can be written in the
form $\left(  x+kb^{\prime},\ y-ka^{\prime}\right)  $ for some $k\in
\mathbb{Z}$. \medskip

[\textbf{Hint:} It is probably easiest to first prove this in the case when
$\gcd\left(  a,b\right)  =1$. In this case, $g=1$ and $a^{\prime}=a$ and
$b^{\prime}=b$.]
\end{exercise}

\begin{exercise}
Let $a\in\mathbb{Z}$. Prove that%
\[
\gcd\left(  a^{n}-1,\ a^{m}-1\right)  =\left\vert a^{\gcd\left(  n,m\right)
}-1\right\vert \ \ \ \ \ \ \ \ \ \ \text{for all }n,m\in\mathbb{N}.
\]

[\textbf{Hint:} Start by showing that if $n>m$, then $a^{n}-1=a^{n-m}\left(
a^{m}-1\right)  +\left(  a^{n-m}-1\right)  $.]
\end{exercise}

\begin{exercise}
\label{exe.fib.gcd1}Prove that the Fibonacci numbers (see Definition
\ref{def.ind.fib}) satisfy%
\[
\gcd\left(  f_{n},f_{n+1}\right)  =1\ \ \ \ \ \ \ \ \ \ \text{for each }%
n\in\mathbb{N}.
\]

\end{exercise}

\subsubsection{The universal property of the gcd}

Bezout's theorem is helpful for proving properties of gcds. Here is the most
important one, which is called the \textbf{universal property of the gcd}:

\begin{theorem}
[universal property of the gcd]\label{thm.ent.gcd.uniprop}Let $a,b,m\in
\mathbb{Z}$. Then, we have the equivalence%
\[
\left(  m\mid a\text{ and }m\mid b\right)  \ \Longleftrightarrow\ \left(
m\mid\gcd\left(  a,b\right)  \right)  .
\]

\end{theorem}

In other words, the common divisors of $a$ and $b$ are precisely the divisors
of $\gcd\left(  a,b\right)  $. In other words, $\gcd\left(  a,b\right)  $ is
not just the greatest among the common divisors of $a$ and $b$ (if $a$ and $b$
are not both $0$), but it also is divisible by all of them.

\begin{proof}
[Proof of Theorem \ref{thm.ent.gcd.uniprop}.]We must prove the two
implications%
\[
\left(  m\mid a\text{ and }m\mid b\right)  \ \Longrightarrow\ \left(
m\mid\gcd\left(  a,b\right)  \right)
\]
and%
\[
\left(  m\mid\gcd\left(  a,b\right)  \right)  \ \Longrightarrow\ \left(  m\mid
a\text{ and }m\mid b\right)  .
\]

The second of these two implications is easy to prove: If $m\mid\gcd\left(
a,b\right)  $, then $m\mid a$ (since $m\mid\gcd\left(  a,b\right)  \mid a$)
and $m\mid b$ (similarly).

It thus remains to prove the first implication: i.e., to prove that%
\[
\left(  m\mid a\text{ and }m\mid b\right)  \ \Longrightarrow\ \left(
m\mid\gcd\left(  a,b\right)  \right)  .
\]

To prove this, we assume that $m\mid a$ and $m\mid b$. We must show that
$m\mid\gcd\left(  a,b\right)  $.

Bezout's theorem (Theorem \ref{thm.ent.gcd.bezout}) tells us that there exist
two integers $x$ and $y$ such that $\gcd\left(  a,b\right)  =xa+yb$. Consider
these $x$ and $y$. Then, $m\mid a\mid xa$, so that $xa$ is a multiple of $m$.
Similarly, $yb$ is a multiple of $m$. Thus, $xa+yb$ is a multiple of $m$ as
well (since a sum of two multiples of $m$ is again a multiple of $m$). But
this is saying that $\gcd\left(  a,b\right)  $ is a multiple of $m$ (since
$\gcd\left(  a,b\right)  =xa+yb$). In other words, $m\mid\gcd\left(
a,b\right)  $. But this is precisely what we wanted to show. Thus, the first
implication is proved, and the proof of Theorem \ref{thm.ent.gcd.uniprop} is complete.
\end{proof}

We note that Theorem \ref{thm.ent.gcd.uniprop} is commonly used in the
\textquotedblleft$\Longrightarrow$\textquotedblright\ direction (since the
\textquotedblleft$\Longleftarrow$\textquotedblright\ direction is trivial).
That is, the following fact is used most of the time:

\begin{corollary}
[universal property of the gcd, forward direction]%
\label{cor.ent.gcd.uniprop.==>}Let $a,b,m\in\mathbb{Z}$. If $m\mid a$ and
$m\mid b$, then $m\mid\gcd\left(  a,b\right)  $.
\end{corollary}

\begin{proof}
This is the \textquotedblleft$\Longrightarrow$\textquotedblright\ direction of
Theorem \ref{thm.ent.gcd.uniprop}.
\end{proof}

\begin{exercise}
\label{exe.ent.gcd-monoton}Let $a_{1},a_{2},b_{1},b_{2}\in\mathbb{Z}$ be
integers satisfying $a_{1}\mid b_{1}$ and $a_{2}\mid b_{2}$. Prove that
$\gcd\left(  a_{1},a_{2}\right)  \mid\gcd\left(  b_{1},b_{2}\right)  $.
\end{exercise}

\subsubsection{Factoring out a common factor from a gcd}

The following theorem has an \textquotedblleft intuitively
obvious\textquotedblright\ feel, but its proof is not as simple as you might suspect:

\begin{theorem}
\label{thm.ent.gcd.sasb}Let $s,a,b\in\mathbb{Z}$. Then,%
\[
\gcd\left(  sa,sb\right)  =\left\vert s\right\vert \cdot\gcd\left(
a,b\right)  .
\]

\end{theorem}

This is saying that when two integers have a common factor $s$, then this
common factor can be pulled out of their gcd. (The caveat is, of course, that
the common factor must be replaced by its absolute value, since a gcd cannot
be negative by definition.)

\begin{proof}
[Proof of Theorem \ref{thm.ent.gcd.sasb}.]Let
\[
g=\gcd\left(  a,b\right)  \ \ \ \ \ \ \ \ \ \ \text{and}%
\ \ \ \ \ \ \ \ \ \ h=\gcd\left(  sa,sb\right)  .
\]
Thus, we must prove that $h=\left\vert s\right\vert \cdot g$. Note that $h$
and $g$ are nonnegative (because Proposition \ref{prop.ent.gcd.basics}
\textbf{(a)} shows that gcds are always nonnegative). Thus, $h=\left\vert
h\right\vert $ and $g=\left\vert g\right\vert $, so that $\left\vert
s\right\vert \cdot g=\left\vert s\right\vert \cdot\left\vert g\right\vert
=\left\vert sg\right\vert $ (since $\left\vert x\right\vert \cdot\left\vert
y\right\vert =\left\vert xy\right\vert $ for any two real numbers $x$ and $y$).

Our goal is to prove that $h=\left\vert s\right\vert \cdot g$. Since
$h=\left\vert h\right\vert $ and $\left\vert s\right\vert \cdot g=\left\vert
sg\right\vert $, this amounts to proving that $\left\vert h\right\vert
=\left\vert sg\right\vert $. So this is our goal now.

One good way to prove that two integers $p$ and $q$ satisfy $\left\vert
p\right\vert =\left\vert q\right\vert $ is by showing that $p\mid q$ and
$q\mid p$. Indeed, from $p\mid q$ and $q\mid p$, it follows that $\left\vert
p\right\vert =\left\vert q\right\vert $ (by Proposition \ref{prop.div.abs}
\textbf{(c)}).

Thus, in order to prove that $\left\vert h\right\vert =\left\vert
sg\right\vert $, it will suffice to show that $h\mid sg$ and $sg\mid h$. Now,
let us do this.

\begin{itemize}
\item \textit{Proof of }$sg\mid h$\textit{:} We have $g=\gcd\left(
a,b\right)  \mid a$. Multiplying both sides by $s$, we thus obtain $sg\mid
sa$\ \ \ \ \footnote{\textquotedblleft Multiplying both sides by
$s$\textquotedblright\ means using the following simple fact: If two integers
$x$ and $y$ satisfy $x\mid y$, then $sx\mid sy$.}. Similarly, $sg\mid sb$.
Hence, Corollary \ref{cor.ent.gcd.uniprop.==>} (applied to $sg$, $sa$ and $sb$
instead of $m$, $a$ and $b$) yields $sg\mid\gcd\left(  sa,sb\right)  $. In
other words, $sg\mid h$ (since $h=\gcd\left(  sa,sb\right)  $).

\item \textit{First proof of }$h\mid sg$\textit{:} If $s=0$, then the claim
$h\mid sg$ is obvious (since $\underbrace{s}_{=0}g=0=h\cdot0$). Thus, let us
consider the case when $s\neq0$.

We have just showed that $sg\mid h$, but we also clearly have $s\mid sg$.
Thus, $s\mid sg\mid h$. Since $s\neq0$, this entails that $\dfrac{h}{s}%
\in\mathbb{Z}$ (by Proposition \ref{prop.div.abs} \textbf{(d)}, applied to $s$
and $h$ instead of $a$ and $b$).

This integer $\dfrac{h}{s}$ satisfies $s\cdot\dfrac{h}{s}=h=\gcd\left(
sa,sb\right)  \mid sa$. Dividing both sides by $s$, we thus obtain $\dfrac
{h}{s}\mid a$\ \ \ \ \footnote{\textquotedblleft Dividing both sides by
$s$\textquotedblright\ means using the following simple fact: If two integers
$x$ and $y$ satisfy $sx\mid sy$, then $x\mid y$. (Note that this relies on
$s\neq0$.)}. Similarly, $\dfrac{h}{s}\mid b$. Hence, Corollary
\ref{cor.ent.gcd.uniprop.==>} (applied to $m=\dfrac{h}{s}$) yields $\dfrac
{h}{s}\mid\gcd\left(  a,b\right)  $. In other words, $\dfrac{h}{s}\mid g$
(since $g=\gcd\left(  a,b\right)  $). Multiplying both sides by $s$, we thus
obtain $s\cdot\dfrac{h}{s}\mid sg$. In other words, $h\mid sg$. Thus, $h\mid
sg$ is proved.

\item \textit{Second proof of }$h\mid sg$\textit{:} We have $h=\gcd\left(
sa,sb\right)  \mid sa$. In other words, $sa=hu$ for some integer $u$.
Similarly, $sb=hv$ for some integer $v$. Consider these integers $u$ and $v$.

However, Bezout's theorem (Theorem \ref{thm.ent.gcd.bezout}) shows that there
exist two integers $x$ and $y$ such that $\gcd\left(  a,b\right)  =xa+yb$.
Consider these $x$ and $y$.

Now, $g=\gcd\left(  a,b\right)  =xa+yb$, so that%
\[
sg=s\left(  xa+yb\right)  =sxa+syb=\underbrace{sa}_{=hu}x+\underbrace{sb}%
_{=hv}y=hux+hvy=h\underbrace{\left(  ux+vy\right)  }_{\text{an integer}}.
\]
This again proves that $h\mid sg$.
\end{itemize}

We have now proved that $h\mid sg$ (proved in two different ways) and $sg\mid
h$. Hence, as explained above, we obtain $\left\vert h\right\vert =\left\vert
sg\right\vert $. As we also explained above, this completes our proof of
Theorem \ref{thm.ent.gcd.sasb}.
\end{proof}

\subsection{Coprime integers}

\subsubsection{Definition and examples}

Greatest common divisors are at their most useful when they are $1$. This is
called \textquotedblleft coprimality\textquotedblright:

\begin{definition}
\label{def.ent.coprime.coprime}Two integers $a$ and $b$ are said to be
\textbf{coprime} (or \textbf{relatively prime}) if $\gcd\left(  a,b\right)
=1$.
\end{definition}

\begin{remark}
This is a symmetric relation: If $a$ and $b$ are coprime, then $b$ and $a$ are
coprime (since $\gcd\left(  b,a\right)  =\gcd\left(  a,b\right)  $).
\end{remark}

\begin{example}
\textbf{(a)} An integer $n$ is coprime to $2$ if and only if $n$ is odd.
Indeed, we know that $\gcd\left(  n,2\right)  $ is a divisor of $2$ and is a
nonnegative integer (since any gcd is a nonnegative integer). Thus,
$\gcd\left(  n,2\right)  $ must be either $1$ or $2$ (since the only
nonnegative divisors of $2$ are $1$ and $2$). Now:

\begin{itemize}
\item If $\gcd\left(  n,2\right)  =2$, then $n$ is even (since $2=\gcd\left(
n,2\right)  \mid n$).

\item If $\gcd\left(  n,2\right)  =1$, then $n$ is odd (because otherwise, $2$
would be a common divisor of $n$ and $2$, but this cannot happen when the
greatest common divisor of $n$ and $2$ is $1$).
\end{itemize}

\textbf{(b)} An integer $n$ is coprime to $3$ if and only if $n$ is not
divisible by $3$. (This can be proved just as part \textbf{(a)}, since the
only nonnegative divisors of $3$ are $1$ and $3$.) \medskip

\textbf{(c)} An integer $n$ is coprime to $4$ if and only if $n$ is odd. (If
you expected \textquotedblleft... if $n$ is not divisible by $4$%
\textquotedblright\ here, then you were wrong. The nonnegative divisors of $4$
are not only $1$ and $4$ but also $2$. Thus, $\gcd\left(  n,4\right)  $ can be
$1$, $2$ or $4$. Specifically:

\begin{itemize}
\item If $\gcd\left(  n,4\right)  =1$, then $n$ is odd (since otherwise, $2$
would be a common divisor of $n$ and $4$, but this cannot happen when the
greatest common divisor is $1$).

\item If $\gcd\left(  n,4\right)  =2$, then $n$ is even (since $2=\gcd\left(
n,4\right)  \mid n$).

\item If $\gcd\left(  n,4\right)  =4$, then $n$ is even as well (since
$2\mid4=\gcd\left(  n,4\right)  \mid n$). )
\end{itemize}

\textbf{(d)} An integer $n$ is coprime to $5$ if and only if $n$ is not
divisible by $5$. (This can be proved just as part \textbf{(a)}, since the
only nonnegative divisors of $5$ are $1$ and $5$.) \medskip

\textbf{(e)} An integer $n$ is coprime to $6$ if and only if $n$ is neither
even nor divisible by $3$. (Indeed, the nonnegative divisors of $6$ are $1$,
$2$, $3$ and $6$. Thus, $\gcd\left(  n,6\right)  $ can be $1$, $2$, $3$ or
$6$. Specifically:

\begin{itemize}
\item If $\gcd\left(  n,6\right)  =1$, then $n$ is odd (since otherwise, $2$
would be a common divisor of $n$ and $6$, but this cannot happen when the
greatest common divisor is $1$) and not divisible by $3$ (for similar reasons
but using $3$ instead of $2$).

\item If $\gcd\left(  n,6\right)  =2$, then $n$ is even (since $2=\gcd\left(
n,6\right)  \mid n$).

\item If $\gcd\left(  n,6\right)  =3$, then $n$ is divisible by $3$ (since
$3=\gcd\left(  n,6\right)  \mid n$).

\item If $\gcd\left(  n,6\right)  =6$, then $n$ is both even and divisible by
$3$.)
\end{itemize}

\textbf{(f)} An integer $n$ is always coprime to $1$. Indeed, the only
nonnegative divisor of $1$ is $1$, so that $\gcd\left(  n,1\right)  $ must
always be $1$.
\end{example}

Informally, I think of coprimality as some sort of \textquotedblleft
unrelatedness\textquotedblright\ or \textquotedblleft
independence\textquotedblright\ or \textquotedblleft
orthogonality\textquotedblright\ or \textquotedblleft
noninterference\textquotedblright\ relation. In other words, two integers $a$
and $b$ are coprime if and only if they have \textquotedblleft nothing to do
with each other\textquotedblright, in some sense. This is nowhere near a
rigorous statement, but it motivates many properties of coprimality, including
the ones we will see below.

\subsubsection{Three theorems about coprimality}

The following three theorems are useful properties of coprime integers:

\begin{theorem}
[coprime divisors theorem]\label{thm.ent.coprime.ab|c}Let $a,b,c\in\mathbb{Z}$
satisfy $a\mid c$ and $b\mid c$. Assume that $a$ and $b$ are coprime. Then,
$ab\mid c$.

(In other words, a product of two coprime divisors of $c$ is again a divisor
of $c$.)
\end{theorem}

\begin{proof}
We have $ab\mid ac$ (since $b\mid c$) and $ba\mid bc$ (because $a\mid c$).
Since $ba=ab$ and $ac=ca$ and $bc=cb$, we can rewrite this as follows: We have
$ab\mid ca$ and $ab\mid cb$. Thus, Corollary \ref{cor.ent.gcd.uniprop.==>}
(applied to $ab$, $ca$ and $cb$ instead of $m$, $a$ and $b$) yields%
\begin{align*}
ab  &  \mid\gcd\left(  ca,cb\right)  =\left\vert c\right\vert \cdot
\underbrace{\gcd\left(  a,b\right)  }_{\substack{=1\\\text{(since }a\text{ is
coprime to }b\text{)}}}\ \ \ \ \ \ \ \ \ \ \left(  \text{by Theorem
\ref{thm.ent.gcd.sasb}}\right) \\
&  =\left\vert c\right\vert .
\end{align*}
Since divisibility does not depend on signs (Proposition \ref{prop.div.abs}
\textbf{(a)}), we thus obtain $ab\mid c$\ \ \ \ \footnote{Here is this
argument in detail: We have just proved that $ab\mid\operatorname*{abs}c$
(where we write $\operatorname*{abs}x$ for $\left\vert x\right\vert $ in order
to avoid confusing absolute-value bars with divisibility symbols). Proposition
\ref{prop.div.abs} \textbf{(a)} shows that we have $ab\mid c$ if and only if
$\operatorname*{abs}\left(  ab\right)  \mid\operatorname*{abs}c$. However, the
same proposition shows that we have $ab\mid\operatorname*{abs}c$ if and only
if $\operatorname*{abs}\left(  ab\right)  \mid\operatorname*{abs}\left(
\operatorname*{abs}c\right)  $. Since $\operatorname*{abs}\left(
\operatorname*{abs}c\right)  =\operatorname*{abs}c$, the latter statement can
be rewritten as $\operatorname*{abs}\left(  ab\right)  \mid\operatorname*{abs}%
c$. Thus, both statements $ab\mid c$ and $ab\mid\operatorname*{abs}c$ are
equivalent to $\operatorname*{abs}\left(  ab\right)  \mid\operatorname*{abs}%
c$, and thus are equivalent to each other. Hence, from $ab\mid
\operatorname*{abs}c$, we obtain $ab\mid c$.}. This proves Theorem
\ref{thm.ent.coprime.ab|c}.
\end{proof}

\begin{example}
We have $4\mid56$ and $7\mid56$. Since $4$ and $7$ are coprime, we can thus
conclude (by Theorem \ref{thm.ent.coprime.ab|c}, applied to $a=4$, $b=7$ and
$c=56$) that $4\cdot7\mid56$.

In contrast, from $6\mid12$ and $4\mid12$, we cannot conclude that
$6\cdot4\mid12$, since $6$ and $4$ are not coprime.
\end{example}

\begin{fineprint}
In terms of our \textquotedblleft coprimality as
independence\textquotedblright\ heuristic, Theorem \ref{thm.ent.coprime.ab|c}
can be made intuitive as follows: If $a$ and $b$ are two coprime divisors of
$c$, then (because $a$ and $b$ are coprime) $a$ and $b$ must divide
\textquotedblleft different parts\textquotedblright\ of $c$, and thus their
product $ab$ is still a divisor of $c$. Of course, the notion of
\textquotedblleft different parts\textquotedblright\ here is not a real thing,
but it is helpful as a mnemonic device.
\end{fineprint}

\begin{theorem}
[coprime removal theorem]\label{thm.ent.coprime.a|c}Let $a,b,c\in\mathbb{Z}$
satisfy $a\mid bc$. Assume that $a$ is coprime to $b$. Then, $a\mid c$.
\end{theorem}

\begin{proof}
We have $a\mid ca$ and $a\mid bc=cb$. Thus, Corollary
\ref{cor.ent.gcd.uniprop.==>} (applied to $a$, $ca$ and $cb$ instead of $m$,
$a$ and $b$) yields
\begin{align*}
a  &  \mid\gcd\left(  ca,cb\right)  =\left\vert c\right\vert \cdot
\underbrace{\gcd\left(  a,b\right)  }_{\substack{=1\\\text{(since }a\text{ is
coprime to }b\text{)}}}\ \ \ \ \ \ \ \ \ \ \left(  \text{by Theorem
\ref{thm.ent.gcd.sasb}}\right) \\
&  =\left\vert c\right\vert .
\end{align*}
Since divisibility does not depend on signs, this means that $a\mid c$. Thus,
Theorem \ref{thm.ent.coprime.a|c} holds.
\end{proof}

\begin{example}
We have $6\mid7\cdot12$, but $6$ is coprime to $7$. Thus, Theorem
\ref{thm.ent.coprime.a|c} (applied to $a=6$, $b=7$ and $c=12$) yields
$6\mid12$ (as if you didn't know this already).

But we cannot obtain $6\mid7$ from $6\mid12\cdot7$, since $6$ is not coprime
to $12$.
\end{example}

\begin{fineprint}
Again, Theorem \ref{thm.ent.coprime.a|c} can be motivated using the
\textquotedblleft independence\textquotedblright\ view on coprimality: If $a$
is coprime to $b$, then $b$ cannot be the \textquotedblleft
reason\textquotedblright\ for the divisibility $a\mid bc$, and thus $b$ can be
removed from this divisibility. Again, this is neither a proof nor even a
rigorous statement, but it makes Theorem \ref{thm.ent.coprime.a|c} look less surprising.
\end{fineprint}

\begin{theorem}
[coprime product theorem]\label{thm.ent.coprime.ab-co-c}Let $a,b,c\in
\mathbb{Z}$. Assume that each of the numbers $a$ and $b$ is coprime to $c$.
Then, $ab$ is also coprime to $c$.
\end{theorem}

\begin{proof}
Let $g=\gcd\left(  ab,c\right)  $. Thus, we must prove that $g=1$.

We have $g=\gcd\left(  ab,c\right)  \mid ab$ and $g=\gcd\left(  ab,c\right)
\mid c\mid ac$. Hence, Corollary \ref{cor.ent.gcd.uniprop.==>} (applied to
$g$, $ab$ and $ac$ instead of $m$, $a$ and $b$) yields
\begin{align*}
g  &  \mid\gcd\left(  ab,ac\right)  =\left\vert a\right\vert \cdot
\underbrace{\gcd\left(  b,c\right)  }_{\substack{=1\\\text{(because }b\text{
is coprime to }c\text{)}}}\ \ \ \ \ \ \ \ \ \ \left(  \text{by Theorem
\ref{thm.ent.gcd.sasb}}\right) \\
&  =\left\vert a\right\vert \cdot1=\left\vert a\right\vert .
\end{align*}
Hence, $g\mid a$ (since divisibility does not depend on signs). Combining this
with $g\mid c$, we obtain $g\mid\gcd\left(  a,c\right)  $ (by Corollary
\ref{cor.ent.gcd.uniprop.==>}, applied to $g$, $a$ and $c$ instead of $m$, $a$
and $b$). However, $\gcd\left(  a,c\right)  =1$ (since $a$ is coprime to $c$),
so we obtain $g\mid\gcd\left(  a,c\right)  =1$.

However, $g$ is a nonnegative integer (since any gcd is a nonnegative
integer). Thus, $g$ is a nonnegative divisor of $1$ (since $g\mid1$). Since
the only nonnegative divisor of $1$ is $1$, we thus conclude that $g=1$.
Hence, $\gcd\left(  ab,c\right)  =g=1$. This shows that $ab$ is coprime to
$c$, and we have proved Theorem \ref{thm.ent.coprime.ab-co-c}.
\end{proof}

\begin{example}
Each of the numbers $3$ and $4$ is coprime to $5$. Thus, Theorem
\ref{thm.ent.coprime.ab-co-c} (applied to $a=3$, $b=4$ and $c=5$) yields that
$3\cdot4$ is coprime to $5$.
\end{example}

\begin{fineprint}
Again, Theorem \ref{thm.ent.coprime.ab-co-c} can be viewed within the
\textquotedblleft independence\textquotedblright\ paradigm: If each of $a$ and
$b$ is coprime to $c$, then so should be $ab$, because any \textquotedblleft
dependence\textquotedblright\ between $ab$ and $c$ should come from $a$ or
from $b$. Alternatively, if you think of coprimality as an analogue of
orthogonality, then you can view Theorem \ref{thm.ent.coprime.ab-co-c} as an
analogue of the fact that if two vectors $\overrightarrow{a}$ and
$\overrightarrow{b}$ are both orthogonal to a given vector $\overrightarrow{c}%
$, then so is their sum $\overrightarrow{a}+\overrightarrow{b}$. Again, none
of these metaphors should be mistaken for a proof of Theorem
\ref{thm.ent.coprime.ab-co-c}.
\end{fineprint}

\begin{exercise}
\label{exe.25hw4.1}Let $a$ and $b$ be two coprime integers. Let $i,j\in
\mathbb{N}$ be arbitrary. Prove that $a^{i}$ and $b^{j}$ are again coprime.
\end{exercise}

\begin{fineprint}
Theorems \ref{thm.ent.coprime.ab|c}, \ref{thm.ent.coprime.a|c} and
\ref{thm.ent.coprime.ab-co-c} can be generalized, dropping some of the
coprimality assumptions (but leading to less memorable results). Here is the
generalization of Theorem \ref{thm.ent.coprime.ab|c}:

\begin{theorem}
\label{thm.ent.coprime.ab|c-gen}Let $a,b,c\in\mathbb{Z}$ satisfy $a\mid c$ and
$b\mid c$. Then, $ab\mid\gcd\left(  a,b\right)  \cdot c$.
\end{theorem}

\begin{proof}
Read our above proof of Theorem \ref{thm.ent.coprime.ab|c} until the point
where it shows that $ab\mid\left\vert c\right\vert \cdot\gcd\left(
a,b\right)  $. Now, observe that $\left\vert c\right\vert $ divides $c$ (since
$\left\vert c\right\vert $ is either $c$ or $-c$), and thus $\left\vert
c\right\vert \cdot\gcd\left(  a,b\right)  $ divides $c\cdot\gcd\left(
a,b\right)  $. Hence,%
\[
ab\mid\left\vert c\right\vert \cdot\gcd\left(  a,b\right)  \mid c\cdot
\gcd\left(  a,b\right)  =\gcd\left(  a,b\right)  \cdot c.
\]
This proves Theorem \ref{thm.ent.coprime.ab|c-gen}.
\end{proof}

Here is the generalization of Theorem \ref{thm.ent.coprime.a|c}:

\begin{theorem}
\label{thm.ent.coprime.a|c-gen}Let $a,b,c\in\mathbb{Z}$ satisfy $a\mid bc$.
Then, $a\mid\gcd\left(  a,b\right)  \cdot c$.
\end{theorem}

\begin{proof}
Read our above proof of Theorem \ref{thm.ent.coprime.a|c} until the point
where it shows that $a\mid\left\vert c\right\vert \cdot\gcd\left(  a,b\right)
$. Now, observe that $\left\vert c\right\vert $ divides $c$ (since $\left\vert
c\right\vert $ is either $c$ or $-c$), and thus $\left\vert c\right\vert
\cdot\gcd\left(  a,b\right)  $ divides $c\cdot\gcd\left(  a,b\right)  $.
Hence,%
\[
a\mid\left\vert c\right\vert \cdot\gcd\left(  a,b\right)  \mid c\cdot
\gcd\left(  a,b\right)  =\gcd\left(  a,b\right)  \cdot c.
\]
This proves Theorem \ref{thm.ent.coprime.a|c-gen}.
\end{proof}
\end{fineprint}

\subsubsection{Reducing a fraction}

Here is one more property of gcds:

\begin{theorem}
\label{thm.ent.gcd.reduce}Let $a$ and $b$ be two integers that are not both
$0$. Let $g=\gcd\left(  a,b\right)  $. Then, the integers $\dfrac{a}{g}$ and
$\dfrac{b}{g}$ are coprime.
\end{theorem}

This theorem is important for understanding rational numbers. Indeed, a ratio
$\dfrac{u}{v}$ of two integers is said to be in \textbf{reduced form} if $u$
and $v$ are coprime. Now, Theorem \ref{thm.ent.gcd.reduce} shows that if we
start with a ratio $\dfrac{a}{b}$ of two integers, and cancel $\gcd\left(
a,b\right)  $ from the numerator and the denominator, then the result will be
a ratio in reduced form. Hence, each rational number can be brought to a
reduced form. For example, $\dfrac{12}{21}=\dfrac{12/3}{21/3}=\dfrac{4}{7}$.

\begin{proof}
[Proof of Theorem \ref{thm.ent.gcd.reduce}.]Since $a$ and $b$ are not both
$0$, we have $\gcd\left(  a,b\right)  \neq0$ (since $0$ cannot divide any
nonzero integer). Since we know that $\gcd\left(  a,b\right)  \in\mathbb{N}$,
we thus conclude that $\gcd\left(  a,b\right)  >0$. In other words, $g>0$
(since $g=\gcd\left(  a,b\right)  $). Thus, $\dfrac{a}{g}$ and $\dfrac{b}{g}$
are well-defined. Also, from $g>0$, we obtain $\left\vert g\right\vert =g$.

Since $g=\gcd\left(  a,b\right)  $, we have $g\mid a$ and $g\mid b$. Hence,
$\dfrac{a}{g}$ and $\dfrac{b}{g}$ are integers. Moreover,
\begin{align*}
g  &  =\gcd\left(  a,b\right)  =\gcd\left(  g\cdot\dfrac{a}{g},\ g\cdot
\dfrac{b}{g}\right)  \ \ \ \ \ \ \ \ \ \ \left(  \text{since }a=g\cdot
\dfrac{a}{g}\text{ and }b=g\cdot\dfrac{b}{g}\right) \\
&  =\underbrace{\left\vert g\right\vert }_{=g}\cdot\gcd\left(  \dfrac{a}%
{g},\dfrac{b}{g}\right)  \ \ \ \ \ \ \ \ \ \ \left(
\begin{array}
[c]{c}%
\text{by Theorem \ref{thm.ent.gcd.sasb},}\\
\text{since }\dfrac{a}{g}\text{ and }\dfrac{b}{g}\text{ are integers}%
\end{array}
\right) \\
&  =g\cdot\gcd\left(  \dfrac{a}{g},\dfrac{b}{g}\right)  .
\end{align*}
Dividing this equality by $g$, we find%
\[
1=\gcd\left(  \dfrac{a}{g},\dfrac{b}{g}\right)  \ \ \ \ \ \ \ \ \ \ \left(
\text{since }g\neq0\right)  .
\]
This shows that $\dfrac{a}{g}$ and $\dfrac{b}{g}$ are coprime. Thus, Theorem
\ref{thm.ent.gcd.reduce} is proven.
\end{proof}

\begin{exercise}
\label{exe.ent.gcd(ak,bk)}Let $a,b\in\mathbb{Z}$ and $k\in\mathbb{N}$. Prove
that%
\[
\gcd\left(  a^{k},b^{k}\right)  =\left(  \gcd\left(  a,b\right)  \right)
^{k}.
\]

\end{exercise}

\begin{exercise}
\label{exe.ent.gcd(ab,c)}Let $a,b,c$ be three integers. \medskip

\textbf{(a)} Prove that%
\[
\gcd\left(  ab,c\right)  \mid\gcd\left(  a,c\right)  \cdot\gcd\left(
b,c\right)  .
\]

\textbf{(b)} Prove that if $a$ and $b$ are coprime, then%
\[
\gcd\left(  ab,c\right)  =\gcd\left(  a,c\right)  \cdot\gcd\left(  b,c\right)
.
\]

[\textbf{Hint:} For part \textbf{(a)}, show first that $\gcd\left(
ab,c\right)  $ divides both $a\gcd\left(  b,c\right)  $ and $c\gcd\left(
b,c\right)  $. For part \textbf{(b)}, begin by arguing that $\gcd\left(
a,c\right)  $ and $\gcd\left(  b,c\right)  $ are themselves coprime.]
\end{exercise}

\begin{exercise}
\label{exe.ent.prod-is-kth-power}Let $a,b,k$ be three positive integers such
that $a$ and $b$ are coprime. Assume that $ab$ is the $k$-th power of an
integer -- say, $ab=c^{k}$ for some $c\in\mathbb{Z}$. Show that $a$ and $b$
are $k$-th powers of integers, too, and in fact we have%
\[
a=\left(  \gcd\left(  a,c\right)  \right)  ^{k}\ \ \ \ \ \ \ \ \ \ \text{and}%
\ \ \ \ \ \ \ \ \ \ b=\left(  \gcd\left(  b,c\right)  \right)  ^{k}.
\]

[\textbf{Hint:} Assume WLOG that $c\geq0$. First show that $a\mid\gcd\left(
a^{k},c^{k}\right)  =\left(  \gcd\left(  a,c\right)  \right)  ^{k}$ and
likewise $b\mid\left(  \gcd\left(  b,c\right)  \right)  ^{k}$. However,
$c=\gcd\left(  ab,c\right)  $ (why?) and thus $ab=c^{k}=\left(  \gcd\left(
ab,c\right)  \right)  ^{k}$. Now use Exercise \ref{exe.ent.gcd(ab,c)}
\textbf{(b)}.]
\end{exercise}

\begin{exercise}
\label{exe.ent.a!b!-div-a+b-1!}Let $a$ and $b$ be two coprime positive
integers. Prove that%
\[
a!\cdot b!\mid\left(  a+b-1\right)  !.
\]

[\textbf{Hint:} Let $g=\dfrac{\left(  a+b-1\right)  !}{a!\cdot b!}$. Prove
that both $ag$ and $bg$ are integers. What then?]
\end{exercise}

\begin{exercise}
Prove that the Fibonacci numbers (see Definition \ref{def.ind.fib}) satisfy%
\[
\gcd\left(  f_{n},f_{m}\right)  =f_{\gcd\left(  n,m\right)  }%
\ \ \ \ \ \ \ \ \ \ \text{for all }n,m\in\mathbb{N}.
\]

[\textbf{Hint:} First show that $f_{n}=f_{n-m-1}f_{m}+f_{n-m}f_{m+1}$ whenever
$n>m$. Then recall Exercise \ref{exe.fib.gcd1}.]
\end{exercise}

\subsection{\label{sec.ent.pri}Prime numbers}

\subsubsection{Definition}

The following is one of the most famous concepts in mathematics:

\begin{definition}
\label{def.ent.prime.prime}An integer $n>1$ is said to be \textbf{prime} (or
\textbf{a prime}) if the only positive divisors of $n$ are $1$ and $n$.
\end{definition}

The first few primes (= prime numbers) are%
\[
2,\ \ 3,\ \ 5,\ \ 7,\ \ 11,\ \ 13,\ \ 17,\ \ 19,\ \ 23,\ \ 29,\ \ 31,\ \ 37,\ \ 41,\ \ 43.
\]

It can be shown that there are infinitely many primes (see Exercise
\ref{exe.23hw4.4} \textbf{(b)} for one proof).

\subsubsection{The friend-or-foe lemma}

The first property of primes that we will show is an important result that we
call the \textbf{friend-or-foe lemma}:

\begin{lemma}
[friend-or-foe lemma]\label{lem.ent.prime.fof}Let $p$ be a prime. Let
$n\in\mathbb{Z}$. Then, $n$ is either divisible by $p$ or coprime to $p$, but
not both.
\end{lemma}

\begin{proof}
The number $p$ is prime, and thus its only positive divisors are $1$ and $p$.
Since $\gcd\left(  n,p\right)  $ is a positive divisor of $p$ (this is easy to
see\footnote{\textit{Proof.} The number $\gcd\left(  n,p\right)  $ is a
divisor of $p$, and thus is nonzero (since $0$ does not divide $p$).
Furthermore, $\gcd\left(  n,p\right)  $ is nonnegative (since any gcd is
nonnegative). Thus, $\gcd\left(  n,p\right)  $ is positive. Hence,
$\gcd\left(  n,p\right)  $ is a positive divisor of $p$.}), we thus conclude
that $\gcd\left(  n,p\right)  $ must be either $1$ or $p$. So we are in one of
the following two cases:

\textit{Case 1:} We have $\gcd\left(  n,p\right)  =1$.

\textit{Case 2:} We have $\gcd\left(  n,p\right)  =p$.

Let us first consider Case 1. In this case, we have $\gcd\left(  n,p\right)
=1$. In other words, $n$ is coprime to $p$. Furthermore, the greatest common
divisor of $n$ and $p$ is $\gcd\left(  n,p\right)  =1$; therefore, $p$ cannot
be a common divisor of $n$ and $p$ (since $p>1$). Thus, $n$ is not divisible
by $p$ (since this would entail that $p$ is a common divisor of $n$ and $p$).
So we have shown that $n$ is coprime to $p$ and not divisible by $p$. Thus,
Lemma \ref{lem.ent.prime.fof} is proved in Case 1.

Let us now consider Case 2. In this case, we have $\gcd\left(  n,p\right)
=p\neq1$. Thus, $n$ is not coprime to $p$. Also, $p=\gcd\left(  n,p\right)
\mid n$ shows that $n$ is divisible by $p$. So we have shown that $n$ is
divisible by $p$ and not coprime to $p$. Hence, Lemma \ref{lem.ent.prime.fof}
is proved in Case 2.

We have now proved Lemma \ref{lem.ent.prime.fof} in both Cases 1 and 2; thus,
Lemma \ref{lem.ent.prime.fof} is fully proved.
\end{proof}

(The moniker \textquotedblleft friend-or-foe lemma\textquotedblright\ is
metaphorical: You can think of integers that are divisible by $p$ as
\textquotedblleft friends of $p$\textquotedblright, and think of integers
coprime to $p$ as \textquotedblleft foes of $p$\textquotedblright. Thus, a
prime number cleanly divides the integers into its \textquotedblleft
friends\textquotedblright\ and its \textquotedblleft foes\textquotedblright.
In contrast, the non-prime number $4$ has a more \textquotedblleft
nuanced\textquotedblright\ relationship with certain integers such as $2$
(since $2$ is neither divisible by $4$ nor coprime to $4$).)

\subsubsection{There are infinitely many primes, and some more exercises}

\begin{exercise}
\label{exe.23hw4.4}Let $\left(  a_{0},a_{1},a_{2},\ldots\right)  $ be a
sequence of integers defined recursively by%
\[
a_{n}=1+a_{0}a_{1}\cdots a_{n-1}\ \ \ \ \ \ \ \ \ \ \text{for all }n\geq0.
\]
(This sequence has been studied in Exercise \ref{exe.23mt1.5}.) \medskip

\textbf{(a)} Prove that $\gcd\left(  a_{n},a_{m}\right)  =1$ for any two
distinct integers $n,m\in\mathbb{N}$. \medskip

For each $n\in\mathbb{N}$, let $p_{n}$ be a prime that divides $a_{n}$. (Such
a prime exists, since $a_{n}=1+\underbrace{a_{0}a_{1}\cdots a_{n-1}}_{\geq
1}\geq1+1>1$. Of course, there will often be several choices. In this case,
just choose one.) \medskip

\textbf{(b)} Prove that the primes $p_{0},p_{1},p_{2},\ldots$ are distinct.
\end{exercise}

Exercise \ref{exe.23hw4.4} \textbf{(b)} shows that there are infinitely many
primes. This is a famous result of Euclid; many other proofs of it are known
(see, e.g., \cite{Conrad-primes}). \medskip

All primes except for $2$ are odd. Thus, the distances between consecutive
primes (except for $2$ and $3$) are always even. Beside this, however, these
distances are rather unpredictable. For instance, the two consecutive primes
$41$ and $43$ are a distance of $2$ apart, whereas the two consecutive primes
$113$ and $127$ are a distance of $14$ apart. Even some very simple-sounding
questions, such as \textquotedblleft are there infinitely many pairs of
consecutive primes at a distance of $2$ from each other?\textquotedblright%
\ (such primes are called \textbf{twin primes}) are so-far unresolved (this
one is called the \textbf{twin primes conjecture}). At least, one can show
that \textbf{three} consecutive primes cannot be at distances of $2$ from each
other except in one very simple case:

\begin{exercise}
\label{exe.23hw4.5}Let $p$ be a prime such that $p-2$ and $p+2$ are also
prime. Prove that $p=5$.

[\textbf{Hint:} Consider the remainders upon division by $6$.]
\end{exercise}

\begin{exercise}
\label{exe.23mt2.2}Let $p$ be a prime larger than $3$. Prove that $p^{2}%
\equiv1\operatorname{mod}24$.

[\textbf{Hint:} Recall some older problems. Also note that the integers $3$
and $8$ are coprime.]
\end{exercise}

The following exercise helps checking whether a given integer $n$ is prime:

\begin{exercise}
\label{exe.21finfin.4rewri}Let $n$ be an integer such that $n>1$ but $n$ is
not a prime. Let $d$ be the smallest divisor of $n$ that is larger than $1$.
Prove that $d^{2}\leq n$. \medskip

(You can use standard properties of inequalities -- e.g., the equivalence
$\left(  u\leq v\right)  \ \Longleftrightarrow\ \left(  u^{2}\leq
v^{2}\right)  $ when $u$ and $v$ are positive.)
\end{exercise}

\subsubsection{Binomial coefficients and primes}

The friend-or-foe lemma has myriad applications. As a first example, recall
Pascal's triangle (which we saw in Section \ref{sec.sums.bin-def}):%
\[
\hspace{-1.0114pc}\hspace{-0.3cm}%
\begin{tabular}
[c]{ll|ccccccccccccccccc}%
\vphantom{$\overset{k=0}{\swarrow}$} &  & $\phantom{20}$ & $\phantom{20}$ &
$\phantom{20}$ & $\phantom{20}$ & $\phantom{20}$ & $\phantom{20}$ &  &  &  &
$\overset{k=0}{\swarrow}$ & $\phantom{20}$ & $\phantom{20}$ & $\phantom{20}$ &
$\phantom{20}$ & $\phantom{20}$ & $\phantom{20}$ & $\phantom{20}$\\
\vphantom{$\overset{k=0}{\swarrow}$}$n=0$ & $\rightarrow$ &  &  &  &  &  &  &
&  & $1$ &  & $\overset{k=1}{\swarrow}$ &  &  &  &  &  & \\
\vphantom{$\overset{k=0}{\swarrow}$}$n=1$ & $\rightarrow$ &  &  &  &  &  &  &
& $1$ &  & $1$ &  & $\overset{k=2}{\swarrow}$ &  &  &  &  & \\
\vphantom{$\overset{k=0}{\swarrow}$}$n=2$ & $\rightarrow$ &  &  &  &  &  &  &
$1$ &  & $2$ &  & $1$ &  & $\overset{k=3}{\swarrow}$ &  &  &  & \\
\vphantom{$\overset{k=0}{\swarrow}$}$n=3$ & $\rightarrow$ &  &  &  &  &  & $1$
&  & $3$ &  & $3$ &  & $1$ &  & $\overset{k=4}{\swarrow}$ &  &  & \\
\vphantom{$\overset{k=0}{\swarrow}$}$n=4$ & $\rightarrow$ &  &  &  &  & $1$ &
& $4$ &  & $6$ &  & $4$ &  & $1$ &  & $\overset{k=5}{\swarrow}$ &  & \\
\vphantom{$\overset{k=0}{\swarrow}$}$n=5$ & $\rightarrow$ &  &  &  & $1$ &  &
$5$ &  & $10$ &  & $10$ &  & $5$ &  & $1$ &  & $\overset{k=6}{\swarrow}$ & \\
\vphantom{$\overset{k=0}{\swarrow}$}$n=6$ & $\rightarrow$ &  &  & $1$ &  & $6$
&  & $15$ &  & $20$ &  & $15$ &  & $6$ &  & $1$ &  & $\overset{k=7}{\swarrow}%
$\\
\vphantom{$\overset{k=0}{\swarrow}$}$n=7$ & $\rightarrow$ &  & $1$ &  & $7$ &
& $21$ &  & $35$ &  & $35$ &  & $21$ &  & $7$ &  & $1$ & \\
\vphantom{$\overset{k=0}{\swarrow}$}$n=8$ & $\rightarrow$ & $1$ &  & $8$ &  &
$28$ &  & $56$ &  & $70$ &  & $56$ &  & $28$ &  & $8$ &  & $1$%
\end{tabular}
\ \ \ \ \ \ \ \ \ \ \ \
\]

One property of Pascal's triangle that you might have already noticed is the
following: All entries in its $n=7$ row except for the two $1$'s (i.e., all
the binomial coefficients $\dbinom{7}{1},\dbinom{7}{2},\ldots,\dbinom{7}{6}$)
are divisible by $7$; all entries in the $n=5$ row except for the two $1$'s
are divisible by $5$; likewise for the $n=3$ and $n=2$ rows. The pattern here
can be generalized to any prime number instead of $7,5,3,2$:

\begin{theorem}
\label{thm.ent.prime.p|(pck)}Let $p$ be a prime. Let $k\in\left\{
1,2,\ldots,p-1\right\}  $. Then, $p\mid\dbinom{p}{k}$.
\end{theorem}

\begin{proof}
Exercise \ref{exe.23hw3.5} \textbf{(a)} (applied to $n=p$) yields%
\[
k\dbinom{p}{k}=p\underbrace{\dbinom{p-1}{k-1}}_{\substack{\text{an
integer}\\\text{(by Theorem \ref{thm.binom.N})}}}.
\]
Thus, $p\mid k\dbinom{p}{k}$.

From $k\in\left\{  1,2,\ldots,p-1\right\}  $, we furthermore obtain $p\nmid k$
(because if we had $p\mid k$, then Proposition \ref{prop.div.abs} \textbf{(b)}
would entail $\left\vert p\right\vert \leq\left\vert k\right\vert $, which
would contradict $\left\vert k\right\vert =k\leq p-1<p=\left\vert p\right\vert
$). In other words, $k$ is not divisible by $p$.

But the friend-or-foe lemma (Lemma \ref{lem.ent.prime.fof}, applied to $n=k$)
says that $k$ is either divisible by $p$ or coprime to $p$. Since $k$ is not
divisible by $p$, we thus conclude that $k$ must be coprime to $p$. In other
words, $p$ is coprime to $k$. Hence, from $p\mid k\dbinom{p}{k}$, we obtain
$p\mid\dbinom{p}{k}$ using the coprime cancellation theorem (Theorem
\ref{thm.ent.coprime.a|c}, applied to $a=p$ and $b=k$ and $c=\dbinom{p}{k}$).
This proves Theorem \ref{thm.ent.prime.p|(pck)}.
\end{proof}

Theorem \ref{thm.ent.prime.p|(pck)} shows that if $p$ is a prime, then all the
binomial coefficients $\dbinom{p}{i}$ in the $p$-th row of Pascal's triangle
are divisible by $p$ (except for the two $1$'s on the borders of the
triangle). The following exercise, in contrast, claims that the binomial
coefficients in the $\left(  p-1\right)  $-st row are alternatingly congruent
to $1$ and to $-1$ modulo $p$:

\begin{exercise}
\label{exe.23mt2.1}Let $p$ be a prime. Prove that%
\[
\dbinom{p-1}{i}\equiv\left(  -1\right)  ^{i}\operatorname{mod}%
p\ \ \ \ \ \ \ \ \ \ \text{for each }i\in\left\{  0,1,\ldots,p-1\right\}  .
\]

[\textbf{Hint:} What connects the three binomial coefficients $\dbinom{p-1}%
{i}$, $\dbinom{p-1}{i-1}$ and $\dbinom{p}{i}$ ?]
\end{exercise}

\subsubsection{Fermat's little theorem}

It is easy to see that every integer $a$ satisfies $a^{2}\equiv
a\operatorname{mod}2$. Indeed, the difference $a^{2}-a=a\left(  a-1\right)  $
is divisible by $2$, since at least one of the two consecutive integers $a$
and $a-1$ must be even and thus contributes a factor of $2$ to the product
$a\left(  a-1\right)  $.

Likewise, every integer $a$ satisfies $a^{3}\equiv a\operatorname{mod}3$,
since the difference $a^{3}-a=\left(  a-1\right)  a\left(  a+1\right)  $ is
divisible by $3$ (because at least one of the three consecutive integers
$a-1$, $a$ and $a+1$ must be divisible by $3$).

This pattern does not persist for $4$: Indeed, $a^{4}\equiv
a\operatorname{mod}4$ does not hold for $a=2$. However, for $5$, the pattern
emerges again: Every integer $a$ satisfies $a^{5}\equiv a\operatorname{mod}5$.
This is not as easy to see as the analogous claims for $a^{2}$ and $a^{3}$
(since $a^{5}-a$ does not factor into linear factors any more), but still can
be checked with a bit of work (there are only $5$ possible values for the
remainder $a\%5$, and each of these values allows us to check $a^{5}\equiv
a\operatorname{mod}5$ by reducing both sides modulo $5$).

The pattern is lost again for $6$ (the congruence $a^{6}\equiv
a\operatorname{mod}6$ fails for $a=2$), but reemerges for $7$.

As you may have guessed, there is a general result here:

\begin{theorem}
[Fermat's Little Theorem]\label{thm.ent.prime.flt}Let $p$ be a prime. Let
$a\in\mathbb{Z}$. Then,%
\[
a^{p}\equiv a\operatorname{mod}p.
\]

\end{theorem}

\begin{proof}
We shall induct on $a$. This will only cover the case $a\geq0$, so we will
have to handle the case $a<0$ by a separate argument afterwards.

\textit{Base case:} The congruence $a^{p}\equiv a\operatorname{mod}p$ clearly
holds for $a=0$ (since $0^{p}=0\equiv0\operatorname{mod}p$).

\textit{Induction step:} Let $a\in\mathbb{N}$. Assume (as the induction
hypothesis) that $a^{p}\equiv a\operatorname{mod}p$. We must prove that
$\left(  a+1\right)  ^{p}\equiv a+1\operatorname{mod}p$.

But the binomial formula (Theorem \ref{thm.binom.bin}) yields%
\begin{align*}
\left(  a+1\right)  ^{p}  &  =\sum_{k=0}^{p}\dbinom{p}{k}a^{k}%
\underbrace{1^{p-k}}_{=1}=\sum_{k=0}^{p}\dbinom{p}{k}a^{k}\\
&  =\underbrace{\dbinom{p}{0}}_{=1}\underbrace{a^{0}}_{=1}+\sum_{k=1}%
^{p-1}\dbinom{p}{k}a^{k}+\underbrace{\dbinom{p}{p}}_{=1}a^{p}\\
&  \ \ \ \ \ \ \ \ \ \ \ \ \ \ \ \ \ \ \ \ \left(
\begin{array}
[c]{c}%
\text{here, we have split off the addends}\\
\text{for }k=0\text{ and for }k=p\text{ from the sum}%
\end{array}
\right) \\
&  =1+\sum_{k=1}^{p-1}\dbinom{p}{k}a^{k}+a^{p}=\sum_{k=1}^{p-1}\dbinom{p}%
{k}a^{k}+a^{p}+1.
\end{align*}
In other words,%
\begin{equation}
\left(  a+1\right)  ^{p}-\left(  a^{p}+1\right)  =\sum_{k=1}^{p-1}\dbinom
{p}{k}a^{k}. \label{pf.thm.ent.prime.flt.3}%
\end{equation}

However, Theorem \ref{thm.ent.prime.p|(pck)} shows that each $k\in\left\{
1,2,\ldots,p-1\right\}  $ satisfies $p\mid\dbinom{p}{k}\mid\dbinom{p}{k}a^{k}%
$. In other words, $\dbinom{p}{k}a^{k}$ is a multiple of $p$ for each
$k\in\left\{  1,2,\ldots,p-1\right\}  $. Hence, $\sum_{k=1}^{p-1}\dbinom{p}%
{k}a^{k}$ is a sum of multiples of $p$, and thus itself a multiple of $p$.
That is, we have $p\mid\sum_{k=1}^{p-1}\dbinom{p}{k}a^{k}$. In view of
(\ref{pf.thm.ent.prime.flt.3}), we can rewrite this as $p\mid\left(
a+1\right)  ^{p}-\left(  a^{p}+1\right)  $. In other words,%
\begin{equation}
\left(  a+1\right)  ^{p}\equiv a^{p}+1\operatorname{mod}p.
\label{pf.thm.ent.prime.flt.4}%
\end{equation}

However, the induction hypothesis says that $a^{p}\equiv a\operatorname{mod}%
p$. Adding the obvious congruence $1\equiv1\operatorname{mod}p$ to this, we
obtain%
\[
a^{p}+1\equiv a+1\operatorname{mod}p.
\]
Combining this congruence with (\ref{pf.thm.ent.prime.flt.4}), we obtain%
\[
\left(  a+1\right)  ^{p}\equiv a^{p}+1\equiv a+1\operatorname{mod}p,
\]
which shows that $\left(  a+1\right)  ^{p}\equiv a+1\operatorname{mod}p$ (by
the transitivity of congruence). This completes the induction step.

Thus, Theorem \ref{thm.ent.prime.flt} is proved for all $a\geq0$. It remains
to prove it for all $a<0$ now. This can be done with a neat trick:

Let $a\in\mathbb{Z}$ satisfy $a<0$. Then, we must prove that $a^{p}\equiv
a\operatorname{mod}p$.

But we already know that $b^{p}\equiv b\operatorname{mod}p$ for all integers
$b\geq0$ (because we have already proved Theorem \ref{thm.ent.prime.flt} for
all $a\geq0$). We can apply this to $b=a\%p$ (since the remainder $a\%p$ is
$\geq0$), and thus obtain%
\[
\left(  a\%p\right)  ^{p}\equiv a\%p\operatorname{mod}p.
\]

However, Proposition \ref{prop.ent.quo-rem.basics} \textbf{(a)} (applied to
$n=a$ and $d=p$) shows that $a\%p\in\left\{  0,1,\ldots,p-1\right\}  $ and
$a\%p\equiv a\operatorname{mod}p$. Taking the congruence $a\%p\equiv
a\operatorname{mod}p$ to the $p$-th power, we obtain $\left(  a\%p\right)
^{p}\equiv a^{p}\operatorname{mod}p$ (we have here used Exercise
\ref{exe.mod.k-th-power}). Therefore, $a^{p}\equiv\left(  a\%p\right)
^{p}\operatorname{mod}p$. Combining all the congruences we have obtained so
far, we obtain%
\[
a^{p}\equiv\left(  a\%p\right)  ^{p}\equiv a\%p\equiv a\operatorname{mod}p,
\]
from which we can conclude that $a^{p}\equiv a\operatorname{mod}p$ (by
transitivity of congruence). Thus, we have proved Theorem
\ref{thm.ent.prime.flt} for $a<0$. This completes the proof of Theorem
\ref{thm.ent.prime.flt}.
\end{proof}

Fermat's Little Theorem has a bunch of applications, some of which we will see
later (Subsection \ref{subsec.ent.crypt.RSA}).

\begin{fineprint}
One wrinkle in the pattern we have discussed above: Theorem
\ref{thm.ent.prime.flt} shows that every prime $p$ satisfies $a^{p}\equiv
a\operatorname{mod}p$ for all $a\in\mathbb{Z}$. But there are some positive
integers $p$ that satisfy this even though they are not prime! The smallest
such integers are $1,\ \ 561,\ \ 1105,\ \ 1729,\ \ 2465$. See
\href{https://en.wikipedia.org/wiki/Carmichael_number}{Carmichael numbers} for
more details.
\end{fineprint}

\subsubsection{Prime divisor separation theorem}

You can think of the primes as \textquotedblleft inseparable\textquotedblright%
\ positive integers: They cannot be written as products of two smaller
positive integers. (Of course, $1$ also has this property but does not count
as a prime. In a way, $1$ is inseparable because there is nothing to separate,
but this is not the \textquotedblleft inseparability\textquotedblright\ that
we are looking for in a prime.)

One useful consequence of this \textquotedblleft
inseparability\textquotedblright\ is that if a prime $p$ divides a product
$ab$ of two integers, then it must divide one of the two factors $a$ and $b$,
since (speaking heuristically) it cannot be \textquotedblleft
separated\textquotedblright\ into a part that divides $a$ and a part that
divides $b$. Never mind that this is not a valid argument, the conclusion is a
true fact:

\begin{theorem}
[prime divisor separation theorem]\label{thm.ent.prime.sep}Let $p$ be a prime.
Let $a,b\in\mathbb{Z}$ be such that $p\mid ab$. Then, $p\mid a$ or $p\mid b$.
\end{theorem}

\begin{proof}
[Proof of Theorem \ref{thm.ent.prime.sep}.]We shall prove the claim of Theorem
\ref{thm.ent.prime.sep} in the following equivalent form: \textquotedblleft If
$p\nmid a$, then $p\mid b$.\textquotedblright

Assume that $p\nmid a$. We must then prove that $p\mid b$.

The friend-or-foe lemma (Lemma \ref{lem.ent.prime.fof}) yields that $a$ is
either divisible by $p$ or coprime to $p$. Thus, $a$ is coprime to $p$ (since
$p\nmid a$). In other words, $p$ is coprime to $a$. Hence, we can use the
coprime cancellation theorem (Theorem \ref{thm.ent.coprime.a|c}, applied to
$p$, $a$ and $b$ instead of $a$, $b$ and $c$) to obtain $p\mid b$ from $p\mid
ab$. This is precisely what we wanted to prove. Theorem
\ref{thm.ent.prime.sep} is thus proved.
\end{proof}

Theorem \ref{thm.ent.prime.sep} shows that if a prime number $p$ divides a
product $ab$, then it must divide $a$ or $b$ (or both). In contrast, a
non-prime number like $4$ can divide a product $ab$ without dividing $a$ or
$b$. For example, $4\mid2\cdot6$ but $4\nmid2$ and $4\nmid6$. \medskip

We can extend Theorem \ref{thm.ent.prime.sep} to products of several factors:

\begin{corollary}
[prime divisor separation theorem for $k$ factors]\label{cor.ent.prime.sepk}%
Let $p$ be a prime. Let $a_{1},a_{2},\ldots,a_{k}\in\mathbb{Z}$ be such that
$p\mid a_{1}a_{2}\cdots a_{k}$. Then, $p\mid a_{i}$ for some $i\in\left\{
1,2,\ldots,k\right\}  $.

(In words: If a prime divides a product of several integers, then it must
divide at least one of the factors.)
\end{corollary}

\begin{proof}
[Proof sketch.]Induct on $k$. In the induction step, use Theorem
\ref{thm.ent.prime.sep}. (The base case is the case $k=0$, in which case
Corollary \ref{cor.ent.prime.sepk} is vacuously true because $p\nmid1$.) (See
\cite[Proposition 2.13.7]{19s} for this proof in detail.)
\end{proof}

The following exercise is another form of Fermat's Little Theorem:

\begin{exercise}
\label{exe.ent.prime.flt-p-1}Let $p$ be a prime. Let $a\in\mathbb{Z}$ be an
integer not divisible by $p$. Prove that $a^{p-1}\equiv1\operatorname{mod}p$.
\end{exercise}

\subsubsection{$p$-valuations: definition}

We will need the following simple lemma:

\begin{lemma}
\label{lem.ent.prime.vpn.wd}Let $p$ be a prime. Let $n$ be a nonzero integer.
Then, there exists a largest $m\in\mathbb{N}$ such that $p^{m}\mid n$.
\end{lemma}

\begin{proof}
The relation $p^{m}\mid n$ means that $\dfrac{n}{p^{m}}\in\mathbb{Z}$. In
other words, it means that we can divide $n$ by $p$ at least $m$ times without
obtaining a non-integer. So the claim of Lemma \ref{lem.ent.prime.vpn.wd} is
saying that there is a largest number of times that we can divide $n$ by $p$
without obtaining a non-integer. But this is clear: Every time we divide $n$
by $p$, the absolute value $\left\vert n\right\vert $ decreases (since $p>1$
and $n\neq0$), and obviously this cannot go on forever without eventually
yielding a non-integer.\footnote{Of course, we are also tacitly using the fact
that $n$ is an integer in the first place, so that $m=0$ does satisfy
$p^{m}\mid n$ (since $p^{0}=1\mid n$).}

(See \cite[Proof of Lemma 2.13.22]{19s} for a more formal proof of Lemma
\ref{lem.ent.prime.vpn.wd}.)
\end{proof}

Lemma \ref{lem.ent.prime.vpn.wd} allows us to make the following definition:

\begin{definition}
\label{def.ent.prime.vpn}Let $p$ be a prime. \medskip

\textbf{(a)} Let $n$ be a nonzero integer. Then, $v_{p}\left(  n\right)  $
shall denote the largest $m\in\mathbb{N}$ such that $p^{m}\mid n$. (This is
well-defined by Lemma \ref{lem.ent.prime.vpn.wd}. Thus, $v_{p}\left(
n\right)  $ is the number of times that you can divide $n$ by $p$ without
getting a non-integer.)

This number $v_{p}\left(  n\right)  $ will be called the $p$%
\textbf{-valuation} (or the $p$\textbf{-adic valuation}) of $n$. \medskip

\textbf{(b)} In order to have $v_{p}\left(  n\right)  $ defined for all
integers $n$ (as opposed to just for nonzero $n$), we also define
$v_{p}\left(  0\right)  $ to be $\infty$ (because $0$ can be divided by $p$ an
arbitrary number of times without any changes). This symbol $\infty$ is not an
actual number, but we shall pretend that it behaves like a number at least in
some regards. In particular, we will eventually add or compare it to other
numbers. In doing so, we shall follow the rules that%
\begin{align*}
k+\infty &  =\infty+k=\infty\ \ \ \ \ \ \ \ \ \ \text{for all }k\in
\mathbb{Z};\\
\infty+\infty &  =\infty;\\
k  &  <\infty\text{ and }\infty>k\ \ \ \ \ \ \ \ \ \ \text{for all }%
k\in\mathbb{Z}.
\end{align*}
Thus, $\infty$ acts like a \textquotedblleft mythical number that is larger
than any actual number\textquotedblright. We can keep up this charade as long
as we only add and compare, but never subtract $\infty$ from anything (since
$1+\infty=\infty$ would turn into $1=0$ if you subtracted $\infty$).
\end{definition}

Here are some examples:%
\begin{align*}
v_{3}\left(  99\right)   &  =2\ \ \ \ \ \ \ \ \ \ \left(  \text{since }%
3^{2}\mid99\text{ but }3^{3}\nmid99\right)  ;\\
v_{3}\left(  98\right)   &  =0\ \ \ \ \ \ \ \ \ \ \left(  \text{since }%
3^{0}\mid98\text{ but }3^{1}\nmid98\right)  ;\\
v_{3}\left(  96\right)   &  =1\ \ \ \ \ \ \ \ \ \ \left(  \text{since }%
3^{1}\mid96\text{ but }3^{2}\nmid96\right)  ;\\
v_{3}\left(  0\right)   &  =\infty.
\end{align*}

We can restate the definition of $v_{p}\left(  n\right)  $ in yet another way:
If $p$ is a prime and $n$ is a positive integer, then $v_{p}\left(  n\right)
$ is the number of zeroes at the end of the base-$p$ representation of the
number $n$. For example, the base-$2$ representation of the number $344$ is
$101011000$, which has three zeroes at its end (the other zeroes don't
count!), so that $v_{2}\left(  344\right)  =3$.

Note that Definition \ref{def.ent.prime.vpn} can be generalized to any
positive integer $p>1$ (prime or not). But most of the useful properties of
$p$-valuations hold only when $p$ is prime.

\subsubsection{$p$-valuations: basic properties}

Let us now discuss some basic properties of $p$-valuations. We begin with a
lemma that is almost trivial, but quite helpful:

\begin{lemma}
\label{lem.ent.vp.equivalence}Let $p$ be a prime. Let $i\in\mathbb{N}$ and
$n\in\mathbb{Z}$. Then, $p^{i}\mid n$ if and only if $v_{p}\left(  n\right)
\geq i$.
\end{lemma}

\begin{proof}
If $n=0$, then this is clear (because in this case, we have both $p^{i}%
\mid0=n$ and $v_{p}\left(  n\right)  =v_{p}\left(  0\right)  =\infty\geq i$).

It remains to deal with the case $n\neq0$. In this case, $v_{p}\left(
n\right)  $ is defined as the largest $m\in\mathbb{N}$ such that $p^{m}\mid
n$. Thus, in this case, we have $p^{i}\mid p^{v_{p}\left(  n\right)  }\mid n$
whenever $i\leq v_{p}\left(  n\right)  $, whereas $p^{i}\nmid n$ whenever
$i>v_{p}\left(  n\right)  $. In other words, we have $p^{i}\mid n$ if and only
if $i\leq v_{p}\left(  n\right)  $. In other words, we have $p^{i}\mid n$ if
and only if $v_{p}\left(  n\right)  \geq i$. Thus, Lemma
\ref{lem.ent.vp.equivalence} is proved in this case as well.
\end{proof}

Recall some standard notations: For any two numbers $x$ and $y$, we let
$\min\left\{  x,y\right\}  $ denote the smaller of these two numbers, and we
let $\max\left\{  x,y\right\}  $ denote the larger of these two numbers. More
generally, if $S$ is a set of numbers, then $\min S$ means the smallest
element of $S$ (if it exists), and $\max S$ means the largest element of $S$
(if it exists). We extend these notations to sets that include $\infty$ in the
obvious way (recalling that $\infty$ is larger than any integer). Thus, in particular,%

\begin{align*}
\max\left\{  \infty,k\right\}   &  =\max\left\{  k,\infty\right\}
=\infty\ \ \ \ \ \ \ \ \ \ \text{for all }k\in\mathbb{Z}\cup\left\{
\infty\right\}  ;\\
\min\left\{  \infty,k\right\}   &  =\min\left\{  k,\infty\right\}
=k\ \ \ \ \ \ \ \ \ \ \text{for all }k\in\mathbb{Z}\cup\left\{  \infty
\right\}  .
\end{align*}

Now, we can state a bunch of rather important properties of $p$-valuations:

\begin{theorem}
[basic properties of $p$-valuations]\label{thm.ent.vp.vpab}Let $p$ be a prime.
Then: \medskip

\textbf{(a)} We have $v_{p}\left(  ab\right)  =v_{p}\left(  a\right)
+v_{p}\left(  b\right)  $ for any $a,b\in\mathbb{Z}$. \medskip

\textbf{(b)} We have $v_{p}\left(  a+b\right)  \geq\min\left\{  v_{p}\left(
a\right)  ,\ v_{p}\left(  b\right)  \right\}  $ for any $a,b\in\mathbb{Z}$.
\medskip

\textbf{(c)} We have $v_{p}\left(  1\right)  =0$. \medskip

\textbf{(d)} We have $v_{p}\left(  p\right)  =1$. \medskip

\textbf{(e)} We have $v_{p}\left(  q\right)  =0$ for any prime $q\neq p$.
\end{theorem}

\begin{proof}
\textbf{(a)} Let $a,b\in\mathbb{Z}$. We must prove that $v_{p}\left(
ab\right)  =v_{p}\left(  a\right)  +v_{p}\left(  b\right)  $.

If $a=0$, then this is saying that $\infty=\infty+v_{p}\left(  b\right)  $,
which follows from our rules for $\infty$ (specifically, from the rules saying
that $\infty+k=\infty$ for all $k\in\mathbb{Z}$ and that $\infty+\infty
=\infty$). Likewise, we can prove our claim if $b=0$.

It thus remains to handle the case when neither $a$ nor $b$ is $0$. So let us
consider this case. Since $a$ and $b$ are nonzero, the numbers $v_{p}\left(
a\right)  $ and $v_{p}\left(  b\right)  $ are nonnegative integers. Let us
give them names: We set%
\[
n=v_{p}\left(  a\right)  \ \ \ \ \ \ \ \ \ \ \text{and}%
\ \ \ \ \ \ \ \ \ \ m=v_{p}\left(  b\right)  .
\]

Thus, $p^{n}\mid a$ and $p^{m}\mid b$. In other words, there are integers $x$
and $y$ such that $a=p^{n}x$ and $b=p^{m}y$. Consider these $x$ and $y$.

If we had $p\mid x$, then we would readily obtain $p^{n+1}\mid a$ (because
$p\mid x$ entails that $x=pz$ for some integer $z$, and thus this integer $z$
must satisfy $a=p^{n}\underbrace{x}_{=pz}=p^{n}pz=p^{n+1}z$) and therefore
$v_{p}\left(  a\right)  \geq n+1$ (by Lemma \ref{lem.ent.vp.equivalence},
applied to $n+1$ and $a$ instead of $i$ and $n$), which would contradict
$v_{p}\left(  a\right)  =n<n+1$. Thus, we cannot have $p\mid x$. For similar
reasons, we cannot have $p\mid y$.

However, multiplying $a=p^{n}x$ with $b=p^{m}y$, we obtain $ab=p^{n}x\cdot
p^{m}y=p^{n+m}xy$, and thus $p^{n+m}\mid ab$. Therefore, $v_{p}\left(
ab\right)  \geq n+m$ (by Lemma \ref{lem.ent.vp.equivalence}, applied to $ab$
and $n+m$ instead of $n$ and $i$).

Now, we shall show that this inequality is an equality. To do so, we must show
that $p^{n+m+1}\nmid ab$.

To prove this, we assume the contrary. Thus, $p^{n+m+1}\mid ab=p^{n+m}xy$.
Dividing both sides of this divisibility by $p^{n+m}$, we obtain $p\mid xy$.

However, the prime divisor separation theorem (Theorem \ref{thm.ent.prime.sep}%
) says that if the prime number $p$ divides a product of two integers, then it
must divide one of these two integers. Therefore, from $p\mid xy$, we obtain
either $p\mid x$ or $p\mid y$ (since $x$ and $y$ are integers). But this
contradicts the fact that we cannot have $p\mid x$ and we cannot have $p\mid
y$. This contradiction shows that our assumption must have been wrong. Thus,
we have shown that $p^{n+m+1}\nmid ab$.

So we know that $p^{n+m}\mid ab$ but $p^{n+m+1}\nmid ab$. In other words, the
largest $i\in\mathbb{N}$ that satisfies $p^{i}\mid ab$ is $n+m$. In other
words, $v_{p}\left(  ab\right)  =n+m$ (by the definition of $v_{p}\left(
ab\right)  $). Since $n=v_{p}\left(  a\right)  $ and $m=v_{p}\left(  b\right)
$, we can rewrite this as $v_{p}\left(  ab\right)  =v_{p}\left(  a\right)
+v_{p}\left(  b\right)  $. This proves Theorem \ref{thm.ent.vp.vpab}
\textbf{(a)}.\medskip

\textbf{(b)} Let $a,b\in\mathbb{Z}$. We must prove that $v_{p}\left(
a+b\right)  \geq\min\left\{  v_{p}\left(  a\right)  ,\ v_{p}\left(  b\right)
\right\}  $.

If $\min\left\{  v_{p}\left(  a\right)  ,\ v_{p}\left(  b\right)  \right\}
=\infty$, then this inequality boils down to $\infty\geq\infty$ (because
$\min\left\{  v_{p}\left(  a\right)  ,\ v_{p}\left(  b\right)  \right\}
=\infty$ yields $v_{p}\left(  a\right)  =\infty$ and $v_{p}\left(  b\right)
=\infty$, so that $a=0$ and $b=0$, and thus $a+b=0$ as well, which in turn
leads to $v_{p}\left(  a+b\right)  =\infty$), which is true.

Thus, it remains to handle the case when $\min\left\{  v_{p}\left(  a\right)
,\ v_{p}\left(  b\right)  \right\}  \neq\infty$. Consider this case. Thus,
$\min\left\{  v_{p}\left(  a\right)  ,\ v_{p}\left(  b\right)  \right\}
\in\mathbb{N}$. Set $k=\min\left\{  v_{p}\left(  a\right)  ,\ v_{p}\left(
b\right)  \right\}  $. Then, $k\leq v_{p}\left(  a\right)  $ and $k\leq
v_{p}\left(  b\right)  $. From $k\leq v_{p}\left(  a\right)  $, we obtain
$v_{p}\left(  a\right)  \geq k$ and thus $p^{k}\mid a$ (by Lemma
\ref{lem.ent.vp.equivalence}, applied to $n=a$ and $i=k$). Similarly,
$p^{k}\mid b$. Thus, $a$ and $b$ are multiples of $p^{k}$. Hence, their sum
$a+b$ is also a multiple of $p^{k}$. In other words, $p^{k}\mid a+b$. Using
Lemma \ref{lem.ent.vp.equivalence} (applied to $n=a+b$ and $i=k$), this in
turn entails $v_{p}\left(  a+b\right)  \geq k=\min\left\{  v_{p}\left(
a\right)  ,\ v_{p}\left(  b\right)  \right\}  $. Thus, Theorem
\ref{thm.ent.vp.vpab} \textbf{(b)} is proved.\medskip

\textbf{(c)} This follows from $p^{0}=1\mid1$ and $p^{1}=p\nmid1$.\medskip

\textbf{(d)} This follows from $p^{1}=p\mid p$ and $p^{2}\nmid p$. \medskip

\textbf{(e)} Let $q\neq p$ be a prime. Then, the only positive divisors of $q$
are $1$ and $q$ (since $q$ is a prime). Hence, $p$ is not a positive divisor
of $q$ (since $p\neq1$ and $p\neq q$). Therefore, $p$ is not a divisor of $q$
(since $p$ is positive). In other words, $p\nmid q$. Now, from $p^{0}=1\mid q$
and $p^{1}=p\nmid q$, we obtain $v_{p}\left(  q\right)  =0$. This proves
Theorem \ref{thm.ent.vp.vpab} \textbf{(e)}.
\end{proof}

\begin{corollary}
\label{cor.ent.vp.prodak}Let $p$ be a prime. Then,
\[
v_{p}\left(  a_{1}a_{2}\cdots a_{k}\right)  =v_{p}\left(  a_{1}\right)
+v_{p}\left(  a_{2}\right)  +\cdots+v_{p}\left(  a_{k}\right)
\]
for any $k$ integers $a_{1},a_{2},\ldots,a_{k}$.
\end{corollary}

\begin{proof}
Induct on $k$. The base case uses $v_{p}\left(  1\right)  =0$. The induction
step relies on Theorem \ref{thm.ent.vp.vpab} \textbf{(a)}.
\end{proof}

Note that Theorem \ref{thm.ent.vp.vpab} \textbf{(a)} would fail if we allowed
$p$ to be non-prime. For instance, $v_{4}\left(  2\cdot2\right)  =1$ but
$v_{4}\left(  2\right)  +v_{4}\left(  2\right)  =0+0=0$.

\subsubsection{Back to Hanoi}

Let us take a closer look at $2$-valuations. The sequence%
\begin{align*}
&  \left(  v_{2}\left(  1\right)  ,\ v_{2}\left(  2\right)  ,\ v_{2}\left(
3\right)  ,\ v_{2}\left(  4\right)  ,\ v_{2}\left(  5\right)  ,\ \ldots\right)
\\
&  =\left(
0,\ 1,\ 0,\ 2,\ 0,\ 1,\ 0,\ 3,\ 0,\ 1,\ 0,\ 2,\ 0,\ 1,\ 0,\ 4,\ \ldots\right)
\end{align*}
is called the
\textbf{\href{https://en.wikipedia.org/wiki/Ruler_function}{\textbf{ruler
sequence}}}, as it resembles the pattern of markings on a ruler (a small
marking at every inch, a slightly larger marking every $2$ inches, an even
larger marking every $4$ inches, and so on):%
\[%
\begin{tikzpicture}[scale=0.8]
\begin{scope}[every edge/.style={draw=black,thick}]
\draw[very thick] (0,0) rectangle (16,-4);
\path[-] (2,0) node[above] {$2$} edge (2,-0.5);
\path[-] (4,0) node[above] {$4$} edge (4,-1);
\path[-] (6,0) node[above] {$6$} edge (6,-0.5);
\path[-] (8,0) node[above] {$8$} edge (8,-2);
\path[-] (10,0) node[above] {$10$} edge (10,-0.5);
\path[-] (12,0) node[above] {$12$} edge (12,-1);
\path[-] (14,0) node[above] {$14$} edge (14,-0.5);
\end{scope}
\end{tikzpicture}%
\]
This sequence tends to appear every once in a while in seemingly unexpected
places. Case in point:

\begin{proposition}
\label{prop.rulerfun.hanoi}Let $n\in\mathbb{N}$.

In Section \ref{sec.ind.hanoi}, we proposed a strategy for solving the Tower
of Hanoi puzzle with $n$ disks. Let $S_{n}$ be this strategy.

Let $k\in\left\{  1,2,\ldots,2^{n}-1\right\}  $. Then, the $k$-th move of the
strategy $S_{n}$ moves the $\left(  v_{2}\left(  k\right)  +1\right)  $-th
smallest disk.
\end{proposition}

Thus, in particular, every odd move (i.e., the $1$-st, the $3$-rd, the $5$-th,
and so on moves) moves the smallest disk (since $v_{2}\left(  k\right)  =0$
when $k$ is odd).

The proof of Proposition \ref{prop.rulerfun.hanoi} relies on the following
lemma about $p$-valuations:

\begin{lemma}
\label{lem.ent.vp.hanoi-lem}Let $p$ be a prime. Let $m\in\mathbb{N}$. Let $k$
be an integer such that $p^{m}\nmid k$. Then, $v_{p}\left(  p^{m}+k\right)
=v_{p}\left(  k\right)  $.
\end{lemma}

\begin{fineprint}
\begin{proof}
[Proof of Lemma \ref{lem.ent.vp.hanoi-lem}.]From $p^{m}\nmid k$, we obtain
$k\neq0$, so that $v_{p}\left(  k\right)  \neq\infty$. In other words,
$v_{p}\left(  k\right)  \in\mathbb{N}$.

Let $i=v_{p}\left(  k\right)  $. Thus, $i\in\mathbb{N}$ (since $v_{p}\left(
k\right)  \in\mathbb{N}$), and the definition of $v_{p}\left(  k\right)  $
shows that $p^{i}\mid k$ and $p^{i+1}\nmid k$.

If we had $m\leq i$, then we would have $p^{m}\mid p^{i}\mid k$, which would
contradict $p^{m}\nmid k$. Thus, we cannot have $m\leq i$. In other words, we
have $i<m$. Thus, $i\leq m-1$ (since $i$ and $m$ are integers), so that
$i+1\leq m$. Therefore, $p^{i+1}\mid p^{m}$.

From the definition of $p$-valuations, it follows easily that $v_{p}\left(
p^{m}\right)  =m$ and $v_{p}\left(  -p^{m}\right)  =m$.

The numbers $p^{m}$ and $k$ are multiples of $p^{i}$ (since $p^{i}\mid
p^{i+1}\mid p^{m}$ and $p^{i}\mid k$). Thus, their sum $p^{m}+k$ is a multiple
of $p^{i}$ as well. In other words, $p^{i}\mid p^{m}+k$.

On the other hand, let us show that $p^{i+1}\nmid p^{m}+k$. Indeed, assume the
contrary. Thus, $p^{i+1}\mid p^{m}+k$.

Therefore, the numbers $p^{m}+k$ and $-p^{m}$ are multiples of $p^{i+1}$
(since $p^{i+1}\mid p^{m}+k$ and $p^{i+1}\mid p^{m}\mid p^{m}\cdot\left(
-1\right)  =-p^{m}$). Hence, their sum $\left(  p^{m}+k\right)  +\left(
-p^{m}\right)  $ is a multiple of $p^{i+1}$ as well. In other words, $k$ is a
multiple of $p^{i+1}$ (since $\left(  p^{m}+k\right)  +\left(  -p^{m}\right)
=k$). But this contradicts $p^{i+1}\nmid k$.

This contradiction shows that our assumption was wrong. Hence, $p^{i+1}\nmid
p^{m}+k$ is proved.

Combining $p^{i}\mid p^{m}+k$ with $p^{i+1}\nmid p^{m}+k$, we see that $i$ is
the largest $j\in\mathbb{N}$ satisfying $p^{j}\mid p^{m}+k$. In other words,
$i=v_{p}\left(  p^{m}+k\right)  $. Hence, $v_{p}\left(  p^{m}+k\right)
=i=v_{p}\left(  k\right)  $. This proves Lemma \ref{lem.ent.vp.hanoi-lem}.
\end{proof}

\begin{proof}
[Proof of Proposition \ref{prop.rulerfun.hanoi}.]We will prove Proposition
\ref{prop.rulerfun.hanoi} by induction on $n$:

\textit{Base case:} If $n=0$, then there exists no $k\in\left\{
1,2,\ldots,2^{n}-1\right\}  $ (since the set $\left\{  1,2,\ldots
,2^{n}-1\right\}  =\left\{  1,2,\ldots,2^{0}-1\right\}  =\left\{
1,2,\ldots,0\right\}  $ is empty in this case). Thus, in this case,
Proposition \ref{prop.rulerfun.hanoi} is vacuously true (i.e., true because it
makes a claim about non-existing objects).

\textit{Induction step:} Let $n$ be a positive integer. Assume (as the
induction hypothesis) that Proposition \ref{prop.rulerfun.hanoi} holds for
$n-1$ instead of $n$. We must now prove that Proposition
\ref{prop.rulerfun.hanoi} holds for $n$ as well.

So let $k\in\left\{  1,2,\ldots,2^{n}-1\right\}  $ be arbitrary. We must prove
that the $k$-th move of the strategy $S_{n}$ moves the $\left(  v_{2}\left(
k\right)  +1\right)  $-th smallest disk.

Recall that the strategy $S_{n}$ was defined recursively: It consists of first
performing the strategy $S_{n-1}$ (but with pegs $2$ and $3$ swapped), then
moving the largest disk (from peg $1$ to peg $3$), and then again performing
the strategy $S_{n-1}$ (but now with pegs $1$ and $2$ swapped). Since strategy
$S_{n-1}$ requires $2^{n-1}-1$ moves in total, we thus conclude that

\begin{enumerate}
\item the first $2^{n-1}-1$ moves of strategy $S_{n}$ are identical with the
corresponding moves of strategy $S_{n-1}$ (except that pegs $2$ and $3$ are swapped);

\item the $2^{n-1}$-th move of strategy $S_{n}$ consists in moving the largest disk;

\item the next $2^{n-1}-1$ moves of strategy $S_{n}$ (that is, the moves
numbered $2^{n-1}+1,\ 2^{n-1}+2,\ \ldots,\ 2^{n}-1$) are identical with the
moves of strategy $S_{n-1}$ (except that pegs $1$ and $2$ are swapped).
\end{enumerate}

Therefore, the $k$-th move of the strategy $S_{n}$

\begin{itemize}
\item moves the same disk as the $k$-th move of $S_{n-1}$ if $k<2^{n-1}$;

\item moves the largest disk if $k=2^{n-1}$;

\item moves the same disk as the $\left(  k-2^{n-1}\right)  $-th move of
$S_{n-1}$ if $k>2^{n-1}$.
\end{itemize}

We thus distinguish between the following three cases:

\textit{Case 1:} We have $k<2^{n-1}$.

\textit{Case 2:} We have $k=2^{n-1}$.

\textit{Case 3:} We have $k>2^{n-1}$.

Let us first consider Case 1. In this case, we have $k<2^{n-1}$. Thus, the
$k$-th move of the strategy $S_{n}$ moves the same disk as the $k$-th move of
$S_{n-1}$ (according to the first of the three bullet points above). But our
induction hypothesis shows that the latter move moves the $\left(
v_{2}\left(  k\right)  +1\right)  $-th smallest disk (since $k\in\left\{
1,2,\ldots,2^{n}-1\right\}  $ and $k<2^{n-1}$ entails $k\in\left\{
1,2,\ldots,2^{n-1}-1\right\}  $). Thus, the former move also moves the
$\left(  v_{2}\left(  k\right)  +1\right)  $-th smallest disk. So the claim we
are trying to prove has been proved in Case 1.

Let us now consider Case 2. In this case, we have $k=2^{n-1}$. Thus, the
$k$-th move of the strategy $S_{n}$ moves the largest disk (according to the
second of the three bullet points above), i.e., the $n$-th smallest disk
(since there are $n$ disks in total, so the largest disk is the $n$-th
smallest). However, we have $k=2^{n-1}$ and thus $v_{2}\left(  k\right)
=v_{2}\left(  2^{n-1}\right)  =n-1$, so that $n=v_{2}\left(  k\right)  +1$.
Thus, the $k$-th move of the strategy $S_{n}$ moves the $\left(  v_{2}\left(
k\right)  +1\right)  $-th smallest disk (because we have shown that it moves
the $n$-th smallest disk). So the claim we are trying to prove has been proved
in Case 2.

Let us finally consider Case 3. In this case, we have $k>2^{n-1}$. Thus, the
$k$-th move of the strategy $S_{n}$ moves the same disk as the $\left(
k-2^{n-1}\right)  $-th move of $S_{n-1}$ (according to the third of the three
bullet points above). But our induction hypothesis (applied to $k-2^{n-1}$
instead of $k$) yields that the latter move moves the $\left(  v_{2}\left(
k-2^{n-1}\right)  +1\right)  $-th smallest disk (since $k\in\left\{
1,2,\ldots,2^{n}-1\right\}  $ and $k>2^{n-1}$ entails $k-2^{n-1}\in\left\{
1,2,\ldots,2^{n-1}-1\right\}  $ quite easily\footnote{Here are the details:
From $k\in\left\{  1,2,\ldots,2^{n}-1\right\}  \subseteq\mathbb{Z}$ and
$k>2^{n-1}$, we see immediately that $k-2^{n-1}$ is a positive integer.
Furthermore, from $k\in\left\{  1,2,\ldots,2^{n}-1\right\}  $, we obtain
$k\leq2^{n}-1=2\cdot2^{n-1}-1=2^{n-1}+2^{n-1}-1$, so that $k-2^{n-1}%
\leq2^{n-1}-1$. Since $k-2^{n-1}$ is a positive integer, this results in
$k-2^{n-1}\in\left\{  1,2,\ldots,2^{n-1}-1\right\}  $.}). Thus, the former
move moves the $\left(  v_{2}\left(  k-2^{n-1}\right)  +1\right)  $-th
smallest disk as well.

However, $k-2^{n-1}>0$ (since $k>2^{n-1}$) and $k-2^{n-1}<2^{n-1}$ (since
$k\leq2^{n}-1<2^{n}=2\cdot2^{n-1}=2^{n-1}+2^{n-1}$). Hence, $2^{n-1}\nmid
k-2^{n-1}$ (because if we had $2^{n-1}\mid k-2^{n-1}$, then Proposition
\ref{prop.div.abs} \textbf{(b)} (applied to $a=2^{n-1}$ and $b=k-2^{n-1}$)
would yield $\operatorname*{abs}\left(  2^{n-1}\right)  \leq
\operatorname*{abs}\left(  k-2^{n-1}\right)  =k-2^{n-1}$ (since $k-2^{n-1}%
>0$), so that $k-2^{n-1}\geq\operatorname*{abs}\left(  2^{n-1}\right)
=2^{n-1}$, which would contradict $k-2^{n-1}<2^{n-1}$). Therefore, Lemma
\ref{lem.ent.vp.hanoi-lem} (applied to $2$, $n-1$ and $k-2^{n-1}$ instead of
$p$, $m$ and $k$) yields%
\[
v_{2}\left(  2^{n-1}+k-2^{n-1}\right)  =v_{2}\left(  k-2^{n-1}\right)  ,
\]
so that%
\begin{equation}
v_{2}\left(  k-2^{n-1}\right)  =v_{2}\left(  \underbrace{2^{n-1}+k-2^{n-1}%
}_{=k}\right)  =v_{2}\left(  k\right)  . \label{pf.prop.rulerfun.hanoi.v2=v2}%
\end{equation}

Now, recall that the $k$-th move of the strategy $S_{n}$ moves the $\left(
v_{2}\left(  k-2^{n-1}\right)  +1\right)  $-th smallest disk (as we showed
above). In view of (\ref{pf.prop.rulerfun.hanoi.v2=v2}), we can restate this
as follows: The $k$-th move of the strategy $S_{n}$ moves the $\left(
v_{2}\left(  k\right)  +1\right)  $-th smallest disk. So the claim we are
trying to prove has been proved in Case 3.

Thus, we have proved our claim in all three Cases 1, 2 and 3. In other words,
we have shown that the $k$-th move of the strategy $S_{n}$ moves the $\left(
v_{2}\left(  k\right)  +1\right)  $-th smallest disk. Hence, we have proved
that Proposition \ref{prop.rulerfun.hanoi} holds for $n$. This completes the
induction step. Thus, Proposition \ref{prop.rulerfun.hanoi} is proved.
\end{proof}
\end{fineprint}

The ruler sequence also has an appearance in data storage:

\begin{remark}
A
\href{https://en.wikipedia.org/wiki/Backup_rotation_scheme}{\textquotedblleft
Tower of Hanoi\textquotedblright\ backup scheme} is a backup scheme where you
have several backup drives for your data. Every odd day, you back up to the
first drive. Every even day that is not divisible by $4$, you back up to the
second drive. Every day that is divisible by $4$ but not by $8$, you back up
to the third drive. And so on. Thus, on the $k$-th day, you back up to the
$\left(  v_{2}\left(  k\right)  +1\right)  $-th drive. This scheme ensures
that at every point in time, you have both a fresh backup and several levels
of older backups available.

(Of course, I only said \textquotedblleft day\textquotedblright\ for
simplicity; you can use any unit of time instead. Note that the first drive
will see the largest traffic and therefore will wear out and need replacement
the fastest.)
\end{remark}

\subsubsection{More exercises}

The following exercise generalizes Theorem \ref{thm.ent.prime.sep}:

\begin{exercise}
\label{exe.24hw4.1}Let $p$ be a prime, and let $m\in\mathbb{N}$. Let $a$ and
$b$ be two integers such that $p^{m}\mid ab$ and $p^{m}\nmid a$. Prove that
$p\mid b$.
\end{exercise}

The following exercise generalizes Theorem \ref{thm.ent.prime.p|(pck)}:

\begin{exercise}
\label{exe.23hw4.6}Let $p$ be a prime. Let $m\in\mathbb{N}$, and let
$k\in\left\{  1,2,\ldots,p^{m}-1\right\}  $. Prove that $p\mid\dbinom{p^{m}%
}{k}$.
\end{exercise}

\begin{exercise}
Let $p$ be a prime such that $p\equiv3\operatorname{mod}4$. Let $n\in
\mathbb{Z}$. Prove that $p\nmid n^{2}+1$.

[\textbf{Hint:} In other words, prove that the congruence $n^{2}%
\equiv-1\operatorname{mod}p$ cannot hold. What happens if you take this
congruence to the $\left(  p-1\right)  /2$-th power?]
\end{exercise}

\begin{exercise}
Let $p$ and $q$ be two distinct primes. Prove that $p^{q}+q^{p}\equiv
p+q\operatorname{mod}pq$ and $p^{q-1}+q^{p-1}\equiv1\operatorname{mod}pq$.

[\textbf{Hint:} First show that $p$ and $q$ are coprime.]
\end{exercise}

\subsubsection{The $p$-valuation of $n!$}

What is the $p$-valuation of a factorial $n!$ ? There turns out to be a nice
formula for this:\footnote{See Definition \ref{def.ent.floor} and Definition
\ref{def.ent.quo-rem} for the notations we are using here. The meaning of the
infinite sums will be discussed in the proof of the theorem.}

\begin{theorem}
[de Polignac's formula]\label{thm.ent.vp.polignac}Let $p$ be a prime. Let
$n\in\mathbb{N}$. Then,%
\begin{align*}
v_{p}\left(  n!\right)   &  =\left\lfloor \dfrac{n}{p^{1}}\right\rfloor
+\left\lfloor \dfrac{n}{p^{2}}\right\rfloor +\left\lfloor \dfrac{n}{p^{3}%
}\right\rfloor +\cdots\\
&  =\left(  n//p^{1}\right)  +\left(  n//p^{2}\right)  +\left(  n//p^{3}%
\right)  +\cdots.
\end{align*}

\end{theorem}

\begin{proof}
[Proof sketch.]First, these sums are infinite sums. Why do they make
sense?\footnote{It is trivially easy to concoct an infinite sum that does not
make sense: for instance, $1+1+1+\cdots$, or $\dfrac{1}{1}+\dfrac{1}{2}%
+\dfrac{1}{3}+\cdots$. In general, \textquotedblleft
infinite\textquotedblright\ operations in mathematics (i.e., \textquotedblleft
do something infinitely many times\textquotedblright) do not usually exist
unless their existence has been justified.}

Because we can discard all the addends that are zero, and then only finitely
many nonzero addends remain. For instance, if $p=2$ and $n=13$, then%
\begin{align*}
&  \left\lfloor \dfrac{n}{p^{1}}\right\rfloor +\left\lfloor \dfrac{n}{p^{2}%
}\right\rfloor +\left\lfloor \dfrac{n}{p^{3}}\right\rfloor +\cdots\\
&  =\left\lfloor \dfrac{13}{2^{1}}\right\rfloor +\left\lfloor \dfrac{13}%
{2^{2}}\right\rfloor +\left\lfloor \dfrac{13}{2^{3}}\right\rfloor +\cdots\\
&  =\left\lfloor 6.5\right\rfloor +\left\lfloor 3.25\right\rfloor
+\left\lfloor 1.625\right\rfloor +\left\lfloor 0.8125\right\rfloor
+\left\lfloor 0.40625\right\rfloor +\cdots\\
&  =6+3+1+\underbrace{0+0+0+0+0+\cdots}_{\text{These are zeroes, thus don't
contribute to the sum}}\\
&  =6+3+1=10,
\end{align*}
which is a well-defined (finite) value. More generally, for any prime $p$ and
any $n\in\mathbb{N}$, the sum $\left\lfloor \dfrac{n}{p^{1}}\right\rfloor
+\left\lfloor \dfrac{n}{p^{2}}\right\rfloor +\left\lfloor \dfrac{n}{p^{3}%
}\right\rfloor +\cdots$ has only finitely many nonzero addends (because for
every $i\geq n$, we have $p^{i}\geq p^{n}>n$ and thus $0\leq\dfrac{n}{p^{i}%
}<1$, so that $\left\lfloor \dfrac{n}{p^{i}}\right\rfloor =0$), and thus
becomes a finite sum once we discard all its addends that are zero; but a
finite sum obviously has a well-defined value.

Moreover, for every positive integer $d$, we have $\left\lfloor \dfrac{n}%
{d}\right\rfloor =n//d$ (by Proposition \ref{prop.ent.quo-rem-floor}). Thus,
the two infinite sums%
\begin{align*}
&  \left\lfloor \dfrac{n}{p^{1}}\right\rfloor +\left\lfloor \dfrac{n}{p^{2}%
}\right\rfloor +\left\lfloor \dfrac{n}{p^{3}}\right\rfloor +\cdots
\ \ \ \ \ \ \ \ \ \ \text{and}\\
&  \left(  n//p^{1}\right)  +\left(  n//p^{2}\right)  +\left(  n//p^{3}%
\right)  +\cdots
\end{align*}
are equal.

It remains to prove that these two sums equal $v_{p}\left(  n!\right)  $. In
other words, we must prove that%
\begin{equation}
v_{p}\left(  n!\right)  =\left\lfloor \dfrac{n}{p^{1}}\right\rfloor
+\left\lfloor \dfrac{n}{p^{2}}\right\rfloor +\left\lfloor \dfrac{n}{p^{3}%
}\right\rfloor +\cdots. \label{pf.thm.ent.vp.polignac.mgoal}%
\end{equation}

We can prove this by induction on $n$:

The \textit{base case} ($n=0$) boils down to $0=0+0+0+\cdots$, which is true.

For the \textit{induction step}, we proceed from $n-1$ to $n$. So we fix a
positive integer $n$, and we assume (as our induction hypothesis) that%
\begin{equation}
v_{p}\left(  \left(  n-1\right)  !\right)  =\left\lfloor \dfrac{n-1}{p^{1}%
}\right\rfloor +\left\lfloor \dfrac{n-1}{p^{2}}\right\rfloor +\left\lfloor
\dfrac{n-1}{p^{3}}\right\rfloor +\cdots, \label{pf.thm.ent.vp.polignac.IH}%
\end{equation}
and we set out to prove that%
\begin{equation}
v_{p}\left(  n!\right)  =\left\lfloor \dfrac{n}{p^{1}}\right\rfloor
+\left\lfloor \dfrac{n}{p^{2}}\right\rfloor +\left\lfloor \dfrac{n}{p^{3}%
}\right\rfloor +\cdots. \label{pf.thm.ent.vp.polignac.IG}%
\end{equation}

We first compare the left hand sides: Let $k=v_{p}\left(  n\right)  $. We know
that $n!=\left(  n-1\right)  !\cdot n$, and therefore%
\begin{align*}
v_{p}\left(  n!\right)   &  =v_{p}\left(  \left(  n-1\right)  !\cdot n\right)
\\
&  =v_{p}\left(  \left(  n-1\right)  !\right)  +\underbrace{v_{p}\left(
n\right)  }_{=k}\ \ \ \ \ \ \ \ \ \ \left(  \text{by Theorem
\ref{thm.ent.vp.vpab} \textbf{(a)}}\right) \\
&  =v_{p}\left(  \left(  n-1\right)  !\right)  +k.
\end{align*}
In other words, the LHS\footnote{The word \textquotedblleft
LHS\textquotedblright\ means \textquotedblleft left hand
side\textquotedblright.
\par
The word \textquotedblleft RHS\textquotedblright\ means \textquotedblleft
right hand side\textquotedblright.} of (\ref{pf.thm.ent.vp.polignac.IG})
equals the LHS of (\ref{pf.thm.ent.vp.polignac.IH}) plus $k$.

Now, we shall show that the RHSs of the two equations differ by $k$ as well.

For each $i\in\left\{  1,2,\ldots,k\right\}  $, we have $p^{i}\mid p^{k}\mid
n$ (since $k=v_{p}\left(  n\right)  $) and therefore%
\[
\left\lfloor \dfrac{n}{p^{i}}\right\rfloor =\left\lfloor \dfrac{n-1}{p^{i}%
}\right\rfloor +1\ \ \ \ \ \ \ \ \ \ \left(  \text{by Corollary
\ref{cor.ent.quo-rem.floor-bday} \textbf{(a)}, applied to }d=p^{i}\right)  .
\]
On the other hand, for each $i\in\left\{  k+1,k+2,k+3,\ldots\right\}  $, we
have $p^{i}\nmid n$ (since $i>k=v_{p}\left(  n\right)  $) and thus%
\[
\left\lfloor \dfrac{n}{p^{i}}\right\rfloor =\left\lfloor \dfrac{n-1}{p^{i}%
}\right\rfloor \ \ \ \ \ \ \ \ \ \ \left(  \text{by Corollary
\ref{cor.ent.quo-rem.floor-bday} \textbf{(b)}, applied to }d=p^{i}\right)  .
\]
These two equalities together yield%
\begin{align*}
&  \left\lfloor \dfrac{n}{p^{1}}\right\rfloor +\left\lfloor \dfrac{n}{p^{2}%
}\right\rfloor +\left\lfloor \dfrac{n}{p^{3}}\right\rfloor +\cdots\\
&  =\left(  \left\lfloor \dfrac{n-1}{p^{1}}\right\rfloor +1\right)  +\left(
\left\lfloor \dfrac{n-1}{p^{2}}\right\rfloor +1\right)  +\cdots+\left(
\left\lfloor \dfrac{n-1}{p^{k}}\right\rfloor +1\right) \\
&  \ \ \ \ \ \ \ \ \ \ +\left\lfloor \dfrac{n-1}{p^{k+1}}\right\rfloor
+\left\lfloor \dfrac{n-1}{p^{k+2}}\right\rfloor +\left\lfloor \dfrac
{n-1}{p^{k+3}}\right\rfloor +\cdots\\
&  =\left(  \left\lfloor \dfrac{n-1}{p^{1}}\right\rfloor +\left\lfloor
\dfrac{n-1}{p^{2}}\right\rfloor +\left\lfloor \dfrac{n-1}{p^{3}}\right\rfloor
+\cdots\right)  +k.
\end{align*}
In other words, the RHS of (\ref{pf.thm.ent.vp.polignac.IG}) equals the RHS of
(\ref{pf.thm.ent.vp.polignac.IH}) plus $k$.

But previously, we have shown the same for the LHSs. Thus, the equality
(\ref{pf.thm.ent.vp.polignac.IG}) is just the equality
(\ref{pf.thm.ent.vp.polignac.IH}) with each side increased by $k$. Since
(\ref{pf.thm.ent.vp.polignac.IH}) holds (by the induction hypothesis), it thus
follows that (\ref{pf.thm.ent.vp.polignac.IG}) also holds. In other words,%
\[
v_{p}\left(  n!\right)  =\left\lfloor \dfrac{n}{p^{1}}\right\rfloor
+\left\lfloor \dfrac{n}{p^{2}}\right\rfloor +\left\lfloor \dfrac{n}{p^{3}%
}\right\rfloor +\cdots.
\]
But this completes the induction step, and thus Theorem
\ref{thm.ent.vp.polignac} is proven.

(For another proof of Theorem \ref{thm.ent.vp.polignac}, see \cite[Exercise
2.17.2 \textbf{(c)}]{19s} or \cite[Theorem 5.3.1]{21f5}.)
\end{proof}

Theorem \ref{thm.ent.vp.polignac} is known as \textbf{de Polignac's formula}
or \textbf{Legendre's formula}. Various uses of this formula can be found in
\cite{21f5}.

\begin{noncompile}
\begin{exercise}
How many $0$'s does the number $50!$ end with when written in decimal notation?
\end{exercise}

\begin{proof}
[Solution idea.]TODO: This is to be completed.

This is asking for $v_{10}\left(  50!\right)  $, but unfortunately $10$ is not
a prime, so our formulas above do not apply.

Nevertheless, de Polignac's formula can be applied. Indeed, $10=2\cdot5$, and
the numbers $2$ and $5$ are two distinct primes (hence coprime). This easily
yields that $2^{k}$ and $5^{k}$ are coprime for any $k\in\mathbb{N}$
(exercise). Thus, for any $k\in\mathbb{N}$, we have the equivalence%
\begin{align*}
\left(  10^{k}\mid n!\right)  \  &  \Longleftrightarrow\ \left(  2^{k}\mid
n!\text{ and }5^{k}\mid n!\right) \\
&  \Longleftrightarrow\ \left(  v_{2}\left(  n!\right)  \geq k\text{ and
}v_{5}\left(  n!\right)  \geq k\right)  .
\end{align*}
So the number of $0$'s at the end of $n!$ is $\min\left\{  v_{2}\left(
n!\right)  ,\ v_{5}\left(  n!\right)  \right\}  $ (because this number of
$0$'s is the largest $k\in\mathbb{N}$ such that $10^{k}\mid n!$). Thus, if we
can find $v_{2}\left(  n!\right)  $ and $v_{5}\left(  n!\right)  $, then we
can also find this number of $0$'s.

For $n=50$, we can do this pretty easily using de Polignac's formula:%
\begin{align*}
v_{2}\left(  50!\right)   &  =\left\lfloor \dfrac{50}{2^{1}}\right\rfloor
+\left\lfloor \dfrac{50}{2^{2}}\right\rfloor +\left\lfloor \dfrac{50}{2^{3}%
}\right\rfloor +\cdots\\
&  =25+12+6+3+1+0+0+0+\cdots\\
&  =47
\end{align*}
and%
\begin{align*}
v_{5}\left(  50!\right)   &  =\left\lfloor \dfrac{50}{5^{1}}\right\rfloor
+\left\lfloor \dfrac{50}{5^{2}}\right\rfloor +\left\lfloor \dfrac{50}{5^{3}%
}\right\rfloor +\cdots\\
&  =10+2+0+0+\cdots\\
&  =12,
\end{align*}
so%
\[
\min\left\{  v_{2}\left(  50!\right)  ,\ v_{5}\left(  50!\right)  \right\}
=\min\left\{  47,12\right\}  =12.
\]
So the number $50!$ ends with $12$ zeroes. And indeed,%
\[
50!=\allowbreak30\,\allowbreak414\,093\,201\,\allowbreak
713\,378\,043\,\allowbreak612\,608\,166\,\allowbreak064\,768\,844\,\allowbreak
377\,641\,568\,\allowbreak960\,512\,000\,\allowbreak000\,000\,000\allowbreak.
\]

\end{proof}
\end{noncompile}

\subsubsection{\label{subsec.ent.primefac}Prime factorization}

We are now ready to prove one of the most important properties of primes: the
fact that every positive integer can be uniquely decomposed into a product of
some primes. For instance,%
\[
200=2\cdot100=2\cdot2\cdot50=2\cdot2\cdot5\cdot10=\underbrace{2\cdot
2\cdot5\cdot2\cdot5}_{\text{a product of primes}}.
\]
The word \textquotedblleft uniquely\textquotedblright\ means here that any two
ways of decomposing a given positive integer $n$ into a product of primes are
equal up to reordering the factors. For example, we can also decompose $200$
as $5\cdot2\cdot2\cdot5\cdot2$, but this is the same product with the factors
in a different order.

Let us state this fact in full generality. First, we introduce a name for
these decompositions:

\begin{definition}
\label{def.ent.primefac}Let $n$ be a positive integer. A \textbf{prime
factorization} of $n$ means a finite list $\left(  p_{1},p_{2},\ldots
,p_{k}\right)  $ of primes (not necessarily distinct) such that%
\[
n=p_{1}p_{2}\cdots p_{k}.
\]

\end{definition}

Thus, $\left(  2,2,5,2,5\right)  $ and $\left(  5,2,2,5,2\right)  $ are prime
factorizations of $200$. Another such is $\left(  2,2,2,5,5\right)  $. There
are more\footnote{In Remark \ref{rmk.enum.primefac}, we will see how many.},
but all of them contain the number $2$ thrice and the number $5$ twice (and no
other numbers), just as we said.

Let us state this as a general claim:\footnote{Note that the integer $1$ has a
prime factorization, too! This factorization is the empty list $\left(
{}\right)  $, and it works because the empty product is $1$ (by definition).}

\begin{theorem}
[Fundamental Theorem of Arithmetic]\label{thm.ent.fta}Let $n$ be a positive
integer. Then: \medskip

\textbf{(a)} There exists a prime factorization of $n$. \medskip

\textbf{(b)} This prime factorization is unique up to reordering its entries.
In other words, if $\left(  p_{1},p_{2},\ldots,p_{k}\right)  $ and $\left(
q_{1},q_{2},\ldots,q_{\ell}\right)  $ are two prime factorizations of $n$,
then $\left(  q_{1},q_{2},\ldots,q_{\ell}\right)  $ can be obtained from
$\left(  p_{1},p_{2},\ldots,p_{k}\right)  $ by reordering the entries.
\medskip

\textbf{(c)} Let $\left(  p_{1},p_{2},\ldots,p_{k}\right)  $ be a prime
factorization of $n$. Let $p$ be any prime. Then, the number of times that $p$
appears in the list $\left(  p_{1},p_{2},\ldots,p_{k}\right)  $ (in other
words, the number of $i\in\left\{  1,2,\ldots,k\right\}  $ satisfying
$p_{i}=p$) is $v_{p}\left(  n\right)  $.
\end{theorem}

\begin{proof}
\textbf{(a)} This is Theorem \ref{thm.ent.primefac.exists}. \medskip

\textbf{(c)} By the definition of a prime factorization, we have $n=p_{1}%
p_{2}\cdots p_{k}$. Thus,%
\begin{align}
v_{p}\left(  n\right)   &  =v_{p}\left(  p_{1}p_{2}\cdots p_{k}\right)
\nonumber\\
&  =v_{p}\left(  p_{1}\right)  +v_{p}\left(  p_{2}\right)  +\cdots
+v_{p}\left(  p_{k}\right)  \label{pf.thm.ent.fta.c.1}%
\end{align}
(by Corollary \ref{cor.ent.vp.prodak}).

The right hand side of this equality is a sum of $k$ addends. Each of these
addends has the form $v_{p}\left(  p_{i}\right)  $ for some $i\in\left\{
1,2,\ldots,k\right\}  $. Each such addend $v_{p}\left(  p_{i}\right)  $ equals
$1$ if $p_{i}=p$ (by Theorem \ref{thm.ent.vp.vpab} \textbf{(d)}) and equals
$0$ if $p_{i}\neq p$ (by Theorem \ref{thm.ent.vp.vpab} \textbf{(e)}).

Thus, our sum $v_{p}\left(  p_{1}\right)  +v_{p}\left(  p_{2}\right)
+\cdots+v_{p}\left(  p_{k}\right)  $ has an addend equal to $1$ for each
$i\in\left\{  1,2,\ldots,k\right\}  $ that satisfies $p_{i}=p$, and an addend
equal to $0$ for each $i$ that doesn't.

Obviously, the addends that are equal to $0$ do not affect the sum. Hence, the
sum equals the number of addends equal to $1$. In other words, the sum equals
the number of $i\in\left\{  1,2,\ldots,k\right\}  $ that satisfy $p_{i}=p$.

In view of (\ref{pf.thm.ent.fta.c.1}), we can restate this as follows:
$v_{p}\left(  n\right)  $ equals the number of $i\in\left\{  1,2,\ldots
,k\right\}  $ that satisfy $p_{i}=p$. In other words, $v_{p}\left(  n\right)
$ equals the number of times that $p$ appears in the list $\left(  p_{1}%
,p_{2},\ldots,p_{k}\right)  $. This proves Theorem \ref{thm.ent.fta}
\textbf{(c)}. \medskip

\textbf{(b)} This follows easily from part \textbf{(c)}. Namely:

Let $\left(  p_{1},p_{2},\ldots,p_{k}\right)  $ and $\left(  q_{1}%
,q_{2},\ldots,q_{\ell}\right)  $ be two prime factorizations of $n$. We must
prove that $\left(  q_{1},q_{2},\ldots,q_{\ell}\right)  $ can be obtained from
$\left(  p_{1},p_{2},\ldots,p_{k}\right)  $ by reordering the entries.

Each prime $p$ appears $v_{p}\left(  n\right)  $ times in the list $\left(
p_{1},p_{2},\ldots,p_{k}\right)  $ (by part \textbf{(c)}), and appears
$v_{p}\left(  n\right)  $ times in the list $\left(  q_{1},q_{2}%
,\ldots,q_{\ell}\right)  $ (similarly). Thus, each prime $p$ appears the same
number of times in either list. Since both lists consist of primes, this shows
that the two lists contain the same numbers the same number of times.
Therefore, $\left(  q_{1},q_{2},\ldots,q_{\ell}\right)  $ can be obtained from
$\left(  p_{1},p_{2},\ldots,p_{k}\right)  $ by reordering the entries. This
proves Theorem \ref{thm.ent.fta} \textbf{(b)}.

(We have used the intuitively obvious fact that if two lists of numbers
contain the same numbers the same number of times, then one can be obtained
from the other by reordering. You are free to trust your intuition on this
one; for a formal proof, see \cite[Lemma 2.13.20]{19s}.)
\end{proof}

\bigskip

Theorem \ref{thm.ent.fta} \textbf{(a)} shows that every positive integer $n$
has a prime factorization. Finding this prime factorization is
\href{https://en.wikipedia.org/wiki/Integer_factorization}{a classical hard
computational problem}. (Quite a few encryption standards rely on its hardness.)

\subsubsection{Applications}

Prime factorizations can be rather useful. The next few exercises provide some examples:

\begin{exercise}
\label{exe.pval.n|m}Let $n$ and $m$ be integers. Prove that $n\mid m$ if and
only if each prime $p$ satisfies $v_{p}\left(  n\right)  \leq v_{p}\left(
m\right)  $. \medskip

[\textbf{Hint:} For the \textquotedblleft if\textquotedblright\ direction,
start by picking prime factorizations of $n$ and $m$ (after first handling the
cases $n=0$ and $m=0$ and getting rid of signs).]
\end{exercise}

\begin{exercise}
\label{exe.pval.e2}Let $k$ be a positive integer. Let $w$ be a rational number
such that $w^{k}$ is an integer. Prove that $w$ is an integer. \medskip

[\textbf{Hint:} Use Exercise \ref{exe.pval.n|m}.]
\end{exercise}

Note that Exercise \ref{exe.pval.e2} shows that (for example) the number
$\sqrt{2}$ is irrational, because it is not an integer (since $1<\sqrt{2}<2$)
but its $2$-nd power $\left(  \sqrt{2}\right)  ^{2}=2$ is an integer.
Likewise, $\sqrt[3]{2}$ and $\sqrt{5}$ and $\sqrt[7]{9}$ and many similar
surds are irrational.

\begin{exercise}
\label{exe.pval.c1}Solve Exercise \ref{exe.ent.prod-is-kth-power} again using
$p$-valuations.
\end{exercise}

\subsection{\label{sec.ent.lcm}Least common multiples}

In Section \ref{sec.ent.gcd}, we have studied greatest common divisors in some
detail. Let me now briefly discuss least common multiples: a kind of
counterpart to greatest common divisors. The greatest common divisor of two
positive integers $a$ and $b$ is usually smaller than both $a$ and $b$; in
contrast, the least common multiple is usually larger than both.

\begin{definition}
Let $a$ and $b$ be two integers. \medskip

\textbf{(a)} The \textbf{common multiples} of $a$ and $b$ are the integers
that are divisible by $a$ and simultaneously divisible by $b$. \medskip

\textbf{(b)} The \textbf{least common multiple} (aka the \textbf{lowest common
multiple}, or just the \textbf{lcm}) of $a$ and $b$ is defined as follows:

\begin{itemize}
\item If $a$ and $b$ are nonzero, then it is the smallest positive common
multiple of $a$ and $b$.

\item Otherwise, it is $0$.
\end{itemize}

It is denoted by $\operatorname{lcm}\left(  a,b\right)  $.
\end{definition}

Some examples:

\begin{itemize}
\item We have $\operatorname{lcm}\left(  3,4\right)  =12$.

\item We have $\operatorname{lcm}\left(  6,4\right)  =12$.

\item We have $\operatorname{lcm}\left(  6,8\right)  =24$.

\item We have $\operatorname{lcm}\left(  2,4\right)  =4$.

\item We have $\operatorname{lcm}\left(  0,5\right)  =0$.

\item We have $\operatorname{lcm}\left(  -2,3\right)  =6$.
\end{itemize}

Note that the lcm of two positive integers is a fairly well-known concept:
When you bring two fractions (of integers) to their lowest common denominator,
this lowest common denominator is actually the lcm of the denominators of the fractions.

Here are some properties of lcms:

\begin{theorem}
\label{thm.ent.lcm.varia}Let $a$ and $b$ be two integers. Then: \medskip

\textbf{(a)} The lcm of $a$ and $b$ exists. \medskip

\textbf{(b)} We have $\operatorname{lcm}\left(  a,b\right)  \in\mathbb{N}$.
\medskip

\textbf{(c)} We have $\operatorname{lcm}\left(  a,b\right)
=\operatorname{lcm}\left(  b,a\right)  $. \medskip

\textbf{(d)} We have $a\mid\operatorname{lcm}\left(  a,b\right)  $ and
$b\mid\operatorname{lcm}\left(  a,b\right)  $. \medskip

\textbf{(e)} We have $\operatorname{lcm}\left(  -a,b\right)
=\operatorname{lcm}\left(  a,b\right)  $ and $\operatorname{lcm}\left(
a,-b\right)  =\operatorname{lcm}\left(  a,b\right)  $.
\end{theorem}

\begin{proof}
[Proof sketch.]Easy consequences of the definitions. (For part \textbf{(a)},
observe that two nonzero integers $a$ and $b$ have at least one positive
common multiple -- namely, $\left\vert ab\right\vert $.)
\end{proof}

Here is a counterpart to the universal property of the gcd (Theorem
\ref{thm.ent.gcd.uniprop}):

\begin{theorem}
[universal property of the lcm]\label{thm.ent.lcm.uniprop}Let $a,b,m\in
\mathbb{Z}$. Then, we have the equivalence%
\[
\left(  a\mid m\text{ and }b\mid m\right)  \ \Longleftrightarrow\ \left(
\operatorname{lcm}\left(  a,b\right)  \mid m\right)  .
\]

\end{theorem}

In other words, the common multiples of two integers $a$ and $b$ are precisely
the multiples of $\operatorname{lcm}\left(  a,b\right)  $.

\begin{proof}
[Proof sketch.](See \cite[Theorem 2.11.7]{19s} for a detailed proof.) \medskip

$\Longleftarrow:$ If $\operatorname{lcm}\left(  a,b\right)  \mid m$, then
$a\mid m$ (since Theorem \ref{thm.ent.lcm.varia} \textbf{(d)} yields
$a\mid\operatorname{lcm}\left(  a,b\right)  \mid m$) and $b\mid m$
(similarly). Thus, the \textquotedblleft$\Longleftarrow$\textquotedblright%
\ direction of the desired equivalence is proved. \medskip

$\Longrightarrow:$ Assume that $a\mid m$ and $b\mid m$. We must show that
$\operatorname{lcm}\left(  a,b\right)  \mid m$.

If one of $a$ and $b$ is $0$, then this is easy (in fact, let's say that
$a=0$; then, $0=a\mid m$, thus $m=0$, and therefore $\operatorname{lcm}\left(
a,b\right)  \mid0=m$). Hence, we need only to consider the case when $a$ and
$b$ are nonzero.

In this case, set $\ell=\operatorname{lcm}\left(  a,b\right)  $. Recall that
$\ell$ is defined as the smallest positive common multiple of $a$ and $b$.
Hence, $\ell$ is a positive integer and is a multiple of $a$ and of $b$. Let
$q$ and $r$ be the quotient and the remainder of the division of $m$ by $\ell
$. Thus,%
\[
q\in\mathbb{Z}\ \ \ \ \ \ \ \ \ \ \text{and}\ \ \ \ \ \ \ \ \ \ r\in\left\{
0,1,\ldots,\ell-1\right\}  \ \ \ \ \ \ \ \ \ \ \text{and}%
\ \ \ \ \ \ \ \ \ \ m=q\ell+r
\]
(by the definition of quotient and remainder). From $r\in\left\{
0,1,\ldots,\ell-1\right\}  $, we obtain $r<\ell$.

From $m=q\ell+r$, we obtain $r=m-q\ell$. Since both $m$ and $\ell$ are
multiples of $a$, we thus conclude that $r$ is a multiple of $a$ as well.
Similarly, $r$ is a multiple of $b$. Thus, $r$ is a common multiple of $a$ and
$b$. But $\ell$ is the smallest positive common multiple of $a$ and $b$. If
$r$ was positive, then $r$ would contradict this minimality (because $r<\ell
$). Hence, $r$ cannot be positive. Since $r\in\left\{  0,1,\ldots
,\ell-1\right\}  $, we conclude that $r$ must be $0$. Hence, $m=q\ell
+\underbrace{r}_{=0}=q\ell$, so that $\ell\mid m$. In other words,
$\operatorname{lcm}\left(  a,b\right)  \mid m$ (since $\ell=\operatorname{lcm}%
\left(  a,b\right)  $). This proves the \textquotedblleft$\Longrightarrow
$\textquotedblright\ direction of the desired equivalence.
\end{proof}

\bigskip

The gcd and the lcm of two integers are connected to each other by the
following formula:

\begin{theorem}
\label{thm.ent.lcm.gcd}Let $a$ and $b$ be two integers. Then,%
\[
\gcd\left(  a,b\right)  \cdot\operatorname{lcm}\left(  a,b\right)  =\left\vert
ab\right\vert .
\]

\end{theorem}

\begin{proof}
[Proof sketch.](See \cite[Theorem 2.11.6]{19s} for a detailed proof.)

The case when $a$ or $b$ is $0$ is easy (since both $\operatorname{lcm}\left(
a,b\right)  $ and $ab$ are $0$ in this case). Thus, we can assume for the rest
of this proof that both $a$ and $b$ are nonzero.

Then, $\dfrac{ab}{\gcd\left(  a,b\right)  }$ is an integer (since $\gcd\left(
a,b\right)  \mid a\mid ab$) and is a common multiple of $a$ and $b$ (indeed,
$\gcd\left(  a,b\right)  \mid b$ yields $\dfrac{b}{\gcd\left(  a,b\right)
}\in\mathbb{Z}$, and therefore $\dfrac{ab}{\gcd\left(  a,b\right)  }%
=a\cdot\dfrac{b}{\gcd\left(  a,b\right)  }$ is a multiple of $a$; an analogous
computation shows that $\dfrac{ab}{\gcd\left(  a,b\right)  }$ is a multiple of
$b$ as well). That is, $a\mid\dfrac{ab}{\gcd\left(  a,b\right)  }$ and
$b\mid\dfrac{ab}{\gcd\left(  a,b\right)  }$. By Theorem
\ref{thm.ent.lcm.uniprop}, this entails that $\operatorname{lcm}\left(
a,b\right)  \mid\dfrac{ab}{\gcd\left(  a,b\right)  }$, so that $\gcd\left(
a,b\right)  \cdot\operatorname{lcm}\left(  a,b\right)  \mid ab$.

On the other hand, Theorem \ref{thm.ent.lcm.uniprop} shows that
$\operatorname{lcm}\left(  a,b\right)  \mid ab$ (since $a\mid ab$ and $b\mid
ab$), so that $\dfrac{ab}{\operatorname{lcm}\left(  a,b\right)  }$ is an
integer (note that $\operatorname{lcm}\left(  a,b\right)  \mid ab\neq0$ yields
$\operatorname{lcm}\left(  a,b\right)  \neq0$, so this fraction $\dfrac
{ab}{\operatorname{lcm}\left(  a,b\right)  }$ is well-defined). This integer
$\dfrac{ab}{\operatorname{lcm}\left(  a,b\right)  }$ divides $a$ (because
$a/\dfrac{ab}{\operatorname{lcm}\left(  a,b\right)  }=\dfrac
{\operatorname{lcm}\left(  a,b\right)  }{b}\in\mathbb{Z}$ (since
$b\mid\operatorname{lcm}\left(  a,b\right)  $)) and divides $b$ (for similar
reasons); therefore, it divides $\gcd\left(  a,b\right)  $ (by Theorem
\ref{thm.ent.gcd.uniprop}, because it divides each of $a$ and $b$). In other
words, $ab\mid\gcd\left(  a,b\right)  \cdot\operatorname{lcm}\left(
a,b\right)  $.

Now, recall that two integers $x$ and $y$ that satisfy $x\mid y$ and $y\mid x$
must satisfy $\left\vert x\right\vert =\left\vert y\right\vert $. Hence, from
$\gcd\left(  a,b\right)  \cdot\operatorname{lcm}\left(  a,b\right)  \mid ab$
and $ab\mid\gcd\left(  a,b\right)  \cdot\operatorname{lcm}\left(  a,b\right)
$, we obtain $\left\vert \gcd\left(  a,b\right)  \cdot\operatorname{lcm}%
\left(  a,b\right)  \right\vert =\left\vert ab\right\vert $. Since
$\gcd\left(  a,b\right)  $ and $\operatorname{lcm}\left(  a,b\right)  $ are
nonnegative, we can drop the absolute values on the left-hand side, so this
simplifies to $\gcd\left(  a,b\right)  \cdot\operatorname{lcm}\left(
a,b\right)  =\left\vert ab\right\vert $.
\end{proof}

\bigskip

Both gcds and lcms have easily computable $p$-valuations:

\begin{theorem}
\label{thm.ent.vp.gcd-lcm}Let $p$ be a prime. Let $a$ and $b$ be two integers.
Then,%
\begin{align*}
v_{p}\left(  \gcd\left(  a,b\right)  \right)   &  =\min\left\{  v_{p}\left(
a\right)  ,\ v_{p}\left(  b\right)  \right\}  \ \ \ \ \ \ \ \ \ \ \text{and}\\
v_{p}\left(  \operatorname{lcm}\left(  a,b\right)  \right)   &  =\max\left\{
v_{p}\left(  a\right)  ,\ v_{p}\left(  b\right)  \right\}  .
\end{align*}

\end{theorem}

\begin{proof}
[Proof sketch.]This is \cite[Exercise 2.13.10]{19s} and is also a particular
case of \cite[Lecture 5, Proposition 5.2.15]{23f}. Anyway, the proof is a nice
exercise in using the universal properties of the gcd and the lcm (and the
definition of $p$-valuation), so you should do it yourself.
\end{proof}

Theorem \ref{thm.ent.vp.gcd-lcm} gives an easy way to compute $\gcd\left(
a,b\right)  $ and $\operatorname{lcm}\left(  a,b\right)  $ if you know prime
factorizations of two positive integers $a$ and $b$. For example, knowing that
$18=2\cdot3^{2}$ and $12=2^{2}\cdot3$, we obtain%
\begin{align*}
\gcd\left(  18,\ 12\right)   &  =2\cdot3=6\ \ \ \ \ \ \ \ \ \ \text{and}\\
\operatorname{lcm}\left(  18,\ 12\right)   &  =2^{2}\cdot3^{2}=36.
\end{align*}

If you don't know the prime factorizations of $a$ and $b$, the quickest way to
find $\operatorname{lcm}\left(  a,b\right)  $ is by using the Euclidean
algorithm to find $\gcd\left(  a,b\right)  $ first, and then solving the
equality $\gcd\left(  a,b\right)  \cdot\operatorname{lcm}\left(  a,b\right)
=\left\vert ab\right\vert $ for $\operatorname{lcm}\left(  a,b\right)  $. This
gives\footnote{Here we are assuming that $a$ and $b$ are nonzero. If $a$ or
$b$ is $0$, then $\operatorname{lcm}\left(  a,b\right)  $ is just $0$.}%
\[
\operatorname{lcm}\left(  a,b\right)  =\dfrac{\left\vert ab\right\vert }%
{\gcd\left(  a,b\right)  }=\left\vert \dfrac{a}{\gcd\left(  a,b\right)  }\cdot
b\right\vert .
\]

\bigskip

Gcds and lcms can be defined for multiple numbers (not just for two numbers).
Their properties are mostly analogous to the case of two numbers, with some
exceptions (i.e., the formula $\gcd\left(  a,b\right)  \cdot\operatorname{lcm}%
\left(  a,b\right)  =\left\vert ab\right\vert $ does not generalize to
$\gcd\left(  a,b,c\right)  \cdot\operatorname{lcm}\left(  a,b,c\right)
=\left\vert abc\right\vert $, but rather to $\gcd\left(  a,b,c\right)
\cdot\operatorname{lcm}\left(  bc,ca,ab\right)  =\left\vert abc\right\vert $).
See \cite[\S 2.11]{19s} for more details.

\subsection{\label{sec.ent.sylv}Sylvester's $xa+yb$ theorem (or the Chicken
McNugget theorem)}

We come to a rather curious (although not overly important) topic in
elementary number theory: the $\mathbb{N}$-linear combinations of two positive integers.

For this entire section, we let $a$ and $b$ be two positive integers.

\begin{definition}
\textbf{(a)} A $\mathbb{Z}$\textbf{-linear combination} (short: $\mathbb{Z}%
$\textbf{-LC}) of $a$ and $b$ will mean a number of the form
\[
xa+yb\ \ \ \ \ \ \ \ \ \ \text{with }x,y\in\mathbb{Z}.
\]
In other words, it means a number of cents that you can pay with $a$-cent
coins and $b$-cent coins if you can get change. \medskip

\textbf{(b)} An $\mathbb{N}$\textbf{-linear combination} (short: $\mathbb{N}%
$\textbf{-LC}) of $a$ and $b$ will mean a number of the form
\[
xa+yb\ \ \ \ \ \ \ \ \ \ \text{with }x,y\in\mathbb{N}.
\]
In other words, it means a number of cents that you can pay with $a$-cent
coins and $b$-cent coins without getting change.
\end{definition}

Thus, Proposition \ref{prop.int.coins-35} is saying that any integer $n\geq8$
is an $\mathbb{N}$-LC of $3$ and $5$. Moreover, as we saw just above that
proposition, the numbers $0,\ 3,\ 5,\ 6$ are $\mathbb{N}$-LCs of $3$ and $5$
as well, whereas the numbers $1,\ 2,\ 4,\ 7$ are not. Thus the complete list
of all $\mathbb{N}$-LCs of $3$ and $5$ is
\[
0,\ 3,\ 5,\ 6,\ \underbrace{8,\ 9,\ 10,\ \ldots}_{\text{all integers }n\geq
8}.
\]

This should prompt us to study $\mathbb{N}$-LCs of $a$ and $b$ in the general
case. We shall begin with the $\mathbb{Z}$-LCs, however, since they are much
easier to describe.

Note that the $\mathbb{N}$-LCs of $a$ and $b$ are always $\geq0$ (because if
$x,y\in\mathbb{N}$, then $\underbrace{x}_{\geq0}\underbrace{a}_{>0}%
+\underbrace{y}_{\geq0}\underbrace{b}_{>0}\geq0$), whereas the $\mathbb{Z}%
$-LCs of $a$ and $b$ can have any sign.

Clearly, any $\mathbb{N}$-LC of $a$ and $b$ is a $\mathbb{Z}$-LC of $a$ and
$b$. However, a $\mathbb{Z}$-LC of $a$ and $b$ doesn't have to be an
$\mathbb{N}$-LC of $a$ and $b$, even if it is $\geq0$. For example, $1$ is a
$\mathbb{Z}$-LC of $3$ and $5$ (since $1=2\cdot3+\left(  -1\right)  \cdot5$),
but not an $\mathbb{N}$-LC of $3$ and $5$.

\bigskip

We can easily describe the $\mathbb{Z}$-LCs of $a$ and $b$:

\begin{proposition}
\label{prop.ent.sylv2Z}The $\mathbb{Z}$-LCs of $a$ and $b$ are exactly the
multiples of $\gcd\left(  a,b\right)  $.
\end{proposition}

\begin{fineprint}
\begin{proof}
We must prove the following two claims:

\begin{statement}
\textit{Claim 1:} Any $\mathbb{Z}$-LC of $a$ and $b$ is a multiple of
$\gcd\left(  a,b\right)  $.
\end{statement}

\begin{statement}
\textit{Claim 2:} Any multiple of $\gcd\left(  a,b\right)  $ is a $\mathbb{Z}%
$-LC of $a$ and $b$.
\end{statement}

But both claims are easy:

\begin{proof}
[Proof of Claim 1.]Let $n$ be a $\mathbb{Z}$-LC of $a$ and $b$. We must show
that $n$ is a multiple of $\gcd\left(  a,b\right)  $.

Indeed, $n$ is a $\mathbb{Z}$-LC of $a$ and $b$, and thus has the form
$n=xa+yb$ for some $x,y\in\mathbb{Z}$. Consider these $x,y$. We have
$\gcd\left(  a,b\right)  \mid a\mid xa$ and $\gcd\left(  a,b\right)  \mid
b\mid yb$. In other words, both numbers $xa$ and $yb$ are multiples of
$\gcd\left(  a,b\right)  $. Hence, their sum $xa+yb$ is a multiple of
$\gcd\left(  a,b\right)  $ as well. In other words, $n$ is a multiple of
$\gcd\left(  a,b\right)  $ (since $n=xa+yb$). This proves Claim 1.
\end{proof}

\begin{proof}
[Proof of Claim 2.]Let $n$ be a multiple of $\gcd\left(  a,b\right)  $. We
must prove that $n$ is a $\mathbb{Z}$-LC of $a$ and $b$.

Bezout's theorem (Theorem \ref{thm.ent.gcd.bezout}) says that there exist two
integers $x$ and $y$ such that $\gcd\left(  a,b\right)  =xa+yb$. Consider
these $x$ and $y$. However, $n$ is a multiple of $\gcd\left(  a,b\right)  $;
in other words, there exists an integer $c$ such that $n=\gcd\left(
a,b\right)  \cdot c$. Consider this $c$. Now,%
\[
n=\underbrace{\gcd\left(  a,b\right)  }_{=xa+yb}\cdot\,c=\left(  xa+yb\right)
\cdot c=xac+ybc=\left(  cx\right)  a+\left(  cy\right)  b.
\]
This shows that $n$ is a $\mathbb{Z}$-LC of $a$ and $b$ (since $cx$ and $cy$
are integers). This proves Claim 2.
\end{proof}

Combining Claim 1 with Claim 2, we conclude that the $\mathbb{Z}$-LCs of $a$
and $b$ are exactly the multiples of $\gcd\left(  a,b\right)  $. Thus,
Proposition \ref{prop.ent.sylv2Z} is proved.
\end{proof}
\end{fineprint}

\bigskip

Now we move on to the $\mathbb{N}$-LCs. What are they? Can we describe them
any better than by their definition?

Let $g=\gcd\left(  a,b\right)  $. Then, $g$ divides each of $a$ and $b$, so
that the numbers $\dfrac{a}{g}$ and $\dfrac{b}{g}$ are positive integers. We
can simplify our problem by replacing $a$ and $b$ with $\dfrac{a}{g}$ and
$\dfrac{b}{g}$. Clearly, the $\mathbb{N}$-LCs of $a$ and $b$ are just the
$\mathbb{N}$-LCs of $\dfrac{a}{g}$ and $\dfrac{b}{g}$, multiplied by $g$. By
Theorem \ref{thm.ent.gcd.reduce}, the two integers $\dfrac{a}{g}$ and
$\dfrac{b}{g}$ are coprime. Thus, understanding the $\mathbb{N}$-LCs of the
original integers $a$ and $b$ is equivalent to understanding the $\mathbb{N}%
$-LCs of the coprime integers $\dfrac{a}{g}$ and $\dfrac{b}{g}$.

Hence, it suffices to solve our problem in the case when $a$ and $b$ are
coprime. In this case, Proposition \ref{prop.ent.sylv2Z} shows that every
integer is a $\mathbb{Z}$-LC of $a$ and $b$ (since every integer is a multiple
of $1=\gcd\left(  a,b\right)  $). The $\mathbb{N}$-LCs are more interesting.
We have already listed the $\mathbb{N}$-LCs of $3$ and $5$ above; let us now
give a somewhat more complicated example: The $\mathbb{N}$-LCs of $5$ and $9$
are%
\[
0,\ 5,\ 9,\ 10,\ 14,\ 15,\ 18,\ 19,\ 20,\ 23,\ 24,\ 25,\ 27,\ 28,\ 29,\ 30,\ \underbrace{32,\ 33,\ 34,\ \ldots
}_{\text{all integers }n\geq32}.
\]
Note that every integer $n\geq32$ is an $\mathbb{N}$-LC of $5$ and $9$. Among
the first $32$ nonnegative integers $0,1,\ldots,31$, exactly half (that is,
$16$) are $\mathbb{N}$-LCs of $5$ and $9$. A similar phenomenon can be seen in
our above example with $3$ and $5$, except that $32$ is replaced by $8$.

This phenomenon generalizes:

\begin{theorem}
[Sylvester's two-coin theorem, or Chicken McNugget theorem]%
\label{thm.ent.sylv2full}Assume that the two positive integers $a$ and $b$ are
coprime. Then: \medskip

\textbf{(a)} Every integer $n>ab-a-b$ is an $\mathbb{N}$-LC of $a$ and $b$.
\medskip

\textbf{(b)} The number $ab-a-b$ is \textbf{not} an $\mathbb{N}$-LC of $a$ and
$b$. \medskip

\textbf{(c)} Among the first $\left(  a-1\right)  \left(  b-1\right)  $
nonnegative integers $0,\ 1,\ \ldots,\ ab-a-b$, exactly half are $\mathbb{N}%
$-LCs of $a$ and $b$. \medskip

\textbf{(d)} Let $n\in\mathbb{Z}$. Then, exactly one of the two numbers $n$
and $ab-a-b-n$ is an $\mathbb{N}$-LC of $a$ and $b$.
\end{theorem}

This theorem was discovered by J. J. Sylvester in 1884, as a side-product of
his work in invariant theory. Its more recent moniker is due to the McDonald's
Chicken McNuggets, which used to be sold in packs of $9$ or $20$, prompting
mathematicians to wonder what numbers of nuggets could be bought.

The theorem stops short of explicitly answering which of the first $\left(
a-1\right)  \left(  b-1\right)  $ nonnegative integers are $\mathbb{N}$-LCs of
$a$ and $b$. There is no \textquotedblleft easy formula\textquotedblright\ for
this answer. But Theorem \ref{thm.ent.sylv2full} \textbf{(a)} gives you all
the information you need to compute all the $\mathbb{N}$-LCs of $a$ and $b$,
since the first $\left(  a-1\right)  \left(  b-1\right)  $ nonnegative
integers can be checked one by one.

The particular case of Theorem \ref{thm.ent.sylv2full} \textbf{(a)} where
$a=p$ and $b=p+1$ was Exercise \ref{exe.23hw3.4}. \medskip

Before we prove Theorem \ref{thm.ent.sylv2full}, we show a basic lemma:

\begin{lemma}
\label{lem.ent.sylv2lim}Assume that the two positive integers $a$ and $b$ are
coprime. Let $n\in\mathbb{Z}$. Then, there exist two integers $u$ and $v$ such
that $0\leq u\leq b-1$ and $ua+vb=n$.
\end{lemma}

\begin{fineprint}
\begin{proof}
[Proof of Lemma \ref{lem.ent.sylv2lim}.]Bezout's theorem (Theorem
\ref{thm.ent.gcd.bezout}) says that there exist two integers $x$ and $y$ such
that $\gcd\left(  a,b\right)  =xa+yb$. Consider these $x$ and $y$. Thus,
$xa+yb=\gcd\left(  a,b\right)  =1$ (since $a$ and $b$ are coprime).

Recall that $b$ is a positive integer. Thus, division with remainder by $b$ is
well-defined (see Definition \ref{def.ent.quo-rem} for the terminology).

Let $q=\left(  nx\right)  //b$ and $r=\left(  nx\right)  \%b$. In other words,
let $q$ and $r$ be the quotient and the remainder of the division of $nx$ by
$b$. By the definition of quotient and remainder, we thus have%
\[
q\in\mathbb{Z}\ \ \ \ \ \ \ \ \ \ \text{and}\ \ \ \ \ \ \ \ \ \ r\in\left\{
0,1,\ldots,b-1\right\}  \ \ \ \ \ \ \ \ \ \ \text{and}%
\ \ \ \ \ \ \ \ \ \ nx=qb+r.
\]

From $r\in\left\{  0,1,\ldots,b-1\right\}  $, we see that $r$ is an integer
satisfying $0\leq r\leq b-1$.

On the other hand, $nxa+nyb=n\underbrace{\left(  xa+yb\right)  }_{=1}=n$, so
that%
\begin{align*}
n  &  =\underbrace{nx}_{=qb+r}a+nyb=\left(  qb+r\right)  a+nyb\\
&  =qba+ra+nyb=ra+\underbrace{qba+nyb}_{=\left(  qa+ny\right)  b}=ra+\left(
qa+ny\right)  b.
\end{align*}
In other words, $ra+\left(  qa+ny\right)  b=n$.

Altogether, we now know that $r$ and $qa+ny$ are two integers satisfying
$0\leq r\leq b-1$ and $ra+\left(  qa+ny\right)  b=n$. Thus, there exist two
integers $u$ and $v$ such that $0\leq u\leq b-1$ and $ua+vb=n$ (namely, $u=r$
and $v=qa+ny$). This proves Lemma \ref{lem.ent.sylv2lim}.
\end{proof}

\begin{proof}
[Proof of Theorem \ref{thm.ent.sylv2full}.]We shall first prove part
\textbf{(b)} and then part \textbf{(d)}. The other two parts will follow quite
easily from these. \medskip

\textbf{(b)} Assume the contrary. Thus, $ab-a-b$ is an $\mathbb{N}$-LC of $a$
and $b$. In other words, there exist nonnegative integers $x$ and $y$ such
that $ab-a-b=xa+yb$. Consider these $x$ and $y$.

From $ab-a-b=xa+yb$, we obtain $ab=xa+yb+a+b=\left(  x+1\right)  a+\left(
y+1\right)  b=a\left(  x+1\right)  +b\left(  y+1\right)  $. Hence,%
\[
b\left(  y+1\right)  =ab-a\left(  x+1\right)  =a\cdot\underbrace{\left(
b-\left(  x+1\right)  \right)  }_{\text{an integer}}.
\]
This shows that $a\mid b\left(  y+1\right)  $. Thus, the coprime removal
theorem (Theorem \ref{thm.ent.coprime.a|c}) yields that $a\mid y+1$ (since $a$
is coprime to $b$). Therefore, $\dfrac{y+1}{a}$ is an integer (since $a\neq
0$). Since $\underbrace{y}_{\geq0}+\,1\geq1>0$ and $a>0$, this integer
$\dfrac{y+1}{a}$ is furthermore positive, and thus is $\geq1$. In other words,
$y+1\geq a$. Hence, $y\geq a-1$. Now,%
\[
ab-a-b=\underbrace{x}_{\geq0}a+\underbrace{y}_{\geq a-1}b\geq0a+\left(
a-1\right)  b=ab-b.
\]
Subtracting $ab-a-b$ from both sides of this inequality, we obtain $0\geq a$,
which contradicts the positivity of $a$. This contradiction shows that our
assumption was false. Thus, Theorem \ref{thm.ent.sylv2full} \textbf{(b)} is
proved. \medskip

\textbf{(d)} Let $m=ab-a-b-n$. Hence, $n+m=ab-a-b$. Thus, $n+m$ is not an
$\mathbb{N}$-LC of $a$ and $b$ (since Theorem \ref{thm.ent.sylv2full}
\textbf{(b)} shows that $ab-a-b$ is not an $\mathbb{N}$-LC of $a$ and $b$).

We shall now prove the following two claims:

\begin{statement}
\textit{Claim 1:} At \textbf{least} one of the two numbers $n$ and $m$ is an
$\mathbb{N}$-LC of $a$ and $b$.
\end{statement}

\begin{statement}
\textit{Claim 2:} At \textbf{most} one of the two numbers $n$ and $m$ is an
$\mathbb{N}$-LC of $a$ and $b$.
\end{statement}

\begin{proof}
[Proof of Claim 1.]Lemma \ref{lem.ent.sylv2lim} shows that there exist two
integers $u$ and $v$ such that $0\leq u\leq b-1$ and $ua+vb=n$. Consider these
$u$ and $v$. Now,%
\begin{align}
\left(  b-1-u\right)  a+\left(  -v-1\right)  b  &  =ba-a-ua-vb-b\nonumber\\
&  =\underbrace{ba}_{=ab}-\,a-b-\underbrace{\left(  ua+vb\right)  }%
_{=n}\nonumber\\
&  =ab-a-b-n=m \label{pf.thm.ent.sylv2full.d.1}%
\end{align}
(by the definition of $m$). We are in one of the following two cases:

\textit{Case 1:} We have $v\geq0$.

\textit{Case 2:} We have $v<0$.

Let us first consider Case 1. In this case, we have $v\geq0$. Thus,
$v\in\mathbb{N}$. Also, $u\in\mathbb{N}$ (since $0\leq u$). Recall that
$ua+vb=n$, so that $n=\underbrace{u}_{\in\mathbb{N}}a+\underbrace{v}%
_{\in\mathbb{N}}b$. This shows that $n$ is an $\mathbb{N}$-LC of $a$ and $b$.
Thus, at least one of the two numbers $n$ and $m$ is an $\mathbb{N}$-LC of $a$
and $b$. So we have proved Claim 1 in Case 1.

Let us next consider Case 2. In this case, we have $v<0$. Hence, $-v>0$, so
that $-v\geq1$ (since $-v$ is an integer) and therefore $-v-1\geq0$. Thus,
$-v-1\in\mathbb{N}$. Moreover, from $u\leq b-1$, we obtain $b-1-u\geq0$, so
that $b-1-u\in\mathbb{N}$. However, (\ref{pf.thm.ent.sylv2full.d.1}) yields%
\[
m=\underbrace{\left(  b-1-u\right)  }_{\in\mathbb{N}}a+\underbrace{\left(
-v-1\right)  }_{\in\mathbb{N}}b.
\]
This shows that $m$ is an $\mathbb{N}$-LC of $a$ and $b$. Thus, at least one
of the two numbers $n$ and $m$ is an $\mathbb{N}$-LC of $a$ and $b$. So we
have proved Claim 1 in Case 2.

Thus, Claim 1 holds in each of Cases 1 and 2. The proof of Claim 1 is
therefore complete.
\end{proof}

\begin{proof}
[Proof of Claim 2.]Assume the contrary. Thus, both numbers $n$ and $m$ are
$\mathbb{N}$-LCs of $a$ and $b$. Therefore, we can write $n$ as $n=xa+yb$ for
some $x,y\in\mathbb{N}$ (since $n$ is an $\mathbb{N}$-LC of $a$ and $b$).
Furthermore, we can write $m$ as $m=za+wb$ for some $z,w\in\mathbb{N}$ (since
$m$ is an $\mathbb{N}$-LC of $a$ and $b$). Consider these $x,y,z,w$. Now,
adding the equalities $n=xa+yb$ and $m=za+wb$ together, we obtain
\[
n+m=\left(  xa+yb\right)  +\left(  za+wb\right)  =\underbrace{\left(
x+z\right)  }_{\in\mathbb{N}}a+\underbrace{\left(  y+w\right)  }%
_{\in\mathbb{N}}b.
\]
This shows that $n+m$ is an $\mathbb{N}$-LC of $a$ and $b$. This contradicts
the fact that $n+m$ is not an $\mathbb{N}$-LC of $a$ and $b$. This
contradiction shows that our assumption was wrong. Hence, Claim 2 is proved.
\end{proof}

Combining Claim 1 with Claim 2, we see that exactly one of the two numbers $n$
and $m$ is an $\mathbb{N}$-LC of $a$ and $b$. In other words, exactly one of
the two numbers $n$ and $ab-a-b-n$ is an $\mathbb{N}$-LC of $a$ and $b$ (since
$m=ab-a-b-n$). This proves Theorem \ref{thm.ent.sylv2full} \textbf{(d)}.
\medskip

\textbf{(a)} Let $n>ab-a-b$. Then, the integer $ab-a-b-n$ is negative, and
thus cannot be an $\mathbb{N}$-LC of $a$ and $b$ (since any $\mathbb{N}$-LC of
$a$ and $b$ is $\geq0$). However, Theorem \ref{thm.ent.sylv2full} \textbf{(d)}
yields that exactly one of the two numbers $n$ and $ab-a-b-n$ is an
$\mathbb{N}$-LC of $a$ and $b$. Since $ab-a-b-n$ cannot be an $\mathbb{N}$-LC
of $a$ and $b$, we thus conclude that $n$ is an $\mathbb{N}$-LC of $a$ and
$b$. This proves Theorem \ref{thm.ent.sylv2full} \textbf{(a)}. \medskip

\textbf{(c)} Consider the following table of integers:%
\[%
\begin{tabular}
[c]{|c|c|c|c|c|c|}\hline
$0$ & $1$ & $2$ & $\cdots\cdots\cdots$ & $ab-a-b-1$ & $ab-a-b$\\\hline
$ab-a-b$ & $ab-a-b-1$ & $ab-a-b-2$ & $\cdots\cdots\cdots$ & $1$ & $0$\\\hline
\end{tabular}
\]
(whose first row is listing the numbers $0,1,2,\ldots,ab-a-b$ in increasing
order, while the second row is listing the same numbers in decreasing order).
This table has $ab-a-b+1=\left(  a-1\right)  \left(  b-1\right)  $ many columns.

Each column of this table contains the numbers $n$ and $ab-a-b-n$ for some
$n\in\left\{  0,1,\ldots,ab-a-b\right\}  $. Thus, each column of this table
contains exactly one $\mathbb{N}$-LC of $a$ and $b$ (by Theorem
\ref{thm.ent.sylv2full} \textbf{(d)}). Hence, in total, exactly $\left(
a-1\right)  \left(  b-1\right)  $ entries of our table are $\mathbb{N}$-LCs of
$a$ and $b$ (since our table has $\left(  a-1\right)  \left(  b-1\right)  $
many columns). Since our table contains each element of the set $\left\{
0,1,\ldots,ab-a-b\right\}  $ exactly twice, this entails that exactly
$\dfrac{\left(  a-1\right)  \left(  b-1\right)  }{2}$ elements of this set are
$\mathbb{N}$-LCs of $a$ and $b$. In other words, among the elements of the set
$\left\{  0,1,\ldots,ab-a-b\right\}  $, exactly half are $\mathbb{N}$-LCs of
$a$ and $b$. But this is precisely the claim of Theorem
\ref{thm.ent.sylv2full} \textbf{(c)}. Thus, Theorem \ref{thm.ent.sylv2full}
\textbf{(c)} is proved.
\end{proof}
\end{fineprint}

Theorem \ref{thm.ent.sylv2full} is one of the deepest results we will see in
this course, but it is only the beginning of a theory! See
\href{https://en.wikipedia.org/wiki/Coin_problem}{the Wikipedia page for
\textquotedblleft Coin problem\textquotedblright} for more general (and
trickier) questions, such as describing the $\mathbb{N}$-LCs of three integers
$a,b,c$. See also
\href{https://www.math.miami.edu/~armstrong/Talks/RCC_FPSAC_17.pdf}{the slides
of Drew Armstrong's talk at FPSAC 2017} for deep connections to algebraic
combinatorics (and a visual proof different from ours).

\subsection{\label{sec.ent.crypt}Digression: An introduction to cryptography}

\begin{fineprint}
In this short section, we shall make a short foray into \textbf{cryptography}
(also known as \textbf{cryptology}): the study of ciphers, i.e., methods of
encrypting data, mostly for the purpose of maintaining secrecy or proving
authenticity. This is a wide field with a several thousand years' long
history; while it is not fully part of mathematics (as it is governed to a
significant extent by real-life limitations and the \textquotedblleft human
factor\textquotedblright), it relies on mathematical concepts and results.

We will see an ancient (Roman) as well as a modern (20th century) cipher. Both
are underlain by elementary number theory. The second is still in use
(occasionally), whereas the former is only used for recreational purposes
(e.g., hiding spoilers in forum posts). Many more ciphers have been invented
over the ages, and much has been learned about how to break them
(\textquotedblleft cryptanalysis\textquotedblright) and how to keep them
secure. As so often, we will only reach skin-deep into the subject. The
interested reader can learn much more from popular introductions such as
\cite{Beutel94} as well as many number theory texts such as \cite{KraWas15}
and \cite{KraWas18}. (Ciphers can also be based on other parts of mathematics,
but the majority use number theory and abstract algebra.)
\end{fineprint}

\subsubsection{Caesar ciphers (alphabet rotation)}

\begin{fineprint}
We begin with an algorithm that was supposedly used by Julius Caesar to
encrypt military communications. We assume that our messages are textual and
are written in the modern Latin alphabet, all in uppercase
letters.\footnote{Other alphabets (and lowercase letters) can be handled
similarly. Note that the Romans had
\href{https://en.wikipedia.org/wiki/Latin_alphabet}{a slightly different Latin
alphabet than we do}, but we shall use the modern one (with its 26 letters)
for the sake of familiarity.}

The modern Latin alphabet has $26$ letters: A, B, $\ldots$, Z. Let us assign a
number to each of these letters in the most natural way:%
\[%
\begin{tabular}
[c]{ccccccccccccc}\hline\hline
\multicolumn{1}{||c}{A} & \multicolumn{1}{|c}{B} & \multicolumn{1}{|c}{C} &
\multicolumn{1}{|c}{D} & \multicolumn{1}{|c}{E} & \multicolumn{1}{|c}{F} &
\multicolumn{1}{|c}{G} & \multicolumn{1}{|c}{H} & \multicolumn{1}{|c}{I} &
\multicolumn{1}{|c}{J} & \multicolumn{1}{|c}{K} & \multicolumn{1}{|c}{L} &
\multicolumn{1}{|c||}{M}\\\hline
\multicolumn{1}{||c}{$0$} & \multicolumn{1}{|c}{$1$} & \multicolumn{1}{|c}{$2$%
} & \multicolumn{1}{|c}{$3$} & \multicolumn{1}{|c}{$4$} &
\multicolumn{1}{|c}{$5$} & \multicolumn{1}{|c}{$6$} & \multicolumn{1}{|c}{$7$}
& \multicolumn{1}{|c}{$8$} & \multicolumn{1}{|c}{$9$} &
\multicolumn{1}{|c}{$10$} & \multicolumn{1}{|c}{$11$} &
\multicolumn{1}{|c||}{$12$}\\\hline\hline
&  &  &  &  &  &  &  &  &  &  &  & \\\hline\hline
\multicolumn{1}{||c}{N} & \multicolumn{1}{|c}{O} & \multicolumn{1}{|c}{P} &
\multicolumn{1}{|c}{Q} & \multicolumn{1}{|c}{R} & \multicolumn{1}{|c}{S} &
\multicolumn{1}{|c}{T} & \multicolumn{1}{|c}{U} & \multicolumn{1}{|c}{V} &
\multicolumn{1}{|c}{W} & \multicolumn{1}{|c}{X} & \multicolumn{1}{|c}{Y} &
\multicolumn{1}{|c||}{Z}\\\hline
\multicolumn{1}{||c}{$13$} & \multicolumn{1}{|c}{$14$} &
\multicolumn{1}{|c}{$15$} & \multicolumn{1}{|c}{$16$} &
\multicolumn{1}{|c}{$17$} & \multicolumn{1}{|c}{$18$} &
\multicolumn{1}{|c}{$19$} & \multicolumn{1}{|c}{$20$} &
\multicolumn{1}{|c}{$21$} & \multicolumn{1}{|c}{$22$} &
\multicolumn{1}{|c}{$23$} & \multicolumn{1}{|c}{$24$} &
\multicolumn{1}{|c||}{$25$}\\\hline\hline
\end{tabular}
\ \ \ .
\]
Thus, each letter corresponds to a unique number in the set $\left\{
0,1,\ldots,25\right\}  $. For instance, the letter F corresponds to the number
$5$, and the letter X corresponds to the number $23$. This gives us a method
to encode letters as numbers (and, conversely, decode numbers back into
letters); this method will be called \textbf{numeric encoding} of letters.

A word is just a finite list of letters: For example, the word
\textquotedblleft KITTEN\textquotedblright\ is the list $\left(
\text{K},\text{ I},\ \text{T},\text{ T},\text{ E},\text{ N}\right)  $. If we
encode each of these six letters numerically, then we obtain the list $\left(
10,\ 8,\ 19,\ 19,\ 4,\ 13\right)  $ (since the letter K corresponds to $10$,
the letter I to $8$, and so on). This way, we can encode any word as a finite
list of numbers (specifically, of numbers in the set $\left\{  0,1,\ldots
,25\right\}  $). Conversely, any finite list of such numbers can be decoded
into a word (although not necessarily a meaningful word): For instance, the
list $\left(  17,\ 0,\ 19\right)  $ decodes as \textquotedblleft
RAT\textquotedblright, since the number $17$ corresponds to the letter R, the
number $0$ to the letter A, and the number $19$ to the letter T.

We can now formulate Caesar's algorithm, which is nowadays known as the
\textquotedblleft Caesar cipher $\operatorname*{ROT}\nolimits_{3}%
$\textquotedblright\ (we will soon see other variants):

\begin{statement}
\textbf{Caesar cipher }$\operatorname*{ROT}\nolimits_{3}$\textbf{:} To encrypt
a word (written in the modern Latin alphabet, all uppercase), proceed as follows:

\begin{enumerate}
\item Encode the word as a finite list of numbers $\left(  a_{1},a_{2}%
,\ldots,a_{n}\right)  $ (using the numeric encoding).

\item Replace each number $a_{i}$ in this list by $\left(  a_{i}+3\right)
\%26$.

\item Decode the resulting list back into a word.
\end{enumerate}
\end{statement}

\begin{example}
Let us encrypt the word \textquotedblleft CRAZY\textquotedblright\ using the
Caesar cipher $\operatorname*{ROT}\nolimits_{3}$. First, we encode it as a
finite list of numbers:%
\[
\text{CRAZY }\longrightarrow\ \left(  2,\ 17,\ 0,\ 25,\ 24\right)  .
\]
Next, we replace each number $a_{i}$ in this list by $\left(  a_{i}+3\right)
\%26$. Thus,

\begin{itemize}
\item we replace the number $2$ by $\left(  2+3\right)  \%26=5$,

\item we replace the number $17$ by $\left(  17+3\right)  \%26=20$,

\item we replace the number $0$ by $\left(  0+3\right)  \%26=3$,

\item we replace the number $25$ by $\left(  25+3\right)  \%26=28\%26=2$, and

\item we replace the number $24$ by $\left(  24+3\right)  \%26=27\%26=1$.
\end{itemize}

Our list $\left(  2,\ 17,\ 0,\ 25,\ 24\right)  $ thus turns into the new list
$\left(  5,\ 20,\ 3,\ 2,\ 1\right)  $. Decoding the latter list back into a
word, we find \textquotedblleft FUDCB\textquotedblright.
\end{example}

An easy way to visualize the Caesar cipher $\operatorname*{ROT}\nolimits_{3}$
is by placing the $26$ letters of the alphabet in the sectors of a
\textquotedblleft$26$-hour clock\textquotedblright\ (an analog clock with $26$
hours instead of the usual $12$), in the order A, B, C, $\ldots$, Z clockwise.
This \textquotedblleft alphabet clock\textquotedblright\ looks as follows:%
\[%
\begin{tikzpicture}[scale=2.5]
\filldraw[fill=yellow!50!white] (0,0) circle (1.8);
\foreach\x in {0,...,25}
{\node(a\x) at ({90 + 360*\x/26}:2) {\symbol{\numexpr90-\x}} ;
\draw({90 + 360*\x/26}:1.8) -- ({90 + 360*\x/26}:1.5); }
\end{tikzpicture}%
\ \ .
\]
Then, $\operatorname*{ROT}\nolimits_{3}$ simply shifts each letter forward by
$3$ \textquotedblleft hours\textquotedblright\ (so A becomes D, whereas B
becomes E, and so on).

Thus, it is clear how we can decrypt a word encrypted using
$\operatorname*{ROT}\nolimits_{3}$: We just need to shift each letter backward
by $3$ \textquotedblleft hours\textquotedblright, i.e., replace each $a_{i}$
by $\left(  a_{i}-3\right)  \%26$. We can denote this operation by
$\operatorname*{ROT}\nolimits_{-3}$.

More generally, we define the operation $\operatorname*{ROT}\nolimits_{k}$ for
any integer $k$ as follows:

\begin{statement}
\textbf{Caesar cipher }$\operatorname*{ROT}\nolimits_{k}$\textbf{ (for a given
integer }$k$\textbf{):} To encrypt a word (written in the modern Latin
alphabet, all uppercase), proceed as follows:

\begin{enumerate}
\item Encode the word as a finite list of numbers $\left(  a_{1},a_{2}%
,\ldots,a_{n}\right)  $ (using the numeric encoding).

\item Replace each number $a_{i}$ in this list by $\left(  a_{i}+k\right)
\%26$.

\item Decode the resulting list back into a word.
\end{enumerate}
\end{statement}

In terms of our \textquotedblleft letter clock\textquotedblright,
$\operatorname*{ROT}\nolimits_{k}$ shifts each letter forward by $k$
\textquotedblleft hours\textquotedblright. It is easy to see that a word
encrypted using $\operatorname*{ROT}\nolimits_{k}$ can be decrypted back using
$\operatorname*{ROT}\nolimits_{-k}$, since $\operatorname*{ROT}\nolimits_{-k}$
shifts each letter backward by $k$ \textquotedblleft hours\textquotedblright.
We can also prove this rigorously using our definition of $\operatorname*{ROT}%
\nolimits_{k}$, using the following simple lemma:

\begin{lemma}
\label{lem.crypto.ROT+-}Let $k$ be an integer. Let $a,b\in\left\{
0,1,\ldots,25\right\}  $ be two numbers satisfying $b=\left(  a+k\right)
\%26$. Then, $a=\left(  b-k\right)  \%26$.
\end{lemma}

\begin{proof}
We have $b=\left(  a+k\right)  \%26\equiv a+k\operatorname{mod}26$ (by
Proposition \ref{prop.ent.quo-rem.basics} \textbf{(a)}, applied to $n=a+k$ and
$d=26$). Subtracting the trivial congruence $k\equiv k\operatorname{mod}26$
from this congruence, we obtain $b-k\equiv a\operatorname{mod}26$, so that
$a\equiv b-k\operatorname{mod}26$. Hence, Proposition
\ref{prop.ent.quo-rem.modd} (applied to $b-k$ and $26$ instead of $b$ and $d$)
yields $a\%26=\left(  b-k\right)  \%26$. However, from $a\in\left\{
0,1,\ldots,25\right\}  $, we see that the remainder $a\%26$ is $a$ itself.
Thus, $a=a\%26=\left(  b-k\right)  \%26$. This proves Lemma
\ref{lem.crypto.ROT+-}.
\end{proof}

Of course, we used nothing special about the number $26$ here; we could just
as well replace it by any fixed positive integer $m$ in Lemma
\ref{lem.crypto.ROT+-}.

Lemma \ref{lem.crypto.ROT+-} shows that the Caesar cipher $\operatorname*{ROT}%
\nolimits_{-k}$ undoes the Caesar cipher $\operatorname*{ROT}\nolimits_{k}$:
Indeed, when we apply $\operatorname*{ROT}\nolimits_{k}$ to a word, each entry
$a$ in the corresponding list of numbers gets replaced by $b:=\left(
a+k\right)  \%26$; then, a subsequent application of $\operatorname*{ROT}%
\nolimits_{-k}$ replaces this new number $b$ by $\left(  b+\left(  -k\right)
\right)  \%26=\left(  b-k\right)  \%26=a$ (by Lemma \ref{lem.crypto.ROT+-}),
which is the original entry before $\operatorname*{ROT}\nolimits_{k}$ was applied.

\begin{exercise}
\label{exe.crypto.ROT.1}\textbf{(a)} Encrypt the word \textquotedblleft
REED\textquotedblright\ using $\operatorname*{ROT}\nolimits_{4}$. \medskip

\textbf{(b)} Encrypt the word \textquotedblleft BOON\textquotedblright\ using
$\operatorname*{ROT}\nolimits_{16}$. \medskip

\textbf{(c)} Some word was encrypted using $\operatorname*{ROT}\nolimits_{10}$
and became \textquotedblleft GSDROB\textquotedblright. Reconstruct the
original word. \medskip

\textbf{(d)} Some (meaningful) word was encrypted using $\operatorname*{ROT}%
\nolimits_{k}$ for some unknown integer $k$ and became \textquotedblleft
WBBSFACGH\textquotedblright. Reconstruct the original word. (This illustrates
why Caesar ciphers are not very secure, to put it mildly. There are quick ways
to solve this without trying \textquotedblleft all\textquotedblright\ possibilities.)
\end{exercise}

Let us make some further observations about Caesar ciphers:

\begin{itemize}
\item The encryption method $\operatorname*{ROT}\nolimits_{0}$ does nothing:
Each word is encrypted as itself (since shifting by $0$ \textquotedblleft
hours\textquotedblright\ on the letter clock changes nothing, or since
$\left(  a+0\right)  \%26=a\%26=a$ for each $a\in\left\{  0,1,\ldots
,25\right\}  $).

\item The encryption method $\operatorname*{ROT}\nolimits_{26}$ also does
nothing: Each word is encrypted as itself (since shifting by $26$
\textquotedblleft hours\textquotedblright\ on the letter clock amounts to a
full revolution, or since $\left(  a+26\right)  \%26=a$ for each $a\in\left\{
0,1,\ldots,25\right\}  $).

\item The encryption method $\operatorname*{ROT}\nolimits_{27}$ does the same
as $\operatorname*{ROT}\nolimits_{1}$ (since $\left(  a+27\right)
\%26=\left(  a+1\right)  \%26$ for each $a\in\left\{  0,1,\ldots,25\right\}  $).

\item More generally, if two integers $u$ and $v$ satisfy $u\equiv
v\operatorname{mod}26$, then $\operatorname*{ROT}\nolimits_{u}%
=\operatorname*{ROT}\nolimits_{v}$. Thus, there are only $26$ distinct Caesar
ciphers, namely%
\[
\operatorname*{ROT}\nolimits_{0},\ \operatorname*{ROT}\nolimits_{1}%
,\ \ldots,\ \operatorname*{ROT}\nolimits_{25}.
\]
Any other $\operatorname*{ROT}\nolimits_{k}$ is just a copy of one of these.
Of these $26$ ciphers, only $25$ are useful, since $\operatorname*{ROT}%
\nolimits_{0}$ does nothing.

\item The cipher $\operatorname*{ROT}\nolimits_{13}$ inverts itself: Any word
encrypted using $\operatorname*{ROT}\nolimits_{13}$ can be decrypted by
applying $\operatorname*{ROT}\nolimits_{13}$ again. Indeed,
$\operatorname*{ROT}\nolimits_{13}$ is undone by $\operatorname*{ROT}%
\nolimits_{-13}$, but $\operatorname*{ROT}\nolimits_{-13}=\operatorname*{ROT}%
\nolimits_{13}$ because $-13\equiv13\operatorname{mod}26$.

\item Encrypting a word using $\operatorname*{ROT}\nolimits_{u}$ (for some
integer $u$) and then encrypting the result using $\operatorname*{ROT}%
\nolimits_{v}$ (for some integer $v$) is the same as encrypting the original
word using $\operatorname*{ROT}\nolimits_{u+v}$.
\end{itemize}

\begin{exercise}
Prove the latter statement rigorously using the description in terms of
remainders. That is, prove the following fact:

If $u$ and $v$ are two integers, and if $a,b,c\in\left\{  0,1,\ldots
,25\right\}  $ are three numbers satisfying $b=\left(  a+u\right)  \%26$ and
$c=\left(  b+v\right)  \%26$, then $c=\left(  a+\left(  u+v\right)  \right)
\%26$.
\end{exercise}

We have so far been encrypting single words. To encrypt an entire text, one
must decide what to do about whitespaces. There are different legitimate
choices: e.g., one can leave them unchanged; one can remove them (at the risk
of making the text hard to read even after decryption); or one can treat them
as a \textquotedblleft$27$th letter\textquotedblright\ of the alphabet (thus
adapting the definition of Caesar ciphers to use $\left(  a+k\right)  \%27$
instead of $\left(  a+k\right)  \%26$). We shall not delve any deeper into
these questions here.
\end{fineprint}

\subsubsection{Keys and ciphers}

\begin{fineprint}
Ciphers such as $\operatorname*{ROT}\nolimits_{k}$ are one-trick ponies: Once
your enemy knows the method, he will be able to decrypt anything you encrypt.

This is true to an extent even if the enemy does \textbf{not} know the $k$,
but only knows that you have used some Caesar cipher. Indeed, there are only
$26$ Caesar ciphers
\[
\operatorname*{ROT}\nolimits_{0},\ \operatorname*{ROT}\nolimits_{1}%
,\ \ldots,\ \operatorname*{ROT}\nolimits_{25}.
\]
Thus, if your enemy finds a text you encrypted using some $\operatorname*{ROT}%
\nolimits_{k}$, he can just try to decrypt it using%
\[
\operatorname*{ROT}\nolimits_{-0},\ \operatorname*{ROT}\nolimits_{-1}%
,\ \ldots,\ \operatorname*{ROT}\nolimits_{-25},
\]
and see which of the results gives a meaningful word/text rather than
gibberish (see Exercise \ref{exe.crypto.ROT.1} \textbf{(d)}).\footnote{The
answer might be non-unique when the word is short (see Exercise
\ref{exe.crypto.ROT.1} \textbf{(b)}), but will practically always be unique
when the word/text is long enough.}

In modern language, this is saying that Caesar ciphers have too small a
\textbf{key size} to be secure. The \textbf{key} here is the number $k$. While
technically there are infinitely many options for $k$, there are only $26$
distinct ciphers obtained, so the \textquotedblleft true\textquotedblright%
\ key is just an element of $\left\{  0,1,\ldots,25\right\}  $. No wonder the
cipher is easily broken.

Another problem with Caesar ciphers is that they are \textquotedblleft too
regular\textquotedblright: e.g., equal letters in the original word remain
equal after encryption. This, too, causes weaknesses that render the cipher
easy to break.

So how can we create a cipher that is harder to break? We need a bigger key
size, and we need \textquotedblleft more chaos\textquotedblright\ (e.g., don't
apply the same rule to each letter). Here are some ciphers that are slightly
better in some of these regards:

\begin{itemize}
\item \textbf{Monoalphabetic substitution:} Here we still do the same thing to
each letter, but this thing is no longer just a shift by $k$ \textquotedblleft
hours\textquotedblright. Instead, we fix \textbf{any} permutation of the
alphabet (i.e., a rule that sends each letter to a different letter) and we
apply this permutation separately to each letter. For instance, we can use the
following permutation:%
\[%
\begin{tabular}
[c]{ccccccccccccc}\hline\hline
\multicolumn{1}{||c}{A} & \multicolumn{1}{|c}{B} & \multicolumn{1}{|c}{C} &
\multicolumn{1}{|c}{D} & \multicolumn{1}{|c}{E} & \multicolumn{1}{|c}{F} &
\multicolumn{1}{|c}{G} & \multicolumn{1}{|c}{H} & \multicolumn{1}{|c}{I} &
\multicolumn{1}{|c}{J} & \multicolumn{1}{|c}{K} & \multicolumn{1}{|c}{L} &
\multicolumn{1}{|c||}{M}\\\hline
\multicolumn{1}{||c}{C} & \multicolumn{1}{|c}{Z} & \multicolumn{1}{|c}{X} &
\multicolumn{1}{|c}{B} & \multicolumn{1}{|c}{N} & \multicolumn{1}{|c}{M} &
\multicolumn{1}{|c}{P} & \multicolumn{1}{|c}{A} & \multicolumn{1}{|c}{D} &
\multicolumn{1}{|c}{T} & \multicolumn{1}{|c}{S} & \multicolumn{1}{|c}{R} &
\multicolumn{1}{|c||}{Q}\\\hline\hline
&  &  &  &  &  &  &  &  &  &  &  & \\\hline\hline
\multicolumn{1}{||c}{N} & \multicolumn{1}{|c}{O} & \multicolumn{1}{|c}{P} &
\multicolumn{1}{|c}{Q} & \multicolumn{1}{|c}{R} & \multicolumn{1}{|c}{S} &
\multicolumn{1}{|c}{T} & \multicolumn{1}{|c}{U} & \multicolumn{1}{|c}{V} &
\multicolumn{1}{|c}{W} & \multicolumn{1}{|c}{X} & \multicolumn{1}{|c}{Y} &
\multicolumn{1}{|c||}{Z}\\\hline
\multicolumn{1}{||c}{K} & \multicolumn{1}{|c}{O} & \multicolumn{1}{|c}{E} &
\multicolumn{1}{|c}{W} & \multicolumn{1}{|c}{Y} & \multicolumn{1}{|c}{U} &
\multicolumn{1}{|c}{I} & \multicolumn{1}{|c}{J} & \multicolumn{1}{|c}{F} &
\multicolumn{1}{|c}{L} & \multicolumn{1}{|c}{G} & \multicolumn{1}{|c}{H} &
\multicolumn{1}{|c||}{V}\\\hline\hline
\end{tabular}
\ \ \ \ \ \ \ .
\]
Then, the word \textquotedblleft KITTEN\textquotedblright\ is encrypted as
\textquotedblleft SDIINK\textquotedblright.

The key size of this encryption method is huge: The number of possible keys is
the number of all permutations of the alphabet, which is (as we will soon
see\footnote{See Corollary \ref{cor.count.select.order-n!} and the discussion
that follows it.}) $26!=\allowbreak403\,291\,461\,\allowbreak
126\,605\,635\,\allowbreak584\,000\,000$. You cannot just try each of these
keys to see which one works. But you can still exploit certain patterns in the
English language (or whatever language the text is written in), such as
frequencies of letters, frequencies of two-letter combinations, and so on. If
you have a ciphertext (i.e., encrypted text) of sufficient length (e.g., a
page, but often a paragraph will be enough), you can break a monoalphabetic
substitution cipher using just statistics and a bit of patience (see, e.g.,
\cite[\S 1.6]{Beutel94} for details). Essentially, this is because the cipher
is \textquotedblleft not chaotic enough\textquotedblright.

\item The \textbf{one-time pad:} Now the key is an infinite (or finite but
sufficiently long) sequence%
\[
\left(  k_{1},k_{2},k_{3},\ldots\right)
\]
of elements of $\left\{  0,1,\ldots,25\right\}  $ (or just of integers). To
encrypt a word, we first encode it as a tuple of integers $\left(  a_{1}%
,a_{2},\ldots,a_{n}\right)  $, and then replace each number $a_{i}$ by
$\left(  a_{i}+k_{i}\right)  \%26$; then, we decode the resulting tuple back
into a word.

This is essentially a generalized Caesar cipher, in which we let each letter
get a different key depending on its position.

This cipher is completely unbreakable (provided the key is truly random), but
it is also very inconvenient: You need an infinitely long key, or at least a
key that is at least as long as the text you want to encrypt.

In many cases, this becomes impractical, so people have tried to
\textquotedblleft cheat\textquotedblright, e.g., by using a periodic sequence
$\left(  k_{1},k_{2},k_{3},\ldots\right)  $ (this is known as the
\textbf{Vigen\`{e}re substitution}) or using a passage from a book as the key
(this is the \textbf{running key cipher}); but such cheats make the cipher
breakable when the ciphertext is sufficiently long (see \cite[\S 2.3]%
{Beutel94}). Likewise, if you reuse the same sequence $\left(  k_{1}%
,k_{2},k_{3},\ldots\right)  $ as a key for several messages, certain
frequency-based patterns will appear in your ciphertexts that will eventually
give away the key.
\end{itemize}

Many different algorithms have been invented through the ages, usually
striking some balance between practicality (ease of use, simplicity, shortness
of the key) and security (unbreakability). See \cite{Singh01} for more
classical algorithms and their history.
\end{fineprint}

\subsubsection{\label{subsec.ent.crypt.RSA}The RSA cipher}

\begin{fineprint}
All ciphers invented until the early 20th century are
\textbf{\href{https://en.wikipedia.org/wiki/Classical_cipher}{\textbf{classical
ciphers}}}: ciphers that can be computed (specifically, encrypted and
decrypted with) by hand, without the use of computers. Practically all these
ciphers can be broken with the help of computers (provided the ciphertext is
long enough), and thus are obsolete in the 21st century. (Actually, breaking
ciphers was one of the earliest uses of computers: The
\href{https://en.wikipedia.org/wiki/Cryptanalysis_of_the_Enigma}{quest} to
break \href{https://en.wikipedia.org/wiki/Enigma_machine}{the German Enigma
cipher} during World War II was a major motivation for the development of
computers in the mid-20th century.)

\textbf{Modern ciphers} are ciphers that require a computer to encrypt and
decrypt (at realistically useful speeds). Using the computational power of
modern electronics, they can afford to rely on much longer calculations with
much larger numbers than $0,1,\ldots,25$. Quite a few modern ciphers are
nowadays considered unbreakable, in the sense that no realistic methods for
breaking them are known (unless the ciphers are used incorrectly), and there
are good reasons to assume that such methods do not exist on any currently
existing hardware.

In this subsection, we will discuss one such modern cipher:
\href{https://en.wikipedia.org/wiki/RSA_(cryptosystem)}{the \textbf{RSA
cipher}}, developed by Rivest, Shamir and Adleman in 1977. Like a Caesar
cipher, it is based upon division with remainder, but in a much less
\textquotedblleft predictable\textquotedblright\ way.

Unlike most classical ciphers, the RSA cipher is surprisingly robust, in the
sense that it can be used in a much less \textquotedblleft
fair-weather\textquotedblright\ situation. For most classical ciphers, the
sender and the receiver must have \textbf{privately} agreed on the key (e.g.,
the $k$ in a Caesar cipher) \textbf{in advance} (ideally at a private
vis-a-vis long before the need for secrecy arises). If the enemy has managed
to eavesdrop on this agreement, he will know the key, and the cipher will be
useless. In contrast, in the RSA cipher, the parties can agree on their keys
\textbf{whenever} they need them, and they can do so even over a completely
\textbf{public} channel (e.g., screaming at each other from the rooftops, or
posting on reddit)!

This rather counterintuitive feature is achieved by using two types of keys:
\textbf{public keys} (which are sent out in plain text, so that every curious
outsider knows them) and \textbf{private keys} (which the sender and the
receiver compute individually, and don't share with anyone -- not even with
each other!). \medskip

Let me describe (while omitting some technicalities) how the RSA cipher works.
Assume that Albert and Julia want to communicate securely over a public
channel (e.g., an internet forum with no private-message functionality). They
cannot hide the fact that they are talking to each other (at least not using
the RSA cipher), but they want to hide the contents of their communications by
encrypting them in a way that no eavesdroppers can decipher. Albert and Julia
have not exchanged keys in advance. What do they do?

First, they need to \textbf{set up the cipher:}

\begin{itemize}
\item Julia tells Albert (over the public channel) that she wants to
communicate, and thus he should start creating keys.

\item Albert generates two distinct large and sufficiently random primes $p$
and $q$.

[What exactly does this mean, and how does he do this? With modern hardware,
\textquotedblleft large\textquotedblright\ means approximately $300$ digits or
more. \textquotedblleft Sufficiently random\textquotedblright\ means
\textquotedblleft%
\href{https://en.wikipedia.org/wiki/Pseudorandomness}{pseudorandom}%
\textquotedblright, e.g., (roughly speaking) practically unpredictable and
devoid of discernible patterns. Large pseudorandom primes can be generated
fairly fast by
\href{https://en.wikipedia.org/wiki/Generation_of_primes}{various algorithms},
with a bit of input from the outside world (e.g., Geiger counters) to generate randomness.]

\item Albert computes the positive integers%
\[
m:=pq\ \ \ \ \ \ \ \ \ \ \text{and}\ \ \ \ \ \ \ \ \ \ \ell:=\left(
p-1\right)  \left(  q-1\right)  .
\]
He makes the number $m$ public (i.e., sends it to Julia over the public
channel), but keeps the number $\ell$ private (even Julia does not need to
know it). Eavesdroppers will thus learn $m$, but will struggle to find $p$ and
$q$, since \href{https://en.wikipedia.org/wiki/Integer_factorization}{no fast
algorithm for factoring numbers into primes is known}. (If anyone finds such
an algorithm, the RSA cipher will be broken.)

\item Albert randomly picks an integer $e\in\left\{  2,3,\ldots,\ell
-1\right\}  $ that is coprime to $\ell$.

[The easiest way to do so is to pick a bunch of numbers in the set $\left\{
2,3,\ldots,\ell-1\right\}  $ at random, and grab the first of them that is
coprime to $\ell$. Coprimality can be checked quickly using the Euclidean
algorithm. If no chosen number is coprime to $\ell$, then roll the dice again.]

\item Albert computes a positive integer $d$ such that $ed\equiv
1\operatorname{mod}\ell$.

[How? Bezout's theorem (Theorem \ref{thm.ent.gcd.bezout}) shows that
$\gcd\left(  e,\ell\right)  =xe+y\ell$ for some integers $x$ and $y$. These
$x$ and $y$ can be computed quickly using the Extended Euclidean Algorithm
(see Subsection \ref{subsec.ent.gcd.bezout}). Having found these integers $x$
and $y$, we have $\gcd\left(  e,\ell\right)  =xe+y\ell\equiv
xe=ex\operatorname{mod}\ell$ and therefore $ex\equiv\gcd\left(  e,\ell\right)
=1\operatorname{mod}\ell$ (since $e$ is coprime to $\ell$). So we just take
$d=x$, or -- if this is negative -- we take $d=x\%\ell$.]

\item Albert publishes the pair $\left(  e,m\right)  $ (so that Julia knows
it, and so does anyone else who cares to listen). This pair is his
\textbf{public key}, whereas the (secret) pair $\left(  d,m\right)  $ is his
\textbf{private key}.
\end{itemize}

\bigskip

\textbf{Encrypting a message:}

Now, assume that Julia wants to send a message to Albert. She encodes this
message as an element $a$ of the set $\left\{  0,1,\ldots,m-1\right\}  $. (If
it does not fit into this set, she just breaks it up into size-$m$ chunks and
encrypts each chunk separately. Note that the encoding has to be agreed on in
advance, but this can be a public method.)

She computes the remainder $a^{e}\%m$ and sends this remainder to Albert.

[Practical issue: To compute $a^{e}\%m$ fast, she should \textbf{not} try to
compute the huge number $a^{e}$, since there is no space in the universe to
store such a huge number. Instead, she can \textquotedblleft work modulo
$m$\textquotedblright, and use
\href{https://en.wikipedia.org/wiki/Exponentiation_by_squaring}{binary
exponentiation}. For example, to compute $a^{190}$, she should not use the
definition $a^{190}=\underbrace{aa\cdots a}_{190\text{ times}}$ but the much
faster formula $a^{190}=\left(  \left(  \left(  \left(  \left(  \left(
a^{2}\right)  ^{2}a\right)  ^{2}a\right)  ^{2}a\right)  ^{2}a\right)
^{2}a\right)  ^{2}$, and moreover, since she only needs the remainder
$a^{190}\%m$, she can reduce each intermediate result modulo $m$ (that is,
replace it by its remainder upon division by $m$), so that no overly large
numbers should appear in the process.]

\bigskip

\textbf{Decrypting a message:}

Albert receives the remainder $b=a^{e}\%m$. To recover Julia's original
message $a$, he just needs to take the $d$-th power and take its remainder
upon division by $m$. In other words,%
\[
a=b^{d}\%m.
\]

[Just like Julia, Albert should use binary exponentiation and work modulo $m$
to compute this efficiently.]

\bigskip

So the encryption algorithm is just \textquotedblleft take the $e$-th power
and then take its remainder when divided by $m$\textquotedblright, whereas the
decryption algorithm is just \textquotedblleft take the $d$-th power and then
take its remainder when divided by $m$\textquotedblright\ (although the
implementation is a bit more complex, in order to be efficient).

Why does this work? Obviously, we need to prove the following proposition:

\begin{proposition}
[correctness of RSA]\label{prop.crypt.RSA-corr}Let $p$ and $q$ be two distinct
primes. Let $m=pq$ and $\ell=\left(  p-1\right)  \left(  q-1\right)  $. Let
$e$ and $d$ be two positive integers such that $ed\equiv1\operatorname{mod}%
\ell$.

Let $a$ and $b$ be two numbers in $\left\{  0,1,\ldots,m-1\right\}  $ such
that $b=a^{e}\%m$. Then, $a=b^{d}\%m$.
\end{proposition}

This is not at all obvious! The RSA cipher might resemble a Caesar cipher in
that it uses remainders, but it is different in that it takes powers instead
of adding/subtracting a fixed $k$.

To prove Proposition \ref{prop.crypt.RSA-corr}, we will need a lemma, which
resembles Fermat's Little Theorem (Theorem \ref{thm.ent.prime.flt}):

\begin{lemma}
\label{lem.crypt.RSA-lem}Let $p$ and $q$ be two distinct primes. Let $N$ be a
positive integer such that $N\equiv1\operatorname{mod}\left(  p-1\right)
\left(  q-1\right)  $. Let $a$ be any integer. Then,%
\[
a^{N}\equiv a\operatorname{mod}pq.
\]

\end{lemma}

\begin{example}
Let $p=3$ and $q=5$ and $N=9$. Then, $N\equiv1\operatorname{mod}\left(
p-1\right)  \left(  q-1\right)  $, so that Lemma \ref{lem.crypt.RSA-lem}
yields that
\[
a^{9}\equiv a\operatorname{mod}15\ \ \ \ \ \ \ \ \ \ \text{for any integer
}a.
\]

\end{example}

\begin{proof}
[Proof of Lemma \ref{lem.crypt.RSA-lem}.]Fermat's little theorem (Theorem
\ref{thm.ent.prime.flt}) says that $a^{p}\equiv a\operatorname{mod}p$ and
$a^{q}\equiv a\operatorname{mod}q$. Our claim looks similar, but not quite the
same. Nevertheless, we are on the right trail.

We must prove that $a^{N}\equiv a\operatorname{mod}pq$. In other words, we
must prove that $pq\mid a^{N}-a$. But $p$ and $q$ are two distinct primes, and
thus are coprime (why?\footnote{\textit{Proof.} The only positive divisors of
$q$ are $1$ and $q$ (since $q$ is prime). Since $p$ is neither $1$ nor $q$, it
thus follows that $p$ is not a positive divisor of $q$. In other words, $q$ is
not a multiple of $p$. Hence, the friend-or-foe lemma (Lemma
\ref{lem.ent.prime.fof}) shows that $q$ is coprime to $p$.}). Hence, $pq\mid
a^{N}-a$ would follow from the coprime divisors theorem (Theorem
\ref{thm.ent.coprime.ab|c}), if we can show that $p\mid a^{N}-a$ and $q\mid
a^{N}-a$.

It thus remains to prove that $p\mid a^{N}-a$ and $q\mid a^{N}-a$. We will
only show $p\mid a^{N}-a$, since $q\mid a^{N}-a$ is analogous.

So we must show that $p\mid a^{N}-a$. In other words, we must show that
$a^{N}\equiv a\operatorname{mod}p$.

However, $p-1\mid\left(  p-1\right)  \left(  q-1\right)  \mid N-1$ (since
$N\equiv1\operatorname{mod}\left(  p-1\right)  \left(  q-1\right)  $). In
other words, $N-1=\left(  p-1\right)  c$ for some integer $c$. Consider this
$c$. It is easy to see that $c\geq0$ (why?\footnote{\textit{Proof.} Assume the
contrary. Thus, $c<0$. But $p>1$ (since $p$ is prime), so that $p-1>0$. Hence,
$\left(  p-1\right)  c<0$ (since $c<0$). Thus, $N-1=\left(  p-1\right)  c<0$,
so that $N<1$. But $N\geq1$ (since $N$ is a positive integer). This is in
obvious contradiction to $N<1$. Hence, our assumption was false, qed.}), so
that $c\in\mathbb{N}$. From $N-1=\left(  p-1\right)  c$, we obtain
\begin{equation}
N=1+\left(  p-1\right)  c. \label{pf.lem.crypt.RSA-lem.N}%
\end{equation}

However, recall that $a^{p}\equiv a\operatorname{mod}p$. Using this fact, we
can easily see that%
\begin{equation}
a^{1+\left(  p-1\right)  k}\equiv a\operatorname{mod}p
\label{pf.lem.crypt.RSA-lem.k}%
\end{equation}
for each $k\in\mathbb{N}$.

[\textit{Proof of (\ref{pf.lem.crypt.RSA-lem.k}):} We can prove this by
induction on $k$:

\textit{Base case} ($k=0$)\textit{:} We have $a^{1+\left(  p-1\right)
0}=a^{1}=a\equiv a\operatorname{mod}p$. Thus, (\ref{pf.lem.crypt.RSA-lem.k})
holds for $k=0$.

\textit{Induction step:} Let $k\in\mathbb{N}$. Assume (as the induction
hypothesis) that (\ref{pf.lem.crypt.RSA-lem.k}) holds for $k$; that is, we
have $a^{1+\left(  p-1\right)  k}\equiv a\operatorname{mod}p$. We must then
prove that (\ref{pf.lem.crypt.RSA-lem.k}) holds for $k+1$ instead of $k$; in
other words, we must prove that $a^{1+\left(  p-1\right)  \left(  k+1\right)
}\equiv a\operatorname{mod}p$.

But $\left(  p-1\right)  \left(  k+1\right)  =\left(  p-1\right)  k+\left(
p-1\right)  $, and therefore
\begin{align*}
a^{1+\left(  p-1\right)  \left(  k+1\right)  }  &  =a^{1+\left(  p-1\right)
k+\left(  p-1\right)  }=a^{1+\left(  p-1\right)  k}a^{p-1}\\
&  \equiv aa^{p-1}\ \ \ \ \ \ \ \ \ \ \left(
\begin{array}
[c]{c}%
\text{here we multiplied the congruence }a^{1+\left(  p-1\right)  k}\equiv
a\operatorname{mod}p\\
\text{by the trivial congruence }a^{p-1}\equiv a^{p-1}\operatorname{mod}p
\end{array}
\right) \\
&  =a^{p}\equiv a\operatorname{mod}p.
\end{align*}
This completes the induction step. Thus, (\ref{pf.lem.crypt.RSA-lem.k}) is proved.]

Now, (\ref{pf.lem.crypt.RSA-lem.k}) (applied to $k=c$) yields $a^{1+\left(
p-1\right)  c}\equiv a\operatorname{mod}p$ (since $c\in\mathbb{N}$). In view
of (\ref{pf.lem.crypt.RSA-lem.N}), we can rewrite this as $a^{N}\equiv
a\operatorname{mod}p$. In other words, $p\mid a^{N}-a$. Similarly, we can show
that $q\mid a^{N}-a$. As explained above, this completes the proof of Lemma
\ref{lem.crypt.RSA-lem}.
\end{proof}

\begin{proof}
[Proof of Proposition \ref{prop.crypt.RSA-corr}.]We have $b=a^{e}\%m\equiv
a^{e}\operatorname{mod}m$ (by Proposition \ref{prop.ent.quo-rem.basics}
\textbf{(a)}, applied to $a^{e}$ and $m$ instead of $n$ and $d$). Taking this
congruence to the $d$-th power (using Exercise \ref{exe.mod.k-th-power}), we
obtain%
\[
b^{d}\equiv\left(  a^{e}\right)  ^{d}=a^{ed}\operatorname{mod}m.
\]
But $ed\equiv1\operatorname{mod}\ell$, that is, $ed\equiv1\operatorname{mod}%
\left(  p-1\right)  \left(  q-1\right)  $ (since $\ell=\left(  p-1\right)
\left(  q-1\right)  $). Hence, Lemma \ref{lem.crypt.RSA-lem} (applied to
$N=ed$) yields%
\[
a^{ed}\equiv a\operatorname{mod}pq.
\]
In other words, $a^{ed}\equiv a\operatorname{mod}m$ (since $pq=m$). Combining
what we have shown, we obtain%
\[
b^{d}\equiv a^{ed}\equiv a\operatorname{mod}m.
\]
Therefore, Proposition \ref{prop.ent.quo-rem.modd} (applied to $b^{d}$, $a$
and $m$ instead of $a$, $b$ and $d$) yields $b^{d}\%m=a\%m=a$ (since
$a\in\left\{  0,1,\ldots,m-1\right\}  $). This proves Proposition
\ref{prop.crypt.RSA-corr}.
\end{proof}

The RSA cipher, as demonstrated above, lets Julia send secret messages to
Albert. If Albert wants to respond secretly, the two can switch roles (i.e.,
now Julia must set up her two primes $p^{\prime}$ and $q^{\prime}$ and her
$m^{\prime}$, $\ell^{\prime}$, $e^{\prime}$ and $d^{\prime}$, publish her
public key $\left(  e^{\prime},m^{\prime}\right)  $, and let Albert encrypt
his message using that public key). \medskip

The RSA cipher is not hard to implement in your favorite programming language,
provided that it supports sufficiently big integers. But there are some
practical considerations:

\begin{itemize}
\item You need sufficiently random primes. (Generally, any cipher requires
something sufficiently random that the eavesdroppers cannot guess.)

\item Certain primes make for bad choices of $p$ and $q$, since they allow
certain tricks for computing $d$. You want to avoid such primes.

\item You want to avoid certain practical \textquotedblleft side
channels\textquotedblright\ (as with any ciphers).

\item You don't want your message $a$ to be much smaller than $m$. If it is,
pad it with random bits.
\end{itemize}

These and many other caveats are discussed on
\href{https://en.wikipedia.org/wiki/RSA_(cryptosystem)#Padding}{the Wikipedia
page for the RSA cipher}, as well as in more serious textbooks on modern
cryptography (e.g., \cite{BaEdHa18}, \cite{Buchma04} or \cite{HoPiSi14}).

The RSA cipher can be used not just for encrypting secret messages, but also
for authentication (i.e., proving that a message is really coming from you).
See \cite[Chapter 4]{HoPiSi14} or \cite[Chapter 12]{Buchma04} for this application.

There are many other modern ciphers, including elliptic curve cryptography
(see \cite[\S 13.2]{Buchma04} or \cite[Chapter 6]{HoPiSi14}), which is also
built on ideas from number theory.

\begin{exercise}
Prove the following three-primes version of Lemma \ref{lem.crypt.RSA-lem}:

Let $p$, $q$ and $r$ be three distinct primes. Let $N$ be a positive integer
such that $N\equiv1\operatorname{mod}\left(  p-1\right)  \left(  q-1\right)
\left(  r-1\right)  $. Let $a$ be any integer. Then,%
\[
a^{N}\equiv a\operatorname{mod}pqr.
\]

\end{exercise}
\end{fineprint}

\newpage

\section{\label{chp.enum}An informal introduction to enumeration}

Enumeration is a fancy word for counting -- i.e., answering questions of the
form \textquotedblleft how many things of a certain type are
there?\textquotedblright. Here are some examples of counting problems:

\begin{itemize}
\item How many ways are there to choose $3$ odd integers between $0$ and $20$,
if the order matters (i.e., we count the choice $1,3,5$ as different from the
choice $3,1,5$)? (The answer is $1000$.)

\item How many ways are there to choose $3$ odd integers between $0$ and $20$,
if the order does not matter? (The answer is $220$.)

\item How many ways are there to choose $3$ distinct odd integers between $0$
and $20$, if the order matters? (The answer is $720$.)

\item How many ways are there to choose $3$ distinct odd integers between $0$
and $20$, if the order does not matter? (The answer is $120$.)

\item How many prime factorizations does $200$ have (where we count different
orderings as distinct)? (The answer is $10$. This is a mix between a number
theory problem and a counting problem.)

\item How many ways are there to tile a $2\times15$-rectangle with dominos
(i.e., rectangles of size $1\times2$ or $2\times1$) ? (The answer is $987$.
For instance, the tiling
\[
\begin{tikzpicture}[scale=0.8]
\draw (0, 0) rectangle (1, 2);
\draw (1, 0) rectangle (3, 1);
\draw (1, 1) rectangle (3, 2);
\draw (3, 0) rectangle (4, 2);
\draw (4, 0) rectangle (5, 2);
\draw (5, 0) rectangle (7, 1);
\draw (5, 1) rectangle (7, 2);
\draw (7, 0) rectangle (9, 1);
\draw (7, 1) rectangle (9, 2);
\draw (9, 0) rectangle (10, 2);
\draw (10, 0) rectangle (12, 1);
\draw (10, 1) rectangle (12, 2);
\draw (12, 0) rectangle (13, 2);
\draw (13, 0) rectangle (14, 2);
\draw (14, 0) rectangle (15, 2);
\end{tikzpicture}
\]
is one of these $987$ ways.)

\item How many addends do you get when you expand the product \newline$\left(
a+b\right)  \left(  c+d+e\right)  \left(  f+g\right)  $ ? (The answer is $12$.)

\item How many different monomials do you get when you expand the product
$\left(  a-b\right)  \left(  a^{2}+ab+b^{2}\right)  $ ? (This one is more of
an algebra problem, but I wanted to list it because it is connected to
counting. The answer is $2$, because $\left(  a-b\right)  \left(
a^{2}+ab+b^{2}\right)  =a^{3}-b^{3}$.)

\item How many positive divisors does $24$ have? (We can actually list them:
$1,\ 2,\ 3,\ 4,\ 6,\ 8,\ 12,\ 24$. This one is again a mix of a counting
problem and a number theory problem.)
\end{itemize}

We will first solve a few basic counting problems informally, and then (in
Chapter \ref{chp.enum2}) make the underlying concepts rigorous.

\subsection{\label{sec.enum.sets}A refresher on sets}

In prerequisite courses, you have seen basic properties of sets, and basic
notations around sets, but let me quickly remind you of them.

Formally, the notion of a set is fundamental and cannot be defined.

Informally, a \textbf{set} is a collection of objects (which can be anything:
numbers, matrices, functions or other sets) that knows which objects it
contains and which objects it doesn't.

That is, if $S$ is a set and $p$ is any object, then $S$ can either contain
$p$ (in which case we write $p\in S$) or not contain $p$ (in which case we
write $p\notin S$). There is no such thing as \textquotedblleft containing $p$
twice\textquotedblright.

The objects that a set $S$ contains are called the \textbf{elements} of $S$;
they are said to \textbf{belong to }$S$ (or \textbf{lie in }$S$, or \textbf{be
contained in }$S$).

A set can be finite or infinite (i.e., contain finitely or infinitely many
elements). It can be empty (i.e., contain nothing) or nonempty (i.e., contain
at least one element).

An example of a set is the set of all odd integers. This is the set that
contains each odd integer and no other objects. Generally, \textquotedblleft
the set of X\textquotedblright\ means the set that contains X and nothing else.

When a set is finite, it can be written by listing all its elements. For
example, the set of all odd integers between $0$ and $10$ can be written as%
\[
\left\{  1,3,5,7,9\right\}  .
\]
The braces $\{$ and $\}$ around the list are there to signal that we mean the
set of all the elements, not the single elements themselves. These braces are
called \textquotedblleft set braces\textquotedblright, and are involved in
several different notations for sets.

Some more examples of finite sets are%
\begin{align*}
&  \left\{  1,2,3,4,5\right\}  ,\\
&  \left\{  1,2\right\}  ,\\
&  \left\{  1\right\}  \ \ \ \ \ \ \ \ \ \ \left(  \text{this is the set that
only contains }1\right)  ,\\
&  \left\{  {}\right\}  \ \ \ \ \ \ \ \ \ \ \left(  \text{the empty set, also
denoted }\varnothing\right)  ,\\
&  \left\{  1,2,\ldots,1000\right\}  \ \ \ \ \ \ \ \ \ \ \left(  \text{you
understand what \textquotedblleft}\ldots\text{\textquotedblright\ means
here}\right)  .
\end{align*}

Some infinite sets can also be written in this form:%
\begin{align*}
&  \left\{  1,2,3,\ldots\right\}  \ \ \ \ \ \ \ \ \ \ \left(  \text{this is
the set of all positive integers}\right)  ,\\
&  \left\{  0,1,2,\ldots\right\}  \ \ \ \ \ \ \ \ \ \ \left(  \text{this is
the set of all nonnegative integers}\right)  ,\\
&  \left\{  4,5,6,\ldots\right\}  \ \ \ \ \ \ \ \ \ \ \left(  \text{this is
the set of all integers }\geq4\right)  ,\\
&  \left\{  -1,-2,-3,\ldots\right\}  \ \ \ \ \ \ \ \ \ \ \left(  \text{this is
the set of all negative integers}\right)  ,\\
&  \left\{  \ldots,-2,-1,0,1,2,\ldots\right\}  \ \ \ \ \ \ \ \ \ \ \left(
\text{this is the set of all integers}\right)  .
\end{align*}
Some others cannot. For example, how would you list all the real numbers? Or
even all the rational numbers? (These questions are not rhetorical!)

Another way to describe a set is just by putting a description of its elements
in set braces. For example:%
\begin{align*}
&  \left\{  \text{all integers}\right\}  \ \ \ \ \ \ \ \ \ \ \left(
\text{this is the set of all integers}\right)  ,\\
&  \left\{  \text{all integers between }3\text{ and }9\text{ inclusive}%
\right\}  ,\\
&  \left\{  \text{all real numbers}\right\}  .
\end{align*}

Often, you want to define a set that contains all objects of a certain type
that satisfy a certain condition. For example, let's say you want the set of
all integers $x$ that satisfy $x^{2}<13$. There is a notation for this:%
\[
\left\{  x\text{ is an integer }\mid\text{ }x^{2}<13\right\}  .
\]
The vertical bar $\mid$ here should be read as \textquotedblleft such
that\textquotedblright\ (don't mistake it for a divisibility or absolute value
bracket). The part before this bar says what type of objects you are
considering (in our case, it is the integers $x$); the part after this bar
imposes a condition (or several) on these objects (in our case, the condition
is $x^{2}<13$). What you get is the set of all objects of the former type that
satisfy the latter condition. For instance,%
\begin{align*}
&  \left\{  x\text{ is an integer }\mid\text{ }x^{2}<13\right\} \\
&  =\left\{  \text{all integers whose square is smaller than }13\right\} \\
&  =\left\{  -3,-2,-1,0,1,2,3\right\}  .
\end{align*}
Some authors write a colon ($:$) instead of the vertical bar $\mid$. Thus,
they write $\left\{  x\text{ is an integer }\mid\text{ }x^{2}<13\right\}  $ as
$\left\{  x\text{ is an integer }:\text{ }x^{2}<13\right\}  $.

Some sets have standard names:%
\begin{align*}
\mathbb{Z}  &  =\left\{  \text{all integers}\right\}  =\left\{  \ldots
,-2,-1,0,1,2,\ldots\right\}  ;\\
\mathbb{N}  &  =\left\{  \text{all nonnegative integers}\right\}  =\left\{
0,1,2,\ldots\right\} \\
&  \ \ \ \ \ \ \ \ \ \ \left(  \text{beware that some authors use }%
\mathbb{N}\text{ for }\left\{  1,2,3,\ldots\right\}  \text{ instead}\right)
;\\
\mathbb{Q}  &  =\left\{  \text{all rational numbers}\right\}  ;\\
\mathbb{R}  &  =\left\{  \text{all real numbers}\right\}
\ \ \ \ \ \ \ \ \ \ \left(  \text{you barely need them in this course}\right)
;\\
\mathbb{C}  &  =\left\{  \text{all complex numbers}\right\}
\ \ \ \ \ \ \ \ \ \ \left(  \text{you don't need them in this course}\right)
;\\
\varnothing &  =\left\{  {}\right\}  \ \ \ \ \ \ \ \ \ \ \left(  \text{this is
the empty set}\right)  .
\end{align*}

Using these notations, we can rewrite
\[
\left\{  x\text{ is an integer }\mid\text{ }x^{2}<13\right\}
\ \ \ \ \ \ \ \ \ \ \text{as}\ \ \ \ \ \ \ \ \ \ \left\{  x\in\mathbb{Z}\text{
}\mid\text{ }x^{2}<13\right\}  .
\]

Yet another way of defining sets is when you let a variable range over a given
set and collect certain derived quantities. For example,%
\[
\left\{  x^{2}+2\ \mid\ x\in\left\{  1,3,5,7,9\right\}  \right\}
\]
means the set whose elements are the numbers $x^{2}+2$ for all $x\in\left\{
1,3,5,7,9\right\}  $. Thus,%
\begin{align*}
\left\{  x^{2}+2\ \mid\ x\in\left\{  1,3,5,7,9\right\}  \right\}   &
=\left\{  1^{2}+2,\ 3^{2}+2,\ 5^{2}+2,\ 7^{2}+2,\ 9^{2}+2\right\} \\
&  =\left\{  3,\ 11,\ 27,\ 51,\ 83\right\}  .
\end{align*}
In general, if $S$ is a given set, then the notation%
\[
\left\{  \text{an expression }\mid\ x\in S\right\}
\]
stands for the set whose elements are the values of the given expression for
all $x\in S$.

Some more examples of this:%
\begin{align*}
\left\{  \dfrac{x+1}{x}\ \mid\ x\in\left\{  1,2,3,4,5\right\}  \right\}   &
=\left\{  \dfrac{1+1}{1},\ \dfrac{2+1}{2},\ \dfrac{3+1}{3},\ \dfrac{4+1}%
{4},\ \dfrac{5+1}{5}\right\} \\
&  =\left\{  2,\ \dfrac{3}{2},\ \dfrac{4}{3},\ \dfrac{5}{4},\ \dfrac{6}%
{5}\right\}
\end{align*}
and%
\begin{align*}
\left\{  x^{2}\%5\ \mid\ x\in\mathbb{N}\right\}   &  =\left\{  0^{2}%
\%5,\ 1^{2}\%5,\ 2^{2}\%5,\ 3^{2}\%5,\ 4^{2}\%5,\ 5^{2}\%5,\ 6^{2}%
\%5,\ \ldots\right\} \\
&  =\left\{
0,\ 1,\ 4,\ 4,\ 1,\ 0,\ 1,\ 4,\ 4,\ 1,\ 0,\ 1,\ 4,\ 4,\ 1,\ 0,\ \ldots
\right\}  .
\end{align*}
Note that the remainders $x^{2}\%5$ repeat every five steps, because every
integer $x$ satisfies $\left(  x+5\right)  ^{2}\equiv x^{2}\operatorname{mod}%
5$ and thus $\left(  x+5\right)  ^{2}\%5=x^{2}\%5$ (by Proposition
\ref{prop.ent.quo-rem.modd}).

Let me stress once again that a set cannot contain an element more than once.
Also, sets do not come with an ordering of their elements. Thus,%
\[
\left\{  1,2\right\}  =\left\{  2,1\right\}  =\left\{  2,1,1\right\}
=\left\{  1,2,1,2,1\right\}  ,
\]
since each of these four sets contains $1$ and $2$ and nothing else. If $S$ is
a set and $p$ is an object, then $S$ either contains $p$ or does not contain
$p$; it cannot \textquotedblleft contain $p$ twice\textquotedblright, nor can
it contain an element \textquotedblleft before\textquotedblright\ another. So
when you write $\left\{  2,1,1\right\}  $, you aren't making a set that
contains $1$ twice; you are just saying twice that it contains $1$, and this
is equivalent to saying the same thing once. Likewise, the sets $\left\{
1,2\right\}  $ and $\left\{  2,1\right\}  $ do not \textquotedblleft contain
$1$ and $2$ in different orders\textquotedblright; you are just saying in
different orders that they contain $1$ and $2$, but the meaning is the same.
So%
\begin{align*}
\left\{  x^{2}\%5\ \mid\ x\in\mathbb{N}\right\}   &  =\left\{
0,\ 1,\ 4,\ 4,\ 1,\ 0,\ 1,\ 4,\ 4,\ 1,\ 0,\ 1,\ 4,\ 4,\ 1,\ 0,\ \ldots\right\}
\\
&  =\left\{  0,\ 1,\ 4\right\}  .
\end{align*}
This is a finite set, even though $\mathbb{N}$ is infinite!

\bigskip

Note that sets can contain any mathematical objects, not just numbers. In
particular, they can contain other sets. Make sure you understand what the
sets%
\[
\left\{  1,2,3\right\}  ,\ \ \ \ \ \ \ \ \ \ \left\{  \left\{  1,2,3\right\}
\right\}  ,\ \ \ \ \ \ \ \ \ \ \left\{  \left\{  1,2\right\}  ,\left\{
3\right\}  \right\}  ,\ \ \ \ \ \ \ \ \ \ \left\{  \left\{  1\right\}
,\left\{  2\right\}  ,\left\{  3\right\}  \right\}
\]
are and why they are different\footnote{Answer:
\par
\begin{itemize}
\item The set $\left\{  1,2,3\right\}  $ contains three elements, namely the
numbers $1$, $2$ and $3$.
\par
\item The set $\left\{  \left\{  1,2,3\right\}  \right\}  $ contains one
element, namely the set $\left\{  1,2,3\right\}  $.
\par
\item The set $\left\{  \left\{  1,2\right\}  ,\left\{  3\right\}  \right\}  $
contains two elements, namely the sets $\left\{  1,2\right\}  $ and $\left\{
3\right\}  $.
\par
\item The set $\left\{  \left\{  1\right\}  ,\left\{  2\right\}  ,\left\{
3\right\}  \right\}  $ contains three elements, namely the sets $\left\{
1\right\}  $, $\left\{  2\right\}  $ and $\left\{  3\right\}  $.
\end{itemize}
}.

\bigskip

Sets can be compared and combined in several ways:

\begin{definition}
Let $A$ and $B$ be two sets. \medskip

\textbf{(a)} We say that $A$ is a \textbf{subset} of $B$ (and we write
$A\subseteq B$) if every element of $A$ is an element of $B$. \medskip

\textbf{(b)} We say that $A$ is a \textbf{superset} of $B$ (and we write
$A\supseteq B$) if every element of $B$ is an element of $A$. This is
tantamount to saying $B\subseteq A$. \medskip

\textbf{(c)} We say that $A=B$ if the sets $A$ and $B$ contain the same
elements. This is tantamount to saying that both $A\subseteq B$ and
$A\supseteq B$ hold. \medskip

\textbf{(d)} We define the \textbf{union} of $A$ and $B$ to be the set%
\begin{align*}
A\cup B:=  &  \left\{  \text{all elements that are contained in }A\text{ or
}B\right\} \\
=  &  \left\{  x\ \mid\ x\in A\text{ or }x\in B\right\}  .
\end{align*}
(The \textquotedblleft or\textquotedblright\ is non-exclusive, as usual. So
this includes the elements that are contained in both $A$ and $B$.) \medskip

\textbf{(e)} We define the \textbf{intersection} of $A$ and $B$ to be the set%
\begin{align*}
A\cap B:=  &  \left\{  \text{all elements that are contained in both }A\text{
and }B\right\} \\
=  &  \left\{  x\ \mid\ x\in A\text{ and }x\in B\right\}  .
\end{align*}

\textbf{(f)} We define the \textbf{set difference} of $A$ and $B$ to be the
set%
\begin{align*}
A\setminus B:=  &  \left\{  \text{all elements that are contained in }A\text{
but not in }B\right\} \\
=  &  \left\{  x\ \mid\ x\in A\text{ and }x\notin B\right\}  =\left\{  x\in
A\ \mid\ x\notin B\right\}  .
\end{align*}
This is also denoted by $A-B$ by certain authors. \medskip

\textbf{(g)} We say that $A$ and $B$ are \textbf{disjoint} if $A\cap
B=\varnothing$ (that is, $A$ and $B$ have no element in common).
\end{definition}

For example,%
\begin{align*}
\left\{  1,3,5\right\}   &  \subseteq\left\{  1,2,3,4,5\right\}  ,\\
\left\{  1,2,3,4,5\right\}   &  \supseteq\left\{  1,3,5\right\}  ,\\
\text{we \textbf{don't} have }\left\{  5,6,7\right\}   &  \subseteq\left\{
1,2,3,4,5\right\}  ,\\
\left\{  1,2,3\right\}   &  =\left\{  3,2,1\right\}  ,\\
\left\{  1,3,5\right\}  \cup\left\{  3,6\right\}   &  =\left\{
1,3,5,3,6\right\}  =\left\{  1,3,5,6\right\}  ,\\
\left\{  1,3,5\right\}  \cap\left\{  3,6\right\}   &  =\left\{  3\right\}  ,\\
\left\{  1,2,4\right\}  \cap\left\{  3,5\right\}   &  =\varnothing
\ \ \ \ \ \ \ \ \ \ \left(  \text{so that the sets }\left\{  1,2,4\right\}
\text{ and }\left\{  3,5\right\}  \text{ are disjoint}\right)  ,\\
\left\{  1,3,5\right\}  \setminus\left\{  3,6\right\}   &  =\left\{
1,5\right\}  ,\\
\left\{  3,6\right\}  \setminus\left\{  1,3,5\right\}   &  =\left\{
6\right\}  ,\\
\mathbb{Z}\setminus\mathbb{N}  &  =\left\{  -1,-2,-3,\ldots\right\}  =\left\{
\text{all negative integers}\right\}  .
\end{align*}

\begin{definition}
Several sets $A_{1},A_{2},\ldots,A_{k}$ are said to be \textbf{disjoint} if
any two of them (not counting a set and itself) are disjoint, i.e., if we have
$A_{i}\cap A_{j}=\varnothing$ for all $i<j$.
\end{definition}

For example, the three sets $\left\{  1,2\right\}  $, $\left\{  5\right\}  $
and $\left\{  0,7\right\}  $ are disjoint. On the other hand, the three sets
$\left\{  1,2\right\}  $, $\left\{  5\right\}  $ and $\left\{  2,3\right\}  $
are not (since $\left\{  1,2\right\}  \cap\left\{  2,3\right\}  =\left\{
2\right\}  \neq\varnothing$).

\begin{remark}
Beware: \textquotedblleft disjoint\textquotedblright\ is not the same thing as
\textquotedblleft distinct\textquotedblright! Two sets (or any other kinds of
objects) are called \textbf{distinct} if they are not equal. The sets
$\left\{  1,2\right\}  $ and $\left\{  2,3\right\}  $ are distinct but not disjoint.
\end{remark}

\begin{remark}
Unions and intersections can be defined not just for two sets, but for any
number of sets. Given $k$ sets $A_{1},A_{2},\ldots,A_{k}$, their union is%
\[
A_{1}\cup A_{2}\cup\cdots\cup A_{k}:=\left\{  x\ \mid\ x\in A_{i}\text{ for
some }i\in\left\{  1,2,\ldots,k\right\}  \right\}  ,
\]
and their intersection (for $k>0$) is%
\[
A_{1}\cap A_{2}\cap\cdots\cap A_{k}:=\left\{  x\ \mid\ x\in A_{i}\text{ for
all }i\in\left\{  1,2,\ldots,k\right\}  \right\}  .
\]
(The union of $0$ sets is thus $\varnothing$, whereas the intersection of $0$
sets is undefined.)
\end{remark}

\subsection{\label{sec.enum.count-informal}Counting, informally}

Now, let us see how the elements of a set can be counted. Formally speaking,
we will define \textquotedblleft counting\textquotedblright\ later, so we will
play around with not-quite-rigorous concepts for now. As long as we are
working with finite sets, your intuitive understanding of \textquotedblleft
counting\textquotedblright\ should not mislead you.

For example, the set of all odd integers between $0$ and $10$ has $5$ elements
($1,\ 3,\ 5,\ 7,\ 9$), and this doesn't change if you write it redundantly as
$\left\{  1,3,5,5,5,5,7,9\right\}  $. In other words, there are $5$ odd
integers between $0$ and $10$.

More generally, I claim:

\begin{proposition}
\label{prop.count.odd-0-n}Let $n\in\mathbb{N}$. Then, there are exactly
$\left(  n+1\right)  //2=\left\lfloor \dfrac{n+1}{2}\right\rfloor $ odd
integers between $0$ and $n$ (inclusive).
\end{proposition}

\begin{proof}
[Informal proof.]The equality $\left(  n+1\right)  //2=\left\lfloor
\dfrac{n+1}{2}\right\rfloor $ follows from Proposition
\ref{prop.ent.quo-rem-floor}. It remains to show that there are exactly
$\left\lfloor \dfrac{n+1}{2}\right\rfloor $ odd integers between $0$ and $n$.
(We shall always understand the word \textquotedblleft
between\textquotedblright\ to be inclusive, so that $n$ itself is counted if
$n$ is odd.)

We prove this by induction on $n$:

\textit{Base case:} For $n=0$, the claim is true, because there are
$0=\left\lfloor \dfrac{0+1}{2}\right\rfloor $ odd integers between $0$ and $0$.

\textit{Induction step:} Let $n$ be a positive integer. Assume (as the
induction hypothesis) that the claim is true for $n-1$. That is, assume that
there are exactly $\left\lfloor \dfrac{n}{2}\right\rfloor $ odd integers
between $0$ and $n-1$. We must show that the claim also holds for $n$, i.e.,
that there are exactly $\left\lfloor \dfrac{n+1}{2}\right\rfloor $ odd
integers between $0$ and $n$.

Let me introduce a shorthand: The symbol \textquotedblleft\#\textquotedblright%
\ shall mean \textquotedblleft number\textquotedblright. Thus, our induction
hypothesis says%
\begin{equation}
\left(  \text{\# of odd integers between }0\text{ and }n-1\right)
=\left\lfloor \dfrac{n}{2}\right\rfloor , \label{pf.prop.count.odd-0-n.IH}%
\end{equation}
and our goal is to prove that%
\[
\left(  \text{\# of odd integers between }0\text{ and }n\right)  =\left\lfloor
\dfrac{n+1}{2}\right\rfloor .
\]

We are in one of the following two cases:

\textit{Case 1:} The number $n$ is even.

\textit{Case 2:} The number $n$ is odd.

Let us consider Case 1 first. In this case, $n$ is even. Thus, $n$ is not odd.
Therefore, the odd integers between $0$ and $n$ are precisely the odd integers
between $0$ and $n-1$ (since the extra integer $n$ does not qualify as odd).
Hence,%
\begin{align}
&  \left(  \text{\# of odd integers between }0\text{ and }n\right) \nonumber\\
&  =\left(  \text{\# of odd integers between }0\text{ and }n-1\right)
\nonumber\\
&  =\left\lfloor \dfrac{n}{2}\right\rfloor \ \ \ \ \ \ \ \ \ \ \left(
\text{by (\ref{pf.prop.count.odd-0-n.IH})}\right)  .
\label{pf.prop.count.odd-0-n.c1.3}%
\end{align}

However, $n+1$ is odd (since $n$ is even), and thus $2\nmid n+1$. Therefore,
Corollary \ref{cor.ent.quo-rem.floor-bday} \textbf{(b)} (applied to $2$ and
$n+1$ instead of $d$ and $n$) yields $\left\lfloor \dfrac{n+1}{2}\right\rfloor
=\left\lfloor \dfrac{\left(  n+1\right)  -1}{2}\right\rfloor =\left\lfloor
\dfrac{n}{2}\right\rfloor $. Comparing this with
(\ref{pf.prop.count.odd-0-n.c1.3}), we find%
\[
\left(  \text{\# of odd integers between }0\text{ and }n\right)  =\left\lfloor
\dfrac{n+1}{2}\right\rfloor .
\]
Thus, we have achieved our goal in Case 1.

Let us now consider Case 2. In this case, $n$ is odd. Thus, the odd integers
between $0$ and $n$ are precisely the odd integers between $0$ and $n-1$ along
with the new odd integer $n$. Hence,%
\begin{align}
&  \left(  \text{\# of odd integers between }0\text{ and }n\right) \nonumber\\
&  =\left(  \text{\# of odd integers between }0\text{ and }n-1\right)
+1\nonumber\\
&  =\left\lfloor \dfrac{n}{2}\right\rfloor +1\ \ \ \ \ \ \ \ \ \ \left(
\text{by (\ref{pf.prop.count.odd-0-n.IH})}\right)  .
\label{pf.prop.count.odd-0-n.c2.3}%
\end{align}

However, $n+1$ is even (since $n$ is odd), and thus $2\mid n+1$. Therefore,
Corollary \ref{cor.ent.quo-rem.floor-bday} \textbf{(a)} (applied to $2$ and
$n+1$ instead of $d$ and $n$) yields $\left\lfloor \dfrac{n+1}{2}\right\rfloor
=\left\lfloor \dfrac{\left(  n+1\right)  -1}{2}\right\rfloor +1=\left\lfloor
\dfrac{n}{2}\right\rfloor +1$. Comparing this with
(\ref{pf.prop.count.odd-0-n.c2.3}), we find%
\[
\left(  \text{\# of odd integers between }0\text{ and }n\right)  =\left\lfloor
\dfrac{n+1}{2}\right\rfloor .
\]
Thus, we have achieved our goal in Case 2.

So the goal has been achieved in either case, and the induction step is
complete. This proves Proposition \ref{prop.count.odd-0-n}.
\end{proof}

\textbf{Note:} We called the above proof \textquotedblleft
informal\textquotedblright\ because we still don't have a rigorous definition
of the size of a set (i.e., of what \textquotedblleft the number
of\textquotedblright\ means). But we will soon see such a definition. Once we
have learnt this definition and its basic properties, the above proof will
become a formal proof with trivial changes.

\bigskip

Incidentally, it is worth stating the formula for the number of integers (not
just odd integers) in a given interval. Before we state it, let us agree that
if $a$ and $b$ are two integers, then the notation%
\[
\left\{  a,\ a+1,\ a+2,\ \ldots,\ b\right\}  \ \ \ \ \ \ \ \ \ \ \left(
\text{or, shorter, }\left\{  a,\ a+1,\ \ldots,\ b\right\}  \right)
\]
stands for the set of all integers between $a$ and $b$ (inclusive), i.e., the
set \newline$\left\{  x\in\mathbb{Z}\ \mid\ a\leq x\leq b\right\}  $. In
particular, this set is just $\left\{  a\right\}  $ if $a=b$, and is empty if
$a>b$. The following proposition gives its size whenever it is nonempty:

\begin{proposition}
\label{prop.count.fencepost}Let $a,b\in\mathbb{Z}$ be such that $a\leq b+1$.

Then, there are exactly $b-a+1$ numbers in the set $\left\{
a,\ a+1,\ a+2,\ \ldots,\ b\right\}  $. In other words, there are exactly
$b-a+1$ integers between $a$ and $b$ (inclusive).
\end{proposition}

\begin{proof}
[Informal proof.]This is intuitively obvious and can be rigorously proved by
induction on $b$.
\end{proof}

The hard part about Proposition \ref{prop.count.fencepost} is not the proof,
but rather remembering the \textquotedblleft$+1$\textquotedblright! If your
intuition comes from calculus, you think of the interval $\left[  a,b\right]
$ as having length $b-a$ (if $b\geq a$). But since we are doing discrete
mathematics, we are computing not the geometric length of this interval, but
rather the number of integers on this interval, including both endpoints; and
this number is $1$ larger than the length. (For example, if $a=b$, then the
geometric interval $\left[  a,b\right]  $ has zero length, but it contains one
integer, namely $a$.)

\bigskip

It is also worth saying that if two integers $a$ and $b$ satisfy $a\leq b-1$,
then there are exactly $b-a-1$ integers between $a$ and $b$ exclusive (meaning
that we count neither $a$ nor $b$).

\bigskip

\begin{convention}
We agree to use the symbol \textquotedblleft\#\textquotedblright\ for
\textquotedblleft number\textquotedblright.
\end{convention}

\subsection{\label{sec.enum.subsets}Counting subsets}

\subsubsection{Counting them all}

Now, let us count something less trivial than numbers.

How many subsets does the set $\left\{  1,2,3\right\}  $ have? These subsets
are%
\begin{align*}
&  \left\{  {}\right\}  ,\ \ \ \ \ \ \ \ \ \ \left\{  1\right\}
,\ \ \ \ \ \ \ \ \ \ \left\{  2\right\}  ,\ \ \ \ \ \ \ \ \ \ \left\{
3\right\}  ,\\
&  \left\{  1,2\right\}  ,\ \ \ \ \ \ \ \ \ \ \left\{  1,3\right\}
,\ \ \ \ \ \ \ \ \ \ \left\{  2,3\right\}  ,\ \ \ \ \ \ \ \ \ \ \left\{
1,2,3\right\}  .
\end{align*}
(Yes, every set $A$ satisfies $A\subseteq A$ and $\left\{  {}\right\}
\subseteq A$.) Thus, there are $8$ subsets of $\left\{  1,2,3\right\}  $ in total.

Likewise,

\begin{itemize}
\item there are $4$ subsets of $\left\{  1,2\right\}  $, namely $\left\{
{}\right\}  ,\ \left\{  1\right\}  ,\ \left\{  2\right\}  ,\ \left\{
1,2\right\}  $.

\item there are $2$ subsets of $\left\{  1\right\}  $, namely $\left\{
{}\right\}  $ and $\left\{  1\right\}  $.

\item there is $1$ subset of $\left\{  {}\right\}  $, namely $\left\{
{}\right\}  $.

\item there are $16$ subsets of $\left\{  1,2,3,4\right\}  $.
\end{itemize}

The pattern here is hard to miss:\footnote{The expression \textquotedblleft%
$\left\{  1,2,\ldots,n\right\}  $\textquotedblright\ should be read as
$\left\{  1,2\right\}  $ if $n=2$, as $\left\{  1\right\}  $ if $n=1$, and as
the empty set $\left\{  {}\right\}  $ if $n=0$.}

\begin{theorem}
\label{thm.count.subsets.all.n}Let $n\in\mathbb{N}$. Then,%
\[
\left(  \text{\# of subsets of }\left\{  1,2,\ldots,n\right\}  \right)
=2^{n}.
\]

\end{theorem}

\begin{proof}
[Informal proof.]We induct on $n$.

The \textit{base case} ($n=0$) is easy: The set $\left\{  1,2,\ldots
,0\right\}  $ is empty, and thus its only subset is $\left\{  {}\right\}  $
itself; hence, the \# of subsets of $\left\{  1,2,\ldots,0\right\}  $ is
$1=2^{0}$.

\textit{Induction step:} We proceed from $n-1$ to $n$. Thus, let $n$ be a
positive integer. We assume (as the induction hypothesis) that Theorem
\ref{thm.count.subsets.all.n} holds for $n-1$ instead of $n$, and we set out
to prove that it holds for $n$.

So our induction hypothesis says that%
\[
\left(  \text{\# of subsets of }\left\{  1,2,\ldots,n-1\right\}  \right)
=2^{n-1}.
\]
Our goal is to prove that%
\[
\left(  \text{\# of subsets of }\left\{  1,2,\ldots,n\right\}  \right)
=2^{n}.
\]

We define

\begin{itemize}
\item a \textbf{red set} to be a subset of $\left\{  1,2,\ldots,n\right\}  $
that contains $n$;

\item a \textbf{green set} to be a subset of $\left\{  1,2,\ldots,n\right\}  $
that does not contain $n$.
\end{itemize}

For example, if $n=3$, then the red sets are%
\[
\left\{  3\right\}  ,\ \ \ \ \ \ \ \ \ \ \left\{  1,3\right\}
,\ \ \ \ \ \ \ \ \ \ \left\{  2,3\right\}  ,\ \ \ \ \ \ \ \ \ \ \left\{
1,2,3\right\}  ,
\]
whereas the green sets are%
\[
\left\{  {}\right\}  ,\ \ \ \ \ \ \ \ \ \ \left\{  1\right\}
,\ \ \ \ \ \ \ \ \ \ \left\{  2\right\}  ,\ \ \ \ \ \ \ \ \ \ \left\{
1,2\right\}  .
\]

Each subset of $\left\{  1,2,\ldots,n\right\}  $ is either red or green, but
not both. Hence,%
\[
\left(  \text{\# of subsets of }\left\{  1,2,\ldots,n\right\}  \right)
=\left(  \text{\# of red sets}\right)  +\left(  \text{\# of green
sets}\right)  .
\]
(This is an instance of a basic counting principle: If some objects are
classified into two types, then we can count these objects by counting the
objects of each type and adding the results. Later we will state this as a
rigorous theorem, called the \textbf{sum rule for two sets}.)

Thus it remains to count the red sets and the green sets separately.

The green sets are easy: They are just the subsets of $\left\{  1,2,\ldots
,n-1\right\}  $. Hence,%
\[
\left(  \text{\# of green sets}\right)  =\left(  \text{\# of subsets of
}\left\{  1,2,\ldots,n-1\right\}  \right)  =2^{n-1}%
\]
(by the induction hypothesis).

Counting the red sets is trickier, but we can reduce the problem to counting
the green sets: Indeed, the red sets are just the green sets with the element
$n$ inserted into them. To be more precise: Each green set can be turned into
a red set by inserting $n$ into it\footnote{For example, if $n=3$, then the
green set $\left\{  2\right\}  $ becomes $\left\{  2,3\right\}  $ in this
way.}. Conversely, each red set can be turned into a green set by removing the
element $n$ from it. These two operations are mutually inverse, and thus set
up a one-to-one correspondence between the green sets and the red
sets.\footnote{For instance, for $n=3$, it looks like this:%
\[%
\begin{tabular}
[c]{|c||c|c|c|c|}\hline
green set & $\left\{  {}\right\}  $ & $\left\{  1\right\}  $ & $\left\{
2\right\}  $ & $\left\{  1,2\right\}  $\\\hline
& $\updownarrow$ & $\updownarrow$ & $\updownarrow$ & $\updownarrow$\\\hline
red set & $\left\{  3\right\}  $ & $\left\{  1,3\right\}  $ & $\left\{
2,3\right\}  $ & $\left\{  1,2,3\right\}  $\\\hline
\end{tabular}
\ \ .
\]
} This reveals that the \# of red sets is the \# of green sets. Thus,%
\[
\left(  \text{\# of red sets}\right)  =\left(  \text{\# of green sets}\right)
=2^{n-1}.
\]

Combining what we have shown, we now obtain
\begin{align*}
\left(  \text{\# of subsets of }\left\{  1,2,\ldots,n\right\}  \right)   &
=\underbrace{\left(  \text{\# of red sets}\right)  }_{=2^{n-1}}%
+\underbrace{\left(  \text{\# of green sets}\right)  }_{=2^{n-1}}\\
&  =2^{n-1}+2^{n-1}=2\cdot2^{n-1}=2^{n}.
\end{align*}
This is precisely what we needed to prove. This completes the induction step,
and thus Theorem \ref{thm.count.subsets.all.n} is proved.
\end{proof}

More generally, we have the following:

\begin{theorem}
\label{thm.count.subsets.all.S}Let $n\in\mathbb{N}$. Let $S$ be an $n$-element
set. Then,%
\[
\left(  \text{\# of subsets of }S\right)  =2^{n}.
\]

\end{theorem}

\begin{proof}
[Informal proof.]This follows from Theorem \ref{thm.count.subsets.all.n},
since we can rename the $n$ elements of $S$ as $1,2,\ldots,n$.
\end{proof}

For example,
\[
\left(  \text{\# of subsets of }\left\{  \text{\textquotedblleft
cat\textquotedblright, \textquotedblleft dog\textquotedblright,
\textquotedblleft bat\textquotedblright}\right\}  \right)  =2^{3}.
\]

\subsubsection{Counting the subsets of a given size}

Let us now refine our question: Instead of counting all subsets of $\left\{
1,2,\ldots,n\right\}  $, we shall only count the ones that have a given size
$k$. Here, the \textbf{size} of a set means the \# of its elements, i.e., how
many distinct elements it has. (For example, the set $\left\{
1,4,1,15\right\}  $ has size $3$, never mind that I needlessly listed one of
its elements twice.) A set of size $k$ is also known as a $k$\textbf{-element
set}. (Soon we will define these concepts rigorously.)

For instance, $\left\{  1,2,3,4\right\}  $ is a $4$-element set. How many
$2$-element subsets does it have? It has six:%
\[
\left\{  1,2\right\}  ,\ \ \ \ \ \ \ \ \ \ \left\{  1,3\right\}
,\ \ \ \ \ \ \ \ \ \ \left\{  1,4\right\}  ,\ \ \ \ \ \ \ \ \ \ \left\{
2,3\right\}  ,\ \ \ \ \ \ \ \ \ \ \left\{  2,4\right\}
,\ \ \ \ \ \ \ \ \ \ \left\{  3,4\right\}  .
\]

More generally, the answer to the question \textquotedblleft how many
$k$-element subsets does a given $n$-element set have\textquotedblright\ turns
out to be the binomial coefficient $\dbinom{n}{k}$. Let us state this as a
theorem and give an informal proof (which will easily become rigorous once we
have the basic concepts of counting pinned down):\footnote{This theorem is
exactly Theorem \ref{thm.binom.comb-1}, which we left unproved a few chapters
ago.}

\begin{theorem}
\label{thm.count.subsets.size.S}Let $n\in\mathbb{N}$, and let $k$ be any
number (not necessarily an integer). Let $S$ be an $n$-element set. Then,%
\[
\left(  \text{\# of }k\text{-element subsets of }S\right)  =\dbinom{n}{k}.
\]
\footnotemark
\end{theorem}

\footnotetext{Note that $k$-element sets exist only when $k\in\mathbb{N}$.
Thus, if $k\notin\mathbb{N}$, then $\left(  \text{\# of }k\text{-element
subsets of }S\right)  =0$. This fits perfectly with our definition of binomial
coefficients, which ensures that $\dbinom{n}{k}=0$ whenever $k\notin%
\mathbb{N}$.}

\begin{proof}
[Informal proof.]We induct on $n$ (without fixing $k$). That is, we use
induction on $n$ to prove the statement%
\[
P\left(  n\right)  :=\left(
\begin{array}
[c]{c}%
\text{\textquotedblleft for any number }k\text{ and any }n\text{-element set
}S,\\
\text{we have }\left(  \text{\# of }k\text{-element subsets of }S\right)
=\dbinom{n}{k}\text{\textquotedblright}%
\end{array}
\right)
\]
for each $n\in\mathbb{N}$.

\textit{Base case:} Let us prove $P\left(  0\right)  $. Let $k$ be any number.
The only $0$-element set is $\varnothing$, and its only subset is
$\varnothing$. Thus, a $0$-element set $S$ necessarily has one $0$-element
subset ($\varnothing$) and no other subsets. Hence, it satisfies%
\[
\left(  \text{\# of }k\text{-element subsets of }S\right)  =%
\begin{cases}
1, & \text{if }k=0;\\
0, & \text{else.}%
\end{cases}
\]
However, we also have%
\[
\dbinom{0}{k}=%
\begin{cases}
1, & \text{if }k=0;\\
0, & \text{else}%
\end{cases}
\]
(this follows easily from the definition of binomial coefficients). By
comparing these two equalities, we see that any $0$-element set $S$ satisfies
\[
\left(  \text{\# of }k\text{-element subsets of }S\right)  =\dbinom{0}{k}.
\]
In other words, $P\left(  0\right)  $ holds.

\textit{Induction step:} Let $n$ be a positive integer. Assume (as the
induction hypothesis) that $P\left(  n-1\right)  $ holds. We must prove that
$P\left(  n\right)  $ holds.

So we consider any number $k$ and any $n$-element set $S$. We must prove that%
\[
\left(  \text{\# of }k\text{-element subsets of }S\right)  =\dbinom{n}{k}.
\]

We rename the $n$ elements of $S$ as $1,2,\ldots,n$, so we must prove that%
\[
\left(  \text{\# of }k\text{-element subsets of }\left\{  1,2,\ldots
,n\right\}  \right)  =\dbinom{n}{k}.
\]

To prove this, we define

\begin{itemize}
\item a \textbf{red set} to be a $k$-element subset of $\left\{
1,2,\ldots,n\right\}  $ that contains $n$;

\item a \textbf{green set} to be a $k$-element subset of $\left\{
1,2,\ldots,n\right\}  $ that does not contain $n$.
\end{itemize}

For instance:

\begin{itemize}
\item For $n=4$ and $k=2$, the red sets are%
\[
\left\{  1,4\right\}  ,\ \ \ \ \ \ \ \ \ \ \left\{  2,4\right\}
,\ \ \ \ \ \ \ \ \ \ \left\{  3,4\right\}  ,
\]
while the green sets are%
\[
\left\{  1,2\right\}  ,\ \ \ \ \ \ \ \ \ \ \left\{  1,3\right\}
,\ \ \ \ \ \ \ \ \ \ \left\{  2,3\right\}  .
\]

\item For $n=5$ and $k=2$, the red sets are%
\[
\left\{  1,5\right\}  ,\ \ \ \ \ \left\{  2,5\right\}  ,\ \ \ \ \ \left\{
3,5\right\}  ,\ \ \ \ \ \left\{  4,5\right\}  ,
\]
while the green sets are%
\[
\left\{  1,2\right\}  ,\ \ \ \ \ \left\{  1,3\right\}  ,\ \ \ \ \ \left\{
1,4\right\}  ,\ \ \ \ \ \left\{  2,3\right\}  ,\ \ \ \ \ \left\{  2,4\right\}
,\ \ \ \ \ \left\{  3,4\right\}  .
\]

\end{itemize}

Each $k$-element subset of $\left\{  1,2,\ldots,n\right\}  $ is either red or
green (but not both). Hence,%
\begin{align}
&  \left(  \text{\# of }k\text{-element subsets of }\left\{  1,2,\ldots
,n\right\}  \right) \nonumber\\
&  =\left(  \text{\# of red sets}\right)  +\left(  \text{\# of green
sets}\right)  . \label{pf.thm.count.subsets.size.S.sum}%
\end{align}

The green sets are just the $k$-element subsets of $\left\{  1,2,\ldots
,n-1\right\}  $. Thus,%
\begin{align*}
\left(  \text{\# of green sets}\right)   &  =\left(  \text{\# of
}k\text{-element subsets of }\left\{  1,2,\ldots,n-1\right\}  \right) \\
&  =\dbinom{n-1}{k}%
\end{align*}
(by the statement $P\left(  n-1\right)  $, which we have assumed to hold).

Now, let's try to count the red sets.

If $T$ is a red set, then $T\setminus\left\{  n\right\}  $ is a $\left(
k-1\right)  $-element subset of $\left\{  1,2,\ldots,n-1\right\}  $.

Let us refer to the $\left(  k-1\right)  $-element subsets of $\left\{
1,2,\ldots,n-1\right\}  $ as \textbf{blue sets}. Thus, if $T$ is a red set,
then $T\setminus\left\{  n\right\}  $ is a blue set. Conversely, if $U$ is a
blue set, then $U\cup\left\{  n\right\}  $ is a red set. This sets up a
one-to-one correspondence between the red sets and the blue sets: We turn red
sets into blue sets by removing the element $n$, and conversely we turn blue
sets red by inserting the element $n$ into the set.\footnote{For instance, for
$n=4$ and $k=2$, this correspondence looks like this:%
\[%
\begin{tabular}
[c]{|c||c|c|c|}\hline
red set & $\left\{  1,4\right\}  $ & $\left\{  2,4\right\}  $ & $\left\{
3,4\right\}  $\\\hline
& $\updownarrow$ & $\updownarrow$ & $\updownarrow$\\\hline
blue set & $\left\{  1\right\}  $ & $\left\{  2\right\}  $ & $\left\{
3\right\}  $\\\hline
\end{tabular}
\ \ \ .
\]
} Hence,%
\begin{align*}
\left(  \text{\# of red sets}\right)   &  =\left(  \text{\# of blue
sets}\right) \\
&  =\left(  \text{\# of }\left(  k-1\right)  \text{-element subsets of
}\left\{  1,2,\ldots,n-1\right\}  \right) \\
&  \ \ \ \ \ \ \ \ \ \ \ \ \ \ \ \ \ \ \ \ \left(  \text{since this is how the
blue sets were defined}\right) \\
&  =\dbinom{n-1}{k-1}%
\end{align*}
(again by the statement $P\left(  n-1\right)  $, but now applied to $k-1$
instead of $k$). Note that we deliberately did not fix $k$ in our induction,
so that we were now able to apply $P\left(  n-1\right)  $ to $k-1$ instead of
$k$.

Now, (\ref{pf.thm.count.subsets.size.S.sum}) becomes%
\begin{align*}
\left(  \text{\# of }k\text{-element subsets of }\left\{  1,2,\ldots
,n\right\}  \right)   &  =\underbrace{\left(  \text{\# of red sets}\right)
}_{=\dbinom{n-1}{k-1}}+\underbrace{\left(  \text{\# of green sets}\right)
}_{=\dbinom{n-1}{k}}\\
&  =\dbinom{n-1}{k-1}+\dbinom{n-1}{k}=\dbinom{n}{k}%
\end{align*}
by Pascal's recurrence (Theorem \ref{thm.binom.rec}). But this is precisely
the equality that we have to prove. This completes the induction step, and
thus Theorem \ref{thm.count.subsets.size.S} is proved.
\end{proof}

The above proof can also be used to write an algorithm that lists all the
$k$-element subsets of $\left\{  1,2,\ldots,n\right\}  $. This algorithm is
recursive and proceeds as follows:

\begin{itemize}
\item If $n=0$, then:

\begin{itemize}
\item if $k=0$, then list $\varnothing$ (i.e., the resulting list will consist
only of $\varnothing$).

\item otherwise, list nothing.
\end{itemize}

\item Otherwise,

\begin{itemize}
\item list the red sets (by listing all the $\left(  k-1\right)  $-element
subsets of $\left\{  1,2,\ldots,n-1\right\}  $, and inserting $n$ into each of them);

\item list the green sets (i.e., the $k$-element subsets of $\left\{
1,2,\ldots,n-1\right\}  $);

\item combine these two lists.
\end{itemize}
\end{itemize}

In Python, this algorithm (or one possible implementation of it) looks as
follows\footnote{Note that lists are enclosed within brackets in Python: e.g.,
a list that we call $\left(  a,b,c\right)  $ would be written \texttt{[a,b,c]}
in Python. Also, Python's notation \texttt{set([a,b,c])} corresponds to our
$\left\{  a,b,c\right\}  $.}:\medskip\Needspace{12pc}

\texttt{\qquad def subsets(n, k):}

\texttt{\qquad\qquad\# listing all subsets of \{1, 2, ..., n\} that have size
k.}

\texttt{\qquad\qquad if n == 0:}

\texttt{\qquad\qquad\qquad if k == 0:}

\texttt{\qquad\qquad\qquad\qquad return [set([])] \# set([]) is the empty set}

\texttt{\qquad\qquad\qquad return [] \# empty list}

\texttt{\qquad\qquad\# Now, the case when n is not 0:}

\texttt{\qquad\qquad green\_sets = subsets(n-1, k)}

\texttt{\qquad\qquad\# This is the list of all green sets.}

\texttt{\qquad\qquad red\_sets = [U.union([n]) for U in subsets(n-1, k-1)]}

\texttt{\qquad\qquad\# This is the list of all red sets. We construct it by}

\texttt{\qquad\qquad\# taking all the (k-1)-element subsets of \{1, 2, ...,
n-1\}}

\texttt{\qquad\qquad\# (i.e., the blue sets), and inserting n into each of}

\texttt{\qquad\qquad\# them.}

\texttt{\qquad\qquad return red\_sets + green\_sets}

\texttt{\qquad\qquad\# In Python, the plus sign can be used to combine two
lists. }\medskip

With this code, \texttt{subsets(4, 2)} yields

\texttt{\qquad\lbrack\{3, 4\}, \{2, 4\}, \{1, 4\}, \{2, 3\}, \{1, 3\}, \{1,
2\}]}

\noindent as an output, and this is indeed a list of all $2$-element subsets
of $\left\{  1,2,\ldots,4\right\}  =\left\{  1,2,3,4\right\}  $. \medskip

Theorem \ref{thm.count.subsets.size.S} is often called the
\textbf{combinatorial interpretation of binomial coefficients}, since it
reveals that the binomial coefficients $\dbinom{n}{k}$ (at least for
$n\in\mathbb{N}$) have a combinatorial meaning (viz., counting $k$-element
subsets of a given $n$-element set). However, it is just one of many such
interpretations, and we will see four others in Chapter \ref{chp.enum2}!

\begin{exercise}
\label{exe.23mt2.3}Let $n\geq2$ be an integer. The symbol \textquotedblleft%
\#\textquotedblright\ means \textquotedblleft number\textquotedblright.
\medskip

\textbf{(a)} Compute the \# of subsets of $\left\{  1,2,\ldots,n\right\}  $
that contain both $1$ and $2$. \medskip

\textbf{(b)} Compute the \# of $3$-element subsets of $\left\{  1,2,\ldots
,n\right\}  $ that contain both $1$ and $2$. \medskip

[To \textquotedblleft compute\textquotedblright\ a number means to find a
closed-form expression for this number (with no summation signs) and to prove
this formula. I expect proofs to be given at the level of detail and rigor
seen in this chapter.]
\end{exercise}

\subsection{\label{sec.enum.tuples}Tuples (aka lists)}

\subsubsection{Definition and disambiguation}

Now that we have mentioned lists, it is time to explain: What is a finite
list? Here is a somewhat awkward definition:

\begin{definition}
\label{def.count.list}A \textbf{finite list} (aka \textbf{tuple}) is a list
consisting of finitely many objects. The objects appear in this list in a
specified order, and they don't have to be distinct.

A finite list is delimited using parentheses: i.e., the list that contains the
objects $a_{1},a_{2},\ldots,a_{n}$ in this order is denoted by $\left(
a_{1},a_{2},\ldots,a_{n}\right)  $.

\textquotedblleft Specified order\textquotedblright\ means that the list has a
well-defined first entry, a well-defined second entry, and so on. Thus, two
lists $\left(  a_{1},a_{2},\ldots,a_{n}\right)  $ and $\left(  b_{1}%
,b_{2},\ldots,b_{m}\right)  $ are considered equal if and only if

\begin{itemize}
\item we have $n=m$, and

\item we have $a_{i}=b_{i}$ for each $i\in\left\{  1,2,\ldots,n\right\}  $.
\end{itemize}
\end{definition}

For example:

\begin{itemize}
\item The lists $\left(  1,2\right)  $ and $\left(  2,1\right)  $ are not
equal (although the sets $\left\{  1,2\right\}  $ and $\left\{  2,1\right\}  $
are equal).

\item The lists $\left(  1,2\right)  $ and $\left(  1,1,2\right)  $ are not
equal (although the sets $\left\{  1,2\right\}  $ and $\left\{  1,1,2\right\}
$ are equal).

\item The lists $\left(  1,2\right)  $ and $\left(  1,2,2\right)  $ are not
equal (although the sets $\left\{  1,2\right\}  $ and $\left\{  1,2,2\right\}
$ are equal).

\item The lists $\left(  1,1,2\right)  $ and $\left(  1,2,2\right)  $ are not
equal (although the sets $\left\{  1,1,2\right\}  $ and $\left\{
1,2,2\right\}  $ are equal).
\end{itemize}

\begin{definition}
\textbf{(a)} The \textbf{length} of a list $\left(  a_{1},a_{2},\ldots
,a_{n}\right)  $ is defined to be the number $n$. \medskip

\textbf{(b)} A list of length $2$ is called a \textbf{pair} (or an
\textbf{ordered pair}). \medskip

\textbf{(c)} A list of length $3$ is called a \textbf{triple}. \medskip

\textbf{(d)} A list of length $4$ is called a \textbf{quadruple}. \medskip

\textbf{(e)} A list of length $n$ is called an $n$\textbf{-tuple}.
\end{definition}

For example, $\left(  1,3,2,2\right)  $ is a list of length $4$ (although it
has only $3$ \textbf{distinct} entries), i.e., a quadruple or a $4$-tuple. For
another example, $\left(  5,8\right)  $ is a pair, i.e., a $2$-tuple.

Note that there is exactly one list of length $0$: the empty list $\left(
{}\right)  $, which contains nothing.

Lists of length $1$ consist of just a single entry. For example, $\left(
3\right)  $ is a list containing only the entry $3$.

\subsubsection{Counting pairs}

Now, let us count some pairs:

\begin{itemize}
\item How many pairs $\left(  a,b\right)  $ are there with $a,b\in\left\{
1,2,3\right\}  $ ? There are nine:%
\begin{align*}
&  \left(  1,1\right)  ,\ \ \ \ \ \ \ \ \ \ \left(  1,2\right)
,\ \ \ \ \ \ \ \ \ \ \left(  1,3\right)  ,\\
&  \left(  2,1\right)  ,\ \ \ \ \ \ \ \ \ \ \left(  2,2\right)
,\ \ \ \ \ \ \ \ \ \ \left(  2,3\right)  ,\\
&  \left(  3,1\right)  ,\ \ \ \ \ \ \ \ \ \ \left(  3,2\right)
,\ \ \ \ \ \ \ \ \ \ \left(  3,3\right)  .
\end{align*}
The fact that there are nine of them is not surprising given how I've laid
them out: They are forming a table with $3$ rows and $3$ columns, where the
row determines the first entry of the pair\footnote{i.e.:
\par
\begin{itemize}
\item The first row contains the pairs that begin with $1$.
\par
\item The second row contains the pairs that begin with $2$.
\par
\item The third row contains the pairs that begin with $3$.
\end{itemize}
} and the column determines the second entry. Thus, their total number is
$3\cdot3=9$.

\item How many pairs $\left(  a,b\right)  $ are there with $a,b\in\left\{
1,2,3\right\}  $ and $a<b$ ? There are three:%
\[
\left(  1,2\right)  ,\ \ \ \ \ \ \ \ \ \ \left(  1,3\right)
,\ \ \ \ \ \ \ \ \ \ \left(  2,3\right)  .
\]

\item How many pairs $\left(  a,b\right)  $ are there with $a,b\in\left\{
1,2,3\right\}  $ and $a=b$ ? Again, three:%
\[
\left(  1,1\right)  ,\ \ \ \ \ \ \ \ \ \ \left(  2,2\right)
,\ \ \ \ \ \ \ \ \ \ \left(  3,3\right)  .
\]

\item How many pairs $\left(  a,b\right)  $ are there with $a,b\in\left\{
1,2,3\right\}  $ and $a>b$ ? Again, three:%
\[
\left(  2,1\right)  ,\ \ \ \ \ \ \ \ \ \ \left(  3,1\right)
,\ \ \ \ \ \ \ \ \ \ \left(  3,2\right)  .
\]

\end{itemize}

Let us generalize this:

\begin{proposition}
\label{prop.count.pairs.nn}Let $n\in\mathbb{N}$. Then: \medskip

\textbf{(a)} The \# of pairs $\left(  a,b\right)  $ with $a,b\in\left\{
1,2,\ldots,n\right\}  $ is $n^{2}$. \medskip

\textbf{(b)} The \# of pairs $\left(  a,b\right)  $ with $a,b\in\left\{
1,2,\ldots,n\right\}  $ and $a<b$ is $1+2+\cdots+\left(  n-1\right)  $.
\medskip

\textbf{(c)} The \# of pairs $\left(  a,b\right)  $ with $a,b\in\left\{
1,2,\ldots,n\right\}  $ and $a=b$ is $n$. \medskip

\textbf{(d)} The \# of pairs $\left(  a,b\right)  $ with $a,b\in\left\{
1,2,\ldots,n\right\}  $ and $a>b$ is $1+2+\cdots+\left(  n-1\right)  $.
\end{proposition}

\begin{proof}
[Informal proof.]\textbf{(a)} These pairs can be arranged in a table with $n$
rows and $n$ columns, where the rows determine the first entry and the columns
determine the second. Here is what this table looks like:%
\[%
\begin{array}
[c]{cccc}%
\left(  1,1\right)  , & \left(  1,2\right)  , & \ldots, & \left(  1,n\right)
,\\
\left(  2,1\right)  , & \left(  2,2\right)  , & \ldots, & \left(  2,n\right)
,\\
\vdots & \vdots & \ddots & \vdots\\
\left(  n,1\right)  , & \left(  n,2\right)  , & \ldots, & \left(  n,n\right)
.
\end{array}
\]
So there are $n\cdot n=n^{2}$ of these pairs. \medskip

\textbf{(b)} In the table we have just shown, a pair $\left(  a,b\right)  $
satisfies $a<b$ if and only if it is placed above the main diagonal (i.e., the
diagonal starting at the northwestern corner and ending at the southeastern
corner of the table). Thus, the \# of such pairs is the \# of cells above the
main diagonal in this table. But this \# is
\[
0+1+2+\cdots+\left(  n-1\right)  ,
\]
because there are $0$ such cells in the first column, $1$ such cell in the
second, $2$ such cells in the third, and so on. Hence,%
\begin{align*}
&  \left(  \text{\# of pairs }\left(  a,b\right)  \text{ with }a,b\in\left\{
1,2,\ldots,n\right\}  \text{ and }a<b\right) \\
&  =0+1+2+\cdots+\left(  n-1\right) \\
&  =1+2+\cdots+\left(  n-1\right)  .
\end{align*}

\textbf{(c)} A pair $\left(  a,b\right)  $ with $a=b$ is just a pair of the
form $\left(  a,a\right)  $, that is, a single element of $\left\{
1,2,\ldots,n\right\}  $ written twice in succession. Counting such pairs is
therefore tantamount to counting single elements of $\left\{  1,2,\ldots
,n\right\}  $; but there are clearly $n$ of them. \medskip

\textbf{(d)} The pairs $\left(  a,b\right)  $ that satisfy $a>b$ are in
one-to-one correspondence with the pairs $\left(  a,b\right)  $ that satisfy
$a<b$: Namely, each former pair becomes a latter pair if we swap its two
entries, and vice versa. Thus, the \# of former pairs equals the \# of latter
pairs. But we have already found (in part \textbf{(b)}) that the \# of latter
pairs is $1+2+\cdots+\left(  n-1\right)  $. Hence, the \# of former pairs is
$1+2+\cdots+\left(  n-1\right)  $ as well.
\end{proof}

Proposition \ref{prop.count.pairs.nn} has a nice consequence: For any
$n\in\mathbb{N}$, we have%
\begin{align*}
n^{2}  &  =\left(  \text{\# of pairs }\left(  a,b\right)  \text{ with }%
a,b\in\left\{  1,2,\ldots,n\right\}  \right)  \ \ \ \ \ \ \ \ \ \ \left(
\text{by Proposition \ref{prop.count.pairs.nn} \textbf{(a)}}\right) \\
&  =\underbrace{\left(  \text{\# of pairs }\left(  a,b\right)  \text{ with
}a,b\in\left\{  1,2,\ldots,n\right\}  \text{ and }a<b\right)  }%
_{\substack{=1+2+\cdots+\left(  n-1\right)  \\\text{(by Proposition
\ref{prop.count.pairs.nn} \textbf{(b)})}}}\\
&  \ \ \ \ \ \ \ \ \ \ +\underbrace{\left(  \text{\# of pairs }\left(
a,b\right)  \text{ with }a,b\in\left\{  1,2,\ldots,n\right\}  \text{ and
}a=b\right)  }_{\substack{=n\\\text{(by Proposition \ref{prop.count.pairs.nn}
\textbf{(c)})}}}\\
&  \ \ \ \ \ \ \ \ \ \ +\underbrace{\left(  \text{\# of pairs }\left(
a,b\right)  \text{ with }a,b\in\left\{  1,2,\ldots,n\right\}  \text{ and
}a>b\right)  }_{\substack{=1+2+\cdots+\left(  n-1\right)  \\\text{(by
Proposition \ref{prop.count.pairs.nn} \textbf{(d)})}}}\\
&  \ \ \ \ \ \ \ \ \ \ \ \ \ \ \ \ \ \ \ \ \left(
\begin{array}
[c]{c}%
\text{since each pair }\left(  a,b\right)  \text{ satisfies either }a<b\text{
or }a=b\text{ or }a>b\text{,}\\
\text{and never more than one of these three conditions}%
\end{array}
\right) \\
&  =\underbrace{\left(  1+2+\cdots+\left(  n-1\right)  \right)  +n}%
_{=1+2+\cdots+n}+\underbrace{\left(  1+2+\cdots+\left(  n-1\right)  \right)
}_{=\left(  1+2+\cdots+n\right)  -n}\\
&  =\left(  1+2+\cdots+n\right)  +\left(  1+2+\cdots+n\right)  -n\\
&  =2\cdot\left(  1+2+\cdots+n\right)  -n.
\end{align*}
Solving this for $1+2+\cdots+n$, we obtain%
\[
1+2+\cdots+n=\dfrac{n^{2}+n}{2}=\dfrac{n\left(  n+1\right)  }{2}.
\]
Thus, we have recovered the Little Gauss formula (Theorem
\ref{thm.ind.1+2+...+n}) by counting pairs. This illustrates the fact that
counting can be used to prove algebraic identities.

\begin{exercise}
How many pairs $\left(  a,b\right)  $ are there with $a\in\left\{
1,2,3\right\}  $ and $b\in\left\{  1,2,3,4,5\right\}  $ ?
\end{exercise}

\begin{proof}
[Solution.]By the same reasoning as in Proposition \ref{prop.count.pairs.nn}
\textbf{(a)}, there are $15$ such pairs, since the pairs can be arranged in a
table with $3$ rows and $5$ columns.
\end{proof}

The same reasoning gives the following more general result:

\begin{theorem}
\label{thm.count.pairs.nm}Let $n,m\in\mathbb{N}$. Let $A$ be an $n$-element
set. Let $B$ be an $m$-element set. Then,%
\[
\left(  \text{\# of pairs }\left(  a,b\right)  \text{ with }a\in A\text{ and
}b\in B\right)  =nm.
\]

\end{theorem}

What about triples?

\begin{theorem}
\label{thm.count.triples.nmp}Let $n,m,p\in\mathbb{N}$. Let $A$ be an
$n$-element set. Let $B$ be an $m$-element set. Let $C$ be a $p$-element set.
Then,%
\[
\left(  \text{\# of triples }\left(  a,b,c\right)  \text{ with }a\in A\text{
and }b\in B\text{ and }c\in C\right)  =nmp.
\]

\end{theorem}

\begin{proof}
[Informal proof.]You can think of these triples as occupying the cells of a
$3$-dimensional table, but this kind of visualization is tricky (and gets even
less reliable when you get to higher dimensions).

A better approach: Re-encode each triple $\left(  a,b,c\right)  $ as a pair
$\left(  \left(  a,b\right)  ,c\right)  $ (a pair whose first entry is itself
a pair). This is a pair whose first entry comes from the set of all pairs
$\left(  a,b\right)  $ with $a\in A$ and $b\in B$, whereas its second entry
comes from $C$. Let $U$ be the set of all pairs $\left(  a,b\right)  $ with
$a\in A$ and $b\in B$. Then, this set $U$ is an $nm$-element set, because
\begin{align*}
\left(  \text{\# of elements of }U\right)   &  =\left(  \text{\# of pairs
}\left(  a,b\right)  \text{ with }a\in A\text{ and }b\in B\right) \\
&  =nm\ \ \ \ \ \ \ \ \ \ \left(  \text{by Theorem \ref{thm.count.pairs.nm}%
}\right)  .
\end{align*}

Now, we have re-encoded each triple $\left(  a,b,c\right)  $ as a pair
$\left(  \left(  a,b\right)  ,c\right)  $ with $\left(  a,b\right)  \in U$ and
$c\in C$. Thus,
\begin{align*}
&  \left(  \text{\# of triples }\left(  a,b,c\right)  \text{ with }a\in
A\text{ and }b\in B\text{ and }c\in C\right) \\
&  =\left(  \text{\# of pairs }\left(  \left(  a,b\right)  ,c\right)  \text{
with }\left(  a,b\right)  \in U\text{ and }c\in C\right) \\
&  =\left(  \text{\# of pairs }\left(  u,c\right)  \text{ with }u\in U\text{
and }c\in C\right) \\
&  =\left(  nm\right)  p
\end{align*}
(by Theorem \ref{thm.count.pairs.nm}, since $U$ is an $nm$-element set while
$C$ is a $p$-element set). In other words,%
\[
\left(  \text{\# of triples }\left(  a,b,c\right)  \text{ with }a\in A\text{
and }b\in B\text{ and }c\in C\right)  =nmp.
\]
This proves Theorem \ref{thm.count.triples.nmp}.
\end{proof}

\subsubsection{Cartesian products}

There is a general notation for sets of pairs:

\begin{definition}
Let $A$ and $B$ be two sets.

The set of all pairs $\left(  a,b\right)  $ with $a\in A$ and $b\in B$ is
denoted by $A\times B$, and is called the \textbf{Cartesian product} (or just
\textbf{product}) of the sets $A$ and $B$.
\end{definition}

For instance, $\left\{  1,2\right\}  \times\left\{  7,8,9\right\}  $ is the
set of all pairs $\left(  a,b\right)  $ with $a\in\left\{  1,2\right\}  $ and
$b\in\left\{  7,8,9\right\}  $. Explicitly, it consists of the following six
pairs:%
\begin{align*}
&  \left(  1,7\right)  ,\ \ \ \ \ \ \ \ \ \ \left(  1,8\right)
,\ \ \ \ \ \ \ \ \ \ \left(  1,9\right)  ,\\
&  \left(  2,7\right)  ,\ \ \ \ \ \ \ \ \ \ \left(  2,8\right)
,\ \ \ \ \ \ \ \ \ \ \left(  2,9\right)  .
\end{align*}
Likewise, the set $\left\{  1,2\right\}  \times\left\{  2,3\right\}  $
consists of the four pairs%
\begin{align*}
&  \left(  1,2\right)  ,\ \ \ \ \ \ \ \ \ \ \left(  1,3\right)  ,\\
&  \left(  2,2\right)  ,\ \ \ \ \ \ \ \ \ \ \left(  2,3\right)  .
\end{align*}

A similar notation exists for sets of triples, of quadruples or of $k$-tuples
in general:

\begin{definition}
Let $A_{1},A_{2},\ldots,A_{k}$ be $k$ sets.

The set of all $k$-tuples $\left(  a_{1},a_{2},\ldots,a_{k}\right)  $ with
$a_{1}\in A_{1}$ and $a_{2}\in A_{2}$ and $\cdots$ and $a_{k}\in A_{k}$ is
denoted by%
\[
A_{1}\times A_{2}\times\cdots\times A_{k},
\]
and is called the \textbf{Cartesian product} (or just \textbf{product}) of the
sets $A_{1},A_{2},\ldots,A_{k}$.
\end{definition}

For example, the set $\left\{  1,2\right\}  \times\left\{  5\right\}
\times\left\{  2,7,6\right\}  $ consists of all triples $\left(  a_{1}%
,a_{2},a_{3}\right)  $ with $a_{1}\in\left\{  1,2\right\}  $ and $a_{2}%
\in\left\{  5\right\}  $ and $a_{3}\in\left\{  2,7,6\right\}  $. One such
triple is $\left(  2,5,2\right)  $; another is $\left(  2,5,6\right)  $. In
total, there are $2\cdot1\cdot3$ such triples (by Theorem
\ref{thm.count.triples.nmp}).

The word \textquotedblleft Cartesian\textquotedblright\ in \textquotedblleft
Cartesian product\textquotedblright\ honors Ren\'{e} Descartes, who observed
that a point in the Euclidean plane can be characterized by its two
coordinates (i.e., a pair of real numbers), whereas a point in space can be
characterized by its three coordinates (i.e., a triple of real numbers). These
two observations allow us to think of the plane as the Cartesian product
$\mathbb{R}\times\mathbb{R}$, and to think of space as the Cartesian product
$\mathbb{R}\times\mathbb{R}\times\mathbb{R}$. \medskip

Using the notation $A\times B$, we can restate Theorem
\ref{thm.count.pairs.nm} as follows:

\begin{theorem}
[product rule for two sets]\label{thm.count.pairs.nm-x}If $A$ is an
$n$-element set, and $B$ is an $m$-element set, then $A\times B$ is an
$nm$-element set.
\end{theorem}

Likewise, we can restate Theorem \ref{thm.count.triples.nmp} as follows:

\begin{theorem}
[product rule for three sets]If $A$ is an $n$-element set, and $B$ is an
$m$-element set, and $C$ is a $p$-element set, then $A\times B\times C$ is an
$nmp$-element set.
\end{theorem}

More generally:

\begin{theorem}
[product rule for $k$ sets]\label{thm.count.triples.nmp-x-n}Let $A_{1}%
,A_{2},\ldots,A_{k}$ be $k$ sets. If each $A_{i}$ is an $n_{i}$-element set,
then $A_{1}\times A_{2}\times\cdots\times A_{k}$ is an $n_{1}n_{2}\cdots
n_{k}$-element set.
\end{theorem}

In other words, when you count $k$-tuples, with each entry coming from a
certain set, the total number is the product of the numbers of options for
each entry.

You can prove Theorem \ref{thm.count.triples.nmp-x-n} by induction on $k$,
using Theorem \ref{thm.count.pairs.nm-x} and the same \textquotedblleft
re-encode a tuple as a nested pair\textquotedblright\ trick that we used in
our proof of Theorem \ref{thm.count.triples.nmp}. We will later come back to
this in more detail.

\begin{remark}
Let $A,B,C$ be three sets. Then, the Cartesian products $A\times B\times C$
and $\left(  A\times B\right)  \times C$ are not literally the same (the
former consists of triples $\left(  a,b,c\right)  $, while the latter consists
of nested pairs $\left(  \left(  a,b\right)  ,c\right)  $), even though they
\textquotedblleft encode the same information\textquotedblright\ (which we
have leveraged in our above proof of Theorem \ref{thm.count.triples.nmp}).
\end{remark}

\begin{exercise}
\label{exe.23mt2.4}Let $n\in\mathbb{N}$. Compute the number of all pairs
$\left(  a,b\right)  \in\left\{  1,2,\ldots,n\right\}  \times\left\{
1,2,\ldots,n\right\}  $ satisfying $a\equiv b\operatorname{mod}2$.

(The answer will depend on whether $n$ is even or odd. You can find a unified
formula using the floor of a number, but you don't have to.)
\end{exercise}

\begin{exercise}
\label{exe.AxB-BxA}Let $A$ and $B$ be two sets. Show that $A\times B=B\times
A$ holds if and only if we have
\[
A=B\text{ or }A=\varnothing\text{ or }B=\varnothing.
\]

\end{exercise}

\subsubsection{Counting strictly increasing tuples (informally)}

In Proposition \ref{prop.count.pairs.nn} \textbf{(b)}, we have seen that for
any given $n\in\mathbb{N}$, the \# of pairs $\left(  a,b\right)  $ of elements
of $\left\{  1,2,\ldots,n\right\}  $ satisfying $a<b$ is%
\[
1+2+\cdots+\left(  n-1\right)  =\dfrac{\left(  n-1\right)  n}{2}=\dbinom{n}%
{2}.
\]

What is the \# of triples $\left(  a,b,c\right)  $ of elements of $\left\{
1,2,\ldots,n\right\}  $ satisfying $a<b<c$ ?

Such a triple $\left(  a,b,c\right)  $ always determines a $3$-element subset
$\left\{  a,b,c\right\}  $ of $\left\{  1,2,\ldots,n\right\}  $ (and yes, this
will really be a $3$-element subset, because $a<b<c$ entails that $a,b,c$ are
distinct). Conversely, any $3$-element subset of $\left\{  1,2,\ldots
,n\right\}  $ becomes a triple $\left(  a,b,c\right)  $ with $a<b<c$ if we
list its elements in increasing order. Thus, the triples $\left(
a,b,c\right)  $ of elements of $\left\{  1,2,\ldots,n\right\}  $ satisfying
$a<b<c$ are just the $3$-element subsets of $\left\{  1,2,\ldots,n\right\}  $
in disguise.\footnote{We are again being informal here. To be more rigorous,
we should be speaking of a one-to-one correspondence between the former
triples and the latter subsets. But it is not yet the time for this pedantry.}
Hence,%
\begin{align*}
&  \left(  \text{\# of triples }\left(  a,b,c\right)  \text{ of elements of
}\left\{  1,2,\ldots,n\right\}  \text{ satisfying }a<b<c\right) \\
&  =\left(  \text{\# of }3\text{-element subsets of }\left\{  1,2,\ldots
,n\right\}  \right) \\
&  =\dbinom{n}{3}\ \ \ \ \ \ \ \ \ \ \left(  \text{by Theorem
\ref{thm.count.subsets.size.S}, applied to }S=\left\{  1,2,\ldots,n\right\}
\text{ and }k=3\right)  .
\end{align*}
More generally, for any $k\in\mathbb{N}$, we have%
\begin{align*}
&  \left(  \text{\# of }k\text{-tuples }\left(  a_{1},a_{2},\ldots
,a_{k}\right)  \text{ of elements of }\left\{  1,2,\ldots,n\right\}  \text{
satisfying }a_{1}<a_{2}<\cdots<a_{k}\right) \\
&  =\dbinom{n}{k}%
\end{align*}
(by a similar argument: these $k$-tuples are just the $k$-element subsets of
$\left\{  1,2,\ldots,n\right\}  $ in disguise). For comparison, if we drop the
\textquotedblleft$a_{1}<a_{2}<\cdots<a_{k}$\textquotedblright\ requirement,
then we have
\begin{align*}
&  \left(  \text{\# of }k\text{-tuples }\left(  a_{1},a_{2},\ldots
,a_{k}\right)  \text{ of elements of }\left\{  1,2,\ldots,n\right\}  \right)
\\
&  =\underbrace{nn\cdots n}_{k\text{ times}}\ \ \ \ \ \ \ \ \ \ \left(
\text{by Theorem \ref{thm.count.triples.nmp-x-n}}\right) \\
&  =n^{k}.
\end{align*}

Other counting problems don't have answers this simple. For instance, it is
not hard to see that every positive integer $k$ satisfies%
\begin{align*}
&  \left(  \text{\# of }k\text{-tuples }\left(  a_{1},a_{2},\ldots
,a_{k}\right)  \text{ of elements of }\left\{  1,2,\ldots,n\right\}  \text{
}\right. \\
&  \ \ \ \ \ \ \ \ \ \ \left.  \text{such that }a_{1}\text{ is the largest
entry}\right) \\
&  =1^{k-1}+2^{k-1}+3^{k-1}+\cdots+n^{k-1},
\end{align*}
but there is no way to express this without a \textquotedblleft$\cdots
$\textquotedblright\ or a $\sum$ sign. For each specific $k$, however, we can
simplify this:%
\begin{align*}
1^{0}+2^{0}+\cdots+n^{0}  &  =\underbrace{1+1+\cdots+1}_{n\text{ times}}=n;\\
1^{1}+2^{1}+\cdots+n^{1}  &  =1+2+\cdots+n=\dfrac{n\left(  n+1\right)  }{2};\\
1^{2}+2^{2}+\cdots+n^{2}  &  =\dfrac{n\left(  n+1\right)  \left(  2n+1\right)
}{6};\\
1^{3}+2^{3}+\cdots+n^{3}  &  =\dfrac{n^{2}\left(  n+1\right)  ^{2}}{4};\\
1^{4}+2^{4}+\cdots+n^{4}  &  =\dfrac{n\left(  2n+1\right)  \left(  n+1\right)
\left(  3n+3n^{2}-1\right)  }{30};\\
&  \ldots.
\end{align*}
Such a closed-form expression for $1^{m}+2^{m}+\cdots+n^{m}$ exists for any
specific value of $m$ (see, e.g., \cite[Lecture 17, Theorem 2.5.3]{22fco} for
how to find it), but it gets messier as $m$ increases. A closed-form
expression for the general case (with $m$ as a variable) would be too much to expect.

\bigskip

In the next two chapters, we will learn what it means for a set to have $n$
elements, and what rules we have actually been using in our above informal
arguments. To do so, we must first get familiar with the concept of
\textbf{maps} (also known as \textbf{functions}).

\newpage

\section{\label{chp.maps}Maps (aka functions)}

\subsection{\label{sec.enum.maps.informal}Functions, informally}

One of the main notions in mathematics is that of a \textbf{function}, aka
\textbf{map}, aka \textbf{mapping}, aka \textbf{transformation}.

Intuitively, a function is a \textquotedblleft black box\textquotedblright%
\ that takes inputs and transforms them into outputs. For example, the
\textquotedblleft$f\left(  t\right)  =t^{2}$\textquotedblright\ function takes
a real number $t$ and outputs its square $t^{2}$.

You can thus think of a function as a rule for producing an output from an
input. This gives the following \textbf{provisional} definition of a function:

\begin{definition}
[Informal definition of a function]\label{def.function.def-informal}Let $X$
and $Y$ be two sets. A \textbf{function} from $X$ to $Y$ is (provisionally) a
rule that transforms each element of $X$ into some element of $Y$.

If this function is called $f$, then the result of applying it to a given
$x\in X$ (that is, the output produced by $f$ when $x$ is the input) will be
called $f\left(  x\right)  $ (or sometimes $fx$).
\end{definition}

This is not a real definition, as it only kicks the can down the road: It
defines \textquotedblleft function\textquotedblright\ in terms of
\textquotedblleft rule\textquotedblright, but what is a rule? But it gives
some good intuition, provided that it is correctly understood. Here are some
comments that should clarify it:

\begin{itemize}
\item A function has to \textquotedblleft work\textquotedblright\ for each
element of $X$. It cannot decline to operate on some elements! Thus,
\textquotedblleft take the reciprocal\textquotedblright\ is not a function
from $\mathbb{R}$ to $\mathbb{R}$, since it does not operate on $0$ (because
$0$ has no reciprocal). However, \textquotedblleft take the
reciprocal\textquotedblright\ is a function from $\mathbb{R}\setminus\left\{
0\right\}  $ to $\mathbb{R}$, since any nonzero real number does have a reciprocal.

\item A function must not be ambiguous. Each input must produce exactly one
output. Thus, \textquotedblleft take your number to some random
power\textquotedblright\ is not a function from $\mathbb{R}$ to $\mathbb{R}$,
since different powers give different results. (There is a \textquotedblleft
multi-valued\textquotedblright\ variant of functions around, but they aren't
called \textquotedblleft functions\textquotedblright.)

\item We write \textquotedblleft$f:X\rightarrow Y$\textquotedblright\ for
\textquotedblleft$f$ is a function from $X$ to $Y$\textquotedblright.

\item Instead of saying \textquotedblleft$f\left(  x\right)  =y$%
\textquotedblright, we can say \textquotedblleft$f$ transforms $x$ into
$y$\textquotedblright\ or \textquotedblleft$f$ sends $x$ to $y$%
\textquotedblright\ or \textquotedblleft$f$ maps $x$ to $y$\textquotedblright%
\ or \textquotedblleft$f$ takes the value $y$ at $x$\textquotedblright\ or
\textquotedblleft$y$ is the value of $f$ at $x$\textquotedblright\ or
\textquotedblleft$y$ is the image of $x$ under $f$\textquotedblright\ or
\textquotedblleft applying $f$ to $x$ yields $y$\textquotedblright\ or
\textquotedblleft$f$ takes $x$ to $y$\textquotedblright\ or \textquotedblleft%
$f:x\mapsto y$\textquotedblright. All of these statements are synonyms.

For instance, if $f$ is the \textquotedblleft take the
square\textquotedblright\ function from $\mathbb{R}$ to $\mathbb{R}$, then
$f\left(  2\right)  =2^{2}=4$, so that $f$ transforms $2$ into $4$, or sends
$2$ to $4$, or takes the value $4$ at $2$, etc., or $f:2\mapsto4$.

Do not confuse the $\rightarrow$ arrow with the $\mapsto$ arrow! The former
arrow is written between the \textbf{sets} $X$ and $Y$, whereas the latter is
written between a specific input and the corresponding output.

\item As the above terminology suggests, the \textbf{value} of a function $f$
at an input $x$ means the corresponding output $f\left(  x\right)  $.

\item The notation%
\begin{align*}
X  &  \rightarrow Y,\\
x  &  \mapsto\left(  \text{some expression involving }x\right)
\end{align*}
(where $X$ and $Y$ are two sets) means \textquotedblleft the function from $X$
to $Y$ that sends each element $x$ of $X$ to the expression to the right of
the \textquotedblleft$\mapsto$\textquotedblright\ symbol\textquotedblright.
Here, the expression can (for example) be $x^{2}$ or $\dfrac{1}{x+4}$ or
$\dfrac{x}{x+2}$.

For example,%
\begin{align*}
\mathbb{R}  &  \rightarrow\mathbb{R},\\
x  &  \mapsto x^{2}%
\end{align*}
is the \textquotedblleft take the square\textquotedblright\ function (sending
each element $x$ of $\mathbb{R}$ to $x^{2}$). For another example,%
\begin{align*}
\mathbb{R}  &  \rightarrow\mathbb{R},\\
x  &  \mapsto\dfrac{x}{\sin x+15}%
\end{align*}
is the function that takes the sine of the input, then adds $15$, then divides
the input by the result. (Note that this is well-defined, since $\sin x+15$ is
never zero and thus the expression $\dfrac{x}{\sin x+15}$ is always
meaningful, so we really get a function from $\mathbb{R}$ to $\mathbb{R}$.)

For yet another example,%
\begin{align*}
\mathbb{R}  &  \rightarrow\mathbb{R},\\
x  &  \mapsto2
\end{align*}
is the function that sends each real number $x$ to $2$; this is an example of
a constant function. (This is a case where our \textquotedblleft expression
involving $x$\textquotedblright\ does not actually contain $x$. This is
perfectly fine; it's just a very simple particular case.)

For yet another example,%
\begin{align*}
\mathbb{Z}  &  \rightarrow\mathbb{Q},\\
x  &  \mapsto2^{x}%
\end{align*}
is a function (sending each integer $x$ to $2^{x}$). Some of its values are
listed in the following table:%
\[%
\begin{tabular}
[c]{|c||c|c|c|c|c|}\hline
$x$ & $-2$ & $-1$ & $0$ & $1$ & $2$\\\hline
$2^{x}$ & $\dfrac{1}{4}$ & $\dfrac{1}{2}$ & $1$ & $2$ & $4$\\\hline
\end{tabular}
\
\]

A more complicated example is the function%
\begin{align*}
\mathbb{Z}  &  \rightarrow\mathbb{Q},\\
x  &  \mapsto%
\begin{cases}
\dfrac{1}{x-1}, & \text{if }x\neq1;\\
5, & \text{if }x=1.
\end{cases}
\end{align*}

\item The notation%
\begin{align*}
f:X  &  \rightarrow Y,\\
x  &  \mapsto\left(  \text{some expression involving }x\right)
\end{align*}
means that we take the function from $X$ to $Y$ that sends each $x\in X$ to
the expression to the right of the \textquotedblleft$\mapsto$%
\textquotedblright\ symbol, and we call this function $f$.

(Or, if a function named $f$ has already been defined, this notation means
that this $f$ \textbf{is} the function from $X$ to $Y$ that sends each $x\in
X$ to the expression to the right of the \textquotedblleft$\mapsto
$\textquotedblright\ symbol.)

For example, if we write%
\begin{align*}
f:\mathbb{R}  &  \rightarrow\mathbb{R},\\
x  &  \mapsto x^{2}+1,
\end{align*}
then $f$ henceforth will denote the function from $\mathbb{R}$ to $\mathbb{R}$
that sends each $x\in\mathbb{R}$ to $x^{2}+1$.

\item If the set $X$ is finite, then a function $f:X\rightarrow Y$ can be
specified by simply listing all its values. For example, I can define a
function $h:\left\{  0,2,4\right\}  \rightarrow\mathbb{N}$ by setting%
\begin{align*}
h\left(  0\right)   &  =92,\\
h\left(  2\right)   &  =20,\\
h\left(  4\right)   &  =92.
\end{align*}
The values here have been chosen at whim, for no particular reason. A function
does not have to be \textquotedblleft natural\textquotedblright\ or
\textquotedblleft meaningful\textquotedblright\ in any way; all it has to do
is transform each element of $X$ into some element of $Y$.

\item If $f$ is a function from $X$ to $Y$, then the sets $X$ and $Y$ are part
of the function. Thus,%
\begin{align*}
g_{1}:\mathbb{Z}  &  \rightarrow\mathbb{Q},\\
x  &  \mapsto2^{x}%
\end{align*}
and%
\begin{align*}
g_{2}:\mathbb{N}  &  \rightarrow\mathbb{Q},\\
x  &  \mapsto2^{x}%
\end{align*}
and%
\begin{align*}
g_{3}:\mathbb{N}  &  \rightarrow\mathbb{N},\\
x  &  \mapsto2^{x}%
\end{align*}
are three distinct functions! We distinguish between them, so that we can
later speak of the \textquotedblleft domain\textquotedblright\ and the
\textquotedblleft target\textquotedblright\ of a function. Namely, the
\textbf{domain} of a function $f:X\rightarrow Y$ is defined to be the set $X$,
whereas the \textbf{target} of a function $f:X\rightarrow Y$ is defined to be
the set $Y$. Thus, the above function $g_{2}$ has target $\mathbb{Q}$, whereas
the function $g_{3}$ has target $\mathbb{N}$. The above function $g_{1}$ has
domain $\mathbb{Z}$, whereas the function $g_{2}$ has domain $\mathbb{N}$.

\item When are two functions equal? In programming, functions are often
understood to be (implemented) algorithms, and two algorithms can be different
even if they compute the same thing. In mathematics, it's different: Only the
domain, the target and the output values matter; the way they are computed
does not (and indeed there might not even be a way to compute them). Two
algorithms that (always) compute the same thing count for one function only.

So when are two functions considered to be equal?

Two functions $f_{1}:X_{1}\rightarrow Y_{1}$ and $f_{2}:X_{2}\rightarrow
Y_{2}$ are said to be \textbf{equal} if and only if%
\begin{align*}
X_{1}  &  =X_{2}\ \ \ \ \ \ \ \ \ \ \text{and}\ \ \ \ \ \ \ \ \ \ Y_{1}%
=Y_{2}\ \ \ \ \ \ \ \ \ \ \text{and}\\
f_{1}\left(  x\right)   &  =f_{2}\left(  x\right)
\ \ \ \ \ \ \ \ \ \ \text{for all }x\in X_{1}.
\end{align*}

An example of two equal functions is%
\begin{align*}
f_{1}:\mathbb{R}  &  \rightarrow\mathbb{R},\\
x  &  \mapsto x^{2}%
\end{align*}
and%
\begin{align*}
f_{2}:\mathbb{R}  &  \rightarrow\mathbb{R},\\
x  &  \mapsto\left\vert x\right\vert ^{2},
\end{align*}
since each $x\in\mathbb{R}$ satisfies $x^{2}=\left\vert x\right\vert ^{2}$.
\end{itemize}

\begin{fineprint}
\begin{itemize}
\item The Caesar ciphers from Section \ref{sec.ent.crypt} can also be viewed
as examples of maps (i.e., functions). Specifically, if $k$ is any integer,
and if we denote the set of all words (including nonsensical ones like
\textquotedblleft OQJCLA\textquotedblright) by $W$, then the Caesar cipher
$\operatorname*{ROT}\nolimits_{k}$ is a map from $W$ to $W$. For instance,
$\operatorname*{ROT}\nolimits_{1}\left(  \text{\textquotedblleft
KITTEN\textquotedblright}\right)  =\ $\textquotedblleft
LJUUFO\textquotedblright.
\end{itemize}
\end{fineprint}

At this point, we have a good idea of what a function is, but the provisional
definition given above (Definition \ref{def.function.def-informal}) wasn't as
precise as we would like. Even worse, the word \textquotedblleft
rule\textquotedblright\ in that definition is still unclear, and prevents us
from dealing with functions that can neither be given by an explicit formula
(such as \textquotedblleft take the square\textquotedblright) nor be specified
by a complete list of values (e.g., since the domain is infinite). Thus, we
need a better definition of a function.

This is what we will do in the present chapter. The trick is to first define
the more general concept of a \textbf{relation}, and then to characterize
functions as relations with a certain property.

\subsection{\label{sec.enum.maps.rel}Relations}

\textbf{Relations} (to be specific: binary relations) are another concept that
you have already seen on myriad examples:

\begin{itemize}
\item The relation $\subseteq$ is a relation between two sets. For example, we
have $\left\{  1,3\right\}  \subseteq\left\{  1,2,3,4\right\}  $ but we don't
have $\left\{  1,5\right\}  \subseteq\left\{  1,2,3,4\right\}  $.

\item The order relations $\leq$ and $<$ and $>$ and $\geq$ are relations
between two integers (or rational numbers, or real numbers). For example,
$1\leq5$ but $1\nleq-1$.

\item The containment relation $\in$ is a relation between an object and a
set. For instance, $3\in\left\{  1,2,3,4\right\}  $ but $5\notin\left\{
1,2,3,4\right\}  $.

\item The divisibility relation $\mid$ is a relation between two integers.

\item The relation \textquotedblleft coprime\textquotedblright\ is a relation
between two integers.

\item In plane geometry, there are lots of relations: \textquotedblleft
parallel\textquotedblright\ (between two lines), \textquotedblleft
perpendicular\textquotedblright\ (between two lines), \textquotedblleft
congruent\textquotedblright\ (between two shapes), \textquotedblleft
similar\textquotedblright, \textquotedblleft directly
similar\textquotedblright, etc.

\item For any given integer $n$, the relation \textquotedblleft congruent
modulo $n$\textquotedblright\ is a relation between two integers. Let me call
it $\overset{n}{\equiv}$. Thus, $a\overset{n}{\equiv}b$ holds if and only if
$a\equiv b\operatorname{mod}n$. For example, $2\overset{3}{\equiv}8$ but
$2\overset{3}{\not \equiv }7$.
\end{itemize}

What do these relations all have in common? They can be applied to pairs of
objects. Applying a relation to a pair of objects gives a statement which can
be true or false. For example, applying the relation \textquotedblleft
coprime\textquotedblright\ to the pair $\left(  5,8\right)  $ yields the
statement \textquotedblleft$5$ is coprime to $8$\textquotedblright, which is
true. Applying it to the pair $\left(  5,10\right)  $ yields the statement
\textquotedblleft$5$ is coprime to $10$\textquotedblright, which is false.

A general relation $R$ relates elements of a set $X$ with elements of a set
$Y$. For any pair $\left(  x,y\right)  \in X\times Y$ (that is, for any pair
consisting of an element $x\in X$ and an element $y\in Y$), we can apply the
relation $R$ to the pair $\left(  x,y\right)  $, obtaining a statement
\textquotedblleft$x\ R\ y$\textquotedblright\ which is either true or false.
To describe this relation $R$, we need to know which pairs $\left(
x,y\right)  \in X\times Y$ do satisfy $x\ R\ y$ and which pairs don't. In
other words, we need to know the \textbf{set} of all pairs $\left(
x,y\right)  \in X\times Y$ that satisfy $x\ R\ y$. For a rigorous definition
of a relation, we simply take the relation $R$ to \textbf{be} this set of
pairs. In other words, we define relations as follows:

\begin{definition}
\label{def.relation}Let $X$ and $Y$ be two sets. A \textbf{relation} from $X$
to $Y$ is a subset of $X\times Y$ (that is, a set of pairs $\left(
x,y\right)  $ with $x\in X$ and $y\in Y$). \medskip

If $R$ is a relation from $X$ to $Y$, and if $\left(  x,y\right)  \in X\times
Y$ is any pair, then

\begin{itemize}
\item we write $x\ R\ y$ if $\left(  x,y\right)  \in R$;

\item we write $x\ \cancel{R} \ y$ if $\left(  x,y\right)  \notin R$.
\end{itemize}
\end{definition}

All the relations we have seen so far can be recast in terms of this definition:

\begin{itemize}
\item The divisibility relation $\mid$ is a subset of $\mathbb{Z}%
\times\mathbb{Z}$, namely the subset%
\begin{align*}
&  \left\{  \left(  x,y\right)  \in\mathbb{Z}\times\mathbb{Z}\ \mid\ x\text{
divides }y\right\} \\
&  =\left\{  \left(  x,y\right)  \in\mathbb{Z}\times\mathbb{Z}\ \mid
\ \text{there exists some }z\in\mathbb{Z}\text{ such that }y=xz\right\} \\
&  =\left\{  \left(  x,xz\right)  \ \mid\ x\in\mathbb{Z}\text{ and }%
z\in\mathbb{Z}\right\}  .
\end{align*}
For instance, the pairs $\left(  2,4\right)  $ and $\left(  3,9\right)  $ and
$\left(  10,20\right)  $ belong to this subset, whereas the pairs $\left(
2,3\right)  $ and $\left(  2,15\right)  $ and $\left(  10,5\right)  $ do not.

\item The coprimality relation (\textquotedblleft coprime to\textquotedblright%
) is a subset of $\mathbb{Z}\times\mathbb{Z}$, namely the subset%
\begin{align*}
&  \left\{  \left(  x,y\right)  \in\mathbb{Z}\times\mathbb{Z}\ \mid\ x\text{
is coprime to }y\right\} \\
&  =\left\{  \left(  x,y\right)  \in\mathbb{Z}\times\mathbb{Z}\ \mid
\ \gcd\left(  x,y\right)  =1\right\}  .
\end{align*}
It contains, for instance, $\left(  2,3\right)  $ and $\left(  7,9\right)  $,
but not $\left(  4,6\right)  $.

\item For any $n\in\mathbb{Z}$, the \textquotedblleft congruent modulo
$n$\textquotedblright\ relation $\overset{n}{\equiv}$ is a subset of
$\mathbb{Z}\times\mathbb{Z}$, namely the subset%
\begin{align*}
&  \left\{  \left(  x,y\right)  \in\mathbb{Z}\times\mathbb{Z}\ \mid\ x\equiv
y\operatorname{mod}n\right\} \\
&  =\left\{  \left(  x,\ x+nz\right)  \ \mid\ x\in\mathbb{Z}\text{ and }%
z\in\mathbb{Z}\right\}
\end{align*}
(because for a given integer $x$, the integers $y$ that satisfy $x\equiv
y\operatorname{mod}n$ are precisely the integers of the form $x+nz$ for
$z\in\mathbb{Z}$).

\item A geometric example: Let $P$ be the set of all points in the plane, and
let $L$ be the set of all lines in the plane. Then, the \textquotedblleft lies
on\textquotedblright\ relation (as in \textquotedblleft a point lies on a
line\textquotedblright) is a subset of $P\times L$, namely the subset%
\[
\left\{  \left(  p,\ell\right)  \in P\times L\ \mid\ \text{the point }p\text{
lies on the line }\ell\right\}  .
\]

\item If $A$ is any set, then the \textbf{equality relation} on $A$ is the
subset $E_{A}$ of $A\times A$ given by%
\begin{align*}
E_{A}  &  =\left\{  \left(  x,y\right)  \in A\times A\ \mid\ x=y\right\} \\
&  =\left\{  \left(  x,x\right)  \ \mid\ x\in A\right\}  .
\end{align*}
Two elements $x$ and $y$ of $A$ satisfy $x\ E_{A}\ y$ if and only if they are equal.

\item We can literally take any subset of $X\times Y$ (where $X$ and $Y$ are
two sets) and it will be a relation from $X$ to $Y$. Just as with functions, a
relation does not have to follow any \textquotedblleft
meaningful\textquotedblright\ rule. For example, here is a relation from
$\left\{  1,2,3\right\}  $ to $\left\{  5,6,7\right\}  $:%
\[
\left\{  \left(  1,6\right)  ,\ \left(  1,7\right)  ,\ \left(  3,5\right)
\right\}  .
\]
Equivalently, it can be specified by the table%
\[%
\begin{tabular}
[c]{c||ccc}
& $5$ & $6$ & $7$\\\hline\hline
$1$ & no & yes & yes\\
$2$ & no & no & no\\
$3$ & yes & no & no
\end{tabular}
\
\]
(where a \textquotedblleft yes\textquotedblright\ in row $x$ and column $y$
means that $\left(  x,y\right)  $ belongs to the relation). If we call this
relation $R$, then we have $1\ R\ 6$ and $1\ R\ 7$ and $3\ R\ 5$ but not
$1\ R\ 5$ or $2\ R\ 6$.
\end{itemize}

A good way to visualize a relation $R$ from a set $X$ to a set $Y$ (at least
when $X$ and $Y$ are finite) is by drawing the sets $X$ and $Y$ as blobs,
drawing their elements as nodes within these blobs, and drawing an arrow from
the $x$-node to the $y$-node for every pair $\left(  x,y\right)  $ that
belongs to the relation $R$. For example, the relation $R$ in our last example
can be visualized as follows:%
\begin{equation}%
\begin{tikzpicture}
\filldraw[fill=cyan!20!white] (0, 1) ellipse (1 and 1.5);
\node(1) at (0, 2) {$1$};
\node(2) at (0, 1) {$2$};
\node(3) at (0, 0) {$3$};
\node(Cname) at (0, 3) {$X$};
\filldraw[fill=cyan!20!white] (3, 1) ellipse (1 and 1.5);
\node(5) at (3, 2) {$5$};
\node(6) at (3, 1) {$6$};
\node(7) at (3, 0) {$7$};
\node(Dname) at (3, 3) {$Y$};
\draw[->, very thick, green!50!black] (1) -- (6);
\draw[->, very thick, green!50!black] (1) -- (7);
\draw[->, very thick, green!50!black] (3) -- (5);
\end{tikzpicture}
\label{eq.relation.123-567}%
\end{equation}

\begin{remark}
\label{rmk.relation.triple}To make the notion of a relation fit our above
(informal) notion of function, we need to tweak its definition somewhat: We
should define a relation from $X$ to $Y$ not as a mere subset $R$ of $X\times
Y$, but rather as a triple $\left(  R,X,Y\right)  $, where $R$ is a subset of
$X\times Y$. By doing so, we ensure that our relation \textquotedblleft
remembers\textquotedblright\ the sets $X$ and $Y$ and not just the pairs
$\left(  x,y\right)  $ that satisfy $x\ R\ y$.

In the following, we will tacitly understand that a relation is defined in
this way (i.e., as a triple $\left(  R,X,Y\right)  $ and not just as a set
$R$). Nevertheless, we will sloppily (but conveniently) refer to the set $R$
as the relation, and rely on the context to make clear what the sets $X$ and
$Y$ are. For example, the set $\left\{  \left(  1,6\right)  ,\ \left(
1,7\right)  ,\ \left(  3,5\right)  \right\}  $ could be interpreted both as a
relation from $\left\{  1,2,3\right\}  $ to $\left\{  5,6,7\right\}  $ and as
a relation from $\left\{  1,3\right\}  $ to $\left\{  5,6,7,9\right\}  $ (and
in many other ways), and all these relations are different; thus, it is not
\textbf{by itself} a relation, but only becomes a relation when the sets $X$
and $Y$ are provided.
\end{remark}

\subsection{\label{sec.enum.maps.formal}Functions, formally}

We can now define functions rigorously:

\begin{definition}
[Rigorous definition of a function]\label{def.function.def-rig}Let $X$ and $Y$
be two sets. A \textbf{function} from $X$ to $Y$ means a relation $R$ from $X$
to $Y$ that has the following property:

\begin{itemize}
\item \textbf{Output uniqueness:} For each $x\in X$, there exists
\textbf{exactly one} $y\in Y$ such that $x\ R\ y$.
\end{itemize}

If $R$ is a function from $X$ to $Y$, and if $x$ is an element of $X$, then
the unique element $y\in Y$ satisfying $x\ R\ y$ will be called $R\left(
x\right)  $.
\end{definition}

In our above example, the relation%
\[
\left\{  \left(  1,6\right)  ,\ \left(  1,7\right)  ,\ \left(  3,5\right)
\right\}
\]
(which we illustrated in (\ref{eq.relation.123-567})) is not a function from
$\left\{  1,2,3\right\}  $ to $\left\{  5,6,7\right\}  $. In fact, it violates
output uniqueness at $x=1$ (since there are two $y\in\left\{  5,6,7\right\}  $
that satisfy $1\ R\ y$) and also violates it at $x=2$ (since there are no
$y\in\left\{  5,6,7\right\}  $ that satisfy $2\ R\ y$). Each of these two
violations is reason enough to disqualify this relation from being a function.

In our above list of relations, only the equality relation $E_{A}$ is a function.

Here is an example of a function from $X=\left\{  1,2,3\right\}  $ to
$Y=\left\{  5,6,7,8\right\}  $: the relation%
\[
\left\{  \left(  1,7\right)  ,\ \left(  2,5\right)  ,\ \left(  3,7\right)
\right\}  .
\]
This relation satisfies output uniqueness and thus is a function. Visualized
by blobs and arrows, it looks as follows:%
\[%
\begin{tikzpicture}
\filldraw[fill=cyan!20!white] (0, 1) ellipse (1 and 1.5);
\node(1) at (0, 2) {$1$};
\node(2) at (0, 1) {$2$};
\node(3) at (0, 0) {$3$};
\node(Cname) at (0, 3) {$X$};
\filldraw[fill=cyan!20!white] (3, 0.5) ellipse (1 and 2);
\node(5) at (3, 2) {$5$};
\node(6) at (3, 1) {$6$};
\node(7) at (3, 0) {$7$};
\node(8) at (3, -1) {$8$};
\node(Dname) at (3, 3) {$Y$};
\draw[->, very thick, green!50!black] (1) -- (7);
\draw[->, very thick, green!50!black] (2) -- (5);
\draw[->, very thick, green!50!black] (3) -- (7);
\end{tikzpicture}%
\ \ .
\]
If we denote this function by $f$, then $f\left(  1\right)  =7$ and $f\left(
2\right)  =5$ and $f\left(  3\right)  =7$.

Our way of visualizing relations by blobs and arrows makes the output
uniqueness property quite intuitive: This property just says that for each
$x\in X$, there is \textbf{exactly one} arrow starting at the $x$-node. In
other words, each node in the $X$-blob has to be the starting point of exactly
one arrow. Thus, a function is a relation whose visual picture has exactly one
arrow coming out of each $X$-node.

\bigskip

Now we have two definitions of a function: the provisional definition
(Definition \ref{def.function.def-informal}) and the rigorous one (Definition
\ref{def.function.def-rig}). These two definitions are equivalent.
Indeed:\footnote{Keep in mind: As we explained in Remark
\ref{rmk.relation.triple}, a relation $R$ from $X$ to $Y$ \textquotedblleft
remembers\textquotedblright\ the sets $X$ and $Y$. Thus, the same is true for
functions in the sense of the rigorous definition.}

\begin{itemize}
\item If $R$ is a function from $X$ to $Y$ in the sense of the rigorous
definition (i.e., a relation from $X$ to $Y$ that satisfies output
uniqueness), then $R$ can also be viewed as a rule that sends each element $x$
of $X$ to some element of $Y$: namely, to the unique $y\in Y$ that satisfies
$x\ R\ y$. Thus, $R$ becomes a function in the provisional sense.

\item Conversely, if $f$ is a function from $X$ to $Y$ in the provisional
sense (i.e., a rule sending elements of $X$ to elements of $Y$), then $f$ can
also be viewed as a function in the rigorous sense (i.e., as a relation from
$X$ to $Y$ that satisfies output uniqueness), as follows: Let $R$ be the set%
\[
\left\{  \left(  x,\ f\left(  x\right)  \right)  \ \mid\ x\in X\right\}  .
\]
This set $R$ is a subset of $X\times Y$, that is, a relation from $X$ to $Y$.
(In a more intuitive language, this relation $R$ is characterized as follows:
Two elements $x\in X$ and $y\in Y$ satisfy $x\ R\ y$ if and only if
$y=f\left(  x\right)  $. That is, roughly speaking, the relation $R$ relates
each input $x\in X$ with the corresponding output value $f\left(  x\right)
\in Y$ and with nothing else.) This relation $R$ satisfies output uniqueness
(because each input $x\in X$ produces exactly one output value $f\left(
x\right)  $), and therefore is a function from $X$ to $Y$ in the rigorous
sense. Thus, $f$ becomes a function in the rigorous sense (namely, the
rigorous function $R$).
\end{itemize}

Therefore, we can translate rigorous functions into provisional ones and vice
versa. We thus shall think of the two concepts as being the same (i.e., we
will regard the rigorous concept as a clarification of the provisional one).
In particular, all the notations we have introduced for provisional functions
will be used for rigorous ones.

\subsection{\label{sec.enum.maps.more-exas}Some more examples of functions}

Let us give some examples of functions as well as some examples of what looks
like functions but are not.

\begin{example}
Consider the function
\[
f_{0}:\left\{  1,2,3,4\right\}  \rightarrow\left\{  1,2,3,4\right\}
\]
that sends $1,2,3,4$ to $3,2,3,3$, respectively. As a rigorous function, it is
the relation $R$ that satisfies%
\[
1\ R\ 3,\ \ \ \ \ \ \ \ \ \ 2\ R\ 2,\ \ \ \ \ \ \ \ \ \ 3\ R\ 3,\ \ \ \ \ \ \ \ \ \ 4\ R\ 3
\]
(and nothing else). In other words, it is the relation%
\[
\left\{  \left(  1,3\right)  ,\ \left(  2,2\right)  ,\ \left(  3,3\right)
,\ \left(  4,3\right)  \right\}  .
\]

\end{example}

\begin{example}
What about the function%
\begin{align*}
f_{1}:\left\{  1,2,3,4\right\}   &  \rightarrow\left\{  1,2,3\right\}  ,\\
n  &  \mapsto n\text{\ \ \ \ ?}%
\end{align*}
Such a function $f_{1}$ does not exist, since it would have to send $4$ to
$4$, but $4$ is not in the target $\left\{  1,2,3\right\}  $.

This is a pedantic issue, but it should be kept in mind: Not every expression
that appears to define a function actually defines a function. Make sure that
the expression to the right of the \textquotedblleft$\mapsto$%
\textquotedblright\ symbol always is an actual element of the target (which,
in this case, is the set $\left\{  1,2,3\right\}  $).
\end{example}

\begin{example}
Consider the function%
\begin{align*}
f_{2}:\left\{  1,2,3,\ldots\right\}   &  \rightarrow\left\{  1,2,3,\ldots
\right\}  ,\\
n  &  \mapsto\left(  \text{the number of positive divisors of }n\right)  .
\end{align*}
As a relation, it is%
\[
\left\{  \left(  1,1\right)  ,\ \left(  2,2\right)  ,\ \left(  3,2\right)
,\ \left(  4,3\right)  ,\ \left(  5,2\right)  ,\ \left(  6,4\right)
,\ \left(  7,2\right)  ,\ \left(  8,4\right)  ,\ \left(  9,3\right)
,\ \ldots\right\}  .
\]
(We cannot list all the pairs, since there are infinitely many.) Thus,
$f_{2}\left(  1\right)  =1$ and $f_{2}\left(  2\right)  =2$ and $f_{2}\left(
3\right)  =2$ and so on.

Note that every prime $p$ satisfies $f_{2}\left(  p\right)  =2$. But the only
$n\in\left\{  1,2,3,\ldots\right\}  $ that satisfies $f_{2}\left(  n\right)
=1$ is $1$.
\end{example}

\begin{example}
What about the function%
\begin{align*}
\widetilde{f}_{2}:\mathbb{Z}  &  \rightarrow\left\{  1,2,3,\ldots\right\}  ,\\
n  &  \mapsto\left(  \text{the number of positive divisors of }n\right)
\ \ \ \ \text{?}%
\end{align*}
There is no such function $\widetilde{f}_{2}$, since $\widetilde{f}_{2}\left(
0\right)  $ would have to be undefined or $\infty$ (because $0$ has infinitely
many positive divisors).

This is the exact same problem that we had with the non-function $f_{1}$ above.
\end{example}

\begin{example}
What about the function%
\begin{align*}
f_{3}:\left\{  1,2,3,\ldots\right\}   &  \rightarrow\left\{  1,2,3,\ldots
\right\}  ,\\
n  &  \mapsto\left(  \text{the smallest prime divisor of }n\right)
\ \ \ \ \text{?}%
\end{align*}
Again, there is no such function $f_{3}$, since $f_{3}\left(  1\right)  $
makes no sense (indeed, the number $1$ has no prime divisors, thus no smallest
prime divisor).

This is essentially the same problem as with the function $\widetilde{f}_{2}$
from the previous example, except that this time the value $f_{3}\left(
1\right)  $ is really undefined (as opposed to just failing to belong to the target).

Note that the function $f_{3}$ \textquotedblleft almost\textquotedblright%
\ exists: There is a relation \textquotedblleft$y$ is the smallest prime
divisor of $x$\textquotedblright\ from $\left\{  1,2,3,\ldots\right\}  $ to
$\left\{  1,2,3,\ldots\right\}  $, but this relation fails the output
uniqueness requirement at $x=1$, and thus is not a function. However, we can
make it into a function by removing the offending element $1$ from its domain.
That is, there is a function%
\begin{align*}
\widetilde{f}_{3}:\left\{  2,3,4,\ldots\right\}   &  \rightarrow\left\{
1,2,3,\ldots\right\}  ,\\
n  &  \mapsto\left(  \text{the smallest prime divisor of }n\right)  .
\end{align*}

\end{example}

\begin{example}
What about the function%
\begin{align*}
f_{4}:\mathbb{Q}  &  \rightarrow\mathbb{Z},\\
\dfrac{a}{b}  &  \mapsto a\ \ \ \ \ \left(  \text{for }a,b\in\mathbb{Z}\text{
with }b\neq0\right)  \text{ ?}%
\end{align*}
Restated in words, this is to be a function that takes a rational number as
input, writes it as a ratio of two integers and outputs the numerator. Is
there such a function?

Again, the answer is \textbf{no}. Again, the problem is a failure of output
uniqueness, but this time, it fails not because the output does not exist (or
does not belong to the target), but rather because the output is non-unique.
For example, if $f_{4}$ were a function, then we would have the two
equalities
\begin{align*}
f_{4}\left(  0.5\right)   &  =f_{4}\left(  \dfrac{1}{2}\right)
=1\ \ \ \ \ \ \ \ \ \ \text{and}\\
f_{4}\left(  0.5\right)   &  =f_{4}\left(  \dfrac{3}{6}\right)  =3,
\end{align*}
which contradict one another. The underlying issue is that a rational number
can be written as a fraction in several different ways, and the numerators of
these fractions will usually \textbf{not be the same}. Thus, if you follow the
rule $\dfrac{a}{b}\mapsto a$ to compute the output of $f_{4}$ for a given
input, your output will depend on how exactly you write your input as a
fraction, and this is a violation of output uniqueness.

(We could fix this issue by requiring in the definition of $f_{4}$ that the
fraction $\dfrac{a}{b}$ be given in reduced form -- specifically, that $b>0$
and $\gcd\left(  a,b\right)  =1$. It is not hard to show that each rational
number can be written in such a form in exactly one way, so that there could
not be any different outputs for the same input.)
\end{example}

\subsection{\label{sec.enum.maps.welldef}Well-definedness}

The issues that we have seen in the last few examples (supposed functions
failing to exist either because their output values make no sense, or because
these values don't lie in $Y$, or because these values are ambiguous) are
known as \textbf{well-definedness} issues. Often, mathematicians say that
\textquotedblleft a function is well-defined\textquotedblright\ when they mean
that its definition does not suffer from such issues (i.e., its definition
really defines a function). So you should read \textquotedblleft This function
is well-defined [or: not well-defined]\textquotedblright\ as \textquotedblleft
The definition we just gave really defines a function [or: does not actually
define a function]\textquotedblright.

For example, as we just saw, the function%
\begin{align*}
f_{4}:\mathbb{Q}  &  \rightarrow\mathbb{Z},\\
\dfrac{a}{b}  &  \mapsto a
\end{align*}
is not well-defined (i.e., there is no such function), but the function%
\begin{align*}
f_{5}:\mathbb{Q}  &  \rightarrow\mathbb{Q},\\
\dfrac{a}{b}  &  \mapsto\dfrac{a^{2}}{b^{2}}%
\end{align*}
is well-defined (because if you write a given rational number as $\dfrac{a}%
{b}$ for different pairs $\left(  a,b\right)  $, the resulting quotients
$\dfrac{a^{2}}{b^{2}}$ will all be equal). The function%
\begin{align*}
f_{1}:\left\{  1,2,3,4\right\}   &  \rightarrow\left\{  1,2,3\right\}  ,\\
n  &  \mapsto n
\end{align*}
is not well-defined (since its supposed output $f_{1}\left(  4\right)  $ fails
to lie in the target $\left\{  1,2,3\right\}  $), whereas the function%
\begin{align*}
f_{6}:\left\{  1,2,3,4\right\}   &  \rightarrow\left\{  1,2,3\right\}  ,\\
n  &  \mapsto1+\left(  n\%3\right)
\end{align*}
is well-defined (since its outputs at $1,2,3,4$ are $2,3,1,2$).

Thus, in order to show that a function specified by an expression is
well-defined (i.e., that the expression really does define a function), we can
follow a three-point checklist:

\begin{enumerate}
\item Prove that the supposed outputs exist (e.g., that the expression
involves no division by zero).

\item Prove that the supposed outputs lie in the target (e.g., that the
expression does not evaluate to $4$ when the target is $\left\{
1,2,3\right\}  $).

\item Prove that the supposed outputs are determined uniquely by the
respective inputs (and not by how the inputs are provided).
\end{enumerate}

For instance, in our above examples, the \textquotedblleft
function\textquotedblright\ $f_{3}$ fails point 1 of this checklist (since
$f_{3}\left(  1\right)  $ does not exist); the \textquotedblleft
function\textquotedblright\ $f_{1}$ fails point 2; and the \textquotedblleft
function\textquotedblright\ $f_{4}$ fails point 3.

\begin{example}
Let us say we want to define a function%
\begin{align*}
f_{7}:\mathbb{Q}  &  \rightarrow\mathbb{Q},\\
x  &  \mapsto%
\begin{cases}
\dfrac{x}{x-1}, & \text{if }x\leq0;\\
\dfrac{x}{x+1}, & \text{if }x\geq0.
\end{cases}
\end{align*}
Is this function $f_{7}$ well-defined (i.e., does it really exist)? To answer
this, we follow the above checklist:

\begin{enumerate}
\item The supposed outputs exist: We need to check that the fraction
$\dfrac{x}{x-1}$ is well-defined for all $x\leq0$, and that the fraction
$\dfrac{x}{x+1}$ is well-defined for all $x\geq0$. This is easy: We just need
to ensure that we don't divide by $0$; but in fact our denominator $x-1$ is
negative for $x\leq0$, while our denominator $x+1$ is positive for $x\geq0$.

\item The supposed outputs lie in the target: We need to check that $\dfrac
{x}{x-1}\in\mathbb{Q}$ when $x\in\mathbb{Q}$ is $\leq0$, and that $\dfrac
{x}{x+1}\in\mathbb{Q}$ when $x\in\mathbb{Q}$ is $\geq0$. But this is easy,
since addition, subtraction and division of rational numbers cannot
\textquotedblleft take us out\textquotedblright\ of the set $\mathbb{Q}$.

\item The supposed outputs are determined uniquely by the respective inputs:
This is not completely obvious here, since the two cases \textquotedblleft%
$x\leq0$\textquotedblright\ and \textquotedblleft$x\geq0$\textquotedblright%
\ have a bit of overlap (the input $0$ falls into them both), and thus we have
provided two different expressions for the output at $x=0$. However, it turns
out that both expressions define the same value (indeed, if $x=0$, then both
$\dfrac{x}{x-1}$ and $\dfrac{x}{x+1}$ are $0$), and thus the output is
uniquely determined even when it is given by two different expressions.
\end{enumerate}

Thus, all three checks come out positive, and it follows that $f_{7}$ exists.
\end{example}

This example is a bit unusual; in most situations, one or more of the items of
our checklist are self-evident. For instance, for the function%
\begin{align*}
f:\mathbb{Q}  &  \rightarrow\mathbb{Q},\\
x  &  \mapsto\dfrac{1}{x^{2}-2},
\end{align*}
items 2 and 3 on the checklist are trivial (a fraction made of rational
numbers always yields a rational number, and clearly there can be no ambiguity
when the input is just taken as a single variable $x$), but item 1 requires
some work (it boils down to checking that $x^{2}-2$ is never $0$ when
$x\in\mathbb{Q}$; this is a consequence of the fact that $\sqrt{2}$ is
irrational\footnote{We observed this fact after Exercise \ref{exe.pval.e2}.}).

\subsection{\label{sec.enum.maps.id}The identity function}

\begin{definition}
For any set $A$, there is an \textbf{identity function} $\operatorname*{id}%
\nolimits_{A}:A\rightarrow A$. This is the function that sends each element
$a\in A$ to $a$ itself. In other words, it is precisely the relation $E_{A}$
defined in Section \ref{sec.enum.maps.rel}.
\end{definition}

Here is the blobs-and-arrows visualization of the identity function
$\operatorname*{id}\nolimits_{A}$ for $A=\left\{  1,2,3\right\}  $:%
\[%
\begin{tikzpicture}
\filldraw[fill=cyan!20!white] (0, 1) ellipse (1 and 1.5);
\node(1) at (0, 2) {$1$};
\node(2) at (0, 1) {$2$};
\node(3) at (0, 0) {$3$};
\node(Cname) at (0, 3) {$X$};
\filldraw[fill=cyan!20!white] (3, 1) ellipse (1 and 1.5);
\node(5) at (3, 2) {$1$};
\node(6) at (3, 1) {$2$};
\node(7) at (3, 0) {$3$};
\node(Dname) at (3, 3) {$Y$};
\draw[->, very thick, green!50!black] (1) -- (5);
\draw[->, very thick, green!50!black] (2) -- (6);
\draw[->, very thick, green!50!black] (3) -- (7);
\end{tikzpicture}%
\]

\subsection{More examples, and multivariate functions}

As we said before, a function $f:X\rightarrow Y$ can be described either by a
rule or by a list of values (if $X$ is finite) or as a relation. For instance,
the \textquotedblleft take the square\textquotedblright\ function on real
numbers is the function%
\begin{align*}
f:\mathbb{R}  &  \rightarrow\mathbb{R},\\
x  &  \mapsto x^{2}.
\end{align*}
As a relation, it is the set%
\[
\left\{  \left(  x,x^{2}\right)  \ \mid\ x\in\mathbb{R}\right\}  .
\]

When the domain of a function $f$ is a Cartesian product of several sets
(i.e., the inputs of $f$ are tuples), $f$ is called a \textbf{multivariate}
function. For instance, the function%
\begin{align*}
f:\mathbb{Z}\times\mathbb{Z}  &  \rightarrow\mathbb{Z},\\
\left(  a,b\right)   &  \mapsto a+b
\end{align*}
(which sends each pair $\left(  a,b\right)  $ of two integers to their sum
$a+b$) is a multivariate function. Its input is a pair of two integers, i.e.,
it really has two inputs ($a$ and $b$). As a relation, it is the subset%
\begin{align*}
&  \left\{  \left(  \left(  a,b\right)  ,\ a+b\right)  \ \mid\ a,b\in
\mathbb{Z}\right\} \\
&  =\left\{  \left(  \left(  a,b\right)  ,c\right)  \ \mid\ a,b,c\in
\mathbb{Z}\text{ such that }c=a+b\right\}
\end{align*}
of $\left(  \mathbb{Z}\times\mathbb{Z}\right)  \times\mathbb{Z}$. Of course,
this function has a name: It is the addition of integers. Other multivariate
functions are
\begin{align*}
\mathbb{Z}\times\mathbb{Z}  &  \rightarrow\mathbb{Z},\\
\left(  a,b\right)   &  \mapsto a-b
\end{align*}
(known as the subtraction of integers) and%
\begin{align*}
\mathbb{Z}\times\mathbb{Z}  &  \rightarrow\mathbb{Z},\\
\left(  a,b\right)   &  \mapsto ab
\end{align*}
(the multiplication of integers), as well as similar functions defined for
other sets of numbers. Keep in mind that there is no \textquotedblleft
division\textquotedblright\ function%
\begin{align*}
\mathbb{Z}\times\mathbb{Z}  &  \rightarrow\mathbb{Z},\\
\left(  a,b\right)   &  \mapsto a/b,
\end{align*}
since $a/b$ is not always an integer (and does not even exist when $b=0$).

When $f$ is a multivariate function whose inputs are $k$-tuples, we commonly
use the shorthand notation $f\left(  a_{1},a_{2},\ldots,a_{k}\right)  $ for
its values $f\left(  \left(  a_{1},a_{2},\ldots,a_{k}\right)  \right)  $.
(That is, we commonly omit the outer pair of parentheses.) For instance, if
$f$ is the addition of integers, then $f\left(  a,b\right)  =f\left(  \left(
a,b\right)  \right)  =a+b$ for all $a,b\in\mathbb{Z}$.

\subsection{\label{sec.enum.maps.fog}Composition of functions}

\subsubsection{\label{subsec.enum.maps.fog-def}Definition}

There are some ways to transform functions into other functions. The most
important one is \textbf{composition}:

\begin{definition}
\label{def.function.fog}Let $X$, $Y$ and $Z$ be three sets. Let
$f:Y\rightarrow Z$ and $g:X\rightarrow Y$ be two functions. Then, $f\circ g$
denotes the function%
\begin{align*}
X  &  \rightarrow Z,\\
x  &  \mapsto f\left(  g\left(  x\right)  \right)  .
\end{align*}
In other words, $f\circ g$ is the function that first applies $g$ and then
applies $f$. This function $f\circ g$ is called the \textbf{composition} of
$f$ with $g$ (and I pronounce it \textquotedblleft$f$ after $g$%
\textquotedblright).
\end{definition}

In terms of relations, if we view $f$ and $g$ as two relations $F$ and $G$ (as
in Definition \ref{def.function.def-rig}), then $f\circ g$ is the relation%
\[
\left\{  \left(  x,z\right)  \ \mid\ \text{there exists }y\in Y\text{ such
that }x\ G\ y\text{ and }y\ F\ z\right\}  \text{ from }X\text{ to }Z.
\]

\begin{example}
Consider the two functions
\begin{align*}
f:\mathbb{R}  &  \rightarrow\mathbb{R},\\
x  &  \mapsto x^{3}%
\end{align*}
and%
\begin{align*}
g:\mathbb{R}  &  \rightarrow\mathbb{R},\\
x  &  \mapsto\dfrac{1}{x^{2}+7}.
\end{align*}
Then, for any real $x\in\mathbb{R}$, we have%
\[
\left(  f\circ g\right)  \left(  x\right)  =f\left(  g\left(  x\right)
\right)  =f\left(  \dfrac{1}{x^{2}+7}\right)  =\left(  \dfrac{1}{x^{2}%
+7}\right)  ^{3}%
\]
whereas%
\[
\left(  g\circ f\right)  \left(  x\right)  =g\left(  f\left(  x\right)
\right)  =g\left(  x^{3}\right)  =\dfrac{1}{\left(  x^{3}\right)  ^{2}%
+7}=\dfrac{1}{x^{6}+7}.
\]
Note that these two results are different. Thus, $f\circ g\neq g\circ f$ in general.
\end{example}

\begin{example}
Consider the two functions $f:\left\{  1,2,3\right\}  \rightarrow\left\{
1,2,3,4\right\}  $ and $g:\left\{  1,2,3,4\right\}  \rightarrow\left\{
1,2,3\right\}  $ given by the following tables of values:%
\[%
\begin{tabular}
[c]{|c||c|c|c|}\hline
$i$ & $1$ & $2$ & $3$\\\hline
$f\left(  i\right)  $ & $1$ & $3$ & $2$\\\hline
\end{tabular}
\ \ \ \ \ \ \ \ \ \ \
\begin{tabular}
[c]{|c||c|c|c|c|}\hline
$i$ & $1$ & $2$ & $3$ & $4$\\\hline
$g\left(  i\right)  $ & $2$ & $1$ & $3$ & $2$\\\hline
\end{tabular}
\ \ \ .
\]

These two functions can be visualized using blobs and arrows, and we can even
reuse the target-blob from $g$ as the domain-blob for $f$:
\[%
\begin{tikzpicture}
\filldraw[fill=cyan!20!white] (0, 0.5) ellipse (1 and 2);
\node(1) at (0, 2) {$1$};
\node(2) at (0, 1) {$2$};
\node(3) at (0, 0) {$3$};
\node(4) at (0, -1) {$4$};
\node(Cname) at (0, 3) {$\set{1,2,3,4}$};
\filldraw[fill=cyan!20!white] (3, 1) ellipse (1 and 1.5);
\node(11) at (3, 2) {$1$};
\node(22) at (3, 1) {$2$};
\node(33) at (3, 0) {$3$};
\node(Dname) at (3, 3) {$\set{1,2,3}$};
\filldraw[fill=cyan!20!white] (6, 0.5) ellipse (1 and 2);
\node(111) at (6, 2) {$1$};
\node(222) at (6, 1) {$2$};
\node(333) at (6, 0) {$3$};
\node(444) at (6, -1) {$4$};
\node(Ename) at (6, 3) {$\set{1,2,3,4}$};
\draw[->, very thick, green!50!black] (11) -- (111);
\draw[->, very thick, green!50!black] (22) -- (333);
\draw[->, very thick, green!50!black] (33) -- (222);
\draw[->, very thick, green!50!black] (1) -- (22);
\draw[->, very thick, green!50!black] (2) -- (11);
\draw[->, very thick, green!50!black] (3) -- (33);
\draw[->, very thick, green!50!black] (4) -- (22);
\node(oname) at (1.5, -0.5) {$\color{green!50!black}{g}$};
\node(pname) at (4.5, -0.5) {$\color{green!50!black}{f}$};
\end{tikzpicture}%
\ \ .
\]
This allows us to visually construct $f\circ g$ by removing the middle blob
and merging each $g$-arrow with the $f$-arrow that starts where the $g$-arrow
ends:
\[%
\begin{tikzpicture}
\filldraw[fill=cyan!20!white] (0, 0.5) ellipse (1 and 2);
\node(1) at (0, 2) {$1$};
\node(2) at (0, 1) {$2$};
\node(3) at (0, 0) {$3$};
\node(4) at (0, -1) {$4$};
\node(Cname) at (0, 3) {$\set{1,2,3,4}$};
\filldraw[fill=cyan!20!white] (6, 0.5) ellipse (1 and 2);
\node(111) at (6, 2) {$1$};
\node(222) at (6, 1) {$2$};
\node(333) at (6, 0) {$3$};
\node(444) at (6, -1) {$4$};
\node(Ename) at (6, 3) {$\set{1,2,3,4}$};
\draw[->, very thick, green!50!black] (1) -- (333);
\draw[->, very thick, green!50!black] (2) -- (111);
\draw[->, very thick, green!50!black] (3) -- (222);
\draw[->, very thick, green!50!black] (4) -- (333);
\node(oname) at (3, -1) {$\color{green!50!black}{f \circ g}$};
\end{tikzpicture}%
\]

\end{example}

\begin{exercise}
\label{exe.mt2.6}For any positive integer $d$, let us define the function%
\begin{align*}
r_{d}:\mathbb{Z}  &  \rightarrow\mathbb{Z},\\
n  &  \mapsto n\%d
\end{align*}
(which sends each integer $n$ to the remainder of the division of $n$ by $d$).
For example, $r_{5}\left(  18\right)  =18\%5=3$ and $r_{6}\left(  18\right)
=18\%6=0$. \medskip

\textbf{(a)} Make a table of the values of the function $r_{2}\circ r_{3}$ on
the inputs $0,1,2,3,4,5$. \medskip

\textbf{(b)} Prove that $r_{2}\circ r_{3}\neq r_{2}$. \medskip

\textbf{(c)} Let $d$ and $e$ be two positive integers such that $d\mid e$.
Prove that $r_{d}\circ r_{e}=r_{d}$.
\end{exercise}

\subsubsection{\label{subsec.enum.maps.fog-props}Basic properties}

Let us recap. In Definition \ref{def.function.fog}, we defined the
\textbf{composition} of two functions\footnote{Recall: \textquotedblleft
Function\textquotedblright\ and \textquotedblleft map\textquotedblright\ mean
the same thing.} $f$ and $g$ to be the function%
\begin{align*}
\left(  \text{domain of }g\right)   &  \rightarrow\left(  \text{target of
}f\right)  ,\\
x  &  \mapsto f\left(  g\left(  x\right)  \right)  .
\end{align*}
This composition is denoted by $f\circ g$.

As we saw, the compositions $f\circ g$ and $g\circ f$ are usually not the same
(in fact, in many cases, one of these is defined and the other isn't). In
other words, composition of functions does not satisfy commutativity. However,
it has a few other nice properties:

\begin{theorem}
[associativity of composition]\label{thm.maps.compos.ass}Let $X,Y,Z,W$ be four
sets. Let $f:Z\rightarrow W$ and $g:Y\rightarrow Z$ and $h:X\rightarrow Y$ be
three functions. Then,%
\[
\left(  f\circ g\right)  \circ h=f\circ\left(  g\circ h\right)  .
\]

\end{theorem}

\begin{proof}
Both $\left(  f\circ g\right)  \circ h$ and $f\circ\left(  g\circ h\right)  $
are functions from $X$ to $W$. Moreover, for each $x\in X$, we have%
\begin{align*}
\left(  f\circ\left(  g\circ h\right)  \right)  \left(  x\right)   &
=f\left(  \left(  g\circ h\right)  \left(  x\right)  \right)
\ \ \ \ \ \ \ \ \ \ \left(  \text{by the definition of }f\circ\left(  g\circ
h\right)  \right) \\
&  =f\left(  g\left(  h\left(  x\right)  \right)  \right)
\ \ \ \ \ \ \ \ \ \ \left(
\begin{array}
[c]{c}%
\text{since the definition of }g\circ h\\
\text{yields }\left(  g\circ h\right)  \left(  x\right)  =g\left(  h\left(
x\right)  \right)
\end{array}
\right)
\end{align*}
and%
\begin{align*}
\left(  \left(  f\circ g\right)  \circ h\right)  \left(  x\right)   &
=\left(  f\circ g\right)  \left(  h\left(  x\right)  \right)
\ \ \ \ \ \ \ \ \ \ \left(  \text{by the definition of }\left(  f\circ
g\right)  \circ h\right) \\
&  =f\left(  g\left(  h\left(  x\right)  \right)  \right)
\ \ \ \ \ \ \ \ \ \ \left(  \text{by the definition of }f\circ g\right)  ,
\end{align*}
so that%
\[
\left(  f\circ\left(  g\circ h\right)  \right)  \left(  x\right)  =f\left(
g\left(  h\left(  x\right)  \right)  \right)  =\left(  \left(  f\circ
g\right)  \circ h\right)  \left(  x\right)  .
\]
Since this holds for each $x\in X$, we conclude that $f\circ\left(  g\circ
h\right)  =\left(  f\circ g\right)  \circ h$ (because two functions $u$ and
$v$ from $X$ to $W$ are equal if and only if the equality $u\left(  x\right)
=v\left(  x\right)  $ holds for each $x\in X$). This proves the theorem.
\end{proof}

Intuitively, the claim of Theorem \ref{thm.maps.compos.ass} is pretty obvious:
It is just saying that if you can do three things (applying $h$, applying $g$
and applying $f$) in succession, then it does not matter whether you view it
as \textquotedblleft first doing $h$ followed by $g$, and then doing
$f$\textquotedblright\ or as \textquotedblleft first doing $h$, and then doing
$g$ followed by $f$\textquotedblright.

Thanks to Theorem \ref{thm.maps.compos.ass}, we can write compositions of
several functions without parentheses: i.e., instead of writing $f\circ\left(
g\circ h\right)  $ or $\left(  f\circ g\right)  \circ h$, we can just write
$f\circ g\circ h$. \medskip

The following property of composition of functions is even easier. We recall
that $\operatorname*{id}\nolimits_{P}$ means the identity map on a given set
$P$; this is the map from $P$ to $P$ that sends each element $p\in P$ to itself.

\begin{theorem}
\label{thm.maps.compos.id}Let $f:X\rightarrow Y$ be a function. Then,%
\[
f\circ\operatorname*{id}\nolimits_{X}=\operatorname*{id}\nolimits_{Y}\circ
f=f.
\]

\end{theorem}

\begin{proof}
For each $x\in X$, we have%
\begin{align*}
\left(  f\circ\operatorname*{id}\nolimits_{X}\right)  \left(  x\right)   &
=f\left(  \operatorname*{id}\nolimits_{X}\left(  x\right)  \right) \\
&  =f\left(  x\right)  \ \ \ \ \ \ \ \ \ \ \left(  \text{since the definition
of }\operatorname*{id}\nolimits_{X}\text{ yields }\operatorname*{id}%
\nolimits_{X}\left(  x\right)  =x\right)  .
\end{align*}
This shows that $f\circ\operatorname*{id}\nolimits_{X}=f$ (since both
$f\circ\operatorname*{id}\nolimits_{X}$ and $f$ are functions from $X$ to
$Y$). A similar computation yields $\operatorname*{id}\nolimits_{Y}\circ f=f$.
Thus, the theorem follows.
\end{proof}

Thanks to Theorem \ref{thm.maps.compos.id}, we can remove identity maps from
compositions: e.g., the composition $f\circ g\circ\operatorname*{id}%
\nolimits_{P}\circ\,h$ (where $P$ is the target of $h$ and the domain of $g$)
can be simplified to $f\circ g\circ h$.

\begin{exercise}
Let $s_{1},\ s_{2},\ s_{3}$ be the three functions from $\left\{
1,2,3,4\right\}  $ to $\left\{  1,2,3,4\right\}  $ defined by the following
tables of values:%
\[%
\begin{tabular}
[c]{|c||c|c|c|c|}\hline
$i$ & $1$ & $2$ & $3$ & $4$\\\hline
$s_{1}\left(  i\right)  $ & $2$ & $1$ & $3$ & $4$\\\hline
$s_{2}\left(  i\right)  $ & $1$ & $3$ & $2$ & $4$\\\hline
$s_{3}\left(  i\right)  $ & $1$ & $2$ & $4$ & $3$\\\hline
\end{tabular}
\ \ .
\]
(That is, each $s_{i}$ is the function that transforms the two numbers $i$ and
$i+1$ into one another while leaving all other inputs unchanged.) \medskip

\textbf{(a)} Make a table of values of the function $s_{2}\circ s_{3}\circ
s_{1}\circ s_{2}$. \medskip

\textbf{(b)} Is $s_{1}\circ s_{3}=s_{3}\circ s_{1}$ ?\medskip

\textbf{(c)} Is $s_{1}\circ s_{2}=s_{2}\circ s_{1}$ ? \medskip

\textbf{(d)} Let $w$ be the function from $\left\{  1,2,3,4\right\}  $ to
$\left\{  1,2,3,4\right\}  $ with the following table of values:%
\[%
\begin{tabular}
[c]{|c||c|c|c|c|}\hline
$i$ & $1$ & $2$ & $3$ & $4$\\\hline
$w\left(  i\right)  $ & $4$ & $2$ & $1$ & $3$\\\hline
\end{tabular}
\ \ .
\]
Write $w$ as a composition of some of the functions $s_{1},s_{2},s_{3}$. (You
can use these functions in any order and any number of times, including none.
For example, \textquotedblleft$s_{2}\circ s_{3}\circ s_{1}\circ s_{2}%
$\textquotedblright\ would be a valid answer if this function was $w$.)
\end{exercise}

\subsection{\label{sec.enum.maps.jective}Jectivities (injectivity,
surjectivity and bijectivity)}

Now we introduce some important properties of functions, which have to do with
how often they attain certain values. There are three of these properties, and
I refer to them as the \textquotedblleft jectivity
properties\textquotedblright, as they are called injectivity, surjectivity and bijectivity.

\begin{definition}
\label{def.maps.jects.def}Let $f:X\rightarrow Y$ be a function. Then: \medskip

\textbf{(a)} We say that $f$ is \textbf{injective} (aka \textbf{one-to-one},
aka an \textbf{injection}) if
\[
\text{for each }y\in Y\text{, there exists \textbf{at most one} }x\in X\text{
such that }f\left(  x\right)  =y.
\]
In other words: We say that $f$ is \textbf{injective} if there are no two
distinct elements $x_{1},x_{2}\in X$ such that $f\left(  x_{1}\right)
=f\left(  x_{2}\right)  $.

In other words: We say that $f$ is \textbf{injective} if any two elements
$x_{1},x_{2}\in X$ satisfying $f\left(  x_{1}\right)  =f\left(  x_{2}\right)
$ must also satisfy $x_{1}=x_{2}$. \medskip

\textbf{(b)} We say that $f$ is \textbf{surjective} (aka \textbf{onto}, aka a
\textbf{surjection}) if
\[
\text{for each }y\in Y\text{, there exists \textbf{at least one} }x\in X\text{
such that }f\left(  x\right)  =y.
\]
In other words: We say that $f$ is \textbf{surjective} if every element of $Y$
is an output value of $f$. \medskip

\textbf{(c)} We say that $f$ is \textbf{bijective} (aka a \textbf{one-to-one
correspondence}, aka a \textbf{bijection}) if
\[
\text{for each }y\in Y\text{, there exists \textbf{exactly one} }x\in X\text{
such that }f\left(  x\right)  =y.
\]
Thus, $f$ is bijective if and only if $f$ is both injective and surjective.
\end{definition}

Here are some examples:

\begin{itemize}
\item The function%
\begin{align*}
f:\mathbb{N}  &  \rightarrow\mathbb{N},\\
k  &  \mapsto k^{2}%
\end{align*}
is injective (because no two distinct nonnegative integers $x_{1},x_{2}$
satisfy $x_{1}^{2}=x_{2}^{2}$) but not surjective (because, e.g., the
nonnegative integer $2\in\mathbb{N}$ is not the square of any nonnegative
integer). Thus, it is not bijective.

\item Let $S=\left\{  0,1,4,9,16,\ldots\right\}  $ be the set of all perfect
squares (i.e., all squares of integers). Then, the function%
\begin{align*}
g:\mathbb{N}  &  \rightarrow S,\\
k  &  \mapsto k^{2}%
\end{align*}
is injective (for the same reason as the $f$ in the previous example) and also
surjective (since every perfect square can be written as $k^{2}$ for some
$k\in\mathbb{N}$). Thus, it is bijective.

Take note: The functions $f$ and $g$ differ only in their choice of target!
Other than that, they are indistinguishable (both have domain $\mathbb{N}$,
and send each element of this domain to its square). But of course, this
little difference matters for the surjectivity, since the surjectivity depends
crucially on the target. No wonder that $g$ is surjective while $f$ is not.

\item Let $S=\left\{  0,1,4,9,16,\ldots\right\}  $ be the set of all perfect
squares again. Consider the function%
\begin{align*}
g_{\mathbb{Z}}:\mathbb{Z}  &  \rightarrow S,\\
k  &  \mapsto k^{2},
\end{align*}
which differs from $g$ only in its domain (it allows all integers rather than
only nonnegative integers as inputs). This function $g_{\mathbb{Z}}$ is not
injective (since $g_{\mathbb{Z}}\left(  1\right)  =g_{\mathbb{Z}}\left(
-1\right)  $), but is still surjective (since each perfect square can be
written as $k^{2}$ for some $k\in\mathbb{Z}$). Since it is not injective, it
cannot be bijective.

\item The function
\begin{align*}
h:\mathbb{N}  &  \rightarrow\mathbb{N},\\
k  &  \mapsto k//2
\end{align*}
(recall that $k//2$ is the quotient of the division of $k$ by $2$) is not
injective (for example, the two distinct elements $0,1\in\mathbb{N}$ satisfy
$h\left(  0\right)  =h\left(  1\right)  $, because both $h\left(  0\right)
=0//2$ and $h\left(  1\right)  =1//2$ are $0$), but is surjective (because for
each $y\in\mathbb{N}$, there exists an $x\in\mathbb{N}$ such that $h\left(
x\right)  =y$, namely for example $x=2y$). Hence, it is not bijective.

\item Let $E=\left\{  0,2,4,6,\ldots\right\}  $ be the set of all even
nonnegative integers. The function%
\begin{align*}
h_{\operatorname*{even}}:E  &  \rightarrow\mathbb{N},\\
k  &  \mapsto k//2
\end{align*}
(note that $k//2=k/2$ here, since $k$ is even) is both injective and
surjective, thus bijective.

\item Let $O=\left\{  1,3,5,7,\ldots\right\}  $ be the set of all odd
nonnegative integers. The function%
\begin{align*}
h_{\operatorname*{odd}}:O  &  \rightarrow\mathbb{N},\\
k  &  \mapsto k//2
\end{align*}
is also injective and surjective, thus bijective.

\item The function%
\begin{align*}
f:\mathbb{Z}\times\mathbb{Z}  &  \rightarrow\mathbb{Z},\\
\left(  a,b\right)   &  \mapsto a+b
\end{align*}
(that is, the addition of integers) is not injective (because, for instance,
$f\left(  0,1\right)  =1=f\left(  1,0\right)  $), but is surjective (since
each $n\in\mathbb{Z}$ satisfies $n=f\left(  n,0\right)  $). In other words,
two pairs of integers can have the same sum, but every integer can be written
as a sum of two integers.
\end{itemize}

The following criterion for injectivity, surjectivity and bijectivity is just
a restatement of Definition \ref{def.maps.jects.def}, but it can be quite
useful for checking these properties:

\begin{remark}
\label{rmk.maps.jects.table}Consider a function $f:X\rightarrow Y$ given by a
table of all its values (possibly an infinite table if $X$ is infinite).
Assume that all possible inputs $x\in X$ appear in the top row (each exactly
once), and the corresponding outputs $f\left(  x\right)  $ appear in the
bottom row, so the table looks as follows:%
\[%
\begin{tabular}
[c]{|c||c|c|c|c|c}\hline
$x$ & $a$ & $b$ & $c$ & $d$ & $\cdots$\\\hline
$f\left(  x\right)  $ & $f\left(  a\right)  $ & $f\left(  b\right)  $ &
$f\left(  c\right)  $ & $f\left(  d\right)  $ & $\cdots$\\\hline
\end{tabular}
\ \ \ \ \ \ .
\]
Then: \medskip

\textbf{(a)} The function $f$ is injective if and only if the bottom row of
this table has no two equal entries. \medskip

\textbf{(b)} The function $f$ is surjective if and only if every element of
$Y$ appears in the bottom row. \medskip

\textbf{(c)} The function $f$ is bijective if and only if every element of $Y$
appears exactly once in the bottom row.
\end{remark}

For example:

\begin{itemize}
\item The function%
\begin{align*}
f:\left\{  1,2,3\right\}   &  \rightarrow\left\{  7,8,9\right\}  ,\\
k  &  \mapsto k+6
\end{align*}
is bijective, as you can see from its table of values:%
\[%
\begin{tabular}
[c]{|c||c|c|c|}\hline
$k$ & $1$ & $2$ & $3$\\\hline
$f\left(  k\right)  $ & $7$ & $8$ & $9$\\\hline
\end{tabular}
\]
(by noticing that every element of $\left\{  7,8,9\right\}  $ appears exactly
once in the bottom row of this table). Of course, this can also be shown
logically (by arguing that $f$ is injective and surjective because adding $6$
can be undone by subtracting $6$).

\item The function\footnote{Recall that $k\%3$ denotes the remainder of the
division of $k$ by $3$.}%
\begin{align*}
f:\left\{  4,6,7\right\}   &  \rightarrow\left\{  0,1,2\right\}  ,\\
k  &  \mapsto k\%3
\end{align*}
is neither injective nor surjective. Indeed, its table of values%
\[%
\begin{tabular}
[c]{|c||c|c|c|}\hline
$k$ & $4$ & $6$ & $7$\\\hline
$f\left(  k\right)  $ & $1$ & $0$ & $1$\\\hline
\end{tabular}
\ \
\]
has the element $1$ appear twice in the bottom row (so $f$ is not injective)
and does not have the element $2$ in its bottom row (so $f$ is not surjective).
\end{itemize}

Here is yet another way to restate Definition \ref{def.maps.jects.def}:

\begin{remark}
If you visualize a function $f:X\rightarrow Y$ as a blobs-and-arrows picture
(as explained in Section \ref{sec.enum.maps.formal}), then

\begin{itemize}
\item the function $f$ is injective if and only if no two arrows hit the same
$Y$-node;

\item the function $f$ is surjective if and only if every node in the $Y$-blob
gets hit by at least one arrow;

\item the function $f$ is bijective if and only if every node in the $Y$-blob
gets hit by exactly one arrow.
\end{itemize}
\end{remark}

This can be illustrated by the following four examples:%
\[%
\begin{tabular}
[c]{|c|c|}\hline
$%
\begin{array}
[c]{c}%
\begin{tikzpicture}
\filldraw[fill=cyan!20!white] (0, 1) ellipse (1 and 1.5);
\node(1) at (0, 2) {$1$};
\node(2) at (0, 1) {$2$};
\node(3) at (0, 0) {$3$};
\node(Cname) at (0, 3) {$X$};
\filldraw[fill=cyan!20!white] (3, 1) ellipse (1 and 1.5);
\node(5) at (3, 2) {$1$};
\node(6) at (3, 1) {$2$};
\node(7) at (3, 0) {$3$};
\node(Dname) at (3, 3) {$Y$};
\draw[->, very thick, green!50!black] (1) -- (6);
\draw[->, very thick, green!50!black] (2) -- (7);
\draw[->, very thick, green!50!black] (3) -- (6);
\end{tikzpicture}%
\\
\text{neither injective nor surjective}\\
\text{(since }1\in Y\text{ is not hit,}\\
\text{while }2\in Y\text{ is hit twice)}%
\end{array}
$ & $%
\begin{array}
[c]{c}%
\begin{tikzpicture}
\filldraw[fill=cyan!20!white] (0, 1.5) ellipse (1 and 1);
\node(1) at (0, 2) {$1$};
\node(2) at (0, 1) {$2$};
\node(Cname) at (0, 3) {$X$};
\filldraw[fill=cyan!20!white] (3, 1) ellipse (1 and 1.5);
\node(5) at (3, 2) {$1$};
\node(6) at (3, 1) {$2$};
\node(7) at (3, 0) {$3$};
\node(Dname) at (3, 3) {$Y$};
\draw[->, very thick, green!50!black] (1) -- (5);
\draw[->, very thick, green!50!black] (2) -- (7);
\end{tikzpicture}%
\\
\text{injective but not surjective}\\
\text{(since }2\in Y\text{ is not hit)}%
\end{array}
$\\\hline
$%
\begin{array}
[c]{c}%
\begin{tikzpicture}
\filldraw[fill=cyan!20!white] (0, 1) ellipse (1 and 1.5);
\node(1) at (0, 2) {$1$};
\node(2) at (0, 1) {$2$};
\node(3) at (0, 0) {$3$};
\node(Cname) at (0, 3) {$X$};
\filldraw[fill=cyan!20!white] (3, 1.5) ellipse (1 and 1);
\node(5) at (3, 2) {$1$};
\node(6) at (3, 1) {$2$};
\node(Dname) at (3, 3) {$Y$};
\draw[->, very thick, green!50!black] (1) -- (5);
\draw[->, very thick, green!50!black] (2) -- (6);
\draw[->, very thick, green!50!black] (3) -- (5);
\end{tikzpicture}%
\\
\text{surjective but not injective}\\
\text{(since }1\in Y\text{ is hit twice)}%
\end{array}
$ & $%
\begin{array}
[c]{c}%
\begin{tikzpicture}
\filldraw[fill=cyan!20!white] (0, 1) ellipse (1 and 1.5);
\node(1) at (0, 2) {$1$};
\node(2) at (0, 1) {$2$};
\node(3) at (0, 0) {$3$};
\node(Cname) at (0, 3) {$X$};
\filldraw[fill=cyan!20!white] (3, 1) ellipse (1 and 1.5);
\node(5) at (3, 2) {$1$};
\node(6) at (3, 1) {$2$};
\node(7) at (3, 0) {$3$};
\node(Dname) at (3, 3) {$Y$};
\draw[->, very thick, green!50!black] (1) -- (6);
\draw[->, very thick, green!50!black] (2) -- (7);
\draw[->, very thick, green!50!black] (3) -- (5);
\end{tikzpicture}%
\\
\text{both injective and surjective}\\
\text{(i.e., bijective)}%
\end{array}
$\\\hline
\end{tabular}
\
\]

\bigskip

\begin{exercise}
For each of the following functions, determine whether it is injective,
surjective and/or bijective: \medskip

\textbf{(a)} The function%
\begin{align*}
f:\mathbb{Z}  &  \rightarrow\mathbb{Z},\\
x  &  \mapsto x^{2}.
\end{align*}

\textbf{(b)} The function%
\begin{align*}
f:\mathbb{Z}  &  \rightarrow\mathbb{Z},\\
x  &  \mapsto x^{3}.
\end{align*}

\textbf{(c)} The function%
\begin{align*}
f:\mathbb{Z}\times\mathbb{Z}  &  \rightarrow\mathbb{Z},\\
\left(  x,y\right)   &  \mapsto x^{2}+y^{2}.
\end{align*}

\textbf{(d)} The function%
\begin{align*}
f:\mathbb{Z}  &  \rightarrow\mathbb{Z},\\
x  &  \mapsto3-x.
\end{align*}

\textbf{(e)} The function%
\begin{align*}
f:\mathbb{Z}  &  \rightarrow\mathbb{Z},\\
x  &  \mapsto3-2x.
\end{align*}

\textbf{(f)} The function%
\begin{align*}
f:\mathbb{N}  &  \rightarrow\mathbb{N},\\
x  &  \mapsto x!.
\end{align*}
(Keep in mind that $0\in\mathbb{N}$.) \medskip

\textbf{(g)} The function%
\begin{align*}
f:\mathbb{Z}\times\mathbb{Z}  &  \rightarrow\mathbb{Z}\times\mathbb{Z},\\
\left(  x,y\right)   &  \mapsto\left(  x+y,\ x-y\right)  .
\end{align*}

\textbf{(h)} The function%
\begin{align*}
f:\mathbb{Z}\times\mathbb{Z}  &  \rightarrow\mathbb{Z}\times\mathbb{Z},\\
\left(  x,y\right)   &  \mapsto\left(  x-y,\ y-x\right)  .
\end{align*}

\textbf{(i)} The function
\begin{align*}
f:\mathbb{Z}\times\mathbb{Z}  &  \rightarrow\mathbb{Z}\times\mathbb{Z},\\
\left(  x,y\right)   &  \mapsto\left(  x+2y,\ x+y\right)  .
\end{align*}

\textbf{(j)} The function%
\begin{align*}
f:\mathbb{Z}\times\mathbb{Z}  &  \rightarrow\mathbb{Z}\times\mathbb{Z},\\
\left(  x,y\right)   &  \mapsto\left(  x+2y,\ 2x+y\right)  .
\end{align*}

\end{exercise}

\subsection{\label{sec.enum.maps.inv}Inverses}

\subsubsection{Definition and examples}

Bijective maps have a special power: They can be \textbf{inverted}. Here is
what this means:

\begin{definition}
Let $f:X\rightarrow Y$ be a function. An \textbf{inverse} of $f$ means a
function $g:Y\rightarrow X$ such that%
\[
f\circ g=\operatorname*{id}\nolimits_{Y}\ \ \ \ \ \ \ \ \ \ \text{and}%
\ \ \ \ \ \ \ \ \ \ g\circ f=\operatorname*{id}\nolimits_{X}.
\]

In other words, an \textbf{inverse} of $f$ means a function $g:Y\rightarrow X$
such that%
\begin{align*}
f\left(  g\left(  y\right)  \right)   &  =y\ \ \ \ \ \ \ \ \ \ \text{for each
}y\in Y,\ \ \ \ \ \ \ \ \ \ \text{and}\\
g\left(  f\left(  x\right)  \right)   &  =x\ \ \ \ \ \ \ \ \ \ \text{for each
}x\in X.
\end{align*}

\end{definition}

Roughly speaking, an inverse of $f$ thus means a map that both undoes $f$ and
is undone by $f$.

Not every function has an inverse. We shall soon see which ones do and which
ones don't; we will also prove that an inverse of $f$ is unique if it exists.
For now, however, let us explore a few examples:

\begin{itemize}
\item Let $f:\left\{  1,2,3\right\}  \rightarrow\left\{  7,8,9\right\}  $ be
the \textquotedblleft add $6$\textquotedblright\ function -- i.e., the
function that sends each $x\in\left\{  1,2,3\right\}  $ to $x+6\in\left\{
7,8,9\right\}  $. Then, $f$ has an inverse: the \textquotedblleft subtract
$6$\textquotedblright\ function (i.e., the function from $\left\{
7,8,9\right\}  $ to $\left\{  1,2,3\right\}  $ that sends each $y$ to $y-6$).
Indeed, if we denote the \textquotedblleft subtract $6$\textquotedblright%
\ function by $g$, then we have%
\begin{align*}
f\left(  g\left(  y\right)  \right)   &  =f\left(  y-6\right)  =\left(
y-6\right)  +6=y\ \ \ \ \ \ \ \ \ \ \text{for each }y\in\left\{
7,8,9\right\}  ,\ \ \ \ \ \ \ \ \ \ \text{and}\\
g\left(  f\left(  x\right)  \right)   &  =g\left(  x+6\right)  =\left(
x+6\right)  -6=x\ \ \ \ \ \ \ \ \ \ \text{for each }x\in\left\{
1,2,3\right\}  .
\end{align*}

\item Let $f:\left\{  1,2,3\right\}  \rightarrow\left\{  7,8,9\right\}  $ be
the \textquotedblleft subtract from $10$\textquotedblright\ function -- i.e.,
the function that sends each $x\in\left\{  1,2,3\right\}  $ to $10-x\in
\left\{  7,8,9\right\}  $. Then, $f$ has an inverse: the \textquotedblleft
subtract from $10$\textquotedblright\ function $g$ from $\left\{
7,8,9\right\}  $ to $\left\{  1,2,3\right\}  $ (which sends each $y$ to
$10-y$). This is because $10-\left(  10-n\right)  =n$ for each $n\in
\mathbb{Z}$. Note that the functions $f$ and $g$ are given by the same formula
($k\mapsto10-k$), but they are not the same function, since they have
different domains (and targets).

\item Let $f:\left\{  1,2,3,4,5\right\}  \rightarrow\left\{
1,2,3,4,5\right\}  $ be the function that sends $1,2,3,4,5$ to $3,4,1,5,2$,
respectively. Then, $f$ has an inverse: namely, the function $g$ that sends
$1,2,3,4,5$ to $3,5,1,2,4$, respectively. We can check that $f\left(  g\left(
y\right)  \right)  =y$ for each $y\in\left\{  1,2,3,4,5\right\}  $. For
example, for $y=3$, this is because $f\left(  g\left(  3\right)  \right)
=f\left(  1\right)  =3$. Similarly we can check that $g\left(  f\left(
x\right)  \right)  =x$ for each $x\in\left\{  1,2,3,4,5\right\}  $.

This is best seen by drawing the blobs-and-arrows diagrams of $f$ and $g$ side
by side:%
\[%
\begin{tabular}
[c]{|c|c|}\hline
$%
\begin{array}
[c]{c}%
\begin{tikzpicture}
\filldraw[fill=cyan!20!white] (0, 0) ellipse (1 and 2.5);
\node(A1) at (0, 2) {$1$};
\node(A2) at (0, 1) {$2$};
\node(A3) at (0, 0) {$3$};
\node(A4) at (0, -1) {$4$};
\node(A5) at (0, -2) {$5$};
\node(Aname) at (0, 3) {$X$};
\filldraw[fill=cyan!20!white] (3, 0) ellipse (1 and 2.5);
\node(B1) at (3, 2) {$1$};
\node(B2) at (3, 1) {$2$};
\node(B3) at (3, 0) {$3$};
\node(B4) at (3, -1) {$4$};
\node(B5) at (3, -2) {$5$};
\node(Aname) at (3, 3) {$Y$};
\draw[->, very thick, green!50!black] (A1) -- (B3);
\draw[->, very thick, green!50!black] (A2) -- (B4);
\draw[->, very thick, green!50!black] (A3) -- (B1);
\draw[->, very thick, green!50!black] (A4) -- (B5);
\draw[->, very thick, green!50!black] (A5) -- (B2);
\end{tikzpicture}%
\\
f
\end{array}
$ & $%
\begin{array}
[c]{c}%
\begin{tikzpicture}
\filldraw[fill=cyan!20!white] (0, 0) ellipse (1 and 2.5);
\node(A1) at (0, 2) {$1$};
\node(A2) at (0, 1) {$2$};
\node(A3) at (0, 0) {$3$};
\node(A4) at (0, -1) {$4$};
\node(A5) at (0, -2) {$5$};
\node(Aname) at (0, 3) {$Y$};
\filldraw[fill=cyan!20!white] (3, 0) ellipse (1 and 2.5);
\node(B1) at (3, 2) {$1$};
\node(B2) at (3, 1) {$2$};
\node(B3) at (3, 0) {$3$};
\node(B4) at (3, -1) {$4$};
\node(B5) at (3, -2) {$5$};
\node(Aname) at (3, 3) {$X$};
\draw[->, very thick, green!50!black] (A3) -- (B1);
\draw[->, very thick, green!50!black] (A4) -- (B2);
\draw[->, very thick, green!50!black] (A1) -- (B3);
\draw[->, very thick, green!50!black] (A5) -- (B4);
\draw[->, very thick, green!50!black] (A2) -- (B5);
\end{tikzpicture}%
\\
g
\end{array}
$\\\hline
\end{tabular}
\ \ \ \ \ .
\]
As you see, there is a \textquotedblleft dual\textquotedblright\ relationship
between these two diagrams: Whenever the diagram of $f$ has an arrow from some
$x\in X$ to some $y\in Y$, the diagram of $g$ has an arrow from $y$ to $x$. In
other words, the diagram of $g$ can be obtained from the diagram of $f$ by
swapping the $X$-blob with the $Y$-blob and reversing the direction of each
arrow. This rule applies not just to our specific two maps $f$ and $g$, but to
any map $f$ that has an inverse. Thus, if you have drawn a blobs-and-arrows
diagram of a function $f$, it is fairly easy to construct its inverse (as long
as such an inverse exists).

This rule can also be restated in terms of tables of values: If you have a
table of all values of a function $f:X\rightarrow Y$, then you can get an
inverse of $f$ by swapping the two rows of this table (provided that the
result is a function). For instance, if $f:\left\{  1,2,3,4,5\right\}
\rightarrow\left\{  1,2,3,4,5\right\}  $ is the function we just showed, then
$f$ has the table of values%
\[%
\begin{tabular}
[c]{|c||c|c|c|c|c|}\hline
$k$ & $1$ & $2$ & $3$ & $4$ & $5$\\\hline
$f\left(  k\right)  $ & $3$ & $4$ & $1$ & $5$ & $2$\\\hline
\end{tabular}
\ \ \ \ \ ,
\]
and thus you can get its inverse $g$ by swapping the two rows:%
\[%
\begin{tabular}
[c]{|c||c|c|c|c|c|}\hline
$k$ & $3$ & $4$ & $1$ & $5$ & $2$\\\hline
$g\left(  k\right)  $ & $1$ & $2$ & $3$ & $4$ & $5$\\\hline
\end{tabular}
\ \ \ \ \ .
\]

\item Let $f:\left\{  1,2,3,4\right\}  \rightarrow\left\{  1,2,3,4\right\}  $
be the function that sends $1,2,3,4$ to $1,2,3,3$, respectively. Then, $f$ has
no inverse. Indeed, if $g$ was an inverse of $f$, then we would have%
\begin{align*}
3  &  =g\left(  f\left(  3\right)  \right)  \ \ \ \ \ \ \ \ \ \ \left(
\text{since }g\left(  f\left(  x\right)  \right)  =x\text{ for each }%
x\in\left\{  1,2,3,4\right\}  \right) \\
&  =g\left(  f\left(  4\right)  \right)  \ \ \ \ \ \ \ \ \ \ \left(
\text{since }f\left(  3\right)  =3=f\left(  4\right)  \right) \\
&  =4\ \ \ \ \ \ \ \ \ \ \left(  \text{since }g\left(  f\left(  x\right)
\right)  =x\text{ for each }x\in\left\{  1,2,3,4\right\}  \right)  ,
\end{align*}
which is absurd.

The same argument shows that more generally, if a function $f:X\rightarrow Y$
is to have an inverse, then $f$ should be injective, because two distinct
elements $x_{1}$ and $x_{2}$ of $X$ satisfying $f\left(  x_{1}\right)
=f\left(  x_{2}\right)  $ would create a contradiction via $x_{1}=g\left(
f\left(  x_{1}\right)  \right)  =g\left(  f\left(  x_{2}\right)  \right)
=x_{2}$.

\item Let $f:\left\{  1,2,3\right\}  \rightarrow\left\{  1,2,3,4\right\}  $ be
the function that sends $1,2,3$ to $1,2,3$, respectively. Then, $f$ has no
inverse. Indeed, if $g$ was an inverse of $f$, then we would have $f\left(
g\left(  4\right)  \right)  =4$, but this is absurd, since $4$ is not an
output of $f$.

The same argument shows that more generally, if a function $f:X\rightarrow Y$
is to have an inverse, then $f$ should be surjective, because each $y\in Y$
will satisfy $y=f\left(  g\left(  y\right)  \right)  $ and thus be an output
value of $f$.
\end{itemize}

\subsubsection{Invertibility is bijectivity by another name}

Combining the morals of the last two examples, we conclude that if a function
$f:X\rightarrow Y$ is to have an inverse, then $f$ should be both injective
and surjective, i.e., should be bijective. In other words, only bijective maps
have a chance at having inverses. This turns out to be sufficient as well: If
a map is bijective, then it has an inverse. Let us summarize this as a theorem:

\begin{theorem}
\label{thm.maps.jects.inv=bij}Let $f:X\rightarrow Y$ be a map between two sets
$X$ and $Y$. Then, $f$ has an inverse if and only if $f$ is bijective.
\end{theorem}

\begin{fineprint}
\begin{proof}
We must prove the logical equivalence%
\begin{equation}
\left(  f\text{ has an inverse}\right)  \Longleftrightarrow\left(  f\text{ is
bijective}\right)  . \label{pf.thm.maps.jects.inv=bij.goal}%
\end{equation}

Let us prove the $\Longrightarrow$ and $\Longleftarrow$ directions separately:

$\Longrightarrow:$ Assume that $f$ has an inverse. We must show that $f$ is
bijective.\footnote{We have already done this in the above examples, but we
repeat it for the sake of completeness.}

We assumed that $f$ has an inverse. Let $g$ be this inverse.

Let us show that $f$ is injective. Let $x_{1},x_{2}\in X$ satisfy $f\left(
x_{1}\right)  =f\left(  x_{2}\right)  $. We must prove that $x_{1}=x_{2}$. Set
$y=f\left(  x_{1}\right)  $; then, $y=f\left(  x_{2}\right)  $ as well (since
$f\left(  x_{1}\right)  =f\left(  x_{2}\right)  $). Since $g$ is an inverse of
$f$, we have $x_{1}=g\left(  f\left(  x_{1}\right)  \right)  =g\left(
y\right)  $ (since $f\left(  x_{1}\right)  =y$) and $x_{2}=g\left(  f\left(
x_{2}\right)  \right)  =g\left(  y\right)  $ (since $f\left(  x_{2}\right)
=y$). Thus, $x_{1}=g\left(  y\right)  =x_{2}$. This completes our proof that
$f$ is injective.

Let us show that $f$ is surjective. Let $y\in Y$. Then, $y=f\left(  g\left(
y\right)  \right)  $ (since $g$ is an inverse of $f$). Therefore, there exists
an $x\in X$ such that $y=f\left(  x\right)  $ (namely, $x=g\left(  y\right)
$). So we have proved for each $y\in Y$ that there exists an $x\in X$ such
that $y=f\left(  x\right)  $. In other words, $f$ is surjective.

So $f$ is both injective and surjective, thus bijective. This proves the
\textquotedblleft$\Longrightarrow$\textquotedblright\ direction of our
equivalence (\ref{pf.thm.maps.jects.inv=bij.goal}).

Let us now prove the \textquotedblleft$\Longleftarrow$\textquotedblright\ direction:

$\Longleftarrow:$ Assume that $f$ is bijective. We must show that $f$ has an inverse.

Since $f$ is bijective, for each $y\in Y$, there exists a \textbf{unique}
$x\in X$ such that $f\left(  x\right)  =y$. Thus, we can define a map%
\[
g:Y\rightarrow X,
\]
which sends each $y\in Y$ to this unique $x$. It is easy to see that $g$ is an
inverse of $f$. Thus, $f$ has an inverse. This proves the \textquotedblleft%
$\Longleftarrow$\textquotedblright\ direction of our equivalence
(\ref{pf.thm.maps.jects.inv=bij.goal}). Thus, the proof of
(\ref{pf.thm.maps.jects.inv=bij.goal}) is complete, i.e., Theorem
\ref{thm.maps.jects.inv=bij} is proved.
\end{proof}
\end{fineprint}

Theorem \ref{thm.maps.jects.inv=bij} says that bijective maps are the same as
invertible maps (i.e., maps that have an inverse). This is a fundamental
result that is used all over mathematics.

\subsubsection{Uniqueness of the inverse}

As we promised, let us now show that an inverse of a map $f$ is unique if it exists:

\begin{theorem}
\label{thm.maps.inverse.unique}Let $f:X\rightarrow Y$ be a function. Then, $f$
has at most one inverse.
\end{theorem}

\begin{proof}
What does \textquotedblleft at most one inverse\textquotedblright\ mean? It
means that $f$ has no two distinct inverses. In other words, it means that any
two inverses of $f$ are identical.

So let us prove this. Let $g_{1}$ and $g_{2}$ be two inverses of $f$. We must
show that $g_{1}=g_{2}$.

Since $g_{1}$ is an inverse of $f$, we have $g_{1}\circ f=\operatorname*{id}%
\nolimits_{X}$ and $f\circ g_{1}=\operatorname*{id}\nolimits_{Y}$.

Since $g_{2}$ is an inverse of $f$, we have $g_{2}\circ f=\operatorname*{id}%
\nolimits_{X}$ and $f\circ g_{2}=\operatorname*{id}\nolimits_{Y}$.

By associativity of composition (Theorem \ref{thm.maps.compos.ass}), the two
maps $\left(  g_{1}\circ f\right)  \circ g_{2}$ and $g_{1}\circ\left(  f\circ
g_{2}\right)  $ are equal. Thus, we can denote both of these maps by
$g_{1}\circ f\circ g_{2}$.

Comparing%
\begin{align*}
g_{1}\circ\underbrace{f\circ g_{2}}_{=\operatorname*{id}\nolimits_{Y}}  &
=g_{1}\circ\operatorname*{id}\nolimits_{Y}=g_{1}%
\ \ \ \ \ \ \ \ \ \ \text{with}\\
\underbrace{g_{1}\circ f}_{=\operatorname*{id}\nolimits_{X}}\circ\,g_{2}  &
=\operatorname*{id}\nolimits_{X}\circ g_{2}=g_{2},
\end{align*}
we find $g_{1}=g_{2}$, qed.
\end{proof}

\begin{definition}
Let $f:X\rightarrow Y$ be a map that has an inverse. Then, this inverse (which
is unique by Theorem \ref{thm.maps.inverse.unique}) is called $f^{-1}$.
\end{definition}

Thus, if $f:X\rightarrow Y$ is a map that has an inverse (i.e., by Theorem
\ref{thm.maps.jects.inv=bij}, a bijective map), then we have%
\[
f^{-1}\circ f=\operatorname*{id}\nolimits_{X}\ \ \ \ \ \ \ \ \ \ \text{and}%
\ \ \ \ \ \ \ \ \ \ f\circ f^{-1}=\operatorname*{id}\nolimits_{Y},
\]
that is,%
\begin{align}
f^{-1}\left(  f\left(  x\right)  \right)   &  =x\ \ \ \ \ \ \ \ \ \ \text{for
each }x\in X,\ \ \ \ \ \ \ \ \ \ \text{and}\label{eq.def.maps.inverse.f-1f}\\
f\left(  f^{-1}\left(  y\right)  \right)   &  =y\ \ \ \ \ \ \ \ \ \ \text{for
each }y\in Y. \label{eq.def.maps.inverse.ff-1}%
\end{align}
These equalities should explain why the notation $f^{-1}$ was chosen for the
inverse of $f$.

\subsubsection{More examples}

Here are some further examples of inverses:

\begin{itemize}
\item Let $E=\left\{  0,2,4,6,\ldots\right\}  $ be the set of all even
nonnegative integers. Consider the function%
\begin{align*}
f:E  &  \rightarrow\mathbb{N},\\
k  &  \mapsto k/2.
\end{align*}
Then, $f$ has an inverse. This inverse is the function%
\begin{align*}
f^{-1}:\mathbb{N}  &  \rightarrow E,\\
k  &  \mapsto2k.
\end{align*}

\item Let $\mathbb{R}_{\geq0}=\left\{  \text{all nonnegative real
numbers}\right\}  $. Then, the function
\begin{align*}
f:\mathbb{R}_{\geq0}  &  \rightarrow\mathbb{R}_{\geq0},\\
x  &  \mapsto x^{2}%
\end{align*}
has an inverse. This inverse is the function%
\begin{align*}
f^{-1}:\mathbb{R}_{\geq0}  &  \rightarrow\mathbb{R}_{\geq0},\\
x  &  \mapsto\sqrt{x}.
\end{align*}

\item In contrast, the function%
\begin{align*}
f:\mathbb{R}  &  \rightarrow\mathbb{R},\\
x  &  \mapsto x^{2}%
\end{align*}
has no inverse. In fact, this function is not injective (since $f\left(
1\right)  =f\left(  -1\right)  $) and not surjective (since $-1$ is not a
square of a real number), so it is certainly not bijective, and thus not invertible.

\item The function%
\begin{align*}
f:\mathbb{R}  &  \rightarrow\mathbb{R},\\
x  &  \mapsto x^{3}%
\end{align*}
has an inverse. This inverse is the function%
\begin{align*}
f^{-1}:\mathbb{R}  &  \rightarrow\mathbb{R},\\
x  &  \mapsto\sqrt[3]{x}.
\end{align*}

\end{itemize}

\begin{fineprint}
\begin{itemize}
\item Another example of inverses comes from cryptography: If $k$ is any
integer, then the Caesar cipher $\operatorname*{ROT}\nolimits_{k}$ (defined in
Section \ref{sec.ent.crypt}, regarded as a map from $W=\left\{  \text{words}%
\right\}  $ to $W$) has an inverse, namely $\operatorname*{ROT}\nolimits_{-k}%
$. This is just saying that any word encrypted with $\operatorname*{ROT}%
\nolimits_{k}$ can be decrypted with $\operatorname*{ROT}\nolimits_{-k}$ and
vice versa.
\end{itemize}
\end{fineprint}

\begin{exercise}
Show that the function%
\begin{align*}
f:\mathbb{Z}\times\mathbb{Z}  &  \rightarrow\mathbb{Z}\times\mathbb{Z},\\
\left(  x,y\right)   &  \mapsto\left(  x+3y,\ 2x+5y\right)
\end{align*}
has an inverse. Give an explicit formula for this inverse (i.e., for
$f^{-1}\left(  \left(  u,v\right)  \right)  $). \medskip

[\textbf{Hint:} This is a linear algebra question, since $f^{-1}\left(
\left(  u,v\right)  \right)  $ should be a pair $\left(  x,y\right)
\in\mathbb{Z}\times\mathbb{Z}$ satisfying $\left(  x+3y,\ 2x+5y\right)
=\left(  u,v\right)  $.]
\end{exercise}

\subsubsection{Inverses of inverses and compositions}

Here are some more general properties of inverses:

\begin{proposition}
Let $X$ be any set. Then, the identity map $\operatorname*{id}\nolimits_{X}%
:X\rightarrow X$ is bijective, and is its own inverse.
\end{proposition}

\begin{proof}
The map $\operatorname*{id}\nolimits_{X}$ is an inverse of itself (since
$\operatorname*{id}\nolimits_{X}\circ\operatorname*{id}\nolimits_{X}%
=\operatorname*{id}\nolimits_{X}$ and $\operatorname*{id}\nolimits_{X}%
\circ\operatorname*{id}\nolimits_{X}=\operatorname*{id}\nolimits_{X}$). Hence,
it has an inverse, and thus is bijective (by Theorem
\ref{thm.maps.jects.inv=bij}).
\end{proof}

\begin{theorem}
\label{thm.maps.inverse.invinv}Let $f:X\rightarrow Y$ be a map that has an
inverse $f^{-1}:Y\rightarrow X$. Then, $f^{-1}$ has an inverse, namely $f$.
\end{theorem}

\begin{proof}
Since $f^{-1}$ is an inverse of $f$, we have $f\circ f^{-1}=\operatorname*{id}%
\nolimits_{Y}$ and $f^{-1}\circ f=\operatorname*{id}\nolimits_{X}$. But the
same two equalities can be read as saying that $f$ is an inverse of $f^{-1}$.
\end{proof}

\begin{theorem}
[socks-and-shoes formula]\label{thm.maps.inverse.socks}Let $X$, $Y$ and $Z$ be
three sets. Let $g:X\rightarrow Y$ and $f:Y\rightarrow Z$ be two bijective
functions. Then, the composition $f\circ g:X\rightarrow Z$ is bijective as
well, and its inverse is%
\[
\left(  f\circ g\right)  ^{-1}=g^{-1}\circ f^{-1}.
\]

\end{theorem}

\begin{proof}
This is obvious from the blobs-and-arrows picture; but let us check this rigorously.

For any $x\in X$, we have%
\[
\left(  g^{-1}\circ f^{-1}\right)  \left(  \left(  f\circ g\right)  \left(
x\right)  \right)  =g^{-1}\left(  \underbrace{f^{-1}\left(  f\left(  g\left(
x\right)  \right)  \right)  }_{=g\left(  x\right)  }\right)  =g^{-1}\left(
g\left(  x\right)  \right)  =x.
\]
For any $z\in Z$, we have%
\[
\left(  f\circ g\right)  \left(  \left(  g^{-1}\circ f^{-1}\right)  \left(
z\right)  \right)  =f\left(  \underbrace{g\left(  g^{-1}\left(  f^{-1}\left(
z\right)  \right)  \right)  }_{=f^{-1}\left(  z\right)  }\right)  =f\left(
f^{-1}\left(  z\right)  \right)  =z.
\]
Thus, $g^{-1}\circ f^{-1}$ is an inverse of $f\circ g$. Hence, $f\circ g$ has
an inverse, and thus is bijective (by Theorem \ref{thm.maps.jects.inv=bij}).
\end{proof}

\begin{remark}
Note that $g^{-1}\circ f^{-1}$ is not the same as $f^{-1}\circ g^{-1}$.
Indeed, $f^{-1}\circ g^{-1}$ might not even exist in Theorem
\ref{thm.maps.inverse.socks}.
\end{remark}

A surprising feature of the socks-and-shoes formula $\left(  f\circ g\right)
^{-1}=g^{-1}\circ f^{-1}$ is that the order in which the inverses $f^{-1}$ and
$g^{-1}$ appear on the right hand side is different from the order in which
$f$ and $g$ appear on the left hand side. However, this is completely natural:
If you want to undo two things you have done in some order, then you should
undo them in the opposite order! For example, if you have put on your socks
and then your shoes in the morning, then you need to first take off the shoes
and then the socks when you go to bed. (The formula owes its moniker to this metaphor.)

\begin{remark}
Part of Theorem \ref{thm.maps.inverse.socks} says that a composition of two
bijective functions is bijective. However, a composition $f\circ g$ of two
non-bijective functions $f$ and $g$ can sometimes also be bijective. Here is
an example:%
\[
\begin{tikzpicture}
\filldraw[fill=cyan!20!white] (0, 0) ellipse (0.5 and 0.5);
\node(1) at (0, 0) {$1$};
\node(Cname) at (0, 1) {$\set{1}$};
\filldraw[fill=cyan!20!white] (3, 0.5) ellipse (1 and 1);
\node(11) at (3, 1) {$1$};
\node(22) at (3, 0) {$2$};
\node(Dname) at (3, 2) {$\set{1,2}$};
\filldraw[fill=cyan!20!white] (6, 0) ellipse (0.5 and 0.5);
\node(111) at (6, 0) {$1$};
\node(Ename) at (6, 1) {$\set{1}$};
\draw[->, very thick, green!50!black] (11) -- (111);
\draw[->, very thick, green!50!black] (22) -- (111);
\draw[->, very thick, green!50!black] (1) -- (11);
\node(oname) at (1.5, -0.5) {$\color{green!50!black}{g}$};
\node(pname) at (4.5, -0.5) {$\color{green!50!black}{f}$};
\end{tikzpicture}
\]

\end{remark}

\subsection{\label{sec.enum.maps.exen}Some exercises on jectivities and
inverses}

\subsubsection{Exercises with solutions}

Here are a few solved exercises on jectivity properties and inverses.

\begin{exercise}
\label{exe.maps.jects.compos}Let $X$, $Y$ and $Z$ be three sets, and
$f:Y\rightarrow Z$ and $g:X\rightarrow Y$ be two maps. Which of the following
are true? \medskip

\textbf{(a)} If $f$ and $g$ are injective, then $f\circ g$ is injective.
\medskip

\textbf{(b)} If $f\circ g$ is injective, then $f$ is injective. \medskip

\textbf{(c)} If $f\circ g$ is injective, then $g$ is injective. \medskip

\textbf{(d)} If $f$ and $g$ are surjective, then $f\circ g$ is surjective.
\medskip

\textbf{(e)} If $f\circ g$ is surjective, then $f$ is surjective. \medskip

\textbf{(f)} If $f\circ g$ is surjective, then $g$ is surjective. \medskip

\textbf{(g)} If $f$ and $g$ are bijective, then $f\circ g$ is bijective.
\medskip

\textbf{(h)} If $f\circ g$ is bijective, then $f$ is bijective. \medskip

\textbf{(i)} If $f\circ g$ is bijective, then $g$ is bijective.
\end{exercise}

\begin{proof}
[Solution.]We shall use the following definitions of \textquotedblleft
injective\textquotedblright, \textquotedblleft surjective\textquotedblright%
\ and \textquotedblleft bijective\textquotedblright\footnote{We gave several
equivalent definitions for \textquotedblleft injective\textquotedblright,
\textquotedblleft surjective\textquotedblright\ and \textquotedblleft
bijective\textquotedblright\ in Definition \ref{def.maps.jects.def}; you can
just as well use any of them instead.}:

\begin{itemize}
\item A map $h:U\rightarrow V$ is \textbf{injective} if and only if it has the
following property: For any $u_{1},u_{2}\in U$ satisfying $h\left(
u_{1}\right)  =h\left(  u_{2}\right)  $, we have $u_{1}=u_{2}$.

\item A map $h:U\rightarrow V$ is \textbf{surjective} if and only if it has
the following property: For any $v\in V$, there exists some $u\in U$ such that
$h\left(  u\right)  =v$.

\item A map $h:U\rightarrow V$ is \textbf{bijective} if and only if $h$ is
both injective and surjective.
\end{itemize}

\textbf{(a)} This is \textbf{true}.

[\textit{Proof:} Assume that $f$ and $g$ are injective. We must prove that
$f\circ g$ is injective.

Let $u_{1},u_{2}\in X$ satisfy $\left(  f\circ g\right)  \left(  u_{1}\right)
=\left(  f\circ g\right)  \left(  u_{2}\right)  $. We shall show that
$u_{1}=u_{2}$.

Indeed, we have $\left(  f\circ g\right)  \left(  u_{1}\right)  =f\left(
g\left(  u_{1}\right)  \right)  $ (by the definition of $f\circ g$), so that%
\[
f\left(  g\left(  u_{1}\right)  \right)  =\left(  f\circ g\right)  \left(
u_{1}\right)  =\left(  f\circ g\right)  \left(  u_{2}\right)  =f\left(
g\left(  u_{2}\right)  \right)
\]
(again by the definition of $f\circ g$). Since $f$ is injective, we thus
conclude that $g\left(  u_{1}\right)  =g\left(  u_{2}\right)  $%
\ \ \ \ \footnote{In some more detail:
\par
We know that $f$ is injective. In other words, for any $v_{1},v_{2}\in Y$
satisfying $f\left(  v_{1}\right)  =f\left(  v_{2}\right)  $, we have
$v_{1}=v_{2}$ (by our definition of \textquotedblleft
injective\textquotedblright). Applying this to $v_{1}=g\left(  u_{1}\right)  $
and $v_{2}=g\left(  u_{2}\right)  $, we obtain $g\left(  u_{1}\right)
=g\left(  u_{2}\right)  $ (since $f\left(  g\left(  u_{1}\right)  \right)
=f\left(  g\left(  u_{2}\right)  \right)  $).}. Since $g$ is injective, we
thus conclude that $u_{1}=u_{2}$.

Forget that we fixed $u_{1},u_{2}$. We thus have shown that for any
$u_{1},u_{2}\in X$ satisfying $\left(  f\circ g\right)  \left(  u_{1}\right)
=\left(  f\circ g\right)  \left(  u_{2}\right)  $, we have $u_{1}=u_{2}$. In
other words, the map $f\circ g$ is injective (by our definition of
\textquotedblleft injective\textquotedblright). This completes our proof.]
\medskip

\textbf{(b)} This is \textbf{false}.

[\textit{Counterexample:} For instance, we can set $X=\left\{  1\right\}  $
and $Y=\left\{  1,2\right\}  $ and $Z=\left\{  1\right\}  $, and let
$f:Y\rightarrow Z$ be the map that sends both elements of $Y$ to $1$, while
$g:X\rightarrow Y$ is the map sending $1$ to $1$. Then, $f\circ g$ is
injective (in fact, $f\circ g$ is the identity map $\operatorname*{id}%
\nolimits_{\left\{  1\right\}  }$), but $f$ is not.] \medskip

\textbf{(c)} This is \textbf{true}.

[\textit{Proof:} Assume that $f\circ g$ is injective. We must prove that $g$
is injective.

Let $u_{1},u_{2}\in X$ satisfy $g\left(  u_{1}\right)  =g\left(  u_{2}\right)
$. We shall show that $u_{1}=u_{2}$.

Indeed, we have $\left(  f\circ g\right)  \left(  u_{1}\right)  =f\left(
g\left(  u_{1}\right)  \right)  $ (by the definition of $f\circ g$) and
$\left(  f\circ g\right)  \left(  u_{2}\right)  =f\left(  g\left(
u_{2}\right)  \right)  $ (likewise). Hence,%
\[
\left(  f\circ g\right)  \left(  u_{1}\right)  =f\left(  \underbrace{g\left(
u_{1}\right)  }_{=g\left(  u_{2}\right)  }\right)  =f\left(  g\left(
u_{2}\right)  \right)  =\left(  f\circ g\right)  \left(  u_{2}\right)  .
\]
Since $f\circ g$ is injective, this entails $u_{1}=u_{2}$.

Forget that we fixed $u_{1},u_{2}$. We thus have shown that for any
$u_{1},u_{2}\in X$ satisfying $g\left(  u_{1}\right)  =g\left(  u_{2}\right)
$, we have $u_{1}=u_{2}$. In other words, the map $g$ is injective (by our
definition of \textquotedblleft injective\textquotedblright). This completes
our proof.] \medskip

\textbf{(d)} This is \textbf{true}.

[\textit{Proof:} Assume that $f$ and $g$ are surjective. We must prove that
$f\circ g$ is surjective.

Let $z\in Z$ be arbitrary. We shall show that there exists some $x\in X$ such
that $\left(  f\circ g\right)  \left(  x\right)  =z$.

Indeed, recall that $f$ is surjective. Thus, there exists some $y\in Y$ such
that $f\left(  y\right)  =z$. Consider this $y$.

Recall now that $g$ is surjective. Thus, there exists some $w\in X$ such that
$g\left(  w\right)  =y$. Consider this $w$.

We have $\left(  f\circ g\right)  \left(  w\right)  =f\left(
\underbrace{g\left(  w\right)  }_{=y}\right)  =f\left(  y\right)  =z$. Hence,
there exists some $x\in X$ such that $\left(  f\circ g\right)  \left(
x\right)  =z$ (namely, $x=w$).

Forget that we fixed $z$. We thus have shown that for any $z\in Z$, there
exists some $x\in X$ such that $\left(  f\circ g\right)  \left(  x\right)
=z$. In other words, the map $f\circ g$ is surjective (by our definition of
\textquotedblleft surjective\textquotedblright). This completes our proof.]
\medskip

\textbf{(e)} This is \textbf{true}.

[\textit{Proof:} Assume that $f\circ g$ is surjective. We must prove that $f$
is surjective.

Let $z\in Z$ be arbitrary. We shall show that there exists some $y\in Y$ such
that $f\left(  y\right)  =z$.

Indeed, recall that $f\circ g$ is surjective. Thus, there exists some $x\in X$
such that $\left(  f\circ g\right)  \left(  x\right)  =z$. Consider this $x$.

Now, $f\left(  g\left(  x\right)  \right)  =\left(  f\circ g\right)  \left(
x\right)  =z$. Hence, there exists some $y\in Y$ such that $f\left(  y\right)
=z$ (namely, $y=g\left(  x\right)  $).

Forget that we fixed $z$. We thus have shown that for any $z\in Z$, there
exists some $y\in Y$ such that $f\left(  y\right)  =z$. In other words, the
map $f$ is surjective (by our definition of \textquotedblleft
surjective\textquotedblright). This completes our proof.] \medskip

\textbf{(f)} This is \textbf{false}.

[\textit{Counterexample:} For instance, we can set $X=\left\{  1\right\}  $
and $Y=\left\{  1,2\right\}  $ and $Z=\left\{  1\right\}  $, and let
$f:Y\rightarrow Z$ be the map that sends both elements of $Y$ to $1$, while
$g:X\rightarrow Y$ is the map sending $1$ to $1$. Then, $f\circ g$ is
surjective (in fact, $f\circ g$ is the identity map $\operatorname*{id}%
\nolimits_{\left\{  1\right\}  }$), but $g$ is not.] \medskip

\textbf{(g)} This is \textbf{true}.

[\textit{Proof:} This is part of Theorem \ref{thm.maps.inverse.socks}. But let
us give a different proof as well: Assume that $f$ and $g$ are bijective.
Thus, $f$ and $g$ are both injective and surjective. Hence, $f\circ g$ is
injective (by Exercise \ref{exe.maps.jects.compos} \textbf{(a)}) and
surjective (by Exercise \ref{exe.maps.jects.compos} \textbf{(d)}). Thus,
$f\circ g$ is bijective.] \medskip

\textbf{(h)} This is \textbf{false}.

[\textit{Counterexample:} For instance, we can set $X=\left\{  1\right\}  $
and $Y=\left\{  1,2\right\}  $ and $Z=\left\{  1\right\}  $, and let
$f:Y\rightarrow Z$ be the map that sends both elements of $Y$ to $1$, while
$g:X\rightarrow Y$ is the map sending $1$ to $1$. Then, $f\circ g$ is
bijective (in fact, $f\circ g$ is the identity map $\operatorname*{id}%
\nolimits_{\left\{  1\right\}  }$), but $f$ is not.] \medskip

\textbf{(i)} This is \textbf{false}.

[\textit{Counterexample:} For instance, we can set $X=\left\{  1\right\}  $
and $Y=\left\{  1,2\right\}  $ and $Z=\left\{  1\right\}  $, and let
$f:Y\rightarrow Z$ be the map that sends both elements of $Y$ to $1$, while
$g:X\rightarrow Y$ is the map sending $1$ to $1$. Then, $f\circ g$ is
bijective (in fact, $f\circ g$ is the identity map $\operatorname*{id}%
\nolimits_{\left\{  1\right\}  }$), but $g$ is not.]
\end{proof}

\begin{exercise}
Let $f:X\rightarrow Y$ be a map that has an inverse $f^{-1}:Y\rightarrow X$.
Let $x\in X$ and $y\in Y$. Prove that we have the logical equivalence%
\[
\left(  f\left(  x\right)  =y\right)  \ \Longleftrightarrow\ \left(
f^{-1}\left(  y\right)  =x\right)  .
\]

\end{exercise}

\begin{proof}
[Solution.]We shall prove the \textquotedblleft$\Longrightarrow$%
\textquotedblright\ and \textquotedblleft$\Longleftarrow$\textquotedblright%
\ parts of this equivalence separately:

$\Longrightarrow:$ If we have $f\left(  x\right)  =y$, then%
\[
f^{-1}\left(  \underbrace{y}_{=f\left(  x\right)  }\right)  =f^{-1}\left(
f\left(  x\right)  \right)  =x\ \ \ \ \ \ \ \ \ \ \left(  \text{by
(\ref{eq.def.maps.inverse.f-1f})}\right)  .
\]
Thus, the \textquotedblleft$\Longrightarrow$\textquotedblright\ part of the
equivalence holds.

$\Longleftarrow:$ If we have $f^{-1}\left(  y\right)  =x$, then%
\[
f\left(  \underbrace{x}_{=f^{-1}\left(  y\right)  }\right)  =f\left(
f^{-1}\left(  y\right)  \right)  =y\ \ \ \ \ \ \ \ \ \ \left(  \text{by
(\ref{eq.def.maps.inverse.ff-1})}\right)  .
\]
Thus, the \textquotedblleft$\Longleftarrow$\textquotedblright\ part of the
equivalence holds.
\end{proof}

\begin{exercise}
Let $A$ and $B$ be two sets. As we know from Exercise \ref{exe.AxB-BxA}, the
two sets $A\times B$ and $B\times A$ are usually not the same. However, I
claim that there is a bijective map from $A\times B$ to $B\times A$. Prove
this (by finding one such map, and showing that it is bijective).
\end{exercise}

\begin{proof}
[Solution.]Consider the map%
\begin{align*}
f:A\times B  &  \rightarrow B\times A,\\
\left(  a,b\right)   &  \mapsto\left(  b,a\right)  .
\end{align*}
This is the map that sends each pair $\left(  a,b\right)  \in A\times B$ to
the pair $\left(  b,a\right)  \in B\times A$; in other words, it swaps the two
entries of the input pair. Likewise, consider the map%
\begin{align*}
g:B\times A  &  \rightarrow A\times B,\\
\left(  b,a\right)   &  \mapsto\left(  a,b\right)
\end{align*}
(which does the same as $f$, but does it to a pair in $B\times A$ instead of a
pair in $A\times B$). Let us show that these two maps $f$ and $g$ are mutually inverse.

Indeed, in order to show this, we must check that $f\circ g=\operatorname*{id}%
\nolimits_{B\times A}$ and $g\circ f=\operatorname*{id}\nolimits_{A\times B}$.

Let us check that $f\circ g=\operatorname*{id}\nolimits_{B\times A}$. This
means checking that $\left(  f\circ g\right)  \left(  y\right)
=\operatorname*{id}\nolimits_{B\times A}\left(  y\right)  $ for each $y\in
B\times A$. So let $y\in B\times A$ be arbitrary. Thus, $y$ is a pair $\left(
b,a\right)  $ with $b\in B$ and $a\in A$. Consider these $b$ and $a$. Hence,
$y=\left(  b,a\right)  $, so that $g\left(  y\right)  =g\left(  \left(
b,a\right)  \right)  =\left(  a,b\right)  $ (by the definition of $g$). By the
definition of $f\circ g$, we have%
\[
\left(  f\circ g\right)  \left(  y\right)  =f\left(  \underbrace{g\left(
y\right)  }_{=\left(  a,b\right)  }\right)  =f\left(  \left(  a,b\right)
\right)  =\left(  b,a\right)
\]
(by the definition of $f$). Comparing this with $\operatorname*{id}%
\nolimits_{B\times A}\left(  y\right)  =y=\left(  b,a\right)  $, we obtain
$\left(  f\circ g\right)  \left(  y\right)  =\operatorname*{id}%
\nolimits_{B\times A}\left(  y\right)  $.

Forget that we fixed $y$. We thus have shown that $\left(  f\circ g\right)
\left(  y\right)  =\operatorname*{id}\nolimits_{B\times A}\left(  y\right)  $
for each $y\in B\times A$. Thus, we have proved the equality $f\circ
g=\operatorname*{id}\nolimits_{B\times A}$. Similarly we can show the equality
$g\circ f=\operatorname*{id}\nolimits_{A\times B}$ (since the maps $f$ and $g$
are constructed in the same way, just with the roles of $A$ and $B$ switched).
These two equalities (together) show that the map $g$ is an inverse of $f$.
Hence, the map $f$ has an inverse, and thus is bijective (by Theorem
\ref{thm.maps.jects.inv=bij}). Thus, there exists a bijective map from
$A\times B$ to $B\times A$ (namely, $f$).
\end{proof}

\subsubsection{More exercises}

Here are some more exercises.

\begin{exercise}
\label{exe.23hw5.1+}Let $A,B,C,D$ be four sets. Let $f:A\rightarrow C$ and
$g:B\rightarrow D$ be two maps. Define a new map $f\ast g:A\times B\rightarrow
C\times D$ by setting%
\[
\left(  f\ast g\right)  \left(  a,b\right)  =\left(  f\left(  a\right)
,\ g\left(  b\right)  \right)  \ \ \ \ \ \ \ \ \ \ \text{for every pair
}\left(  a,b\right)  \in A\times B.
\]

Prove the following: \medskip

\textbf{(a)} If $f$ and $g$ are injective, then $f\ast g$ is injective.
\medskip

\textbf{(b)} If $f$ and $g$ are surjective, then $f\ast g$ is surjective.
\medskip

\textbf{(c)} If $f\ast g$ is injective and the sets $A$ and $B$ are nonempty,
then $f$ and $g$ are injective. \medskip

\textbf{(d)} If $f\ast g$ is surjective and the sets $C$ and $D$ are nonempty,
then $f$ and $g$ are surjective.
\end{exercise}

\begin{exercise}
\label{exe.hw5.2}Let $X$ and $Y$ be two sets. Let $f:X\rightarrow Y$ be a map.
\medskip

A \textbf{left inverse} of $f$ means a map $g:Y\rightarrow X$ that satisfies
$g\circ f=\operatorname*{id}\nolimits_{X}$ (but not necessarily $f\circ
g=\operatorname*{id}\nolimits_{Y}$). \medskip

A \textbf{right inverse} of $f$ means a map $g:Y\rightarrow X$ that satisfies
$f\circ g=\operatorname*{id}\nolimits_{Y}$ (but not necessarily $g\circ
f=\operatorname*{id}\nolimits_{X}$). \medskip

\textbf{(a)} Prove that $f$ has a right inverse if and only if $f$ is
surjective. (Readers aware of the technicalities of logic will recognize that
this proof requires the axiom of choice.) \medskip

\textbf{(b)} Assume that $X\neq\varnothing$. Prove that $f$ has a left inverse
if and only if $f$ is injective. \medskip

\textbf{(c)} Find two distinct left inverses of the map%
\[%
\begin{tikzpicture}
\filldraw[fill=cyan!20!white] (0, 1.5) ellipse (1 and 1);
\node(1) at (0, 2) {$1$};
\node(2) at (0, 1) {$2$};
\node(Cname) at (0, 3) {$X$};
\filldraw[fill=cyan!20!white] (3, 1) ellipse (1 and 1.5);
\node(5) at (3, 2) {$1$};
\node(6) at (3, 1) {$2$};
\node(7) at (3, 0) {$3$};
\node(Dname) at (3, 3) {$Y$};
\draw[->, very thick, green!50!black] (1) -- (5);
\draw[->, very thick, green!50!black] (2) -- (7);
\end{tikzpicture}%
\ \ .
\]

\textbf{(d)} Find two distinct right inverses of the map%
\[%
\begin{tikzpicture}
\filldraw[fill=cyan!20!white] (0, 1) ellipse (1 and 1.5);
\node(1) at (0, 2) {$1$};
\node(2) at (0, 1) {$2$};
\node(3) at (0, 0) {$3$};
\node(Cname) at (0, 3) {$X$};
\filldraw[fill=cyan!20!white] (3, 1.5) ellipse (1 and 1);
\node(5) at (3, 2) {$1$};
\node(6) at (3, 1) {$2$};
\node(Dname) at (3, 3) {$Y$};
\draw[->, very thick, green!50!black] (1) -- (5);
\draw[->, very thick, green!50!black] (2) -- (6);
\draw[->, very thick, green!50!black] (3) -- (5);
\end{tikzpicture}%
\ \ .
\]

\end{exercise}

\begin{exercise}
For each of the following functions, determine whether it is injective,
surjective and/or bijective: \medskip

\textbf{(a)} The function%
\begin{align*}
f:\mathbb{Q}  &  \rightarrow\mathbb{Q},\\
x  &  \mapsto\dfrac{x}{1+x^{2}}.
\end{align*}

\textbf{(b)} The function%
\begin{align*}
f:\mathbb{Z}  &  \rightarrow\mathbb{Q},\\
x  &  \mapsto\dfrac{x}{1+x^{2}}.
\end{align*}

\textbf{(c)} The function%
\begin{align*}
f:\left\{  \text{finite nonempty subsets of }\mathbb{Z}\right\}   &
\rightarrow\mathbb{Z},\\
S  &  \mapsto\min S
\end{align*}
(which sends each set to its smallest element). \medskip

\textbf{(d)} The function
\begin{align*}
f:\mathbb{Z}\times\mathbb{Z}  &  \rightarrow\mathbb{Z},\\
\left(  x,y\right)   &  \mapsto2x+3y.
\end{align*}

\textbf{(e)} The function%
\begin{align*}
f:\mathbb{N}\times\mathbb{N}  &  \rightarrow\mathbb{N},\\
\left(  x,y\right)   &  \mapsto2x+3y.
\end{align*}

\end{exercise}

\subsection{\label{sec.enum.maps.iso}Isomorphic sets}

As an application of inverses, we can define the concept of isomorphic sets:

\begin{definition}
Let $X$ and $Y$ be two sets. We say that these two sets $X$ and $Y$ are
\textbf{isomorphic as sets} (or, for short, \textbf{isomorphic}, or \textbf{in
bijection}, or \textbf{in one-to-one correspondence}, or \textbf{equinumerous}%
) if there exists a bijective map from $X$ to $Y$. We use the notation
\textquotedblleft$X\cong Y$\textquotedblright\ for \textquotedblleft$X$ and
$Y$ are isomorphic as sets\textquotedblright.
\end{definition}

Note that this relation \textquotedblleft isomorphic as sets\textquotedblright%
\ is symmetric (i.e., if $X$ and $Y$ are isomorphic, then $Y$ and $X$ are
isomorphic). This is because if $f:X\rightarrow Y$ is a bijective map, then
$f$ has an inverse $f^{-1}$ (by Theorem \ref{thm.maps.jects.inv=bij}), and
this inverse $f^{-1}$ is again bijective (since Theorem
\ref{thm.maps.inverse.invinv} shows that $f^{-1}$ again has an inverse).

Some examples:

\begin{itemize}
\item The sets $\left\{  1,2\right\}  $ and $\left\{  1,2,3\right\}  $ are
\textbf{not} isomorphic. In fact, there is no surjective map $f:\left\{
1,2\right\}  \rightarrow\left\{  1,2,3\right\}  $ (since, informally, a map
from $\left\{  1,2\right\}  $ to $\left\{  1,2,3\right\}  $ has only two
arrows, but two arrows cannot hit all three elements of $\left\{
1,2,3\right\}  $). Thus, there is no bijective map $f:\left\{  1,2\right\}
\rightarrow\left\{  1,2,3\right\}  $ either.

\item The sets $\left\{  1,2,3\right\}  $ and $\left\{  6,7,8\right\}  $ are
isomorphic. In fact, the map%
\begin{align*}
\left\{  1,2,3\right\}   &  \rightarrow\left\{  6,7,8\right\}  ,\\
k  &  \mapsto k+5
\end{align*}
(that is, the \textquotedblleft add $5$\textquotedblright\ map) is bijective
(and its inverse sends $k\mapsto k-5$).

\item The sets $\left\{  1,2,3\right\}  $ and $\left\{  3,8,19\right\}  $ are
isomorphic. In fact, the map $f:\left\{  1,2,3\right\}  \rightarrow\left\{
3,8,19\right\}  $ with the table of values%
\[%
\begin{tabular}
[c]{|c||c|c|c|}\hline
$x$ & $1$ & $2$ & $3$\\\hline
$f\left(  x\right)  $ & $3$ & $8$ & $19$\\\hline
\end{tabular}
\]
is bijective.

\item The sets $\left\{  1,2,3\right\}  $ and $\left\{  1,3,5\right\}  $ are
isomorphic. In fact, the map%
\begin{align*}
\left\{  1,2,3\right\}   &  \rightarrow\left\{  1,3,5\right\}  ,\\
k  &  \mapsto2k-1
\end{align*}
is a bijection.

\item The sets $\mathbb{N}$ and $E:=\left\{  \text{all even nonnegative
integers}\right\}  $ are isomorphic, since the map%
\begin{align*}
\mathbb{N}  &  \rightarrow E,\\
n  &  \mapsto2n
\end{align*}
is a bijection.

\item The sets $\mathbb{N}$ and $O:=\left\{  \text{all odd nonnegative
integers}\right\}  $ are isomorphic, since the map%
\begin{align*}
\mathbb{N}  &  \rightarrow O,\\
n  &  \mapsto2n+1
\end{align*}
is a bijection.

\item The sets $\mathbb{N}$ and $\mathbb{Z}$ are isomorphic, since there is a
bijection from $\mathbb{N}$ to $\mathbb{Z}$ that sends%
\begin{align*}
&  0,1,2,3,4,5,6,7,8,\ldots\ \ \ \ \ \ \ \ \ \ \text{to}\\
&  0,1,-1,2,-2,3,-3,4,-4,\ldots,\ \ \ \ \ \ \ \ \ \ \text{respectively.}%
\end{align*}
Explicitly, this bijection $f$ can be defined by the following formula:%
\[
f\left(  n\right)  =%
\begin{cases}
-n/2, & \text{if }n\text{ is even};\\
\left(  n+1\right)  /2, & \text{if }n\text{ is odd}%
\end{cases}
\ \ \ \ \ \ \ \ \ \ \text{for each }n\in\mathbb{N}.
\]
(This formula ensures that the values $f\left(  0\right)  ,f\left(  2\right)
,f\left(  4\right)  ,f\left(  6\right)  ,\ldots$ cover exactly the integers
$0,-1,-2,-3,\ldots$ that are $\leq0$, whereas the values $f\left(  1\right)
,f\left(  3\right)  ,f\left(  5\right)  ,f\left(  7\right)  ,\ldots$ cover
exactly the positive integers $1,2,3,4,\ldots$.)

The inverse $f^{-1}$ of $f$ is the map $\widetilde{f}:\mathbb{Z}%
\rightarrow\mathbb{N}$ defined by%
\[
\widetilde{f}\left(  n\right)  =%
\begin{cases}
-2n, & \text{if }n\leq0;\\
2n-1, & \text{if }n>0
\end{cases}
\ \ \ \ \ \ \ \ \ \ \text{for each }n\in\mathbb{Z}.
\]
It is a nice exercise to check that the maps $f$ and $\widetilde{f}$ really
are inverses of each other (and well-defined!), thus are bijections.

There are, of course, many other bijections from $\mathbb{N}$ to $\mathbb{Z}$.

\item The sets $\mathbb{N}$ and $\mathbb{Q}$ are isomorphic, since there is a
bijection from $\mathbb{N}$ to $\mathbb{Q}$ that sends%
\begin{align*}
&  0,1,2,3,4,5,6,7,8,9,\ldots\ \ \ \ \ \ \ \ \ \ \text{to}\\
&  \underbrace{\dfrac{-1}{1},\dfrac{0}{1},\dfrac{1}{1}}_{\substack{\text{all
reduced fractions}\\\text{whose numerator and}\\\text{denominator are }%
\leq1\\\text{in absolute value}\\\text{(ordered from smallest}\\\text{to
largest)}}},\ \underbrace{\dfrac{-2}{1},\dfrac{-1}{2},\dfrac{1}{2},\dfrac
{2}{1}}_{\substack{\text{all reduced fractions}\\\text{whose numerator
and}\\\text{denominator are }\leq2\\\text{but not }\leq1\\\text{in absolute
value}\\\text{(ordered from smallest}\\\text{to largest)}}%
},\ \underbrace{\dfrac{-3}{1},\dfrac{-3}{2},\dfrac{-2}{3},\dfrac{-1}{3}%
,\dfrac{1}{3},\dfrac{2}{3},\dfrac{3}{2},\dfrac{3}{1}}_{\substack{\text{all
reduced fractions}\\\text{whose numerator and}\\\text{denominator are }%
\leq3\\\text{but not }\leq2\\\text{in absolute value}\\\text{(ordered from
smallest}\\\text{to largest)}}},\ldots,
\end{align*}
respectively. (To be precise, we must only allow \textbf{fully reduced}
fractions -- i.e., fractions $\dfrac{a}{b}$ with $a\in\mathbb{Z}$ and
$b\in\left\{  1,2,3,\ldots\right\}  $ and $\gcd\left(  a,b\right)  =1$ -- in
order to avoid having the same rational number appear twice.)

\item The sets $\mathbb{N}$ and $\mathbb{N}\times\mathbb{N}$ are isomorphic,
since there is a bijection $f$ from $\mathbb{N}$ to $\mathbb{N}\times
\mathbb{N}$ that sends%
\begin{align*}
&  0,1,2,3,4,5,6,7,8,9,\ldots\ \ \ \ \ \ \ \ \ \ \text{to}\\
&  \underbrace{\left(  0,0\right)  }_{\substack{\text{all pairs}\\\text{whose
entries}\\\text{sum to }0}},\ \underbrace{\left(  0,1\right)  ,\ \left(
1,0\right)  }_{\substack{\text{all pairs}\\\text{whose entries}\\\text{sum to
}1\\\text{(ordered by}\\\text{increasing}\\\text{first entry)}}%
},\ \underbrace{\left(  0,2\right)  ,\ \left(  1,1\right)  ,\ \left(
2,0\right)  }_{\substack{\text{all pairs}\\\text{whose entries}\\\text{sum to
}2\\\text{(ordered by}\\\text{increasing}\\\text{first entry)}}%
},\ \underbrace{\left(  0,3\right)  ,\ \left(  1,2\right)  ,\ \left(
2,1\right)  ,\ \left(  3,0\right)  }_{\substack{\text{all pairs}\\\text{whose
entries}\\\text{sum to }3\\\text{(ordered by}\\\text{increasing}\\\text{first
entry)}}},\ \ldots,
\end{align*}
respectively. The inverse $f^{-1}:\mathbb{N}\times\mathbb{N}\rightarrow
\mathbb{N}$ of this bijection can actually be described by an explicit
formula:%
\[
f^{-1}\left(  n,m\right)  =\dfrac{\left(  n+m\right)  \left(  n+m+1\right)
}{2}+n
\]
(nice and not-too-easy exercise: prove this!). This is the so-called
\href{https://en.wikipedia.org/wiki/Pairing_function}{Cantor pairing function}.

\item The sets $\mathbb{N}$ and $\mathbb{Z}\times\mathbb{Z}$ are isomorphic.
The easiest way to find a bijection from $\mathbb{N}$ to $\mathbb{Z}%
\times\mathbb{Z}$ is by combining some of our bijections constructed above:
Namely, pick a bijection $h$ from $\mathbb{N}$ to $\mathbb{N}\times\mathbb{N}$
and a bijection $g$ from $\mathbb{N}$ to $\mathbb{Z}$ (we already have
constructed such bijections). Consider the map $g\ast g:\mathbb{N}%
\times\mathbb{N}\rightarrow\mathbb{Z}\times\mathbb{Z}$ that sends each pair
$\left(  i,j\right)  \in\mathbb{N}\times\mathbb{N}$ to the pair $\left(
g\left(  i\right)  ,g\left(  j\right)  \right)  \in\mathbb{Z}\times\mathbb{Z}$
(this is an instance of the definition of $f\ast g$ in Exercise
\ref{exe.23hw5.1+}). This map $g\ast g$ is injective (by Exercise
\ref{exe.23hw5.1+} \textbf{(a)}) and surjective (by Exercise
\ref{exe.23hw5.1+} \textbf{(b)}), hence bijective. Thus, the composition
$\left(  g\ast g\right)  \circ h:\mathbb{N}\rightarrow\mathbb{Z}%
\times\mathbb{Z}$ is also bijective (by Exercise \ref{exe.maps.jects.compos}
\textbf{(g)}, since $h$ is also bijective). Hence, we have found a bijection
from $\mathbb{N}$ to $\mathbb{Z}\times\mathbb{Z}$, so that the sets
$\mathbb{N}$ and $\mathbb{Z}\times\mathbb{Z}$ are isomorphic.

\item The sets $\mathbb{N}$ and $\mathbb{N}\times\mathbb{N}\times\mathbb{N}$
are isomorphic. Indeed, we can pick a bijection $h$ from $\mathbb{N}$ to
$\mathbb{N}\times\mathbb{N}$ (for example, the one we constructed above) and
use it to construct a bijection from $\mathbb{N}$ to $\mathbb{N}%
\times\mathbb{N}\times\mathbb{N}$ as follows: We define the map%
\begin{align*}
\mathbb{N}\times\mathbb{N}\times\mathbb{N}  &  \rightarrow\mathbb{N},\\
\left(  a,b,c\right)   &  \mapsto h^{-1}\left(  h^{-1}\left(  a,b\right)
,c\right)  .
\end{align*}
It is easy to see that this map is a bijection (Exercise
\ref{exe.maps.jects.compos} \textbf{(g)} comes in handy here!), and thus its
inverse is a bijection from $\mathbb{N}$ to $\mathbb{N}\times\mathbb{N}%
\times\mathbb{N}$, which shows that the sets $\mathbb{N}$ and $\mathbb{N}%
\times\mathbb{N}\times\mathbb{N}$ are isomorphic.

A nicer way to express the same argument is as follows (some details are left
to the reader): Show that if $A$, $B$ and $C$ are three sets such that $A\cong
B$ and $B\cong C$, then $A\cong C$ (this is called \textbf{transitivity of
isomorphism}, and follows easily from Exercise \ref{exe.maps.jects.compos}
\textbf{(g)}). Show that if $A$, $B$ and $C$ are three sets such that $A\cong
B$, then $A\times C\cong B\times C$ (this is easy using Exercise
\ref{exe.23hw5.1+}). Hence, from $\mathbb{N}\cong\mathbb{N}\times\mathbb{N}$,
we obtain $\mathbb{N}\times\mathbb{N}\cong\left(  \mathbb{N}\times
\mathbb{N}\right)  \times\mathbb{N}$. But $\left(  \mathbb{N}\times
\mathbb{N}\right)  \times\mathbb{N}\cong\mathbb{N}\times\mathbb{N}%
\times\mathbb{N}$ using the \textquotedblleft unpacking\textquotedblright%
\ bijection
\begin{align*}
\left(  \mathbb{N}\times\mathbb{N}\right)  \times\mathbb{N}  &  \rightarrow
\mathbb{N}\times\mathbb{N}\times\mathbb{N},\\
\left(  \left(  a,b\right)  ,c\right)   &  \mapsto\left(  a,b,c\right)  .
\end{align*}
Hence, altogether, $\mathbb{N}\cong\mathbb{N}\times\mathbb{N}\cong\left(
\mathbb{N}\times\mathbb{N}\right)  \times\mathbb{N}\cong\mathbb{N}%
\times\mathbb{N}\times\mathbb{N}$, so that $\mathbb{N}\cong\mathbb{N}%
\times\mathbb{N}\times\mathbb{N}$ (by transitivity of isomorphism).

Likewise, it can be shown that $\mathbb{N}$ is isomorphic to any Cartesian
product of the form $\underbrace{\mathbb{N}\times\mathbb{N}\times\cdots
\times\mathbb{N}}_{k\text{ times}}$ of $\mathbb{N}$'s for $k>0$.

\item The sets $\mathbb{N}$ and $\mathbb{R}$ are \textbf{not} isomorphic,
i.e., there exists no bijection from $\mathbb{N}$ to $\mathbb{R}$. Informally
speaking, this is because there are \textquotedblleft a lot
more\textquotedblright\ real numbers than there are nonnegative integers. This
is not a proof at all (after all, $\mathbb{N}$ and $\mathbb{Q}$ are
isomorphic, despite the rational numbers seemingly outnumbering the
nonnegative integers!); an actual proof can be found (e.g.) in \cite[Theorem
6.2.21]{Newste23} or in \cite[Corollary 8.1.17]{LeLeMe16}.
\end{itemize}

Sets that are isomorphic to $\mathbb{N}$ are said to be \textbf{countably
infinite}. Thus, our above examples show that the sets $\mathbb{Z}$,
$\mathbb{Q}$, $\mathbb{N}\times\mathbb{N}$, $\mathbb{Z}\times\mathbb{Z}$ and
$\mathbb{N}\times\mathbb{N}\times\mathbb{N}$ are all countably infinite, but
the set $\mathbb{R}$ is not. More about countably infinite sets can be found
in \cite[\S 6.2]{Newste23} or \cite[\S 4.4]{Mileti22}.

\newpage

\section{\label{chp.enum2}Enumeration revisited}

\subsection{\label{sec.enum.count-formal}Counting, formally}

\subsubsection{\label{subsec.enum.count-formal.def}Definition}

As you might have noticed, isomorphic sets (at least when they are finite)
have the same number of elements -- i.e., the same size. We shall now use this
to \textbf{define} the size of a set!

First, some notations:

\begin{definition}
\label{def.sets.ive-n}\textbf{(a)} If $n\in\mathbb{N}$, then $\left[
n\right]  $ shall mean the set $\left\{  1,2,\ldots,n\right\}  $.

For example, $\left[  3\right]  =\left\{  1,2,3\right\}  $ and $\left[
7\right]  =\left\{  1,2,3,4,5,6,7\right\}  $ and $\left[  0\right]
=\varnothing$ and $\left[  1\right]  =\left\{  1\right\}  $. \medskip

\textbf{(b)} If $a,b\in\mathbb{Z}$, then $\left[  a,b\right]  $ shall mean the
set
\begin{align*}
\left\{  a,a+1,a+2,\ldots,b\right\}   &  =\left\{  \text{all integers }x\text{
satisfying }a\leq x\leq b\right\} \\
&  =\left\{  x\in\mathbb{Z}\ \mid\ a\leq x\leq b\right\}  .
\end{align*}
If $a>b$, then this is understood to be the empty set.

For example, $\left[  2,6\right]  =\left\{  2,3,4,5,6\right\}  $ and $\left[
3,3\right]  =\left\{  3\right\}  $ and $\left[  4,2\right]  =\varnothing$.
\end{definition}

Now, let us define the size of a finite set:

\begin{definition}
\label{def.sets.size.1}Let $n\in\mathbb{N}$. A set $S$ is said to have
\textbf{size }$n$ if $S$ is isomorphic to $\left[  n\right]  $ (that is, if
there exists a bijection from $S$ to $\left[  n\right]  $).
\end{definition}

For example:

\begin{itemize}
\item The set $\left\{  \text{\textquotedblleft cat\textquotedblright%
},\ \text{\textquotedblleft dog\textquotedblright},\ \text{\textquotedblleft
rat\textquotedblright}\right\}  $ has size $3$, since the map%
\begin{align*}
\left\{  \text{\textquotedblleft cat\textquotedblright}%
,\ \text{\textquotedblleft dog\textquotedblright},\ \text{\textquotedblleft
rat\textquotedblright}\right\}   &  \rightarrow\left[  3\right]  ,\\
\text{\textquotedblleft cat\textquotedblright}  &  \mapsto1,\\
\text{\textquotedblleft dog\textquotedblright}  &  \mapsto2,\\
\text{\textquotedblleft rat\textquotedblright}  &  \mapsto3
\end{align*}
is a bijection.

\item The set $\left[  4,7\right]  =\left\{  4,5,6,7\right\}  $ has size $4$,
since the map%
\begin{align*}
\left[  4,7\right]   &  \rightarrow\left[  4\right]  ,\\
k  &  \mapsto k-3
\end{align*}
is a bijection.

\item The set $\mathbb{N}$ is infinite, so there is no bijection from
$\mathbb{N}$ to $\left[  n\right]  $ for any $n\in\mathbb{N}$. Thus,
$\mathbb{N}$ does not have size $n$ for any $n\in\mathbb{N}$.
\end{itemize}

Here is another equivalent definition of size:

\begin{definition}
\label{def.sets.size.2}We define the notion of a \textquotedblleft set of size
$n$\textquotedblright\ recursively as follows: \medskip

\textbf{(a)} A set $S$ is said to have \textbf{size }$0$ if and only if it is
empty. \medskip

\textbf{(b)} Let $n$ be a positive integer. A set $S$ is said to have
\textbf{size }$n$ if and only if there exists an $s\in S$ such that
$S\setminus\left\{  s\right\}  $ has size $n-1$.
\end{definition}

In other words, a set has size $n$ (for $n>0$) if and only if we can remove a
single element from it and obtain a set of size $n-1$. This is a recursive
definition, as it reduces the question \textquotedblleft what is a set of size
$n$\textquotedblright\ to the (simpler) question \textquotedblleft what is a
set of size $n-1$\textquotedblright.

The following fact is not obvious, but can be proved:

\begin{theorem}
\label{thm.sets.size.1=2}\textbf{(a)} The above two definitions of size
(Definition \ref{def.sets.size.1} and Definition \ref{def.sets.size.2}) are
equivalent. \medskip

\textbf{(b)} The size of a finite set is determined uniquely -- i.e., a set
cannot have two different sizes at the same time.
\end{theorem}

Now, we are ready to introduce some notations for sizes of sets:

\begin{definition}
\textbf{(a)} An $n$\textbf{-element set} (for some $n\in\mathbb{N}$) means a
set of size $n$. \medskip

\textbf{(b)} A set is said to be \textbf{finite} if it has size $n$ for some
$n\in\mathbb{N}$. \medskip

\textbf{(c)} If $S$ is a finite set, then $\left\vert S\right\vert $ shall
denote the size of $S$ (which is unique because of Theorem
\ref{thm.sets.size.1=2} \textbf{(b)}). \medskip

\textbf{(d)} We also refer to $\left\vert S\right\vert $ as the
\textbf{cardinality} of $S$, or as the \textbf{number} of elements of $S$. In
particular, the \textbf{number} of some things means the size of the set of
these things.
\end{definition}

Thus, our examples above show that
\[
\left\vert \left\{  \text{\textquotedblleft cat\textquotedblright%
},\ \text{\textquotedblleft dog\textquotedblright},\ \text{\textquotedblleft
rat\textquotedblright}\right\}  \right\vert =3\ \ \ \ \ \ \ \ \ \ \text{and}%
\ \ \ \ \ \ \ \ \ \ \left\vert \left[  4,7\right]  \right\vert =4.
\]
The number of odd integers between $4$ and $10$ is the size of the set
\[
\left\{  \text{odd integers between }4\text{ and }10\right\}  =\left\{
5,7,9\right\}  ,
\]
and thus equals $3$.

(Don't forget that sets cannot \textquotedblleft contain an element more than
once\textquotedblright. Thus, $\left\vert \left\{  5,6,5\right\}  \right\vert
$ is $2$, not $3$.)

\subsubsection{Rules for sizes of finite sets}

We have defined the size $\left\vert S\right\vert $ of a finite set $S$ in
Subsection \ref{subsec.enum.count-formal.def}. Let us now state some rules for
these sizes that make them easier to compute. We will not prove these rules,
as they are all dictated by common sense and their rigorous proofs would
reasonably belong into a text on formalized foundations of mathematics.
\medskip

The most important rule is the following:

\begin{theorem}
[Bijection Principle]\label{thm.count.rules.bij}Let $A$ and $B$ be two finite
sets. Then, $\left\vert A\right\vert =\left\vert B\right\vert $ if and only if
there exists a bijection from $A$ to $B$.
\end{theorem}

(As we recall, a \textquotedblleft bijection\textquotedblright\ means a
bijective map. By Theorem \ref{thm.maps.jects.inv=bij}, this is the same as a
map that has an inverse. We also note that the claim of Theorem
\ref{thm.count.rules.bij} can be restated as \textquotedblleft$\left\vert
A\right\vert =\left\vert B\right\vert $ if and only if $A\cong B$%
\textquotedblright.)

The intuition behind Theorem \ref{thm.count.rules.bij} is that a bijection
between two finite sets $A$ and $B$ is a way to \textquotedblleft pair
up\textquotedblright\ each element of $A$ with an element of $B$ (its image
under the bijection), and thus ensures that $A$ has equally many elements as
$B$. Conversely, if two finite sets $A$ and $B$ have equally many elements --
say, $n$ elements --, then we can list the elements of $A$ as $a_{1}%
,a_{2},\ldots,a_{n}$ and list the elements of $B$ as $b_{1},b_{2},\ldots
,b_{n}$, and then construct a bijection from $A$ to $B$ by sending each
$a_{i}$ to the corresponding $b_{i}$. These justifications fall short of the
requirements of a rigorous proof, but they should make Theorem
\ref{thm.count.rules.bij} intuitively transparent. \medskip

The next rule is obvious from one of our definitions of size:

\begin{theorem}
\label{thm.count.rules.n}For each $n\in\mathbb{N}$, we have $\left\vert
\left[  n\right]  \right\vert =n$.
\end{theorem}

Here, we recall that $\left[  n\right]  $ means the set $\left\{
1,2,\ldots,n\right\}  $ consisting of the first $n$ positive integers.
\medskip

The next rule classifies sets of small size:

\begin{theorem}
\label{thm.count.rules.012}Let $S$ be a set. Then: \medskip

\textbf{(a)} We have $\left\vert S\right\vert =0$ if and only if
$S=\varnothing$ (that is, if $S$ is empty). \medskip

\textbf{(b)} We have $\left\vert S\right\vert =1$ if and only if $S=\left\{
s\right\}  $ for a single element $s$. \medskip

\textbf{(c)} We have $\left\vert S\right\vert =2$ if and only if $S=\left\{
s,t\right\}  $ for two distinct elements $s$ and $t$.
\end{theorem}

The next rule says that inserting a new element into a finite set increases
the size of this set by $1$:

\begin{theorem}
\label{thm.count.rules.+1}Let $S$ be a finite set. Let $t$ be an object such
that $t\notin S$ (that is, $t$ does not belong to $S$). Then,%
\[
\left\vert S\cup\left\{  t\right\}  \right\vert =\left\vert S\right\vert +1.
\]

\end{theorem}

Here are some more substantial facts:

\begin{theorem}
[Sum rule for two sets]\label{thm.count.rules.sum-2}Let $A$ and $B$ be two
disjoint finite sets. (Recall that \textquotedblleft
disjoint\textquotedblright\ means $A\cap B=\varnothing$.) Then, the set $A\cup
B$ is again finite, and has size%
\[
\left\vert A\cup B\right\vert =\left\vert A\right\vert +\left\vert
B\right\vert .
\]

\end{theorem}

\begin{theorem}
[Sum rule for $k$ sets]\label{thm.count.rules.sum-k}Let $A_{1},A_{2}%
,\ldots,A_{k}$ be $k$ disjoint finite sets. (Recall that \textquotedblleft
disjoint\textquotedblright\ for $k$ sets\ means that any two of them are
disjoint, i.e., that $A_{i}\cap A_{j}=\varnothing$ for any $i<j$.) Then, the
set $A_{1}\cup A_{2}\cup\cdots\cup A_{k}$ is finite, and has size%
\[
\left\vert A_{1}\cup A_{2}\cup\cdots\cup A_{k}\right\vert =\left\vert
A_{1}\right\vert +\left\vert A_{2}\right\vert +\cdots+\left\vert
A_{k}\right\vert .
\]

\end{theorem}

\begin{theorem}
[Difference rule]\label{thm.count.rules.diff}Let $T$ be a subset of a finite
set $S$. Then: \medskip

\textbf{(a)} The set $T$ is finite, and its size $\left\vert T\right\vert $
satisfies $\left\vert T\right\vert \leq\left\vert S\right\vert $. \medskip

\textbf{(b)} We have $\left\vert S\setminus T\right\vert =\left\vert
S\right\vert -\left\vert T\right\vert $. \medskip

\textbf{(c)} If $\left\vert T\right\vert =\left\vert S\right\vert $, then
$T=S$.
\end{theorem}

The following theorem has been previously stated (without using the
\textquotedblleft size\textquotedblright\ terminology) as Theorem
\ref{thm.count.pairs.nm-x}:

\begin{theorem}
[Product rule for two sets]\label{thm.count.rules.prod-2}Let $A$ and $B$ be
any finite sets. Then, the set%
\[
A\times B=\left\{  \text{all pairs }\left(  a,b\right)  \text{ with }a\in
A\text{ and }b\in B\right\}
\]
is again finite and has size%
\[
\left\vert A\times B\right\vert =\left\vert A\right\vert \cdot\left\vert
B\right\vert .
\]

\end{theorem}

Likewise, the following theorem was Theorem \ref{thm.count.triples.nmp-x-n}:

\begin{theorem}
[Product rule for $k$ sets]\label{thm.count.rules.prod-k}Let $A_{1}%
,A_{2},\ldots,A_{k}$ be any $k$ finite sets. Then, the set%
\[
A_{1}\times A_{2}\times\cdots\times A_{k}=\left\{  \text{all }k\text{-tuples
}\left(  a_{1},a_{2},\ldots,a_{k}\right)  \text{ with }a_{i}\in A_{i}\text{
for each }i\in\left[  k\right]  \right\}
\]
is again finite and has size%
\[
\left\vert A_{1}\times A_{2}\times\cdots\times A_{k}\right\vert =\left\vert
A_{1}\right\vert \cdot\left\vert A_{2}\right\vert \cdot\cdots\cdot\left\vert
A_{k}\right\vert .
\]

\end{theorem}

All the above theorems are foundational, and are perhaps the reason why the
arithmetic operations $+$, $-$ and $\cdot$ on nonnegative integers were
introduced some millennia ago. Nevertheless, they can be rigorously proved,
but this is not something we will do here.\footnote{For instance, Theorem
\ref{thm.count.rules.prod-2} can be proved by induction on $\left\vert
B\right\vert $ using Theorem \ref{thm.count.rules.sum-k} and Theorem
\ref{thm.count.rules.bij}, whereas Theorem \ref{thm.count.rules.prod-k} can be
proved by induction on $k$ using Theorem \ref{thm.count.rules.prod-2}.}

The above theorems are known as \textquotedblleft basic counting
rules\textquotedblright\ or \textquotedblleft counting
principles\textquotedblright. There are a few more counting principles, which
we might state later on.

\subsubsection{$A\cup B$ and $A\cap B$ revisited}

As a first application of these rules, let us derive the following
\textquotedblleft generalized sum rule for two sets\textquotedblright:

\begin{theorem}
\label{thm.count.rules.AuB-AnB}Let $A$ and $B$ be two finite sets (not
necessarily disjoint). Then, the set $A\cup B$ is finite and has size%
\[
\left\vert A\cup B\right\vert =\left\vert A\right\vert +\left\vert
B\right\vert -\left\vert A\cap B\right\vert .
\]

\end{theorem}

\begin{proof}
[Partial proof.]We shall take for granted that $A\cup B$ is finite, and only
prove the equality $\left\vert A\cup B\right\vert =\left\vert A\right\vert
+\left\vert B\right\vert -\left\vert A\cap B\right\vert $ here.

We first claim that%
\begin{equation}
\left(  A\cup B\right)  \setminus A=B\setminus\left(  A\cap B\right)  .
\label{pf.thm.count.rules.AuB-AnB.1}%
\end{equation}
This equality is obvious using Venn diagrams, but let us prove it rigorously
using \textquotedblleft element chasing\textquotedblright:\footnote{It is
worth noting that both sides of (\ref{pf.thm.count.rules.AuB-AnB.1}) are equal
to $B\setminus A$. However, we will not need this fact.}

\begin{proof}
[Proof of (\ref{pf.thm.count.rules.AuB-AnB.1}).]Let us first prove $\left(
A\cup B\right)  \setminus A\subseteq B\setminus\left(  A\cap B\right)  $. In
order to do so, we must show that each $x\in\left(  A\cup B\right)  \setminus
A$ belongs to $B\setminus\left(  A\cap B\right)  $. Let us do this: Let
$x\in\left(  A\cup B\right)  \setminus A$. Thus, $x\in A\cup B$ but $x\notin
A$. From $x\in A\cup B$, we see that $x\in A$ or $x\in B$. But the first of
these two possibilities is impossible (since $x\notin A$). Thus, the second
possibility must hold. In other words, we have $x\in B$. Furthermore, we have
$x\notin A\cap B$ (since $x\in A\cap B$ would entail $x\in A\cap B\subseteq
A$, which would contradict $x\notin A$). Combining $x\in B$ with $x\notin
A\cap B$, we obtain $x\in B\setminus\left(  A\cap B\right)  $.

Forget that we fixed $x$. We thus have shown that each $x\in\left(  A\cup
B\right)  \setminus A$ belongs to $B\setminus\left(  A\cap B\right)  $. In
other words, $\left(  A\cup B\right)  \setminus A\subseteq B\setminus\left(
A\cap B\right)  $.

Next, let us prove that $B\setminus\left(  A\cap B\right)  \subseteq\left(
A\cup B\right)  \setminus A$. To do so, we must show that each $x\in
B\setminus\left(  A\cap B\right)  $ belongs to $\left(  A\cup B\right)
\setminus A$. We do this as follows: Let $x\in B\setminus\left(  A\cap
B\right)  $. Thus, $x\in B$ but $x\notin A\cap B$. If we had $x\in A$, then we
would have $x\in A\cap B$ (since $x\in A$ and $x\in B$), which would
contradict $x\notin A\cap B$. Hence, we cannot have $x\in A$. Thus, $x\notin
A$. Also, $x\in B\subseteq A\cup B$. Combining this with $x\notin A$, we find
$x\in\left(  A\cup B\right)  \setminus A$.

Forget that we fixed $x$. We thus have shown that each $x\in B\setminus\left(
A\cap B\right)  $ belongs to $\left(  A\cup B\right)  \setminus A$. In other
words, $B\setminus\left(  A\cap B\right)  \subseteq\left(  A\cup B\right)
\setminus A$.

Now, combining the two inclusions%
\[
\left(  A\cup B\right)  \setminus A\subseteq B\setminus\left(  A\cap B\right)
\ \ \ \ \ \ \ \ \ \ \text{and}\ \ \ \ \ \ \ \ \ \ B\setminus\left(  A\cap
B\right)  \subseteq\left(  A\cup B\right)  \setminus A,
\]
we obtain $\left(  A\cup B\right)  \setminus A=B\setminus\left(  A\cap
B\right)  $. Thus, (\ref{pf.thm.count.rules.AuB-AnB.1}) is proved.]
\end{proof}

We now step to the counting. Taking sizes on both sides of
(\ref{pf.thm.count.rules.AuB-AnB.1}), we obtain%
\begin{equation}
\left\vert \left(  A\cup B\right)  \setminus A\right\vert =\left\vert
B\setminus\left(  A\cap B\right)  \right\vert .
\label{pf.thm.count.rules.AuB-AnB.2}%
\end{equation}

But $A$ is a subset of $A\cup B$. Thus, the difference rule (Theorem
\ref{thm.count.rules.diff} \textbf{(b)}, applied to $S=A\cup B$ and $T=A$)
yields%
\begin{equation}
\left\vert \left(  A\cup B\right)  \setminus A\right\vert =\left\vert A\cup
B\right\vert -\left\vert A\right\vert . \label{pf.thm.count.rules.AuB-AnB.3}%
\end{equation}

Also, $A\cap B$ is a subset of $B$. Thus, the difference rule (Theorem
\ref{thm.count.rules.diff} \textbf{(b)}, applied to $S=B$ and $T=A\cap B$)
yields%
\begin{equation}
\left\vert B\setminus\left(  A\cap B\right)  \right\vert =\left\vert
B\right\vert -\left\vert A\cap B\right\vert .
\label{pf.thm.count.rules.AuB-AnB.4}%
\end{equation}

But we know from (\ref{pf.thm.count.rules.AuB-AnB.2}) that the left hand sides
of the two equalities (\ref{pf.thm.count.rules.AuB-AnB.3}) and
(\ref{pf.thm.count.rules.AuB-AnB.4}) are equal. Thus, their right hand sides
are also equal. In other words,%
\[
\left\vert A\cup B\right\vert -\left\vert A\right\vert =\left\vert
B\right\vert -\left\vert A\cap B\right\vert .
\]
Solving this for $\left\vert A\cup B\right\vert $, we find%
\[
\left\vert A\cup B\right\vert =\left\vert A\right\vert +\left\vert
B\right\vert -\left\vert A\cap B\right\vert .
\]
This proves Theorem \ref{thm.count.rules.AuB-AnB}.

(A nicer proof can be given using finite sums; this is done, e.g., in
\cite[Lecture 19, \S 2.7]{22fco}.)
\end{proof}

Note that Theorem \ref{thm.count.rules.AuB-AnB} has an analogue for three
sets: If $A,B,C$ are three finite sets, then%
\[
\left\vert A\cup B\cup C\right\vert =\left\vert A\right\vert +\left\vert
B\right\vert +\left\vert C\right\vert -\left\vert A\cap B\right\vert
-\left\vert A\cap C\right\vert -\left\vert B\cap C\right\vert +\left\vert
A\cap B\cap C\right\vert .
\]
More generally, such a formula can be stated for any $k$ finite sets, and is
known as the \textquotedblleft principle of inclusion and
exclusion\textquotedblright\ or \textquotedblleft Sylvester's sieve
formula\textquotedblright. See \cite[Lecture 19, \S 2.7]{22fco} or any
textbook on combinatorics.

\begin{exercise}
Let $A$ and $B$ be two finite sets. Prove that%
\[
\left\vert A\cup B\right\vert =\left\vert A\setminus B\right\vert +\left\vert
B\setminus A\right\vert +\left\vert A\cap B\right\vert .
\]

\end{exercise}

\subsubsection{Finiteness principles}

It is often clear when a set in combinatorics is finite, so we will usually
omit the rigorous proof of finiteness. However, let us state the principles
that are used in such proofs:

\begin{theorem}
[Finiteness principles]\label{thm.count.finite}\textbf{(a)} Any subset of a
finite set is finite. \medskip

\textbf{(b)} If there exists a bijection between two sets $A$ and $B$, and if
one of these two sets $A$ and $B$ is finite, then so is the other. \medskip

\textbf{(c)} If $A_{1},A_{2},\ldots,A_{n}$ are finitely many finite sets, then
their union $A_{1}\cup A_{2}\cup\cdots\cup A_{n}$ is finite. \medskip

\textbf{(d)} If $A_{1},A_{2},\ldots,A_{n}$ are finitely many finite sets, then
their Cartesian product $A_{1}\times A_{2}\times\cdots\times A_{n}$ is finite.
\end{theorem}

\subsection{Redoing some proofs rigorously}

Previously (in Section \ref{sec.enum.count-informal}), we have proved some
results using informal counting arguments. Let us now revisit these results
and make these arguments rigorous.

\subsubsection{Integers in an interval}

We recall that the notation $\left[  a,b\right]  $ means the set%
\[
\left\{  a,\ a+1,\ a+2,\ \ldots,\ b\right\}  =\left\{  x\in\mathbb{Z}%
\ \mid\ a\leq x\leq b\right\}
\]
whenever $a$ and $b$ are two integers. In particular, $\left[  n\right]
=\left[  1,n\right]  $ for every $n\in\mathbb{N}$.

We begin with Proposition \ref{prop.count.fencepost} (rewritten using the
notation $\left[  a,b\right]  $):

\begin{proposition}
\label{prop.count.fencepost-repeat}Let $a,b\in\mathbb{Z}$ be such that $a\leq
b+1$.

Then, there are exactly $b-a+1$ numbers in the set $\left[  a,b\right]  $. In
other words, there are exactly $b-a+1$ integers between $a$ and $b$ (inclusive).
\end{proposition}

Back in Section \ref{sec.enum.count-informal}, we proved this informally by
inducting on $b$. This proof can be trivially made rigorous; the induction
step relies on Theorem \ref{thm.count.rules.+1} (because $\left[
a,b+1\right]  =\left[  a,b\right]  \cup\left\{  b+1\right\}  $ and
$b+1\notin\left[  a,b\right]  $).

But there is also a more direct proof:

\begin{proof}
[Second proof of Proposition \ref{prop.count.fencepost-repeat}.]Consider the
map%
\begin{align*}
f:\underbrace{\left[  b-a+1\right]  }_{=\left\{  1,2,\ldots,b-a+1\right\}  }
&  \rightarrow\underbrace{\left[  a,b\right]  }_{=\left\{  a,a+1,\ldots
,b\right\}  },\\
i  &  \mapsto i+\left(  a-1\right)  .
\end{align*}
This map $f$ just adds $a-1$ to its input. (Informally, we can view it as
moving numbers to the right by $a-1$ units on the number line.)

It is easy to see that this map $f$ has an inverse: Namely, the map%
\begin{align*}
\left[  a,b\right]   &  \rightarrow\left[  b-a+1\right]  ,\\
j  &  \mapsto j-\left(  a-1\right)
\end{align*}
is an inverse of $f$ (since subtraction undoes addition). Thus, the map $f$ is
bijective (by Theorem \ref{thm.maps.jects.inv=bij}), i.e., is a bijection.
Hence, there is a bijection from $\left[  b-a+1\right]  $ to $\left[
a,b\right]  $ (namely, $f$). The bijection principle (Theorem
\ref{thm.count.rules.bij}, applied to $A=\left[  b-a+1\right]  $ and
$B=\left[  a,b\right]  $) thus yields
\[
\left\vert \left[  b-a+1\right]  \right\vert =\left\vert \left[  a,b\right]
\right\vert .
\]
Hence,%
\[
\left\vert \left[  a,b\right]  \right\vert =\left\vert \left[  b-a+1\right]
\right\vert =b-a+1
\]
(by Theorem \ref{thm.count.rules.n}, since $a\leq b+1$ yields $b-a+1\in
\mathbb{N}$). In other words, there are exactly $b-a+1$ numbers in the set
$\left[  a,b\right]  $. This proves Proposition
\ref{prop.count.fencepost-repeat} again.
\end{proof}

We could also reprove Proposition \ref{prop.count.odd-0-n} rigorously, but
(again) the proof we gave was already rigorous enough; we just need to rewrite
it using Theorem \ref{thm.count.rules.+1}.

\subsubsection{Counting all subsets}

Now, we recall Theorem \ref{thm.count.subsets.all.n} (but shorten it using the
notation $\left[  n\right]  $ for $\left\{  1,2,\ldots,n\right\}  $):

\begin{theorem}
\label{thm.count.subsets.all.n-repeat}Let $n\in\mathbb{N}$. Then,%
\[
\left(  \text{\# of subsets of }\left[  n\right]  \right)  =2^{n}.
\]

\end{theorem}

The proof we gave in Section \ref{sec.enum.subsets} had some informal steps;
let us now make it rigorous:\footnote{Most of the proof below is copied
verbatim from Section \ref{sec.enum.subsets}.}

\begin{proof}
[Rigorous proof of Theorem \ref{thm.count.subsets.all.n-repeat}.]We induct on
$n$.

The \textit{base case} ($n=0$) is easy: The set $\left[  0\right]  $ is empty,
and thus its only subset is $\left\{  {}\right\}  $ itself; hence, the \# of
subsets of $\left[  0\right]  $ is $1=2^{0}$. In other words, Theorem
\ref{thm.count.subsets.all.n-repeat} holds for $n=0$.

\textit{Induction step:} We proceed from $n-1$ to $n$. Thus, let $n$ be a
positive integer. We assume (as the induction hypothesis) that Theorem
\ref{thm.count.subsets.all.n-repeat} holds for $n-1$ instead of $n$, and we
set out to prove that it holds for $n$.

So our induction hypothesis says that%
\[
\left(  \text{\# of subsets of }\left[  n-1\right]  \right)  =2^{n-1}.
\]
Our goal is to prove that%
\[
\left(  \text{\# of subsets of }\left[  n\right]  \right)  =2^{n}.
\]

We define

\begin{itemize}
\item a \textbf{red set} to be a subset of $\left[  n\right]  $ that contains
$n$;

\item a \textbf{green set} to be a subset of $\left[  n\right]  $ that does
not contain $n$.
\end{itemize}

For example, if $n=3$, then the red sets are%
\[
\left\{  3\right\}  ,\ \ \ \ \ \ \ \ \ \ \left\{  1,3\right\}
,\ \ \ \ \ \ \ \ \ \ \left\{  2,3\right\}  ,\ \ \ \ \ \ \ \ \ \ \left\{
1,2,3\right\}  ,
\]
whereas the green sets are%
\[
\left\{  {}\right\}  ,\ \ \ \ \ \ \ \ \ \ \left\{  1\right\}
,\ \ \ \ \ \ \ \ \ \ \left\{  2\right\}  ,\ \ \ \ \ \ \ \ \ \ \left\{
1,2\right\}  .
\]

A set cannot be red and green at the same time. In other words, the sets%
\[
\left\{  \text{red sets}\right\}  \ \ \ \ \ \ \ \ \ \ \text{and}%
\ \ \ \ \ \ \ \ \ \ \left\{  \text{green sets}\right\}
\]
are disjoint\footnote{Keep in mind: The notation \textquotedblleft$\left\{
\text{red sets}\right\}  $\textquotedblright\ stands for \textbf{the set of}
all red sets. For example, if $n=3$, then%
\[
\left\{  \text{red sets}\right\}  =\left\{  \left\{  3\right\}  ,\ \ \left\{
1,3\right\}  ,\ \ \left\{  2,3\right\}  ,\ \ \left\{  1,2,3\right\}  \right\}
.
\]
}. Hence, the sum rule for two sets (Theorem \ref{thm.count.rules.sum-2},
applied to $A=\left\{  \text{red sets}\right\}  $ and $B=\left\{  \text{green
sets}\right\}  $) yields%
\[
\left\vert \left\{  \text{red sets}\right\}  \cup\left\{  \text{green
sets}\right\}  \right\vert =\left\vert \left\{  \text{red sets}\right\}
\right\vert +\left\vert \left\{  \text{green sets}\right\}  \right\vert .
\]
(This is just a formal way to say \textquotedblleft the \# of all sets that
are red or green equals the \# of red sets plus the \# of green
sets\textquotedblright. Indeed, the notation $\left\{  \text{red
sets}\right\}  $ means the \textbf{set} of all red sets, and thus the
expression $\left\vert \left\{  \text{red sets}\right\}  \right\vert $ means
the \textbf{size} of the set of all red sets, i.e., the \# of all red sets.)

Furthermore, each subset of $\left[  n\right]  $ is either red or green (and
conversely, each red or green set is a subset of $\left[  n\right]  $). Hence,%
\[
\left\{  \text{subsets of }\left[  n\right]  \right\}  =\left\{  \text{red
sets}\right\}  \cup\left\{  \text{green sets}\right\}  .
\]
Therefore,%
\begin{align*}
\left\vert \left\{  \text{subsets of }\left[  n\right]  \right\}  \right\vert
&  =\left\vert \left\{  \text{red sets}\right\}  \cup\left\{  \text{green
sets}\right\}  \right\vert \\
&  =\left\vert \left\{  \text{red sets}\right\}  \right\vert +\left\vert
\left\{  \text{green sets}\right\}  \right\vert
\end{align*}
(as we have proved above). In other words,%
\begin{equation}
\left(  \text{\# of subsets of }\left[  n\right]  \right)  =\left(  \text{\#
of red sets}\right)  +\left(  \text{\# of green sets}\right)  .
\label{pf.thm.count.subsets.all.n-repeat.sum}%
\end{equation}

Thus it remains to count the red sets and the green sets separately.

The green sets are easy: They are just the subsets of $\left[  n-1\right]  $.
Hence,%
\[
\left(  \text{\# of green sets}\right)  =\left(  \text{\# of subsets of
}\left[  n-1\right]  \right)  =2^{n-1}%
\]
(by the induction hypothesis).

Counting the red sets is trickier. In Section \ref{sec.enum.subsets}, we did
this by setting up a one-to-one correspondence between the red sets and the
green sets. Formally, a one-to-one correspondence is just a bijection. Thus,
our one-to-one correspondence should become a bijection from $\left\{
\text{green sets}\right\}  $ to $\left\{  \text{red sets}\right\}  $ (i.e.,
from the set of all green sets to the set of all red sets).

As we recall, we obtained this correspondence as follows: To turn a green set
red, we insert $n$ into it; conversely, to turn a red set green, we remove $n$
from it. Rigorously, this means that we define two maps%
\begin{align*}
\operatorname*{ins}\nolimits_{n}:\left\{  \text{green sets}\right\}   &
\rightarrow\left\{  \text{red sets}\right\}  ,\\
G  &  \mapsto G\cup\left\{  n\right\}
\end{align*}
and%
\begin{align*}
\operatorname*{rem}\nolimits_{n}:\left\{  \text{red sets}\right\}   &
\rightarrow\left\{  \text{green sets}\right\}  ,\\
R  &  \mapsto R\setminus\left\{  n\right\}  .
\end{align*}
It is easy to see that both of these maps $\operatorname*{ins}\nolimits_{n}$
and $\operatorname*{rem}\nolimits_{n}$ are well-defined\footnote{Indeed, we
need to show that
\par
\begin{itemize}
\item if $G$ is a green set, then $G\cup\left\{  n\right\}  $ is a red set;
\par
\item if $R$ is a red set, then $R\setminus\left\{  n\right\}  $ is a green
set.
\end{itemize}
\par
Both of these claims are very easy. For instance, if $G$ is a green set, then
$G$ is a subset of $\left[  n\right]  $, and thus $G\cup\left\{  n\right\}  $
is a subset of $\left[  n\right]  $ as well (since $n\in\left[  n\right]  $),
and furthermore is red (since $n\in\left\{  n\right\}  \subseteq G\cup\left\{
n\right\}  $).}. A little bit of set-theoretic computation shows that%
\[
\operatorname*{ins}\nolimits_{n}\left(  \operatorname*{rem}\nolimits_{n}%
\left(  R\right)  \right)  =R\ \ \ \ \ \ \ \ \ \ \text{for every red set }R
\]
(because if $R$ is a red set, then%
\begin{align*}
\operatorname*{ins}\nolimits_{n}\left(  \operatorname*{rem}\nolimits_{n}%
\left(  R\right)  \right)   &  =\underbrace{\operatorname*{rem}\nolimits_{n}%
\left(  R\right)  }_{\substack{=R\setminus\left\{  n\right\}  \\\text{(by the
definition of }\operatorname*{rem}\nolimits_{n}\text{)}}}\cup\left\{
n\right\}  \ \ \ \ \ \ \ \ \ \ \left(  \text{by the definition of
}\operatorname*{ins}\nolimits_{n}\right) \\
&  =\left(  R\setminus\left\{  n\right\}  \right)  \cup\left\{  n\right\}
=R\ \ \ \ \ \ \ \ \ \ \left(  \text{since }n\in R\text{ (because }R\text{ is
red)}\right)
\end{align*}
). Similarly,%
\[
\operatorname*{rem}\nolimits_{n}\left(  \operatorname*{ins}\nolimits_{n}%
\left(  G\right)  \right)  =G\ \ \ \ \ \ \ \ \ \ \text{for every green set
}G.
\]
These two equalities show that the map $\operatorname*{rem}\nolimits_{n}$ is
an inverse of $\operatorname*{ins}\nolimits_{n}$. Hence, the map
$\operatorname*{ins}\nolimits_{n}$ has an inverse, i.e., is bijective (by
Theorem \ref{thm.maps.jects.inv=bij}). In other words, $\operatorname*{ins}%
\nolimits_{n}$ is a bijection. Hence, there exists a bijection from $\left\{
\text{green sets}\right\}  $ to $\left\{  \text{red sets}\right\}  $ (namely,
$\operatorname*{ins}\nolimits_{n}$). Thus, the bijection principle yields%
\[
\left\vert \left\{  \text{green sets}\right\}  \right\vert =\left\vert
\left\{  \text{red sets}\right\}  \right\vert .
\]
In other words,%
\[
\left(  \text{\# of green sets}\right)  =\left(  \text{\# of red sets}\right)
,
\]
and thus%
\[
\left(  \text{\# of red sets}\right)  =\left(  \text{\# of green sets}\right)
=2^{n-1}.
\]

Combining what we have shown, we now obtain
\begin{align*}
\left(  \text{\# of subsets of }\left[  n\right]  \right)   &
=\underbrace{\left(  \text{\# of red sets}\right)  }_{=2^{n-1}}%
+\underbrace{\left(  \text{\# of green sets}\right)  }_{=2^{n-1}}\\
&  =2^{n-1}+2^{n-1}=2\cdot2^{n-1}=2^{n}.
\end{align*}
This is precisely what we needed to prove. This completes the induction step,
and thus Theorem \ref{thm.count.subsets.all.n-repeat} is proved.
\end{proof}

Theorem \ref{thm.count.subsets.all.S} says the following:

\begin{theorem}
\label{thm.count.subsets.all.S-repeat}Let $n\in\mathbb{N}$. Let $S$ be an
$n$-element set. Then,%
\[
\left(  \text{\# of subsets of }S\right)  =2^{n}.
\]

\end{theorem}

\begin{proof}
[Rigorous proof.]Informally, we derived this from Theorem
\ref{thm.count.subsets.all.n-repeat} by renaming the elements of $S$ as
$1,2,\ldots,n$ (so that $S$ became the set $\left[  n\right]  $).

Rigorously, this means setting up a one-to-one correspondence between the
subsets of $S$ and the subsets of $\left[  n\right]  $, and then using the
bijection principle to argue that the \# of the former equals the \# of the latter.

How do we get this correspondence? First, we set up a one-to-one
correspondence between the \textbf{elements} of $S$ and the elements of
$\left[  n\right]  $. (This is what the \textquotedblleft
renaming\textquotedblright\ in our informal proof was secretly doing.)
Formally, this can be done as follows:

The set $S$ is an $n$-element set, i.e., has size $n$. Hence, by the
definition of size, the set $S$ is isomorphic to $\left[  n\right]  $. In
other words, there is a bijection $\alpha:S\rightarrow\left[  n\right]  $.
Consider this $\alpha$. Being a bijection, the map $\alpha$ has an inverse
$\alpha^{-1}$ (by Theorem \ref{thm.maps.jects.inv=bij}).

Now, define a map%
\begin{align*}
\alpha_{\ast}:\left\{  \text{subsets of }S\right\}   &  \rightarrow\left\{
\text{subsets of }\left[  n\right]  \right\}  ,\\
T  &  \mapsto\left\{  \alpha\left(  t\right)  \ \mid\ t\in T\right\}  .
\end{align*}
Explicitly, this map $\alpha_{\ast}$ sends every subset $\left\{  s_{1}%
,s_{2},\ldots,s_{k}\right\}  $ of $S$ to the subset $\left\{  \alpha\left(
s_{1}\right)  ,\alpha\left(  s_{2}\right)  ,\ldots,\alpha\left(  s_{k}\right)
\right\}  $ of $\left[  n\right]  $; that is, it applies $\alpha$ to every
element of the input subset. (For example, if $n=3$ and $S=\left\{
\text{\textquotedblleft cat\textquotedblright},\ \text{\textquotedblleft
dog\textquotedblright},\text{ \textquotedblleft rat\textquotedblright%
}\right\}  $ and if $\alpha\left(  \text{\textquotedblleft
cat\textquotedblright}\right)  =1$ and $\alpha\left(  \text{\textquotedblleft
dog\textquotedblright}\right)  =2$ and $\alpha\left(  \text{\textquotedblleft
rat\textquotedblright}\right)  =3$, then $\alpha_{\ast}\left(  \left\{
\text{\textquotedblleft cat\textquotedblright},\ \text{\textquotedblleft
rat\textquotedblright}\right\}  \right)  =\left\{  1,3\right\}  $.)

Conversely, we can define a map%
\begin{align*}
\left(  \alpha^{-1}\right)  _{\ast}:\left\{  \text{subsets of }\left[
n\right]  \right\}   &  \rightarrow\left\{  \text{subsets of }S\right\}  ,\\
T  &  \mapsto\left\{  \alpha^{-1}\left(  t\right)  \ \mid\ t\in T\right\}  .
\end{align*}
(This map $\left(  \alpha^{-1}\right)  _{\ast}$ is defined in the same way as
$\alpha_{\ast}$, but using the map $\alpha^{-1}$ instead of $\alpha$. For
example, if $n=3$ and $S=\left\{  \text{\textquotedblleft
cat\textquotedblright},\ \text{\textquotedblleft dog\textquotedblright},\text{
\textquotedblleft rat\textquotedblright}\right\}  $ and if $\alpha\left(
\text{\textquotedblleft cat\textquotedblright}\right)  =1$ and $\alpha\left(
\text{\textquotedblleft dog\textquotedblright}\right)  =2$ and $\alpha\left(
\text{\textquotedblleft rat\textquotedblright}\right)  =3$, then $\left(
\alpha^{-1}\right)  _{\ast}\left(  \left\{  2,3\right\}  \right)  =\left\{
\text{\textquotedblleft dog\textquotedblright},\text{ \textquotedblleft
rat\textquotedblright}\right\}  $.)

It is easy to see that the map $\left(  \alpha^{-1}\right)  _{\ast}$ is an
inverse of $\alpha_{\ast}$ (because applying $\alpha$ to each element of a
given set and then applying $\alpha^{-1}$ to the results will recover the
original set, and likewise if you apply $\alpha^{-1}$ first and then $\alpha
$). Thus, the map $\alpha_{\ast}$ has an inverse, i.e., is a bijection (by
Theorem \ref{thm.maps.jects.inv=bij}). Thus, we have found a bijection from
$\left\{  \text{subsets of }S\right\}  $ to $\left\{  \text{subsets of
}\left[  n\right]  \right\}  $ (namely, $\alpha_{\ast}$). Hence, the bijection
principle (Theorem \ref{thm.count.rules.bij}) yields
\[
\left\vert \left\{  \text{subsets of }S\right\}  \right\vert =\left\vert
\left\{  \text{subsets of }\left[  n\right]  \right\}  \right\vert .
\]
In other words,%
\[
\left(  \text{\# of subsets of }S\right)  =\left(  \text{\# of subsets of
}\left[  n\right]  \right)  =2^{n}%
\]
(by Theorem \ref{thm.count.subsets.all.n-repeat}).
\end{proof}

\subsubsection{Counting all $k$-element subsets}

We move on to counting subsets of a given size.

Theorem \ref{thm.count.subsets.size.S} says:

\begin{theorem}
\label{thm.count.subsets.size.S-repeat}Let $n\in\mathbb{N}$, and let $k$ be
any number (not necessarily an integer). Let $S$ be an $n$-element set. Then,%
\[
\left(  \text{\# of }k\text{-element subsets of }S\right)  =\dbinom{n}{k}.
\]

\end{theorem}

\begin{proof}
[Rigorous proof.]We induct on $n$ (without fixing $k$). That is, we use
induction on $n$ to prove the statement%
\[
P\left(  n\right)  :=\left(
\begin{array}
[c]{c}%
\text{\textquotedblleft for any number }k\text{ and any }n\text{-element set
}S,\\
\text{we have }\left(  \text{\# of }k\text{-element subsets of }S\right)
=\dbinom{n}{k}\text{\textquotedblright}%
\end{array}
\right)
\]
for each $n\in\mathbb{N}$.

\textit{Base case:} Let $k$ be any number. The only $0$-element set is
$\varnothing$, and its only subset is $\varnothing$. Thus, a $0$-element set
$S$ necessarily has one $0$-element subset ($\varnothing$) and no other
subsets. Hence, it satisfies%
\[
\left(  \text{\# of }k\text{-element subsets of }S\right)  =%
\begin{cases}
1, & \text{if }k=0;\\
0, & \text{else.}%
\end{cases}
\]
However, we also have%
\[
\dbinom{0}{k}=%
\begin{cases}
1, & \text{if }k=0;\\
0, & \text{else}%
\end{cases}
\]
(this follows easily from the definition of binomial coefficients). By
comparing these two equalities, we see that any $0$-element set $S$ satisfies
\[
\left(  \text{\# of }k\text{-element subsets of }S\right)  =\dbinom{0}{k}.
\]
In other words, $P\left(  0\right)  $ holds.

\textit{Induction step:} Let $n$ be a positive integer. Assume (as the
induction hypothesis) that $P\left(  n-1\right)  $ holds. We must prove that
$P\left(  n\right)  $ holds.

So we consider any number $k$ and any $n$-element set $S$. We must prove that%
\[
\left(  \text{\# of }k\text{-element subsets of }S\right)  =\dbinom{n}{k}.
\]

We rename the $n$ elements of $S$ as $1,2,\ldots,n$ (this corresponds formally
to constructing a bijection $\alpha:S\rightarrow\left[  n\right]  $ and
applying it elementwise to subsets of $S$, as we did in the proof of Theorem
\ref{thm.count.subsets.all.S-repeat}), so we must prove that%
\[
\left(  \text{\# of }k\text{-element subsets of }\left[  n\right]  \right)
=\dbinom{n}{k}.
\]

To prove this, we define

\begin{itemize}
\item a \textbf{red set} to be a $k$-element subset of $\left[  n\right]  $
that contains $n$;

\item a \textbf{green set} to be a $k$-element subset of $\left[  n\right]  $
that does not contain $n$.
\end{itemize}

For instance, for $n=4$ and $k=2$, the red sets are%
\[
\left\{  1,4\right\}  ,\ \ \ \ \ \ \ \ \ \ \left\{  2,4\right\}
,\ \ \ \ \ \ \ \ \ \ \left\{  3,4\right\}  ,
\]
while the green sets are%
\[
\left\{  1,2\right\}  ,\ \ \ \ \ \ \ \ \ \ \left\{  1,3\right\}
,\ \ \ \ \ \ \ \ \ \ \left\{  2,3\right\}  .
\]

Each $k$-element subset of $\left[  n\right]  $ is either red or green (but
not both). Hence, using the sum rule for two sets, we find%
\begin{align}
&  \left(  \text{\# of }k\text{-element subsets of }\left[  n\right]  \right)
\nonumber\\
&  =\left(  \text{\# of red sets}\right)  +\left(  \text{\# of green
sets}\right)  . \label{pf.thm.count.subsets.size.S-repeat.sum}%
\end{align}
(This is proved just as we proved (\ref{pf.thm.count.subsets.all.n-repeat.sum}%
) in the rigorous proof of Theorem \ref{thm.count.subsets.all.n-repeat}.)

The green sets are just the $k$-element subsets of $\left[  n-1\right]  $.
Thus,%
\begin{align*}
\left(  \text{\# of green sets}\right)   &  =\left(  \text{\# of
}k\text{-element subsets of }\left[  n-1\right]  \right) \\
&  =\dbinom{n-1}{k}%
\end{align*}
(by the statement $P\left(  n-1\right)  $, which we have assumed to hold).

Now, let's try to count the red sets.

Let us refer to the $\left(  k-1\right)  $-element subsets of $\left[
n-1\right]  $ as \textbf{blue sets}. If $R$ is a red set, then $R\setminus
\left\{  n\right\}  $ is a blue set\footnote{\textit{Proof.} Let $R$ be a red
set. Then, $R$ is a $k$-element set (by the definition of a red set), so that
$\left\vert R\right\vert =k$. Moreover, $n\in R$ (by the definition of a red
set), so that $\left\{  n\right\}  \subseteq R$. Hence, the difference rule
(Theorem \ref{thm.count.rules.diff} \textbf{(b)}, applied to $S=R$ and
$T=\left\{  n\right\}  $) yields $\left\vert R\setminus\left\{  n\right\}
\right\vert =\underbrace{\left\vert R\right\vert }_{=k}-\underbrace{\left\vert
\left\{  n\right\}  \right\vert }_{=1}=k-1$. Hence, $R\setminus\left\{
n\right\}  $ is a $\left(  k-1\right)  $-element set. Since $R\setminus
\left\{  n\right\}  $ is furthermore a subset of $\left[  n-1\right]  $
(because $R$ is a subset of $\left[  n\right]  $, and we are removing $n$ from
it), we thus conclude that $R\setminus\left\{  n\right\}  $ is a $\left(
k-1\right)  $-element subset of $\left[  n-1\right]  $, that is, a blue set.}.
Thus, we obtain a map%
\begin{align*}
\operatorname*{rem}\nolimits_{n}:\left\{  \text{red sets}\right\}   &
\rightarrow\left\{  \text{blue sets}\right\}  ,\\
R  &  \mapsto R\setminus\left\{  n\right\}  .
\end{align*}
Conversely, if $B$ is a blue set, then $B\cup\left\{  n\right\}  $ is a red
set\footnote{\textit{Proof.} Let $B$ be a blue set. Then, $B$ is a $\left(
k-1\right)  $-element subset of $\left[  n-1\right]  $ (by the definition of
\textquotedblleft blue set\textquotedblright). In other words, $\left\vert
B\right\vert =k-1$ and $B\subseteq\left[  n-1\right]  $. From $B\subseteq
\left[  n-1\right]  $, we obtain $n\notin B$ (since $n\notin\left[
n-1\right]  $). Hence, Theorem \ref{thm.count.rules.+1} (applied to $S=B$ and
$t=n$) yields $\left\vert B\cup\left\{  n\right\}  \right\vert =\left\vert
B\right\vert +1=k$ (since $\left\vert B\right\vert =k-1$). Hence,
$B\cup\left\{  n\right\}  $ is a $k$-element set. Furthermore, $B\cup\left\{
n\right\}  $ is a subset of $\left[  n\right]  $ (since $B\subseteq\left[
n-1\right]  \subseteq\left[  n\right]  $ and $\left\{  n\right\}
\subseteq\left[  n\right]  $) that contains $n$ (since $n\in\left\{
n\right\}  \subseteq B\cup\left\{  n\right\}  $). Thus, $B\cup\left\{
n\right\}  $ is a $k$-element subset of $\left[  n\right]  $ that contains
$n$. In other words, $B\cup\left\{  n\right\}  $ is a red set (by the
definition of \textquotedblleft red set\textquotedblright).}. Thus, we obtain
a map%
\begin{align*}
\operatorname*{ins}\nolimits_{n}:\left\{  \text{blue sets}\right\}   &
\rightarrow\left\{  \text{red sets}\right\}  ,\\
B  &  \mapsto B\cup\left\{  n\right\}  .
\end{align*}
These two maps $\operatorname*{rem}\nolimits_{n}$ and $\operatorname*{ins}%
\nolimits_{n}$ are mutually inverse\footnote{This can be proved just as in our
above proof of Theorem \ref{thm.count.subsets.all.n-repeat}.}. Thus, the map
$\operatorname*{rem}\nolimits_{n}$ has an inverse, i.e., is bijective (by
Theorem \ref{thm.maps.jects.inv=bij}). Hence, we have found a bijection from
$\left\{  \text{red sets}\right\}  $ to $\left\{  \text{blue sets}\right\}  $
(namely, $\operatorname*{rem}\nolimits_{n}$). The bijection principle
therefore yields%
\[
\left\vert \left\{  \text{red sets}\right\}  \right\vert =\left\vert \left\{
\text{blue sets}\right\}  \right\vert .
\]
In other words,%
\begin{align*}
\left(  \text{\# of red sets}\right)   &  =\left(  \text{\# of blue
sets}\right) \\
&  =\left(  \text{\# of }\left(  k-1\right)  \text{-element subsets of
}\left[  n-1\right]  \right) \\
&  \ \ \ \ \ \ \ \ \ \ \ \ \ \ \ \ \ \ \ \ \left(  \text{since this is how the
blue sets were defined}\right) \\
&  =\dbinom{n-1}{k-1}%
\end{align*}
(again by the statement $P\left(  n-1\right)  $, but now applied to $k-1$
instead of $k$). Note that we deliberately formulated $P\left(  n\right)  $ as
a \textquotedblleft for any $k$\textquotedblright\ statement (rather than
fixing $k$ at the onset of our proof), so that we were now able to apply
$P\left(  n-1\right)  $ to $k-1$ instead of $k$.

Now, (\ref{pf.thm.count.subsets.size.S-repeat.sum}) becomes%
\begin{align*}
\left(  \text{\# of }k\text{-element subsets of }\left[  n\right]  \right)
&  =\underbrace{\left(  \text{\# of red sets}\right)  }_{=\dbinom{n-1}{k-1}%
}+\underbrace{\left(  \text{\# of green sets}\right)  }_{=\dbinom{n-1}{k}}\\
&  =\dbinom{n-1}{k-1}+\dbinom{n-1}{k}=\dbinom{n}{k}%
\end{align*}
by Pascal's recurrence (Theorem \ref{thm.binom.rec}). But this is precisely
the equality that we have to prove. This completes the induction step, and
thus Theorem \ref{thm.count.subsets.size.S-repeat} is proved.
\end{proof}

\begin{remark}
Our above proof of Theorem \ref{thm.count.subsets.size.S-repeat} can be
simplified: There is no need to \textquotedblleft rename\textquotedblright%
\ the elements of $S$ as $1,2,\ldots,n$ in the induction step. Instead, we
could have just as well picked an arbitrary element $t$ of $S$ (such an
element exists, since $\left\vert S\right\vert =n>0$ entails that $S$ is
nonempty) and defined

\begin{itemize}
\item a \textbf{red set} to be a $k$-element subset of $S$ that contains $t$;

\item a \textbf{green set} to be a $k$-element subset of $S$ that does not
contain $t$.
\end{itemize}

Then, a simple application of Theorem \ref{thm.count.rules.diff} \textbf{(b)}
would have shown that $S\setminus\left\{  t\right\}  $ is an $\left(
n-1\right)  $-element set, so we could apply our induction hypothesis
$P\left(  n-1\right)  $ to it. Thus, the above argument could be made using
$S$, $t$ and $S\setminus\left\{  t\right\}  $ instead of $\left[  n\right]  $,
$n$ and $\left[  n-1\right]  $. In particular, the green sets would be
precisely the $k$-element subsets of $S\setminus\left\{  t\right\}  $, whereas
the red sets would be in one-to-one correspondence (i.e., bijection) with the
$\left(  k-1\right)  $-element subsets of $S\setminus\left\{  t\right\}  $
(and the bijection would be given by removing/inserting $t$). This argument
would be not only shorter but also more conceptual than the one we gave above.

However, I chose to give the proof I gave because it has the advantage of
familiarity (the set $\left[  n\right]  =\left\{  1,2,\ldots,n\right\}  $ is
easier to visualize than an arbitrary $n$-element set), and in order to
illustrate how the bijection principle can be used to rename the elements of a
given set in a convenient way.

Likewise, Theorem \ref{thm.count.subsets.all.S-repeat} could also be proved
more directly: Instead of deducing it from Theorem
\ref{thm.count.subsets.all.n-repeat} via \textquotedblleft
renaming\textquotedblright, we could have proved it by induction, again
picking an element $t$ of $S$ in the induction step, defining red and green
sets, and counting both kinds of sets using the induction hypothesis (applied
to the $\left(  n-1\right)  $-element set $S\setminus\left\{  t\right\}  $).
\end{remark}

Let us derive a nice, if simple, corollary from our last few theorems:

\begin{corollary}
\label{cor.count.subsets.2n=}Let $n\in\mathbb{N}$. Then,%
\[
\sum_{k=0}^{n}\dbinom{n}{k}=2^{n}.
\]

\end{corollary}

\begin{proof}
Consider the $n$-element set $\left[  n\right]  =\left\{  1,2,\ldots
,n\right\}  $. This set has size $n$, so each subset of $\left[  n\right]  $
must have size $\leq n$ (by Theorem \ref{thm.count.rules.diff} \textbf{(a)}).
Hence, each subset of $\left[  n\right]  $ has size $0$ or size $1$ or size
$2$ or $\cdots$ or size $n$. Thus, we can write the set%
\[
\left\{  \text{subsets of }\left[  n\right]  \right\}
\]
as a union%
\begin{align*}
&  \left\{  0\text{-element subsets of }\left[  n\right]  \right\} \\
&  \ \ \ \ \ \ \ \ \ \ \cup\left\{  1\text{-element subsets of }\left[
n\right]  \right\} \\
&  \ \ \ \ \ \ \ \ \ \ \cup\left\{  2\text{-element subsets of }\left[
n\right]  \right\} \\
&  \ \ \ \ \ \ \ \ \ \ \cup\cdots\\
&  \ \ \ \ \ \ \ \ \ \ \cup\left\{  n\text{-element subsets of }\left[
n\right]  \right\}  .
\end{align*}
Furthermore, this union is a union of disjoint sets (since a subset of
$\left[  n\right]  $ cannot have several distinct sizes at once). Therefore,
the sum rule for $k$ sets (Theorem \ref{thm.count.rules.sum-k}) yields%
\begin{align*}
&  \left\vert \left\{  \text{subsets of }\left[  n\right]  \right\}
\right\vert \\
&  =\left\vert \left\{  0\text{-element subsets of }\left[  n\right]
\right\}  \right\vert \\
&  \ \ \ \ \ \ \ \ \ \ +\left\vert \left\{  1\text{-element subsets of
}\left[  n\right]  \right\}  \right\vert \\
&  \ \ \ \ \ \ \ \ \ \ +\left\vert \left\{  2\text{-element subsets of
}\left[  n\right]  \right\}  \right\vert \\
&  \ \ \ \ \ \ \ \ \ \ +\cdots\\
&  \ \ \ \ \ \ \ \ \ \ +\left\vert \left\{  n\text{-element subsets of
}\left[  n\right]  \right\}  \right\vert .
\end{align*}
\footnote{In more details:
\par
The $n+1$ sets%
\begin{align*}
&  \left\{  0\text{-element subsets of }\left[  n\right]  \right\}  ,\\
&  \left\{  1\text{-element subsets of }\left[  n\right]  \right\}  ,\\
&  \left\{  2\text{-element subsets of }\left[  n\right]  \right\}  ,\\
&  \ldots,\\
&  \left\{  n\text{-element subsets of }\left[  n\right]  \right\}
\end{align*}
are finite (since $\left[  n\right]  $ has only finitely many subsets) and
disjoint (since a subset of $\left[  n\right]  $ cannot have several distinct
sizes at once). Thus, the sum rule for $k$ sets (Theorem
\ref{thm.count.rules.sum-k}) yields%
\begin{align*}
&  \Big|\ \left\{  0\text{-element subsets of }\left[  n\right]  \right\}
\cup\left\{  1\text{-element subsets of }\left[  n\right]  \right\} \\
&  \ \ \ \ \ \ \ \ \ \ \cup\left\{  2\text{-element subsets of }\left[
n\right]  \right\}  \cup\cdots\cup\left\{  n\text{-element subsets of }\left[
n\right]  \right\}  \Big|\\
&  =\left\vert \left\{  0\text{-element subsets of }\left[  n\right]
\right\}  \right\vert \\
&  \ \ \ \ \ \ \ \ \ \ +\left\vert \left\{  1\text{-element subsets of
}\left[  n\right]  \right\}  \right\vert \\
&  \ \ \ \ \ \ \ \ \ \ +\left\vert \left\{  2\text{-element subsets of
}\left[  n\right]  \right\}  \right\vert \\
&  \ \ \ \ \ \ \ \ \ \ +\cdots\\
&  \ \ \ \ \ \ \ \ \ \ +\left\vert \left\{  n\text{-element subsets of
}\left[  n\right]  \right\}  \right\vert .
\end{align*}
\par
Since
\[
\left\{  \text{subsets of }\left[  n\right]  \right\}
\]
is the union%
\begin{align*}
&  \left\{  0\text{-element subsets of }\left[  n\right]  \right\} \\
&  \ \ \ \ \ \ \ \ \ \ \cup\left\{  1\text{-element subsets of }\left[
n\right]  \right\} \\
&  \ \ \ \ \ \ \ \ \ \ \cup\left\{  2\text{-element subsets of }\left[
n\right]  \right\} \\
&  \ \ \ \ \ \ \ \ \ \ \cup\cdots\\
&  \ \ \ \ \ \ \ \ \ \ \cup\left\{  n\text{-element subsets of }\left[
n\right]  \right\}  ,
\end{align*}
we can rewrite this equality as
\begin{align*}
&  \left\vert \left\{  \text{subsets of }\left[  n\right]  \right\}
\right\vert \\
&  =\left\vert \left\{  0\text{-element subsets of }\left[  n\right]
\right\}  \right\vert \\
&  \ \ \ \ \ \ \ \ \ \ +\left\vert \left\{  1\text{-element subsets of
}\left[  n\right]  \right\}  \right\vert \\
&  \ \ \ \ \ \ \ \ \ \ +\left\vert \left\{  2\text{-element subsets of
}\left[  n\right]  \right\}  \right\vert \\
&  \ \ \ \ \ \ \ \ \ \ +\cdots\\
&  \ \ \ \ \ \ \ \ \ \ +\left\vert \left\{  n\text{-element subsets of
}\left[  n\right]  \right\}  \right\vert .
\end{align*}
}

In other words,%
\begin{align*}
&  \left(  \text{\# of subsets of }\left[  n\right]  \right) \\
&  =\left(  \text{\# of }0\text{-element subsets of }\left[  n\right]  \right)
\\
&  \ \ \ \ \ \ \ \ \ \ +\left(  \text{\# of }1\text{-element subsets of
}\left[  n\right]  \right) \\
&  \ \ \ \ \ \ \ \ \ \ +\left(  \text{\# of }2\text{-element subsets of
}\left[  n\right]  \right) \\
&  \ \ \ \ \ \ \ \ \ \ +\cdots\\
&  \ \ \ \ \ \ \ \ \ \ +\left(  \text{\# of }n\text{-element subsets of
}\left[  n\right]  \right) \\
&  =\sum_{k=0}^{n}\underbrace{\left(  \text{\# of }k\text{-element subsets of
}\left[  n\right]  \right)  }_{\substack{=\dbinom{n}{k}\\\text{(by Theorem
\ref{thm.count.subsets.size.S-repeat}, applied to }S=\left[  n\right]
\text{)}}}\\
&  =\sum_{k=0}^{n}\dbinom{n}{k}.
\end{align*}
Thus,%
\[
\sum_{k=0}^{n}\dbinom{n}{k}=\left(  \text{\# of subsets of }\left[  n\right]
\right)  =2^{n}%
\]
(by Theorem \ref{thm.count.subsets.all.n-repeat}).
\end{proof}

Corollary \ref{cor.count.subsets.2n=} can also be easily obtained from the
binomial formula (this was part of Exercise \ref{exe.binom.2^n} \textbf{(a)}).
Our above proof, however, reveals its combinatorial meaning: It comes from the
comparison of two different ways to count one and the same thing (viz., the
subsets of $\left[  n\right]  $). This technique of proving equalities is
called \emph{double counting}, and has multiple other applications (see, e.g.,
\cite[\S 8.1]{Newste23}).

\subsubsection{Recounting pairs}

Proposition \ref{prop.count.pairs.nn} says:

\begin{proposition}
\label{prop.count.pairs.nn-repeat}Let $n\in\mathbb{N}$. Then: \medskip

\textbf{(a)} The \# of pairs $\left(  a,b\right)  $ with $a,b\in\left[
n\right]  $ is $n^{2}$. \medskip

\textbf{(b)} The \# of pairs $\left(  a,b\right)  $ with $a,b\in\left[
n\right]  $ and $a<b$ is $1+2+\cdots+\left(  n-1\right)  $. \medskip

\textbf{(c)} The \# of pairs $\left(  a,b\right)  $ with $a,b\in\left[
n\right]  $ and $a=b$ is $n$. \medskip

\textbf{(d)} The \# of pairs $\left(  a,b\right)  $ with $a,b\in\left[
n\right]  $ and $a>b$ is $1+2+\cdots+\left(  n-1\right)  $.
\end{proposition}

Let us reprove part \textbf{(b)} of this proposition rigorously:

\begin{proof}
[Rigorous proof of Proposition \ref{prop.count.pairs.nn-repeat} \textbf{(b)}
(sketched).]If $\left(  a,b\right)  $ is a pair with $a,b\in\left[  n\right]
$ and $a<b$, then the first entry of this pair (that is, the number $a$) must
be one of the numbers $1,2,\ldots,n-1$ (because $a<b\leq n$ forces $a$ to be
$\leq n-1$). Thus, by the sum rule for $k$ sets (Theorem
\ref{thm.count.rules.sum-k}), we have%
\begin{align*}
&  \left(  \text{\# of pairs }\left(  a,b\right)  \text{ with }a,b\in\left[
n\right]  \text{ and }a<b\right) \\
&  =\sum_{k=1}^{n-1}\underbrace{\left(  \text{\# of pairs }\left(  a,b\right)
\text{ with }a,b\in\left[  n\right]  \text{ and }a<b\text{ and }a=k\right)
}_{\substack{=n-k\\\text{(because these pairs are }\left(  k,\ k+1\right)
,\ \left(  k,\ k+2\right)  ,\ \ldots,\ \left(  k,n\right)  \\\text{(strictly
speaking, this argument is an application of}\\\text{the bijection
principle))}}}\\
&  =\sum_{k=1}^{n-1}\left(  n-k\right)  =\left(  n-1\right)  +\left(
n-2\right)  +\cdots+\left(  n-\left(  n-1\right)  \right) \\
&  =\left(  n-1\right)  +\left(  n-2\right)  +\cdots+1\\
&  =1+2+\cdots+\left(  n-1\right)  ,
\end{align*}
and thus Proposition \ref{prop.count.pairs.nn-repeat} \textbf{(b)} is proven.
\end{proof}

\begin{exercise}
Let $n\in\mathbb{N}$. Consider the set $\left[  2n\right]  =\left\{
1,2,\ldots,2n\right\}  $.

A set of integers will be called \textbf{parity-ambivalent} if it contains at
least one even element and at least one odd element. (For instance, $\left\{
2,4,5\right\}  $ is parity-ambivalent, but $\left\{  2,4,10\right\}  $ is not.)

Compute the \# of all parity-ambivalent subsets of $\left[  2n\right]  $.
\medskip

[\textbf{Hint:} How many subsets of $\left[  2n\right]  $ contain \textbf{no}
even element? How many contain \textbf{no} odd element? How many contain neither?]
\end{exercise}

\begin{exercise}
\label{exe.23hw5.4}Let $n\in\mathbb{N}$. Compute the \# of pairs $\left(
A,B\right)  $ of subsets of $\left[  n\right]  $ that satisfy $A\cap
B=\varnothing$.

(For example, if $n=2$, then this \# is $9$, since there are $9$ such pairs:%
\begin{align*}
&  \left(  \varnothing,\ \varnothing\right)  ,\qquad\left(  \varnothing
,\ \left\{  1\right\}  \right)  ,\qquad\left(  \varnothing,\ \left\{
2\right\}  \right)  ,\qquad\left(  \varnothing,\ \left\{  1,2\right\}
\right)  ,\\
&  \left(  \left\{  1\right\}  ,\ \varnothing\right)  ,\qquad\left(  \left\{
1\right\}  ,\ \left\{  2\right\}  \right)  ,\qquad\left(  \left\{  2\right\}
,\ \varnothing\right)  ,\qquad\left(  \left\{  2\right\}  ,\ \left\{
1\right\}  \right)  ,\\
&  \left(  \left\{  1,2\right\}  ,\ \varnothing\right)  .
\end{align*}
)
\end{exercise}

\begin{exercise}
Let $n\in\mathbb{N}$. Compute the \# of pairs $\left(  A,B\right)  $ of
subsets of $\left[  n\right]  $ that satisfy $A\subseteq B$.

(For example, if $n=2$, then this \# is $9$, since there are $9$ such pairs:%
\begin{align*}
&  \left(  \varnothing,\ \varnothing\right)  ,\qquad\left(  \varnothing
,\ \left\{  1\right\}  \right)  ,\qquad\left(  \varnothing,\ \left\{
2\right\}  \right)  ,\qquad\left(  \varnothing,\ \left\{  1,2\right\}
\right)  ,\\
&  \left(  \left\{  1\right\}  ,\ \left\{  1\right\}  \right)  ,\qquad\left(
\left\{  1\right\}  ,\ \left\{  1,2\right\}  \right)  ,\qquad\left(  \left\{
2\right\}  ,\ \left\{  2\right\}  \right)  ,\qquad\left(  \left\{  2\right\}
,\ \left\{  1,2\right\}  \right)  ,\\
&  \left(  \left\{  1,2\right\}  ,\ \left\{  1,2\right\}  \right)  .
\end{align*}
)
\end{exercise}

\subsection{\label{sec.enum.intermed}Where do we stand now?}

Recall the introductory counting problems from the start of Chapter
\ref{chp.enum} (before Section \ref{sec.enum.count-informal}). We can now
answer some of these:

\begin{itemize}
\item How many ways are there to choose $3$ odd integers between $0$ and $20$,
if the order matters (i.e., we count the choice $1,3,5$ as different from the
choice $3,1,5$)? (The answer is $1000$.)

\textbf{We can solve this now:} To choose $3$ odd integers between $0$ and
$20$, if the order matters, amounts to choosing a $3$-tuple $\left(
a,b,c\right)  $ where $a,b,c\in\left\{  1,3,5,\ldots,19\right\}  $. Since this
set $\left\{  1,3,5,\ldots,19\right\}  $ is a $10$-element set (because
Proposition \ref{prop.count.odd-0-n} yields that the \# of odd integers
between $0$ and $20$ is $\left(  20+1\right)  //2=10$), the \# of these
$3$-tuples is $10\cdot10\cdot10=1000$ (by Theorem \ref{thm.count.triples.nmp}).

\item How many ways are there to choose $3$ odd integers between $0$ and $20$,
if the order does not matter? (The answer is $220$.)

\textbf{We cannot solve this yet}, at least not if the values $3$ and $20$ are
generalized to $k$ and $n$. This will be done in Theorem
\ref{thm.count.select.uotup-mulset}.

\item How many ways are there to choose $3$ distinct odd integers between $0$
and $20$, if the order matters? (The answer is $720$.)

\textbf{We cannot solve this yet}, at least not if the values $3$ and $20$ are
generalized to $k$ and $n$. This will be done in Theorem
\ref{thm.count.select.injtup}.

\item How many ways are there to choose $3$ distinct odd integers between $0$
and $20$, if the order does not matter? (The answer is $120$.)

\textbf{We can solve this now:} This amounts to counting the $3$-element
subsets of $\left\{  1,3,5,\ldots,19\right\}  $; but Theorem
\ref{thm.count.subsets.size.S-repeat} answers such questions. Since the set
$\left\{  1,3,5,\ldots,19\right\}  $ has size $10$, its number of $3$-element
subsets is $\dbinom{10}{3}=\dfrac{10\cdot9\cdot8}{3!}=120$.

\item How many prime factorizations does $200$ have (where we count different
orderings as distinct)? (The answer is $10$. This is a mix between a number
theory problem and a counting problem.)

\textbf{We can solve this now, at least for }$200$\textbf{:} We know that
$200=2\cdot2\cdot2\cdot5\cdot5$. Thus, by the fundamental theorem of
arithmetic, all prime factorizations of $200$ consist of five factors, three
of which are $2$'s and two of which are $5$'s. The only freedom is in choosing
where to place the three $2$'s among the five positions (of course, the two
$5$'s will then have to occupy the remaining positions). There are $5$ factors
in total, so $5$ positions, and we have to choose $3$ of these $5$ positions
to put our three $2$'s in. This is tantamount to choosing a $3$-element subset
of $\left[  5\right]  $ (the subset of the positions in which we put the
$2$'s), and the \# of ways to do this is $\dbinom{5}{3}=10$ (by Theorem
\ref{thm.count.subsets.size.S-repeat}). Thus, $200$ has $10$ prime
factorizations (if we count different orderings as distinct).

However, it is trickier to extend this reasoning to prime factorizations of
$150$. Indeed, $150=2\cdot3\cdot5\cdot5$, so a prime factorization of $150$
has one $2$, one $3$ and two $5$'s. How many ways are there to place one $2$,
one $3$ and two $5$'s in altogether four positions? I'll leave this one to you
for now, but we will come back to this later.

\item How many ways are there to tile a $2\times15$-rectangle with dominos
(i.e., rectangles of size $1\times2$ or $2\times1$) ? (The answer is $987$.)

\textbf{We cannot solve this yet.} But we will outline a solution in
Subsection \ref{subsec.enum.lac.domino}.

\item How many addends do you get when you expand the product \newline$\left(
a+b\right)  \left(  c+d+e\right)  \left(  f+g\right)  $ ? (The answer is $12$.)

\textbf{We can solve this now:} Each addend consists of exactly one of $a$ and
$b$, exactly one of $c$, $d$ and $e$, and exactly one of $f$ and $g$. So the
addends are in one-to-one correspondence with the triples $\left(
x,y,z\right)  $ where $x\in\left\{  a,b\right\}  $ and $y\in\left\{
c,d,e\right\}  $ and $z\in\left\{  f,g\right\}  $. Thus, their \# is
$2\cdot3\cdot2$ (since $\left\{  a,b\right\}  $ is a $2$-element set,
$\left\{  c,d,e\right\}  $ is a $3$-element set, and $\left\{  f,g\right\}  $
is a $2$-element set).

Note that we are using the fact that the addends all end up distinct, so they
don't cancel or combine.

\item How many different monomials do you get when you expand the product
$\left(  a-b\right)  \left(  a^{2}+ab+b^{2}\right)  $ ? (This one is more of
an algebra problem, but I wanted to list it because it is connected to
counting. The answer is $2$, because $\left(  a-b\right)  \left(
a^{2}+ab+b^{2}\right)  =a^{3}-b^{3}$.)

\textbf{This is not a combinatorics problem:} The answer is $2$, because we
have $\left(  a-b\right)  \left(  a^{2}+ab+b^{2}\right)  =a^{3}-b^{3}$. The
other addends all cancel out, so you get an answer much less than $6$.

In general, problems like this (where you count addends after cancellation and
combination) cannot be solved combinatorially; you have to actually expand and collect.

\item How many positive divisors does $24$ have? (We can actually list them:
$1,\ 2,\ 3,\ 4,\ 6,\ 8,\ 12,\ 24$. This one is again a mix of a counting
problem and a number theory problem.)

Okay, but let us generalize: How many positive divisors does a given positive
integer $n$ have? \textbf{We cannot solve this yet.} In this course, we will
not get to solve it, but it is not too hard to solve using the methods we have
learned (see \cite[\S 2.18.1]{19s} or \cite[Lecture 5, Exercise 5.3.3
\textbf{(a)}]{23f}).
\end{itemize}

\subsection{\label{sec.enum.lac}Lacunar subsets}

\subsubsection{\label{subsec.enum.lac.def}Definition}

Another type of objects that can be counted are the so-called \textbf{lacunar
subsets} (also known as \textbf{sparse subsets} to some authors). Here is
their definition:

\begin{definition}
A set $S$ of integers is said to be \textbf{lacunar} if it contains no two
consecutive integers (i.e., if there is no integer $i$ such that both $i$ and
$i+1$ belong to $S$).
\end{definition}

The word \textquotedblleft lacunar\textquotedblright\ comes from Latin
\textquotedblleft lacuna\textquotedblright\ (= \textquotedblleft
gap\textquotedblright). The idea is that a lacunar set has a \textquotedblleft
gap\textquotedblright\ (or \textquotedblleft buffer zone\textquotedblright)
between any two distinct elements.

For example, the set $\left\{  2,4,7\right\}  $ is lacunar, but the set
$\left\{  2,4,5\right\}  $ is not (since $4$ and $5$ are consecutive
integers). Any $1$-element set of integers is lacunar, and so is the empty set.

Now we can ask ourselves some natural questions: For given $n\in\mathbb{N}$,

\begin{enumerate}
\item how many lacunar subsets does the set $\left[  n\right]  =\left\{
1,2,\ldots,n\right\}  $ have?

\item how many $k$-element lacunar subsets does $\left[  n\right]  $ have for
a given $k\in\mathbb{N}$?

\item what is the largest size of a lacunar subset of $\left[  n\right]  $ ?
\end{enumerate}

We shall answer all these three questions in this section.

\subsubsection{\label{subsec.enum.lac.max}The maximum size of a lacunar
subset}

We start with the third question, as it is the easiest one to answer. Recall
the floor notation (Definition \ref{def.ent.floor}).

\begin{proposition}
\label{prop.count.lacs.max-size}Let $n\in\mathbb{N}$. Then, the maximum size
of a lacunar subset of $\left[  n\right]  $ is $\left\lfloor \dfrac{n+1}%
{2}\right\rfloor $.
\end{proposition}

\begin{proof}
The set
\begin{align*}
\left\{  \text{all odd numbers in }\left[  n\right]  \right\}   &  =\left\{
\text{all odd integers between }1\text{ and }n\text{ (inclusive)}\right\} \\
&  =\left\{  \text{all odd integers between }0\text{ and }n\text{
(inclusive)}\right\} \\
&  =\left\{  1,3,5,\ldots\right\}  \cap\left[  n\right]
\end{align*}
is a lacunar subset of $\left[  n\right]  $, and has size $\left\lfloor
\dfrac{n+1}{2}\right\rfloor $ (by Proposition \ref{prop.count.odd-0-n}). Thus,
the size $\left\lfloor \dfrac{n+1}{2}\right\rfloor $ is attainable (for a
lacunar subset of $\left[  n\right]  $).

It remains to show that this size is the largest possible -- i.e., that if $L$
is a lacunar subset of $\left[  n\right]  $, then
\[
\left\vert L\right\vert \leq\left\lfloor \dfrac{n+1}{2}\right\rfloor .
\]

So let $L$ be a lacunar subset of $\left[  n\right]  $. Our goal is to prove
that $\left\vert L\right\vert \leq\left\lfloor \dfrac{n+1}{2}\right\rfloor $.

We shall first prove that $\left\vert L\right\vert \leq\dfrac{n+1}{2}$.

Here are two different ways to prove this (each way illustrates a nice technique):

\begin{proof}
[First proof of $\left\vert L\right\vert \leq\dfrac{n+1}{2}$.]Let $\ell
_{1},\ell_{2},\ldots,\ell_{k}$ be the elements of $L$, listed in increasing
order, so that $L=\left\{  \ell_{1},\ell_{2},\ldots,\ell_{k}\right\}  $ and
$\ell_{1}<\ell_{2}<\cdots<\ell_{k}$. Thus, $\left\vert L\right\vert =k$.

Now, we assume (for the moment) that $k>0$. Thus, $k\geq1$ (since $k$ is an
integer). We have $\ell_{1}\in L\subseteq\left[  n\right]  $, so that
$\ell_{1}\geq1$. Moreover, the elements $\ell_{1}$ and $\ell_{2}$ of $L$
satisfy $\ell_{1}<\ell_{2}$ and $\ell_{2}\neq\ell_{1}+1$ (since $L$ is
lacunar), so that $\ell_{2}\geq\underbrace{\ell_{1}}_{\geq1}+\,2\geq1+2=3$.
Furthermore, the elements $\ell_{2}$ and $\ell_{3}$ of $L$ satisfy $\ell
_{2}<\ell_{3}$ and $\ell_{3}\neq\ell_{2}+1$ (since $L$ is lacunar), so that
$\ell_{3}\geq\underbrace{\ell_{2}}_{\geq3}+\,2\geq3+2=5$. Proceeding in the
same way, we find that%
\begin{equation}
\ell_{i}\geq2i-1\ \ \ \ \ \ \ \ \ \ \text{for each }i\in\left[  k\right]  .
\label{pf.prop.count.lacs.max-size.2i-1}%
\end{equation}
(Strictly speaking, this can be proved by induction on $i$. The base case
follows from $\ell_{1}\geq1=2\cdot1-1$, whereas the induction step requires
deriving $\ell_{i+1}\geq2\left(  i+1\right)  -1$ from $\ell_{i}\geq2i-1$,
which can be done by observing that $L$ is lacunar and therefore $\ell
_{i+1}\geq\underbrace{\ell_{i}}_{\geq2i-1}+\,2\geq2i-1+2=2\left(  i+1\right)
-1$.)

Now, we can apply (\ref{pf.prop.count.lacs.max-size.2i-1}) to $i=k$, and thus
obtain $\ell_{k}\geq2k-1$. However, $\ell_{k}\in L\subseteq\left[  n\right]
$, so that $\ell_{k}\leq n$. Thus, $n\geq\ell_{k}\geq2k-1$, so that
$n+1\geq2k$ and thus $\dfrac{n+1}{2}\geq k$. We have proved this under the
assumption that $k>0$, but this also holds in the opposite case (because if
$k\leq0$, then $\dfrac{n+1}{2}\geq0\geq k$). Thus, we always have $\dfrac
{n+1}{2}\geq k$ (independently of any assumptions). In other words, we have
$\dfrac{n+1}{2}\geq\left\vert L\right\vert $ (since $k=\left\vert L\right\vert
$). In other words, we have $\left\vert L\right\vert \leq\dfrac{n+1}{2}$.
\end{proof}

\begin{proof}
[Second proof of $\left\vert L\right\vert \leq\dfrac{n+1}{2}$.]Define a new
set
\[
L^{+}:=\left\{  \ell+1\ \mid\ \ell\in L\right\}  .
\]
This set $L^{+}$ consists of each element of $L$, incremented by $1$. For
example, if $L=\left\{  3,5,9\right\}  $, then $L^{+}=\left\{  4,6,10\right\}
$. Another way to view $L^{+}$ is as follows:
\[
L^{+}=\left\{  i\in\mathbb{Z}\ \mid\ i-1\in L\right\}
\]
(because an integer $i$ satisfies $i-1\in L$ if and only if it has the form
$\ell+1$ for some $\ell\in L$).

The set $L^{+}$ is just $L$ with each element incremented by $1$. Thus,
$\left\vert L^{+}\right\vert =\left\vert L\right\vert $.

Moreover, since $L$ is a subset of $\left[  n\right]  =\left\{  1,2,\ldots
,n\right\}  $, we conclude that $L^{+}$ is a subset of $\left\{
2,3,\ldots,n+1\right\}  $. Hence, both sets $L$ and $L^{+}$ are subsets of
$\left[  n+1\right]  $. Their union $L\cup L^{+}$ is thus a subset of $\left[
n+1\right]  $ as well. Therefore (by Theorem \ref{thm.count.rules.diff}
\textbf{(a)}, applied to $S=\left[  n+1\right]  $ and $T=L\cup L^{+}$), we
conclude that
\[
\left\vert L\cup L^{+}\right\vert \leq\left\vert \left[  n+1\right]
\right\vert =n+1.
\]

If the sets $L$ and $L^{+}$ had an element $j$ in common, then both $j-1$ and
$j$ would belong to $L$ (indeed, $j\in L^{+}=\left\{  i\in\mathbb{Z}%
\ \mid\ i-1\in L\right\}  $ would entail $j-1\in L$), which would contradict
the fact that $L$ is lacunar (since $j-1$ and $j$ are two consecutive
integers). Thus, the sets $L$ and $L^{+}$ have no element in common. In other
words, they are disjoint. Hence, by the sum rule (Theorem
\ref{thm.count.rules.sum-2}, applied to $A=L$ and $B=L^{+}$), we have
$\left\vert L\cup L^{+}\right\vert =\left\vert L\right\vert
+\underbrace{\left\vert L^{+}\right\vert }_{=\left\vert L\right\vert
}=\left\vert L\right\vert +\left\vert L\right\vert =2\cdot\left\vert
L\right\vert $. Hence,%
\[
2\cdot\left\vert L\right\vert =\left\vert L\cup L^{+}\right\vert \leq n+1.
\]
In other words, $\left\vert L\right\vert \leq\dfrac{n+1}{2}$.
\end{proof}

We have now proved (in two different ways) that $\left\vert L\right\vert
\leq\dfrac{n+1}{2}$. Now, recall the definition of the floor of a real number:
If $x$ is a real number, then $\left\lfloor x\right\rfloor $ is the largest
integer that is $\leq x$. Hence, $\left\lfloor \dfrac{n+1}{2}\right\rfloor $
is the largest integer that is $\leq\dfrac{n+1}{2}$. Therefore, any integer
that is $\leq\dfrac{n+1}{2}$ must also be $\leq\left\lfloor \dfrac{n+1}%
{2}\right\rfloor $. Applying this to the integer $\left\vert L\right\vert $,
we conclude that $\left\vert L\right\vert \leq\left\lfloor \dfrac{n+1}%
{2}\right\rfloor $ (since $\left\vert L\right\vert \leq\dfrac{n+1}{2}$). As
explained above, this completes the proof of Proposition
\ref{prop.count.lacs.max-size}.
\end{proof}

\subsubsection{\label{subsec.enum.lac.count-all}Counting all lacunar subsets
of $\left[  n\right]  $}

Now let us count the lacunar subsets of $\left[  n\right]  $. We shall first
count them all, then count the ones of a given size $k$.

First, a few words about how to find answers to counting questions like this.
For any specific value of $n$, finding the \# of lacunar subsets of $\left[
n\right]  $ is a \textquotedblleft finite problem\textquotedblright: You can
just count them all. Or, better, you can have your computer do this. In
\href{https://www.sagemath.org/}{SageMath} (a computer algebra system, one of
the best suited to combinatorial questions), this takes just a few
lines:\footnote{SageMath is built on top of the Python programming language,
so you will recognize a lot of Python syntax. Actually, the only piece of
non-Python code in the following code snippet is the \texttt{Subsets(n)} part.
If you want to use (pure) Python instead of SageMath, you can replace
\texttt{sum(1\ for\ S\ in\ Subsets(n)\ if\ is\_lacunar(S))} by
\texttt{sum(1\ for\ i in range(n+1) for S\ in\ combinations(range(1, n+1),
i)\ if\ is\_lacunar(S))}, after first importing the \texttt{combinations}
function from the \texttt{itertools} package (using \texttt{from itertools
import combinations}).} \medskip
\begin{verbatim}
  def is_lacunar(S): # test if the set S is lacunar
    return all(i+1 not in S for i in S)

  def num_lacs(n): # number of lacunar subsets of [n]
    return sum(1 for S in Subsets(n) if is_lacunar(S))

  for n in range(10):
    print("For n = " + str(n) + ", the number is " + str(num_lacs(n)))
\end{verbatim}

\medskip The first two lines here speak for themselves (once you know that
\texttt{all} is the universal quantifier). The function \texttt{Subsets}
computes the set of all subsets of a given set, or (if we provide it an
integer $n$ as input) all subsets of $\left[  n\right]  $. The \texttt{sum(1
for S in SomeSet)} construction is just a slick way of counting the elements
of \texttt{SomeSet}, exploiting the fact that a sum of the form $1+1+\cdots+1$
equals the number of its addends. The last two lines are prompting SageMath to
compute the \# of lacunar subsets of $\left[  n\right]  $ for each
$n\in\left[  0,9\right]  $ (note that \texttt{range(a, b)} means the integer
interval $\left[  a,b-1\right]  $ in SageMath) and to output these $10$
numbers. I refer to \cite[\S 1.4.3]{19fco} for more hints on the use of
SageMath, and to \href{https://doc.sagemath.org/}{its documentation} for a
more systematic introduction. Note that you can use
\href{https://sagecell.sagemath.org/}{SageMathCell} to easily call SageMath
from your browser (although the computations you call are limited by 30
seconds each, since they happen on the server).

The answers we get from SageMath are interesting:%
\[%
\begin{tabular}
[c]{|c||c|c|c|c|c|c|c|c|c|c|}\hline
$n$ & $0$ & $1$ & $2$ & $3$ & $4$ & $5$ & $6$ & $7$ & $8$ & $9$\\\hline
\# of lacunar subsets of $\left[  n\right]  $ & $1$ & $2$ & $3$ & $5$ & $8$ &
$13$ & $21$ & $34$ & $55$ & $89$\\\hline
\end{tabular}
\ \ \ \ .
\]
Haven't we seen these numbers before?

Yes, we have: In Definition \ref{def.ind.fib}, we defined the
\textbf{Fibonacci sequence}. This is the sequence $\left(  f_{0},f_{1}%
,f_{2},\ldots\right)  $ of nonnegative integers defined recursively by setting%
\begin{align*}
f_{0}  &  =0,\ \ \ \ \ \ \ \ \ \ f_{1}=1,\ \ \ \ \ \ \ \ \ \ \text{and}\\
f_{n}  &  =f_{n-1}+f_{n-2}\ \ \ \ \ \ \ \ \ \ \text{for each }n\geq2.
\end{align*}
Its first few entries are%
\[%
\begin{tabular}
[c]{|c||c|c|c|c|c|c|c|c|c|c|c|c|c|c|}\hline
$n$ & $0$ & $1$ & $2$ & $3$ & $4$ & $5$ & $6$ & $7$ & $8$ & $9$ & $10$ & $11$
& $12$ & $13$\\\hline
$f_{n}$ & $0$ & $1$ & $1$ & $2$ & $3$ & $5$ & $8$ & $13$ & $21$ & $34$ & $55$
& $89$ & $144$ & $233$\\\hline
\end{tabular}
\ \ \ \ \ .
\]

The two above tables have the same entries, if you discount the fact that the
first two Fibonacci numbers $f_{0}=0$ and $f_{1}=1$ are missing from the
former table. So we have good reasons to suspect that%
\[
\left(  \text{\# of lacunar subsets of }\left[  n\right]  \right)  =f_{n+2}%
\]
for each $n\in\mathbb{N}$. And indeed, this is true:

\begin{theorem}
\label{thm.count.lacs.all}For any integer $n\geq-1$, we have%
\[
\left(  \text{\# of lacunar subsets of }\left[  n\right]  \right)  =f_{n+2}.
\]
Here, we agree that $\left[  -1\right]  :=\varnothing$. More generally, we
agree that $\left[  k\right]  :=\varnothing$ for any $k\leq0$.
\end{theorem}

\begin{example}
The lacunar subsets of $\left[  4\right]  $ are
\[
\varnothing,\ \ \left\{  1\right\}  ,\ \ \left\{  2\right\}  ,\ \ \left\{
3\right\}  ,\ \ \left\{  4\right\}  ,\ \ \left\{  1,3\right\}  ,\ \ \left\{
1,4\right\}  ,\ \ \left\{  2,4\right\}  .
\]
So there are $8$ of them, as predicted by Theorem \ref{thm.count.lacs.all}
(since $f_{4+2}=f_{6}=8$).
\end{example}

(Are you wondering why we are allowing $n$ to be $-1$ in Theorem
\ref{thm.count.lacs.all}? The answer is \textquotedblleft because we
can\textquotedblright, and more precisely \textquotedblleft because it will
make our induction easier\textquotedblright. The case $n=-1$ is not
interesting by itself; the claim of Theorem \ref{thm.count.lacs.all} in this
case is just that the \# of lacunar subsets of $\varnothing$ is $1$.)

\begin{proof}
[Proof of Theorem \ref{thm.count.lacs.all}.]For any integer $n\geq-1$, let us
set%
\[
\ell_{n}:=\left(  \text{\# of lacunar subsets of }\left[  n\right]  \right)
.
\]
Thus, we must prove that
\begin{equation}
\ell_{n}=f_{n+2}\ \ \ \ \ \ \ \ \ \ \text{for each }n\geq-1.
\label{pf.thm.count.lacs.all.goal}%
\end{equation}

We have $\ell_{-1}=1$ (since the set $\left[  -1\right]  =\varnothing$ has
only one lacunar subset, namely $\varnothing$ itself) and $f_{-1+2}=f_{1}=1$.
Hence, $\ell_{-1}=1=f_{-1+2}$. In other words,
(\ref{pf.thm.count.lacs.all.goal}) holds for $n=-1$. A similar computation
shows that (\ref{pf.thm.count.lacs.all.goal}) holds for $n=0$.

Let us next show the following:

\begin{statement}
\textit{Claim 1:} We have $\ell_{n}=\ell_{n-1}+\ell_{n-2}$ for each integer
$n\geq1$.
\end{statement}

\begin{proof}
[Proof of Claim 1.]Let $n\geq1$ be an integer. We shall call a subset of
$\left[  n\right]  $

\begin{itemize}
\item \textbf{red} if it contains $n$, and

\item \textbf{green} if it does not contain $n$.
\end{itemize}

Then, the definition of $\ell_{n}$ shows that%
\begin{align*}
\ell_{n}  &  =\left(  \text{\# of lacunar subsets of }\left[  n\right]
\right) \\
&  =\left(  \text{\# of red lacunar subsets of }\left[  n\right]  \right)
+\left(  \text{\# of green lacunar subsets of }\left[  n\right]  \right)
\end{align*}
(by the sum rule, since each lacunar subset of $\left[  n\right]  $ is either
red or green but cannot be both at the same time\footnote{This is the same
argument that has been used in the proof of Theorem
\ref{thm.count.subsets.all.n-repeat}.}).

The green lacunar subsets of $\left[  n\right]  $ are just the lacunar subsets
of $\left[  n-1\right]  $ (since \textquotedblleft green\textquotedblright%
\ means \textquotedblleft does not contain $n$\textquotedblright). Thus,%
\begin{align*}
&  \left(  \text{\# of green lacunar subsets of }\left[  n\right]  \right) \\
&  =\left(  \text{\# of lacunar subsets of }\left[  n-1\right]  \right)
=\ell_{n-1}%
\end{align*}
(by the definition of $\ell_{n-1}$).

Counting the red lacunar subsets is trickier. We shall show that their \# is
$\ell_{n-2}$.

If $R$ is a red lacunar subset of $\left[  n\right]  $, then $R$ contains $n$
(by the definition of \textquotedblleft red\textquotedblright), so that $R$
does not contain $n-1$ (by lacunarity), and therefore $R\setminus\left\{
n\right\}  $ is a lacunar subset of $\left[  n-2\right]  $ (since
$R\setminus\left\{  n\right\}  $ contains neither $n$ nor $n-1$). Thus, we
obtain a map%
\begin{align*}
\operatorname*{rem}\nolimits_{n}:\left\{  \text{red lacunar subsets of
}\left[  n\right]  \right\}   &  \rightarrow\left\{  \text{lacunar subsets of
}\left[  n-2\right]  \right\}  ,\\
R  &  \mapsto R\setminus\left\{  n\right\}  .
\end{align*}

Conversely, if $L$ is a lacunar subset of $\left[  n-2\right]  $, then
$L\cup\left\{  n\right\}  $ is a lacunar subset of $\left[  n\right]  $
(indeed, the integer $n-1$ is a \textquotedblleft buffer
zone\textquotedblright\ between the elements of $L$ and the new element $n$,
so that the lacunarity of $L$ is preserved when we insert $n$ into the set),
and is red (since $n\in\left\{  n\right\}  \subseteq L\cup\left\{  n\right\}
$). Thus, we obtain a map%
\begin{align*}
\operatorname*{ins}\nolimits_{n}:\left\{  \text{lacunar subsets of }\left[
n-2\right]  \right\}   &  \rightarrow\left\{  \text{red lacunar subsets of
}\left[  n\right]  \right\}  ,\\
L  &  \mapsto L\cup\left\{  n\right\}  .
\end{align*}

It is easy to see (just as in the proof of Theorem
\ref{thm.count.subsets.all.n-repeat}) that the map $\operatorname*{rem}%
\nolimits_{n}$ is an inverse of $\operatorname*{ins}\nolimits_{n}$. Thus, the
map $\operatorname*{ins}\nolimits_{n}$ has an inverse, i.e., is bijective (by
Theorem \ref{thm.maps.jects.inv=bij}). Hence, we have found a bijection%
\[
\text{from }\left\{  \text{lacunar subsets of }\left[  n-2\right]  \right\}
\text{ to }\left\{  \text{red lacunar subsets of }\left[  n\right]  \right\}
\]
(namely, $\operatorname*{ins}\nolimits_{n}$). Therefore, by the bijection
principle, we have%
\[
\left\vert \left\{  \text{lacunar subsets of }\left[  n-2\right]  \right\}
\right\vert =\left\vert \left\{  \text{red lacunar subsets of }\left[
n\right]  \right\}  \right\vert .
\]
In other words,%
\[
\left(  \text{\# of lacunar subsets of }\left[  n-2\right]  \right)  =\left(
\text{\# of red lacunar subsets of }\left[  n\right]  \right)  .
\]
Thus,%
\[
\left(  \text{\# of red lacunar subsets of }\left[  n\right]  \right)
=\left(  \text{\# of lacunar subsets of }\left[  n-2\right]  \right)
=\ell_{n-2}%
\]
(by the definition of $\ell_{n-2}$).

Altogether,%
\begin{align*}
\ell_{n}  &  =\underbrace{\left(  \text{\# of red lacunar subsets of }\left[
n\right]  \right)  }_{=\ell_{n-2}}+\underbrace{\left(  \text{\# of green
lacunar subsets of }\left[  n\right]  \right)  }_{=\ell_{n-1}}\\
&  =\ell_{n-2}+\ell_{n-1}=\ell_{n-1}+\ell_{n-2}.
\end{align*}
This proves Claim 1.
\end{proof}

Now we still need to prove (\ref{pf.thm.count.lacs.all.goal}). In other words,
we need to prove that the two sequences $\left(  \ell_{-1},\ell_{0},\ell
_{1},\ldots\right)  $ and $\left(  f_{1},f_{2},f_{3},\ldots\right)  $ are
identical. But at this point, this is very easy: These two sequences

\begin{itemize}
\item have the same two starting entries $\ell_{-1}=f_{1}$ and $\ell_{0}%
=f_{2}$ (this can be easily checked directly),

\item and satisfy the same recursive equation: namely, each entry of either
sequence is the sum of the preceding two entries (since Claim 1 yields
$\ell_{n}=\ell_{n-1}+\ell_{n-2}$, whereas the definition of the Fibonacci
numbers yields $f_{n+2}=f_{n+1}+f_{n}$).
\end{itemize}

Since a recursively defined sequence is uniquely determined by its starting
entries and its recursive equation, we thus conclude that the two sequences
$\left(  \ell_{-1},\ell_{0},\ell_{1},\ldots\right)  $ and $\left(  f_{1}%
,f_{2},f_{3},\ldots\right)  $ are identical. Thus,
(\ref{pf.thm.count.lacs.all.goal}) follows. This slightly informal argument
can be formalized as a straightforward strong
induction\footnote{\textit{Proof.} Let us prove
(\ref{pf.thm.count.lacs.all.goal}) by strong induction on $n$:
\par
\textit{Base case:} We have already checked that
(\ref{pf.thm.count.lacs.all.goal}) holds for $n=-1$.
\par
\textit{Induction step:} Let $n\geq0$ be an integer. Assume (as the induction
hypothesis) that the claim (\ref{pf.thm.count.lacs.all.goal}) holds for each
of $-1,\ 0,\ 1,\ \ldots,\ n-1$ instead of $n$. We must prove that
(\ref{pf.thm.count.lacs.all.goal}) holds for $n$ as well, i.e., that we have
$\ell_{n}=f_{n+2}$.
\par
If $n=0$, then this follows from the fact (observed above) that
(\ref{pf.thm.count.lacs.all.goal}) holds for $n=0$. It thus remains to
consider the case when $n\neq0$. So let us assume that $n\neq0$. Since
$n\geq0$, we thus obtain $n\geq1$, so that $n-1\geq0$ and $n-2\geq-1$.
\par
In particular, $n-1\geq0\geq-1$. Hence, our induction hypothesis yields that
the claim (\ref{pf.thm.count.lacs.all.goal}) holds for $n-1$ instead of $n$.
In other words, we have $\ell_{n-1}=f_{\left(  n-1\right)  +2}=f_{n+1}$.
\par
Also, our induction hypothesis yields that the claim
(\ref{pf.thm.count.lacs.all.goal}) holds for $n-2$ instead of $n$ (since
$n-2\geq-1$). In other words, we have $\ell_{n-2}=f_{\left(  n-2\right)
+2}=f_{n}$.
\par
Now, Claim 1 yields $\ell_{n}=\underbrace{\ell_{n-1}}_{=f_{n+1}}%
+\underbrace{\ell_{n-2}}_{=f_{n}}=f_{n+1}+f_{n}$. But the recursive definition
of the Fibonacci sequence also yields $f_{n+2}=f_{n+1}+f_{n}$. Comparing these
two equalities, we find $\ell_{n}=f_{n+2}$. In other words,
(\ref{pf.thm.count.lacs.all.goal}) holds for $n$. This completes the induction
step. Thus, (\ref{pf.thm.count.lacs.all.goal}) is proved.}.

Thus we have proved (\ref{pf.thm.count.lacs.all.goal}). In other words, we
have proved Theorem \ref{thm.count.lacs.all} (because we have $\ell
_{n}=\left(  \text{\# of lacunar subsets of }\left[  n\right]  \right)  $).
\end{proof}

\subsubsection{\label{subsec.enum.lac.count-k}Counting all $k$-element lacunar
subsets of $\left[  n\right]  $}

Let us now address the remaining question about lacunar subsets: counting
$k$-element lacunar subsets of $\left[  n\right]  $ for given $n$ and $k$.

Again, we start by asking SageMath for some data:
\begin{verbatim}
  def is_lacunar(S): # test if the set S is lacunar
    return all(i+1 not in S for i in S)

  def num_lacs(n, k): # number of k-element lacunar subsets of [n]
    return sum(1 for S in Subsets(n, k) if is_lacunar(S))

  for n in range(10):
    print("For n = " + str(n) + ", the numbers are " + \
          str([num_lacs(n, k) for k in range(n+1)]))
\end{verbatim}

We obtain the following table:%
\[%
\begin{tabular}
[c]{|c||c|c|c|c|c|c|}\hline
& $k=0$ & $k=1$ & $k=2$ & $k=3$ & $k=4$ & $k=5$\\\hline\hline
$n=0$ & $1$ &  &  &  &  & \\\hline
$n=1$ & $1$ & $1$ &  &  &  & \\\hline
$n=2$ & $1$ & $2$ &  &  &  & \\\hline
$n=3$ & $1$ & $3$ & $1$ &  &  & \\\hline
$n=4$ & $1$ & $4$ & $3$ &  &  & \\\hline
$n=5$ & $1$ & $5$ & $6$ & $1$ &  & \\\hline
$n=6$ & $1$ & $6$ & $10$ & $4$ &  & \\\hline
$n=7$ & $1$ & $7$ & $15$ & $10$ & $1$ & \\\hline
$n=8$ & $1$ & $8$ & $21$ & $20$ & $5$ & \\\hline
$n=9$ & $1$ & $9$ & $28$ & $35$ & $15$ & $1$\\\hline
\end{tabular}
\]
(where each entry is the \# of lacunar $k$-element subsets of $\left[
n\right]  $ for the corresponding values of $n$ and $k$, and where an empty
box means that the corresponding \# is $0$). The many $0$'s are unsurprising
(they are predicted by Proposition \ref{prop.count.lacs.max-size}), and
likewise the values for $k=0$ and $k=1$ are clear (since every subset that has
size $\leq1$ is lacunar). But staring at the table for a bit longer reveals
something subtler: It is a sheared Pascal's triangle! For example, the $n=7$
row contains the numbers $1,7,15,10,1$, which appear along a diagonal in
Pascal's triangle. All the entries are binomial coefficients, and a bit of
work reveals the exact formula:

\begin{theorem}
\label{thm.count.lacs.size}Let $n\in\mathbb{Z}$ and $k\in\mathbb{N}$ be such
that $k\leq n+1$. Then,%
\[
\left(  \text{\# of }k\text{-element lacunar subsets of }\left[  n\right]
\right)  =\dbinom{n+1-k}{k}.
\]

\end{theorem}

For instance, for $n=7$ and $k=3$, this yields
\[
\left(  \text{\# of }3\text{-element lacunar subsets of }\left[  7\right]
\right)  =\dbinom{7+1-3}{3}=\dbinom{5}{3}=10,
\]
which agrees with our above table.

Note that the condition $k\leq n+1$ in Theorem \ref{thm.count.lacs.size} is
needed. If $k>n+1$ and $n\in\mathbb{N}$, then the \# of $k$-element lacunar
subsets of $\left[  n\right]  $ is $0$ (since a subset of $\left[  n\right]  $
cannot have more than $n$ elements, let alone more than $n+1$ elements, and
even less so when it is lacunar), but the binomial coefficient $\dbinom
{n+1-k}{k}$ is nonzero (since the $n+1-k$ on its top is negative).

You can prove Theorem \ref{thm.count.lacs.size} by induction on $n$, using a
similar red/green coloring as in our above proof of Theorem
\ref{thm.count.lacs.all} (and carefully checking that the condition $k\leq
n+1$ is satisfied whenever you apply the induction hypothesis\footnote{This
necessitates a bit of casework.}). Such a proof can be found in \cite[Exercise
3 \textbf{(a)}]{17f-hw2s}\footnote{To be \textbf{very} pedantic:
\cite[Exercise 3 \textbf{(a)}]{17f-hw2s} only states Theorem
\ref{thm.count.lacs.size} in the case when $n\in\mathbb{N}$. But the remaining
case is trivial (since $k\leq n+1$ leads to $k=0$ when $n$ is negative, and
thus we have to count $0$-element subsets of an empty set, which is not a deep
question).}.

There is, however, a nicer proof, which proceeds by constructing a bijection%
\begin{align*}
&  \text{from }\left\{  k\text{-element lacunar subsets of }\left[  n\right]
\right\} \\
&  \text{to }\left\{  k\text{-element subsets of }\left[  n+1-k\right]
\right\}  ,
\end{align*}
and observing that the \# of $k$-element subsets of $\left[  n+1-k\right]  $
is $\dbinom{n+1-k}{k}$ (by Theorem \ref{thm.count.subsets.size.S-repeat}).
Such a proof has the advantage of not just proving Theorem
\ref{thm.count.lacs.size} but also explaining \textquotedblleft
why\textquotedblright\ it holds (at least if you consider it as a given that
binomial coefficients count $k$-element subsets).

This second proof rests upon a basic feature of finite sets of integers:

\begin{proposition}
\label{prop.intset.inclist-k} Let $k\in\mathbb{N}$. Let $S$ be a $k$-element
set of integers. Then, there exists a unique $k$-tuple $\left(  s_{1}%
,s_{2},\ldots,s_{k}\right)  $ of integers satisfying
\[
\left\{  s_{1},s_{2},\ldots,s_{k}\right\}  =S\ \ \ \ \ \ \ \ \ \ \text{and}%
\ \ \ \ \ \ \ \ \ \ s_{1}<s_{2}<\cdots<s_{k}.
\]

\end{proposition}

This proposition is just saying that if you are given a $k$-element set $S$ of
integers, then there is a unique way to list the elements of $S$ in increasing
order (with no repetitions). Intuitively, this is clear (just write down the
smallest element of $S$, then the second-smallest element, then the
third-smallest, and so on, until you run out of elements; it's not like you
have any other options!). But intuition is not proof. Nevertheless, we will
not stoop down to this low a foundational level here\footnote{A boring and
detailed (but ultimately very simple) proof of Proposition
\ref{prop.intset.inclist-k} can be found in \cite[proof of Theorem
2.46]{detnotes}.}, and just take Proposition \ref{prop.intset.inclist-k} for granted.

In connection with Proposition \ref{prop.intset.inclist-k}, we introduce a notation:

\begin{convention}
\label{conv.intset.inclist-setbrace}Let $s_{1},s_{2},\ldots,s_{k}$ be some
integers. Then, the notation \textquotedblleft$\left\{  s_{1}<s_{2}%
<\cdots<s_{k}\right\}  $\textquotedblright\ shall mean the set $\left\{
s_{1},s_{2},\ldots,s_{k}\right\}  $ and additionally signify that the chain of
inequalities $s_{1}<s_{2}<\cdots<s_{k}$ holds.
\end{convention}

Thus, for example, $\left\{  2<4<5\right\}  $ is the set $\left\{
2,4,5\right\}  $, whereas the expression $\left\{  4<2<5\right\}  $ is meaningless.

Proposition \ref{prop.intset.inclist-k} can now be restated as follows: If
$k\in\mathbb{N}$, then any $k$-element set of integers can be written in the
form $\left\{  s_{1}<s_{2}<\cdots<s_{k}\right\}  $ for a unique $k$-tuple
$\left(  s_{1},s_{2},\ldots,s_{k}\right)  $ of integers.

We are now ready to prove Theorem \ref{thm.count.lacs.size}:

\begin{proof}
[Proof of Theorem \ref{thm.count.lacs.size}.]Let $m:=n+1-k$. Then,
$m=n+1-k\geq0$ (since $k\leq n+1$), so that $\left[  m\right]  $ is an
$m$-element set. Also, $m=n+1-k=n-\left(  k-1\right)  $, so that $m+\left(
k-1\right)  =n$.

Now, if $S=\left\{  s_{1}<s_{2}<\cdots<s_{k}\right\}  $ is a $k$-element
lacunar subset of $\left[  n\right]  $, then $\overleftarrow{S}$ shall mean
the set%
\[
\left\{  s_{i}-\left(  i-1\right)  \ \mid\ i\in\left[  k\right]  \right\}
=\left\{  s_{1},\ s_{2}-1,\ s_{3}-2,\ \ldots,\ s_{k}-\left(  k-1\right)
\right\}  .
\]
This set $\overleftarrow{S}$ is obtained from $S$ by the following process:

\begin{itemize}
\item Leave the smallest element of $S$ unchanged.

\item Decrease the second-smallest element of $S$ by $1$.

\item Decrease the third-smallest element of $S$ by $2$.

\item And so on, until eventually decreasing the largest (= $k$-th-smallest)
element of $S$ by $k-1$.
\end{itemize}

We refer to this process as the \textbf{compression process}, as it causes the
elements of $S$ to come closer together (in such a way that the distance
between any two \textquotedblleft positionally adjacent\textquotedblright%
\ elements\footnote{We call two elements $i$ and $j$ of $S$ \textquotedblleft%
\textbf{positionally adjacent}\textquotedblright\ if they satisfy $i<j$ but
there are no other elements of $S$ lying between them (i.e., there are no
elements $s\in S$ satisfying $i<s<j$). For example, the elements $4$ and $6$
of the set $\left\{  2,4,6,8\right\}  $ are positionally adjacent, but the
elements $4$ and $6$ of the set $\left\{  2,3,4,5,6\right\}  $ are not (since
the element $5$ lies between them).} of $S$ shrinks by $1$). Consequently, we
call the resulting set $\overleftarrow{S}$ the \textbf{compression} of $S$.
For example, if $S=\left\{  3<5<9<11\right\}  $, then $\overleftarrow{S}%
=\left\{  3<4<7<8\right\}  $. Let us illustrate this example graphically:%
\[%
\begin{tikzpicture}
\filldraw[fill=cyan!20!white] (0.5, 1.5) rectangle (12, 2.5);
\node(A3) at (3, 2) {$3$};
\node(A5) at (5, 2) {$5$};
\node(A9) at (9, 2) {$9$};
\node(A11) at (11, 2) {$11$};
\node(Aname) at (0, 2) {$S$};
\filldraw[fill=yellow!80!white] (0.5, -0.5) rectangle (9, 0.5);
\node(B3) at (3, 0) {$3$};
\node(B4) at (4, 0) {$4$};
\node(B7) at (7, 0) {$7$};
\node(B8) at (8, 0) {$8$};
\node(Bname) at (0, 0) {$\overleftarrow{S}$};
\draw[->, very thick, red] (A3) -- (B3);
\draw[->, very thick, red] (A5) -- (B4);
\draw[->, very thick, red] (A9) -- (B7);
\draw[->, very thick, red] (A11) -- (B8);
\end{tikzpicture}%
\]
(note that each of the red arrows is slightly more horizontal than the
previous one).

We note the following properties of compression: If $S=\left\{  s_{1}%
<s_{2}<\cdots<s_{k}\right\}  $ is a $k$-element lacunar subset of $\left[
n\right]  $, then its compression $\overleftarrow{S}$ is still a $k$-element
set (i.e., the compression process does not cause any two distinct elements to
\textquotedblleft collide\textquotedblright) and can be written as%
\[
\left\{  s_{1}<s_{2}-1<s_{3}-2<\cdots<s_{k}-\left(  k-1\right)  \right\}
\]
(since $S$ is lacunar, so that any two \textquotedblleft positionally
adjacent\textquotedblright\ elements $s_{i}$ and $s_{i+1}$ of $S$ satisfy
$s_{i}<s_{i+1}-1$ and thus $s_{i}-\left(  i-1\right)  <\left(  s_{i+1}%
-1\right)  -\left(  i-1\right)  =s_{i+1}-i$). Furthermore, $\overleftarrow{S}$
is a subset of $\left[  m\right]  $ (because its smallest
element\footnote{Talking of the \textquotedblleft smallest
element\textquotedblright\ is somewhat risky, since $k$ can be $0$, in which
case $S$ has no elements at all. But this kind of problem appears only for
$k=0$, in which case everything we are trying to prove is trivial or very
easy. Thus, we will be cavalier about the case $k=0$ in this proof.} is
$s_{1}\geq1$ (since $s_{1}\in S\subseteq\left[  n\right]  $), whereas its
largest element is $\underbrace{s_{k}}_{\substack{\leq n\\\text{(since }%
s_{k}\in S\subseteq\left[  n\right]  \text{)}}}-\left(  k-1\right)  \leq
n-\left(  k-1\right)  =m$). Thus, we can define a map%
\begin{align*}
\operatorname*{compress}:\left\{  k\text{-element lacunar subsets of }\left[
n\right]  \right\}   &  \rightarrow\left\{  k\text{-element subsets of
}\left[  m\right]  \right\}  ,\\
S  &  \mapsto\overleftarrow{S}.
\end{align*}

Conversely, if $T=\left\{  t_{1}<t_{2}<\cdots<t_{k}\right\}  $ is a
$k$-element subset of $\left[  m\right]  $, then $\overrightarrow{T}$ shall
mean the set%
\[
\left\{  t_{i}+\left(  i-1\right)  \ \mid\ i\in\left[  k\right]  \right\}
=\left\{  t_{1},\ t_{2}+1,\ t_{3}+2,\ \ldots,\ t_{k}+\left(  k-1\right)
\right\}  .
\]
This set $\overrightarrow{T}$ is obtained from $T$ by the following process:

\begin{itemize}
\item Leave the smallest element of $T$ unchanged.

\item Increase the second-smallest element of $T$ by $1$.

\item Increase the third-smallest element of $T$ by $2$.

\item And so on, until eventually increasing the largest (= $k$-th-smallest)
element of $T$ by $k-1$.
\end{itemize}

We refer to this process as the \textbf{expansion process}, as it causes the
elements of $T$ to drift further apart (in such a way that the distance
between any two \textquotedblleft positionally adjacent\textquotedblright%
\ elements of $T$ increases by $1$). Consequently, we call the resulting set
$\overrightarrow{T}$ the \textbf{expansion} of $T$. For example, if
$T=\left\{  3<4<7<8\right\}  $, then $\overrightarrow{T}=\left\{
3<5<9<11\right\}  $. Let us illustrate this example graphically:%
\[%
\begin{tikzpicture}
\filldraw[fill=yellow!80!white] (0.5, 1.5) rectangle (9, 2.5);
\node(A3) at (3, 2) {$3$};
\node(A4) at (4, 2) {$4$};
\node(A7) at (7, 2) {$7$};
\node(A8) at (8, 2) {$8$};
\node(Aname) at (0, 2) {$T$};
\filldraw[fill=cyan!20!white] (0.5, -0.5) rectangle (12, 0.5);
\node(B3) at (3, 0) {$3$};
\node(B5) at (5, 0) {$5$};
\node(B9) at (9, 0) {$9$};
\node(B11) at (11, 0) {$11$};
\node(Bname) at (0, 0) {$\overrightarrow{T}$};
\draw[->, very thick, red] (A3) -- (B3);
\draw[->, very thick, red] (A4) -- (B5);
\draw[->, very thick, red] (A7) -- (B9);
\draw[->, very thick, red] (A8) -- (B11);
\end{tikzpicture}%
\]
(note that each of the red arrows is slightly more horizontal than the
previous one).

We note the following properties of expansion: If $T=\left\{  t_{1}%
<t_{2}<\cdots<t_{k}\right\}  $ is a $k$-element subset of $\left[  m\right]
$, then its expansion $\overrightarrow{T}$ is still a $k$-element set (i.e.,
the expansion process does not cause any two distinct elements to
\textquotedblleft collide\textquotedblright) and can be written as%
\[
\left\{  t_{1}<t_{2}+1<t_{3}+2<\cdots<t_{k}+\left(  k-1\right)  \right\}
\]
(since each $i\in\left[  k-1\right]  $ satisfies $t_{i}<t_{i+1}$ and thus
$t_{i}+\left(  i-1\right)  <t_{i+1}+\left(  i-1\right)  <t_{i+1}+i$).
Furthermore, $\overrightarrow{T}$ is a subset of $\left[  n\right]  $ (because
its smallest element is $t_{1}\geq1$ (since $t_{1}\in T\subseteq\left[
m\right]  $), whereas its largest element is $\underbrace{t_{k}}%
_{\substack{\leq m\\\text{(since }t_{k}\in T\subseteq\left[  m\right]
\text{)}}}+\left(  k-1\right)  \leq m+\left(  k-1\right)  =n$), and is lacunar
(since the expansion process ensures that the distance between any two
\textquotedblleft positionally adjacent\textquotedblright\ elements of $T$ has
been increased by $1$ in $\overrightarrow{T}$, so they can no longer be
consecutive integers). Thus, we can define a map%
\begin{align*}
\operatorname{expand}:\left\{  k\text{-element subsets of }\left[  m\right]
\right\}   &  \rightarrow\left\{  k\text{-element lacunar subsets of }\left[
n\right]  \right\}  ,\\
T  &  \mapsto\overrightarrow{T}.
\end{align*}

It is easy to see that the map $\operatorname{expand}$ is an inverse of
$\operatorname*{compress}$\ \ \ \ \footnote{In fact, each $k$-element subset
$T$ of $\left[  m\right]  $ satisfies $\operatorname*{compress}\left(
\operatorname{expand}T\right)  =T$, because if we write $T$ as $T=\left\{
t_{1}<t_{2}<\cdots<t_{k}\right\}  $, then%
\[
\operatorname{expand}T=\operatorname{expand}\left(  \left\{  t_{1}%
<t_{2}<\cdots<t_{k}\right\}  \right)  =\left\{  t_{1}<t_{2}+1<t_{3}%
+2<\cdots<t_{k}+\left(  k-1\right)  \right\}
\]
and therefore%
\begin{align*}
\operatorname*{compress}\left(  \operatorname{expand}T\right)   &
=\operatorname*{compress}\left(  \left\{  t_{1}<t_{2}+1<t_{3}+2<\cdots
<t_{k}+\left(  k-1\right)  \right\}  \right) \\
&  =\left\{  t_{1}<\left(  t_{2}+1\right)  -1<\left(  t_{3}+2\right)
-2<\cdots<\left(  t_{k}+\left(  k-1\right)  \right)  -\left(  k-1\right)
\right\} \\
&  =\left\{  t_{1}<t_{2}<\cdots<t_{k}\right\}  =T.
\end{align*}
A similar argument shows that any $k$-element lacunar subset $S$ of $\left[
n\right]  $ satisfies $\operatorname{expand}\left(  \operatorname*{compress}%
S\right)  =S$.}. Hence, the map $\operatorname*{compress}$ has an inverse,
i.e., is bijective. Thus, it is a bijection from $\left\{  k\text{-element
lacunar subsets of }\left[  n\right]  \right\}  $ to $\left\{  k\text{-element
subsets of }\left[  m\right]  \right\}  $. Hence, the bijection principle
yields%
\begin{align*}
&  \left(  \text{\# of }k\text{-element lacunar subsets of }\left[  n\right]
\right) \\
&  =\left(  \text{\# of }k\text{-element subsets of }\left[  m\right]  \right)
\\
&  =\dbinom{m}{k}\ \ \ \ \ \ \ \ \ \ \left(
\begin{array}
[c]{c}%
\text{by Theorem \ref{thm.count.subsets.size.S-repeat}}\\
\text{(applied to }m\text{ and }\left[  m\right]  \text{ instead of }n\text{
and }S\text{)}%
\end{array}
\right) \\
&  =\dbinom{n+1-k}{k}\ \ \ \ \ \ \ \ \ \ \left(  \text{since }m=n+1-k\right)
.
\end{align*}
This proves Theorem \ref{thm.count.lacs.size}.
\end{proof}

\subsubsection{\label{subsec.enum.lac.binom-sum}A corollary}

Combining Theorem \ref{thm.count.lacs.size} with Theorem
\ref{thm.count.lacs.all}, we obtain a curious formula for the Fibonacci
numbers in terms of binomial coefficients:

\begin{corollary}
\label{cor.fib.sum-binom}Let $n\in\mathbb{N}$. Then, the Fibonacci number
$f_{n+1}$ is%
\[
f_{n+1}=\sum_{k=0}^{n}\dbinom{n-k}{k}=\dbinom{n-0}{0}+\dbinom{n-1}{1}%
+\cdots+\dbinom{n-n}{n}.
\]

\end{corollary}

\begin{example}
For $n=6$, Corollary \ref{cor.fib.sum-binom} says that%
\begin{align*}
f_{7}  &  =\dbinom{6-0}{0}+\dbinom{6-1}{1}+\dbinom{6-2}{2}+\dbinom{6-3}{3}\\
&  \ \ \ \ \ \ \ \ \ \ +\dbinom{6-4}{4}+\dbinom{6-5}{5}+\dbinom{6-6}{6}\\
&  =\dbinom{6}{0}+\dbinom{5}{1}+\dbinom{4}{2}+\dbinom{3}{3}+\dbinom{2}%
{4}+\dbinom{1}{5}+\dbinom{0}{6}\\
&  =1+5+6+1+0+0+0,
\end{align*}
which is indeed true. Of course, the three summands that are $0$ could just as
well be excluded from the sum, and the sum $\sum_{k=0}^{n}\dbinom{n-k}{k}$ in
Corollary \ref{cor.fib.sum-binom} could be replaced by the smaller sum
$\sum_{k=0}^{\left\lfloor n/2\right\rfloor }\dbinom{n-k}{k}$ (since
$\dbinom{n-k}{k}=0$ whenever $\left\lfloor n/2\right\rfloor <k\leq n$); but I
find it more important to keep the sum simple than to minimize the number of
its addends.
\end{example}

\begin{proof}
[Proof of Corollary \ref{cor.fib.sum-binom}.]It is easy to see that any subset
of $\left[  n-1\right]  $ has a size between $0$ and $n$
(inclusive)\footnote{\textit{Proof.} Let $T$ be a subset of $\left[
n-1\right]  $. We must show that $T$ has a size between $0$ and $n$
(inclusive). In other words, we must prove that $\left\vert T\right\vert
\in\left\{  0,1,\ldots,n\right\}  $.
\par
However, we have $T\subseteq\left[  n-1\right]  \subseteq\left[  n\right]  $
and therefore $\left\vert T\right\vert \leq\left\vert \left[  n\right]
\right\vert $ (by Theorem \ref{thm.count.rules.diff} \textbf{(a)}, applied to
$S=\left[  n\right]  $). Hence, $\left\vert T\right\vert \leq\left\vert
\left[  n\right]  \right\vert =n$. Since $\left\vert T\right\vert $ is a
nonnegative integer, we thus obtain $\left\vert T\right\vert \in\left\{
0,1,\ldots,n\right\}  $, as desired.}. (Actually, it cannot have size $n$
unless $n=0$, but I find it more convenient to nevertheless include the
\textquotedblleft unnecessary\textquotedblright\ value $n$ among the
theoretically possible sizes; I am not saying that all of these sizes actually
are achievable.)

Now, from $n\in\mathbb{N}$, we obtain $n\geq0$, thus $n-1\geq-1$. Hence,
Theorem \ref{thm.count.lacs.all} (applied to $n-1$ instead of $n$) yields%
\[
\left(  \text{\# of lacunar subsets of }\left[  n-1\right]  \right)
=f_{\left(  n-1\right)  +2}=f_{n+1}.
\]
Therefore,%
\begin{align*}
f_{n+1}  &  =\left(  \text{\# of lacunar subsets of }\left[  n-1\right]
\right) \\
&  =\left(  \text{\# of lacunar subsets of }\left[  n-1\right]  \text{ having
size }0\right) \\
&  \ \ \ \ \ \ \ \ \ \ +\left(  \text{\# of lacunar subsets of }\left[
n-1\right]  \text{ having size }1\right) \\
&  \ \ \ \ \ \ \ \ \ \ +\left(  \text{\# of lacunar subsets of }\left[
n-1\right]  \text{ having size }2\right) \\
&  \ \ \ \ \ \ \ \ \ \ +\cdots\\
&  \ \ \ \ \ \ \ \ \ \ +\left(  \text{\# of lacunar subsets of }\left[
n-1\right]  \text{ having size }n\right) \\
&  \ \ \ \ \ \ \ \ \ \ \ \ \ \ \ \ \ \ \ \ \left(
\begin{array}
[c]{c}%
\text{by the sum rule (Theorem \ref{thm.count.rules.sum-k}), since any}\\
\text{subset of }\left[  n-1\right]  \text{ has a size between }0\text{ and
}n\text{ (inclusive)}%
\end{array}
\right) \\
&  =\sum_{k=0}^{n}\underbrace{\left(  \text{\# of lacunar subsets of }\left[
n-1\right]  \text{ having size }k\right)  }_{\substack{=\left(  \text{\# of
}k\text{-element lacunar subsets of }\left[  n-1\right]  \right)
=\dbinom{\left(  n-1\right)  +1-k}{k}\\\text{(by Theorem
\ref{thm.count.lacs.size}, applied to }n-1\text{ instead of }n\\\text{(since
}k\leq n=\left(  n-1\right)  +1\text{))}}}\\
&  =\sum_{k=0}^{n}\dbinom{\left(  n-1\right)  +1-k}{k}=\sum_{k=0}^{n}%
\dbinom{n-k}{k}\ \ \ \ \ \ \ \ \ \ \left(  \text{since }\left(  n-1\right)
+1=n\right) \\
&  =\dbinom{n-0}{0}+\dbinom{n-1}{1}+\cdots+\dbinom{n-n}{n}.
\end{align*}
This proves Corollary \ref{cor.fib.sum-binom}.
\end{proof}

\subsubsection{\label{subsec.enum.lac.domino}The domino tilings connection}

At the beginning of Chapter \ref{chp.enum}, I asked for the \# of ways to tile
a $2\times15$-rectangle with dominos (i.e., rectangles of size $1\times2$ or
$2\times1$), such as the following:%
\[
\begin{tikzpicture}[scale=0.8]
\draw (0, 0) rectangle (1, 2);
\draw (1, 0) rectangle (3, 1);
\draw (1, 1) rectangle (3, 2);
\draw (3, 0) rectangle (4, 2);
\draw (4, 0) rectangle (5, 2);
\draw (5, 0) rectangle (7, 1);
\draw (5, 1) rectangle (7, 2);
\draw (7, 0) rectangle (9, 1);
\draw (7, 1) rectangle (9, 2);
\draw (9, 0) rectangle (10, 2);
\draw (10, 0) rectangle (12, 1);
\draw (10, 1) rectangle (12, 2);
\draw (12, 0) rectangle (13, 2);
\draw (13, 0) rectangle (14, 2);
\draw (14, 0) rectangle (15, 2);
\end{tikzpicture}\ \ .
\]
Of course, the same problem can be asked for $n\times m$-rectangles for
arbitrary $n$ and $m$, but we shall focus on the case $n=2$ (that is, a
rectangle of height $2$). (See \cite[\S 1.1]{19fco} for some references on the
much harder cases when $n>2$.)

It turns out that the ways to tile a $2\times m$-rectangle with dominos are in
bijection with the lacunar subsets of $\left[  m-1\right]  $. Indeed, if
$\mathcal{T}$ is a way to tile the $2\times m$-rectangle, then we let
$C\left(  \mathcal{T}\right)  $ be the set of all columns (counted from the
left) in which horizontal dominos of $\mathcal{T}$ start (where we say that a
\textbf{horizontal domino} is a domino of height $1$ and width $2$, and it
\textbf{starts} in the leftmost of the two columns that it spans). For
example, if $\mathcal{T}$ is the tiling shown above, then $C\left(
\mathcal{T}\right)  =\left\{  2,6,8,11\right\}  $. Now, it is not hard to see
(but not completely obvious; see \cite[\S 1.4.4, Second proof of Proposition
1.4.9]{19fco}) that the map%
\begin{align*}
\left\{  \text{ways to tile a }2\times m\text{-rectangle with dominos}%
\right\}   &  \rightarrow\left\{  \text{lacunar subsets of }\left[
m-1\right]  \right\}  ,\\
\mathcal{T}  &  \mapsto C\left(  \mathcal{T}\right)
\end{align*}
is a bijection, and therefore the bijection principle yields%
\begin{align*}
&  \left(  \text{\# of ways to tile a }2\times m\text{-rectangle with
dominos}\right) \\
&  =\left(  \text{\# of lacunar subsets of }\left[  m-1\right]  \right)
=f_{m+1}%
\end{align*}
(by Theorem \ref{thm.count.lacs.all}, applied to $n=m-1$). In particular, for
$m=15$, we obtain%
\[
\left(  \text{\# of ways to tile a }2\times15\text{-rectangle with
dominos}\right)  =f_{15+1}=f_{16}=987.
\]

\begin{exercise}
\label{exe.23hw5.5}A set $S$ of integers will be called \textbf{pseudolacunar}
if it has the property that no two elements $s,t$ of $S$ satisfy $\left\vert
s-t\right\vert =2$. For instance, the set $\left\{  2,5,6\right\}  $ is
pseudolacunar, but the set $\left\{  2,5,7\right\}  $ is not (since
$\left\vert 5-7\right\vert =2$).

For each $n\in\mathbb{N}$, let $p_{n}$ be the \# of pseudolacunar subsets of
$\left[  n\right]  $.

Prove that
\[
p_{n}=p_{n-1}+p_{n-3}+p_{n-4}\ \ \ \ \ \ \ \ \ \ \text{for each }n\geq4.
\]

[\textbf{Hint:} To each pseudolacunar subset, assign one of three colors.]
\end{exercise}

\begin{exercise}
\label{exe.count.subsets.self-start}A set $S$ of integers shall be called
\textbf{self-starting} if its size $\left\vert S\right\vert $ is also its
smallest element. (For example, $\left\{  3,5,6\right\}  $ is self-starting,
while $\left\{  2,3,4\right\}  $ and $\left\{  3\right\}  $ are not.)

Let $n\in\mathbb{N}$. \medskip

\textbf{(a)} For any $k\in\left[  n\right]  $, find the number of
self-starting subsets of $\left[  n\right]  $ having size $k$. \medskip

\textbf{(b)} Find the number of all self-starting subsets of $\left[
n\right]  $.
\end{exercise}

\subsection{Compositions and weak compositions}

Two other useful objects to count are \textbf{compositions} and \textbf{weak
compositions}.

\subsubsection{Compositions}

How many ways are there to write the integer $5$ as a sum of $3$ positive
integers, if the order matters? Since $5$ and $3$ are not very large numbers,
we can just list all these ways:%
\begin{align*}
5  &  =2+2+1=2+1+2=1+2+2\\
&  =3+1+1=1+3+1=1+1+3.
\end{align*}
So there are $6$ such ways.

What if we replace $5$ and $3$ by arbitrary nonnegative integers $n$ and $k$ ?
So we want to count the $k$-tuples $\left(  a_{1},a_{2},\ldots,a_{k}\right)  $
of positive integers satisfying $a_{1}+a_{2}+\cdots+a_{k}=n$. These tuples
have a name:

\begin{definition}
\label{def.count.comps.comps}\textbf{(a)} If $n\in\mathbb{N}$, then a
\textbf{composition of }$n$ shall mean a tuple (i.e., finite list) of positive
integers whose sum is $n$. \medskip

\textbf{(b)} If $n,k\in\mathbb{N}$, then a \textbf{composition of }$n$\textbf{
into }$k$\textbf{ parts} shall mean a $k$-tuple of positive integers whose sum
is $n$. \medskip

(The word \textquotedblleft composition\textquotedblright\ here is completely
unrelated to the notion of composition of two functions.)
\end{definition}

\begin{example}
\textbf{(a)} The compositions of $5$ into $3$ parts are%
\begin{align*}
&  \left(  2,2,1\right)  ,\ \ \ \ \ \ \ \ \ \ \left(  2,1,2\right)
,\ \ \ \ \ \ \ \ \ \ \left(  1,2,2\right)  ,\\
&  \left(  3,1,1\right)  ,\ \ \ \ \ \ \ \ \ \ \left(  1,3,1\right)
,\ \ \ \ \ \ \ \ \ \ \left(  1,1,3\right)  .
\end{align*}
These are exactly the $6$ ways we found above (but written as $3$-tuples).
\medskip

\textbf{(b)} The compositions of $3$ are%
\[
\left(  1,1,1\right)  ,\ \ \ \ \ \ \ \ \ \ \left(  2,1\right)
,\ \ \ \ \ \ \ \ \ \ \left(  1,2\right)  ,\ \ \ \ \ \ \ \ \ \ \left(
3\right)  .
\]

\textbf{(c)} The only composition of $0$ is the empty list $\left(  {}\right)
$, which is a $0$-tuple. It is a composition into $0$ parts.
\end{example}

Let us now count compositions of $n$ into $k$ parts. (Later, we will count all
compositions of $n$.) Again, the answer turns out to be a binomial coefficient:

\begin{theorem}
\label{thm.count.comps.comps-nk}Let $n,k\in\mathbb{N}$. Then,%
\begin{equation}
\left(  \text{\# of compositions of }n\text{ into }k\text{ parts}\right)
=\dbinom{n-1}{n-k}. \label{eq.thm.count.comps.comps-nk.n-k}%
\end{equation}
If $n>0$, then we furthermore have%
\begin{equation}
\left(  \text{\# of compositions of }n\text{ into }k\text{ parts}\right)
=\dbinom{n-1}{k-1}. \label{eq.thm.count.comps.comps-nk.k-1}%
\end{equation}

\end{theorem}

\begin{proof}
[Proof sketch.]The proof is straightforward in the case when $n=0$. (Indeed,
if $n=0$, then the only composition of $n$ is the empty list $\left(
{}\right)  $, and this is a composition of $n$ into $0$ parts. Thus, if $n=0$,
then we have
\[
\left(  \text{\# of compositions of }n\text{ into }k\text{ parts}\right)  =%
\begin{cases}
1, & \text{if }k=0;\\
0, & \text{if }k\neq0;
\end{cases}
\]
but we also have%
\[
\dbinom{n-1}{n-k}=\dbinom{0-1}{0-k}=\dbinom{-1}{-k}=%
\begin{cases}
1, & \text{if }k=0;\\
0, & \text{if }k\neq0
\end{cases}
\ \ \ \ \ \ \ \ \ \ \left(  \text{check this!}\right)
\]
in this case, and we obtain (\ref{eq.thm.count.comps.comps-nk.n-k}) by
comparing these two equalities. Thus, Theorem \ref{thm.count.comps.comps-nk}
holds for $n=0$ (because the equality (\ref{eq.thm.count.comps.comps-nk.k-1})
is claimed for $n>0$ only).)

Thus, we only need to consider the case when $n\neq0$. Let us thus focus on
this case. From $n\neq0$, we obtain $n\geq1$ (since $n\in\mathbb{N}$), thus
$n-1\in\mathbb{N}$. Furthermore, for any composition $a=\left(  a_{1}%
,a_{2},\ldots,a_{k}\right)  $ of $n$ into $k$ parts, we must have $k\geq1$
(since the sum of the entries of this composition is $n\neq0$, and therefore
it must have at least one entry), therefore $k-1\in\mathbb{N}$.

For any composition $a=\left(  a_{1},a_{2},\ldots,a_{k}\right)  $ of $n$ into
$k$ parts, we define the \textbf{partial sum set} $C\left(  a\right)  $ to be
the set%
\begin{align*}
&  \ \left\{  a_{1},\ \ a_{1}+a_{2},\ \ a_{1}+a_{2}+a_{3},\ \ \ldots
,\ \ a_{1}+a_{2}+\cdots+a_{k-1}\right\} \\
=  &  \ \left\{  a_{1}+a_{2}+\cdots+a_{i}\ \mid\ i\in\left[  k-1\right]
\right\}  .
\end{align*}
This set $C\left(  a\right)  $ consists of all the \textquotedblleft partial
sums\textquotedblright\ $a_{1}+a_{2}+\cdots+a_{i}$ of the sum $a_{1}%
+a_{2}+\cdots+a_{k}$, except for the empty partial sum $a_{1}+a_{2}%
+\cdots+a_{0}$ (which is $0$ by definition) and the full sum $a_{1}%
+a_{2}+\cdots+a_{k}$ (which is $n$, since $a$ is a composition of $n$). Thus,
all elements of $C\left(  a\right)  $ are integers between $0$ and $n$
(exclusive) (since they have more addends than the empty partial sum, but
fewer than the full sum\footnote{and since all these addends are positive
(because a composition has positive entries)}). In other words, $C\left(
a\right)  $ is a subset of $\left\{  1,2,\ldots,n-1\right\}  =\left[
n-1\right]  $.

We can visualize the partial sum set $C\left(  a\right)  $ of a composition
$a=\left(  a_{1},a_{2},\ldots,a_{k}\right)  $ as follows: The interval
$\left[  0,n\right]  _{\mathbb{R}}:=\left\{  x\in\mathbb{R}\ \mid\ 0\leq x\leq
n\right\}  $ on the real line has length $n$. If we split this interval into
blocks of lengths $a_{1},a_{2},\ldots,a_{k}$ (from left to right), then the
elements of $C\left(  a\right)  $ are precisely the endpoints of these blocks
(i.e., the points at which one block ends and the next begins), except for the
leftmost endpoint $0$ and the rightmost endpoint $n$. See this picture:%
\[%
\begin{tikzpicture}[scale=1.5]
\draw[->, gray] (-0.3,0) -- (7.3,0);
\node(0) at (0 cm, 0 cm) {};
\node(3) at (3 cm, 0 cm) {};
\node(4) at (4 cm, 0 cm) {};
\node(6) at (6 cm, 0 cm) {};
\node(7) at (7 cm, 0 cm) {};
\node[blue] (a) at (1.5 cm, 0.5 cm) {$a_1$};
\node[blue] (b) at (3.5 cm, 0.5 cm) {$a_2$};
\node[blue] (c) at (5.0 cm, 0.5 cm) {$\cdots$};
\node[blue] (d) at (6.5 cm, 0.5 cm) {$a_k$};
\fill(0) circle (1pt) node[below] {$0$};
\fill(3) circle (1pt) node[below] {$s_1$};
\fill(4) circle (1pt) node[below] {$s_2$};
\node(u) at (5 cm, -0.2 cm) {$\cdots$};
\fill(6) circle (1pt) node[below] {$s_{k-1}$};
\fill(7) circle (1pt) node[below] {$n$};
\draw[->, blue] (0.1,0.5) -- (1.35,0.5);
\draw[->, blue] (2.9,0.5) -- (1.65,0.5);
\draw[->, blue] (3.1,0.5) -- (3.35,0.5);
\draw[->, blue] (3.9,0.5) -- (3.65,0.5);
\draw[->, blue] (4.1,0.5) -- (4.75,0.5);
\draw[->, blue] (5.9,0.5) -- (5.25,0.5);
\draw[->, blue] (6.1,0.5) -- (6.35,0.5);
\draw[->, blue] (6.9,0.5) -- (6.65,0.5);
\draw[very thick, red] (0,0) -- (0,1);
\draw[very thick, red] (3,0) -- (3,1);
\draw[very thick, red] (4,0) -- (4,1);
\draw[very thick, red] (6,0) -- (6,1);
\draw[very thick, red] (7,0) -- (7,1);
\end{tikzpicture}%
\]
(on which the partial sums $a_{1}+a_{2}+\cdots+a_{i}$ are denoted by $s_{i}$).

It is thus easy to see that if $a$ is a composition of $n$ into $k$ parts,
then $C\left(  a\right)  $ is a $\left(  k-1\right)  $-element subset of
$\left[  n-1\right]  $. Thus, we obtain a map%
\begin{align*}
C:\left\{  \text{compositions of }n\text{ into }k\text{ parts}\right\}   &
\rightarrow\left\{  \left(  k-1\right)  \text{-element subsets of }\left[
n-1\right]  \right\}  ,\\
a  &  \mapsto C\left(  a\right)  .
\end{align*}
Furthermore, it is not hard to see that this map $C$ has an
inverse\footnote{This is easiest to see using the visual description of
$C\left(  a\right)  $ that we showed above: Given a $\left(  k-1\right)
$-element subset $I$ of $\left[  n-1\right]  $, we can use the elements of $I$
to subdivide the interval $\left[  0,n\right]  _{\mathbb{R}}$ into $k$ blocks.
The lengths of these blocks (listed from left to right) form a composition $a$
of $n$ into $k$ parts, and this composition satisfies $C\left(  a\right)  =I$.
Moreover, this composition is the only one with this property. Thus, the map
that sends each $\left(  k-1\right)  $-element subset $I$ of $\left[
n-1\right]  $ to the corresponding composition $a$ (whose construction we just
explained) is an inverse map of $C$.
\par
Rigorously, this can be restated as follows: For each $\left(  k-1\right)
$-element subset $I=\left\{  i_{1}<i_{2}<\cdots<i_{k-1}\right\}  $ of $\left[
n-1\right]  $ (where we are using Convention
\ref{conv.intset.inclist-setbrace} again), we can define a composition%
\[
A\left(  I\right)  :=\left(  i_{1}-i_{0},\ i_{2}-i_{1},\ i_{3}-i_{2}%
,\ \ldots,\ i_{k-1}-i_{k-2},\ i_{k}-i_{k-1}\right)  ,
\]
where we set $i_{0}:=0$ and $i_{k}:=n$. Then, the map%
\begin{align*}
A:\left\{  \left(  k-1\right)  \text{-element subsets of }\left[  n-1\right]
\right\}   &  \rightarrow\left\{  \text{compositions of }n\text{ into }k\text{
parts}\right\}  ,\\
I  &  \mapsto A\left(  I\right)
\end{align*}
is easily seen to be an inverse map of $C$. A detailed proof can be found in
\cite[solution to Exercise 1 \textbf{(b)}]{19f-hw0s} (except that the latter
solution does not pay attention to the size of the subset).}, and thus is a
bijection. Hence, the bijection principle yields%
\begin{align*}
&  \left(  \text{\# of compositions of }n\text{ into }k\text{ parts}\right) \\
&  =\left(  \text{\# of }\left(  k-1\right)  \text{-element subsets of
}\left[  n-1\right]  \right) \\
&  =\dbinom{n-1}{k-1}\ \ \ \ \ \ \ \ \ \ \left(
\begin{array}
[c]{c}%
\text{by Theorem \ref{thm.count.subsets.size.S-repeat}}\\
\text{(applied to }k-1\text{, }n-1\text{ and }\left[  n-1\right]  \text{
instead of }k\text{, }n\text{ and }S\text{)}%
\end{array}
\right) \\
&  =\dbinom{n-1}{\left(  n-1\right)  -\left(  k-1\right)  }%
\ \ \ \ \ \ \ \ \ \ \left(
\begin{array}
[c]{c}%
\text{by the symmetry of Pascal's triangle}\\
\text{(Theorem \ref{thm.binom.sym}), since }n-1\in\mathbb{N}%
\end{array}
\right) \\
&  =\dbinom{n-1}{n-k}\ \ \ \ \ \ \ \ \ \ \left(  \text{since }\left(
n-1\right)  -\left(  k-1\right)  =n-k\right)  .
\end{align*}
Thus, both (\ref{eq.thm.count.comps.comps-nk.n-k}) and
(\ref{eq.thm.count.comps.comps-nk.k-1}) have been proved. This completes the
proof of Theorem \ref{thm.count.comps.comps-nk}.
\end{proof}

We can also count all compositions of a given $n$:

\begin{theorem}
\label{thm.count.comps.comps-n}Let $n$ be a positive integer. Then, the \# of
all compositions of $n$ is $2^{n-1}$.
\end{theorem}

\begin{proof}
[Proof sketch.]This can be proved using a similar argument as in Theorem
\ref{thm.count.comps.comps-nk} (but now we need to count all subsets of
$\left[  n-1\right]  $). See \cite[Exercise 1 \textbf{(b)}]{19f-hw0s} for details.
\end{proof}

Note that Theorem \ref{thm.count.comps.comps-n} does not hold for $n=0$ (since
$0$ has $1$ composition, but $2^{0-1}=\dfrac{1}{2}$).

\subsubsection{Weak compositions}

One particularly useful variant of compositions are the so-called \textbf{weak
compositions}. These are defined as tuples of nonnegative integers (i.e., they
differ from compositions in that their entries are allowed to be $0$). In
other words:

\begin{definition}
\label{def.count.wcomps.wcomps}\textbf{(a)} If $n\in\mathbb{N}$, then a
\textbf{weak composition of }$n$ shall mean a tuple of nonnegative integers
whose sum is $n$. \medskip

\textbf{(b)} If $n,k\in\mathbb{N}$, then a \textbf{weak composition of }%
$n$\textbf{ into }$k$\textbf{ parts} shall mean a $k$-tuple of nonnegative
integers whose sum is $n$.
\end{definition}

For instance:

\begin{itemize}
\item The weak compositions of $2$ into $3$ parts are%
\begin{align*}
&  \left(  1,1,0\right)  ,\ \ \ \ \ \ \ \ \ \ \left(  1,0,1\right)
,\ \ \ \ \ \ \ \ \ \ \left(  0,1,1\right)  ,\\
&  \left(  2,0,0\right)  ,\ \ \ \ \ \ \ \ \ \ \left(  0,2,0\right)
,\ \ \ \ \ \ \ \ \ \ \left(  0,0,2\right)  .
\end{align*}

\item The weak compositions of $2$ into $2$ parts are%
\[
\left(  2,0\right)  ,\ \ \ \ \ \ \ \ \ \ \left(  1,1\right)
,\ \ \ \ \ \ \ \ \ \ \left(  0,2\right)  .
\]
(Note that any composition is a weak composition, but there are usually more
weak compositions than that.)

\item The weak compositions of $1$ are all tuples of the form
\[
\left(  \underbrace{0,0,\ldots,0}_{\text{any number of zeroes}}%
,1,\underbrace{0,0,\ldots,0}_{\text{any number of zeroes}}\right)  .
\]
Here, \textquotedblleft any number\textquotedblright\ allows for the
possibility of \textquotedblleft none\textquotedblright, and in particular the
$1$-tuple $\left(  1\right)  $ is a weak composition of $1$.
\end{itemize}

Counting all weak compositions of a given $n$ is no longer possible, since
there are infinitely many (as we just saw). But we can still count all weak
compositions of a given $n$ into $k$ parts for a given $k$.

\begin{theorem}
\label{thm.count.wcomps.wcomps-nk}Let $n,k\in\mathbb{N}$. Then,%
\[
\left(  \text{\# of weak compositions of }n\text{ into }k\text{ parts}\right)
=\dbinom{n+k-1}{n}.
\]

Moreover, if $n+k>0$ (that is, if $n$ and $k$ are not both $0$), then%
\[
\left(  \text{\# of weak compositions of }n\text{ into }k\text{ parts}\right)
=\dbinom{n+k-1}{k-1}.
\]

\end{theorem}

\begin{proof}
We shall deduce this from Theorem \ref{thm.count.comps.comps-nk}.

Indeed, if $b$ is a nonnegative integer, then $b+1$ is a positive integer.
Thus, if $\left(  a_{1},a_{2},\ldots,a_{k}\right)  $ is a weak composition of
$n$ into $k$ parts, then the $k$-tuple $\left(  a_{1}+1,\ a_{2}+1,\ \ldots
,\ a_{k}+1\right)  $ is a composition of $n+k$ into $k$ parts (since the sum
of its entries is%
\begin{align*}
\left(  a_{1}+1\right)  +\left(  a_{2}+1\right)  +\cdots+\left(
a_{k}+1\right)   &  =\underbrace{\left(  a_{1}+a_{2}+\cdots+a_{k}\right)
}_{\substack{=n\\\text{(since }\left(  a_{1},a_{2},\ldots,a_{k}\right)  \text{
is a weak composition of }n\text{)}}}+\,k\\
&  =n+k
\end{align*}
). Thus, the map%
\begin{align*}
\left\{  \text{weak compositions of }n\text{ into }k\text{ parts}\right\}   &
\rightarrow\left\{  \text{compositions of }n+k\text{ into }k\text{
parts}\right\}  ,\\
\left(  a_{1},a_{2},\ldots,a_{k}\right)   &  \mapsto\left(  a_{1}%
+1,\ a_{2}+1,\ \ldots,\ a_{k}+1\right)
\end{align*}
is well-defined. Similarly, the map
\begin{align*}
\left\{  \text{compositions of }n+k\text{ into }k\text{ parts}\right\}   &
\rightarrow\left\{  \text{weak compositions of }n\text{ into }k\text{
parts}\right\}  ,\\
\left(  a_{1},a_{2},\ldots,a_{k}\right)   &  \mapsto\left(  a_{1}%
-1,\ a_{2}-1,\ \ldots,\ a_{k}-1\right)
\end{align*}
is well-defined. These two maps are clearly inverses of each other (since
adding $1$ and subtracting $1$ are inverse operations). Therefore, they are
bijections. The bijection principle thus yields%
\begin{align*}
&  \left(  \text{\# of weak compositions of }n\text{ into }k\text{
parts}\right) \\
&  =\left(  \text{\# of compositions of }n+k\text{ into }k\text{ parts}\right)
\\
&  =\dbinom{n+k-1}{n+k-k}\ \ \ \ \ \ \ \ \ \ \left(  \text{by
(\ref{eq.thm.count.comps.comps-nk.n-k}), applied to }n+k\text{ instead of
}n\right) \\
&  =\dbinom{n+k-1}{n}.
\end{align*}
If $n+k>0$, then $n+k\geq1$ and thus $n+k-1\in\mathbb{N}$, so that this
becomes%
\begin{align*}
&  \left(  \text{\# of weak compositions of }n\text{ into }k\text{
parts}\right) \\
&  =\dbinom{n+k-1}{n}\\
&  =\dbinom{n+k-1}{\left(  n+k-1\right)  -n}\ \ \ \ \ \ \ \ \ \ \left(
\begin{array}
[c]{c}%
\text{by the symmetry of Pascal's triangle}\\
\text{(Theorem \ref{thm.binom.sym}), since }n+k-1\in\mathbb{N}%
\end{array}
\right) \\
&  =\dbinom{n+k-1}{k-1}\ \ \ \ \ \ \ \ \ \ \text{in this case.}%
\end{align*}
Thus, Theorem \ref{thm.count.wcomps.wcomps-nk} is fully proved.
\end{proof}

\begin{exercise}
\label{exe.binom.compos.12}Let $n\in\mathbb{N}$.

A $\left\{  1,2\right\}  $\textbf{-composition of} $n$ shall mean a
composition $\left(  a_{1},a_{2},\ldots,a_{k}\right)  $ of $n$ such that
$a_{1},a_{2},\ldots,a_{k}\in\left\{  1,2\right\}  $.

For example, the $\left\{  1,2\right\}  $-compositions of $5$ are%
\begin{align*}
&  \left(  1,1,1,1,1\right)  ,\ \ \ \ \ \ \ \ \ \ \left(  1,1,1,2\right)
,\ \ \ \ \ \ \ \ \ \ \left(  1,1,2,1\right)  ,\ \ \ \ \ \ \ \ \ \ \left(
1,2,1,1\right)  ,\\
&  \left(  2,1,1,1\right)  ,\ \ \ \ \ \ \ \ \ \ \left(  2,2,1\right)
,\ \ \ \ \ \ \ \ \ \ \left(  2,1,2\right)  ,\ \ \ \ \ \ \ \ \ \ \left(
1,2,2\right)  .
\end{align*}

\textbf{(a)} Prove that%
\[
\left(  \text{\# of }\left\{  1,2\right\}  \text{-compositions of }n\right)
=f_{n+1}%
\]
(where $\left(  f_{0},f_{1},f_{2},\ldots\right)  $ denotes the Fibonacci
sequence, as defined in Definition \ref{def.ind.fib}). \medskip

\textbf{(b)} Let $k\in\mathbb{N}$. A $\left\{  1,2\right\}  $%
\textbf{-composition of }$n$\textbf{ into }$k$\textbf{ parts} shall mean a
composition $\left(  a_{1},a_{2},\ldots,a_{k}\right)  $ of $n$ into $k$ parts
such that $a_{1},a_{2},\ldots,a_{k}\in\left\{  1,2\right\}  $.

Prove that
\[
\left(  \text{\# of }\left\{  1,2\right\}  \text{-compositions of }n\text{
into }k\text{ parts}\right)  =\dbinom{k}{n-k}.
\]

\end{exercise}

\subsection{Selections}

We now come back to a class of problems that we have posed at the start of
Chapter \ref{chp.enum} (before Section \ref{sec.enum.sets}) but haven't fully
answered yet: counting the ways to select a bunch of elements from a given set.

To be more specific, these are problems of the following form: Given an
$n$-element set $S$, how many ways are there to select $k$ elements from $S$
(where $n$ and $k$ are fixed nonnegative integers)?

The words \textquotedblleft$k$ elements\textquotedblright\ in this question
are ambiguous, as they allow for several interpretations:

\begin{enumerate}
\item Do we want $k$ arbitrary elements or $k$ distinct elements?

\item Does the order of these $k$ elements matter or not? (In other words,
would \textquotedblleft$1,2$\textquotedblright\ and \textquotedblleft%
$2,1$\textquotedblright\ count as two different selections?)
\end{enumerate}

In total, these decisions leave you with $4$ options, leading to $4$ different
problems. In this section, we shall address them all.

\subsubsection{Unordered selections without repetition (= without
replacement)}

Let us begin with the case when we want to select $k$ distinct elements, and
the order does not matter. This just means selecting a $k$-element subset of
$S$. We already know how to count these subsets (Theorem
\ref{thm.count.subsets.size.S-repeat}):

\begin{theorem}
\label{thm.count.select.subset}Let $n\in\mathbb{N}$, and let $k$ be any
number. Let $S$ be an $n$-element set. Then,%
\[
\left(  \text{\# of }k\text{-element subsets of }S\right)  =\dbinom{n}{k}.
\]

\end{theorem}

In other words, the \# of ways to choose $k$ distinct elements from a given
$n$-element set $S$, if the order does not matter, is $\dbinom{n}{k}$.

\subsubsection{Ordered selections without repetition (= without replacement)}

Now, let us consider the case when the order does matter. Thus, we are looking
not for subsets, but for $k$-tuples. But these $k$-tuples are not arbitrary
$k$-tuples; they are $k$-tuples of \textbf{distinct} elements. We shall call
such $k$-tuples \textbf{injective} (in analogy to injective functions):

\begin{definition}
\label{def.count.select.inj-k-tuple}Let $k\in\mathbb{N}$. A $k$-tuple $\left(
i_{1},i_{2},\ldots,i_{k}\right)  $ is said to be \textbf{injective} if its $k$
entries $i_{1},i_{2},\ldots,i_{k}$ are distinct (i.e., if we have $i_{a}\neq
i_{b}$ for all $a\neq b$).
\end{definition}

For example, the $3$-tuple $\left(  6,1,2\right)  $ is injective, but $\left(
2,1,2\right)  $ is not.

Note that injective $k$-tuples and injective functions are closely related: A
function $f:\left[  k\right]  \rightarrow S$ (for a set $S$ and a number
$k\in\mathbb{N}$) is injective if and only if the $k$-tuple $\left(  f\left(
1\right)  ,\ f\left(  2\right)  ,\ \ldots,\ f\left(  k\right)  \right)  $ is injective.

Next, we introduce another convenient notation:

\begin{definition}
\label{def.count.cartesian-power}Let $S$ be any set, and let $k\in\mathbb{N}$.
Then, $S^{k}$ shall mean the Cartesian product%
\begin{align*}
\underbrace{S\times S\times\cdots\times S}_{k\text{ times}}  &  =\left\{
\left(  a_{1},a_{2},\ldots,a_{k}\right)  \ \mid\ a_{1},a_{2},\ldots,a_{k}\in
S\right\} \\
&  =\left\{  k\text{-tuples whose all entries belong to }S\right\}  .
\end{align*}

\end{definition}

For example, $\left\{  5,6\right\}  ^{3}$ is the set
\begin{align*}
&  \left\{  5,6\right\}  \times\left\{  5,6\right\}  \times\left\{
5,6\right\} \\
&  =\left\{  \left(  5,5,5\right)  ,\ \left(  5,5,6\right)  ,\ \left(
5,6,5\right)  ,\ \left(  5,6,6\right)  ,\ \left(  6,5,5\right)  ,\ \left(
6,5,6\right)  ,\ \left(  6,6,5\right)  ,\ \left(  6,6,6\right)  \right\}  .
\end{align*}
None of the $3$-tuples (i.e., triples) in this set is injective, but it is
easy to find an example where injective $k$-tuples do appear: For instance,
the set $\left\{  1,2,3,4\right\}  ^{3}$ contains both injective $3$-tuples
such as $\left(  1,4,3\right)  $ and non-injective $3$-tuples such as $\left(
3,3,1\right)  $.

Now, we can define rigorously what we are looking for: A way to select $k$
distinct elements from a given set $S$, if the order matters, is the same as
an injective $k$-tuple in $S^{k}$. We shall now count such ways:

\begin{theorem}
\label{thm.count.select.injtup}Let $n\in\mathbb{N}$ and $k\in\mathbb{N}$. Let
$S$ be an $n$-element set. Then,%
\[
\left(  \text{\# of injective }k\text{-tuples in }S^{k}\right)  =n\left(
n-1\right)  \left(  n-2\right)  \cdots\left(  n-k+1\right)  .
\]

\end{theorem}

\begin{example}
Applying Theorem \ref{thm.count.select.injtup} to $n=5$, $k=3$ and $S=\left\{
1,2,3,4,5\right\}  $, we find that%
\[
\left(  \text{\# of injective }3\text{-tuples in }\left\{  1,2,3,4,5\right\}
^{3}\right)  =5\left(  5-1\right)  \left(  5-2\right)  =5\cdot4\cdot3=60.
\]
And indeed, there are $60$ injective $3$-tuples in $\left\{
1,2,3,4,5\right\}  ^{3}$. For example, $\left(  2,5,4\right)  $ and $\left(
5,3,2\right)  $ are two of them.
\end{example}

Note that the right hand side in Theorem \ref{thm.count.select.injtup} is
precisely the numerator in the definition of the binomial coefficient
$\dbinom{n}{k}$ (Definition \ref{def.binom}), and thus can be rewritten as
$k!\cdot\dbinom{n}{k}$ (since $k!$ is the denominator). Thus, the claim of
Theorem \ref{thm.count.select.injtup} can be restated as
\[
\left(  \text{\# of injective }k\text{-tuples in }S^{k}\right)  =k!\cdot
\dbinom{n}{k}.
\]

Now, how do we prove the theorem? Let us first give an informal proof:

\begin{proof}
[Informal proof of Theorem \ref{thm.count.select.injtup}.]Let us look at an
example (which is representative of the general case): We let $n=5$ and $k=3$
and $S=\left\{  a,b,c,d,e\right\}  $. How many injective $k$-tuples are there
in $S^{k}$ ? In other words (since $k=3$): How many injective $3$-tuples are
there in $S^{3}$ ?

Such a $3$-tuple has the form $\left(  x,y,z\right)  $, where $x,y,z$ are
three distinct elements of $S$. Let us see how such a $3$-tuple can be chosen:

\begin{enumerate}
\item First, we choose its first entry $x$. There are $5$ options for this,
since $S$ has $5$ elements (and $x$ can be any of these $5$).

\item Then, we choose its second entry $y$. There are $4$ options for it,
since $y$ can be any of the $5$ elements of $S$ except for $x$ (because the
injectivity of $\left(  x,y,z\right)  $ demands $y$ to be distinct from $x$).

\item Finally, we choose its third entry $z$. There are $3$ options for it,
since $z$ can be any of the $5$ elements of $S$ except for $x$ and $y$
(because the injectivity of $\left(  x,y,z\right)  $ demands $z$ to be
distinct from $x$ and $y$) and since $x$ and $y$ are already distinct.
\end{enumerate}

Altogether, we have $5$ options at the first step, then $4$ options at the
second step (no matter which option has been chosen at the first step), and
finally $3$ options at the third step. Altogether, we can therefore choose our
$3$-tuple in $5\cdot4\cdot3$ many different ways, because the numbers of
options multiply. Here, we have used a counting rule called \textquotedblleft%
\textbf{dependent product rule}\textquotedblright, which informally says that
if we perform a multi-step construction, and we have

\begin{itemize}
\item exactly $n_{1}$ options in step $1$,

\item exactly $n_{2}$ options in step $2$,

\item $\ldots$,

\item exactly $n_{k}$ options in step $k$,
\end{itemize}

\noindent then the entire construction can be performed in $n_{1}n_{2}\cdots
n_{k}$ many different ways. We shall not formalize this rule (let alone prove
it); the reader can find rigorous versions of this rule in \cite[\S 1.8]%
{Loehr-BC} and in \cite[Theorem 8.1.19]{Newste23}. However, we shall next give
a more rigorous proof of Theorem \ref{thm.count.select.injtup}, which uses
induction on $k$ instead of this \textquotedblleft dependent product
rule\textquotedblright\ (although the underlying idea is the same).
\end{proof}

\begin{proof}
[Rigorous proof of Theorem \ref{thm.count.select.injtup}.]Forget that we fixed
$S$ and $n$. We thus must prove the statement%
\[
P\left(  k\right)  :=\left(
\begin{array}
[c]{c}%
\text{\textquotedblleft for all }n\in\mathbb{N}\text{ and all }n\text{-element
sets }S\text{, we have}\\
\left(  \text{\# of injective }k\text{-tuples in }S^{k}\right)  =n\left(
n-1\right)  \left(  n-2\right)  \cdots\left(  n-k+1\right)
\text{\textquotedblright}%
\end{array}
\right)
\]
for each $k\in\mathbb{N}$. We shall prove this by induction on $k$.

\textit{Base case:} We must prove that $P\left(  0\right)  $ holds. In other
words, we must prove that for all $n\in\mathbb{N}$ and all $n$-element sets
$S$, we have%
\[
\left(  \text{\# of injective }0\text{-tuples in }S^{0}\right)  =n\left(
n-1\right)  \left(  n-2\right)  \cdots\left(  n-0+1\right)  .
\]
But this is an easy exercise in understanding emptiness: Let $n\in\mathbb{N}$,
and let $S$ be an $n$-element set. The only $0$-tuple in $S^{0}$ is $\left(
{}\right)  $, and this $0$-tuple is injective. Thus,%
\[
\left(  \text{\# of injective }0\text{-tuples in }S^{0}\right)  =1\text{.}%
\]
Comparing this with%
\[
n\left(  n-1\right)  \left(  n-2\right)  \cdots\left(  n-0+1\right)  =\left(
\text{empty product}\right)  =1,
\]
we obtain $\left(  \text{\# of injective }0\text{-tuples in }S^{0}\right)
=n\left(  n-1\right)  \left(  n-2\right)  \cdots\left(  n-0+1\right)  $. Thus,
$P\left(  0\right)  $ is proved. This completes the base case.

\textit{Induction step:} Let $k$ be a positive integer. Assume (as the
induction hypothesis) that $P\left(  k-1\right)  $ holds. Our goal is to prove
$P\left(  k\right)  $.

We have assumed that $P\left(  k-1\right)  $ holds. In other words, for all
$n\in\mathbb{N}$ and all $n$-element sets $S$, we have%
\begin{align}
&  \left(  \text{\# of injective }\left(  k-1\right)  \text{-tuples in
}S^{k-1}\right) \nonumber\\
&  =n\left(  n-1\right)  \left(  n-2\right)  \cdots\left(  n-\left(
k-1\right)  +1\right)  . \label{pf.thm.count.select.injtup.IH}%
\end{align}

Now, let us focus on proving $P\left(  k\right)  $. Thus, we fix an
$n\in\mathbb{N}$ and an $n$-element set $S$. Our goal is then to prove that%
\[
\left(  \text{\# of injective }k\text{-tuples in }S^{k}\right)  \overset{?}{=}%
n\left(  n-1\right)  \left(  n-2\right)  \cdots\left(  n-k+1\right)  .
\]
(Again, the question mark atop the equality sign reminds us that this is not
proved yet.)

Let $s_{1},s_{2},\ldots,s_{n}$ be the $n$ elements of $S$ (listed without
repetition). Then, any $k$-tuple in $S^{k}$ ends\footnote{We say that a
$k$-tuple \textbf{ends} with a given element $b$ if $b$ is the last entry of
this $k$-tuple. Note that every $k$-tuple does indeed have a last entry, since
$k$ is positive.} with exactly one of $s_{1},s_{2},\ldots,s_{n}$. Hence, by
the sum rule, we have%
\begin{align}
&  \left(  \text{\# of injective }k\text{-tuples in }S^{k}\right) \nonumber\\
&  =\left(  \text{\# of injective }k\text{-tuples in }S^{k}\text{ that end
with }s_{1}\right) \nonumber\\
&  \ \ \ \ \ \ \ \ \ \ +\left(  \text{\# of injective }k\text{-tuples in
}S^{k}\text{ that end with }s_{2}\right) \nonumber\\
&  \ \ \ \ \ \ \ \ \ \ +\cdots\nonumber\\
&  \ \ \ \ \ \ \ \ \ \ +\left(  \text{\# of injective }k\text{-tuples in
}S^{k}\text{ that end with }s_{n}\right) \nonumber\\
&  =\sum_{i=1}^{n}\left(  \text{\# of injective }k\text{-tuples in }%
S^{k}\text{ that end with }s_{i}\right)  .
\label{pf.thm.count.select.injtup.sum}%
\end{align}
Now, we shall compute the addends in this sum.

Fix any $i\in\left[  n\right]  $. An injective $k$-tuple in $S^{k}$ that ends
with $s_{i}$ must have the form%
\[
\left(  \ldots,s_{i}\right)  ,
\]
where the \textquotedblleft$\ldots$\textquotedblright\ are $k-1$ distinct
elements of $S\setminus\left\{  s_{i}\right\}  $ (not merely of $S$, but
actually of $S\setminus\left\{  s_{i}\right\}  $, because if any of them was
$s_{i}$, then our $k$-tuple would contain the entry $s_{i}$ twice and thus
fail to be injective). In other words, an injective $k$-tuple in $S^{k}$ that
ends with $s_{i}$ is an injective $\left(  k-1\right)  $-tuple in $\left(
S\setminus\left\{  s_{i}\right\}  \right)  ^{k-1}$ followed by the entry
$s_{i}$. Thus, we obtain a map%
\begin{align*}
\left\{  \text{injective }k\text{-tuples in }S^{k}\text{ that end with }%
s_{i}\right\}   &  \rightarrow\left\{  \text{injective }\left(  k-1\right)
\text{-tuples in }\left(  S\setminus\left\{  s_{i}\right\}  \right)
^{k-1}\right\}  ,\\
\left(  \ldots,s_{i}\right)   &  \mapsto\left(  \ldots\right)
\end{align*}
(which removes the last entry from our $k$-tuple and leaves the other entries
as they are)\footnote{For example, if $k=4$, then this map sends a $k$-tuple
$\left(  a,b,c,s_{i}\right)  $ to $\left(  a,b,c\right)  $.}. Conversely, we
have a map%
\begin{align*}
\left\{  \text{injective }\left(  k-1\right)  \text{-tuples in }\left(
S\setminus\left\{  s_{i}\right\}  \right)  ^{k-1}\right\}   &  \rightarrow
\left\{  \text{injective }k\text{-tuples in }S^{k}\text{ that end with }%
s_{i}\right\}  ,\\
\left(  \ldots\right)   &  \mapsto\left(  \ldots,s_{i}\right)
\end{align*}
(which inserts an $s_{i}$ after the end of a $\left(  k-1\right)  $-tuple; the
result is still injective\footnote{\textit{Proof.} We must show that if we
insert an $s_{i}$ after the end of an injective $\left(  k-1\right)  $-tuple
in $\left(  S\setminus\left\{  s_{i}\right\}  \right)  ^{k-1}$, then the
result is still injective.
\par
Indeed, the only way this could fail is if the newly inserted entry $s_{i}$
would already appear in the original $\left(  k-1\right)  $-tuple. However,
this is impossible, since the original $\left(  k-1\right)  $-tuple belongs to
$\left(  S\setminus\left\{  s_{i}\right\}  \right)  ^{k-1}$ and thus cannot
contain the entry $s_{i}$.})\footnote{For example, if $k=4$, then this map
sends a $\left(  k-1\right)  $-tuple $\left(  a,b,c\right)  $ to $\left(
a,b,c,s_{i}\right)  $.}. These two maps are clearly inverses of each
other\footnote{because a $k$-tuple that ends with $s_{i}$ stays unchanged if
we replace its last entry with $s_{i}$}, and thus are bijections. Hence, the
bijection principle yields%
\begin{align*}
&  \left(  \text{\# of injective }k\text{-tuples in }S^{k}\text{ that end with
}s_{i}\right) \\
&  =\left(  \text{\# of injective }\left(  k-1\right)  \text{-tuples in
}\left(  S\setminus\left\{  s_{i}\right\}  \right)  ^{k-1}\right)  .
\end{align*}

However, recall our induction hypothesis (\ref{pf.thm.count.select.injtup.IH}%
). We have $\left\vert S\right\vert =n$ (since $S$ is an $n$-element set).
Since $s_{i}$ is an element of $S$, the set $\left\{  s_{i}\right\}  $ is a
subset of $S$. Thus, the difference rule (Theorem \ref{thm.count.rules.diff}
\textbf{(b)}) yields%
\[
\left\vert S\setminus\left\{  s_{i}\right\}  \right\vert
=\underbrace{\left\vert S\right\vert }_{=n}-\underbrace{\left\vert \left\{
s_{i}\right\}  \right\vert }_{=1}=n-1,
\]
so that $S\setminus\left\{  s_{i}\right\}  $ is an $\left(  n-1\right)
$-element set, and we have $n-1=\left\vert S\setminus\left\{  s_{i}\right\}
\right\vert \in\mathbb{N}$. Hence, we can apply
(\ref{pf.thm.count.select.injtup.IH}) to $n-1$ and $S\setminus\left\{
s_{i}\right\}  $ instead of $n$ and $S$ (because
(\ref{pf.thm.count.select.injtup.IH}) is a \textquotedblleft for all
$n\in\mathbb{N}$\textquotedblright\ statement, not just a statement about the
specific $n$ that we have fixed right now!). As a result, we obtain%
\begin{align*}
&  \left(  \text{\# of injective }\left(  k-1\right)  \text{-tuples in
}\left(  S\setminus\left\{  s_{i}\right\}  \right)  ^{k-1}\right) \\
&  =\left(  n-1\right)  \underbrace{\left(  \left(  n-1\right)  -1\right)
}_{=n-2}\underbrace{\left(  \left(  n-1\right)  -2\right)  }_{=n-3}%
\cdots\underbrace{\left(  \left(  n-1\right)  -\left(  k-1\right)  +1\right)
}_{=n-k+1}\\
&  =\left(  n-1\right)  \left(  n-2\right)  \left(  n-3\right)  \cdots\left(
n-k+1\right)  .
\end{align*}
Combining what we have found, we obtain%
\begin{align*}
&  \left(  \text{\# of injective }k\text{-tuples in }S^{k}\text{ that end with
}s_{i}\right) \\
&  =\left(  \text{\# of injective }\left(  k-1\right)  \text{-tuples in
}\left(  S\setminus\left\{  s_{i}\right\}  \right)  ^{k-1}\right) \\
&  =\left(  n-1\right)  \left(  n-2\right)  \left(  n-3\right)  \cdots\left(
n-k+1\right)  .
\end{align*}

Now, forget that we fixed $i$. We have thus proved that%
\begin{align}
&  \left(  \text{\# of injective }k\text{-tuples in }S^{k}\text{ that end with
}s_{i}\right) \nonumber\\
&  =\left(  n-1\right)  \left(  n-2\right)  \left(  n-3\right)  \cdots\left(
n-k+1\right)  \label{pf.thm.count.select.injtup.si}%
\end{align}
for \textbf{every} $i\in\left[  n\right]  $. Therefore,
(\ref{pf.thm.count.select.injtup.sum}) becomes%
\begin{align*}
&  \left(  \text{\# of injective }k\text{-tuples in }S^{k}\right) \\
&  =\sum_{i=1}^{n}\underbrace{\left(  \text{\# of injective }k\text{-tuples in
}S^{k}\text{ that end with }s_{i}\right)  }_{\substack{=\left(  n-1\right)
\left(  n-2\right)  \left(  n-3\right)  \cdots\left(  n-k+1\right)
\\\text{(by (\ref{pf.thm.count.select.injtup.si}))}}}\\
&  =\sum_{i=1}^{n}\left(  n-1\right)  \left(  n-2\right)  \left(  n-3\right)
\cdots\left(  n-k+1\right) \\
&  =n\cdot\left(  n-1\right)  \left(  n-2\right)  \left(  n-3\right)
\cdots\left(  n-k+1\right) \\
&  \ \ \ \ \ \ \ \ \ \ \ \ \ \ \ \ \ \ \ \ \left(  \text{since }\sum_{i=1}%
^{n}a=na\text{ for any number }a\right) \\
&  =n\left(  n-1\right)  \left(  n-2\right)  \cdots\left(  n-k+1\right)  .
\end{align*}

Forget that we fixed $n$ and $S$. We thus have proved that for all
$n\in\mathbb{N}$ and all $n$-element sets $S$, we have%
\[
\left(  \text{\# of injective }k\text{-tuples in }S^{k}\right)  =n\left(
n-1\right)  \left(  n-2\right)  \cdots\left(  n-k+1\right)  .
\]
In other words, we have proved $P\left(  k\right)  $. Thus, the induction step
is complete, and Theorem \ref{thm.count.select.injtup} is proved.
\end{proof}

\subsubsection{Intermezzo: Listing $n$ elements}

Theorem \ref{thm.count.select.injtup} tells us that if $S$ is an $n$-element
set, then the \# of ways to choose $k$ distinct elements from $S$, if the
order matters, is%
\[
n\left(  n-1\right)  \left(  n-2\right)  \cdots\left(  n-k+1\right)
=k!\cdot\dbinom{n}{k}.
\]

In particular, applying this to $k=n$, we conclude that the \# of ways to
choose $n$ distinct elements from $S$, if the order matters, is%
\[
n\left(  n-1\right)  \left(  n-2\right)  \cdots\left(  n-n+1\right)
=n!\cdot\underbrace{\dbinom{n}{n}}_{\substack{=1\\\text{(by Corollary
\ref{cor.binom.nn})}}}=n!.
\]

Of course, when we are choosing $n$ distinct elements from an $n$-element set,
we are not actually choosing the elements (since all elements have to be
chosen\footnote{This follows from Theorem \ref{thm.count.rules.diff}
\textbf{(c)}.}); we are only choosing the order in which we list them. So what
we have just shown (if somewhat informally) is the following result:

\begin{corollary}
\label{cor.count.select.order-n!}Let $n\in\mathbb{N}$. Let $S$ be an
$n$-element set. Then, the \# of ways to list the $n$ elements of $S$ in some
order (that is, the \# of $n$-tuples that contain each element of $S$ exactly
once) is $n!$.
\end{corollary}

\begin{example}
Applying Corollary \ref{cor.count.select.order-n!} to $n=3$ and $S=\left\{
1,2,3\right\}  $, we see that the \# of ways to list the $3$ numbers $1,2,3$
in some order (i.e., the \# of $3$-tuples that contain each of the numbers
$1,2,3$ exactly once) is $3!=6$. And indeed, here are these $6$ ways:%
\[
\left(  1,2,3\right)  ,\ \ \ \ \ \ \ \ \ \ \left(  1,3,2\right)
,\ \ \ \ \ \ \ \ \ \ \left(  2,1,3\right)  ,\ \ \ \ \ \ \ \ \ \ \left(
2,3,1\right)  ,\ \ \ \ \ \ \ \ \ \ \left(  3,1,2\right)
,\ \ \ \ \ \ \ \ \ \ \left(  3,2,1\right)  .
\]

\end{example}

Corollary \ref{cor.count.select.order-n!} is one of the reasons why factorials
are ubiquitous in combinatorics. The $n!$ ways to list the $n$ elements of a
given $n$-element set $S$ are sometimes called the \textquotedblleft
permutations\textquotedblright\ of $S$, but this name is more frequently used
for the bijective maps from $S$ to $S$. (The \# of the latter maps is also
$n!$, and the two concepts are closely related. For details, see
\cite[\S 1.7.4 in Lecture 13]{22fco}. See also \cite[Lectures 26--28]{22fco}
for much more about permutations.)

\begin{exercise}
\label{exe.simple-counts}\textbf{(a)} How many $7$-digit numbers are there? (A
\textquotedblleft$k$-digit number\textquotedblright\ means a nonnegative
integer that has $k$ digits when written in the decimal system (without
leading zeroes). For example, $3902$ is a $4$-digit number, not a $5$-digit
number.) \medskip

\textbf{(b)} How many $7$-digit numbers are there that have no two equal
digits? \medskip

\textbf{(c)} How many $7$-digit numbers have an even sum of digits? \medskip

\textbf{(d)} How many $7$-digit numbers are palindromes? (A \textquotedblleft%
\textbf{palindrome}\textquotedblright\ is a number such that reading its
digits from right to left yields the same number. For example, $5$ and $1331$
and $49094$ are palindromes.) \medskip

[If your answer is a product or power, \textbf{you do not need to simplify it
to a number}.]
\end{exercise}

\subsubsection{Ordered selections with repetition (= with replacement)}

We have now solved two variants of our \textquotedblleft select $k$ out of
$n$\textquotedblright\ counting question. We have two more variants to go: the
ones where the $k$ elements are arbitrary (not necessarily distinct). Again,
we have the choice of caring or not caring about their order.

If we care about their order, then we are just counting all $k$-tuples in
$S^{k}$. The answer to this question is simple:

\begin{theorem}
\label{thm.count.select.tup}Let $n\in\mathbb{N}$ and $k\in\mathbb{N}$. Let $S$
be an $n$-element set. Then,%
\[
\left(  \text{\# of all }k\text{-tuples in }S^{k}\right)  =n^{k}.
\]

\end{theorem}

\begin{proof}
The set $S$ is an $n$-element set; in other words, $\left\vert S\right\vert
=n$. Now,%
\begin{align*}
&  \left(  \text{\# of all }k\text{-tuples in }S^{k}\right) \\
&  =\left\vert S^{k}\right\vert =\left\vert \underbrace{S\times S\times
\cdots\times S}_{k\text{ times}}\right\vert \ \ \ \ \ \ \ \ \ \ \left(
\text{since }S^{k}\text{ is defined to be }\underbrace{S\times S\times
\cdots\times S}_{k\text{ times}}\right) \\
&  =\underbrace{\left\vert S\right\vert \cdot\left\vert S\right\vert
\cdot\cdots\cdot\left\vert S\right\vert }_{k\text{ times}}%
\ \ \ \ \ \ \ \ \ \ \left(
\begin{array}
[c]{c}%
\text{by the product rule for }k\text{ sets}\\
\text{(Theorem \ref{thm.count.rules.prod-k})}%
\end{array}
\right) \\
&  =\left\vert S\right\vert ^{k}=n^{k}\ \ \ \ \ \ \ \ \ \ \left(  \text{since
}\left\vert S\right\vert =n\right)  .
\end{align*}
This proves Theorem \ref{thm.count.select.tup}.
\end{proof}

\begin{exercise}
\label{exe.23hw5.3}Let $n\in\mathbb{N}$. Compute the \# of $4$-tuples $\left(
a,b,c,d\right)  \in\left[  n\right]  ^{4}$ that satisfy $a\leq b<c\leq d$.
(Not a typo: the second sign is a $<$, not a $\leq$.)

(Recall that $\left[  n\right]  ^{4}=\left[  n\right]  \times\left[  n\right]
\times\left[  n\right]  \times\left[  n\right]  $, so that a $4$-tuple
$\left(  a,b,c,d\right)  \in\left[  n\right]  ^{4}$ means a $4$-tuple of
integers $a,b,c,d\in\left\{  1,2,\ldots,n\right\}  $.)
\end{exercise}

\subsubsection{Unordered selections with repetition (= with replacement)}

Now only one question remains: What is the \# of ways to choose $k$ arbitrary
elements from an $n$-element set $S$ if we \textbf{don't} care about their order?

There are several equivalent ways to rigorously define what this means:

\begin{enumerate}
\item We can define the notion of a \textbf{multiset}, which is
\textquotedblleft like a finite set but allowing an element to be contained
multiple times\textquotedblright. This is done, e.g., in \cite[\S 2.9
(Lectures 21--22)]{22fco} or (in more detail) in \cite[\S 2.11]{19fco}. Then,
a selection of $k$ arbitrary elements from a set $S$, disregarding the order,
can be formalized as a size-$k$ multisubset of the set $S$.

\item Alternatively, we can define the notion of an \textbf{unordered }%
$k$\textbf{-tuple}, which is \textquotedblleft a $k$-tuple up to reordering
its entries\textquotedblright. Formally, these unordered $k$-tuples are
defined as the equivalence classes of usual (i.e., ordered) $k$-tuples with
respect to a certain equivalence relation. (See, e.g., \cite[Example
3.3.24]{19fco} for the details.) Then, a selection of $k$ arbitrary elements
from a set $S$, disregarding the order, can be formalized as an unordered
$k$-tuple of elements of $S$.

\item Finally, if we restrict ourselves to the case when $S=\left[  n\right]
$ (which case is sufficient for all practical purposes, since we can otherwise
rename the elements of $S$ as $1,2,\ldots,n$), then the following
\textquotedblleft low-tech\textquotedblright\ solution becomes available: We
say that a $k$-tuple $\left(  i_{1},i_{2},\ldots,i_{k}\right)  \in S^{k}$ is
\textbf{weakly increasing} (aka \textbf{sorted in weakly increasing order}) if
it satisfies $i_{1}\leq i_{2}\leq\cdots\leq i_{k}$. Now, a selection of $k$
arbitrary elements from $S=\left[  n\right]  $, disregarding the order, can be
defined as a weakly increasing $k$-tuple in $S^{k}$ (because if we don't care
about the order of our $k$ elements, then we can just as well sort them in
increasing order, and the result of such a sorting operation is clearly
unique\footnote{Clearly if you believe in common sense. Not so clearly if you
want a formal proof. See, e.g., \cite[Exercise 2.11.2]{19fco} for such a
proof.}).
\end{enumerate}

These three definitions yield different objects, but these objects are
equivalent, in the sense that there are bijections from each to each. In
particular, the \# of selections of $k$ arbitrary elements from $S$ (without
regard for their order) does not depend on which way we define these
selections. Thus, when it comes to counting them, we can pick whatever
definition we prefer.

\begin{noncompile}
Now, instead of choosing $k$ elements of $S$ up to order, we shall choose a
\textbf{weakly increasing} (aka \textbf{sorted}) $k$-tuple in $S^{k}$: that
is, a $k$-tuple $\left(  i_{1},i_{2},\ldots,i_{k}\right)  \in S^{k}$ such that
$i_{1}\leq i_{2}\leq\cdots\leq i_{k}$. This is equivalent, because if you
don't care about the order of the entries of a $k$-tuple, then you can just
assume that they are arranged in increasing order. Any $k$-tuple of elements
of $S$ can be sorted into weakly increasing order, and two $k$-tuples are
equal up to order if and only if sorting makes them literally equal.
\end{noncompile}

Now that all the requisite warnings and disclaimers have been said, we can
finally count these selections:

\begin{theorem}
\label{thm.count.select.uotup-mulset}Let $n\in\mathbb{N}$ and $k\in\mathbb{N}%
$. Let $S$ be an $n$-element set. Then,%
\begin{align*}
&  \left(  \text{\# of all ways to select }k\text{ elements from }S\text{ (if
order does not matter)}\right) \\
&  =\dbinom{k+n-1}{k}%
\end{align*}
(where our $k$ elements don't have to be distinct).
\end{theorem}

\begin{example}
Applying Theorem \ref{thm.count.select.uotup-mulset} to $n=5$ and $k=2$ and
$S=\left[  5\right]  =\left\{  1,2,3,4,5\right\}  $, we obtain%
\begin{align*}
&  \left(  \text{\# of all ways to select }2\text{ elements from }\left[
5\right]  \text{ (if order does not matter)}\right) \\
&  =\dbinom{2+5-1}{2}=\dbinom{6}{2}=15.
\end{align*}
And indeed, here are these $15$ ways:%
\[%
\begin{array}
[c]{ccccc}%
\left(  1,1\right)  , & \left(  1,2\right)  , & \left(  1,3\right)  , &
\left(  1,4\right)  , & \left(  1,5\right)  ,\\
& \left(  2,2\right)  , & \left(  2,3\right)  , & \left(  2,4\right)  , &
\left(  2,5\right)  ,\\
&  & \left(  3,3\right)  , & \left(  3,4\right)  , & \left(  3,5\right)  ,\\
&  &  & \left(  4,4\right)  , & \left(  4,5\right)  ,\\
&  &  &  & \left(  5,5\right)  .
\end{array}
\]
Here, we have represented each of these selections as a weakly increasing
$k$-tuple in $S^{k}$ (as explained above).
\end{example}

\begin{proof}
[Informal proof of Theorem \ref{thm.count.select.uotup-mulset} (sketched).]For
the sake of simplicity, we assume that $S=\left[  n\right]  $ (since
otherwise, we can rename the $n$ elements of $S$ as $1,2,\ldots,n$). Then, as
we said above, a selection of $k$ arbitrary elements from $S=\left[  n\right]
$ (disregarding the order) can be defined as a weakly increasing $k$-tuple in
$S^{k}$. But a weakly increasing $k$-tuple in $S^{k}$ must always look as
follows:%
\[
\left(  \underbrace{1,1,\ldots,1}_{a_{1}\text{ many }1\text{'s}}%
,\underbrace{2,2,\ldots,2}_{a_{2}\text{ many }2\text{'s}},\ldots
,\underbrace{n,n,\ldots,n}_{a_{n}\text{ many }n\text{'s}}\right)
\]
for some numbers $a_{1},a_{2},\ldots,a_{n}\in\mathbb{N}$ (in particular, each
$a_{i}$ can be $0$, which means that $i$ does not appear in our $k$-tuple)
that satisfy $a_{1}+a_{2}+\cdots+a_{n}=k$ (because we want a $k$-tuple). Such
a $k$-tuple is uniquely determined by these numbers $a_{1},a_{2},\ldots,a_{n}%
$, and conversely, any choice of these numbers $a_{1},a_{2},\ldots,a_{n}$
leads to a different $k$-tuple.

Thus, there is a bijection
\begin{align*}
&  \text{from }\left\{  \text{weakly increasing }k\text{-tuples in }%
S^{k}\right\} \\
&  \text{to }\left\{  n\text{-tuples }\left(  a_{1},a_{2},\ldots,a_{n}\right)
\in\mathbb{N}^{n}\text{ satisfying }a_{1}+a_{2}+\cdots+a_{n}=k\right\}  .
\end{align*}
Hence, the bijection principle yields%
\begin{align*}
&  \left(  \text{\# of weakly increasing }k\text{-tuples in }S^{k}\right) \\
&  =\left(  \text{\# of }n\text{-tuples }\left(  a_{1},a_{2},\ldots
,a_{n}\right)  \in\mathbb{N}^{n}\text{ satisfying }a_{1}+a_{2}+\cdots
+a_{n}=k\right) \\
&  =\left(  \text{\# of weak compositions of }k\text{ into }n\text{
parts}\right) \\
&  \ \ \ \ \ \ \ \ \ \ \ \ \ \ \ \ \ \ \ \ \left(
\begin{array}
[c]{c}%
\text{since the }n\text{-tuples }\left(  a_{1},a_{2},\ldots,a_{n}\right)
\in\mathbb{N}^{n}\\
\text{satisfying }a_{1}+a_{2}+\cdots+a_{n}=k\\
\text{are precisely the weak compositions of }k\text{ into }n\text{ parts}%
\end{array}
\right) \\
&  =\dbinom{k+n-1}{k}\ \ \ \ \ \ \ \ \ \ \left(
\begin{array}
[c]{c}%
\text{by Theorem \ref{thm.count.wcomps.wcomps-nk},}\\
\text{applied to }k\text{ and }n\text{ instead of }n\text{ and }k
\end{array}
\right)  .
\end{align*}
This proves Theorem \ref{thm.count.select.uotup-mulset} (because these weakly
increasing $k$-tuples in $S^{k}$ are the ways to select $k$ elements from $S$
(if order does not matter)).

For a rigorous proof, see \cite[Corollary 2.11.3]{19fco} (but note that the
meanings of the letters $n$ and $k$ are switched in \cite[Corollary
2.11.3]{19fco}).
\end{proof}

Theorem \ref{thm.count.select.uotup-mulset} is our fifth combinatorial
interpretation of binomial coefficients so far! Previously, we have seen that
they count subsets (Theorem \ref{thm.count.subsets.size.S-repeat}), lacunar
subsets (Theorem \ref{thm.count.lacs.size}), compositions (Theorem
\ref{thm.count.comps.comps-nk}) and weak compositions (Theorem
\ref{thm.count.wcomps.wcomps-nk}). This all is not too surprising, since we
proved four of these five theorems using the bijection principle (reducing
them to previously proved theorems), but it is impressive to see so many
counting problems answered by the same family of numbers. \medskip

We have now solved all our four selection problems. We now come to a different
counting problem.

\subsection{Anagrams and multinomial coefficients}

\subsubsection{Counting anagrams}

An \textbf{anagram} of a given word $w$ means a word that consists of the same
letters as $w$ but possibly in a different order. For example:

\begin{itemize}
\item The anagrams of the word \textquotedblleft cat\textquotedblright\ are
\textquotedblleft act\textquotedblright, \textquotedblleft
atc\textquotedblright, \textquotedblleft cat\textquotedblright,
\textquotedblleft cta\textquotedblright, \textquotedblleft
tac\textquotedblright\ and \textquotedblleft tca\textquotedblright.

\item The word \textquotedblleft labl\textquotedblright\ is an anagram of
\textquotedblleft ball\textquotedblright\ (and so are several others).
\end{itemize}

As you see here, we make no distinction between meaningful and meaningless
words. (Also, being logically coherent at the expense of common sense, we
consider each word $w$ to be an anagram of itself.) \medskip

Now, we can take a given word $w$ and ask how many anagrams $w$ has. For instance:

\begin{itemize}
\item How many anagrams does the word \textquotedblleft cat\textquotedblright\ have?

It has six (and we have just listed them above). In fact, we can put the three
letters in any order, and there are $6$ possible orders (by Corollary
\ref{cor.count.select.order-n!}).

\item How many anagrams does the word \textquotedblleft dud\textquotedblright\ have?

It has three (\textquotedblleft dud\textquotedblright, \textquotedblleft
ddu\textquotedblright, \textquotedblleft udd\textquotedblright). Note that the
answer does not directly follow from Corollary \ref{cor.count.select.order-n!}%
, since two of the three letters are equal.

\item How many anagrams does the word \textquotedblleft ball\textquotedblright\ have?

It has $12$ of them: In fact, if the two \textquotedblleft l\textquotedblright
s were two different letters, then it would have $24$ anagrams (again by
Corollary \ref{cor.count.select.order-n!}), but since the two
\textquotedblleft l\textquotedblright s are the same, these $24$ anagrams
merge into pairs of equal words (you get \textquotedblleft
ball\textquotedblright\ twice, you get \textquotedblleft
blal\textquotedblright\ twice, etc.), so the answer is $\dfrac{24}{2}=12$.

(Not convinced? Good; it's worth being skeptical about arguments like this.
Still, this argument can be made precise and rigorous. See \cite[first proof
of Theorem 1.46]{Loehr-BC} for this.)

\item How many anagrams does the word \textquotedblleft
bookkeeper\textquotedblright\ have?

Too many to list by brute force, and the \textquotedblleft divide by
$2$\textquotedblright\ technique from the previous example gets muddled
somewhat as there are several equal letters\footnote{Actually, the technique
can be salvaged, but this requires some carefulness that I am too lazy for
right now. (Once again, see \cite[first proof of Theorem 1.46]{Loehr-BC}.)}.

Thus, let us try a new strategy. The word \textquotedblleft
bookkeeper\textquotedblright\ has $10$ letters. Hence, any anagram of it is a
$10$-letter word as well. Its letters are%
\[
1\text{ \textquotedblleft b\textquotedblright,}\ \ \ \ \ \ \ \ \ \ 3\text{
\textquotedblleft e\textquotedblright s,}\ \ \ \ \ \ \ \ \ \ 2\text{
\textquotedblleft k\textquotedblright s,}\ \ \ \ \ \ \ \ \ \ 2\text{
\textquotedblleft o\textquotedblright s,}\ \ \ \ \ \ \ \ \ \ 1\text{
\textquotedblleft p\textquotedblright}\ \ \ \ \ \ \ \ \ \ \text{and }1\text{
\textquotedblleft r\textquotedblright.}%
\]
In order to choose an anagram of \textquotedblleft
bookkeeper\textquotedblright, we have to distribute all these letters into
$10$ positions. In other words, we have to choose which position the $1$
\textquotedblleft b\textquotedblright\ will occupy, which positions the $3$
\textquotedblleft e\textquotedblright s will occupy, and so on. Let us do this
step by step:

\begin{itemize}
\item We first choose the position of the $1$ \textquotedblleft
b\textquotedblright. There are $\dbinom{10}{1}$ many options for this, since
we need to choose a $1$-element subset of the set of all $10$ positions.

\item We then choose the positions of the $3$ \textquotedblleft
e\textquotedblright s. There are $\dbinom{9}{3}$ many options for this, since
we need to choose a $3$-element subset of the set of all $9$ positions not
already occupied.

\item We then choose the positions of the $2$ \textquotedblleft
k\textquotedblright s. There are $\dbinom{6}{2}$ many options for this, since
we need to choose a $2$-element subset of the set of all $6$ positions not
already occupied.

\item We then choose the positions of the $2$ \textquotedblleft
o\textquotedblright s. There are $\dbinom{4}{2}$ many options for this, since
we need to choose a $2$-element subset of the set of all $4$ positions not
already occupied.

\item We then choose the positions of the $1$ \textquotedblleft
p\textquotedblright. There are $\dbinom{2}{1}$ many options for this, since we
need to choose a $1$-element subset of the set of all $2$ positions not
already occupied.

\item We then choose the positions of the $1$ \textquotedblleft
r\textquotedblright. There are $\dbinom{1}{1}$ many options for this, since we
need to choose a $1$-element subset of the set of all $1$ positions not
already occupied.
\end{itemize}

By the dependent product rule (see the informal proof of Theorem
\ref{thm.count.select.injtup} above), the total \# of ways to perform this
construction is therefore%
\begin{align*}
&  \dbinom{10}{1}\cdot\dbinom{9}{3}\cdot\dbinom{6}{2}\cdot\dbinom{4}{2}%
\cdot\dbinom{2}{1}\cdot\dbinom{1}{1}\\
&  =\dfrac{10!}{1!\cdot9!}\cdot\dfrac{9!}{3!\cdot6!}\cdot\dfrac{6!}{2!\cdot
4!}\cdot\dfrac{4!}{2!\cdot2!}\cdot\dfrac{2!}{1!\cdot1!}\cdot\dfrac{1!}%
{1!\cdot0!}\\
&  \ \ \ \ \ \ \ \ \ \ \ \ \ \ \ \ \ \ \ \ \left(  \text{by the factorial
formula (Theorem \ref{thm.binom.fac})}\right) \\
&  =\dfrac{10!}{1!\cdot3!\cdot2!\cdot2!\cdot1!\cdot1!\cdot0!}%
\ \ \ \ \ \ \ \ \ \ \left(  \text{by cancellations}\right) \\
&  =\dfrac{10!}{1!\cdot3!\cdot2!\cdot2!\cdot1!\cdot1!}%
\ \ \ \ \ \ \ \ \ \ \left(  \text{since }0!=1\right) \\
&  =151\,200.
\end{align*}
Thus, the word \textquotedblleft bookkeeper\textquotedblright\ has
$151\ 200=\dfrac{10!}{1!\cdot3!\cdot2!\cdot2!\cdot1!\cdot1!}$ anagrams.

\item How many anagrams does the word \textquotedblleft
anteater\textquotedblright\ have?

By the same logic as we just used, it has%
\[
\dfrac{8!}{2!\cdot2!\cdot1!\cdot1!\cdot2!}=5\ 040\text{ anagrams.}%
\]

\end{itemize}

The same argument works in the general case:

\begin{theorem}
\label{thm.count.anagrams}Let $s_{1},s_{2},\ldots,s_{n}$ be $n$ distinct
objects, and let $a_{1},a_{2},\ldots,a_{n}$ be $n$ nonnegative integers. Then,
the \# of tuples that consist of%
\begin{align*}
&  a_{1}\text{ copies of }s_{1},\\
&  a_{2}\text{ copies of }s_{2},\\
&  \ldots,\\
&  a_{n}\text{ copies of }s_{n}%
\end{align*}
is%
\[
\dfrac{\left(  a_{1}+a_{2}+\cdots+a_{n}\right)  !}{a_{1}!\cdot a_{2}%
!\cdot\cdots\cdot a_{n}!}=\prod_{k=1}^{n}\dbinom{a_{k}+a_{k+1}+\cdots+a_{n}%
}{a_{k}}.
\]

\end{theorem}

\begin{proof}
[Informal proof (sketched).]Follow the same logic as we used for
\textquotedblleft bookkeeper\textquotedblright\ above. To construct such a
tuple, we

\begin{itemize}
\item first choose the positions for the $a_{1}$ many $s_{1}$'s among its
entries (there are $\dbinom{a_{1}+a_{2}+\cdots+a_{n}}{a_{1}}$ many options for this);

\item then choose the positions for the $a_{2}$ many $s_{2}$'s among its
entries (there are $\dbinom{a_{2}+a_{3}+\cdots+a_{n}}{a_{2}}$ many options for this);

\item then choose the positions for the $a_{3}$ many $s_{3}$'s among its
entries (there are $\dbinom{a_{3}+a_{4}+\cdots+a_{n}}{a_{3}}$ many options for this);

\item and so on, until finally choosing the positions for the $a_{n}$ many
$s_{n}$'s among its entries (there are $\dbinom{a_{n}}{a_{n}}$ many options
for this).
\end{itemize}

By the dependent product rule, the total \# of such tuples is therefore%
\begin{align*}
&  \dbinom{a_{1}+a_{2}+\cdots+a_{n}}{a_{1}}\dbinom{a_{2}+a_{3}+\cdots+a_{n}%
}{a_{2}}\dbinom{a_{3}+a_{4}+\cdots+a_{n}}{a_{3}}\cdots\dbinom{a_{n}}{a_{n}}\\
&  =\prod_{k=1}^{n}\dbinom{a_{k}+a_{k+1}+\cdots+a_{n}}{a_{k}}\\
&  =\prod_{k=1}^{n}\dfrac{\left(  a_{k}+a_{k+1}+\cdots+a_{n}\right)  !}%
{a_{k}!\left(  \left(  a_{k}+a_{k+1}+\cdots+a_{n}\right)  -a_{k}\right)
!}\ \ \ \ \ \ \ \ \ \ \left(
\begin{array}
[c]{c}%
\text{by the factorial formula}\\
\text{(Theorem \ref{thm.binom.fac})}%
\end{array}
\right) \\
&  =\prod_{k=1}^{n}\dfrac{\left(  a_{k}+a_{k+1}+\cdots+a_{n}\right)  !}%
{a_{k}!\left(  a_{k+1}+a_{k+2}+\cdots+a_{n}\right)  !}\\
&  =\dfrac{\prod_{k=1}^{n}\left(  a_{k}+a_{k+1}+\cdots+a_{n}\right)
!}{\left(  \prod_{k=1}^{n}a_{k}!\right)  \left(  \prod_{k=1}^{n}\left(
a_{k+1}+a_{k+2}+\cdots+a_{n}\right)  !\right)  }\\
&  =\dfrac{\left(  a_{1}+a_{2}+\cdots+a_{n}\right)  !\cdot\left(  a_{2}%
+a_{3}+\cdots+a_{n}\right)  !\cdot\cdots\cdot a_{n}!}{\left(  \prod_{k=1}%
^{n}a_{k}!\right)  \left(  \left(  a_{2}+a_{3}+\cdots+a_{n}\right)
!\cdot\left(  a_{3}+a_{4}+\cdots+a_{n}\right)  !\cdot\cdots\cdot a_{n}%
!\cdot0!\right)  }\\
&  =\dfrac{\left(  a_{1}+a_{2}+\cdots+a_{n}\right)  !}{\left(  \prod_{k=1}%
^{n}a_{k}!\right)  \cdot0!}\ \ \ \ \ \ \ \ \ \ \left(
\begin{array}
[c]{c}%
\text{here, we have cancelled factors that appear}\\
\text{both in the numerator and the denominator}%
\end{array}
\right) \\
&  =\dfrac{\left(  a_{1}+a_{2}+\cdots+a_{n}\right)  !}{\prod_{k=1}^{n}a_{k}%
!}\ \ \ \ \ \ \ \ \ \ \left(  \text{since }0!=1\right) \\
&  =\dfrac{\left(  a_{1}+a_{2}+\cdots+a_{n}\right)  !}{a_{1}!\cdot a_{2}%
!\cdot\cdots\cdot a_{n}!}.
\end{align*}
This proves Theorem \ref{thm.count.anagrams}.

(For a rigorous proof, see \cite[Proposition 2.12.13]{19fco}. Note that the
objects $s_{1},s_{2},\ldots,s_{n}$ are required to be $1,2,\ldots,n$ in
\cite[Proposition 2.12.13]{19fco}, but this makes no serious difference, since
we can always rename them at will.)
\end{proof}

\begin{remark}
\label{rmk.enum.primefac}We can now answer the question \textquotedblleft how
many prime factorizations does a given number have?\textquotedblright\ from
Subsection \ref{subsec.ent.primefac}. For example, consider the number
$600=2^{3}\cdot3\cdot5^{2}$. A prime factorization of $600$ is a tuple that
consists of three $2$'s, one $3$ and two $5$'s, in an arbitrary order. Thus,
the \# of such prime factorizations is $\dfrac{6!}{3!\cdot1!\cdot2!}$ (by
Theorem \ref{thm.count.anagrams}). Similarly, we can proceed for any positive
integer instead of $600$. For example, $200=2^{3}\cdot5^{2}$ has $\dfrac
{5!}{3!\cdot2!}=10$ prime factorizations, whereas $150=2\cdot3\cdot5^{2}$ has
$\dfrac{4!}{1!\cdot1!\cdot2!}=12$ prime factorizations.
\end{remark}

\subsubsection{Multinomial coefficients}

The number%
\[
\dfrac{\left(  a_{1}+a_{2}+\cdots+a_{n}\right)  !}{a_{1}!\cdot a_{2}%
!\cdot\cdots\cdot a_{n}!}%
\]
in Theorem \ref{thm.count.anagrams} has a name: It is called a
\textbf{multinomial coefficient}. By Theorem \ref{thm.count.anagrams}, it is
an integer (since it counts something), and can be rewritten as $\prod
_{k=1}^{n}\dbinom{a_{k}+a_{k+1}+\cdots+a_{n}}{a_{k}}$. Note that for $n=2$, it
becomes a binomial coefficient:%
\[
\dfrac{\left(  a+b\right)  !}{a!\cdot b!}=\dbinom{a+b}{a}.
\]

Multinomial coefficients have some further properties. There is a standard
notation for them: Namely, if $a_{1},a_{2},\ldots,a_{n}\in\mathbb{N}$ are any
nonnegative integers, and if we set $b=a_{1}+a_{2}+\cdots+a_{n}$, then the
multinomial coefficient%
\[
\dfrac{\left(  a_{1}+a_{2}+\cdots+a_{n}\right)  !}{a_{1}!\cdot a_{2}%
!\cdot\cdots\cdot a_{n}!}=\dfrac{b!}{a_{1}!\cdot a_{2}!\cdot\cdots\cdot
a_{n}!}%
\]
is denoted by%
\[
\dbinom{b}{a_{1},a_{2},\ldots,a_{n}}.
\]

As already mentioned, multinomial coefficients generalize the binomial
coefficients that are found in Pascal's triangle: With our new notation, a
binomial coefficient $\dbinom{n}{k}$ with $n\in\mathbb{N}$ and $k\in\left\{
0,1,\ldots,n\right\}  $ equals the multinomial coefficient $\dbinom
{n}{k,\ n-k}$. Pascal's identity (Theorem \ref{thm.binom.rec}, at least for
$n>0$ and $k\in\left\{  0,1,\ldots,n\right\}  $) thus can be rewritten as%
\begin{align*}
\dbinom{b}{a_{1},a_{2}}  &  =\dbinom{b-1}{a_{1}-1,\ a_{2}}+\dbinom{b-1}%
{a_{1},\ a_{2}-1}\\
&  \ \ \ \ \ \ \ \ \ \ \text{for }b>0\text{ and }a_{1},a_{2}\in\mathbb{N}%
\text{ with }a_{1}+a_{2}=b,
\end{align*}
where we agree to interpret a multinomial coefficient with a negative number
at the bottom to mean $0$. An analogue of this identity holds for multinomial
coefficients with more parameters:

\begin{theorem}
[Recurrence of the multinomial coefficients]\label{thm.multinom.rec}Let
$b\in\mathbb{N}$ and $a_{1},a_{2},\ldots,a_{n}\in\mathbb{N}$ be such that
$a_{1}+a_{2}+\cdots+a_{n}=b>0$. Then,%
\[
\dbinom{b}{a_{1},a_{2},\ldots,a_{n}}=\sum_{i=1}^{n}\underbrace{\dbinom
{b-1}{a_{1},\ldots,a_{i-1},a_{i}-1,a_{i+1},\ldots,a_{n}}}_{\text{This should
be interpreted as }0\text{ if }a_{i}=0}.
\]

\end{theorem}

\begin{proof}
Nice and fairly easy exercise! (See \cite[Exercise 2.12.6]{19fco} for a proof.)
\end{proof}

Just like the binomial coefficients $\dbinom{n}{k}$ with $n\in\mathbb{N}$ and
$k\in\left\{  0,1,\ldots,n\right\}  $ can be arranged into Pascal's triangle,
the multinomial coefficients $\dbinom{b}{a_{1},a_{2},\ldots,a_{n}}$ (for a
given $n$) can be arranged into an $n$-dimensional analogue of Pascal's
triangle, called \textbf{Pascal's simplex} (or, for $n=3$, \textbf{Pascal's
pyramid}). Theorem \ref{thm.multinom.rec} then says that each entry in this
simplex (except for the $1$ at the apex) is the sum of its $n$ adjacent
entries just above it.

Multinomial coefficients owe their name to another fundamental property they
satisfy: a generalization of the binomial formula, called the
\textbf{multinomial formula}:

\begin{theorem}
[the multinomial formula]\label{thm.multinom.mulfor}Let $x_{1},x_{2}%
,\ldots,x_{n}$ be $n$ numbers. Let $b\in\mathbb{N}$. Then,%
\[
\left(  x_{1}+x_{2}+\cdots+x_{n}\right)  ^{b}=\sum_{\substack{\left(
a_{1},a_{2},\ldots,a_{n}\right)  \in\mathbb{N}^{n};\\a_{1}+a_{2}+\cdots
+a_{n}=b}}\dbinom{b}{a_{1},a_{2},\ldots,a_{n}}x_{1}^{a_{1}}x_{2}^{a_{2}}\cdots
x_{n}^{a_{n}}.
\]

\end{theorem}

\begin{proof}
See \cite[Theorem 2.12.17]{19fco} (which gives two references). Here is the
simplest proof in a nutshell:

We expand $\left(  x_{1}+x_{2}+\cdots+x_{n}\right)  ^{b}$ and collect equal
terms. For instance, if $n=2$ and $b=3$, then%
\begin{align*}
&  \left(  x_{1}+x_{2}+\cdots+x_{n}\right)  ^{b}\\
&  =\left(  x_{1}+x_{2}\right)  ^{3}\\
&  =\left(  x_{1}+x_{2}\right)  \left(  x_{1}+x_{2}\right)  \left(
x_{1}+x_{2}\right) \\
&  =x_{1}x_{1}x_{1}+x_{1}x_{1}x_{2}+x_{1}x_{2}x_{1}+x_{1}x_{2}x_{2}+x_{2}%
x_{1}x_{1}+x_{2}x_{1}x_{2}+x_{2}x_{2}x_{1}+x_{2}x_{2}x_{2}\\
&  =x_{1}^{3}+3x_{1}^{2}x_{2}+3x_{1}x_{2}^{2}+x_{2}^{3}.
\end{align*}

What terms do we get for general $n$ and $b$ ? Well, if we expand the product%
\begin{align*}
&  \left(  x_{1}+x_{2}+\cdots+x_{n}\right)  ^{b}\\
&  =\underbrace{\left(  x_{1}+x_{2}+\cdots+x_{n}\right)  \left(  x_{1}%
+x_{2}+\cdots+x_{n}\right)  \cdots\left(  x_{1}+x_{2}+\cdots+x_{n}\right)
}_{b\text{ times}},
\end{align*}
then we obtain the sum of all $n^{b}$ possible products of the form%
\[
x_{i_{1}}x_{i_{2}}\cdots x_{i_{b}}\text{ with }i_{1},i_{2},\ldots,i_{b}%
\in\left[  n\right]  .
\]
Each such product can be rewritten as the monomial $x_{1}^{a_{1}}x_{2}^{a_{2}%
}\cdots x_{n}^{a_{n}}$, where $a_{1}$ is the \# of $1$'s in the $b$-tuple
$\left(  i_{1},i_{2},\ldots,i_{b}\right)  $, where $a_{2}$ is the \# of $2$'s
in this $b$-tuple, and so on. Moreover, this monomial satisfies $a_{1}%
+a_{2}+\cdots+a_{n}=b$, since the total \# of entries of the $b$-tuple
$\left(  i_{1},i_{2},\ldots,i_{b}\right)  $ is $b$.

Thus, expanding $\left(  x_{1}+x_{2}+\cdots+x_{n}\right)  ^{b}$, we obtain a
sum of monomials of the form $x_{1}^{a_{1}}x_{2}^{a_{2}}\cdots x_{n}^{a_{n}}$
with $a_{1}+a_{2}+\cdots+a_{n}=b$, but each such monomial can appear several
times in this sum. The total \# of copies of a given monomial $x_{1}^{a_{1}%
}x_{2}^{a_{2}}\cdots x_{n}^{a_{n}}$ that appear in this sum equals the \# of
all $b$-tuples that consist of%
\begin{align*}
&  a_{1}\text{ copies of }1,\\
&  a_{2}\text{ copies of }2,\\
&  \ldots,\\
&  a_{n}\text{ copies of }n
\end{align*}
(because of the previous paragraph). But this latter \# equals
\begin{align*}
&  \dfrac{\left(  a_{1}+a_{2}+\cdots+a_{n}\right)  !}{a_{1}!\cdot a_{2}%
!\cdot\cdots\cdot a_{n}!}\ \ \ \ \ \ \ \ \ \ \left(  \text{by Theorem
\ref{thm.count.anagrams}}\right) \\
&  =\dfrac{b!}{a_{1}!\cdot a_{2}!\cdot\cdots\cdot a_{n}!}%
\ \ \ \ \ \ \ \ \ \ \left(
\begin{array}
[c]{c}%
\text{since }a_{1}+a_{2}+\cdots+a_{n}=b\text{ (because}\\
\text{our }b\text{-tuple }\left(  i_{1},i_{2},\ldots,i_{b}\right)  \text{ has
}b\text{ entries in total)}%
\end{array}
\right) \\
&  =\dbinom{b}{a_{1},a_{2},\ldots,a_{n}}\ \ \ \ \ \ \ \ \ \ \left(  \text{by
the definition of }\dbinom{b}{a_{1},a_{2},\ldots,a_{n}}\right)  .
\end{align*}
Thus, each monomial $x_{1}^{a_{1}}x_{2}^{a_{2}}\cdots x_{n}^{a_{n}}$ with
$a_{1}+a_{2}+\cdots+a_{n}=b$ appears exactly $\dbinom{b}{a_{1},a_{2}%
,\ldots,a_{n}}$ times in the sum that we obtain by expanding $\left(
x_{1}+x_{2}+\cdots+x_{n}\right)  ^{b}$. Collecting all copies of each monomial
in this expansion, we thus obtain%
\[
\left(  x_{1}+x_{2}+\cdots+x_{n}\right)  ^{b}=\sum_{\substack{\left(
a_{1},a_{2},\ldots,a_{n}\right)  \in\mathbb{N}^{n};\\a_{1}+a_{2}+\cdots
+a_{n}=b}}\dbinom{b}{a_{1},a_{2},\ldots,a_{n}}x_{1}^{a_{1}}x_{2}^{a_{2}}\cdots
x_{n}^{a_{n}}.
\]
This proves Theorem \ref{thm.multinom.mulfor}.
\end{proof}

We note that this yields a new proof of the binomial formula (Theorem
\ref{thm.binom.bin}), since the latter formula is the particular case of
Theorem \ref{thm.multinom.mulfor} for $n=2$.

\begin{fineprint}
\begin{remark}
We note that Theorem \ref{thm.multinom.rec} can be used to give a second proof
of Theorem \ref{thm.count.anagrams}. Here is a rough outline of this proof:

A tuple that consists of%
\begin{align*}
&  a_{1}\text{ copies of }s_{1},\\
&  a_{2}\text{ copies of }s_{2},\\
&  \ldots,\\
&  a_{n}\text{ copies of }s_{n}%
\end{align*}
will be called an $\left(
\begin{array}
[c]{cccc}%
s_{1} & s_{2} & \cdots & s_{n}\\
a_{1} & a_{2} & \cdots & a_{n}%
\end{array}
\right)  $\textbf{-tuple}. Thus, Theorem \ref{thm.count.anagrams} is claiming
that the \# of $\left(
\begin{array}
[c]{cccc}%
s_{1} & s_{2} & \cdots & s_{n}\\
a_{1} & a_{2} & \cdots & a_{n}%
\end{array}
\right)  $-tuples is $\dbinom{b}{a_{1},a_{2},\ldots,a_{n}}$, where
$b:=a_{1}+a_{2}+\cdots+a_{n}$. We shall now prove this by induction on $b$.
The base case ($b=0$) is trivial (since $b=0$ entails $a_{1}=a_{2}%
=\cdots=a_{n}=0$, so we are counting $0$-tuples). In the induction step (from
$b-1$ to $b$), we separate the $\left(
\begin{array}
[c]{cccc}%
s_{1} & s_{2} & \cdots & s_{n}\\
a_{1} & a_{2} & \cdots & a_{n}%
\end{array}
\right)  $-tuples according to their last entry (just as in our above rigorous
proof of Theorem \ref{thm.count.select.injtup}). This last entry is either
$s_{1}$ or $s_{2}$ or $\cdots$ or $s_{n}$. Hence, the sum rule yields%
\begin{align*}
&  \left(  \text{\# of }\left(
\begin{array}
[c]{cccc}%
s_{1} & s_{2} & \cdots & s_{n}\\
a_{1} & a_{2} & \cdots & a_{n}%
\end{array}
\right)  \text{-tuples}\right) \\
&  =\sum_{i=1}^{n}\underbrace{\left(  \text{\# of }\left(
\begin{array}
[c]{cccc}%
s_{1} & s_{2} & \cdots & s_{n}\\
a_{1} & a_{2} & \cdots & a_{n}%
\end{array}
\right)  \text{-tuples that end with }s_{i}\right)  }_{\substack{=\left(
\text{\# of }\left(
\begin{array}
[c]{cccccccc}%
s_{1} & s_{2} & \cdots & s_{i-1} & s_{i} & s_{i+1} & \cdots & s_{n}\\
a_{1} & a_{2} & \cdots & a_{i-1} & a_{i}-1 & a_{i+1} & \cdots & a_{n}%
\end{array}
\right)  \text{-tuples}\right)  \\\text{(by a bijection argument, just as in
the proof of Theorem \ref{thm.count.select.injtup},}\\\text{using the
bijection that removes the last entry from a tuple)}}}\\
&  =\sum_{i=1}^{n}\underbrace{\left(  \text{\# of }\left(
\begin{array}
[c]{cccccccc}%
s_{1} & s_{2} & \cdots & s_{i-1} & s_{i} & s_{i+1} & \cdots & s_{n}\\
a_{1} & a_{2} & \cdots & a_{i-1} & a_{i}-1 & a_{i+1} & \cdots & a_{n}%
\end{array}
\right)  \text{-tuples}\right)  }_{\substack{=\dbinom{b-1}{a_{1}%
,\ldots,a_{i-1},a_{i}-1,a_{i+1},\ldots,a_{n}}\\\text{(by the induction
hypothesis if }a_{i}>0\text{, and for obvious reasons if }a_{i}=0\text{)}}}\\
&  =\sum_{i=1}^{n}\dbinom{b-1}{a_{1},\ldots,a_{i-1},a_{i}-1,a_{i+1}%
,\ldots,a_{n}}\\
&  =\dbinom{b}{a_{1},a_{2},\ldots,a_{n}}\ \ \ \ \ \ \ \ \ \ \left(  \text{by
Theorem \ref{thm.multinom.rec}}\right)  ,
\end{align*}
which completes the induction step. This proof is less conceptual than the
proof we sketched above, but it is easier to formalize, since it does not use
the dependent product rule.
\end{remark}
\end{fineprint}

\subsection{More counting problems}

Recall the notion of a left inverse, as defined in Exercise \ref{exe.hw5.2}.

\begin{exercise}
\label{exe.23mt3.2}Let $n,m\in\mathbb{N}$. Let $X$ be an $n$-element set. Let
$Y$ be an $m$-element set. Let $f:X\rightarrow Y$ be an injective map. Prove
that $f$ has exactly $n^{m-n}$ many left inverses.
\end{exercise}

If $S$ is any set, and $n$ is any nonnegative integer, then the Cartesian
product $\underbrace{S\times S\times\cdots\times S}_{n\text{ times}}$ is
denoted by $S^{n}$. For example, $S^{3}=S\times S\times S$.

Recall that a $k$-tuple $\left(  i_{1},i_{2},\ldots,i_{k}\right)  $ is called
\textbf{injective} if its $k$ entries $i_{1},i_{2},\ldots,i_{k}$ are all
distinct (i.e., if $i_{a}\neq i_{b}$ for all $a\neq b$).

\begin{exercise}
\label{exe.23mt3.3}Let $n\in\mathbb{N}$. How many injective $\left(
2n\right)  $-tuples $\left(  i_{1},i_{2},\ldots,i_{2n}\right)  \in\left[
2n\right]  ^{2n}$ are there such that all of the first $n$ entries
$i_{1},i_{2},\ldots,i_{n}$ are even?

(For instance, for $n=2$, there are $4$ such tuples: $\left(  2,4,1,3\right)
$, $\left(  2,4,3,1\right)  $, $\left(  4,2,1,3\right)  $ and $\left(
4,2,3,1\right)  $.)
\end{exercise}

\begin{exercise}
\label{exe.23mt3.4}Let $n\geq2$ be an integer. \medskip

\textbf{(a)} How many injective $n$-tuples $\left(  i_{1},i_{2},\ldots
,i_{n}\right)  \in\left[  n\right]  ^{n}$ begin with the entry $2$ ? \medskip

\textbf{(b)} How many injective $n$-tuples $\left(  i_{1},i_{2},\ldots
,i_{n}\right)  \in\left[  n\right]  ^{n}$ contain the entry $1$ before the
entry $2$ ? (\textquotedblleft Before\textquotedblright\ means
\textquotedblleft somewhere to the left of\textquotedblright, not necessarily
\textquotedblleft immediately before\textquotedblright. For instance, for
$n=4$, the $4$-tuple $\left(  1,3,2,4\right)  $ qualifies, but the $4$-tuple
$\left(  2,3,1,4\right)  $ does not.) \medskip

\textbf{(c)} How many injective $n$-tuples $\left(  i_{1},i_{2},\ldots
,i_{n}\right)  \in\left[  n\right]  ^{n}$ contain the entry $1$ immediately
preceding the entry $2$ ? (Here, $\left(  1,3,2,4\right)  $ no longer
qualifies, but $\left(  4,1,2,3\right)  $ does.)
\end{exercise}

If $h:S\rightarrow S$ is any map from a set to itself, then a \textbf{fixed
point} of $h$ means an element $s\in S$ satisfying $h\left(  s\right)  =s$.
The set of all fixed points of $h$ will be called $\operatorname*{Fix}h$.

\begin{exercise}
\label{exe.23mt3.5}Let $X$ and $Y$ be two finite sets (not necessarily of the
same size).

Let $f:X\rightarrow Y$ and $g:Y\rightarrow X$ be two maps. Prove that%
\[
\left\vert \operatorname*{Fix}\left(  f\circ g\right)  \right\vert =\left\vert
\operatorname*{Fix}\left(  g\circ f\right)  \right\vert .
\]

[\textbf{Hint:} Show that $f\left(  x\right)  \in\operatorname*{Fix}\left(
f\circ g\right)  $ for each $x\in\operatorname*{Fix}\left(  g\circ f\right)
$. Thus, there is a map%
\begin{align*}
f^{\prime}:\operatorname*{Fix}\left(  g\circ f\right)   &  \rightarrow
\operatorname*{Fix}\left(  f\circ g\right)  ,\\
x  &  \mapsto f\left(  x\right)  .
\end{align*}
Construct a similar map $g^{\prime}$ in the opposite direction. Prove that
these two maps $f^{\prime}$ and $g^{\prime}$ are inverse to each other.]
\end{exercise}

Now, recall Exercise \ref{exe.23hw5.5}. In that exercise, we decided to call a
set $S$ of integers \textbf{pseudolacunar} if no two elements $s,t$ of $S$
satisfy $\left\vert s-t\right\vert =2$. We denoted the \# of pseudolacunar
subsets of $\left[  n\right]  $ (for a given $n\in\mathbb{N}$) by $p_{n}$.

Recall also the Fibonacci sequence $\left(  f_{0},f_{1},f_{2},\ldots\right)  $
that we introduced in Definition \ref{def.ind.fib}, and the floor function
introduced in Definition \ref{def.ent.floor}.

\begin{exercise}
\label{exe.mt3.6}Prove that
\[
p_{n}=f_{\left\lfloor \left(  n+1\right)  /2\right\rfloor +2}\cdot
f_{\left\lfloor n/2\right\rfloor +2}\ \ \ \ \ \ \ \ \ \ \text{for each }%
n\geq2.
\]

[\textbf{Hint:} What does the pseudolacunarity of a set $S$ mean for the even
elements of $S$ ? What does it mean for the odd elements of $S$ ?]
\end{exercise}

\begin{exercise}
\label{exe.count.slow-bitstr}Let $n\in\mathbb{N}$. An $n$\textbf{-bitstring}
shall mean an $n$-tuple $\left(  a_{1},a_{2},\ldots,a_{n}\right)  \in\left\{
0,1\right\}  ^{n}$ (that is, an $n$-tuple of $0$'s and $1$'s). The product
rule shows that there are $2^{n}$ many $n$-bitstrings. (For example, $\left(
1,1,0,1\right)  $ is a $4$-bitstring.)

\begin{enumerate}
\item[\textbf{(a)}] An $n$-bitstring $\left(  a_{1},a_{2},\ldots,a_{n}\right)
$ is said to be \textbf{lacunar} if it contains no two consecutive $1$'s (that
is, there exists no $i\in\left\{  1,2,\ldots,n-1\right\}  $ such that
$a_{i}=a_{i+1}=1$). How many lacunar $n$-bitstrings are there?

[\textbf{Example:} The bitstring $\left(  0,1,0,0,1\right)  $ is lacunar, but
the bitstring $\left(  0,0,1,1,0\right)  $ is not.]

\item[\textbf{(b)}] An $n$-bitstring $\left(  a_{1},a_{2},\ldots,a_{n}\right)
$ is said to be \textbf{slow} if it contains no entry that differs from both
its neighbors (i.e., there exists no $i\in\left\{  2,3,\ldots,n-1\right\}  $
such that $a_{i}$ is distinct from both $a_{i-1}$ and $a_{i+1}$). How many
slow $n$-bitstrings are there?

[\textbf{Example:} The bitstring $\left(  0,0,1,1,0\right)  $ is slow, but the
bitstring $\left(  0,0,1,0,0\right)  $ is not.]
\end{enumerate}
\end{exercise}

\subsection{The pigeonhole principles}

While studying maps, you might have observed an intuitively obvious fact: A
map $f:X\rightarrow Y$ between two finite sets cannot be injective if
$\left\vert X\right\vert >\left\vert Y\right\vert $ (since there are
\textquotedblleft too many arrows\textquotedblright\ to hit each element of
$Y$ only once), and cannot be surjective if $\left\vert X\right\vert
<\left\vert Y\right\vert $ (since there are \textquotedblleft not enough
arrows\textquotedblright\ to hit each element of $Y$). This is indeed true,
and a bit more can be said:

\begin{theorem}
[pigeonhole principles for maps]\label{thm.count.pigeon}Let $X$ and $Y$ be two
finite sets. Let $f:X\rightarrow Y$ be a map. Then: \medskip

\textbf{(a)} If $\left\vert X\right\vert >\left\vert Y\right\vert $, then $f$
cannot be injective. \medskip

\textbf{(b)} If $f$ is injective and $\left\vert X\right\vert =\left\vert
Y\right\vert $, then $f$ is bijective. \medskip

\textbf{(c)} If $\left\vert X\right\vert <\left\vert Y\right\vert $, then $f$
cannot be surjective. \medskip

\textbf{(d)} If $f$ is surjective and $\left\vert X\right\vert =\left\vert
Y\right\vert $, then $f$ is bijective.
\end{theorem}

Theorem \ref{thm.count.pigeon} is known as the \textbf{pigeonhole principle}
(or \textbf{principles}) because of a traditional way to state it in terms of
pigeons and pigeonholes. For example, part \textbf{(a)} says that if $n$
pigeons are placed in $m$ pigeonholes where $n>m$, then there are (at least)
two pigeons in the same hole. (Here, the pigeons are the elements of $X$, the
pigeonholes are the elements of $Y$, and the assignment of a hole to each
pigeon is the map $f$.) Similarly, the other three parts of Theorem
\ref{thm.count.pigeon} can be reformulated. All parts of Theorem
\ref{thm.count.pigeon} are intuitively obvious\footnote{But beware of
extending your intuition to infinite sets! It is easy to construct an
injective but not surjective map $f:\mathbb{N}\rightarrow\mathbb{N}$.}, but
surprisingly useful (see, e.g., \cite[Worksheet 3]{23f} for multiple
applications.) \medskip

Thus ends our introduction to \emph{enumerative combinatorics} (the study of
finite sets and their sizes, i.e., counting). A further-reaching introduction
can be found in \cite{22fco} (see also \cite{19fco} for a more detailed first
few chapters). Those interested in a deeper immersion should look at Loehr's
\cite{Loehr-BC}, B\'{o}na's \cite{Bona07} or Cameron's \cite{Cameron-NOC}, or
some of the other resources listed in
\url{https://math.stackexchange.com/questions/1454339} .

\newpage


















\end{document}